# ELEMENTARY FUZZY MATRIX THEORY AND FUZZY MODELS FOR SOCIAL SCIENTISTS

**W. B. Vasantha Kandasamy**
**Florentin Smarandache**
**K. Ilanthenral**

**2007**

# ELEMENTARY FUZZY MATRIX THEORY AND FUZZY MODELS FOR SOCIAL SCIENTISTS

**W. B. Vasantha Kandasamy**
e-mail: vasanthakandasamy@gmail.com
web: http://mat.iitm.ac.in/~wbv
www.vasantha.net

**Florentin Smarandache**
e-mail: smarand@unm.edu

**K. Ilanthenral**
e-mail: ilanthenral@gmail.com

# CONTENTS

# PREFACE

This book aims to assist social scientists to analyze their problems using fuzzy models. The basic and essential fuzzy matrix theory is given. The book does not promise to give the complete properties of basic fuzzy theory or basic fuzzy matrices. Instead, the authors have only tried to give those essential basically needed to develop the fuzzy model. The authors do not present elaborate mathematical theories to work with fuzzy matrices; instead they have given only the needed properties by way of examples. The authors feel that the book should mainly help social scientists who are interested in finding out ways to emancipate the society. Everything is kept at the simplest level and even difficult definitions have been omitted. Another main feature of this book is the description of each fuzzy model using examples from real-world problems. Further, this book gives lots of references so that the interested reader can make use of them.

This book has two chapters. In Chapter One, basic concepts about fuzzy matrices are introduced. Basic notions of matrices are given in section one in order to make the book self-contained. Section two gives the properties of fuzzy matrices. Since the data need to be transformed into fuzzy models, some elementary properties of graphs are given. Further, this section provides details of how to prepare a linguistic questionnaire to make use of in these fuzzy models when the data related with the problem is unsupervised.

Chapter Two has six sections. Section one deals with basic fuzzy matrix theory and can be used in a simple and effective way for analyzing supervised or unsupervised data. The simple elegant graphs related with this model can be understood even by a layman. The notion of Fuzzy Cognitive Maps (FCMs) model is introduced in the second section. This model is illustrated by a few examples. It can give the hidden pattern of the problem under analysis. The generalization of the FCM models, which are known as Fuzzy Relational Maps (FRMs), come handy when the attributes related with the problem can be divided into two disjoint sets. This model comes handy when the number of attributes under study is large. This is described in section three. This also gives a pair of fixed points or limit cycle which happens to be the hidden pattern of the dynamical system. Bidirectional Associative Memories (BAM) model is described in the fourth section of this chapter. They are time or period dependent and are defined in real intervals. One can make use of them when the change or solution is time-dependent. This is also illustrated using real-world problems. The fifth section deals with Fuzzy Associative Memories (FAM) model and the model comes handy when one wants the gradations of each and every attribute under study. This model is also described and its working is shown through examples. The last section of this chapter deals with the Fuzzy Relational Equations (FRE) model. This model is useful when there are a set of predicted results and the best solution can be constructed very close, or at times, even equal to the predicted results. The working of this model is also given. Thus the book describes simple but powerful and accurate models that can be used by social scientists. We thank Dr. K. Kandasamy and Meena without their unflinching support this book would have never been possible.

W.B.VASANTHA KANDASAMY
FLORENTIN SMARANDACHE
ILANTHENRAL. K

## Chapter One

# BASIC MATRIX THEORY AND FUZZY MATRIX THEORY

This chapter has three sections. In section one; we give some basic matrix theory. Section two recalls some fundamentals of fuzzy matrix theory. Section three gives the use of mean and standard deviation in matrices.

## 1.1 Basic Matrix Theory

In this section we give some basic matrix theory essential to make the book a self contained one. However the book of Paul Horst on Matrix Algebra for social scientists [92] would be a boon to social scientists who wish to make use of matrix theory in their analysis. We give some very basic matrix algebra which is need for the development of fuzzy matrix theory and the related fuzzy model used for the analysis of socio-economic and psychological problems. However these fuzzy models have been used by applied mathematicians, to study social and psychological problems. These models are very much used by doctors, engineers, scientists, industrialists and statisticians. Here we proceed on to give some basic properties of matrix theory.

A matrix is a table of numbers with a finite number of rows and finite number of columns. The following

$$\begin{matrix} 3 & 0 & 1 & 1 \\ 0 & 5 & 8 & 9 \\ 8 & 1 & 3 & 7 \end{matrix}$$

is an example of a matrix with three rows and four columns.

Since data is always put in the table form it is very easy to consider the very table as matrix. This is very much seen when we use fuzzy matrix as the fuzzy models or fuzzy dynamical systems. So we are mainly going to deal in fuzzy model data matrices which are got from feelings, not always concrete numbers.

When we speak of the table we have the rows and columns clearly marked out so the table by removing the lines can become a matrix with rows and columns. The horizontal entries of the table form the rows and the vertical entries forms the columns of a matrix. Since a table by any technical person (an analyst or a statistician) can have only real entries from the set of reals, thus we have any matrix to be a rectangular array of numbers.

We can interchange the rows or columns i.e., the concepts or the entities or the attributes are interchanged.

Let

$$A = \begin{bmatrix} 2 & 4 & 3 & 1 & 0 \\ 5 & 7 & 8 & -9 & 4 \\ 7 & -6 & 5 & 4 & -3 \\ 8 & 2 & 1 & 0 & 1 \\ 10 & 1 & 0 & 4 & -7 \\ -2 & 2 & 5 & 9 & 8 \end{bmatrix}$$

be a 6 × 5 matrix i.e., a matrix with six rows and five columns, suppose we wish to interchange the sixth row with the second row we get the interchanged matrix as A'

$$A' = \begin{bmatrix} 2 & 4 & 3 & 1 & 0 \\ -2 & 2 & 5 & 9 & 8 \\ 7 & -6 & 5 & 4 & -3 \\ 8 & 2 & 1 & 0 & 1 \\ 10 & 1 & 0 & 4 & -7 \\ 5 & 7 & 8 & -9 & 4 \end{bmatrix}.$$

Now if one wants to interchange the first and fourth column one gets

$$A'' = \begin{bmatrix} 1 & 4 & 3 & 2 & 0 \\ 9 & 2 & 5 & -2 & 8 \\ 4 & -6 & 5 & 7 & -3 \\ 0 & 2 & 1 & 8 & 1 \\ 4 & 1 & 0 & 10 & -7 \\ -9 & 7 & 8 & 5 & 4 \end{bmatrix}.$$

Some times the interchange of row or column in a matrix may not be allowed in certain fuzzy models.

If a matrix has only one column but any number of rows then the matrix is called as a column matrix.

$$C = \begin{bmatrix} 9 \\ 12 \\ 10 \\ 3 \\ 7 \\ -2 \end{bmatrix}$$

is a column matrix. This is a special case of a matrix and sometimes known as the column vector. If a matrix has only one row then we call such a matrix to be a row matrix or a row vector it is also a special case of a matrix.

$$R = [8\ 9\ 12\ 14\ -17\ 10\ 1\ -2\ 5],$$

R is a row matrix or a row vector.

If a matrix has more number of rows than columns then we call it as a vertical matrix.

$$X = \begin{bmatrix} 3 & 1 & 2 \\ 0 & 5 & 9 \\ -10 & 2 & 3 \\ 11 & 7 & -8 \\ 6 & 9 & 10 \\ 7 & -1 & 16 \end{bmatrix}$$

is a vertical matrix.

A matrix which has more number of rows than columns will be known as the horizontal matrix.

$$Y = \begin{bmatrix} 3 & 0 & 12 & 7 & -9 & 8 & -10 & 9 \\ 1 & 9 & 10 & -15 & 7 & 1 & 7 & 3 \end{bmatrix}$$

is a horizontal matrix. The number of rows and columns in a matrix is called the order of a matrix.

$$A = \begin{bmatrix} 3 & 0 & -7 & 8 & 9 & 10 & 13 \\ 9 & 8 & 9 & -11 & 6 & 5 & -9 \\ 11 & 0 & 11 & 2 & 1 & 0 & 8 \\ 12 & 3 & -7 & 6 & 5 & -4 & 3 \\ 7 & -5 & 2 & 1 & -2 & 3 & 4 \\ -9 & 6 & -5 & 6 & 7 & 8 & -9 \end{bmatrix}$$

is $6 \times 7$ matrix.

In any matrix

$$A = \begin{bmatrix} 3 & -1 & 6 & 1 & 8 \\ 8 & 9 & 3 & 2 & 1 \\ 4 & 2 & 1 & 4 & 6 \\ 5 & 6 & 5 & 6 & 5 \end{bmatrix}$$

which is a $4 \times 5$ matrix has 20 elements in total, is a rectangular matrix.

A matrix which has the number of columns is equal to the number of rows is known as a square matrix.

The matrix

$$B = \begin{bmatrix} 8 & -4 & 5 & 3 & -7 \\ 0 & 3 & -2 & 1 & 6 \\ 1 & 8 & 5 & 6 & 0 \\ -2 & 0 & 7 & 8 & 1 \\ 6 & 6 & 1 & 0 & -4 \end{bmatrix}$$

is a square matrix.

Let

$$A = \begin{bmatrix} a_{11} & a_{12} & \cdots & a_{1m} \\ a_{21} & a_{22} & \cdots & a_{2m} \\ \vdots & \vdots & & \vdots \\ a_{n1} & a_{n2} & \cdots & a_{nm} \end{bmatrix}$$

is a n × m matrix where the elements $a_{ij}$ are all real numbers $1 \leq i \leq n$ and $1 \leq j \leq m$ .

Consider a 4 × 4 square matrix A,

$$A = \begin{bmatrix} 3 & 1 & 2 & -4 \\ 1 & 7 & 8 & 9 \\ 2 & 8 & 6 & 0 \\ -4 & 9 & 0 & 7 \end{bmatrix}$$

is a symmetric square matrix.

$$A = \begin{bmatrix} a_{11} & a_{12} & a_{13} & a_{14} \\ a_{21} & a_{22} & a_{23} & a_{24} \\ a_{31} & a_{32} & a_{33} & a_{34} \\ a_{41} & a_{42} & a_{43} & a_{44} \end{bmatrix};$$

A is called the symmetric matrix if $a_{12} = a_{21}$, $a_{13} = a_{31}$, $a_{14} = a_{41}$, $a_{23} = a_{32}$, $a_{24} = a_{42}$ and $a_{34} = a_{43}$.

Let

$$A = \begin{bmatrix} 3 & 0 & 0 & 0 & 0 \\ 0 & 8 & 0 & 0 & 0 \\ 0 & 0 & 7 & 0 & 0 \\ 0 & 0 & 0 & 9 & 0 \\ 0 & 0 & 0 & 0 & 11 \end{bmatrix}.$$

We call A to be a diagonal matrix i.e., if it is a square matrix and all elements other than the main diagonal elements are zero. Thus all diagonal matrices are symmetric but a symmetric matrix in general is not a diagonal matrix. A diagonal matrix is the scalar matrix if all diagonal elements are equal.

$$A = \begin{bmatrix} 5 & 0 & 0 & 0 \\ 0 & 5 & 0 & 0 \\ 0 & 0 & 5 & 0 \\ 0 & 0 & 0 & 5 \end{bmatrix}$$

is a scalar diagonal matrix.
The matrix

$$I_n = \begin{bmatrix} 1 & 0 & 0 & 0 & 0 \\ 0 & 1 & 0 & 0 & 0 \\ 0 & 0 & 1 & 0 & 0 \\ 0 & 0 & 0 & 1 & 0 \\ 0 & 0 & 0 & 0 & 1 \end{bmatrix}$$

is a scalar diagonal matrix which is defined as the identity matrix that is a scalar diagonal matrix is defined to be an identity matrix if all of its diagonal terms are equal and is equal to 1.

A matrix is called a zero or null matrix if every entry of it is zero i.e.,

$$A = \begin{bmatrix} 0 & 0 & 0 & 0 \\ 0 & 0 & 0 & 0 \\ 0 & 0 & 0 & 0 \\ 0 & 0 & 0 & 0 \end{bmatrix}$$

is the zero(null ) matrix. A is a 5 × 4 null matrix.

A column vector or a row vector is said to be a unit vector of every entry in it is one. i.e., P = (1 1 1 1 1 1) is a row vector which is a unit row vector. Similarly

$$Q = \begin{bmatrix} 1 \\ 1 \\ 1 \\ 1 \\ 1 \end{bmatrix}$$

is a column vector which is a unit vector. The column or row vector which has its entries to be zero or ones is defined to be a zero-one vector. T = [0 1 0 1 1 0 0 0 1] is a row zero one vector.

$$S = \begin{bmatrix} 1 \\ 0 \\ 1 \\ 1 \\ 1 \\ 0 \end{bmatrix}$$

is a column zero one vector. Clearly a column or a row vector which has all its entries to be zero is called as a zero vector.

$$[0\ \ 0\ \ 0\ \ 0\ \ 0\ \ 0\ \ 0] \text{ and } \begin{bmatrix} 0 \\ 0 \\ 0 \\ 0 \\ 0 \\ 0 \\ 0 \end{bmatrix}$$

are zero row vector and zero column vector. Two matrices

$$A = \begin{bmatrix} a_{11} & a_{12} & a_{13} & a_{14} \\ a_{21} & a_{22} & a_{23} & a_{24} \\ a_{31} & a_{32} & a_{33} & a_{34} \end{bmatrix} \text{ and } B = \begin{bmatrix} b_{11} & b_{12} & b_{13} & b_{14} \\ b_{21} & b_{22} & b_{23} & b_{24} \\ b_{31} & b_{32} & b_{33} & b_{34} \end{bmatrix}$$

are said to be equal if and only if both A and B (they) are of same order and $a_{11}$ is equal to $b_{11}$, $a_{12}$ equal to $b_{12}$, $a_{13}$ equal to $b_{13}$, $a_{14}$ equal to $b_{14}$, $a_{21}$ equal to $b_{21}$, $a_{22}$ equal to $b_{22}$ and so on.

Thus two matrices A and B of different orders can never be equal.

$$A = \begin{bmatrix} 3 & 1 & 0 \\ 1 & 4 & 5 \\ 8 & 9 & 7 \\ 6 & 3 & 2 \end{bmatrix} \text{ and } B = \begin{bmatrix} 3 & 4 & 9 & 1 & 9 \\ 1 & 5 & 2 & 2 & 7 \\ 2 & 7 & 3 & 6 & 8 \end{bmatrix}$$

are never equal, for order A is 4 × 3 where as order of B is 3 × 5 matrix. A is said to be triangular if

$$A = \begin{bmatrix} 9 & 0 & 0 & 0 \\ 6 & 2 & 0 & 0 \\ 7 & 5 & 3 & 0 \\ 1 & 4 & 2 & 1 \end{bmatrix}.$$

A matrix

$$B = \begin{bmatrix} 9 & 0 & 0 \\ 2 & 4 & 0 \\ 2 & 1 & 9 \\ 5 & 2 & 1 \\ 6 & 1 & 5 \end{bmatrix}$$

is defined to be a partial triangular matrix.

Now we proceed on to define the notion of a transpose of a matrix which we will be using in most of our fuzzy models. Consider the matrix

$$A = \begin{bmatrix} 3 & 1 & 8 & 9 \\ 4 & 0 & 9 & -1 \\ 5 & 2 & 7 & 2 \end{bmatrix}.$$

The transpose of A is

$$A^T = \begin{bmatrix} 3 & 4 & 5 \\ 1 & 0 & 2 \\ 8 & 9 & 7 \\ 9 & -1 & 2 \end{bmatrix}.$$

We denote the transpose of A by $A^T$ or by A'.

So if

$$A = \begin{bmatrix} a_{11} & a_{12} & \cdots & a_{1n} \\ a_{21} & a_{22} & \cdots & a_{2n} \\ \vdots & \vdots & & \vdots \\ a_{m1} & a_{m2} & \cdots & a_{mn} \end{bmatrix}$$

transpose of A denoted by

$$A^T = A' = \begin{bmatrix} a_{11} & a_{21} & \cdots & a_{m1} \\ a_{12} & a_{22} & \cdots & a_{m2} \\ \vdots & \vdots & & \vdots \\ a_{1n} & a_{2n} & \cdots & a_{nm} \end{bmatrix}$$

Clearly the transpose of an identity matrix is itself ie. $A^T$ = A' = A. The transpose of a symmetric matrix is also itself A' = A.

If

$$A = \begin{bmatrix} 3 & 1 & 0 & 1 & 5 \\ 1 & 7 & 8 & 9 & 1 \\ 0 & 8 & 3 & 2 & 4 \\ 1 & 9 & 2 & 9 & 6 \\ 5 & 1 & 4 & 6 & 1 \end{bmatrix}$$

be the symmetric matrix we see $A^T = A$.

Also the transpose of a diagonal matrix is itself. For if

$$A = \begin{bmatrix} 3 & 0 & 0 & 0 & 0 \\ 0 & 2 & 0 & 0 & 0 \\ 0 & 0 & 1 & 0 & 0 \\ 0 & 0 & 0 & -1 & 0 \\ 0 & 0 & 0 & 0 & -2 \end{bmatrix}$$

we get $A^T = A$.

The transpose of a row vector is a column vector and the transpose of a column vector is a row vector.

For if A = [0 1 2 3 5] be a row vector then

$$A' = \begin{bmatrix} 0 \\ 1 \\ 2 \\ 3 \\ 5 \end{bmatrix}$$

which is a column vector.

Suppose

$$B = \begin{bmatrix} 7 \\ 5 \\ 0 \\ 1 \\ 2 \\ 6 \end{bmatrix}$$

be the column vector then $B^T$ = [7 5 0 1 2 6], is the transpose of B which is a row vector. Two matrices of some order alone can be added i.e., if

$$A = \begin{bmatrix} 3 & 4 & 5 \\ 0 & 1 & 2 \\ -1 & 4 & 0 \\ 7 & 0 & 6 \\ 6 & 8 & 9 \end{bmatrix}$$

be a 5 × 3 matrix and B be another 5 × 3 matrix given as

$$B = \begin{bmatrix} 9 & 2 & 1 \\ -1 & 8 & 6 \\ -7 & 0 & 1 \\ 0 & 8 & -2 \\ 2 & -4 & 0 \end{bmatrix}.$$

Then we can add A with B or B with A.
As

$$A + B = \begin{bmatrix} 3 & 4 & 5 \\ 0 & 1 & 2 \\ -1 & 4 & 0 \\ 7 & 0 & 6 \\ 6 & 8 & 9 \end{bmatrix} + \begin{bmatrix} 9 & 2 & 1 \\ -1 & 8 & 6 \\ -7 & 0 & 1 \\ 0 & 8 & -2 \\ 2 & -4 & 0 \end{bmatrix}$$

$$= \begin{bmatrix} 3+9 & 4+2 & 5+1 \\ 0+(-1) & 1+8 & 2+6 \\ -1+(-7) & 4+0 & 0+1 \\ 7+0 & 0+8 & 6+(-2) \\ 6+2 & 8+(-4) & 9+0 \end{bmatrix} = \begin{bmatrix} 12 & 6 & 6 \\ -1 & 9 & 8 \\ -8 & 4 & 1 \\ 7 & 8 & 4 \\ 8 & 4 & 9 \end{bmatrix}.$$

A 2 × 4 matrix A can be added with another 2 × 4 matrix B, where

$$A = \begin{bmatrix} 2 & 0 & 4 & -1 \\ 8 & 2 & 0 & 5 \end{bmatrix} \text{ and } B = \begin{bmatrix} -2 & 1 & -3 & 4 \\ -7 & 5 & 6 & -3 \end{bmatrix}.$$

Now

$$A + B = \begin{bmatrix} 2 & 0 & 4 & -1 \\ 8 & 2 & 0 & 5 \end{bmatrix} + \begin{bmatrix} -2 & 1 & -3 & 4 \\ -7 & 5 & 6 & -3 \end{bmatrix}$$

$$= \begin{bmatrix} 2+(-2) & 0+1 & 4+(-3) & (-1)+4 \\ 8+(-7) & 2+5 & 0+6 & 5+(-3) \end{bmatrix} = \begin{bmatrix} 0 & 1 & 1 & 3 \\ 1 & 7 & 6 & 2 \end{bmatrix}.$$

We can also add a 3 × 3 diagonal matrix D with a 3 × 3 matrix A where

$$D = \begin{bmatrix} 7 & 0 & 0 \\ 0 & 41 & 0 \\ 0 & 0 & 2 \end{bmatrix} \text{ and } A = \begin{bmatrix} 7 & 1 & 2 \\ 5 & -40 & 20 \\ 4 & 7 & 8 \end{bmatrix}$$

$$D + A = \begin{bmatrix} 7 & 0 & 0 \\ 0 & 41 & 0 \\ 0 & 0 & 2 \end{bmatrix} + \begin{bmatrix} 7 & 1 & 2 \\ 5 & -40 & 20 \\ 4 & 7 & 8 \end{bmatrix}$$

$$= \begin{bmatrix} 7+7 & 0+1 & 0+2 \\ 0+5 & 41+(-40) & 0+20 \\ 0+4 & 0+7 & 2+8 \end{bmatrix} = \begin{bmatrix} 14 & 1 & 2 \\ 5 & 1 & 20 \\ 4 & 7 & 10 \end{bmatrix}.$$

Also a column vector

$$A = \begin{bmatrix} 3 \\ 1 \\ 4 \\ 5 \end{bmatrix}$$

of order 4 × 1 can be added to a column vector B of 4 × 1 order where

$$B = \begin{bmatrix} -7 \\ -2 \\ 0 \\ 4 \end{bmatrix}.$$

$$A + B = \begin{bmatrix} 3 \\ 1 \\ 4 \\ 5 \end{bmatrix} + \begin{bmatrix} -7 \\ -2 \\ 0 \\ 4 \end{bmatrix} = \begin{bmatrix} 3+(-7) \\ 1+(-2) \\ 4+0 \\ 5+4 \end{bmatrix} = \begin{bmatrix} -4 \\ -1 \\ 4 \\ 9 \end{bmatrix}.$$

We say addition is commutative if a + b = b + a. We see addition of numbers is always commutative. We can see even the addition of matrices A with B is commutative i.e., A + B = B + A. It is left as a very simple exercise for the reader to prove matrix addition is commutative. We can also add A + B + C where all the 3 matrices A, B and C are matrices of same order. For take

$$A = \begin{bmatrix} 7 & 1 & 2 & 0 \\ 0 & 3 & 4 & 5 \\ 1 & 2 & 3 & 1 \end{bmatrix}, B = \begin{bmatrix} 1 & 4 & 1 & 0 \\ 0 & 1 & -1 & 2 \\ 7 & -2 & 2 & -1 \end{bmatrix}$$

and

$$C = \begin{bmatrix} -6 & 1 & 0 & 1 \\ 2 & -3 & 1 & 2 \\ -7 & 2 & 0 & 1 \end{bmatrix}$$

to be 3 × 4 matrices we can add

$$(A + B) + C = \left[\begin{bmatrix} 7 & 1 & 2 & 0 \\ 0 & 3 & 4 & 5 \\ 1 & 2 & 3 & 1 \end{bmatrix} + \begin{bmatrix} 1 & 4 & 1 & 0 \\ 0 & 1 & -1 & 2 \\ 7 & -2 & 2 & -1 \end{bmatrix}\right] + \begin{bmatrix} -6 & 1 & 0 & 1 \\ 2 & -3 & 1 & 2 \\ -7 & 2 & 0 & 1 \end{bmatrix}$$

$$= \begin{bmatrix} 7+1 & 1+4 & 2+1 & 0+0 \\ 0+0 & 3+1 & 4+(-1) & 5+2 \\ 1+7 & -2+2 & 3+2 & 1+(-1) \end{bmatrix} + \begin{bmatrix} -6 & 1 & 0 & 1 \\ 2 & -3 & 1 & 2 \\ -7 & 2 & 0 & 1 \end{bmatrix}$$

$$= \begin{bmatrix} 8 & 5 & 3 & 0 \\ 0 & 4 & 3 & 7 \\ 8 & 0 & 5 & 0 \end{bmatrix} + \begin{bmatrix} -6 & 1 & 0 & 1 \\ 2 & -3 & 1 & 2 \\ -7 & 2 & 0 & 1 \end{bmatrix}$$

$$= \begin{bmatrix} 8+(-6) & 5+1 & 3+0 & 0+1 \\ 0+2 & 4+(-3) & 1+3 & 7+2 \\ 8+(-7) & 0+2 & 5+0 & 0+1 \end{bmatrix} = \begin{bmatrix} 2 & 6 & 3 & 1 \\ 2 & 1 & 4 & 9 \\ 1 & 2 & 5 & 1 \end{bmatrix}.$$

Now consider A + (B + C)

$$= \begin{bmatrix} 7 & 1 & 2 & 0 \\ 0 & 3 & 4 & 5 \\ 1 & 2 & 3 & 1 \end{bmatrix} + \left\{ \begin{bmatrix} 1 & 4 & 1 & 0 \\ 0 & 1 & -1 & 2 \\ 7 & -2 & 2 & -1 \end{bmatrix} + \begin{bmatrix} -6 & 1 & 0 & 1 \\ 2 & -3 & 1 & 2 \\ -7 & 2 & 0 & 1 \end{bmatrix} \right\}$$

$$= \begin{bmatrix} 7 & 1 & 2 & 0 \\ 0 & 3 & 4 & 5 \\ 1 & 2 & 3 & 1 \end{bmatrix} + \begin{bmatrix} 1-(-6) & 4+1 & 1+0 & 0+1 \\ 0+2 & 1+(-3) & (-1)+1 & 2+2 \\ 7+(-7) & -2+2 & 2+0 & -1+1 \end{bmatrix}$$

$$= \begin{bmatrix} 7 & 1 & 2 & 0 \\ 0 & 3 & 4 & 5 \\ 1 & 2 & 3 & 1 \end{bmatrix} + \begin{bmatrix} -5 & 5 & 1 & 1 \\ 2 & -2 & 0 & 4 \\ 0 & 0 & 2 & 0 \end{bmatrix}$$

$$= \begin{bmatrix} 7+(-5) & 1+5 & 2+1 & 0+1 \\ 0+2 & 3+(-2) & 4+0 & 5+4 \\ 1+0 & 2+0 & 3+2 & 1+0 \end{bmatrix} = \begin{bmatrix} 2 & 6 & 3 & 1 \\ 2 & 1 & 4 & 9 \\ 1 & 2 & 5 & 1 \end{bmatrix}$$

we see (A+B) + C = A + (B+C). We call this rule (a + b) + c = a + (b + c) as the associative law. The real number under addition satisfy the associative law. Once again we leave it as a simple exercise for the reader to verify that matrix addition follows the associative law. We next show the subtraction of elementary matrices. Just like matrix addition we see matrix subtraction is also compatible only when the two matrices are of same order. Let

$$A = \begin{bmatrix} 3 & 1 & 2 & 3 \\ 5 & 2 & -1 & 0 \\ 3 & 2 & 4 & 6 \\ 0 & 2 & 1 & 7 \\ 1 & 0 & 3 & 1 \end{bmatrix}$$

be a 5 × 4 matrix.

$$B = \begin{bmatrix} -2 & 0 & 1 & 5 \\ 7 & 2 & -7 & 0 \\ 6 & 0 & 2 & 1 \\ -1 & 4 & 5 & 2 \\ 0 & 3 & 1 & 4 \end{bmatrix}$$

be another matrix of 5 × 4 order.

$$A - B = \begin{bmatrix} 3 & 1 & 2 & 3 \\ 5 & 2 & -1 & 0 \\ 3 & 2 & 4 & 6 \\ 0 & 2 & 1 & 7 \\ 1 & 0 & 3 & 1 \end{bmatrix} - \begin{bmatrix} -2 & 0 & 1 & 5 \\ 7 & 2 & -7 & 0 \\ 6 & 0 & 2 & 1 \\ -1 & 4 & 5 & 2 \\ 0 & 3 & 1 & 4 \end{bmatrix}$$

$$= \begin{bmatrix} 3-(-2) & 1-0 & 2-1 & 3-5 \\ 5-7 & 2-2 & -1-(-7) & 0-0 \\ 3-6 & 2-0 & 4-2 & 6-1 \\ 0-(-1) & 2-4 & 1-5 & 7-2 \\ 1-0 & 0-3 & 3-1 & 1-4 \end{bmatrix}$$

$$= \begin{bmatrix} 5 & 1 & 1 & -2 \\ -2 & 0 & 6 & 0 \\ -3 & 2 & 2 & 5 \\ 1 & -2 & -4 & 5 \\ 1 & -3 & 2 & -3 \end{bmatrix}.$$

Now find

$$B - A = \begin{bmatrix} -2 & 0 & 1 & 5 \\ 7 & 2 & -7 & 0 \\ 6 & 0 & 2 & 1 \\ -1 & 4 & 5 & 2 \\ 0 & 3 & 1 & 4 \end{bmatrix} - \begin{bmatrix} 3 & 1 & 2 & 3 \\ 5 & 2 & -1 & 0 \\ 3 & 2 & 4 & 6 \\ 0 & 2 & 1 & 7 \\ 1 & 0 & 3 & 1 \end{bmatrix}$$

$$= \begin{bmatrix} -2-3 & 0-1 & 1-2 & 5-3 \\ 7-5 & 2-2 & -7-(-1) & 0-0 \\ 6-3 & 0-2 & 2-4 & 1-6 \\ -1-0 & 4-2 & 5-1 & 2-7 \\ 0-1 & 3-0 & 1-3 & 4-1 \end{bmatrix}$$

$$= \begin{bmatrix} -5 & -1 & -1 & 2 \\ 2 & 0 & -6 & 0 \\ 6 & -2 & -2 & -5 \\ -1 & 2 & 4 & -5 \\ -1 & 3 & -2 & 3 \end{bmatrix}.$$

We see A – B ≠ B – A. Thus matrix subtraction is non commutative. It can also be verified by any interested reader that the operation subtraction of matrices in general is not associative i.e., A – (B – C) ≠ (A – B) – C.

Now having seen the simple matrix operations we proceed on to define matrix multiplication we start with the product of two vectors. Before we say any thing, if the product AB is to be defined we need the only condition that the number of columns of A must be equal to the number of rows of B. BA for the matrix B with A may or may not be defined.

Consider the row vector A = [1 0 5 6 1] and the column vector,

$$B = \begin{bmatrix} 1 \\ 2 \\ 0 \\ 1 \\ 5 \end{bmatrix},$$

AB the product of A with B is well defined for A is a 1 × 5 matrix and B is a 5 × 1 matrix. So A B is a 1 × 1 matrix given by

$$AB = [1\ 0\ 5\ 6\ 1]\begin{bmatrix}1\\2\\0\\1\\5\end{bmatrix}$$

$$= (1 \times 1) + (0 \times 2) + (5 \times 0) + (6 \times 1) + (1 \times 5)$$
$$= 1 + 0 + 0 + 6 + 5$$
$$= 12.$$

BA is 5 × 1 matrix with a 1 × 5 matrix given by

$$BA = \begin{bmatrix}1\\2\\0\\1\\5\end{bmatrix}[1\ 0\ 5\ 6\ 1] = \begin{bmatrix}1 & 0 & 5 & 6 & 1\\2 & 0 & 10 & 12 & 2\\0 & 0 & 0 & 0 & 0\\1 & 0 & 5 & 6 & 1\\5 & 0 & 25 & 30 & 5\end{bmatrix}$$

which is a 5 × 5 matrix. Thus in this case we see both AB and BA are defined however we see AB is a singleton i.e., a 1 × 1 matrix where as BA is a 5 × 5 matrix.

Thus AB ≠ BA.

Now we consider yet another example where the product AB alone is defined and BA remains undefined.

Let

$$A = \begin{bmatrix}3 & 1 & 0\\2 & 0 & 1\\5 & 6 & 0\\0 & 1 & 2\\1 & 1 & 1\\0 & 0 & 7\\1 & 1 & 0\end{bmatrix}$$

and

$$B = \begin{bmatrix} 0 & 3 & 6 & 1 \\ 1 & 4 & 1 & 0 \\ 2 & 5 & 0 & 1 \end{bmatrix}$$

We can find AB as A is a 7 × 3 matrix and B is a 3 × 4 matrix.

Thus

$$AB = \begin{bmatrix} 3 & 1 & 0 \\ 2 & 0 & 1 \\ 5 & 6 & 0 \\ 0 & 1 & 2 \\ 1 & 1 & 1 \\ 0 & 0 & 7 \\ 1 & 1 & 0 \end{bmatrix} \begin{bmatrix} 0 & 3 & 6 & 1 \\ 1 & 4 & 1 & 0 \\ 2 & 5 & 0 & 1 \end{bmatrix}$$

$$= \begin{bmatrix} 1 & 13 & 19 & 3 \\ 2 & 11 & 12 & 3 \\ 6 & 39 & 36 & 5 \\ 5 & 14 & 1 & 2 \\ 3 & 12 & 7 & 2 \\ 14 & 35 & 0 & 7 \\ 1 & 7 & 7 & 1 \end{bmatrix}$$

is a 7 × 4 matrix.

Now BA is not even defined as B is a 3 × 4 matrix and A is only a 7 × 3 matrix. Thus we see in case of matrices the product of the matrix A with B may be defined by the product of B with A may not be defined. Further case both AB and BA are defined we may not have AB to be equal to BA.

But it is important and interesting to note that if A is any matrix and $A^T$ the transpose of A then both $AA^T$ and $A^TA$ are well defined. For consider the matrix A = [1 0 1 2 3 4 0 1]. Now

$$A^T = \begin{bmatrix} 1 \\ 0 \\ 1 \\ 2 \\ 3 \\ 4 \\ 0 \\ 1 \end{bmatrix}, \quad AA^T = [1\ 0\ 1\ 2\ 3\ 4\ 0\ 1\ ] \begin{bmatrix} 1 \\ 0 \\ 1 \\ 2 \\ 3 \\ 4 \\ 0 \\ 1 \end{bmatrix}$$

= $(1 \times 1) + (0 \times 0) + (1 \times 1) + (2 \times 2) + (3 \times 3) + (4 \times 4) + (0 \times 0) + (1 \times 1)$

= $1 + 0 + 1 + 4 + 9 + 16 + 0 + 1 = 32.$

$$A^TA = \begin{bmatrix} 1 \\ 0 \\ 1 \\ 2 \\ 3 \\ 4 \\ 0 \\ 1 \end{bmatrix} [1\ 0\ 1\ 2\ 3\ 4\ 0\ 1]$$

$$= \begin{bmatrix} 1 & 0 & 1 & 2 & 3 & 4 & 0 & 1 \\ 0 & 0 & 0 & 0 & 0 & 0 & 0 & 0 \\ 1 & 0 & 1 & 2 & 3 & 4 & 0 & 1 \\ 2 & 0 & 2 & 4 & 6 & 8 & 0 & 2 \\ 3 & 0 & 3 & 6 & 9 & 12 & 0 & 3 \\ 4 & 0 & 4 & 8 & 12 & 16 & 0 & 4 \\ 0 & 0 & 0 & 0 & 0 & 0 & 0 & 0 \\ 1 & 0 & 1 & 2 & 3 & 4 & 0 & 1 \end{bmatrix}.$$

Clearly both $AA^T$ and $A^TA$ are symmetric matrices for $AA^T$ being a singleton is a trivially symmetric matrix $A^TA$ is a $8 \times 8$

symmetric matrix. Now we consider a horizontal row vector A and the product of A with $A^T$, and $A^T$ with A.

Let

$$A = \begin{bmatrix} 3 & 1 & 2 & 1 & 1 & 1 & 0 \\ 0 & 1 & 0 & 1 & 0 & 2 & 1 \\ 1 & 0 & 1 & 0 & 0 & 0 & 2 \\ 1 & 0 & 0 & 1 & 1 & 1 & 0 \end{bmatrix};$$

be a 4 × 7 matrix. Now

$$A^T = \begin{bmatrix} 3 & 0 & 1 & 1 \\ 1 & 1 & 0 & 0 \\ 2 & 0 & 1 & 0 \\ 1 & 1 & 0 & 1 \\ 1 & 0 & 0 & 1 \\ 1 & 2 & 0 & 1 \\ 0 & 1 & 2 & 0 \end{bmatrix}.$$

The product

$$AA^T = \begin{bmatrix} 3 & 1 & 2 & 1 & 1 & 1 & 0 \\ 0 & 1 & 0 & 1 & 0 & 2 & 1 \\ 1 & 0 & 1 & 0 & 0 & 0 & 2 \\ 1 & 0 & 0 & 1 & 1 & 1 & 0 \end{bmatrix} \begin{bmatrix} 3 & 0 & 1 & 1 \\ 1 & 1 & 0 & 0 \\ 2 & 0 & 1 & 0 \\ 1 & 1 & 0 & 1 \\ 1 & 0 & 0 & 1 \\ 1 & 2 & 0 & 1 \\ 0 & 1 & 2 & 0 \end{bmatrix}$$

$$= \begin{bmatrix} 17 & 4 & 5 & 6 \\ 4 & 7 & 2 & 3 \\ 5 & 2 & 6 & 1 \\ 6 & 3 & 1 & 4 \end{bmatrix}$$

is a symmetric 4 × 4 matrix. Now we find

$$A^TA = \begin{bmatrix} 3 & 0 & 1 & 1 \\ 1 & 1 & 0 & 0 \\ 2 & 0 & 1 & 0 \\ 1 & 1 & 0 & 1 \\ 1 & 0 & 0 & 1 \\ 1 & 2 & 0 & 1 \\ 0 & 1 & 2 & 0 \end{bmatrix} \begin{bmatrix} 3 & 1 & 2 & 1 & 1 & 1 & 0 \\ 0 & 1 & 0 & 1 & 0 & 2 & 1 \\ 1 & 0 & 1 & 0 & 0 & 0 & 2 \\ 1 & 0 & 0 & 1 & 1 & 1 & 0 \end{bmatrix}$$

$$= \begin{bmatrix} 11 & 3 & 7 & 4 & 4 & 4 & 2 \\ 3 & 2 & 2 & 2 & 1 & 3 & 1 \\ 7 & 2 & 5 & 2 & 2 & 2 & 2 \\ 4 & 2 & 2 & 3 & 2 & 4 & 1 \\ 4 & 1 & 2 & 2 & 2 & 2 & 0 \\ 4 & 3 & 2 & 4 & 2 & 6 & 2 \\ 2 & 1 & 2 & 1 & 0 & 2 & 5 \end{bmatrix}$$

is a symmetric 7 × 7 matrix.

Further the identity matrix $I_n$ serves as the identity for any product with a compatible multiplication provided both of them are square matrices.

Let

$$A = \begin{bmatrix} 3 & 1 & 5 & 6 & 0 \\ 2 & 0 & 1 & 1 & 7 \\ 5 & 9 & 0 & 8 & 1 \\ 8 & 10 & 7 & 6 & 2 \\ 11 & 6 & 2 & 0 & 1 \end{bmatrix}$$

be a 5 × 5 matrix. Take

$$I_5 = \begin{bmatrix} 1 & 0 & 0 & 0 & 0 \\ 0 & 1 & 0 & 0 & 0 \\ 0 & 0 & 1 & 0 & 0 \\ 0 & 0 & 0 & 1 & 0 \\ 0 & 0 & 0 & 0 & 1 \end{bmatrix}.$$

We see $AI_5 = I_5A = A$;

$$\begin{bmatrix} 3 & 1 & 5 & 6 & 0 \\ 2 & 0 & 1 & 1 & 7 \\ 5 & 9 & 0 & 8 & 1 \\ 8 & 10 & 7 & 6 & 2 \\ 11 & 6 & 2 & 0 & 1 \end{bmatrix} \begin{bmatrix} 1 & 0 & 0 & 0 & 0 \\ 0 & 1 & 0 & 0 & 0 \\ 0 & 0 & 1 & 0 & 0 \\ 0 & 0 & 0 & 1 & 0 \\ 0 & 0 & 0 & 0 & 1 \end{bmatrix}$$

$$= \begin{bmatrix} 1 & 0 & 0 & 0 & 0 \\ 0 & 1 & 0 & 0 & 0 \\ 0 & 0 & 1 & 0 & 0 \\ 0 & 0 & 0 & 1 & 0 \\ 0 & 0 & 0 & 0 & 1 \end{bmatrix} \begin{bmatrix} 3 & 1 & 5 & 6 & 0 \\ 2 & 0 & 1 & 1 & 7 \\ 5 & 9 & 0 & 8 & 1 \\ 8 & 10 & 7 & 6 & 2 \\ 11 & 6 & 2 & 0 & 1 \end{bmatrix}$$

$$= \begin{bmatrix} 3 & 1 & 5 & 6 & 0 \\ 2 & 0 & 1 & 1 & 7 \\ 5 & 9 & 0 & 8 & 1 \\ 8 & 10 & 7 & 6 & 2 \\ 11 & 6 & 2 & 0 & 1 \end{bmatrix}.$$

Also we see $I^T = I$ as we have already mentioned that I is a symmetric matrix. It is to be verified by the reader that in general the product of (AB) C = A (BC). Let

$$A = \begin{bmatrix} 3 & 4 & 1 & 0 \\ 1 & 2 & 4 & 6 \\ 7 & 4 & 1 & 2 \\ 1 & 0 & 1 & 0 \end{bmatrix}, \quad B = \begin{bmatrix} 1 & 0 & 2 & 4 \\ 0 & 1 & 5 & 1 \\ 2 & 1 & 0 & 0 \\ 7 & 0 & 0 & 5 \end{bmatrix}$$

and

$$C = \begin{bmatrix} 0 & 1 & 2 & 0 \\ 0 & 0 & 1 & 0 \\ 1 & 1 & 0 & 0 \\ 0 & 1 & 0 & 0 \end{bmatrix}.$$

Consider (A × B) × C

$$= \left[\begin{bmatrix} 3 & 4 & 1 & 0 \\ 1 & 2 & 4 & 6 \\ 7 & 4 & 1 & 2 \\ 1 & 0 & 1 & 0 \end{bmatrix} \times \begin{bmatrix} 1 & 0 & 2 & 4 \\ 0 & 1 & 5 & 1 \\ 2 & 1 & 0 & 0 \\ 7 & 0 & 0 & 5 \end{bmatrix}\right] \times \begin{bmatrix} 0 & 1 & 2 & 0 \\ 0 & 0 & 1 & 0 \\ 1 & 1 & 0 & 0 \\ 0 & 1 & 0 & 0 \end{bmatrix}$$

$$= \begin{bmatrix} 5 & 5 & 26 & 16 \\ 51 & 6 & 12 & 36 \\ 23 & 5 & 34 & 42 \\ 3 & 1 & 2 & 4 \end{bmatrix} \times \begin{bmatrix} 0 & 1 & 2 & 0 \\ 0 & 0 & 1 & 0 \\ 1 & 1 & 0 & 0 \\ 0 & 1 & 0 & 0 \end{bmatrix} = \begin{bmatrix} 26 & 47 & 15 & 0 \\ 12 & 99 & 108 & 0 \\ 34 & 99 & 51 & 0 \\ 2 & 9 & 7 & 0 \end{bmatrix}$$

Now A × (B × C)

$$= \begin{bmatrix} 3 & 4 & 1 & 0 \\ 1 & 2 & 4 & 6 \\ 7 & 4 & 1 & 2 \\ 1 & 0 & 1 & 0 \end{bmatrix} \times \left[\begin{bmatrix} 1 & 0 & 2 & 4 \\ 0 & 1 & 5 & 1 \\ 2 & 1 & 0 & 0 \\ 7 & 0 & 0 & 5 \end{bmatrix} \times \begin{bmatrix} 0 & 1 & 2 & 0 \\ 0 & 0 & 1 & 0 \\ 1 & 1 & 0 & 0 \\ 0 & 1 & 0 & 0 \end{bmatrix}\right]$$

$$= \begin{bmatrix} 3 & 4 & 1 & 0 \\ 1 & 2 & 4 & 6 \\ 7 & 4 & 1 & 2 \\ 1 & 0 & 1 & 0 \end{bmatrix} \times \begin{bmatrix} 2 & 7 & 2 & 0 \\ 5 & 6 & 1 & 0 \\ 0 & 2 & 5 & 0 \\ 0 & 12 & 14 & 0 \end{bmatrix} = \begin{bmatrix} 26 & 47 & 15 & 0 \\ 12 & 99 & 108 & 0 \\ 34 & 79 & 51 & 0 \\ 2 & 9 & 7 & 0 \end{bmatrix}$$

$$= A \times (B \times C).$$

Clearly (A × B) × C = A × (B × C).

Thus we see the matrix product is associative. Apart from the matrix addition and multiplication one can also work with the operation called max or min which is similar or analogous to addition and yet another type of operation analogous to max {min}. Here also as in case of addition of matrices to be defined we must have the same order like wise we have for min or max rule to be applied to two matrices A and B. We need both A and B must be of the same order.

Let

$$A = \begin{bmatrix} 3 & 9 & 0 & 2 & 1 \\ 5 & 12 & 20 & 9 & 8 \\ 7 & 9 & 7 & 3 & 19 \end{bmatrix}$$

be a 3 × 5 matrix. Let

$$B = \begin{bmatrix} 9 & 0 & 11 & 3 & 1 \\ 2 & 17 & 12 & 0 & 9 \\ 1 & 9 & 14 & 1 & 21 \end{bmatrix}$$

be another 3 × 5 matrix.

Max {A, B} = (max {$a_{ij}$, $b_{ij}$})=

$$\begin{bmatrix} \max(3,9) & \max(9,0) & \max(0,11) & \max(2,3) & \max(1,1) \\ \max(5,2) & \max(12,17) & \max(20,12) & \max(9,0) & \max(8,9) \\ \max(7,1) & \max(9,9) & \max(7,14) & \max(3,1) & \max(19,21) \end{bmatrix}$$

$$= \begin{bmatrix} 9 & 9 & 11 & 3 & 1 \\ 5 & 17 & 20 & 9 & 9 \\ 7 & 9 & 14 & 3 & 21 \end{bmatrix} = \max (B, A).$$

Yet we can find min operation of A and B. Min {A, B} =

$$\begin{bmatrix} \min(3,9) & \min(9,0) & \min(0,11) & \min(2,3) & \min(1,1) \\ \min(5,2) & \min(12,17) & \min(20,12) & \min(9,0) & \min(8,9) \\ \min(7,1) & \min(9,9) & \min(7,14) & \min(3,1) & \min(19,21) \end{bmatrix}$$

$$= \begin{bmatrix} 3 & 0 & 0 & 2 & 1 \\ 2 & 12 & 12 & 0 & 8 \\ 1 & 9 & 7 & 1 & 19 \end{bmatrix} = \min (B, A).$$

Now we just show min or max function like addition of matrices is both commutative and associative.

Let us now define for any matrix A and B max {min ($a_{ij}$, $b_{ij}$)}. For max min function to be defined on two matrices A and B, we need the number of columns of A equal to the number of rows of B; otherwise max {min ($a_{ij}$, $b_{ij}$)} will not be defined.

Let

$$A = \begin{bmatrix} 3 & 7 & 5 & 2 & 0 \\ 1 & 2 & 3 & 4 & 9 \\ 8 & 1 & 0 & 2 & 5 \\ 4 & 0 & 1 & 1 & 2 \end{bmatrix} \text{ and } B = \begin{bmatrix} 6 & 1 & 7 & 1 & 9 \\ 3 & 0 & 2 & 0 & 2 \\ 1 & 3 & 0 & 1 & 3 \\ 2 & 1 & 5 & 1 & 4 \\ 0 & 5 & 7 & 1 & 5 \end{bmatrix}$$

Max {min ($a_{ij}$, $b_{ij}$)} is found as follows. Suppose max {min ($a_{ij}$, $b_{ij}$)}

$$= \begin{bmatrix} C_{11} & C_{12} & C_{13} & C_{14} & C_{15} \\ C_{21} & C_{22} & C_{23} & C_{24} & C_{25} \\ C_{31} & C_{32} & C_{33} & C_{34} & C_{35} \\ C_{41} & C_{42} & C_{43} & C_{44} & C_{45} \end{bmatrix}$$

$C_{11}$ = max {min {3, 6}, min {7, 3}, min {5, 1}, min {2, 2}, min {0, 0}}
= max {3, 3, 1, 2, 0}
= 3.

$C_{12}$ = max {min {3, 1}, min {7, 0}, min {5, 3}, min {2, 1}, min {0, 5}}
= max {1, 0, 3, 1, 0}
= 3.

$C_{13}$ = max {min {3, 7}, min P7, 2}, min {5, 0}, min {2, 5}, min {0, 7}}
= max {3, 2, 0, 2, 0}
= 3.

Thus we get

| | | | | | | | | |
|---|---|---|---|---|---|---|---|---|
| $C_{14}$ | = | 1, | $C_{15}$ | = | 3, | $C_{21}$ | = | 2, |
| $C_{22}$ | = | 5, | $C_{23}$ | = | 7, | $C_{24}$ | = | 1, |
| $C_{25}$ | = | 5. | $C_{31}$ | = | 6, | $C_{32}$ | = | 5, |
| $C_{33}$ | = | 7, | $C_{34}$ | = | 1, | $C_{35}$ | = | 8, |

$C_{41} = 4$, $C_{42} = 2$, $C_{43} = 4$,
$C_{44} = 1$, $C_{45} = 4$.

Thus

$$\max\{\min(a_{ij}, b_{ij})\} = \begin{bmatrix} 3 & 3 & 3 & 1 & 3 \\ 2 & 5 & 7 & 1 & 5 \\ 6 & 5 & 7 & 1 & 8 \\ 4 & 2 & 4 & 1 & 4 \end{bmatrix}.$$

Now suppose we find max min ($b_{ij}$, $a_{ij}$}. We see max {min {$b_{ij}$, $a_{ij}$}) is not defined. Thus we see even if we use max min operation on the matrices A and B and if this is defined in general max min operation of the matrix B and A may not be defined. Now we see this operation will be used in almost all fuzzy models.

Yet one more type of operation which we perform in the fuzzy models is illustrated below.

Let

$$A = \begin{bmatrix} 1 & 4 & 2 & 3 & 0 & 1 \\ 6 & 0 & 1 & 2 & 1 & 0 \\ 3 & 1 & 0 & 5 & 2 & 6 \\ 1 & 2 & 6 & 1 & 0 & 1 \\ 0 & 0 & 1 & 0 & 0 & 4 \\ 1 & 2 & 1 & 0 & 1 & 0 \end{bmatrix}$$

be a 6 × 6 matrix and let X = [0 1 0 0 1 1] be a row vector we calculate

$XA = [7, 2, 3, 2, 2, 4] = X_1$
$X_1A = (a_1, a_2, a_3, a_4, a_5, a_6) = X_2$
$X_2A = X_3$ and so on.

This operation is the usual product max min {X, A} = $x_{11}$, $x_{12}$, $x_{13}$ $x_{14}$, $x_{15}$, $x_{16}$)

= max {min {0, 1}, min {1, 6} min {0, 3} min (0, 1), min (1, 0), min {1, 1}}

= max {0, 1, 0, 0, 0, 1}.
$x_{11} = 1$, $x_{12} = 1$, $x_{13} = 1$, $x_{14} = 1$, $x_{15} = 1$, $x_{16} = 1$,

Thus max min {X, A} = {1 1 1 1 1 1) = $X_1$.

$X_1A = X_2$ and so on. This operation is frequently used in fuzzy models. In the next section we proceed on to describe some of the basic properties of fuzzy matrices.

## 1.2 Basic Concepts on fuzzy matrices

Here for the self containment of this book we recall some of the basic properties about fuzzy matrices and operations using them. Further we recall the fuzzy matrix representation of directed graphs.

Throughout this book [0, 1] denotes the unit interval. We say $x \in [0, 1]$ if $0 \leq x \leq 1$. We also call the unit interval as a fuzzy interval. We say [a, b] is a fuzzy sub interval of the fuzzy interval [0, 1] if $0 \leq a < b \leq 1$: we denote this by $[a, b] \subseteq [0, 1]$. We also use the convention of calling [–1, 1] to be also a fuzzy interval. $x \in [-1, 1]$ if $-1 \leq x \leq 1$. Thus we have $\{x \mid x \in [0, 1]$ i.e., $0 \leq x \leq 1\}$ is uncountable; hence [0, 1] is an infinite set as [0, 1] is an uncountable set. Let X be any universal set. The characteristic function maps elements of X to elements of the set to elements of the set {0, 1}, which is formally expressed by $\chi_A: X \to [0, 1]$. Set A is defined by its characeristic function $\chi_A$. To be more non technical a fuzzy set can be defined mathematically by assigning to each possible individual in the universe of discourse a value representing its grade of membership in the fuzzy set. For Zadeh introduced a theory whose objects fuzzy sets are set with boundaries that are not precise. The membership in a fuzzy set in not a matter of affirmation or denial but rather a matter of a degree. The significance of Zadeh's contribution was that it challenged not only probability theory as a sole agent for uncertainty; but the very foundations upon which the probability theory is based, Aristotelian two – valued logic. For when A is a fuzzy set and x is a relevant object the proposition x is a member of A is not necessarily either true or false as required by two valued logic, but it may be true only to some degree the degree to which x is actually a member of A. It is most common, but not required to express degrees of membership in the fuzzy sets as well as degrees of truth of the associated propositions by numbers in the closed unit interval [0, 1]. The extreme values in this interval 0 and 1, then represent respectively, the total denial and affirmation of the membership in a given fuzzy set as well as the falsity or truth of the associated proposition.

The capability of fuzzy sets to express gradual transitions to membership to non membership and vice versa has a broad utility. This not only helps in the representation of the measurement of uncertainties but also gives a meaningful representation of vague concepts in a simple natural language.

For example a worker wants to find the moods of his master, he cannot say cent percent in mood or 0 percent of mood or depending on the person who is going to meet the boss he can say some 20% in mood or 50% in mood or 1% in mood or 98% in mood. So even 98% in mood or 50% in mood the worker can meet with some confidence. It 1% in mood the worker may not meet his boss. 20% in mood may have fear while meeting him. Thus we see however that this definition eventually leads us to accept all degrees of mood of his boss as in mood no matter how gloomy the boss mood is! In order to resolve this paradox the term mood may introduce vagueness by allowing some sort of gradual transition from degrees of good mood that are considered to be in mood and those that are not. This is infant precisely the basic concept of the fuzzy set, a concept that is both simple and intuitively pleasing and that forms in essence, a generalization of the classical or crisp set. The crisp set is defined in such a way as to dichotomize the individuals in some given universe of discourse into two groups: members (those that certainly belong to the set) and non members (those that certainly do not). Because full membership and full non membership in a fuzzy set can still be indicated by the values of 1 and 0, respectively. The function can be generalized such that the values assigned to the elements of the universal set fall within a specified range and indicate the membership grade of these elements in the set in question. Larger values denote higher degrees of set membership. Such a function is called a membership function and the set defined by it a fuzzy set. As we are not going to work with fuzzy sets or its related graphs we advice the interested reader to refer to [106]. These books supply a lot of information about fuzzy sets. We have just defined fuzzy sets and the need of its in our study. Now we proceed on to define various types of fuzzy matrices without going deep into their structure. Only what is needed for the models are given, for more properties refer [106].

Let X = (0.6, 0.7, 0, 0.3, 1, 0.2, 0.004, 0.0031, 1, 0.102). X is a 1 × 10 row vector. We see every entry in this row vector is from the unit interval [0, 1]. More generally if X = ($x_1$, $x_2$, …, $x_n$), $x_i \in [0, 1]$; $1 \leq i \leq n$ then we call X a 1 × n fuzzy row matrix / vector. Thus A = (0.3, 1, 0) is a fuzzy row vector. Y = (1 1 1 1) is a fuzzy row vector T = (1 0 1 1 0 0) is a fuzzy row vector. Also W = (0.3, 0.2, 0.7, 0.4, 1, 0.1) is a 1 × 6 fuzzy row vector.

Let

$$V = \begin{bmatrix} 0.3 \\ 0.1 \\ 0.201 \\ 0.11 \\ 0.31 \\ 0 \\ 0.12 \end{bmatrix},$$

V is a column matrix or a column vector. But V enjoys an additional property viz. all the entries in V are from the unit interval [0, 1] i.e., from the fuzzy set [0, 1] hence V is a fuzzy column vector / matrix. V is a 7 × 1 fuzzy column matrix. Take

$$W = \begin{bmatrix} 1 \\ 1 \\ 1 \\ 1 \\ 1 \end{bmatrix};$$

W is a fuzzy row vector. If

$$B = \begin{bmatrix} 1 \\ 0 \\ 1 \\ 0 \\ 1 \\ 1 \\ 0 \end{bmatrix}$$

then B is also a fuzzy column vector whose entries are from the fuzzy crisp set {0, 1}. We call

$$A = \begin{bmatrix} 0.3 & 0.1 & 0.2 & 0.4 & 0 \\ 0.6 & 0.01 & 0.5 & 0.8 & 1 \\ 0 & 0.02 & 0.1 & 0.5 & 0.7 \\ 0.4 & 0.12 & 0.7 & 0.91 & 0.6 \end{bmatrix}$$

to be fuzzy matrix of order 4 × 5. We see A is also a 4 × 5 matrix. Thus we can still see that all fuzzy matrices are matrices but every matrix in general is not a fuzzy matrix. For if

$$B = \begin{bmatrix} 3 & 0.6 & 1 & 9 & 12 \\ 4 & 9 & 8 & 2 & 0 \\ 1 & 6 & 7 & 9 & 0.2 \\ 2 & 8 & 5 & 1 & 0.3 \end{bmatrix}$$

B is a 4 × 5 matrix but B is not a fuzzy matrix; We see the fuzzy interval i.e., the unit interval is a subset of reals. Thus a matrix in general is not a fuzzy matrix since the unit interval [0, 1] is contained in the set of reals we see all fuzzy matrices are matrices. The matrix

$$P = \begin{bmatrix} 0 & 1 & 0 & 0 & 1 \\ 1 & 0 & 1 & 1 & 0 \\ 0 & 1 & 0 & 1 & 0 \\ 0 & 0 & 0 & 0 & 1 \end{bmatrix}$$

is also a fuzzy matrix. Now we have seen fuzzy matrices. The big question is can we add two fuzzy matrices A and B and get the sum of them to be fuzzy matrix. The answer in general is not possible for the sum of two fuzzy matrices may turn out to be a matrix which is not a fuzzy matrix.

For consider the fuzzy matrices A and B where A = (0.7, 0.8, 0.9, 1, 0.3, 0.4, 0.1) and B = (0.5, 0.3, 0.2, 0.7, 0.9, .1, 0).

Now A + B where + is the usual addition of number gives A + B = (0.7, 0.8, 0.9, 1, 0.3, 0.4, 0.1) + (0.5, 0.3. 0.2, 0.7, 0.9, 1, 0) = (1.2, 1.1, 1.1, 1.7, 1.2, 1.4, 0.1); we see all the entries in A + B are not in [0, 1] hence A + B is only just a matrix and not a fuzzy matrix.

So only in case of fuzzy matrices the max or min operations are defined. Clearly under the max or min operations the resultant matrix is again a fuzzy matrix which is in someway analogous to our usual addition.

Let

X = (0.3, 0.7, 0.8, 0.9, 0.1) and
Y = (1, 0.2, 0.6, 0.5, 0.4).

Min{X, Y} = {min (0.3, 1), min (0.7, 0.2), min (0.8, 0.6) min (0.9, 0.5), min (0.1, 0.4)

= (0.3, 0.2, 0.6, 0.5, 0.1)
= T.

One can also say infimum of (X, Y) instead of saying as min {X, Y}. We have to make clear for min operation to be defined on the fuzzy matrices, we demand both the matrices must be of the same order.

Next we take the same X = (0.3, 0.7, 0.8, 0.9, 0.1) and Y = (1, 0.2, 0.6, 0.5, 0.4). Now max (X, Y) = {max (0.3, 1), max (0.7, 0.2), max (0.8, 0.6), max (0.9, 0.5), max (0.1, 0.4)} = (1, 0.7, 0.8, 0.9, 0.4) = V.

Clearly min {X, Y} ≠ max {X, Y} i.e., T ≠ V. We see both operations can be used depending on the model and the problem at hand. Like the min operation on the fuzzy matrices for max operation to be defined we needed both the matrices must be of same order.

Let

$$X = \begin{bmatrix} 0.3 & 0.7 & 0.8 & 1 & 0.5 & 0.4 \\ 0.4 & 0.5 & 1 & 0.3 & 0.8 & 0.5 \\ 0.6 & 0.1 & 0.4 & 0.8 & 0 & 0.2 \\ 0.9 & 0.4 & 0.6 & 1 & 0.3 & 0 \end{bmatrix}$$

and

$$Y = \begin{bmatrix} 1 & 0.2 & 0.3 & 0.4 & 0.5 & 0 \\ 0.8 & 0.5 & 0.2 & 0.1 & 0.1 & 1 \\ 0.5 & 1 & 0.8 & 1 & 0 & 0.3 \\ 0.2 & 0.7 & 1 & 0.5 & 0.6 & 0.2 \end{bmatrix}$$

be any two 4 × 6 fuzzy matrices.

To find max {X, Y} and min {X, Y}

Max {X, Y} =

$$\left[\begin{matrix} \max(0.3, 1) & \max(0.7,\ 0.2) & \max(0.8,\ 0.3) \\ \max(0.4,\ 0.8) & \max(0.5,\ 0.5) & \max(1,\ 0.2) \\ \max(0.6,\ 0.5) & \max(0.1, 1) & \max(0.4,\ 0.8) \\ \max(0.9,\ 0.2) & \max(0.4,\ 0.7) & \max(0.6, 1) \end{matrix}\right.$$

$$\left.\begin{matrix} \max(1, 0.4) & \max(0.5, 0.5) & \max(0.4, 0) \\ \max(0.3,\ 0.1) & \max(0.8, 0.1) & \max(0.5, 1) \\ \max(0.8,\ 1) & \max(0, 0) & \max(0.2, 0.3) \\ \max(1,\ 0.5) & \max(0.3, 0.6) & \max(0, 0.2) \end{matrix}\right]$$

$$= \begin{bmatrix} 1 & 0.7 & 0.8 & 1 & 0.5 & 0.4 \\ 0.8 & 0.5 & 1 & 0.3 & 0.8 & 1 \\ 0.6 & 1 & 0.8 & 1 & 0 & 0.3 \\ 0.9 & 0.7 & 1 & 1 & 0.6 & 0.2 \end{bmatrix}$$

$= P = (p_{ij}),\ 1 \le i \le 4,\ 1 \le j \le 6$

Min (X, Y) =

$$\begin{bmatrix} \min(0.3,1) & \min(0.7,0.2) & \min(0.8,0.3) \\ \min(0.4,0.8) & \min(0.5,0.5) & \min(1,0.2) \\ \min(0.6,0.5) & \min(0.1,1) & \min(0.4,0.8) \\ \min(0.9,0.2) & \min(0.4,0.7) & \min(0.6,1) \end{bmatrix.$$

$$\begin{matrix} \min(1,0.4) & \min(0.5,0.5) & \min(0.4,0) \\ \min(0.3,0.1) & \min(0.8,0.1) & \min(0.5,1) \\ \min(0.8,1) & \min(0,0) & \min(0.2,0.3) \\ \min(1,0.5) & \min(0.3,0.6) & \min(0,0.2) \end{bmatrix}$$

$$= \begin{bmatrix} 0.3 & 0.2 & 0.3 & 0.4 & 0.5 & 0 \\ 0.4 & 0.5 & 0.2 & 0.1 & 0.1 & 0.5 \\ 0.5 & 0.1 & 0.4 & 0 & 0 & 0.2 \\ 0.2 & 0.4 & 0.6 & 0.5 & 0.3 & 0 \end{bmatrix}.$$

$= V = (v_{ij})$, $1 \le i \le 4$ and $1 \le j \le 6$.

If $X = (x_{ij})$ and $Y = (y_{ij})$, $1 \le i \le 4$ and $1 \le j \le 6$ we see

$$v_{ij} \le x_{ij} \le p_{ij} \text{ and } v_{ij} \le y_{ij} \le p_{ij}$$

for $1 \le i \le 4$ and $1 \le j \le 6$.

Now when we wish to find the product of two fuzzy matrices X with Y where X and Y are compatible under multiplication; that is the number of column of X equal to the number of row of Y; still we may not have the product XY to be a fuzzy matrix. For take

$$X = \begin{bmatrix} 0.5 & 0 & 0.4 & 1 \\ 0.3 & 0.8 & 1 & 0.9 \\ 0.1 & 1 & 0.9 & 0.2 \\ 1 & 0.7 & 0.5 & 0 \end{bmatrix}$$

and

$$Y = \begin{bmatrix} 0.9 & 0.4 \\ 0.8 & 0.2 \\ 1 & 0.9 \\ 0.7 & 1 \end{bmatrix}$$

be the fuzzy matrices. We see the under the usual matrix product

$$XY = \begin{bmatrix} 0.5 & 0 & 0.4 & 1 \\ 0.3 & 0.8 & 1 & 0.9 \\ 0.1 & 1 & 0.9 & 0.2 \\ 1 & 0.7 & 0.5 & 0 \end{bmatrix} \begin{bmatrix} 0.9 & 0.4 \\ 0.8 & 0.2 \\ 1 & 0.9 \\ 0.7 & 1 \end{bmatrix}$$

$$= \begin{bmatrix} 1.55 & 1.56 \\ 2.54 & 2.08 \\ 1.93 & 1.25 \\ 2.12 & 0.99 \end{bmatrix};$$

XY is not a fuzzy matrix. Thus we see the product of two fuzzy matrices under usual matrix multiplication is not a fuzzy matrix. So we need to define a compatible operation analogous to product so that the product again happens to be a fuzzy matrix. However even for this new operation if the product XY is to be defined we need the number of columns of X is equal to the number of rows of Y. The two types of operations which we can have are max-min operation and min-max operation.

Let

$$X = \begin{bmatrix} 0.3 & 1 & 0.7 & 0.2 & 0.5 \\ 1 & 0.9 & 0 & 0.8 & 0.1 \\ 0.8 & 0.2 & 0.3 & 1 & 0.4 \\ 0.5 & 1 & 0.6 & 0.7 & 0.8 \end{bmatrix}$$

be a $4 \times 5$ fuzzy matrix and let

$$Y = \begin{bmatrix} 0.8 & 0.3 & 1 \\ 0.7 & 0 & 0.2 \\ 1 & 0.7 & 1 \\ 0.5 & 0.4 & 0.5 \\ 0.4 & 0 & 0.7 \end{bmatrix}$$

be a 5 × 3 fuzzy matrix.
XY defined using max. min function.

$$XY = \begin{bmatrix} C_{11} & C_{12} & C_{13} \\ C_{21} & C_{22} & C_{23} \\ C_{31} & C_{32} & C_{33} \\ C_{41} & C_{42} & C_{43} \end{bmatrix}$$

where,

$C_{11}$ = max {min (0.3, 0.8), min (1, 0.7), min (0.7, 1), min (0.2, 0.5), min (0.5, 0.4)}
= max {0.3, 0.7, 0.7, 0.2, 0.4}
= 0.7.

$C_{12}$ = max {min (0.3, 0.3), min (1, 0), min (0.7, 0.7), min (0.2, 0.4), min (0.5, 0)}
= max {0.3, 0, 0.7, 0.2, 0}
= 0.7 and so on.

$$XY = \begin{bmatrix} 0.7 & 0.7 & 0.7 \\ 0.8 & 0.4 & 1 \\ 0.8 & 0.4 & 0.8 \\ 0.7 & 0.6 & 0.7 \end{bmatrix}$$

is a 4 × 3 matrix.

Now suppose for the same X and Y we adopt the operation as min. max operation we get

$$D = \begin{bmatrix} D_{11} & D_{12} & D_{13} \\ D_{21} & D_{22} & D_{23} \\ D_{31} & D_{32} & D_{33} \\ D_{41} & D_{42} & D_{43} \end{bmatrix}$$

$D_{11}$ = min {max (0.3, 0.8), max (1, 0.7), max (0.7, 1) max (0.2, 0.5) max (0.5, 0.4)}
= min {0.8, 1, 1, 0.5, 0.5}
= 0.5.

$D_{12}$ = min {max (0.3, 0.3), max (1, 0), max (0.7, 0.7), max (0.2, 0.4), max (0.5, 0)}
= min {0.3, 1, 0.7, 0.4, 0.5}
= 0.3 and so on .

Thus we have

$$D = \begin{bmatrix} 0.5 & 0.3 & 0.5 \\ 0.4 & 0.1 & 0.7 \\ 0.4 & 0.4 & 0.2 \\ 0.6 & 0.5 & 0.7 \end{bmatrix}.$$

We see C $\neq$ D some experts may wish to work with minimum value and some with the maximum value and accordingly they can adopt for the same.

We need to make the following important observation we see XY is defined but YX may not be defined. In this case we say YX is not defined.

We also find products using max min function as follows. Let

$$A = \begin{bmatrix} 0.3 & 1 & 0.8 & 0.2 & 0 & 1 \\ 0.2 & 0 & 1 & 0.7 & 0.6 & 0.4 \\ 0.5 & 0.3 & 0.2 & 1 & 0 & 0.2 \\ 0.9 & 0.7 & 1 & 0.3 & 0.1 & 0.9 \\ 0.1 & 0.8 & 0 & 0.8 & 1 & 0.7 \\ 0.2 & 1 & 1 & 0.6 & 0 & 0.4 \end{bmatrix}.$$

Let

$$Y = (0.8, 0.1, 0.7, 0, 0.9, 1).$$

Now we can find YA using the max. min function.

$$[0.8 \quad 0.1 \quad 0.7 \quad 0 \quad 0.9 \quad 1]\begin{bmatrix} 0.3 & 1 & 0.8 & 0.2 & 0 & 1 \\ 0.2 & 0 & 1 & 0.7 & 0.6 & 0.4 \\ 0.5 & 0.3 & 0.2 & 1 & 0 & 0.2 \\ 0.9 & 0.7 & 1 & 0.3 & 0.1 & 0.9 \\ 0.1 & 0.8 & 0 & 0.8 & 1 & 0.7 \\ 0.2 & 1 & 1 & 0.6 & 0 & 0.4 \end{bmatrix}$$

= (max{min(0.8, 0.3), min(0.1, 0.2), min(0.7, 0.5), min(0, 0.9), min (0.9, 0.1), min (1, 0.2)}, max {min (0.8, 1), min (0.1, 0), min (0.7, 0.3), min (0, 0.7), min (0.9, 0.8), min (1, 1)} …., max {min (0.8, 1), min (0.1, 0.4), min (0.7, 0.2), min (0, 0.9), min (0.9, 0.7), min (1, 0.4)}).

= (max {0.3, 0.1, 0.5, 0, 0.1, 0.2}, max (0.8, 0, 0.3, 0, 0.9, 1}, max {0.8, 0.1, 0.2, 0, 0, 1}, max {0.2, 0.1, 0.7, 0, 0.8, 0.6}, max (0, 0.1, 0, 0, 0.9, 0}, max {0.8, 0.1, 0.2, 0, 0.7, 0.4})

= (0.5, 1, 1, 0.8, 0.9, 0.8) = $Y_1$.

Now one can calculate $Y_1A$. and so on till the desired number of steps. We give yet another product using max. min operation. Let

$$A = \begin{bmatrix} 0.2 & 1 & 0.5 & 0.5 & 0 \\ 0.9 & 0.8 & 0 & 1 & 0.7 \\ 0.5 & 1 & 0.6 & 0.8 & 1 \\ 0.1 & 0.7 & 0.9 & 0 & 0.5 \end{bmatrix}.$$

Let

$$Y = \begin{bmatrix} 0.7 \\ 0.2 \\ 0.9 \\ 1 \\ 0.4 \end{bmatrix}$$

Now A.Y using is the max. min function (or operation) is as follows.

$$A.Y = \begin{bmatrix} 0.2 & 1 & 0.5 & 0.5 & 0 \\ 0.9 & 0.8 & 0 & 1 & 0.7 \\ 0.5 & 1 & 0.6 & 0.8 & 1 \\ 0.1 & 0.7 & 0.9 & 0 & 0.5 \end{bmatrix} \begin{bmatrix} 0.7 \\ 0.2 \\ 0.9 \\ 1 \\ 0.4 \end{bmatrix}$$

= (max {min (0.2, 0.7), min (1, 0.2), min (0.6, 0.9), min (0.5, 1), min (0, 0.4)}, max {min (0.9, 0.7), min (0.8, 0.2), min (0, 0.9), min (1, 1), min (0.7, 0.4)}, max {min (0.5, 0.7), min (1, 0.2), min (0.6, 0.9), min (0.8, 1), min (1, 0.4)}, max {min (0.1, 0.7), min (0.7, 0.2), min (0.9, 0.9), min (0, 1), min (0.5, 0.4)}

= max (0.2, 0.2, 0.6, 0.5, 0), max (0.7, 0.2, 0.1, 0.4), max (0.5, 0.2, 0.6, 0.8, 0.4), max (0.1, 0.2, 0.9, 0, 0.4))

= (0.6, 1, 0.8, 0.9) = X.

Now we find X.A using max. min operation.

$$(0.6, 1, 0.8, 0.9) \begin{bmatrix} 0.2 & 1 & 0.6 & 0.5 & 0 \\ 0.9 & 0.8 & 0 & 1 & 0.7 \\ 0.5 & 1 & 0.6 & 0.8 & 1 \\ 0.1 & 0.7 & 0.9 & 0 & 0.5 \end{bmatrix}$$

= (max {min (0.6, 0.2), min (1, 0.9), min (0.8, 0.5), min (0.9, 0.1)}, max {min (0.6, 1), min (1, 0.8), min (0.8, 1), min (0.9, 0.7)}, max {min (0.6, 0.6), min (1, 0), min (0.8, 0.6), min (0.9, 0.9)}, max {(0.6, 0.5), min (1, 1), min (0.8, 0.8), min (0.9, 0)}, max {min (0.6, 0), min (1, 0.7), min (0.8, 1), min (0.9, 0.5)})

= (max {0.2, 0.9, 0.5, 0.1}, max {0.6, 0.8, 0.8, 0.7}, max{0.6, 0, 0.6, 0.9}, max {0.5, 1, 0.8, 0}, max {0, 0.7, 0.8, 0.5})

= (0.9, 0.8, 0.9, 1, 0.8) = $Y_1$

If one needs one can find $A.Y_1$ and so on.

Next we show how the product of X.Y where X is a fuzzy row vector and Y a fuzzy column vector is obtained using a max. min product. Let X = (0.9, 1, 0.7, 0, 0.4, 0.6) and

$$Y = \begin{bmatrix} 0.4 \\ 0.7 \\ 1 \\ 0.2 \\ 0 \\ 0.5 \end{bmatrix}.$$

Now we find X.Y using the max min product.

$$X.Y = [0.9, 1, 0.7, 0, 0.4, 0.6] \begin{bmatrix} 0.4 \\ 0.7 \\ 1 \\ 0.2 \\ 0 \\ 0.5 \end{bmatrix}$$

= max {min (0.9, 0.4), min (1, 0.7), min (0.7, 1), min (0, 0.2), min (0.4, 0), min (0.6, 0.5)}

= max {0.4, 0.7, 0.7, 0, 0, 0.5} = 0.7 .

X.Y is a singleton fuzzy matrix.

Now using the same max. min function one can find

$$Y.X = \begin{bmatrix} 0.4 \\ 0.7 \\ 1 \\ 0.2 \\ 0 \\ 0.5 \end{bmatrix} [0.9, 1, 0.7, 0, 0.4, 0.6]$$

$$= \begin{bmatrix} \max\min(0.4,9) & \max\min(0.4,1) & \cdots & \max\min(0.4,0.6) \\ \max\min(0.7,0.9) & \max\min(0.7,1) & \cdots & \max\min(0.7,0.6) \\ \vdots & \vdots & & \vdots \\ \max\min(0.5,0.4) & \max\min(0.5,1) & \cdots & \max\min(0.5,0.6) \end{bmatrix}$$

$$= \begin{bmatrix} 0.4 & 0.4 & 0.4 & 0 & 0.4 & 0.4 \\ 0.7 & 0.7 & 0.7 & 0 & 0.4 & 0.6 \\ 0.9 & 1 & 0.7 & 0 & 0.4 & 0.6 \\ 0.2 & 0.2 & 0.2 & 0 & 0.2 & 0.2 \\ 0 & 0 & 0 & 0 & 0 & 0 \\ 0.5 & 0.5 & 0.5 & 0 & 0.4 & 0.5 \end{bmatrix}.$$

Here also we see XY ≠ YX. Suppose we have a fuzzy square matrix A;

$$A = \begin{bmatrix} -1 & 0 & 1 & 0 & 1 & 1 & 1 \\ 1 & 0 & 1 & 1 & 0 & 0 & 1 \\ 0 & 1 & 0 & 1 & 0 & 1 & 0 \\ -1 & 0 & 1 & 0 & 0 & 0 & 1 \\ 0 & 1 & 0 & 1 & 0 & -1 & 0 \\ 1 & 0 & 1 & 0 & 1 & 0 & 1 \\ 1 & 1 & 0 & 1 & 0 & 1 & 0 \end{bmatrix}$$

Let Y = [1 1 1 1 1 1 1],

$$AY = \begin{bmatrix} -1 & 0 & 1 & 0 & 1 & 1 & 1 \\ 1 & 0 & 1 & 1 & 0 & 0 & 1 \\ 0 & 1 & 0 & 1 & 0 & 1 & 0 \\ -1 & 0 & 1 & 0 & 0 & 0 & 1 \\ 0 & 1 & 0 & 1 & 0 & -1 & 0 \\ 1 & 0 & 1 & 0 & 1 & 0 & 1 \\ 1 & 1 & 0 & 1 & 0 & 1 & 0 \end{bmatrix} \begin{bmatrix} 1 \\ 1 \\ 1 \\ 1 \\ 1 \\ 1 \\ 1 \end{bmatrix} = \begin{bmatrix} 3 \\ 4 \\ 3 \\ 1 \\ 1 \\ 4 \\ 4 \end{bmatrix}.$$

We see that multiplication of A by the column unit matrix Y gives us the row sums of the fuzzy matrix like if we multiply X = [1 1 1 1 1 1 1] the row unit matrix with the fuzzy matrix A will give the column sums.

$$[1\ 1\ 1\ 1\ 1\ 1\ 1]\begin{bmatrix} -1 & 0 & 1 & 0 & 1 & 1 & 1 \\ 1 & 0 & 1 & 1 & 0 & 0 & 1 \\ 0 & 1 & 0 & 1 & 0 & 1 & 0 \\ -1 & 0 & 1 & 0 & 0 & 0 & 1 \\ 0 & 1 & 0 & 1 & 0 & -1 & 0 \\ 1 & 0 & 1 & 0 & 1 & 0 & 1 \\ 1 & 1 & 0 & 1 & 0 & 1 & 0 \end{bmatrix}$$

$$= \quad [1\ 3\ 4\ 4\ 2\ 2\ 4].$$

This mode of getting the row sums and column sums would be used in the fuzzy models in the later chapter.

## 1.3 Basic Concepts on Graphs

In this section we give some of the basic notions of graphs and its relation to matrices and how mean and SD are used in matrices is illustrated. A brief way of obtaining the linguistic questionnaire is also discussed.

We now illustrate how the mean and Standard Deviation are used in matrices.

Mean of n terms $x_1, x_2, \ldots, x_n$ is given by

$$\frac{x_1 + x_2 + \ldots + x_n}{n} = \bar{x}$$

and

Standard Deviation of $x_1, \ldots, x_n$ is given by $\sqrt{\dfrac{\sum (x - \bar{x})^2}{n-1}}$ .

Initial raw data table relating the age groups and symptoms of the cardio vascular problem

| Years/ Symptoms | $S_1$ | $S_2$ | $S_3$ | $S_4$ | $S_5$ | $S_6$ | $S_7$ | $S_8$ |
|---|---|---|---|---|---|---|---|---|
| 20-30 | 23 | 18 | 24 | 16 | 29 | 10 | 16 | 10 |
| 31-43 | 38 | 32 | 38 | 31 | 35 | 18 | 33 | 10 |
| 44-65 | 22 | 21 | 21 | 22 | 20 | 11 | 20 | 5 |

Now this matrix is just the table made into a matrix.

$$\begin{bmatrix} 23 & 18 & 24 & 16 & 29 & 10 & 16 & 10 \\ 38 & 32 & 38 & 31 & 35 & 18 & 33 & 10 \\ 22 & 21 & 21 & 22 & 20 & 11 & 20 & 5 \end{bmatrix}$$

Now this matrix is made uniform as follows:

Every element of the first row is divided by 11. Every element of second row by 13. Every element of third row is divided by 22, which happens to be the length of the class interval of the above table.

This matrix is known as the Average Time Dependent data (ATD) matrix.

The ATD Matrix of Cardio Vascular problem

$$\begin{bmatrix} 2.09 & 1.64 & 2.18 & 1.46 & 2.64 & 0.91 & 1.46 & 0.91 \\ 2.92 & 2.46 & 2.92 & 2.39 & 2.69 & 1.39 & 2.54 & 0.77 \\ 1 & 0.95 & 0.95 & 1 & 0.91 & 0.5 & 0.91 & 0.23 \end{bmatrix}$$

Now for this matrix the average of the 1st column is

$$\frac{2.09+2.92+1}{3} = \frac{6.01}{3} = 2$$

Like wise the SD is given by the formulae $\sqrt{\frac{\sum (x-\bar{x})^2}{n-1}}$ where $\bar{x}$ is the average and the SD of the first column is found to be 0.96.

The Average and Standard Deviation of the above ATD matrix

| | | | | | | | | |
|---|---|---|---|---|---|---|---|---|
| Average | 2.00 | 1.68 | 2.02 | 1.62 | 2.08 | 0.93 | 1.64 | 0.64 |
| SD | 0.96 | 0.76 | 0.995 | 0.71 | 1.01 | 0.45 | 0.83 | 0.36 |

This technique will be used in section 2.1 of chapter 2.

Here we recall some of the basic results from Graph Theory. This is mainly done to make this book a self contained one. We just give only the definition and results on graph theory which we have used. We have taken the results from [84, 161-2, 234]. It is no coincidence that graph theory has been independently discovered many times, since it may quite properly be regarded as an area of applied mathematics.

Euler (1707-1782) became the father of graph theory. In 1847 Kirchoff developed the theory of trees, in order to solve the system of simultaneous linear equations, which give the current in each branch and each circuit of an electric network. Thus in effect Kirchoff replaced each electrical network by its underlying graph and showed that it is not necessary to consider every cycle in the graph of an electric network separately in order to solve the system of equations.

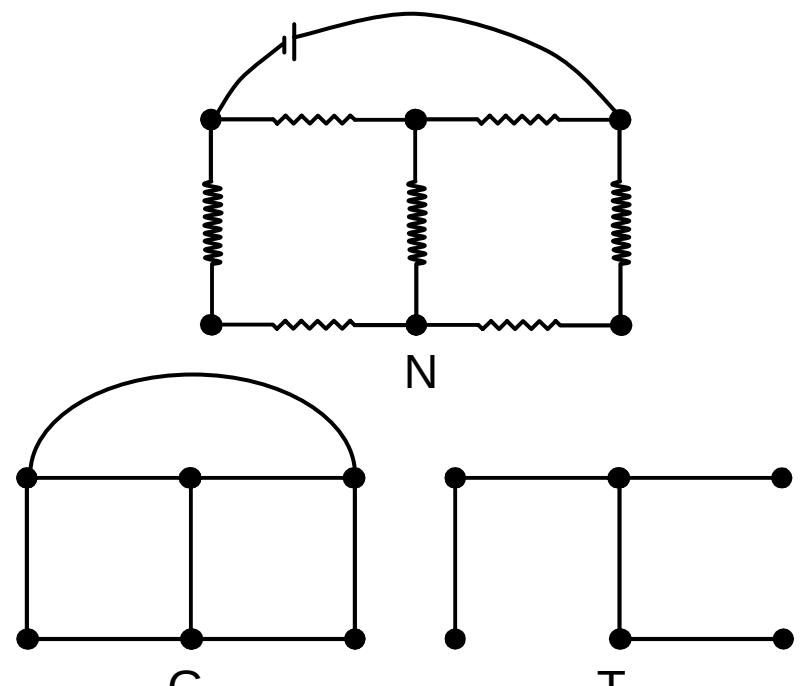

FIGURE: 1.3.1

Instead he pointed out by a simple but powerful construction, which has since become the standard procedure that the independent cycles of a graph determined by any of its spanning trees will suffice. A contrived electrical network N, its under lying graph G and a spanning tree T are shown in the following figure.

In 1857 Cayley discovered the important class of graphs called trees by considering the changes of variables in the differential calculus. Later he was engaged in enumerating the isomers of the saturated hydrocarbons $C_n\ H_{2n+2}$, with a given number n of carbon atoms as shown in the Figure.

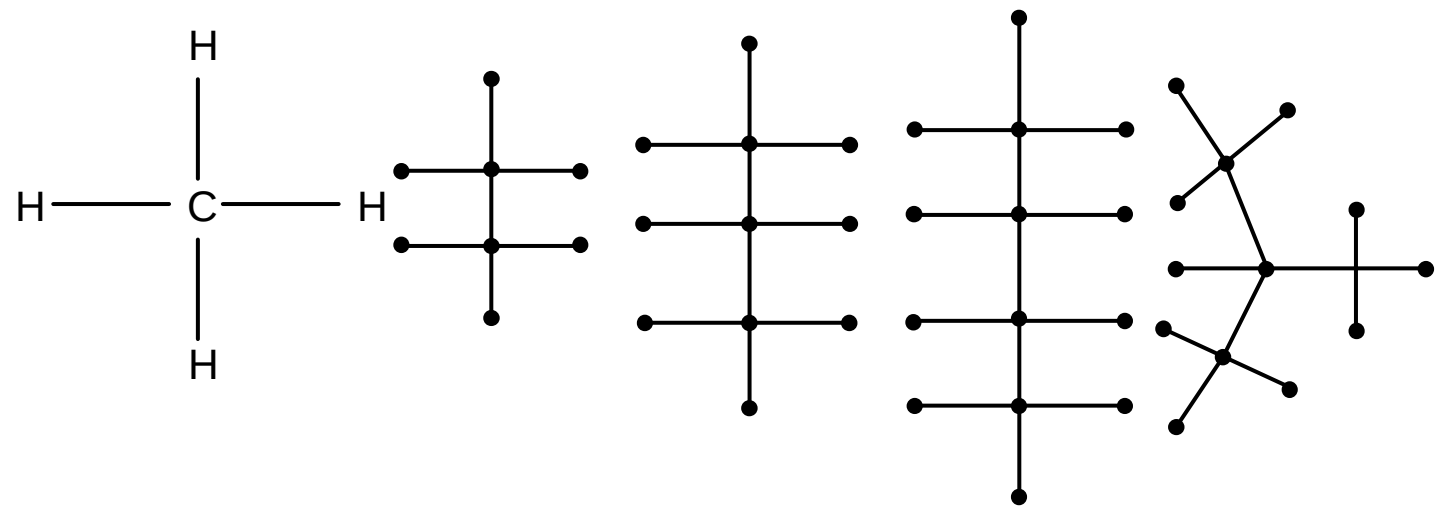

FIGURE: 1.3.2

Cayley restated the problem abstractly, find the number of tress with p points in which every point has degree 1 or 4. He did not immediately succeed in solving this and so he altered the problem until he was able to enumerate; rooted trees (in which one point is distinguished from the others), with point of degree at most 4 and finally the chemical problem of trees in which every point has degree 1 or 4.

Jordan in 1869 independently discovered trees as a purely mathematical discipline and Sylvester 1882 wrote that Jordan did so without having any suspicion of its bearing on modern chemical doctrine.

The most famous problem in graph theory and perhaps in all of mathematics is the celebrated four color conjecture. The remarkable problem can be explained in five minutes by any mathematician to the so called man in the street. At the end of the explanation both will understand the problem but neither will be able to solve it.

The following quotation from the historical article which state the Four color conjecture and describe its role.

Any map on a plane or the surface of a sphere can be colored with only four colors so that no two adjacent countries have the same color.

Each country must consists of a single connected region and adjacent countries are those having a boundary line (not merely a single point) in common.

The conjecture has acted as a catalyst in the branch of mathematics known as combinatorial topology and is closely related to the currently fashionable field of graph theory. Although the computer oriented proof settled the conjecture in 1976 and has stood a test of time, a theoretical proof of the four colour problem is still to be found.

Lewin the psychologist proposed in 1936 that the life span of an individual be represented by a planar map. In such a map the regions would represent the various activities of a person such as his work environment, his home and his hobbies. It was pointed out that Lewin was actually dealing with graphs as indicated by the following figure:

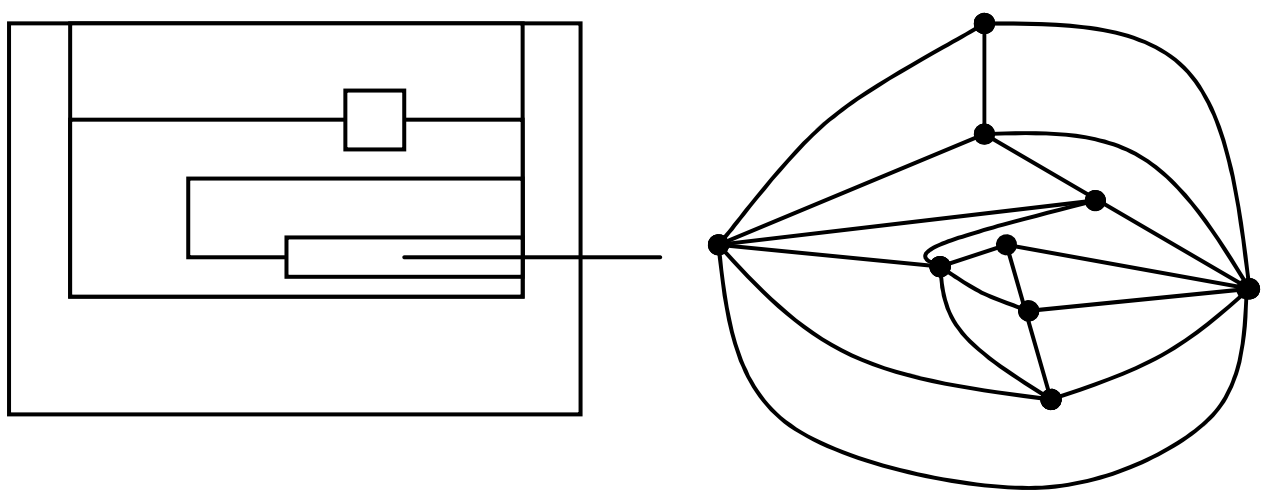

FIGURE: 1.3.3

This viewpoint led the psychologists at the Research center for Group Dynamics to another psychological interpretation of a graph in which people are represented by points and interpersonal relations by lines. Such relations include love, hate, communication and power. In fact it was precisely this approach which led the author to a personal discovery of graph theory, aided and abetted by psychologists L. Festinger and D. Cartwright.

The world of Theoretical physics discovered graph theory for its own purposes more than once. In the study of statistical mechanics by Uhlenbeck the points stand for molecules and two adjacent points indicate nearest neighbor interaction of some physical kind, for example magnetic attraction or repulsion. In a similar interpretation by Lee and Yang the points stand for small cubes in Euclidean space where each cube may or may not be occupied by a molecule.

Then two points are adjacent whenever both spaces are occupied. Another aspect of physics employs graph theory rather as pictorial device. Feynmann proposed the diagram in which the points represent physical particles and the lines represent paths of the particles after collisions.

The study of Markov chains in probability theory involves directed graphs in the sense that events are represented by points and a directed line from one point to another indicates a positive probability of direct succession of these two events. This is made explicit in which Markov chain is defined as a network with the sum of the values of the directed lines from each point equal to 2. A similar representation of a directed graph arises in that part of numerical analysis involving matrix inversion and the calculation of eigen values.

A square matrix is given preferable sparse and a directed graph is associated with it in the following ways. The points denote the index of the rows and columns of the given matrix and there is a directed line from point i to point j whenever the i, j entry of the matrix is nonzero. The similarity between this approach and that for Markov chains in immediate.

Thus finally in the $21^{st}$ century the graph theory has been fully exploited by fuzzy theory. The causal structure of fuzzy cognitive maps from sample data [108-112] mainly uses the notion of fuzzy signed directed graphs with feedback. Thus the use of graph theory especially in the field of applications of fuzzy theory is a grand one for most of analysis of unsupervised data are very successfully carried out by the use of Fuzzy Cognitive Maps (FCMs) which is one of the very few tools which can give the hidden pattern of the dynamical system.

The study of Combined Fuzzy Cognitive Maps (CFCMs) mainly uses the concept of digraphs.

The directed graphs or the diagraphs are used in the representation of Binary relations on a single set. Thus one of the forms of representation of a fuzzy relation R(X, X) is represented by the digraph.

Let X = {1, 2, 3, 4, 5} and R (X, X) the binary relation on X defined by the following membership matrix

| | 1 | 2 | 3 | 4 | 5 |
|---|---|---|---|---|---|
| 1 | 0.2 | 0 | 0 | 0.1 | 0.6 |
| 2 | 0 | 0.8 | 0.4 | 0 | 0 |
| 3 | 0.3 | 0 | 0.9 | 0.2 | 0 |
| 4 | 0 | 0 | 0.2 | 0.9 | 0 |
| 5 | 0 | 0.8 | 0 | 0.5 | 0 |

The related graph for the binary relation.

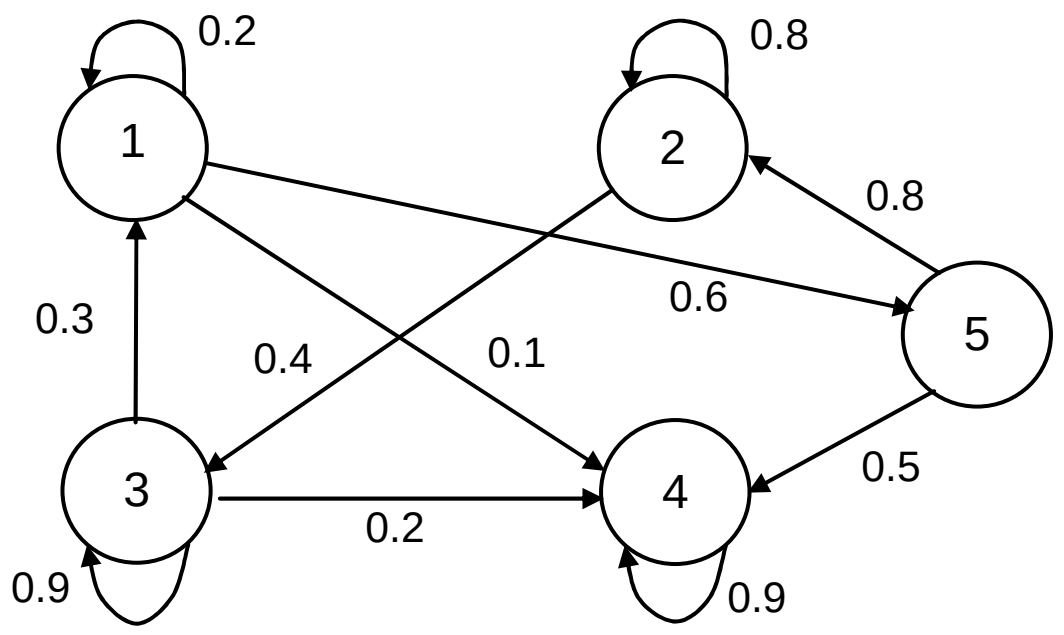

FIGURE: 1.3.4

Further in the description of the fuzzy compatibility relations also the graphs were used.

The notion of graph theory was used in describing the fuzzy ordering relations. In fact one can say that the graph theory method was more simple and an easy representation even by a lay man.

Thus we can say whenever the data had a fuzzy matrix representation it is bound to get the digraph representation. Also in the application side graph theory has been scrupulously used in the description of automaton and semi automaton i.e., in finite machines we do not study in this direction in this book.

Here we recall the definition of graphs and some of its properties. [84, 234]

**DEFINITION 1.3.1:** *A graph G consists of a finite non empty set V = V (G) of p points. (vertex, node, junction O-simplex elements) together with a prescribed set X of q unordered pairs of distinct points of V.*

*Each pair x = {u, v} of points in X is a line (edge, arc, branch, 1-simplex elements) of G and x is said to join u and v.*

*We write x = uv and say that u and v are adjacent points (some times denoted as u adj v); point u and line x are incident with each other as arc v and x. E(G) will denote the edges or lines of G.*

*If two distinct lines x and y are incident with a common point, then they are adjacent lines. A graph with p points and q lines is called a (p, q) graph.*

*Clearly (1, 0) graph is trivial. A graph is represented always by a diagram and we refer to it as the graph.*

The graph in figure 1.3.5 is totally disconnected.

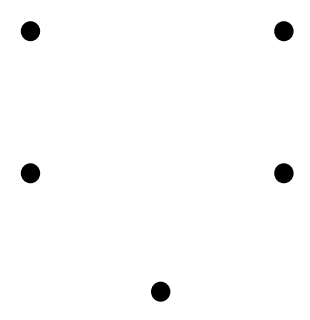

FIGURE: 1.3.5

The graph in figure 1.3.6 is disconnected.

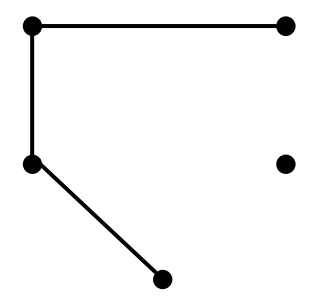

FIGURE: 1.3.6

The graph with four lines in figure 1.3.7 is a path.

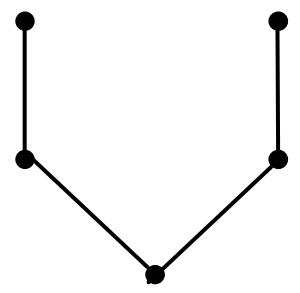

FIGURE: 1.3.7

The graph in figure 1.3.8 is a connected graph.

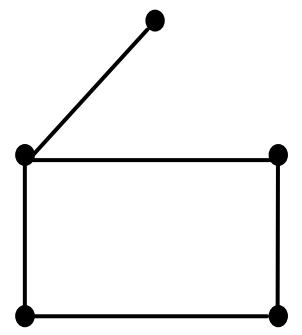

FIGURE: 1.3.8

The graph in figure 1.3.9 is a complete graph.

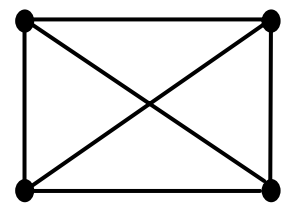

FIGURE: 1.3.9

The graph in figure 1.3.10 with four lines is a cycle.

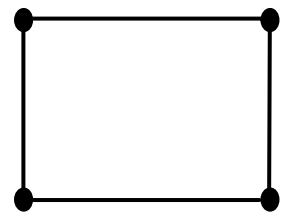

FIGURE:1.3.10

It is important to note that in a graph if any two lines intersect it is not essential that their intersection is a point of the graph.

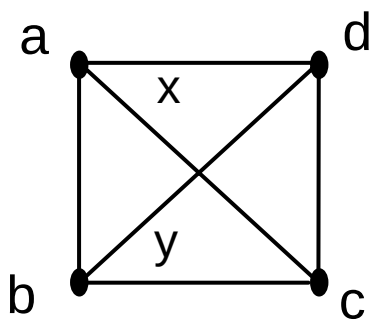

FIGURE: 1.3.11

i.e., the lines x and y intersect in the diagram their intersection is not a point of the graph.

Recall in the graph in figure 1.3.12.

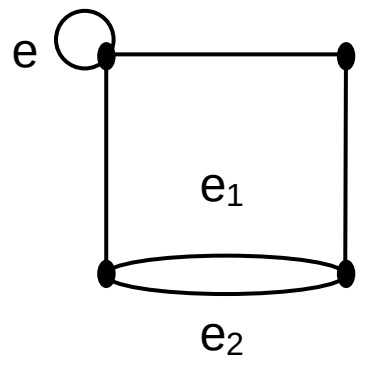

FIGURE: 1.3.12

e is loop and {$e_1$, $e_2$} is a set of multiple edges. Thus a graph with loops and multiple edges will be known by some authors as pseudo graphs. But we shall specify graphs with multiple edges and graphs with loops distinctly. A graph is simple if it has no loops and multiple edges.

Now we proceed on to give the application of graphs to fuzzy models.

Graphs have been the basis for several of the fuzzy models like binary fuzzy relations, sagittal diagrams, fuzzy compatibility relations fuzzy partial ordering relations, fuzzy morphisms fuzzy cognitive maps and fuzzy relational map models. Most of these models also basically rely on the under lying matrix. Both square and rectangular matrices are used. For more about these please refer [7, 19-20, 84].

Just for the sake of completeness we at each stage illustrate each of these by giving a brief definition and by example.

**DEFINITION 1.3.2:** *Let X and Y be two finite sets contrary to functions from X to Y, binary relations R (X, Y) may assign to each element X two or more elements of Y some basic operations on functions such as composition or inverse may also be applicable to binary fuzzy relations.*

*The fuzzy relation R (X, Y) from the domain set X to the range set Y depicted by bi partite graph with edge weights.*

***Example 1.3.1:*** Let X = {$x_1$, $x_2$, $x_3$ $x_4$ $x_5$} and Y = {$y_1$, $y_2$, $y_3$, $y_4$, $y_5$, $y_6$} suppose R (X, Y) fuzzy membership relation.

The related bipartite graph with associated edge weights.

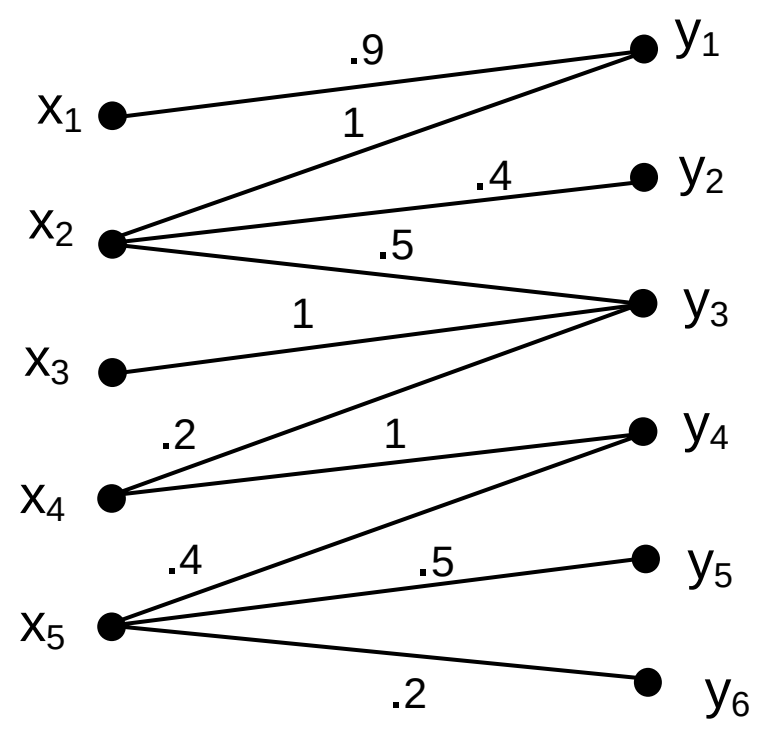

FIGURE: 1.3.13

The relational membership matrix

$$R\ (X, Y) = R = \begin{array}{c} \\ x_1 \\ x_2 \\ x_3 \\ x_4 \\ x_5 \end{array} \begin{array}{c} \begin{array}{cccccc} y_1 & y_2 & y_3 & y_4 & y_5 & y_6 \end{array} \\ \begin{bmatrix} .9 & 0 & 0 & 0 & 0 & 0 \\ 1 & .4 & .5 & 0 & 0 & 0 \\ 0 & 0 & 1 & 0 & 0 & 0 \\ 0 & 0 & .2 & 1 & 0 & 0 \\ 0 & 0 & 0 & .4 & .5 & .2 \end{bmatrix} \end{array}.$$

The extension of two binary relation is a ternary relation denoted by R (X, Y, Z) where P (X, Y) and Q (Y, Z) are fuzzy relations from X to Y and Y to Z respectively.

Suppose X = $\{x_1\ x_2\ x_3\ x_4\}$ Y = $\{y_1\ y_2\ y_3\}$ and Z = $\{Z_1\ Z_2\}$ The representation of the composition of the two binary relations given by R (X, Y, Z) is given by the following graph with the desired edge weights

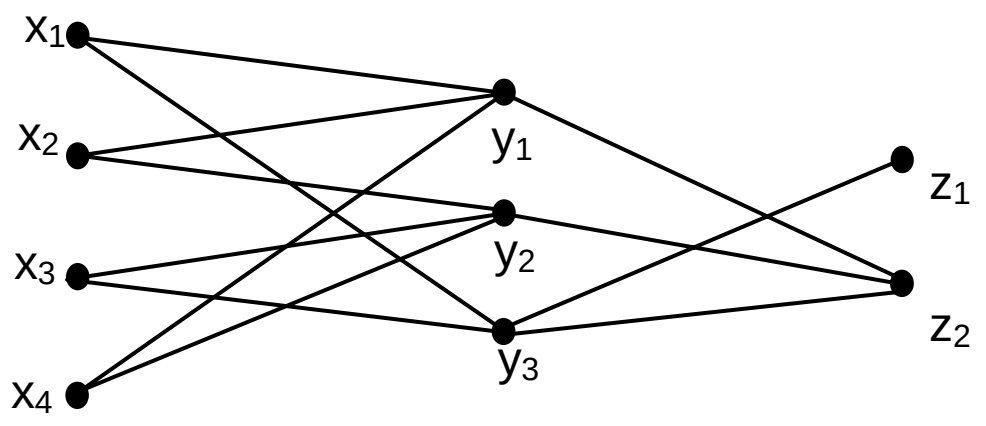

FIGURE: 1.3.14

***Example 1.3.2***: The Binary relation on a single set. The related weighted graph or the sagittal diagram is follows. Let X = {$x_1$, $x_2$ $x_3$ $x_4$ $x_5$}.

Membership matrix

$$\begin{array}{c} \\ x_1 \\ x_2 \\ x_3 \\ x_4 \\ x_5 \end{array} \begin{array}{c} \begin{array}{ccccc} x_1 & x_2 & x_3 & x_4 & x_5 \end{array} \\ \begin{bmatrix} .8 & .3 & 0 & 0 & 0 \\ 0 & .7 & 0 & .8 & 0 \\ 0 & 0 & 0 & 0 & .3 \\ 0 & 0 & 0 & .4 & .3 \\ 0 & .1 & .6 & 0 & 0 \end{bmatrix} \end{array}.$$

The related bipartite graph

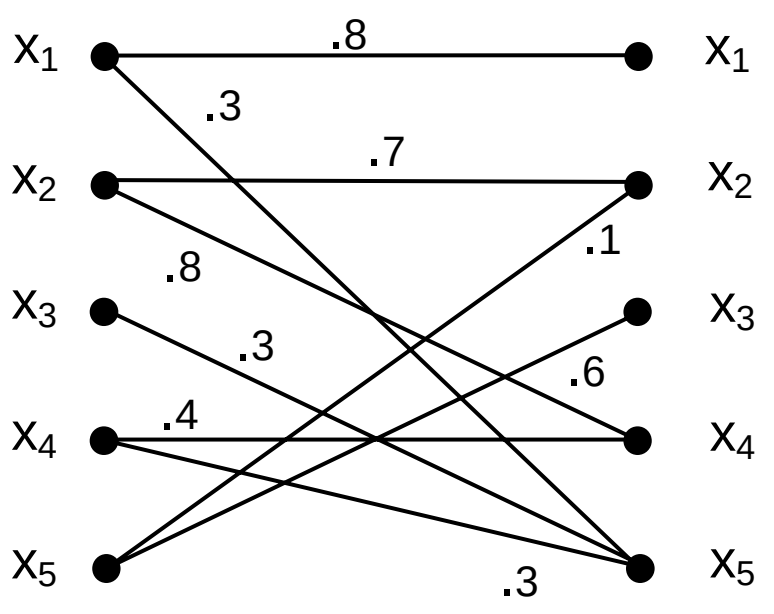

FIGURE: 1.3.15

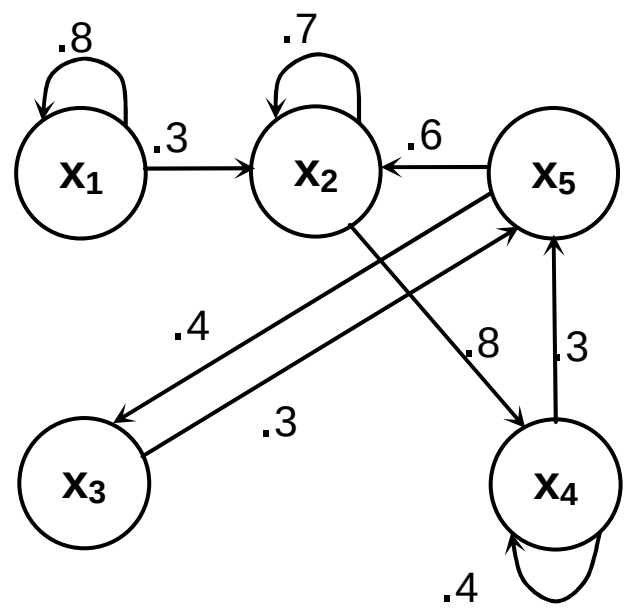

FIGURE: 1.3.16

The graph representation by a single diagram.

Study in this direction using this graph will be carried out for any appropriate models which is under investigation.

**DEFINITION 1.3.3:** *A binary relation R (X, X) that is reflexive and symmetric is usually called a compatibility relation or tolerance relation.*

***Example 1.3.3:*** Consider a fuzzy relation R (X, X) defined on X = {$x_1, x_2, \ldots, x_7$} by the following membership matrix.

$$\begin{array}{c} \\ x_1 \\ x_2 \\ x_3 \\ x_4 \\ x_5 \\ x_6 \\ x_7 \end{array}\begin{array}{c} \begin{array}{ccccccc} x_1 & x_2 & x_3 & x_4 & x_5 & x_6 & x_7 \end{array} \\ \begin{bmatrix} 1 & .3 & 0 & 0 & 0 & 0 & .2 \\ .3 & 1 & 0 & 0 & 0 & 0 & 0 \\ 0 & 0 & 1 & 1 & 0 & .6 & 0 \\ 0 & 0 & 1 & 1 & 0 & .1 & 0 \\ 0 & 0 & 0 & 0 & 1 & .4 & 0 \\ 0 & 0 & .6 & .1 & .4 & 1 & 0 \\ .2 & 0 & .0 & 0 & 0 & 0 & 1 \end{bmatrix} \end{array}.$$

Graph of the compatibility relation or the compatibility relation graph.

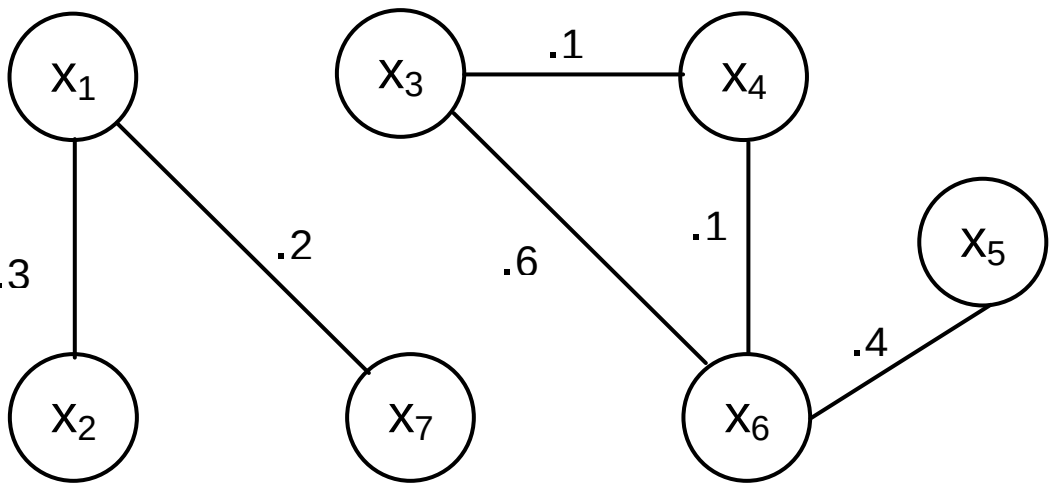

FIGURE: 1.3.17

Now we proceed on to define fuzzy partial ordering and its graphical representation.

**DEFINITION 1.3.4:** *A fuzzy binary relation R on a set X is a fuzzy partial ordering if and only if it is reflexive antisymmetric and transitive under some form of fuzzy transitivity.*

***Example 1.3.4:*** Let X = {a, b, c, d, e, f}. P Q and Q denote the crisp partial ordering on the set X which are defined by their membership matrices and their graphical representation.

$$P = \begin{array}{c} \\ a \\ b \\ c \\ d \\ e \\ f \end{array}\!\!\begin{array}{c} \begin{array}{cccccc} a & b & c & d & e & f \end{array} \\ \begin{bmatrix} 1 & 0 & 0 & 0 & 0 & 0 \\ \underline{1} & 1 & 0 & 0 & 0 & 0 \\ 1 & \underline{1} & 1 & 0 & 0 & 0 \\ 1 & 1 & \underline{1} & 1 & 0 & 0 \\ 1 & 1 & 1 & \underline{1} & 1 & 0 \\ 1 & 1 & 1 & 1 & \underline{1} & 1 \end{bmatrix} \end{array}.$$

The graph

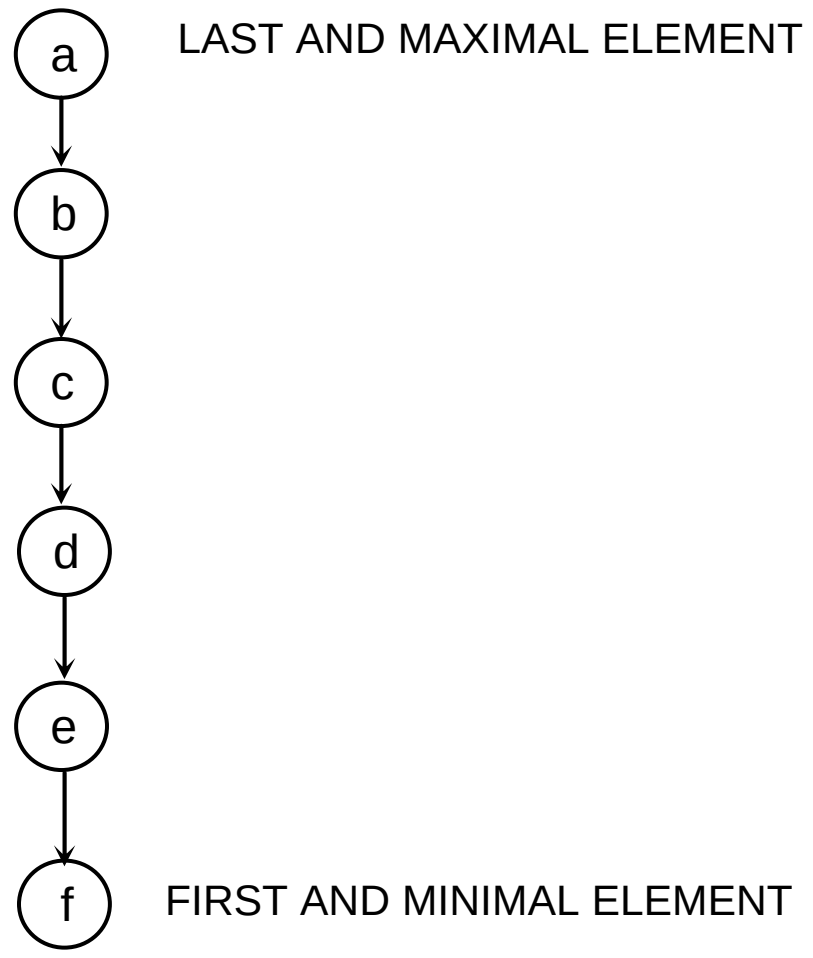

FIGURE: 1.3.18

The related matrix Q and its graphical representation

$$Q = \begin{array}{c} \\ a \\ b \\ c \\ d \\ r \\ f \end{array}\!\!\begin{array}{c} \begin{array}{cccccc} a & b & c & d & e & f \end{array} \\ \begin{bmatrix} 1 & 0 & 0 & 0 & 0 & 0 \\ \underline{1} & 1 & 0 & 0 & 0 & 0 \\ 1 & \underline{1} & 1 & 0 & 0 & 0 \\ 1 & \underline{1} & 0 & 1 & 0 & 0 \\ 1 & 1 & \underline{1} & \underline{1} & 1 & 0 \\ 1 & 1 & 1 & 1 & \underline{1} & 1 \end{bmatrix} \end{array}.$$

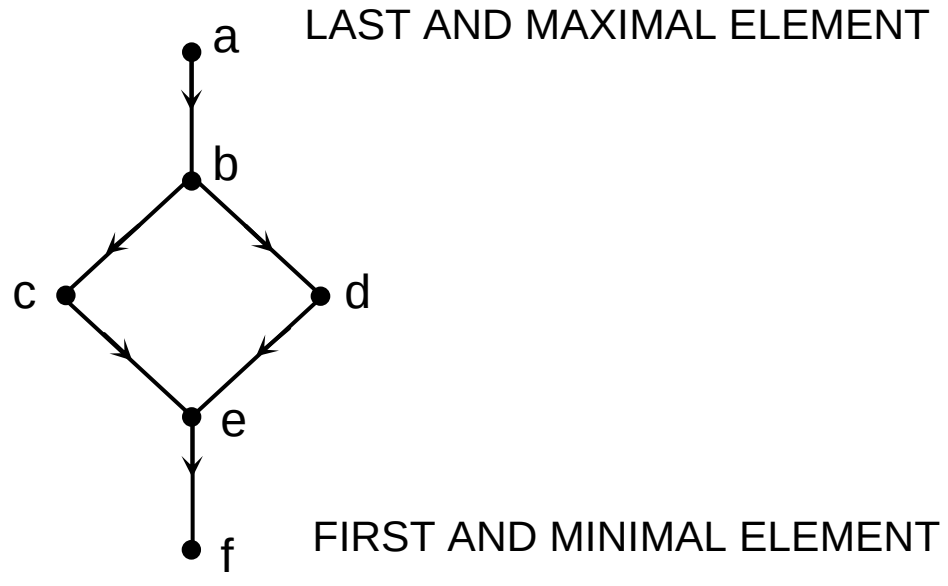

FIGURE: 1.3.19

The membership matrix of R

$$R = \begin{array}{c} \\ a \\ b \\ c \\ d \\ r \\ f \end{array} \begin{array}{c} \begin{array}{cccccc} a & b & c & d & e & f \end{array} \\ \begin{bmatrix} 1 & 0 & 0 & 0 & 0 & 0 \\ \underline{1} & 1 & 0 & 0 & 0 & 0 \\ 1 & \underline{1} & 1 & 0 & 0 & 0 \\ \underline{1} & 0 & 0 & 1 & 0 & 0 \\ 1 & 1 & \underline{1} & \underline{1} & 1 & 0 \\ 1 & 1 & 1 & 1 & \underline{1} & 1 \end{bmatrix} \end{array}.$$

The related graph

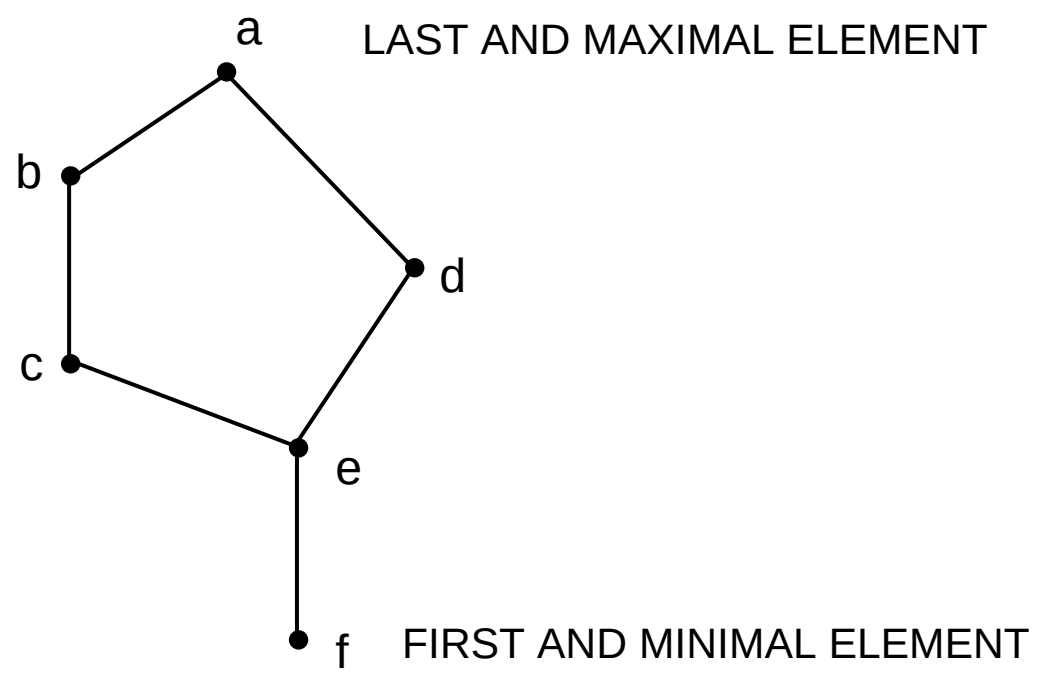

FIGURE: 1.3.20

Next we recall the definition of morphism and there graphical representation.

**DEFINITION 1.3.5**: *Let R (X, X) and Q (Y, Y) be fuzzy binary relations on the sets X and Y respectively A function h : X → Y is said to be a fuzzy homomorphism from (X, R) to (Y Q) if R ($x_1$, $x_2$)≤ Q (h ($x_1$), h ($x_2$) for all $x_1$, $x_2$ ∈ X and their images h ($x_1$), h ($x_2$)∈ Y.*

*Thus, the strength of relation between two elements under R is equated or excepted by the strength of relation between their homomorphic images under Q.*

It is possible for a relation to exist under Q between the homomorphic images of two elements that are themselves unrelated under R.

When this is never the case under a homomorphic function h, the function is called a strong homomorphism. It satisfies the two implications ($x_1$, $x_2$) ∈ R implies (h ($x_1$), h ($x_2$)) ∈ Q for all $x_1$, $x_2$ ∈ X and ($y_1$, $y_2$) ∈ Q implies ($x_1$, $x_2$) ∈ R for all $y_1$ ,$y_2$ ∈ Y where $x_1 \in h^{-1}$ ($y_1$) and $x_2 \in h^{-1}$ ($y_2$).

***Example 1.3.5:*** Let X = {a, b, c} y = {α, β, γ, δ} be sets with the following membership matrices which represent the fuzzy relations R (X, X) and Q (Y, Y)

$$R = \begin{array}{c} \\ a \\ b \\ c \end{array} \overset{\begin{matrix}\alpha & \beta & \gamma & \delta\end{matrix}}{\begin{bmatrix} 0 & .5 & 0 & 0 \\ 0 & 0 & .9 & 0 \\ 1 & 0 & 0 & .5 \\ 0 & .6 & 0 & 0 \end{bmatrix}}$$

and

$$Q = \begin{array}{c} \alpha \\ \beta \\ \gamma \\ \delta \end{array} \overset{\begin{matrix}\alpha & \beta & \gamma & \delta\end{matrix}}{\begin{bmatrix} 0 & .5 & 0 & 0 \\ 0 & 0 & .9 & 0 \\ 1 & 0 & 0 & .5 \\ 0 & .6 & 0 & 0 \end{bmatrix}}.$$

The graph of the ordinary fuzzy homomorphism.

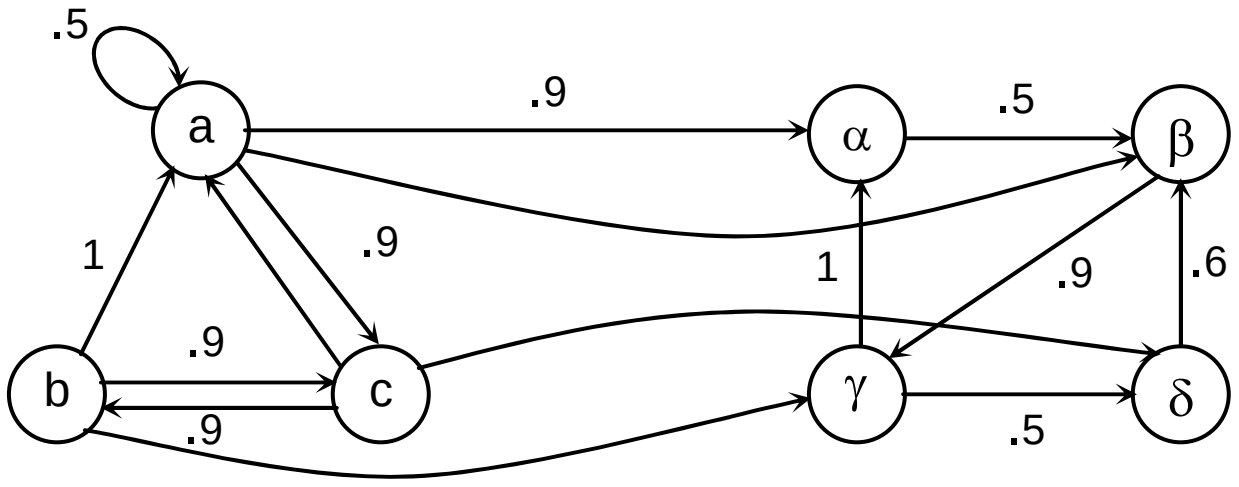

FIGURE: 1.3.21

Now we leave it for the reader to construct a strong fuzzy homomorphism and illustrate it by a graph. Next we just recall the definition of Fuzzy Cognitive maps (FCMs) and illustrate the directed graph of a FCM by an example.

**DEFINITION 1.3.6:** *An FCM is a directed graph with concepts like policies, events etc as nodes and causalities as edges. It represents causal relationship between concepts.*

The problem studied in this case is for a fixed source S, a fixed destination D and a unique route from the source to the destination, with the assumption that all the passengers travel in the same route, we identify the preferences in the regular services at the peak hour of a day.

We have considered only the peak-hour since the passenger demand is very high only during this time period, where the transport sector caters to the demands of the different groups of people like the school children, the office goers, the vendors etc.

We have taken a total of eight characteristic of the transit system, which includes the level of service and the convenience factors.

We have the following elements, Frequency of the service, in-vehicle travel time, the travel fare along the route, the speed of the vehicle, the number of intermediate points, the waiting time, the number of transfers and the crowd in the bus or equivalently the congestion in the service.

Before defining the cognitive structure of the relationship, we give notations to the concepts involved in the analysis as below.

| | | |
|---|---|---|
| $C_1$ | - | Frequency of the vehicles along the route |
| $C_2$ | - | In-vehicle travel time along the route |
| $C_3$ | - | Travel fare along the route |
| $C_4$ | - | Speed of the vehicles along the route |
| $C_5$ | - | Number of intermediate points in the route |
| $C_6$ | - | Waiting time |
| $C_7$ | - | Number of transfers in the route |
| $C_8$ | - | Congestion in the vehicle. |

The graphical representation of the inter-relationship between the nodes is given in the form of directed graph given in Figure: 1.3.22.

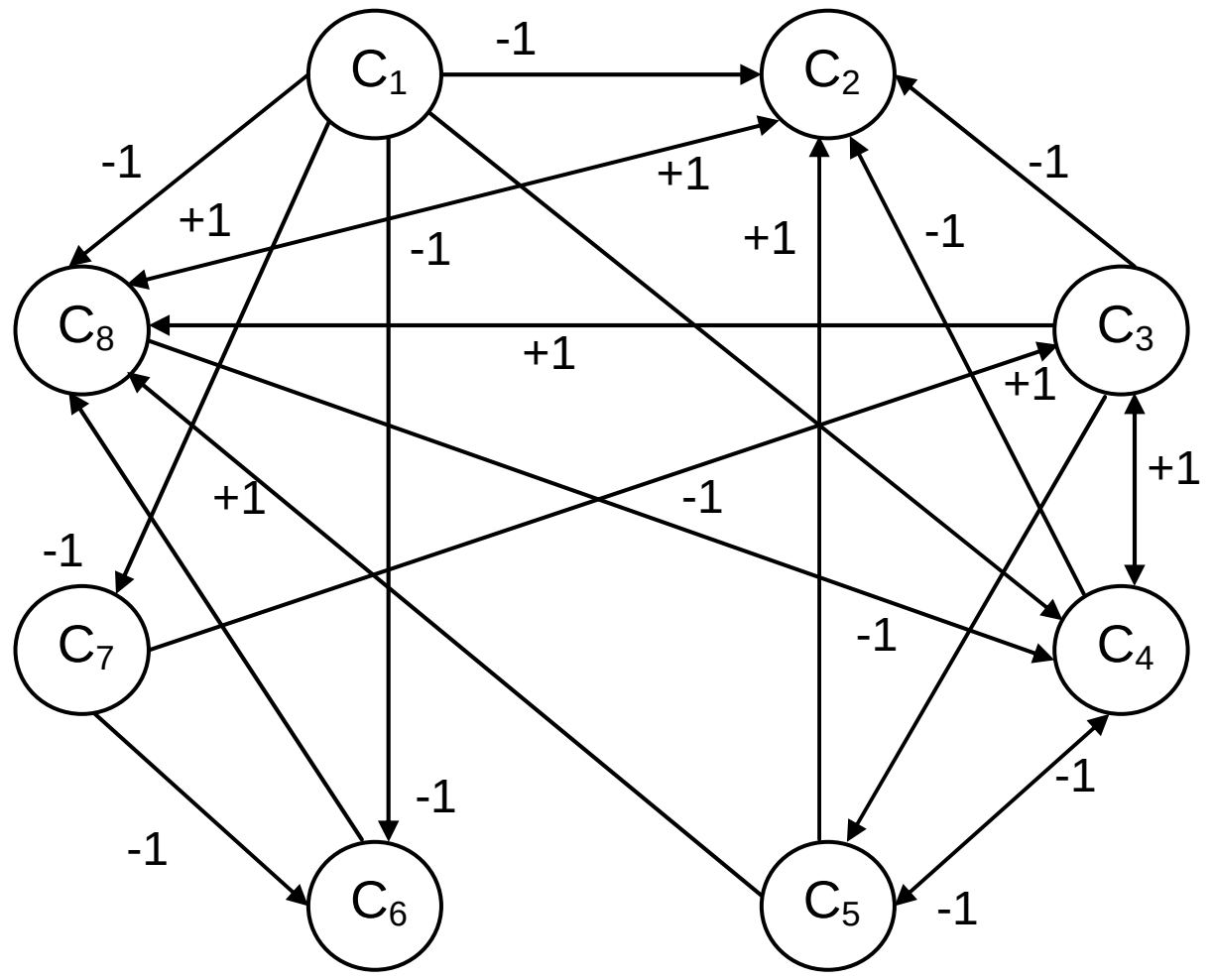

FIGURE: 1.3.22

From the above signed directed graph, we obtain a connection matrix E, since the number of concepts used here are eight, the connection matrix is a $8 \times 8$ matrix.

Thus we have $E = [A_y]_{8\times8}$

$$E = \begin{bmatrix} 0 & -1 & 0 & 1 & 0 & -1 & -1 & -1 \\ 0 & 0 & 0 & 0 & 0 & 0 & 0 & 1 \\ 0 & -1 & 0 & 1 & -1 & 0 & 0 & -1 \\ 0 & -1 & 1 & 0 & -1 & 0 & 0 & 0 \\ 0 & 1 & 0 & -1 & 0 & 0 & 0 & 1 \\ 0 & 0 & 0 & 0 & 0 & 0 & 0 & 1 \\ 0 & 0 & 1 & 0 & 0 & -1 & 0 & 0 \\ 0 & 1 & 0 & -1 & 0 & 0 & 0 & 0 \end{bmatrix}.$$

Now we just recall the definition of Fuzzy Relational Maps (FRMs).

**DEFINITION 1.3.7:** *A FRM is a directed graph or a map from D to R with concepts like policies or events etc as nodes and causalities as edges. It represents causal relation between spaces D and R.*

The employee-employer relationship is an intricate one. For, the employers expect to achieve performance in quality and production in order to earn profit, on the other hand employees need good pay with all possible allowances.

Here we have taken three experts opinion in the study of Employee and Employer model.

The three experts whose opinions are taken are the Industry Owner, Employees' Association Union Leader and an Employee. The data and the opinion are taken only from one industry.

Using the opinion we obtain the hidden patterns. The following concepts are taken as the nodes relative to the employee.

We can have several more nodes and also several experts' opinions for it a clearly evident theory which professes that more the number of experts the better is the result.

We have taken as the concepts / nodes of domain only 8 notions which pertain to the employee.

| | | |
|---|---|---|
| $D_1$ | – | Pay with allowances and bonus to the employee |
| $D_2$ | – | Only pay to the employee |
| $D_3$ | – | Pay with allowances (or bonus) to the employee |
| $D_4$ | – | Best performance by the employee |
| $D_5$ | – | Average performance by the employee |
| $D_6$ | – | Poor performance by the employee |
| $D_7$ | – | Employee works for more number for hours |
| $D_8$ | – | Employee works for less number of hours. |

$D_1, D_2, \ldots, D_8$ are elements related to the employee space which is taken as the domain space.

We have taken only 5 nodes / concepts related to the employer in this study.

These concepts form the range space which is listed below.

| | | |
|---|---|---|
| $R_1$ | – | Maximum profit to the employer |
| $R_2$ | – | Only profit to the employer |
| $R_3$ | – | Neither profit nor loss to the employer |
| $R_4$ | – | Loss to the employer |
| $R_5$ | – | Heavy loss to the employer |

The directed graph as given by the employer is given in Figure 1.3.23.

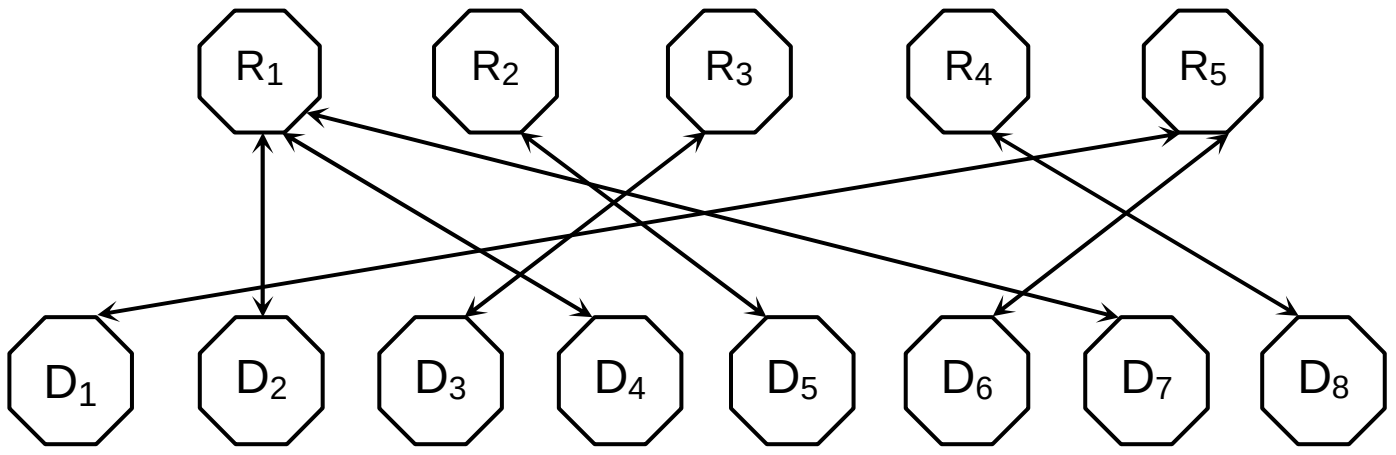

FIGURE: 1.3.23

The associated relational matrix $E_1$ of the employer as given by following;

$$E_1 = \begin{bmatrix} 0 & 0 & 0 & 0 & 1 \\ 1 & 0 & 0 & 0 & 0 \\ 0 & 0 & 1 & 0 & 0 \\ 1 & 0 & 0 & 0 & 0 \\ 0 & 1 & 0 & 0 & 0 \\ 0 & 0 & 0 & 0 & 1 \\ 1 & 0 & 0 & 0 & 0 \\ 0 & 0 & 0 & 1 & 0 \end{bmatrix}.$$

Thus we give several of the fuzzy models to make use of the graphs.

The formatting of Linguistic questionnaire for making use of fuzzy models is briefly given here [220-2].

When we have the data in hand i.e., the related statistics we can easily make use of any statistical methods or operational research and analyse the data.

The analysis can be from the simple calculation of mean and standard deviation or any other higher statistical technique appropriate to analyse the problem in hand.

This is the case when we have the supervised data i.e., the statistical data. This is not always possible for many a times when the need for study of any particular problem arises we may not always be guaranteed of a statistical data related with it.

For instance one is interested in the study of the mental problems undergone by migrant labourers affected with HIV/AIDS or for in general the social problems faced by a HIV/AIDS infected person. This may vary from person to person, depending on gender, their economic status, educational qualification and above all their place of living, urban, semi urban or rural and so on.

The humiliation faced by them in the midst of their relatives will also be distinct for all these cannot be captured by number so only, we are forced to use fuzzy models. To capture these feelings we make use of a linguistic questionnaire which may give proper expressions to the gradation of feelings for instance the mental stress felt by a female or youth may be very different from the stress felt by a considerably old man in his late forties or by a house wife infected by her husband.

These vagueness can be given some weightage by means of these linguistic questionnaire. For the sample questionnaire please refer [220-2, 233].

So even before the formation of or making the linguistic questionnaire the concerned socio-scientist is expected to make a pilot survey by meeting the persons from the category from which he wants to study and put forth some questions.

The first suggestion is to allow them to speak out their experience, problems, feelings, their expected solutions and so on. For this alone will show so much of weight to our problem at hand; so that we can from these persons views and problem formulate more questions with lot of relevance to our problem.

Basically this collection of data done by using a pilot survey will enable the socio-scientists or the experts who wish to make this study to form a good linguistic questionnaire.

Specially by a linguistic questionnaire we only mean a usual questionnaire in which gradations are given; for instance, suppose we want to study the feelings of a student about his teacher, it cannot be just very good or very bad or excellent or worst, their can be many feelings or qualities good or bad in between these extreme feelings like bad, just bad neutral okay just good, rude, temperamental, very partial, partial at times and so on.

These gradations of feelings can be easily got only after the pilot survey.

It may be at times a facial change, just silence or just a laughter or a smile or so on. Depending on the cost that can be spent on any investigation, at times we may suggest the analyzer to use both audio and video for this pilot survey, if one wishes to make the study as accurate as possible, because at times the student even may become restless panic or agitated when one talks about school or teacher or at times pretend to cover his / her feelings so on and so forth.

Thus a sincere pilot survey alone can contribute for a very balanced, sensitive useful linguistic questionnaire. Once the linguistic questionnaire is ready still it is always not possible to get information by asking them to fill it for the person may be literate may be illiterate even if he is literate may not easily follow the questionnaire, so depending on the category of persons or public for study it is advisable to use tapes to record their feelings.

When we interview them if enough finance is available they can also video tape the interview provided the candidate is not uneasy with these equipments. Even taping should be done in such a way they feel free before the gadget or with permission done in an invisible way after giving in the full guarantee that their identity will be kept as confident.

Only on their primary wishes and if need be only by them it would be revealed! For several HIV/AIDS patients or sometimes people working in government or private sectors or say even students may not wish to reveal their identity. So in general they should be given confidence that their identity would never be revealed.

Even after collection of data, it is not an easy task to apply the model, for this data needs experts views what they wish to achieve from their study, is it only time dependent, do they expect a gradation in their solution or do they want the hidden pattern, or do they want sets of solution. Or if they have predicted results and are interested for a feasible solution so on and so forth.

Only after knowing this they can choose their fuzzy models. Also at times the study and analysis of the linguistic questionnaire by some fuzzy model experts can throw the light to choose an appropriate fuzzy model for the problem. Thus we see we have to use a lot of caution and at the same time expertise in making the choice of the model. Once such a choice is made certainly one is guaranteed of a proper sensitive and accurate solution for these unsupervised data. For more about questionnaire refer [220-2, 233].

We have given some linguistic questionnaire in appendix of the books [220-2, 233]. For instance the Fuzzy Cognitive Maps (FCMs) and Fuzzy Relational Maps (FRM) are the models that will give the hidden pattern. The fuzzy matrix model discussed in section 2.1 of this book will give time dependent solutions, where as the Fuzzy Associative Memories (FAMs) model will give a solution closer to the expected solution and so on.

**Chapter Two**

# DESCRIPTION OF SIMPLE FUZZY MODELS AND THEIR APPLICATIONS TO REAL WORLD PROBLEMS

In this chapter we introduce 6 simple fuzzy models using fuzzy matrices and give their applications to the real world problems. In section one of this chapter the notion of fuzzy models using the concept of mean and Standard Deviation (SD) of the columns of the real data matrices obtained using the real statistical data is described. Here the statistical data is formed into a matrix. Using each column of this matrix the mean and standard deviation of the columns are obtained. This mean and SD of the columns are used to convert the statistical data matrix into a fuzzy matrix. This fuzzy matrix model so formed is used in the analysis of the data. Description of the Fuzzy Cognitive Maps (FCMs) is given in section two and their applications to the real world problems are also illustrated. Section three gives the generalization and description of the FCMs; known as Fuzzy Relational Maps (FRMs). Illustrations using real world problems are also given in this section.

The concept of Bidirectional Associative memories (BAM) model is recalled in section four and its applications are also given. Section five gives the concept of Fuzzy Associative Memories (FAM) model and their applications to real social world problems is described. In the final section the notion of Fuzzy Relations Equations (FREs) model is given. We give simple models to illustrate them, so that a socio scientist is able to work with them with ease.

## 2.1 Description of Simple Fuzzy Matrix Model

In this section we describe a simple fuzzy matrix model when we have a raw data in hand. We describe how this raw data is transformed into a fuzzy matrix model using the techniques of average (mean) and standard deviation. The steps involved are; first the raw data at hand is converted or transformed into a time dependent matrix.

After obtaining the time dependent matrix using the techniques of average and standard deviation the time dependent data matrix is converted into an Average Time Dependent Data matrix (ATD-matrix).

In the next stage i.e., in the third stage, to make the calculations simple and easy we make use of simple average technique, to convert the ATD matrix into a fuzzy matrix with entries $e_{ij}$, where $e_{ij} \in \{-1, 0, 1\}$. We call this matrix as the Refined Time Dependent Data matrix (RTD-matrix). The value of $e_{ij}$ is obtained in a special way. At the fourth stage we get the Combined Effect Time Dependent Data Matrix (CETD matrix) which gives the cumulative effect of all these entries. In the final stage we obtain the row sums using C-program or Java program. Graphs are drawn taking the row sum of the CETD matrix along the y axis and time scale along the x-axis. These simple graphs are understandable even by a layman. Hence this method is very effective at the same time simple.

Here the real data collected from the agricultural labourers suffering from health hazards due to chemical pollution is studied using this fuzzy matrix model. Using this model we estimate the maximum age group in which the agricultural labourers suffer health hazards due to chemical pollution. In this section, we give an algebraic approach to the health hazards faced by agriculture labourers due to chemical pollution in Chengalpet District in Tamil Nadu. First by the term chemical pollution we mean the pollution due to the spray of pesticides and insecticides, also the pollution of the grain due to the use of fertilizer, which are mainly chemicals. From our interviews and the fieldwork undertaken we saw that agriculture labourer who are free of any tensions or mental problems suffer from these symptoms, which were mainly due to chemical pollution, this was spelt out by 95% of the interviewed agricultural labourers. They all described their health problems was acute during and after the spray of pesticides, insecticides and manuring. It is pertinent to mention that they informed us that two persons died in the agriculture field itself while they were spraying the pesticides.

They said these victims fainted in the field and before they could be taken to the hospital they died of suffocation due to the spray. Further, the older people were very angry at this event. In fact, they expressed in those days the agriculture methods used by them were really eco friendly. Due to modernization nowadays the pesticides and insecticides are sprayed on the fields using helicopters; which has largely affected the health conditions of the agriculture labourers. For the total atmosphere is polluted by these methods. This study is significant because most of the villages in India adopt the same procedure. To the best of our knowledge no one has ever cared to study the health hazards suffered by these people due to chemical pollution. Thus our study can be adopted to any agriculture field/village in India. We approach the problem of pollution by determining the peak age group in which they are maximum affected by pollution. By knowing this age group the government at least can take steps to treat them and rehabilitate them and give medicines which will antinode the chemicals due to which they are suffering these health problems. We analyze these problems using fuzzy matrix. In the first stage, the raw data is given the matrix representation. Entries corresponding to the intersection of rows and columns, are values corresponding to a live network. The raw data, as it is transformed into a raw time dependent data matrix by taking along the rows the age group and along the columns the health problems the agricultural labourer or coolie suffers. Using the raw data matrix, we convert it into the Average Time Dependent Data (ATD) matrix ($a_{ij}$) by dividing each entry of the raw data matrix by the number of years that is, the time period. This matrix represents a data, which is totally uniform. At the third stage, the average or mean and the Standard Deviation (S.D) of every column in the ATD matrix, are determined. Using the average $\mu_j$ of each $j^{th}$ column and $\sigma_j$ the S.D of each $j^{th}$ column, a parameter $\alpha$ from the interval [0, 1] is chosen and the Refined Time Dependent Data matrix (RTD matrix) ($e_{ij}$) is formed using the formula:

$$\text{if } a_{ij} \leq (u_j - \alpha * \sigma_j) \text{ then } e_{ij} = -1$$
$$\text{else if } a_{ij} \in (u_j - \alpha * \sigma_j, u_j + \alpha * \sigma_j) \text{ then } e_{ij} = 0$$
$$\text{else if } a_{ij} \geq (u_j + \alpha * \sigma_j) \text{ then } e_{ij} = 1, \text{ where, } a_{ij}\text{'s are}$$

the entries of the ATD matrix.

The ATD matrix is thus converted into the Refined Time Dependent Data Matrix. This matrix is also at times termed as the fuzzy matrix as the entries are 1, 0, and –1. Now, the row sum of this matrix gives the maximum age group, who are prone to heath hazards. One can combine these matrices by varying the parameter $\alpha \in [0, 1]$, so that the Combined Effective Time Dependent Data (CETD) matrix is obtained. The row sum is found out for the CETD matrix and conclusions are derived based on the row sums. All these are represented by graphs and graphs play a vital role in exhibiting the data by the simplest means that can be even understood by a layman.

The raw data, which we have obtained from the 110 agriculture labours, is transformed in to time dependent matrices. After obtaining the time dependent matrices using the techniques of average and standard deviation we identify the peak age group in which they suffer the maximum health hazard. Identification of the maximum age group will play a vital role in improving their health conditions by providing them the best health facilities, like medicine, good food and better hygiene.

To the best of our knowledge such study has not been mathematically carried out by anyone. The study in general using these fuzzy matrix have been carried out in the analysis of transportation of passenger problems and in the problem of migrant labourers affected by HIV/AIDS [225, 253]. Apart from this we do not have any study using these techniques. The raw data (health problems) under investigation has been classified under four broad heads:

1. Cardiovascular problems,
2. Digestive problems,
3. Nervous problems and
4. Respiratory problems

These four major heads are further divided into subheads. These four major heads with the introduction and conclusion forms the subsections of this section. Now each of these four broad heads are dealt separately. For each broad head is divided into 8 or more or less subheads and fuzzy matrix /RTD matrix model just described is adopted. These subheads forms the columns of the matrix and the age of the agriculture labourers grouped in varying intervals forms the rows of the RTD matrix. Estimation of the maximum age group is a five stage process. This section has five subsections and each of them deal with different types of symptoms and diseases suffered by agriculture labourer in all age groups. We have varied the age group as well as the parameter $\alpha \in [0, 1]$ to obtain the most sensitive result. We have grouped the data under three sets of groups. We have formed in all the cases the matrix as well as the graphs. It is pertinent to mention here that such type of study among agriculture labourers and their health problem is for the first time analyzed by these methods. Further, we wish to state that this is the most simple and very effective method of analysis for even the lay man looking at the graph can conclude under which age group they suffer which type of disease. At each and every stage we give the sketch of the conclusions at the end of the graph so that the reader is able to get a clear grasp of the working of this model.

### 2.1.1 Estimation of the maximum age group of the agricultural labourers having cardio vascular problem due to chemical pollution, using matrices

One of the major and broad heads of study of the health problems faced by the agricultural labourers was the cardio-vascular problem under this disease the agriculture labourers suffered eight symptoms viz.

$S_1$ - Chest pain,
$S_2$ - Pain at the rib's sides,
$S_3$ - Back pain,
$S_4$ - Shoulder pain,
$S_5$ - Left arm and leg pain,
$S_6$ - Swollen limbs,
$S_7$ - Burning chest and
$S_8$ - Blood pressure (low or high B.P)

which are taken as the columns of the initial raw data matrix. The age group in years 20-30, 31-43 and 44-65 are taken as the rows of the matrix. The estimation of the maximum age group is a five-stage process. In the first stage we give the matrix representation of the raw data. Entries corresponding to the intersection of rows and columns are values corresponding to a live network. The $3 \times 8$ matrix is not uniform i.e. the number of individual years in each interval may not be the same. So in the second stage we in order to obtain an unbiased uniform effect on each and every data so collected, transform this initial matrix into an Average Time Dependent Data (ATD) matrix. To make the calculations easier and simpler we in the third stage using the simple average techniques convert the above average time dependent data matrix in to a matrix with entries $e_{ij} \in \{-1,0,1\}$.

We name this matrix as the Refined Time Dependent data matrix (RTD matrix) or as the fuzzy matrix. The value of $e_{ij}$ corresponding to each entry is determined in a special way described in section 2.1 of this book page 74. At the fourth stage using the fuzzy matrices we obtain the Combined Effect Time Dependent Data matrix (CETD matrix), which gives the cumulative effect of all these entries. In the final stage we obtain the row sums of the CETD matrix. The tables given are self-explanatory at each stage. The graphs of the RTD matrix and CETD matrix are given.

#### 2.1.1.1 Estimation of maximum age group, using 3 × 8 matrices

Initial raw data matrix of cardio vascular problem of order 3 × 8

| Years | $S_1$ | $S_2$ | $S_3$ | $S_4$ | $S_5$ | $S_6$ | $S_7$ | $S_8$ |
|---|---|---|---|---|---|---|---|---|
| 20-30 | 23 | 18 | 24 | 16 | 29 | 10 | 16 | 10 |
| 31-43 | 38 | 32 | 38 | 31 | 35 | 18 | 33 | 10 |
| 44-65 | 22 | 21 | 21 | 22 | 20 | 11 | 20 | 5 |

The ATD Matrix of cardio vascular problem of order 3 × 8

| Years | $S_1$ | $S_2$ | $S_3$ | $S_4$ | $S_5$ | $S_6$ | $S_7$ | $S_8$ |
|---|---|---|---|---|---|---|---|---|
| 20-30 | 2.09 | 1.64 | 2.18 | 1.46 | 2.64 | 0.91 | 1.46 | 0.91 |
| 31-43 | 2.92 | 2.46 | 2.92 | 2.39 | 2.69 | 1.39 | 2.54 | 0.77 |
| 44-65 | 1 | 0.95 | 0.95 | 1 | 0.91 | 0.5 | 0.91 | 0.23 |

The Average and Standard Deviation of the above ATD matrix

| | | | | | | | | |
|---|---|---|---|---|---|---|---|---|
| Average | 2.00 | 1.68 | 2.02 | 1.62 | 2.08 | 0.93 | 1.64 | 0.64 |
| SD | 0.96 | 0.76 | 0.995 | 0.71 | 1.01 | 0.45 | 0.83 | 0.36 |

RTD matrix for $\alpha = 0.15$ Row sum matrix

$$\begin{bmatrix} 0 & 0 & 1 & -1 & 1 & 0 & -1 & 1 \\ 1 & 1 & 1 & 1 & 1 & 1 & 1 & 1 \\ -1 & -1 & -1 & -1 & -1 & -1 & -1 & -1 \end{bmatrix} \qquad \begin{bmatrix} 1 \\ 8 \\ -8 \end{bmatrix}$$

RTD matrix for $\alpha = 0.35$ Row sum matrix

$$\begin{bmatrix} 0 & 0 & 0 & 0 & 1 & 0 & 0 & 1 \\ 1 & 1 & 1 & 1 & 1 & 1 & 1 & 0 \\ -1 & -1 & -1 & -1 & -1 & -1 & -1 & -1 \end{bmatrix} \qquad \begin{bmatrix} 2 \\ 7 \\ -8 \end{bmatrix}$$

RTD matrix for α = 0.45 Row sum matrix

$$\begin{bmatrix} 0 & 0 & 0 & 0 & 1 & 0 & 0 & 1 \\ 1 & 1 & 1 & 1 & 1 & 1 & 1 & 0 \\ -1 & -1 & -1 & -1 & -1 & -1 & -1 & -1 \end{bmatrix} \qquad \begin{bmatrix} 2 \\ 7 \\ -8 \end{bmatrix}$$

RTD matrix for α = 0.75 Row sum matrix

$$\begin{bmatrix} 0 & 0 & 0 & 0 & 0 & 0 & 0 & 0 \\ 1 & 1 & 1 & 1 & 0 & 1 & 1 & 0 \\ -1 & -1 & -1 & -1 & -1 & -1 & -1 & -1 \end{bmatrix} \qquad \begin{bmatrix} 0 \\ 6 \\ -8 \end{bmatrix}$$

Now we give the graphs:

Graph 1: Graph depicting maximum age group of labourers with pollution-related Cardio-Vascular Problem (α = 0.15)

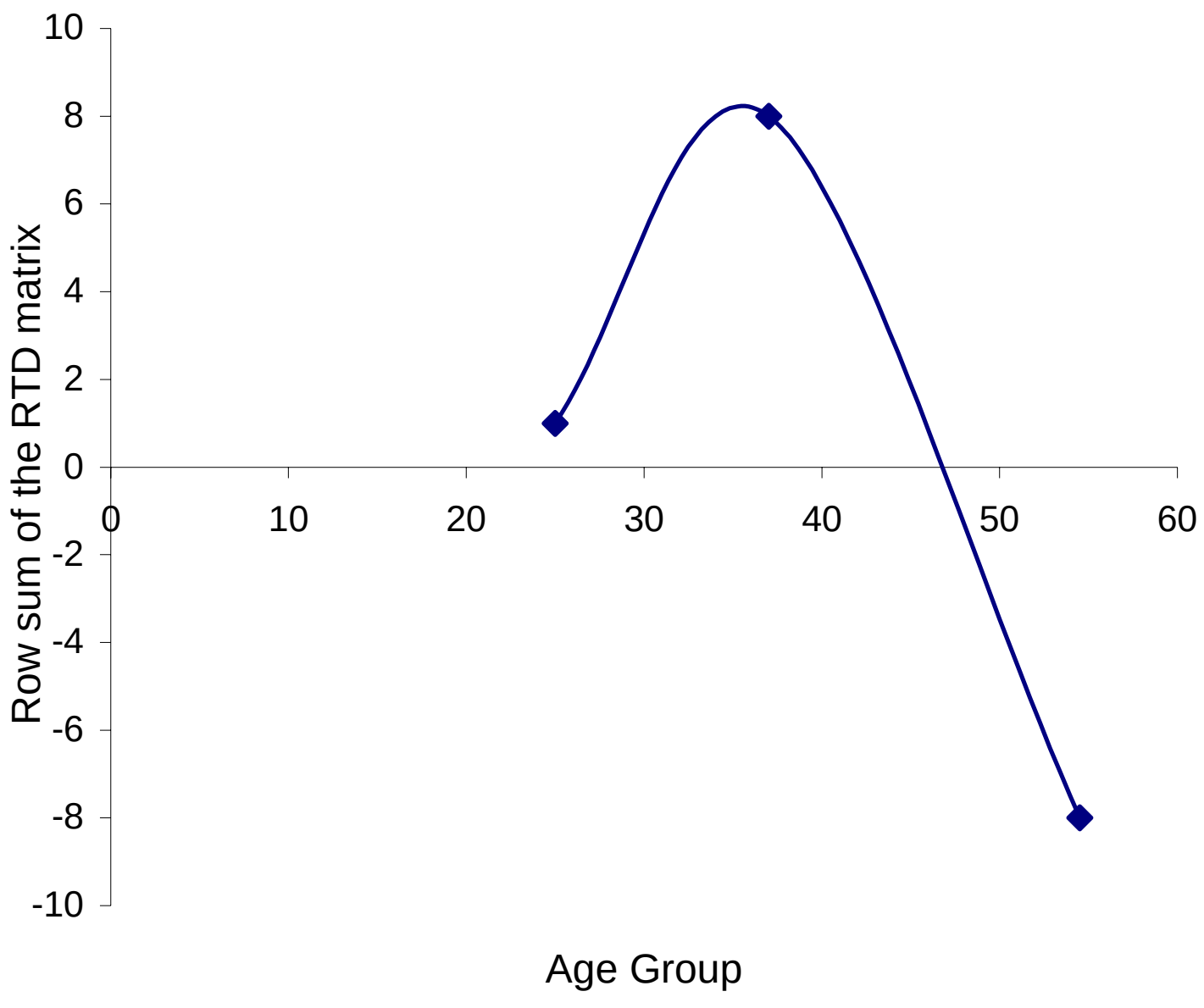

Graph 2: Graph depicting maximum age group of labourers with pollution- related Cardio-Vascular Problem ($\alpha$ = 0.35)

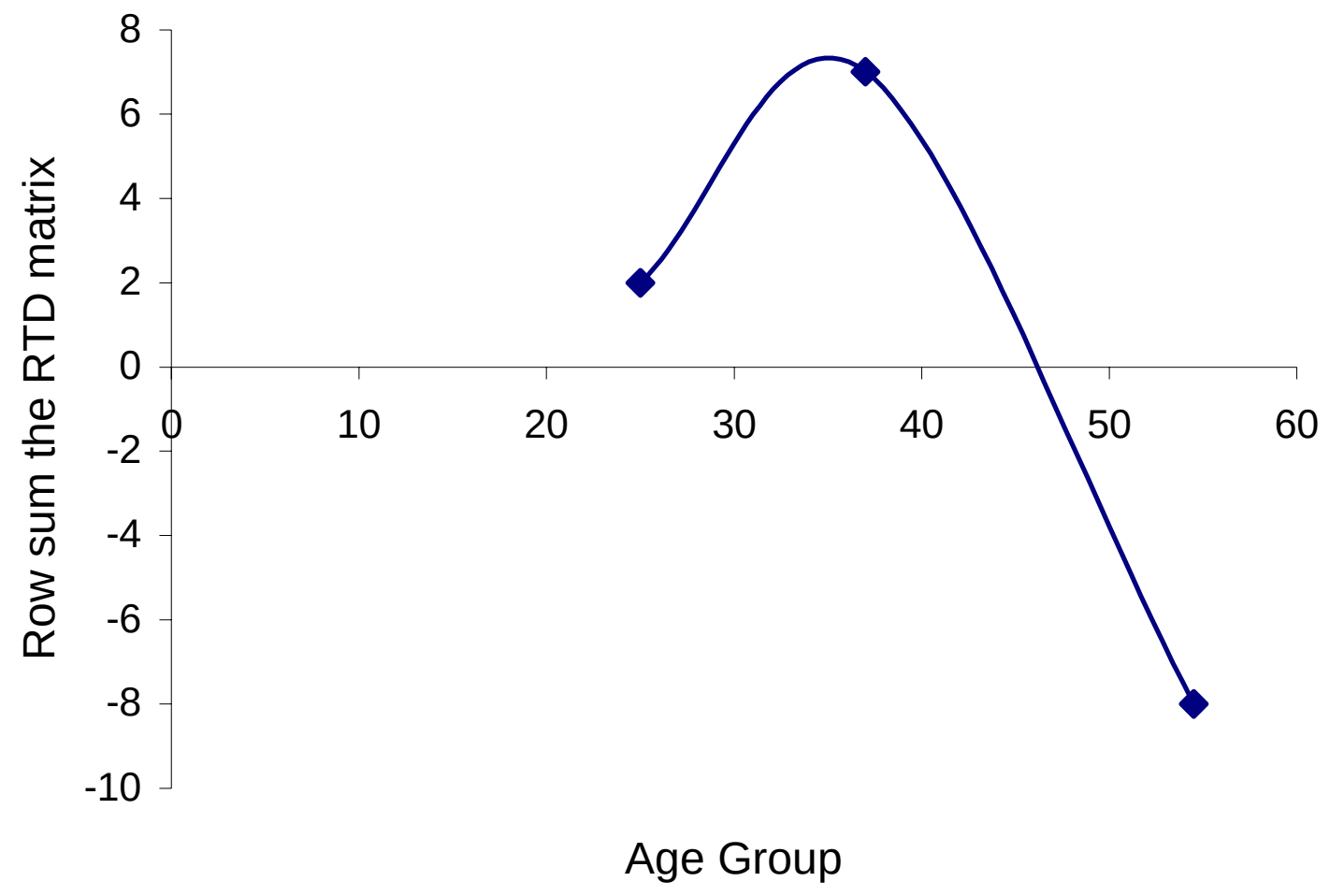

Graph 3: Graph depicting maximum age group of labourers with pollution-related Cardio- Vascular Problem ($\alpha$ = 0.45)

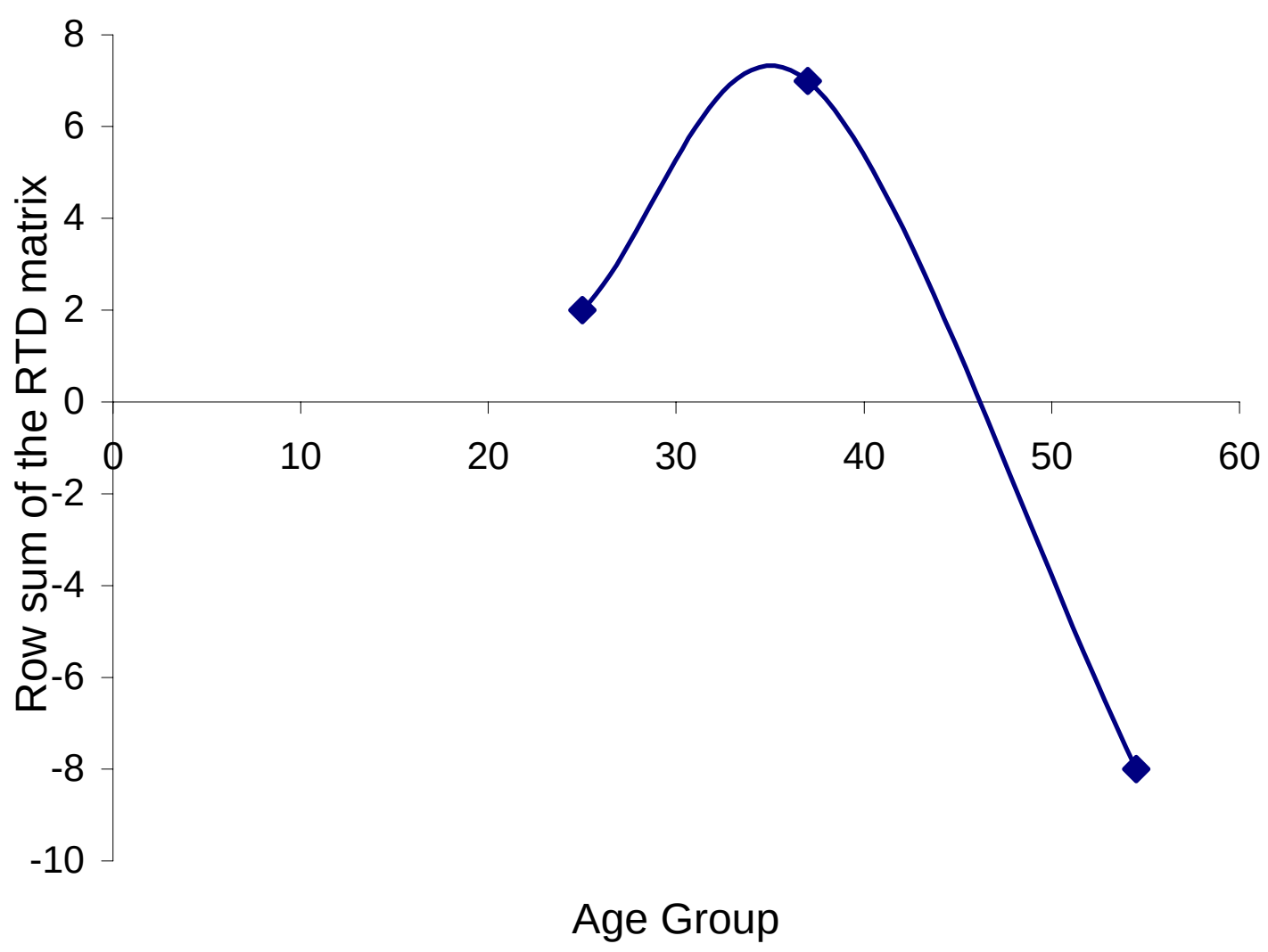

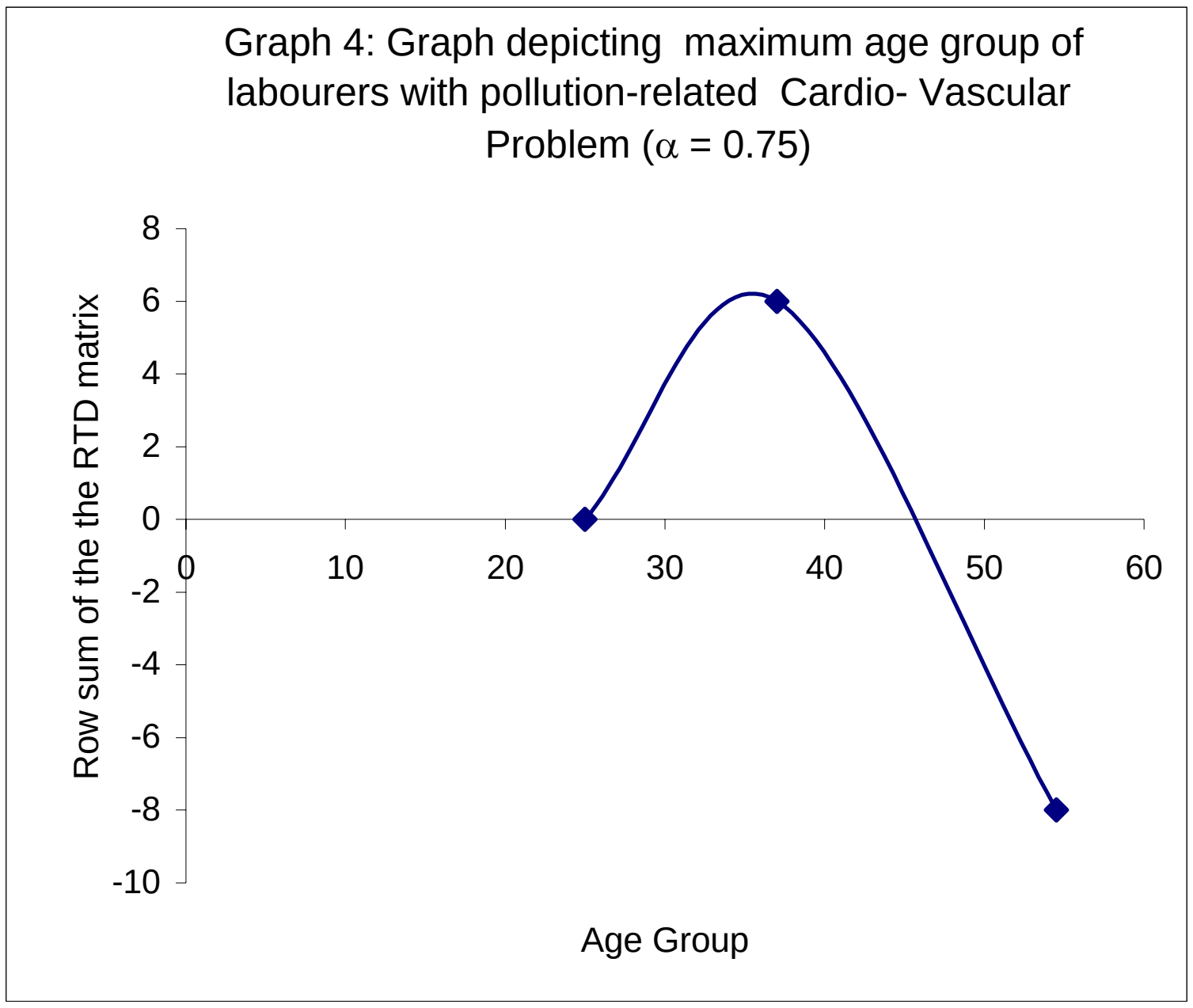
Graph 4: Graph depicting maximum age group of labourers with pollution-related Cardio- Vascular Problem (α = 0.75)
Row sum of the the RTD matrix
Age Group
8
6
4
2
0
-2
-4
-6
-8
-10
0
10
20
30
40
50
60

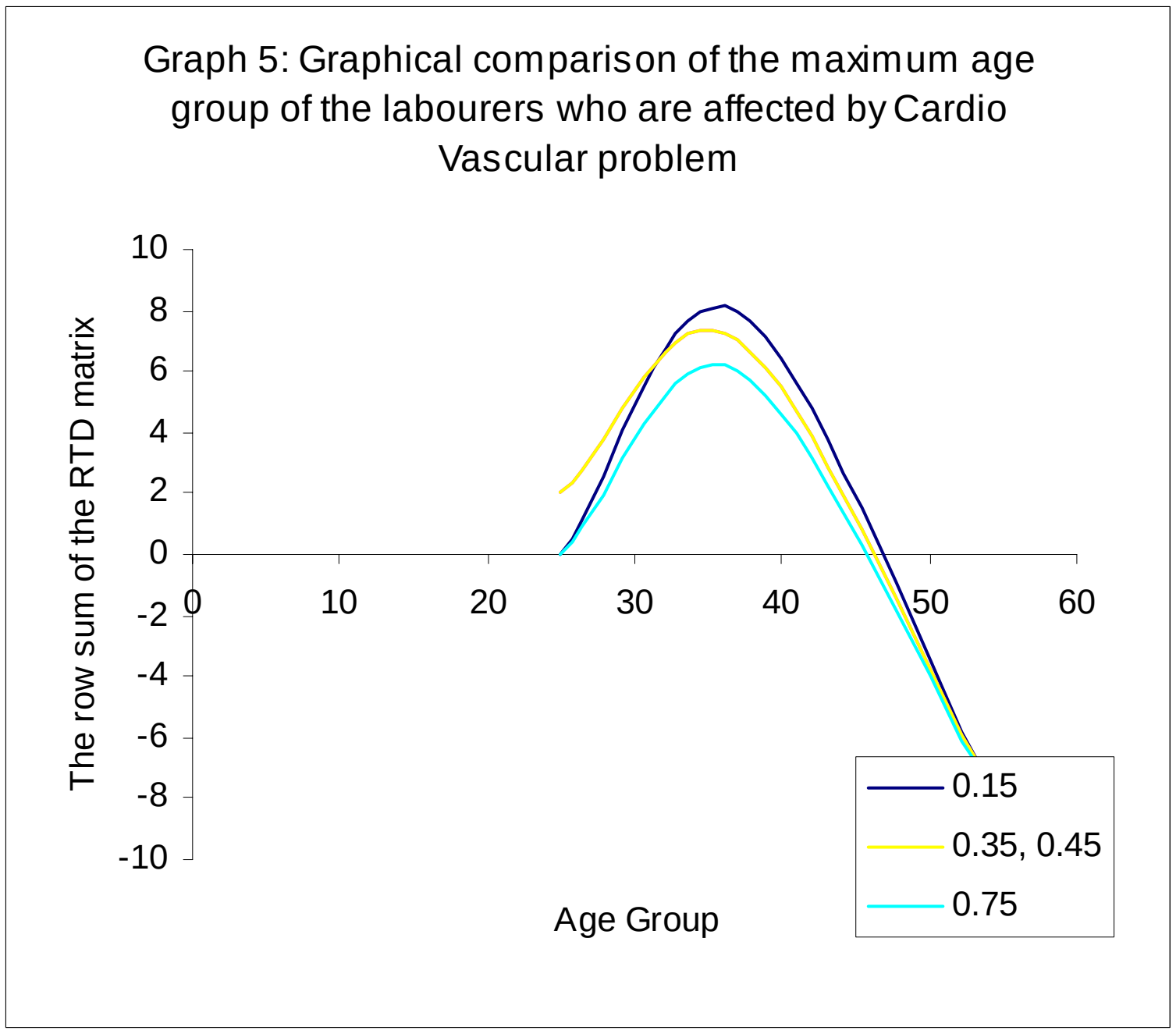
Graph 5: Graphical comparison of the maximum age group of the labourers who are affected by Cardio Vascular problem
The row sum of the RTD matrix
Age Group
10
8
6
4
2
0
-2
-4
-6
-8
-10
0
10
20
30
40
50
60
0.15
0.35, 0.45
0.75

From the above analysis, we observe that the maximum age group getting cardiovascular problem has not changed with the change in the value of the parameter from 0 to 1. The mathematical inference is that the maximum age group of agricultural labourer to have cardio vascular problem is 35-37, and the Combined Effect Time Dependent data matrix also confirm the same result.

The CETD Matrix — Row sum matrix

$$\begin{bmatrix} 0 & 0 & 0 & -1 & 3 & 0 & -1 & 3 \\ 4 & 4 & 4 & 4 & 4 & 4 & 4 & 1 \\ -4 & -4 & -4 & -4 & -4 & -4 & -4 & -4 \end{bmatrix} \qquad \begin{bmatrix} 4 \\ 29 \\ -32 \end{bmatrix}$$

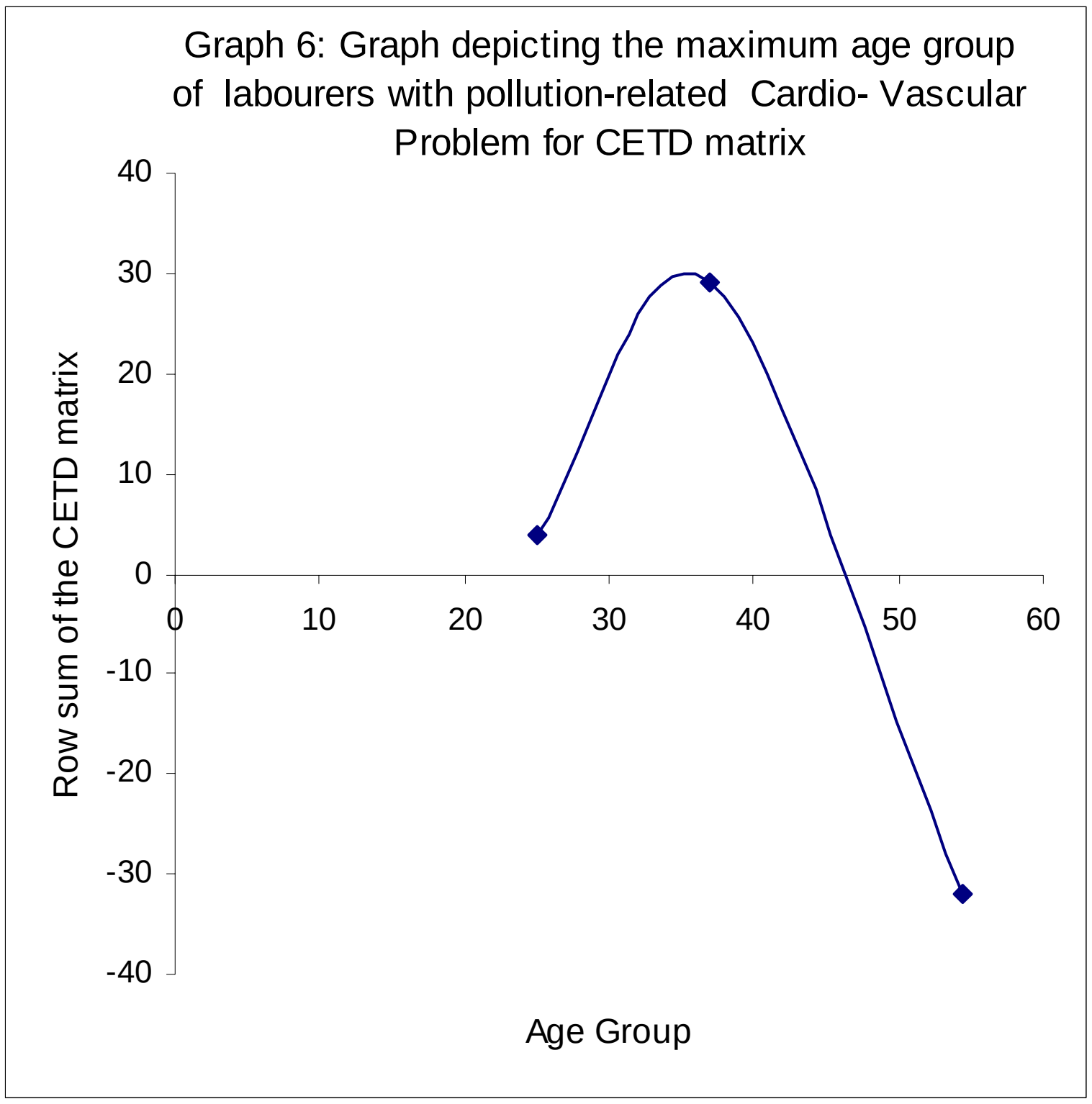

### 2.1.1.2 Estimation of maximum age group, using 4 × 8 matrices

Now to make the study more sensitive we increase the number of rows by 4 and see whether the decision arrived is more sensitive to the earlier one we have discussed. Thus we give the raw data of 4 × 8 matrix.

Initial raw data matrix of cardio vascular problem of order 4 × 8

| Years | $S_1$ | $S_2$ | $S_3$ | $S_4$ | $S_5$ | $S_6$ | $S_7$ | $S_8$ |
|---|---|---|---|---|---|---|---|---|
| 20-30 | 23 | 18 | 24 | 16 | 29 | 10 | 16 | 10 |
| 31-36 | 18 | 15 | 20 | 14 | 17 | 6 | 13 | 5 |
| 37-43 | 20 | 17 | 18 | 17 | 18 | 12 | 20 | 5 |
| 44-65 | 22 | 21 | 21 | 22 | 20 | 11 | 20 | 5 |

The ATD Matrix for cardio vascular problem of order 4 × 8

| Years | $S_1$ | $S_2$ | $S_3$ | $S_4$ | $S_5$ | $S_6$ | $S_7$ | $S_8$ |
|---|---|---|---|---|---|---|---|---|
| 20-30 | 2.09 | 1.64 | 2.18 | 1.46 | 2.64 | 0.91 | 1.46 | 0.91 |
| 31-36 | 3 | 2.5 | 3.33 | 2.33 | 2.83 | 1 | 2.17 | 0.83 |
| 37-43 | 2.86 | 2.43 | 2.57 | 2.43 | 2.57 | 1.71 | 2.86 | 0.71 |
| 44-65 | 1 | 0.95 | 0.95 | 1 | 0.91 | 0.5 | 0.91 | 0.23 |

The Average and Standard Deviation for the above ATD matrix

| | | | | | | | | |
|---|---|---|---|---|---|---|---|---|
| Average | 2.24 | 1.88 | 2.26 | 1.81 | 2.24 | 1.03 | 1.85 | 0.67 |
| SD | 0.92 | 0.73 | 0.99 | 0.69 | 0.89 | 0.50 | 0.85 | 0.30 |

RTD matrix for $\alpha = 0.15$ Row sum matrix

$$\begin{bmatrix} -1 & -1 & 0 & -1 & 1 & -1 & -1 & 1 \\ 1 & 1 & 1 & 1 & 1 & 0 & 1 & 1 \\ 1 & 1 & 1 & 1 & 1 & 1 & 1 & 0 \\ -1 & -1 & -1 & -1 & -1 & -1 & -1 & -1 \end{bmatrix} \qquad \begin{bmatrix} -3 \\ 7 \\ 7 \\ -8 \end{bmatrix}$$

RTD matrix for $\alpha = 0.35$ Row sum matrix

$$\begin{bmatrix} 0 & 0 & 0 & -1 & 1 & 0 & 1 & 1 \\ 1 & 1 & 1 & 1 & 1 & 0 & 1 & 1 \\ 1 & 1 & 0 & 1 & 1 & 1 & 1 & 0 \\ 1 & -1 & -1 & -1 & -1 & -1 & 0 & -1 \end{bmatrix} \qquad \begin{bmatrix} 2 \\ 7 \\ 6 \\ -5 \end{bmatrix}$$

RTD matrix for $\alpha = 0.45$ Row sum matrix

$$\begin{bmatrix} 0 & 0 & 0 & -1 & 0 & 0 & 0 & 1 \\ 1 & 1 & 1 & 1 & 1 & 0 & 0 & 1 \\ 1 & 1 & 0 & 1 & 0 & 1 & 1 & 0 \\ 1 & -1 & -1 & -1 & -1 & -1 & 0 & -1 \end{bmatrix} \qquad \begin{bmatrix} 0 \\ 6 \\ 5 \\ -5 \end{bmatrix}$$

RTD matrix for $\alpha = 0.75$ Row sum matrix

$$\begin{bmatrix} 0 & 0 & 0 & 0 & 0 & 0 & 0 & 1 \\ 1 & 1 & 1 & 0 & 0 & 0 & 0 & 0 \\ 0 & 0 & 0 & 1 & 0 & 1 & 1 & 0 \\ -1 & -1 & -1 & -1 & -1 & -1 & 0 & -1 \end{bmatrix} \qquad \begin{bmatrix} 1 \\ 3 \\ 3 \\ -7 \end{bmatrix}$$

Graph 7: Graph depicting maximum age group of labourers with pollution related Cardio Vascular Problem for $\alpha = 0.15$

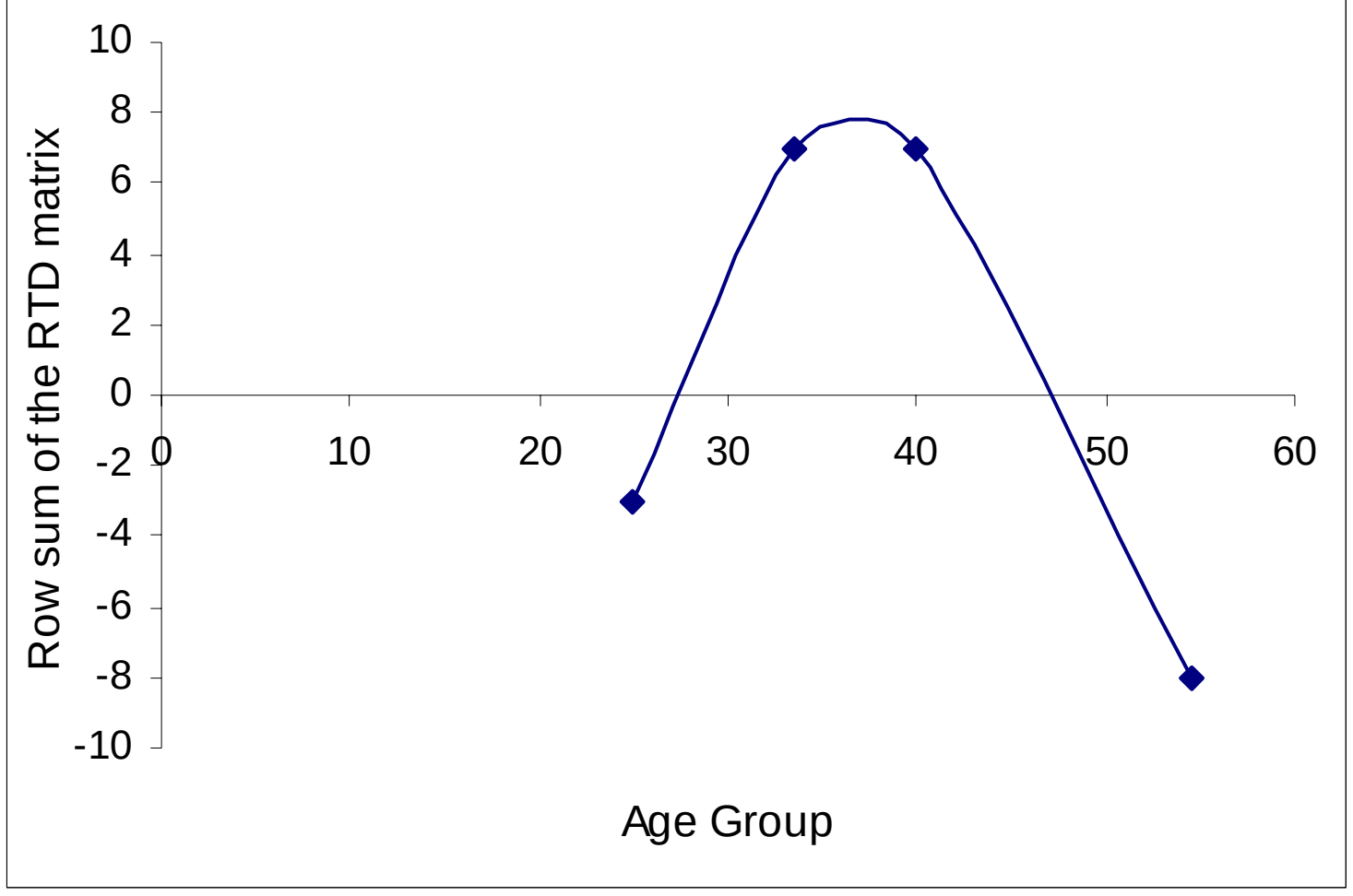

Graph 8: Graph depicting maximum age group of labourers with pollution related Cardio Vascular Problem for $\alpha = 0.35$

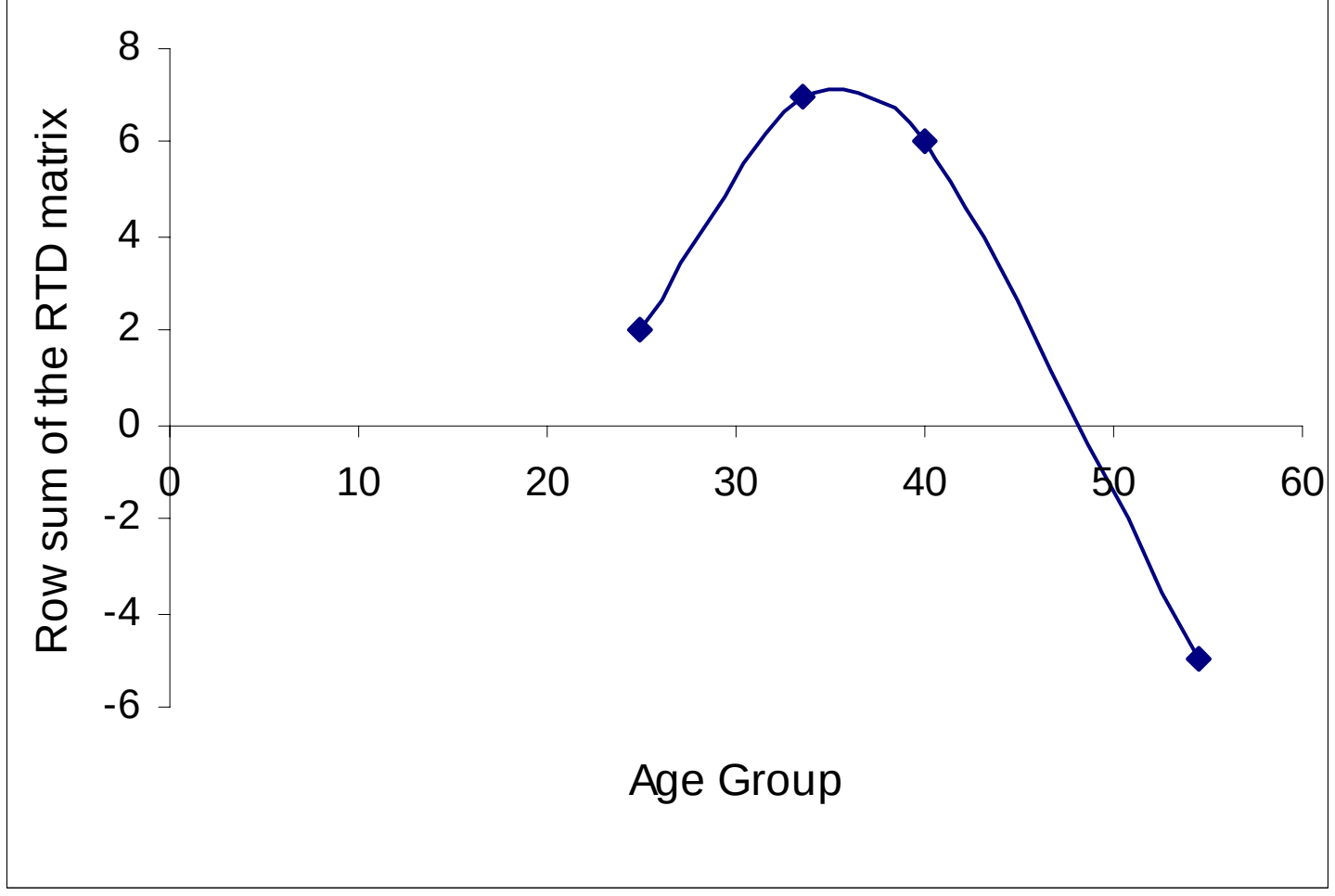

Graph 9 : Graph depicting maximum age group of labourers with pollution related Cardio Vascular Problem for $\alpha = 0.45$

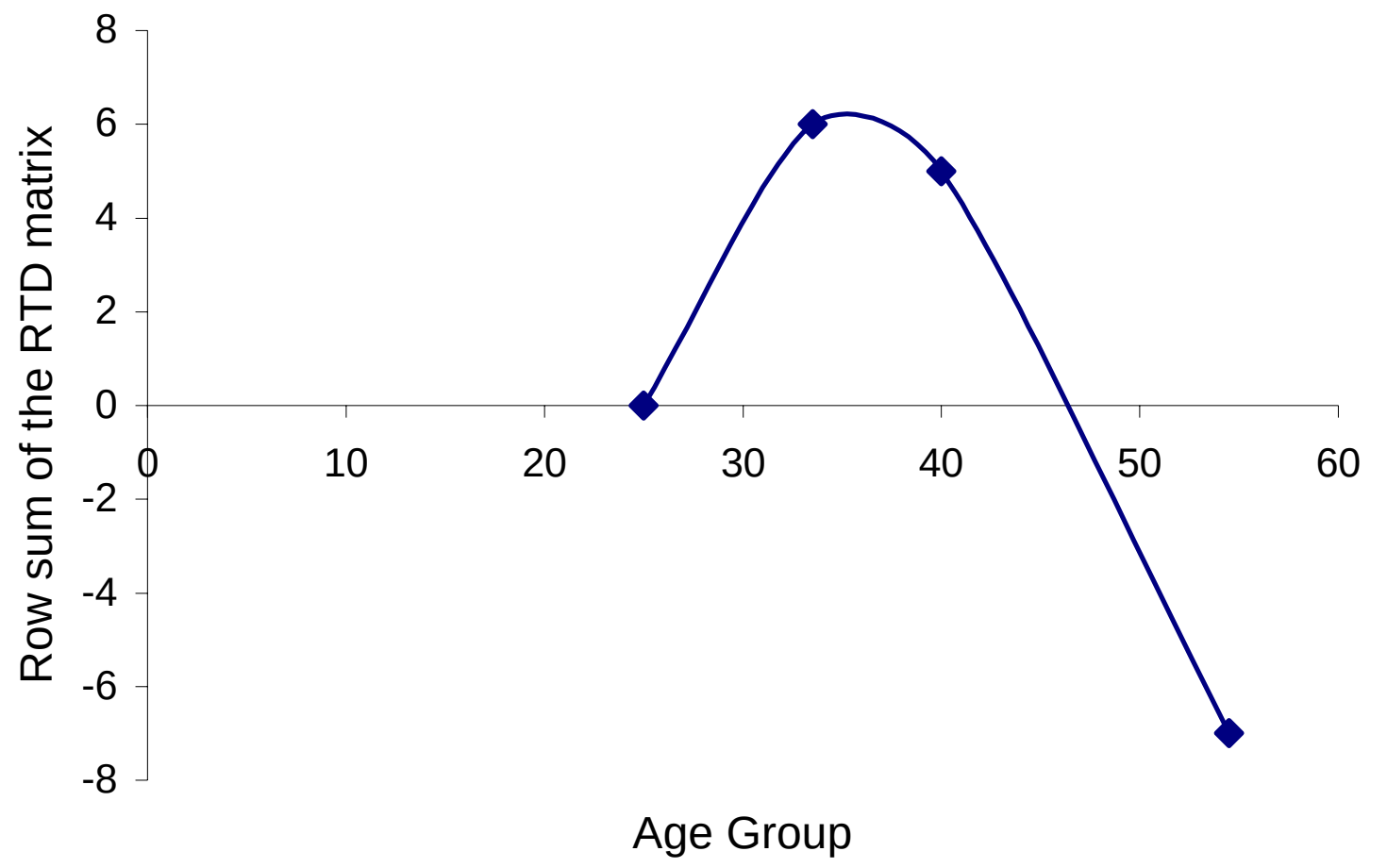

Graph 10 : Graph depicting maximum age group of labourers with pollution related Cardio Vascular Problem for $\alpha = 0.75$

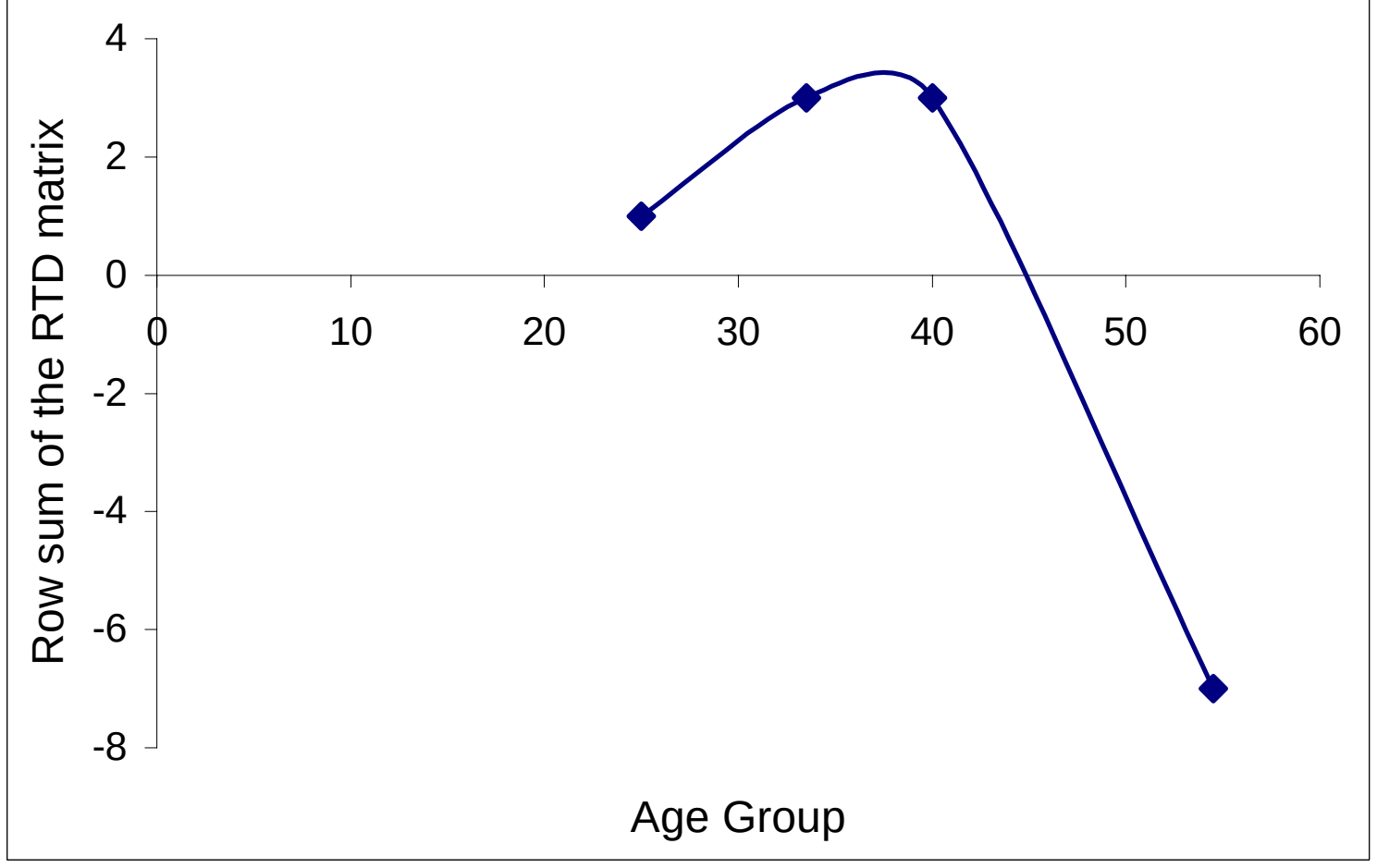

We now give the comparative graph of the maximum age group of labourers with pollution related Cardio Vascular Problem.

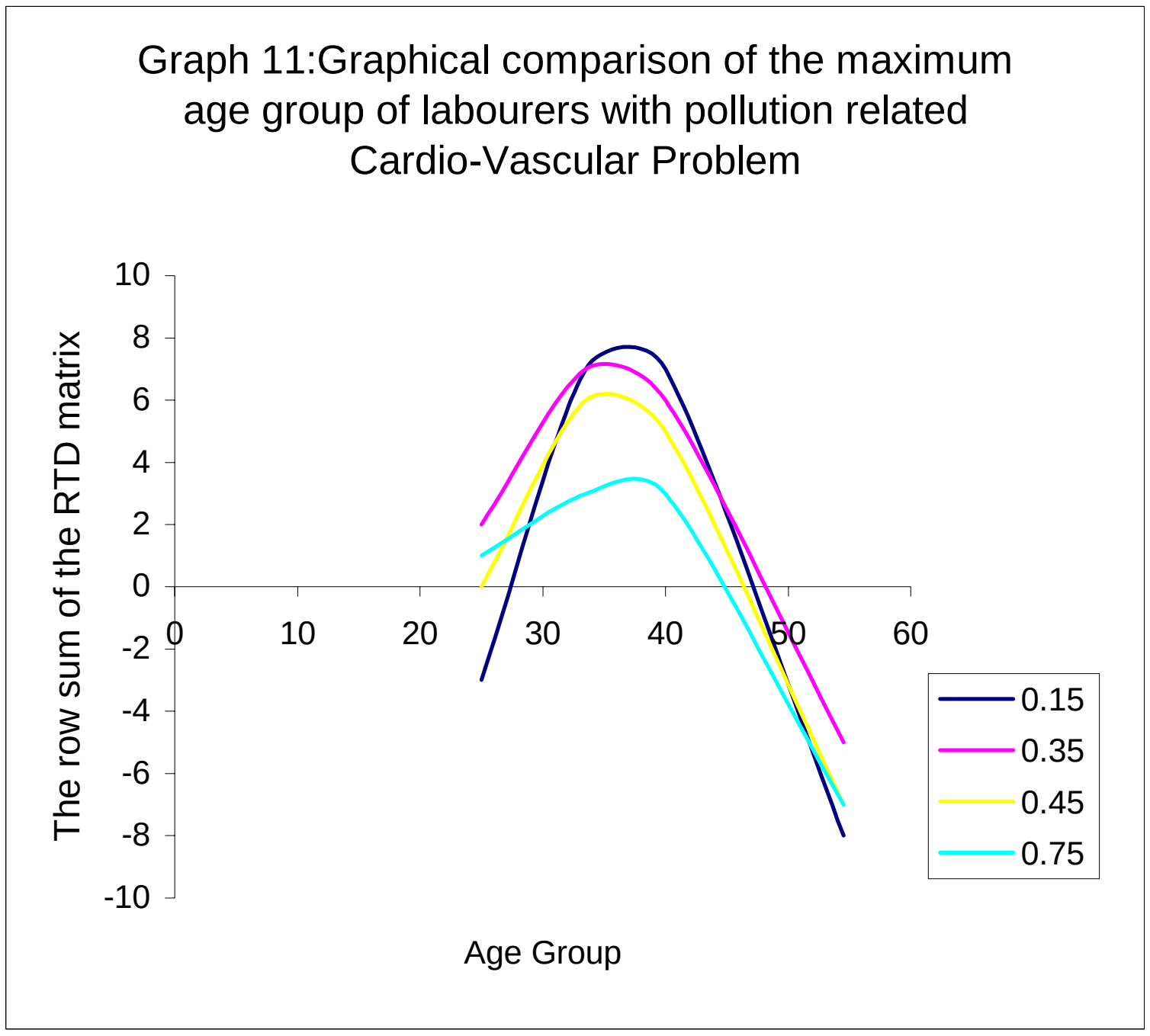

The CETD matrix and the row sum matrix are given here:

CETD Matrix　　　　Row sum matrix

$$\begin{bmatrix} -1 & -1 & 0 & -3 & 2 & -1 & 0 & 4 \\ 4 & 4 & 4 & 3 & 3 & 0 & 2 & 3 \\ 3 & 3 & 1 & 4 & 2 & 4 & 4 & 0 \\ -2 & -4 & -4 & -4 & -4 & -4 & -1 & -2 \end{bmatrix} \qquad \begin{bmatrix} 0 \\ 23 \\ 21 \\ -25 \end{bmatrix}$$

Now we give the related graph for the CETD matrix :

Graph 12: The graph depicting the maximum age group of labourers who are affected by Cardio Vascular Problem for CETD matrix

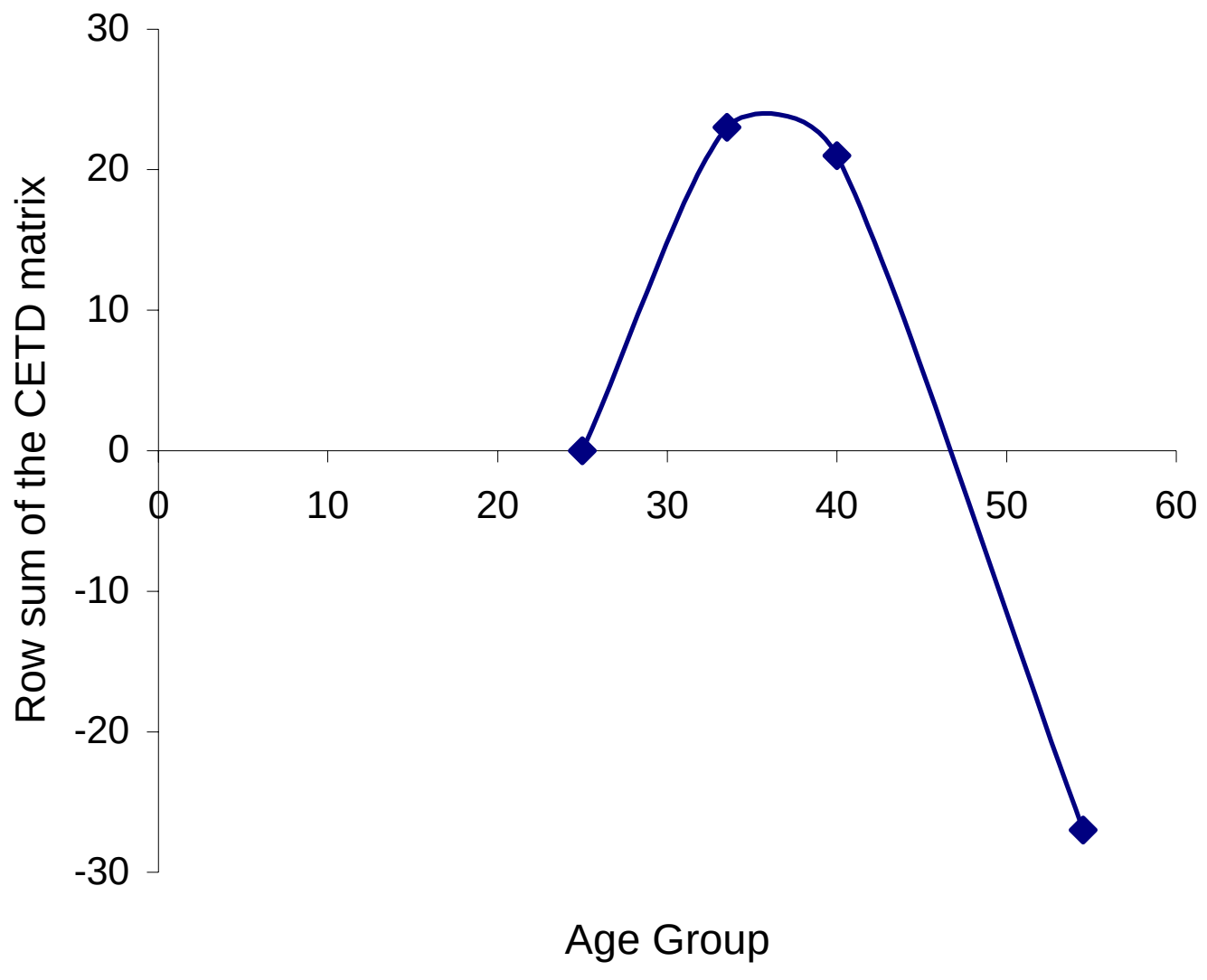

We observe from the above graph that

1. The cardio vascular problem starts only at the age of 25
2. The maximum age for getting cardio vascular problem is 33-36
3. The peak period for a heart problem is 35
4. The above three results also confirmed from the CETD matrix

**2.1.1.3 Estimation of maximum age group of agriculture labourer with pollution related cardio vascular problem, using 5 × 8 matrices**

Now to make the study more sensitive we increase the number of rows by 5 and see whether the decision arrived is more sensitive to the earlier one.

Initial Raw Data Matrix of
Cardio Vascular Problem of order $5 \times 8$

| Years | $S_1$ | $S_2$ | $S_3$ | $S_4$ | $S_5$ | $S_6$ | $S_7$ | $S_8$ |
|---|---|---|---|---|---|---|---|---|
| 20-24 | 2 | 1 | 3 | 1 | 3 | 0 | 2 | 1 |
| 25-30 | 21 | 17 | 21 | 15 | 26 | 10 | 14 | 9 |
| 31-36 | 18 | 15 | 20 | 14 | 17 | 6 | 13 | 5 |
| 37-43 | 20 | 17 | 18 | 17 | 18 | 12 | 20 | 5 |
| 44-65 | 22 | 21 | 21 | 22 | 20 | 11 | 20 | 5 |

The ATD Matrix for
Cardio Vascular problem of order $5 \times 8$

| Years | $S_1$ | $S_2$ | $S_3$ | $S_4$ | $S_5$ | $S_6$ | $S_7$ | $S_8$ |
|---|---|---|---|---|---|---|---|---|
| 20-24 | 0.4 | 0.2 | 0.6 | 0.2 | 0.6 | 0 | 0.4 | 0.2 |
| 25-30 | 3.5 | 2.83 | 3.5 | 2.5 | 4.33 | 1.67 | 2.33 | 1.5 |
| 31-36 | 3 | 2.5 | 3.33 | 2.33 | 2.83 | 1 | 2.17 | 0.83 |
| 37-43 | 2.86 | 2.43 | 2.57 | 2.43 | 2.57 | 1.71 | 2.86 | 0.71 |
| 44-65 | 1 | 0.95 | 0.95 | 1 | 0.91 | 0.5 | 0.91 | 0.23 |

Average and the Standard Deviation
of the above given ATD matrix

| | | | | | | | | |
|---|---|---|---|---|---|---|---|---|
| Average | 2.15 | 1.78 | 2.19 | 1.69 | 2.25 | 0.98 | 1.73 | 0.69 |
| SD | 1.36 | 1.14 | 1.34 | 1.03 | 1.52 | 1.74 | 1.03 | 0.53 |

RTD matrix for $\alpha = 0.1$ Row sum matrix

$$\begin{bmatrix} -1 & -1 & -1 & -1 & -1 & -1 & -1 & -1 \\ 1 & 1 & 1 & 1 & 1 & -1 & 1 & 1 \\ 1 & 1 & 1 & 1 & 1 & 0 & 1 & 1 \\ 1 & 1 & 1 & 1 & 1 & 1 & 1 & 0 \\ -1 & 1 & -1 & -1 & -1 & -1 & -1 & -1 \end{bmatrix} \qquad \begin{bmatrix} -8 \\ 6 \\ 7 \\ 7 \\ -6 \end{bmatrix}$$

RTD matrix for $\alpha = 0.15$ Row sum matrix

$$\begin{bmatrix} -1 & -1 & -1 & -1 & -1 & -1 & -1 & -1 \\ 1 & 1 & 1 & 1 & 1 & 1 & 1 & 1 \\ 1 & 1 & 1 & 1 & 1 & 0 & 1 & 1 \\ 1 & 1 & 1 & 1 & 1 & 1 & 1 & 0 \\ -1 & 0 & -1 & -1 & -1 & -1 & -1 & -1 \end{bmatrix} \qquad \begin{bmatrix} -8 \\ 8 \\ 7 \\ 7 \\ -7 \end{bmatrix}$$

RTD matrix for $\alpha = 0.2$ Row sum matrix

$$\begin{bmatrix} -1 & -1 & -1 & -1 & -1 & -1 & -1 & -1 \\ 1 & -1 & 1 & 1 & 1 & 1 & 1 & 1 \\ 1 & 1 & 1 & 1 & 1 & 0 & 1 & 1 \\ 1 & 1 & 1 & 1 & 1 & 1 & 1 & 0 \\ -1 & -1 & -1 & -1 & -1 & -1 & -1 & -1 \end{bmatrix} \qquad \begin{bmatrix} -8 \\ 6 \\ 7 \\ 7 \\ -8 \end{bmatrix}$$

RTD matrix for $\alpha = 0.35$ Row sum matrix

$$\begin{bmatrix} -1 & -1 & -1 & -1 & -1 & -1 & -1 & -1 \\ 1 & 1 & 1 & 1 & 1 & 1 & 1 & 1 \\ 1 & 1 & 1 & 1 & 1 & 0 & 1 & 0 \\ 1 & 1 & 0 & 1 & 0 & 1 & 1 & 0 \\ -1 & -1 & -1 & -1 & -1 & -1 & -1 & -1 \end{bmatrix} \qquad \begin{bmatrix} -8 \\ 8 \\ 6 \\ 5 \\ -8 \end{bmatrix}$$

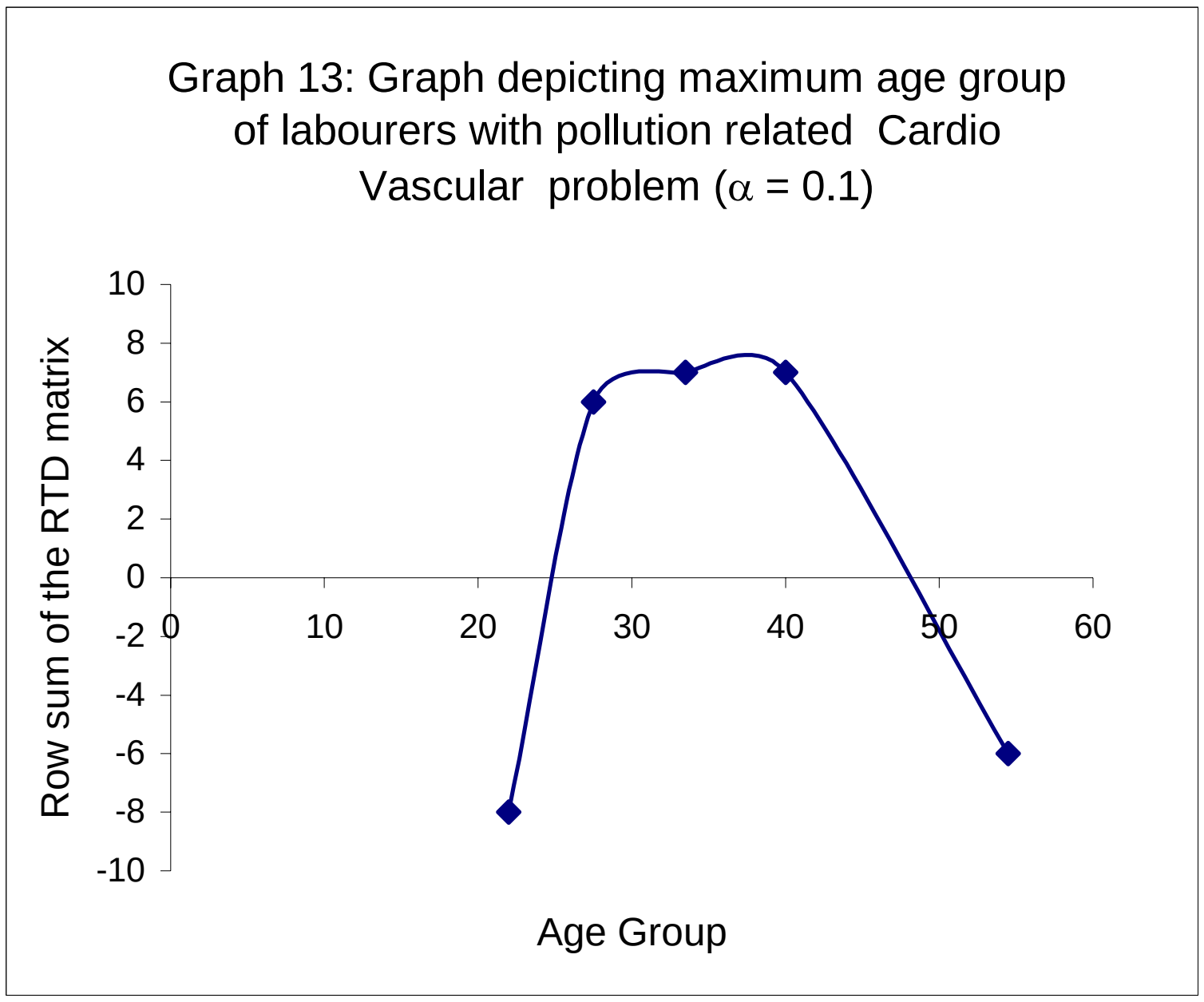
Graph 13: Graph depicting maximum age group of labourers with pollution related Cardio Vascular problem (α = 0.1)
Row sum of the RTD matrix
Age Group
10
8
6
4
2
0
-2
-4
-6
-8
-10
0
10
20
30
40
50
60

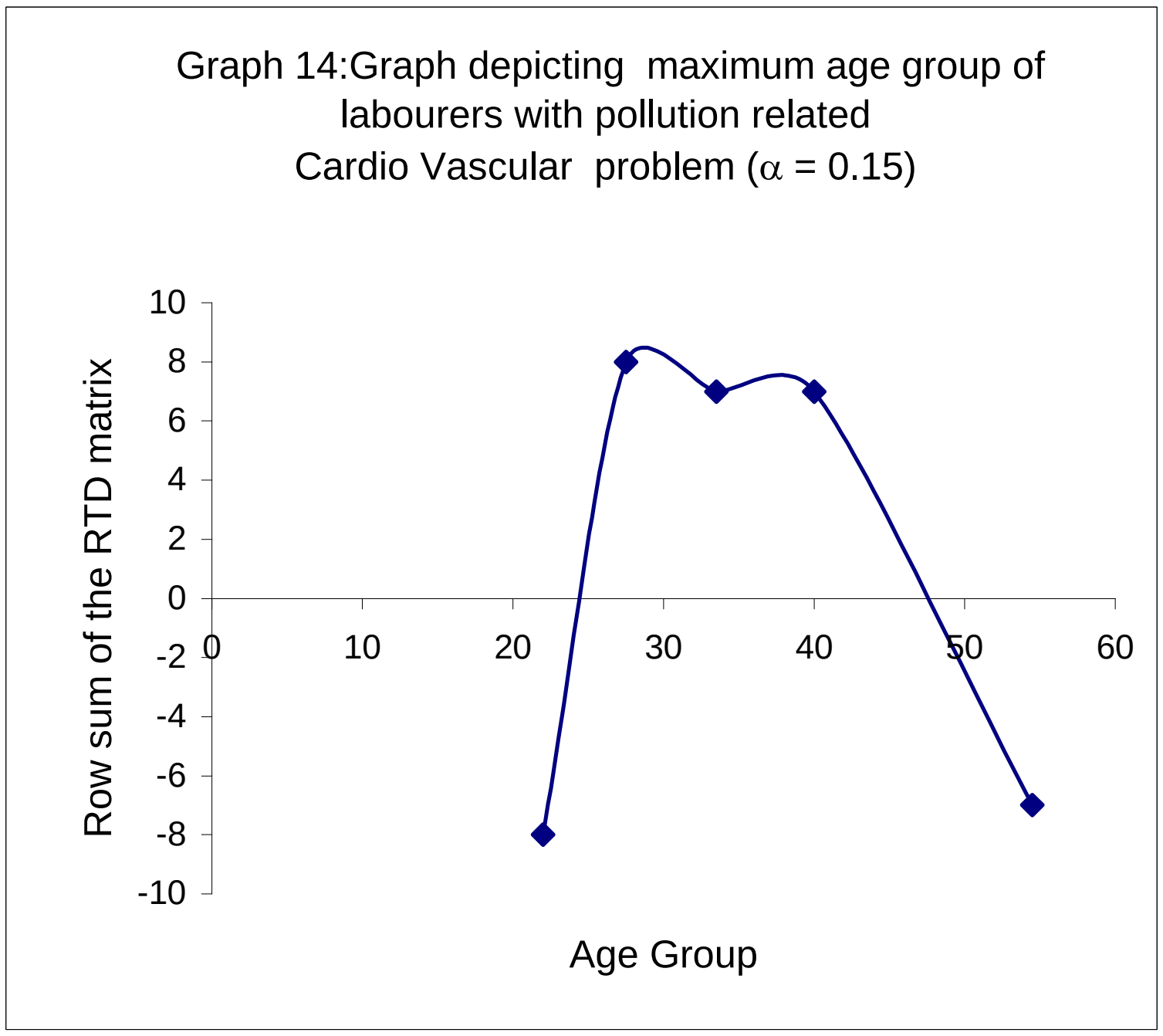
Graph 14:Graph depicting maximum age group of labourers with pollution related Cardio Vascular problem (α = 0.15)
Row sum of the RTD matrix
Age Group
10
8
6
4
2
0
-2
-4
-6
-8
-10
0
10
20
30
40
50
60

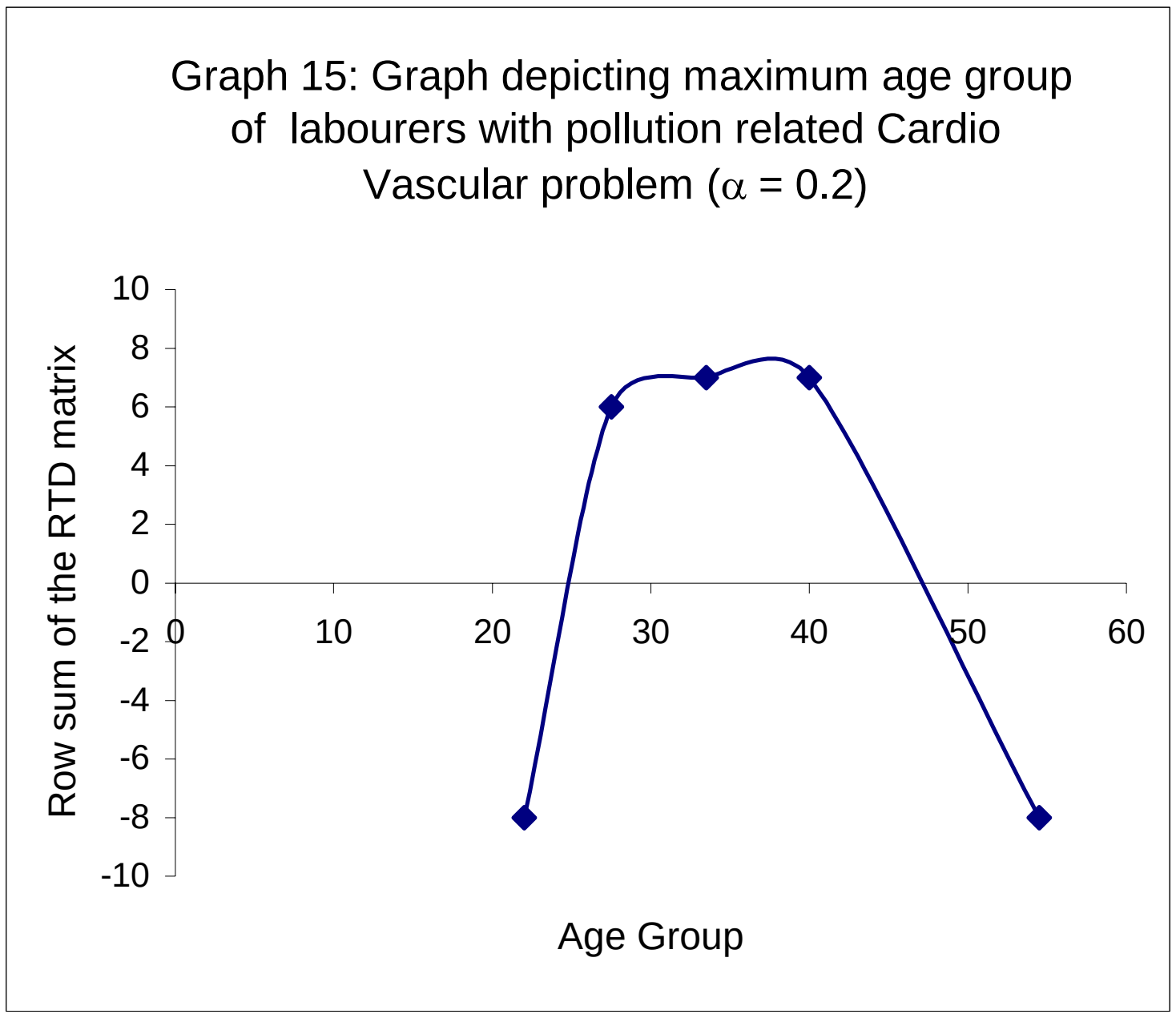
Graph 15: Graph depicting maximum age group of labourers with pollution related Cardio Vascular problem (α = 0.2)
Row sum of the RTD matrix
Age Group
10
8
6
4
2
0
-2
-4
-6
-8
-10
0
10
20
30
40
50
60

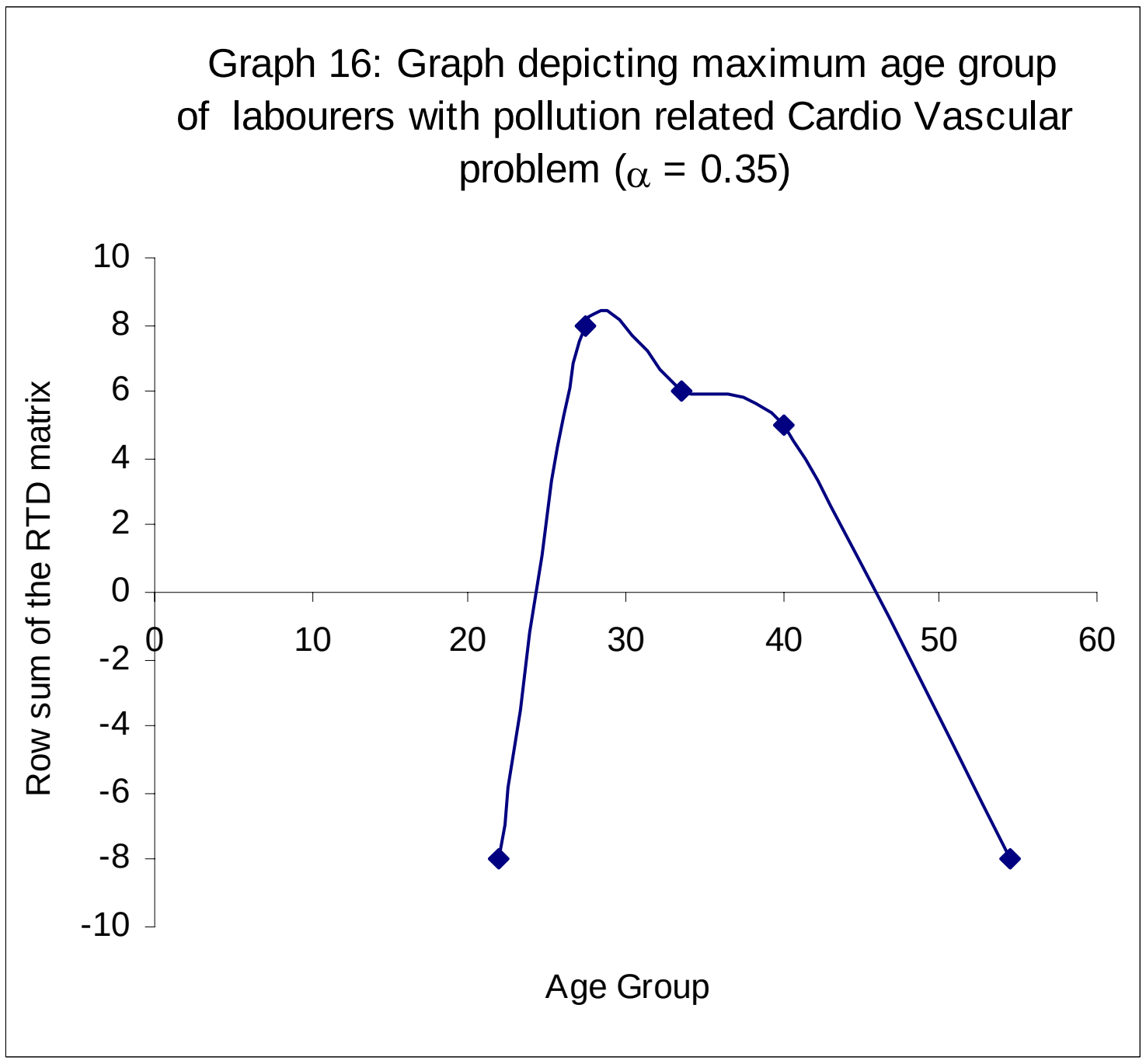
Graph 16: Graph depicting maximum age group of labourers with pollution related Cardio Vascular problem (α = 0.35)
Row sum of the RTD matrix
Age Group
10
8
6
4
2
0
-2
-4
-6
-8
-10
0
10
20
30
40
50
60

The comparative graph of the maximum age group of labourers who are affected by Cardio Vascular Problem is given below:

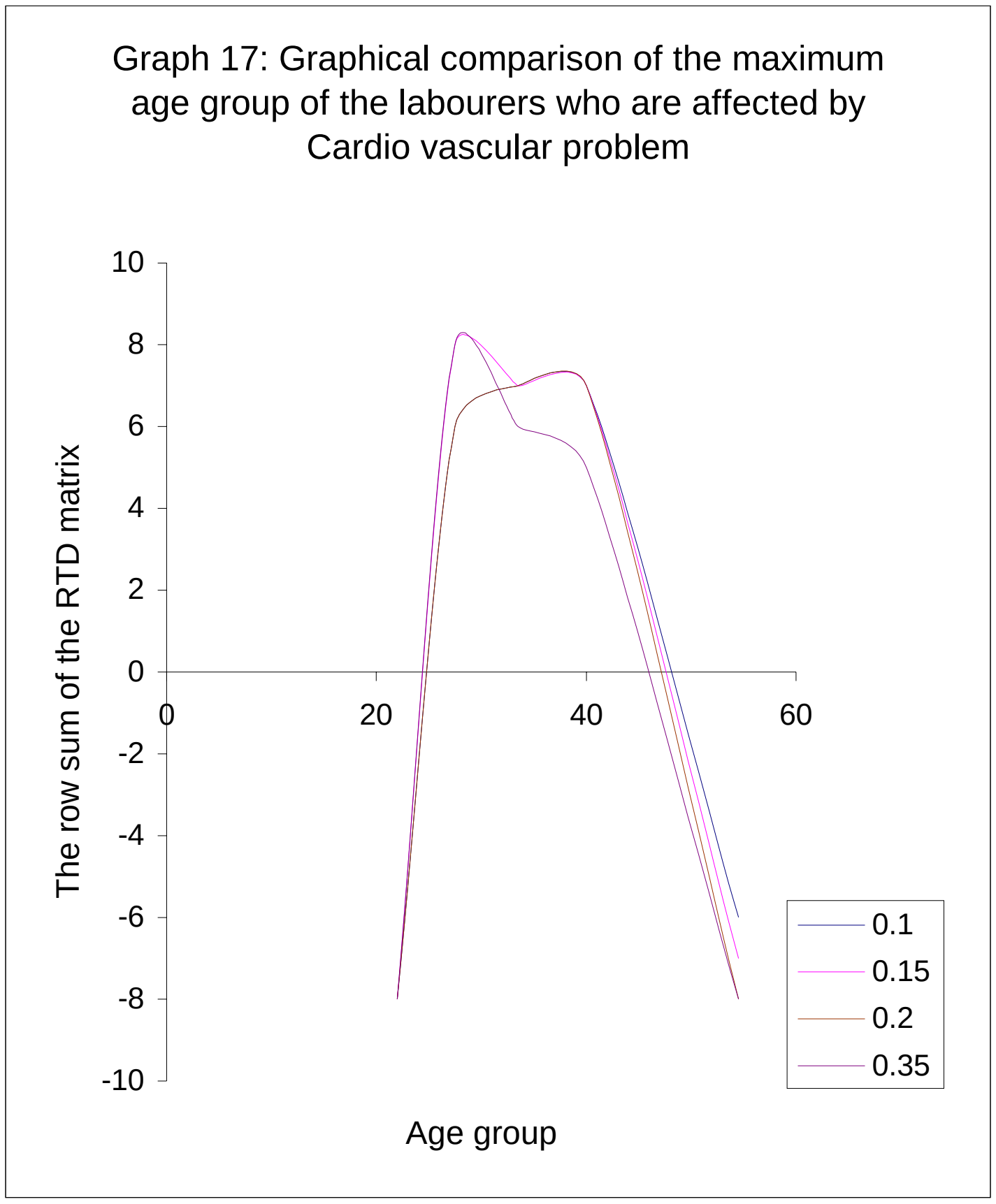

CETD matrix

$$\begin{bmatrix} -4 & -4 & -4 & -4 & -4 & -4 & -4 & -4 \\ 4 & 2 & 4 & 4 & 4 & 2 & 4 & 4 \\ 4 & 3 & 4 & 4 & 4 & 0 & 4 & 3 \\ 4 & 4 & 3 & 4 & 3 & 4 & 4 & 0 \\ -3 & -2 & -4 & -4 & -4 & -4 & -4 & -4 \end{bmatrix}$$

Row sum matrix

$$\begin{bmatrix} -32 \\ 28 \\ 26 \\ 26 \\ -29 \end{bmatrix}$$

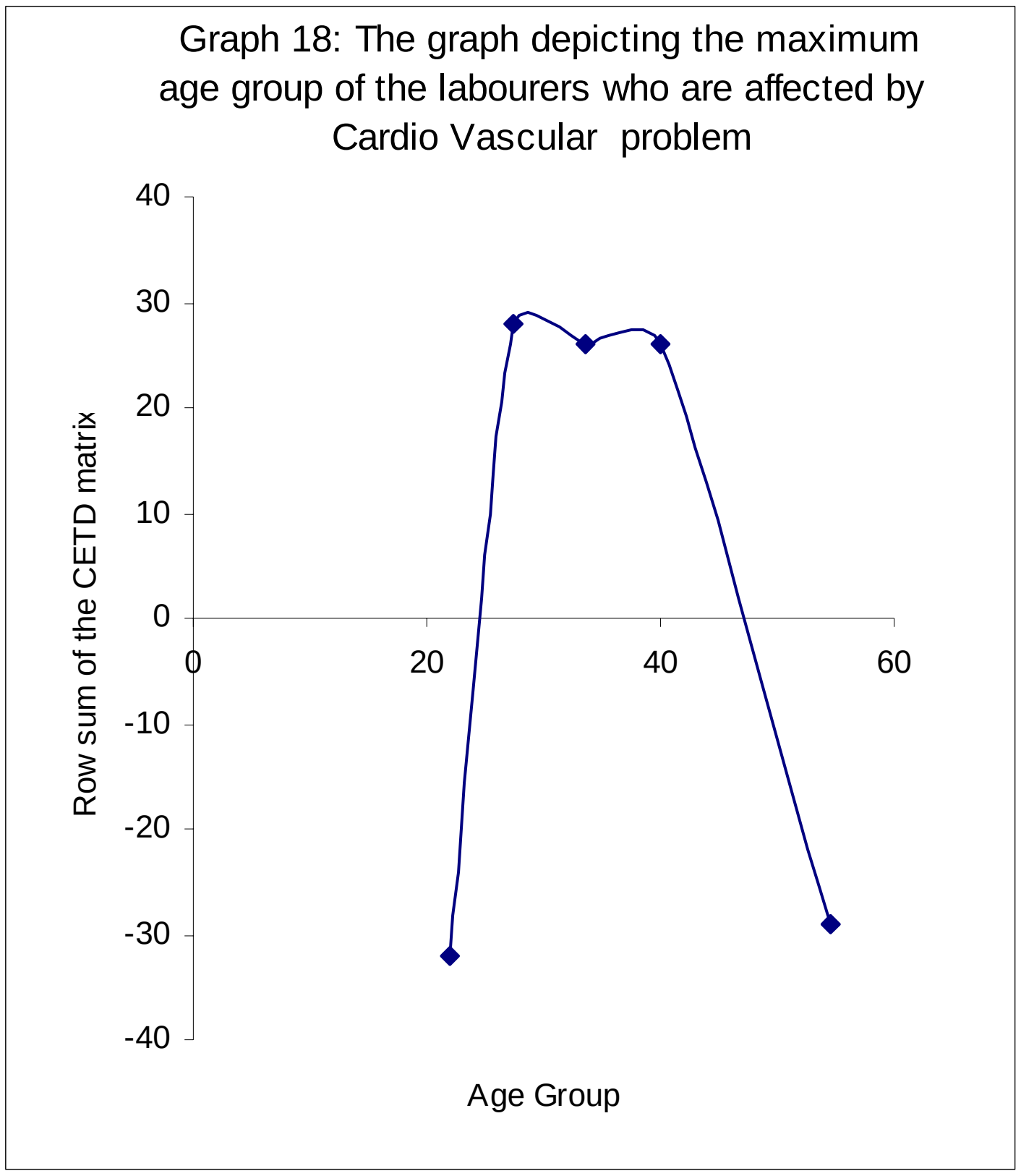

### 2.1.1.4 Conclusions

First contrary to the natural happening that a person after 40 has more chances of getting the cardio vascular disease, we see in case of these agricultural labourers who toil in the sun from day to dawn and who have no other mental tensions or hyper tensions become victims of all types of cardio vascular symptoms mainly due to the evil effects of chemicals used as the pesticides and insecticides. It is unfortunate to state that most of the agriculture labour suffer symptom of cardio vascular disease in the very early age group of 35-37 which is surprising. For there is no natural rhyme or reason for it for even after work hours they do not have traces of tensions or humiliation to be more precise they after a long days physical labour can have a peaceful sleep but on the contrary they suffer serious symptoms like B.P., chest pain burning heart/chest etc. Of the 110 agriculture labourers interviewed we saw 83 of them suffered from the cardio vascular symptom of which 38 of them were from the age group 31 to 43. Only 22 from the age group 41 to 65 which is against the natural medical laws for they should be in large number the two reasons for this have already been attributed. Also 33 of them in the age group 31 to 43 had burning chest symptoms also. They all attributed this mainly due to the pollution by chemicals. For uniformly all of them felt that after the spray of pesticides they suffered more of the symptoms say for at least a fortnight or more, due to the evil effects of the chemicals used in the pesticides and insecticides. It is unfortunate to state that most of the agriculture labourers suffer the symptom of cardio vascular disease in the very early age group of 35-37, which is surprising. As people who live after sixty are 2% and in fifties are less than the 10% we see people suffering from cardio vascular disease and the related symptoms in the age group 44-65 are negligible seen from the large negative deviation. Also the age group 20-24 is negative for this age group forms the migrants due to agriculture failure they are saved from the clutches of the chemical pollution. This is easily seen from the graphs.

### 2.1.2 Estimation of the Maximum Age Group of the Agriculture Labourers having Digestive Problems due to Chemical Pollution, using Matrices

Digestive problem is taken under six symptom diseases viz.

$S_1$ - No appetite,
$S_2$ - Vomiting,
$S_3$ - Diarrhea,
$S_4$ - Constipation,
$S_5$ - Mouth-Stomach Ulcer and
$S_6$ - Indigestion .

The symptoms are taken along the rows and ages of the labourers are taken along the columns.

Initial Raw Data Matrix of Digestive Problem of order $3 \times 6$

| Years | $S_1$ | $S_2$ | $S_3$ | $S_4$ | $S_5$ | $S_6$ |
|---|---|---|---|---|---|---|
| 20-30 | 17 | 11 | 8 | 16 | 16 | 17 |
| 31-43 | 29 | 26 | 15 | 17 | 20 | 20 |
| 44-65 | 22 | 15 | 11 | 9 | 12 | 15 |

ATD Matrix of Digestive Problem of order $3 \times 6$

| Years | $S_1$ | $S_2$ | $S_3$ | $S_4$ | $S_5$ | $S_6$ |
|---|---|---|---|---|---|---|
| 20-30 | 1.55 | 1 | 0.73 | 1.45 | 1.45 | 1.55 |
| 31-43 | 2.23 | 2 | 1.15 | 1.31 | 1.54 | 1.54 |
| 44-65 | 1 | 0.68 | 0.5 | 0.4 | 0.55 | 0.68 |

The Average and Standard Deviation of the above ATD matrix

| | | | | | | |
|---|---|---|---|---|---|---|
| Average | 1.59 | 1.23 | 0.79 | 1.05 | 1.18 | 1.26 |
| Standard Deviation | 0.62 | 0.69 | 0.33 | 0.57 | 0.55 | 0.50 |

RTD matrix for $\alpha = 0.15$     Row sum matrix

$$\begin{bmatrix} 0 & -1 & -1 & 1 & 1 & 1 \\ 1 & 1 & 1 & 1 & 1 & 1 \\ -1 & -1 & -1 & -1 & -1 & -1 \end{bmatrix} \qquad \begin{bmatrix} 1 \\ 6 \\ -6 \end{bmatrix}$$

RTD matrix for $\alpha = 0.35$     Row sum matrix

$$\begin{bmatrix} 0 & 0 & 0 & 1 & 1 & 1 \\ 1 & 1 & 1 & 1 & 1 & 1 \\ -1 & -1 & -1 & -1 & -1 & -1 \end{bmatrix} \qquad \begin{bmatrix} 3 \\ 6 \\ -6 \end{bmatrix}$$

RTD matrix for $\alpha = 0.45$ — Row sum matrix

$$\begin{bmatrix} 0 & 0 & 0 & 1 & 1 & 1 \\ 1 & 1 & 1 & 0 & 1 & 1 \\ -1 & -1 & -1 & -1 & -1 & -1 \end{bmatrix} \qquad \begin{bmatrix} 3 \\ 5 \\ -6 \end{bmatrix}$$

RTD matrix for $\alpha = 0.75$ — Row sum matrix

$$\begin{bmatrix} 0 & 0 & 0 & 0 & 0 & 0 \\ 1 & 1 & 1 & 0 & 0 & 0 \\ -1 & -1 & -1 & -1 & -1 & -1 \end{bmatrix} \qquad \begin{bmatrix} 0 \\ 3 \\ -6 \end{bmatrix}$$

Graph 19: Graph depicting the maximum age group of labourers with pollution related digestive problem for $\alpha = 0.15$

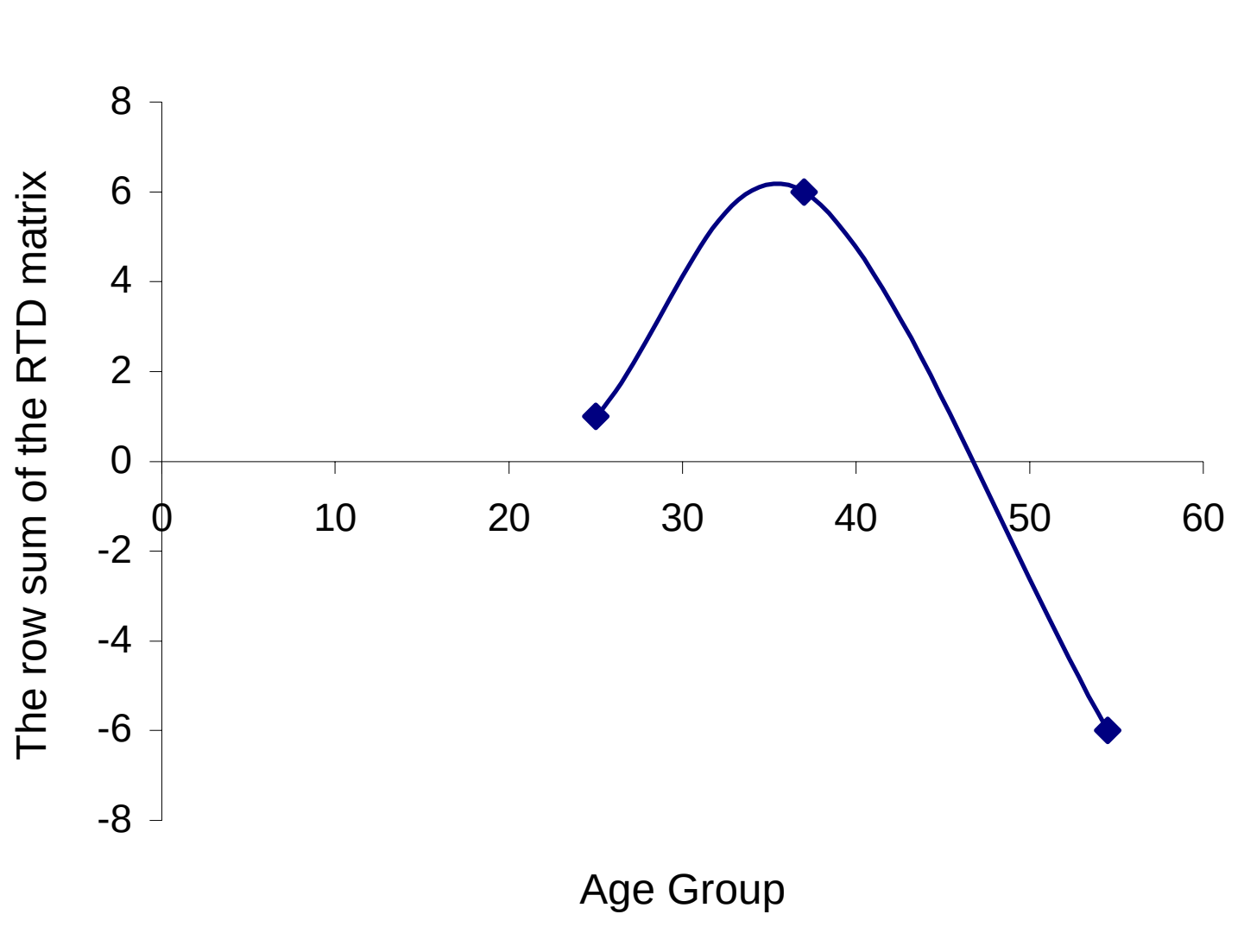

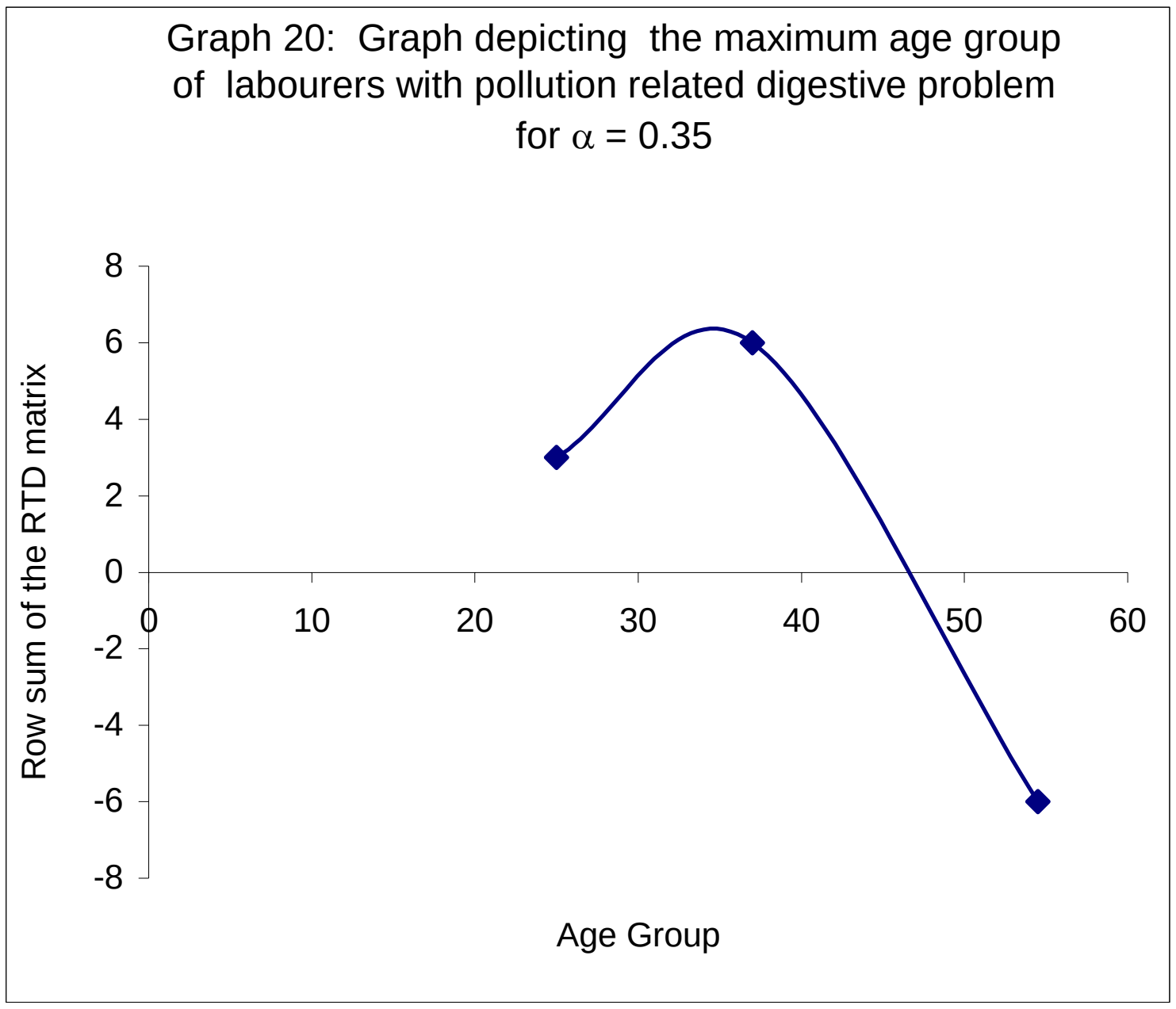

Graph 20: Graph depicting the maximum age group of labourers with pollution related digestive problem for α = 0.35
8
6
4
2
0
-2
-4
-6
-8
0
10
20
30
40
50
60
Row sum of the RTD matrix
Age Group

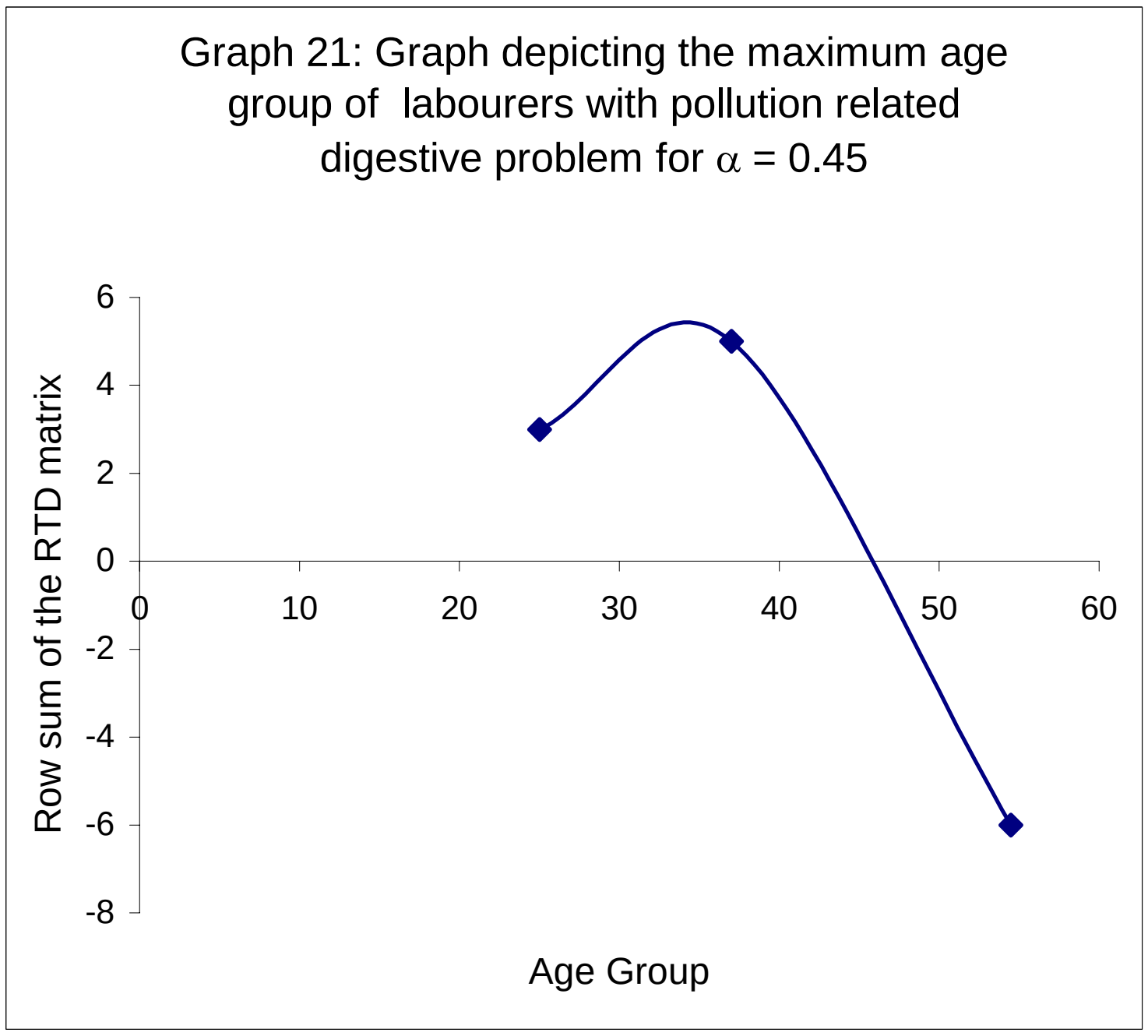

Graph 21: Graph depicting the maximum age group of labourers with pollution related digestive problem for α = 0.45
6
4
2
0
-2
-4
-6
-8
0
10
20
30
40
50
60
Row sum of the RTD matrix
Age Group

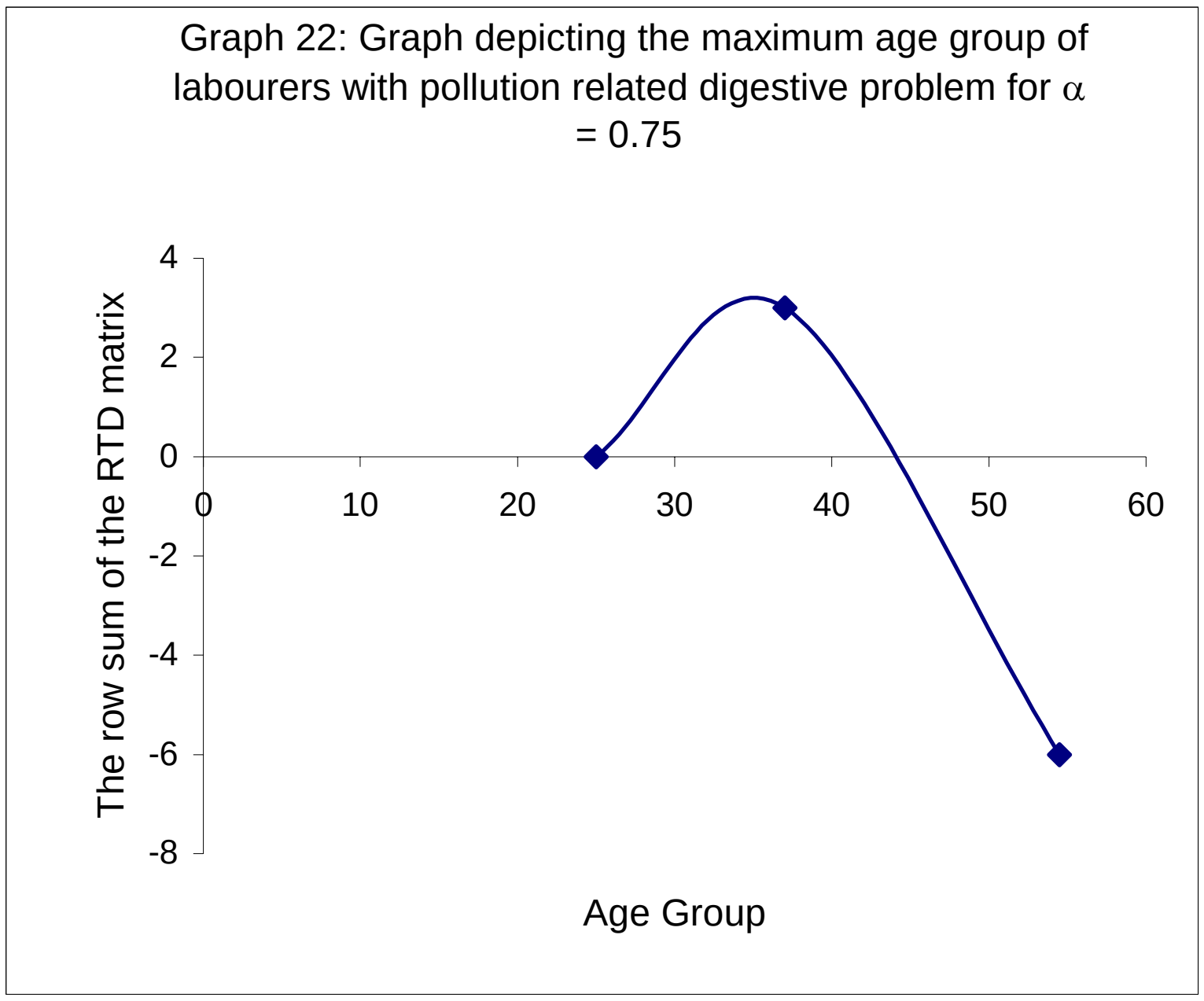
Graph 22: Graph depicting the maximum age group of labourers with pollution related digestive problem for α = 0.75
The row sum of the RTD matrix
Age Group
4
2
0
-2
-4
-6
-8
0
10
20
30
40
50
60

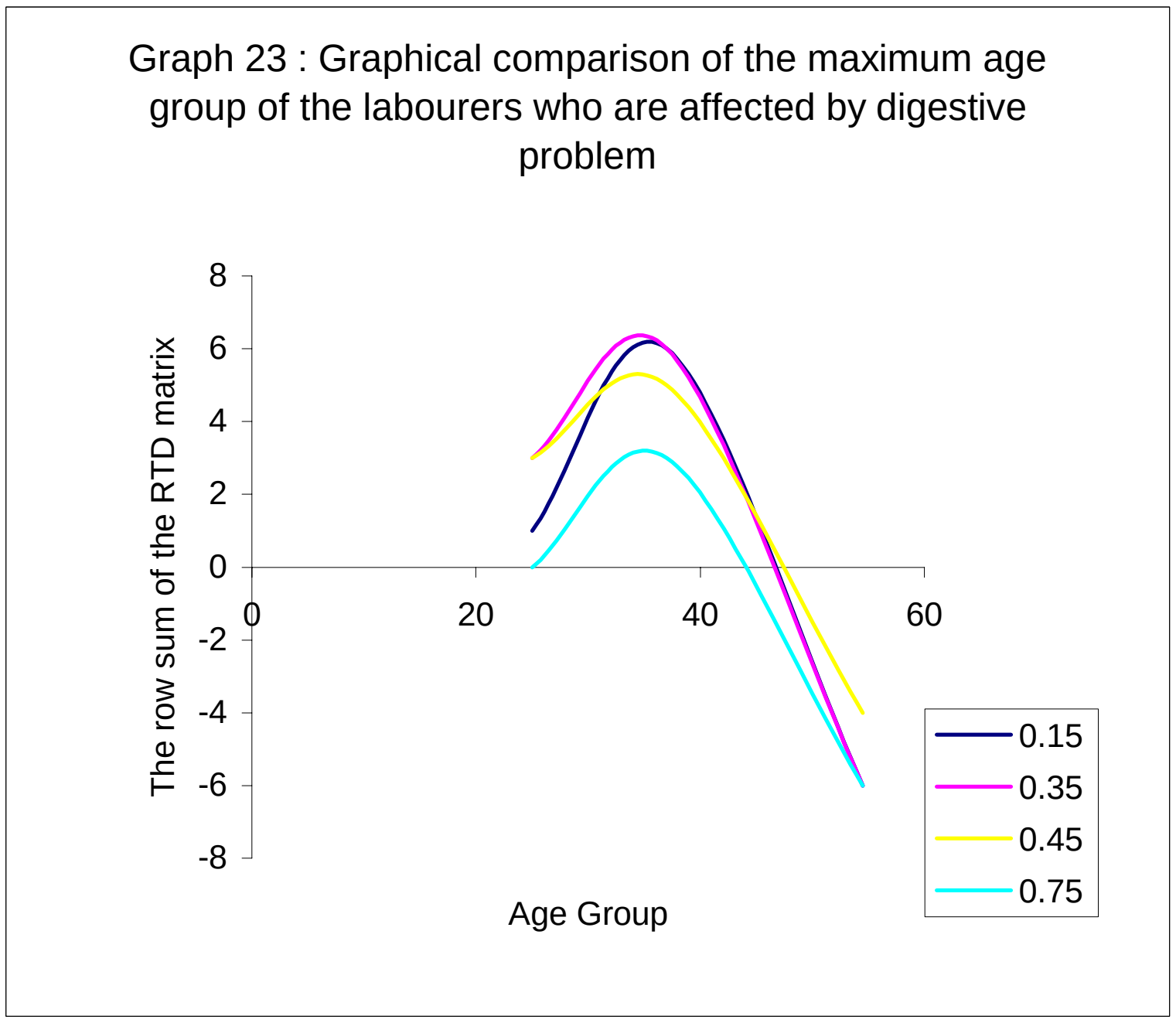
Graph 23 : Graphical comparison of the maximum age group of the labourers who are affected by digestive problem
The row sum of the RTD matrix
Age Group
8
6
4
2
0
-2
-4
-6
-8
0
20
40
60
0.15
0.35
0.45
0.75

CETD matrix Row sum matrix

$$\begin{bmatrix} 0 & -1 & -1 & 3 & 3 & 3 \\ 4 & 4 & 4 & 2 & 3 & 3 \\ -4 & -4 & -4 & -4 & -4 & -4 \end{bmatrix} \qquad \begin{bmatrix} 7 \\ 20 \\ -24 \end{bmatrix}$$

Graph 24 : The graph depicting the maximum age group of the labourers who are affected by the digestive problem for the CETD matrix

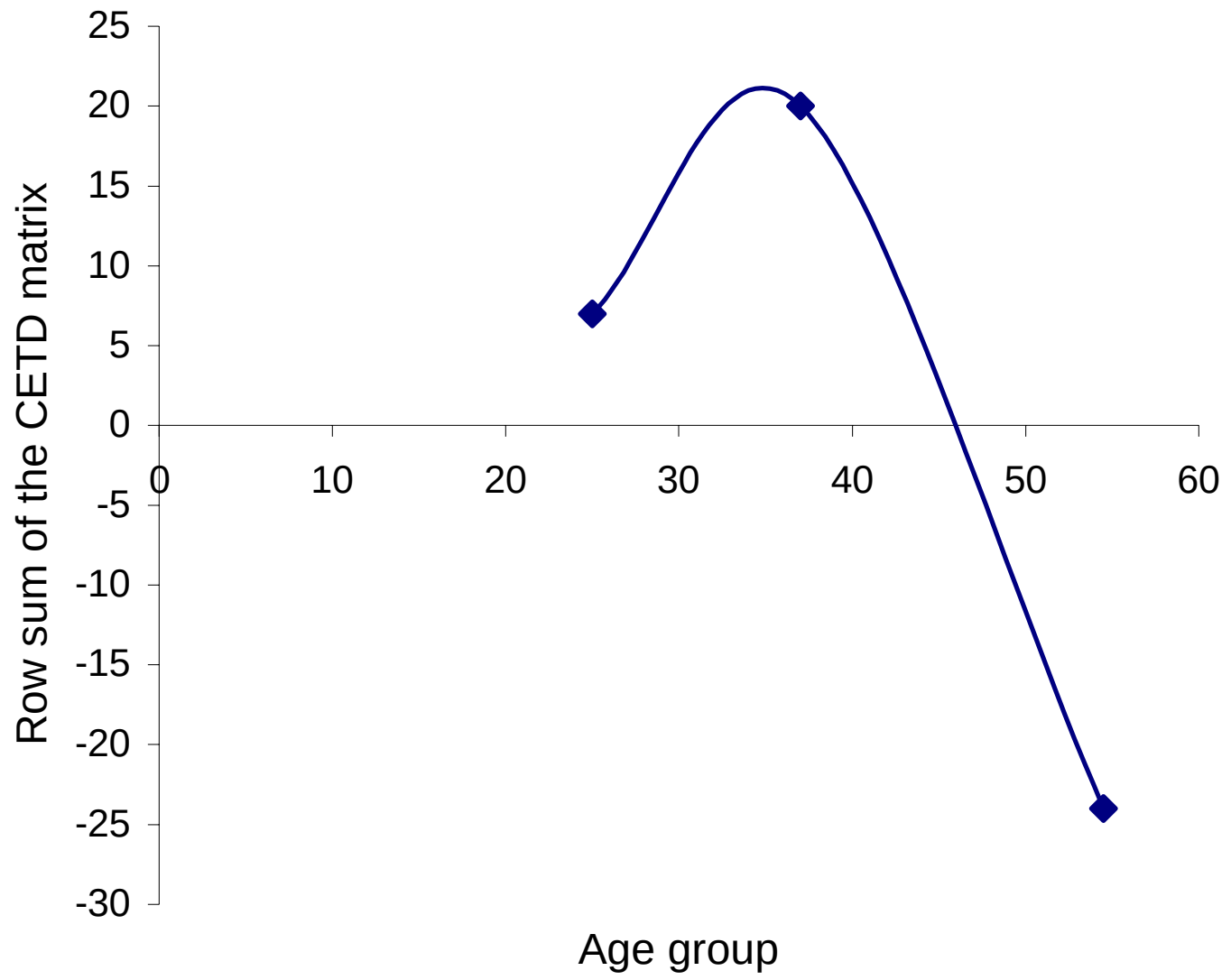

### 2.1.2.2 Estimation of maximum age group of agricultural labourers with pollution related digestive problems using 4 × 6 matrices

Now to make the study more sensitive, by increasing number of rows by 4 and see whether the decision arrived is more sensitive to the earlier one we have discussed. Thus we give the raw data of 4 × 6 matrix.

Initial Raw Data Matrix of digestive problem of order 4 × 6

| Years | $S_1$ | $S_2$ | $S_3$ | $S_4$ | $S_5$ | $S_6$ |
|---|---|---|---|---|---|---|
| 20-30 | 17 | 11 | 8 | 16 | 16 | 17 |
| 31-36 | 14 | 11 | 4 | 9 | 11 | 11 |
| 37-43 | 15 | 15 | 11 | 8 | 9 | 9 |
| 44-65 | 22 | 15 | 11 | 9 | 12 | 15 |

ATD Matrix of digestive problem of order 4 × 6

| Years | $S_1$ | $S_2$ | $S_3$ | $S_4$ | $S_5$ | $S_6$ |
|---|---|---|---|---|---|---|
| 20-30 | 1.55 | 1 | 0.73 | 1.45 | 1.45 | 1.55 |
| 31-36 | 2.33 | 1.83 | 0.67 | 1.5 | 1.83 | 1.83 |
| 37-43 | 2.14 | 2.14 | 1.57 | 1.14 | 1.29 | 1.29 |
| 44-65 | 1 | 0.68 | 0.5 | 0.4 | 0.55 | 0.68 |

The Average and the Standard Deviation

| | | | | | | |
|---|---|---|---|---|---|---|
| Average | 1.76 | 1.41 | 0.87 | 1.12 | 1.28 | 1.34 |
| SD | 0.6 | 0.69 | 0.48 | 0.51 | 0.54 | 0.49 |

RTD matrix for $\alpha = 0.15$ Row sum matrix

$$\begin{bmatrix} -1 & -1 & -1 & 1 & 1 & 1 \\ 1 & 1 & -1 & 1 & 1 & 1 \\ 1 & 1 & 1 & 0 & 0 & 0 \\ -1 & -1 & -1 & -1 & -1 & -1 \end{bmatrix} \qquad \begin{bmatrix} 0 \\ 4 \\ 3 \\ -6 \end{bmatrix}$$

RTD matrix for $\alpha = 0.35$ Row sum matrix

$$\begin{bmatrix} 0 & -1 & 0 & 1 & 0 & 1 \\ 1 & 1 & -1 & 1 & 1 & 1 \\ 1 & 1 & 1 & 0 & 0 & 0 \\ -1 & -1 & -1 & -1 & -1 & -1 \end{bmatrix} \qquad \begin{bmatrix} 1 \\ 4 \\ 3 \\ -6 \end{bmatrix}$$

RTD matrix for α = 0.45 Row sum matrix

$$\begin{bmatrix} 0 & -1 & 0 & 1 & 0 & 0 \\ 1 & 1 & 0 & 1 & 1 & 1 \\ 1 & 1 & 1 & 0 & 0 & 0 \\ -1 & -1 & -1 & -1 & -1 & -1 \end{bmatrix} \qquad \begin{bmatrix} 0 \\ 5 \\ 3 \\ -6 \end{bmatrix}$$

RTD matrix for α = 0.75 Row sum matrix

$$\begin{bmatrix} 0 & 0 & 0 & 0 & 0 & 0 \\ 1 & 0 & 0 & 0 & 1 & 1 \\ 0 & 1 & 1 & 0 & 0 & 0 \\ -1 & -1 & -1 & -1 & -1 & -1 \end{bmatrix} \qquad \begin{bmatrix} 0 \\ 3 \\ 2 \\ -6 \end{bmatrix}$$

Graph 25: Graph depicting maximum age group of labourers with pollution related digestive problem for α = 0.15

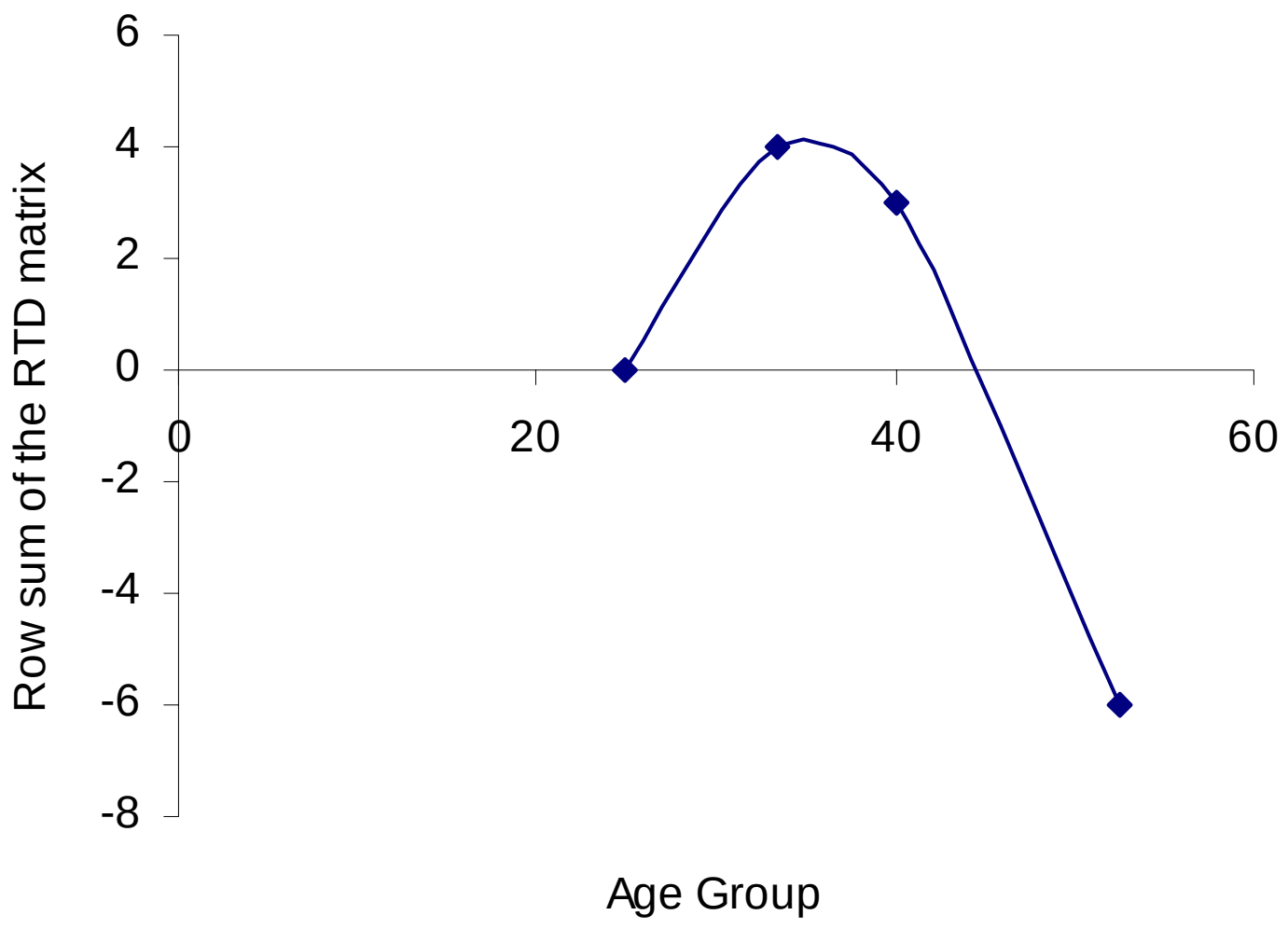

Graph 26: Graph depicting maximum age group of labourers with pollution related digestive problem for $\alpha = 0.35$

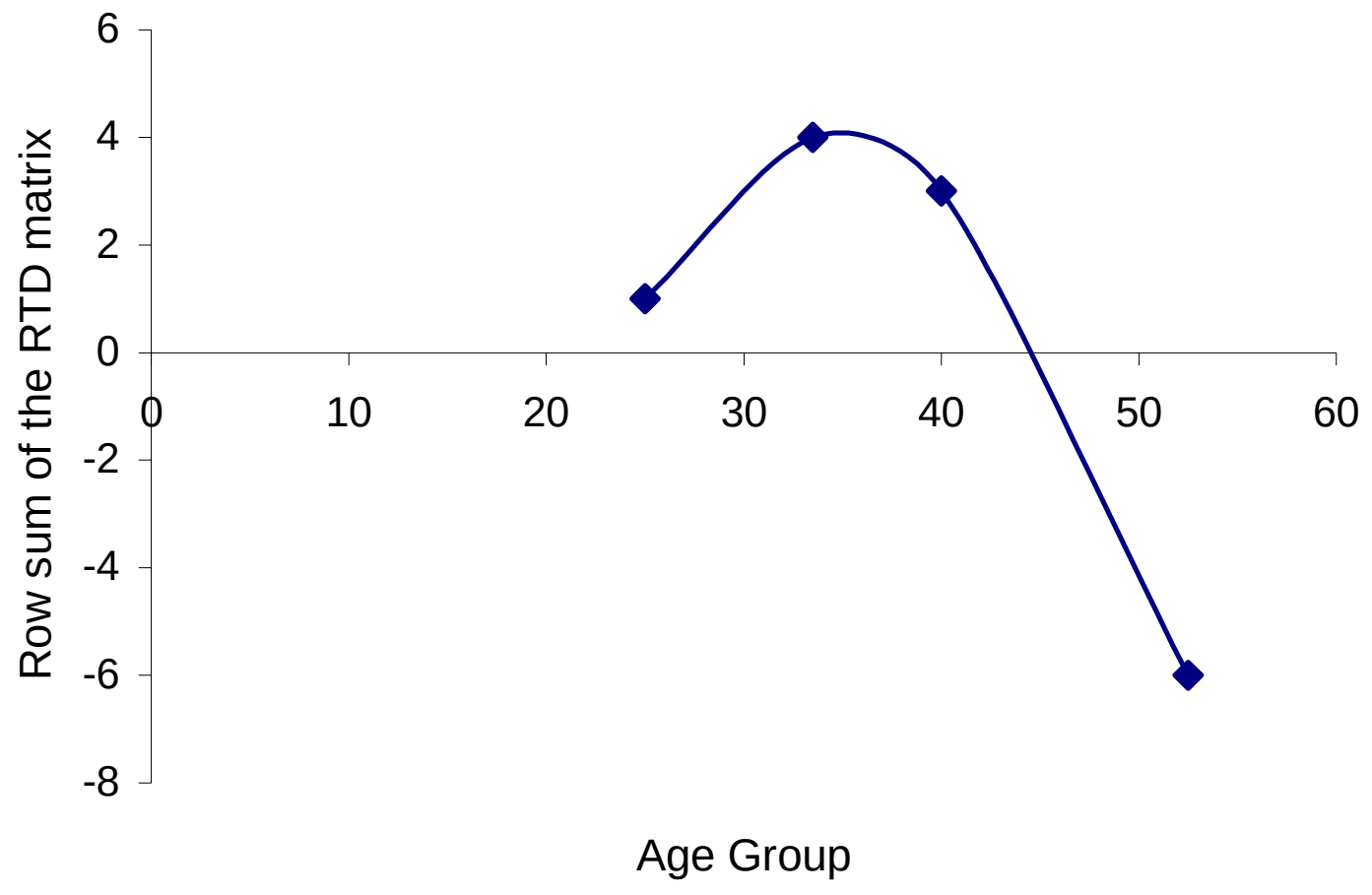

Graph 27: Graph depicting maximum age group of labourers with pollution related digestive problem for $\alpha = 0.45$

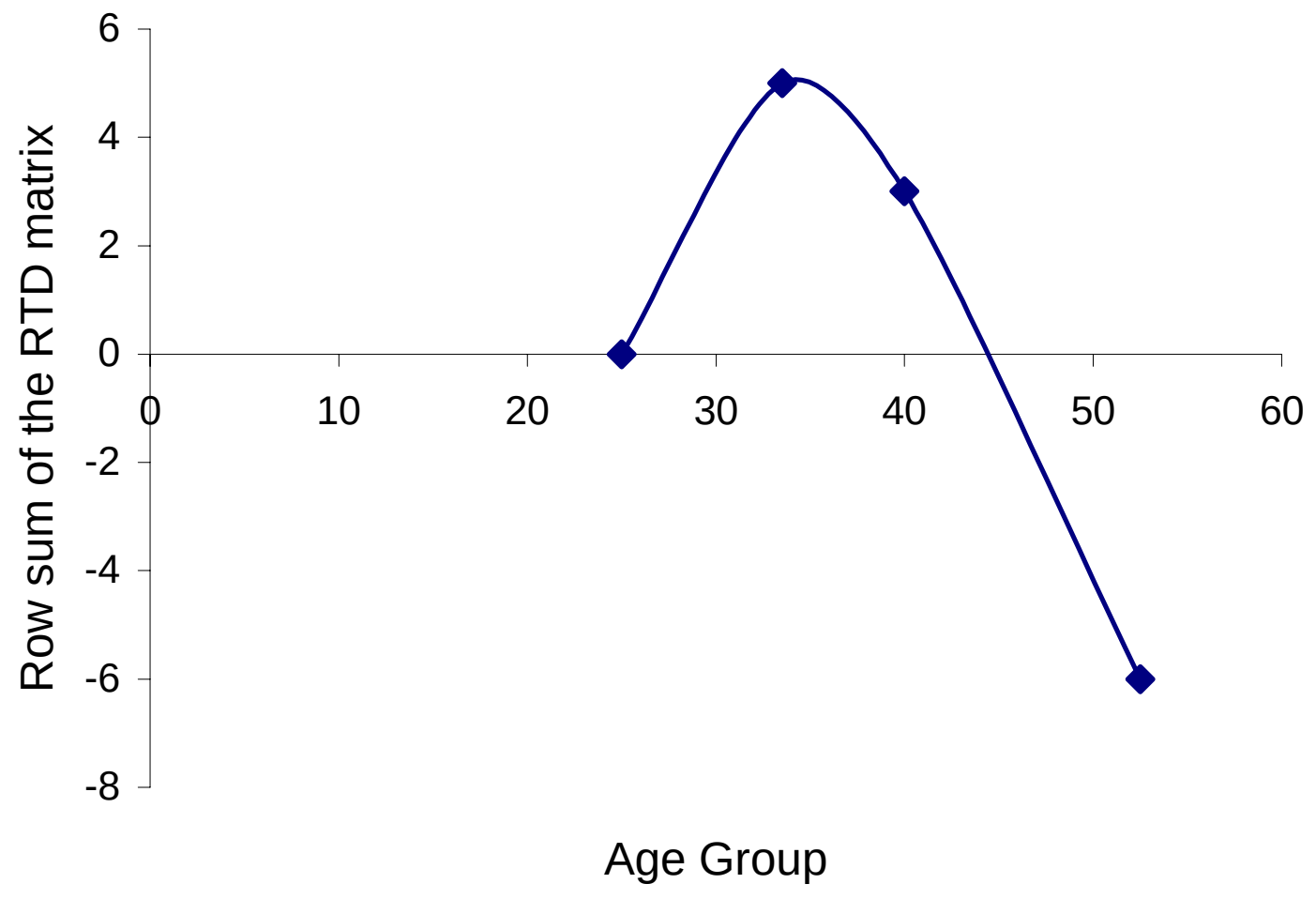

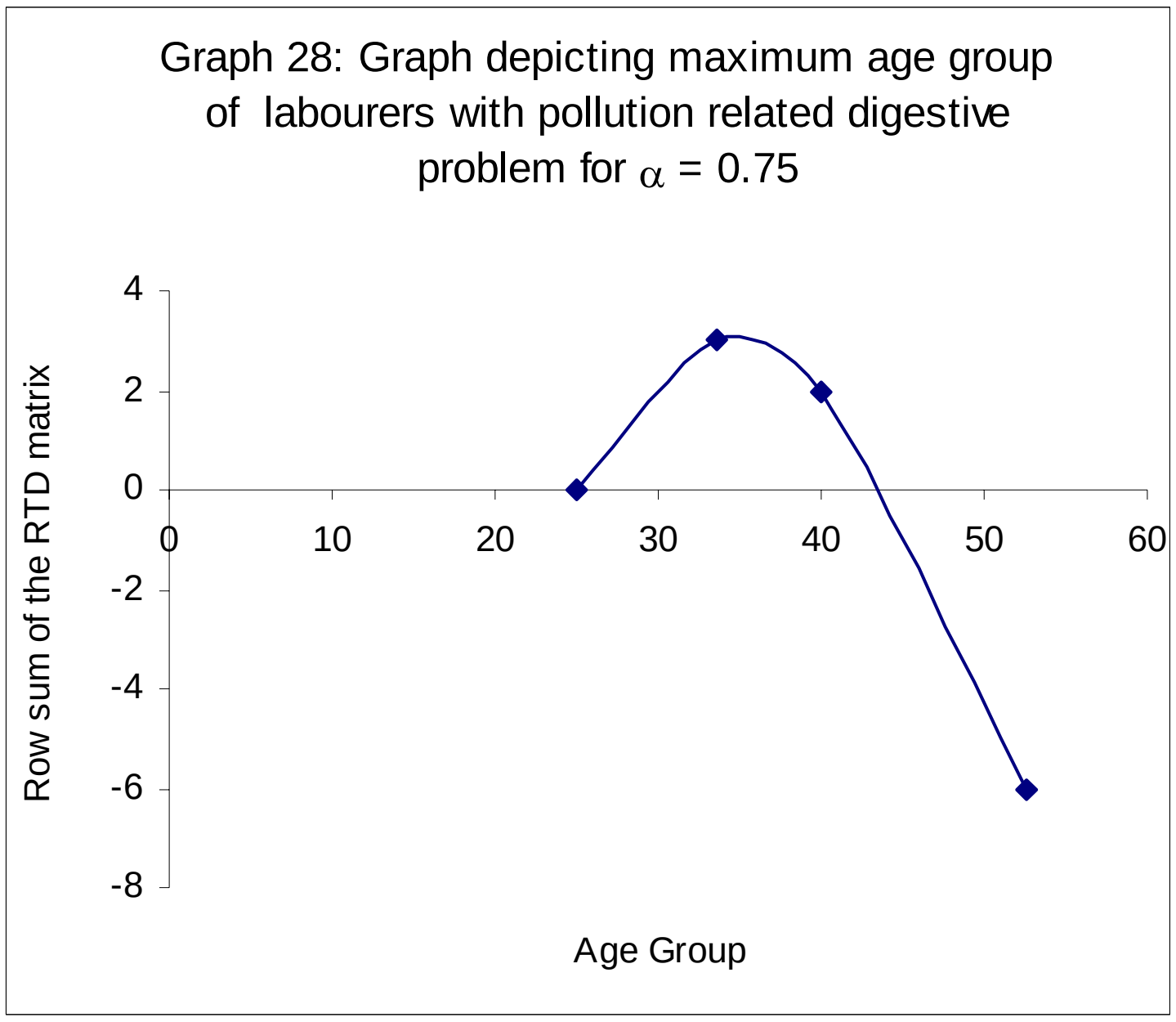
Graph 28: Graph depicting maximum age group of labourers with pollution related digestive problem for α = 0.75
Row sum of the RTD matrix
Age Group
4
2
0
-2
-4
-6
-8
0
10
20
30
40
50
60

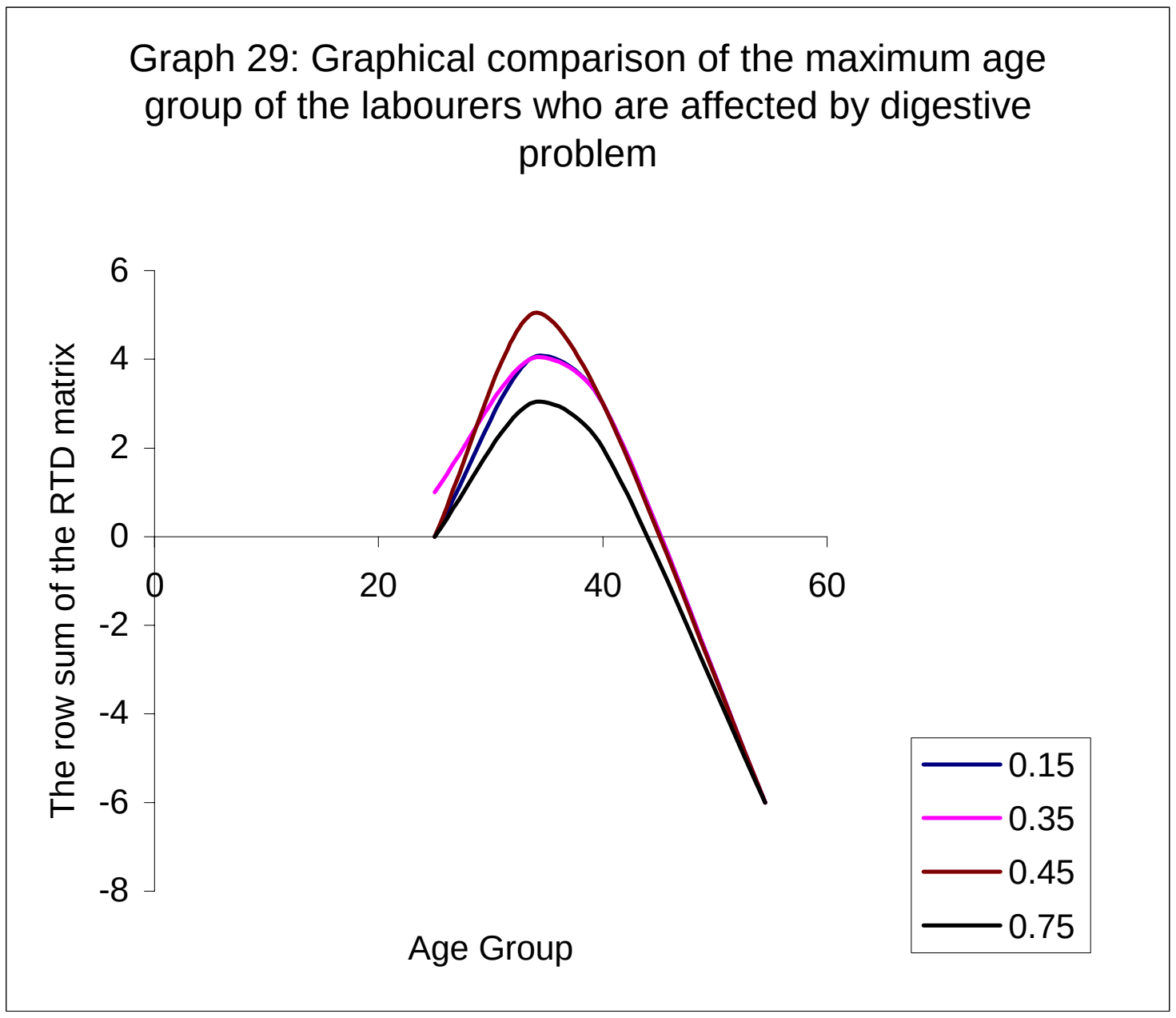
Graph 29: Graphical comparison of the maximum age group of the labourers who are affected by digestive problem
The row sum of the RTD matrix
Age Group
6
4
2
0
-2
-4
-6
-8
0
20
40
60
0.15
0.35
0.45
0.75

The CETD matrix and row sum matrix are given below:

CETD matrix

$$\begin{bmatrix} -1 & -3 & -1 & 3 & 1 & 2 \\ 4 & 3 & -2 & 3 & 4 & 4 \\ 3 & 4 & 4 & 0 & 0 & 0 \\ -4 & -4 & -4 & -4 & -4 & -4 \end{bmatrix}$$

Row Sum matrix

$$\begin{bmatrix} 1 \\ 16 \\ 11 \\ -24 \end{bmatrix}$$

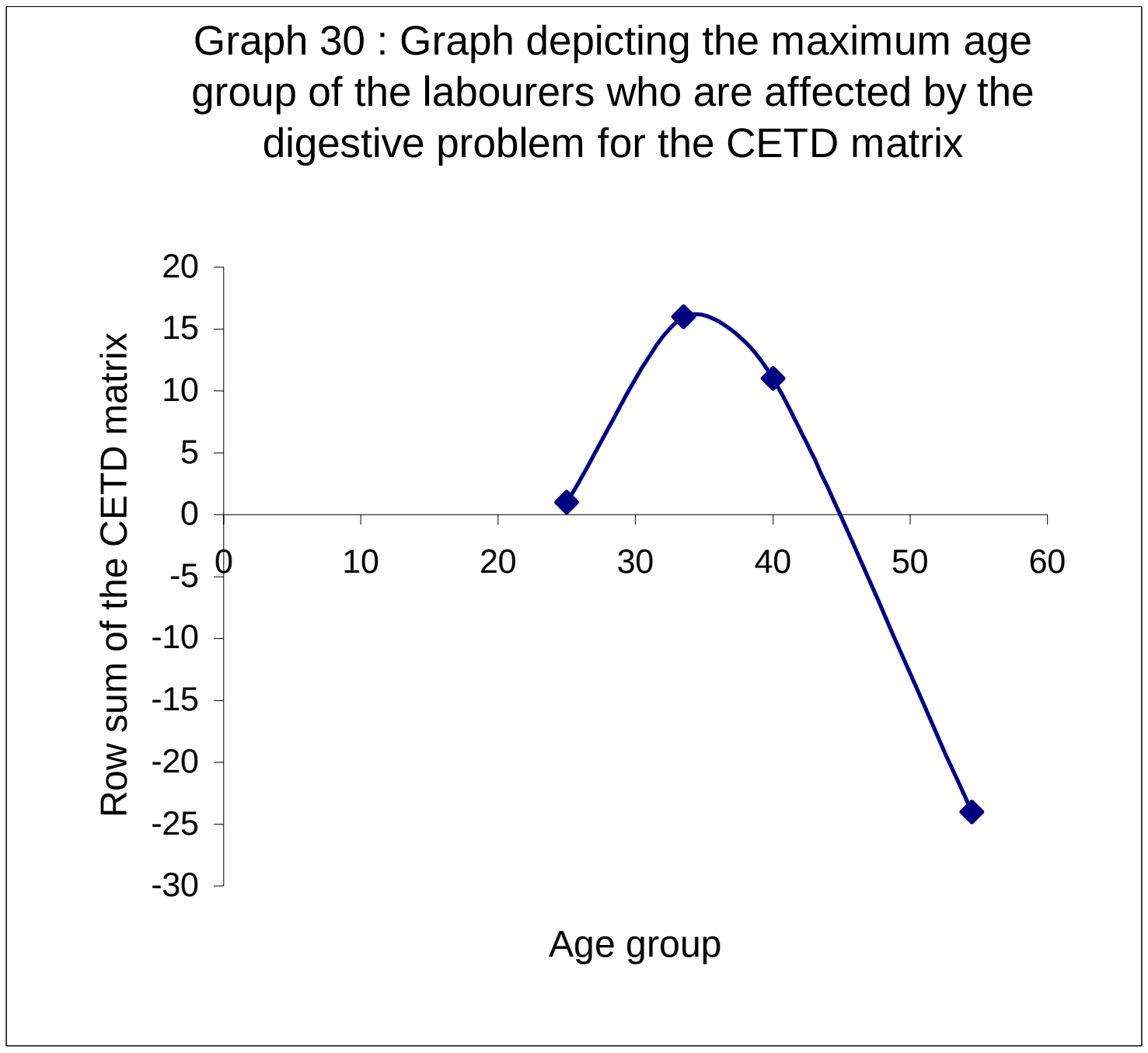

### 2.1.2.3 Estimation of maximum age group of agriculture labourer with pollution related digestive problem, using 5 × 6 matrix

Now to make the study still more sensitive by increasing the number of rows by 5 and see whether the decision arrived is more sensitive to the earlier one we have discussed. Thus we give the raw data of 5 × 6 matrix.

Initial Raw Data Matrix of digestive problem of order 5 × 6

| Years | $S_1$ | $S_2$ | $S_3$ | $S_4$ | $S_5$ | $S_6$ |
|---|---|---|---|---|---|---|
| 20-24 | 3 | 3 | 2 | 5 | 2 | 3 |
| 25-30 | 14 | 8 | 6 | 11 | 14 | 14 |
| 31-36 | 14 | 11 | 4 | 9 | 11 | 11 |
| 37-43 | 15 | 15 | 11 | 8 | 9 | 9 |
| 44-65 | 22 | 15 | 11 | 9 | 12 | 15 |

The ATD Matrix of digestive problem of order 5 × 6

| Years | $S_1$ | $S_2$ | $S_3$ | $S_4$ | $S_5$ | $S_6$ |
|---|---|---|---|---|---|---|
| 20-24 | 0.6 | 0.6 | 0.4 | 1 | 0.4 | 0.6 |
| 25-30 | 2.33 | 1.33 | 1 | 1.83 | 2.33 | 2.33 |
| 31-36 | 2.33 | 1.83 | 0.67 | 1.5 | 1.83 | 1.83 |
| 37-43 | 2.14 | 2.14 | 1.57 | 1.14 | 1.29 | 1.29 |
| 44-65 | 1 | 0.68 | 0.5 | 0.4 | 0.55 | 0.68 |

The Average and Standard Deviation of the above ATD matrix

| | | | | | | |
|---|---|---|---|---|---|---|
| Average | 1.68 | 1.32 | 0.83 | 1.17 | 1.28 | 1.35 |
| Standard Deviation | 0.82 | 0.68 | 0.47 | 0.54 | 0.82 | 0.74 |

RTD matrix for α = 0.1 — Row sum matrix

$$\begin{bmatrix} -1 & -1 & -1 & -1 & -1 & -1 \\ 1 & 0 & 1 & 1 & 1 & 1 \\ 1 & 1 & 0 & 1 & 1 & 1 \\ 1 & 1 & 1 & 0 & 0 & 0 \\ -1 & -1 & -1 & -1 & -1 & -1 \end{bmatrix} \qquad \begin{bmatrix} -6 \\ 5 \\ 5 \\ 3 \\ -6 \end{bmatrix}$$

RTD matrix for α = 0.15 — Row sum matrix

$$\begin{bmatrix} -1 & -1 & -1 & -1 & -1 & -1 \\ 1 & 0 & 1 & 1 & 1 & 1 \\ 1 & 1 & -1 & 1 & 1 & 1 \\ 1 & 1 & 1 & 0 & 0 & 0 \\ -1 & -1 & -1 & -1 & -1 & -1 \end{bmatrix} \qquad \begin{bmatrix} -6 \\ 5 \\ 4 \\ 3 \\ -6 \end{bmatrix}$$

RTD matrix for α = 0.2 — Row sum matrix

$$\begin{bmatrix} -1 & -1 & -1 & -1 & -1 & -1 \\ 1 & 0 & 1 & 1 & 1 & 1 \\ 1 & 1 & -1 & 1 & 1 & 1 \\ 1 & 1 & 1 & 0 & 0 & 0 \\ -1 & -1 & -1 & -1 & -1 & -1 \end{bmatrix} \qquad \begin{bmatrix} -6 \\ 5 \\ 4 \\ 3 \\ -6 \end{bmatrix}$$

RTD matrix for α = 0.35 — Row sum matrix

$$\begin{bmatrix} -1 & -1 & -1 & 0 & -1 & -1 \\ 1 & 0 & 1 & 1 & 1 & 1 \\ 1 & 1 & 0 & 1 & 1 & 1 \\ 1 & 1 & 1 & 0 & 0 & 0 \\ -1 & -1 & -1 & -1 & -1 & -1 \end{bmatrix} \qquad \begin{bmatrix} -5 \\ 5 \\ 5 \\ 3 \\ -6 \end{bmatrix}$$

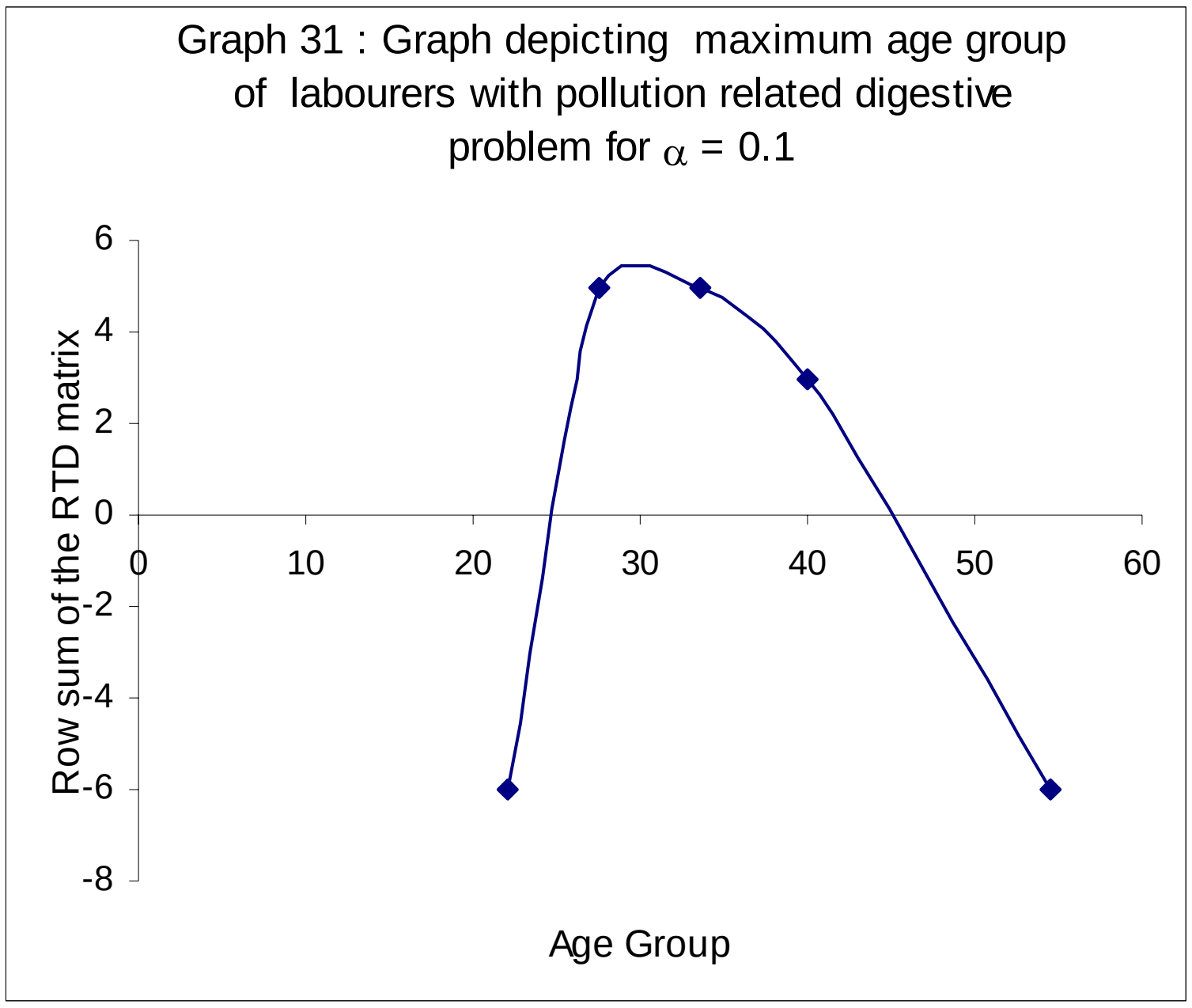
Graph 31 : Graph depicting maximum age group of labourers with pollution related digestive problem for α = 0.1
Row sum of the RTD matrix
6
4
2
0
-2
-4
-6
-8
0
10
20
30
40
50
60
Age Group

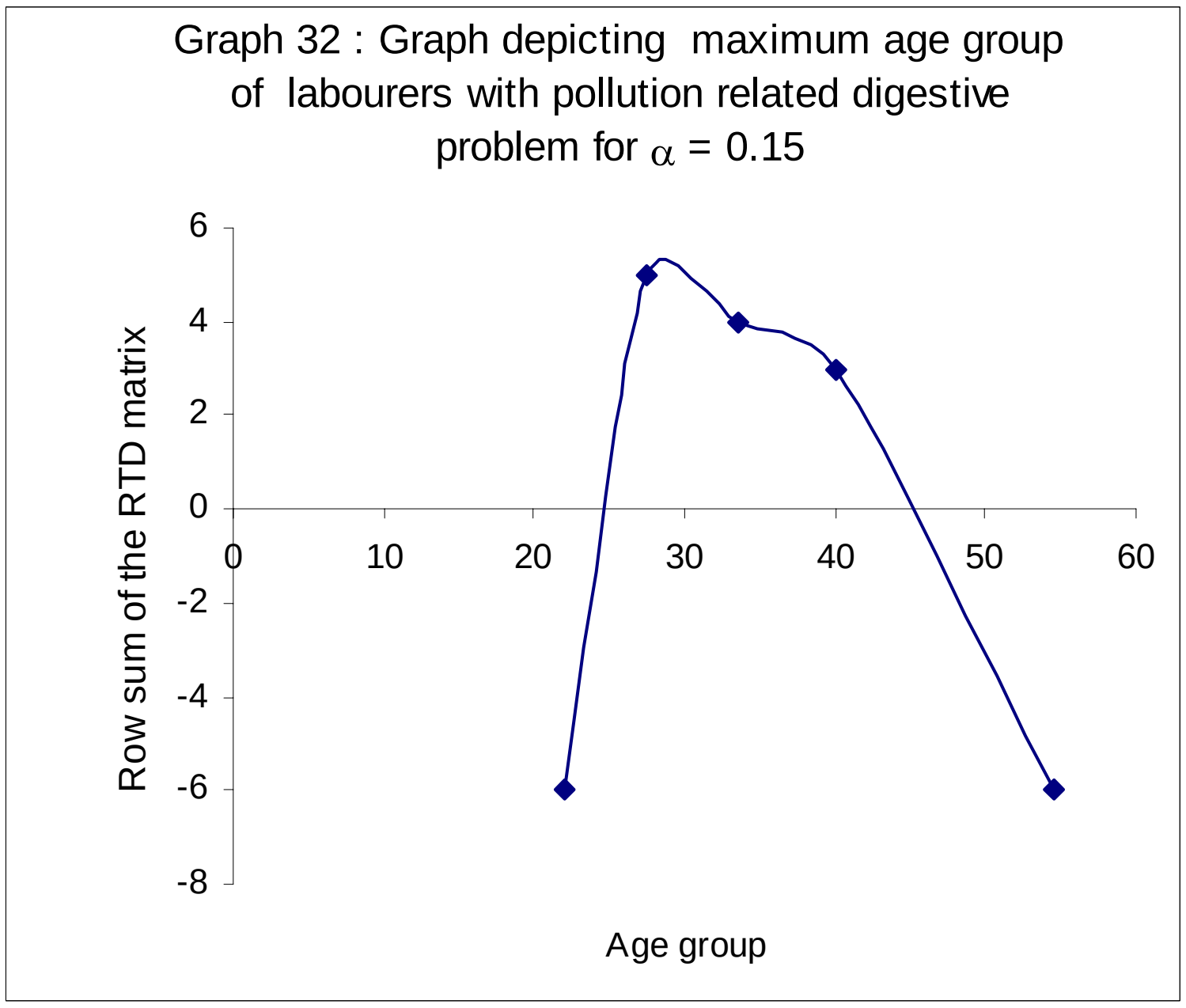
Graph 32 : Graph depicting maximum age group of labourers with pollution related digestive problem for α = 0.15
Row sum of the RTD matrix
6
4
2
0
-2
-4
-6
-8
0
10
20
30
40
50
60
Age group

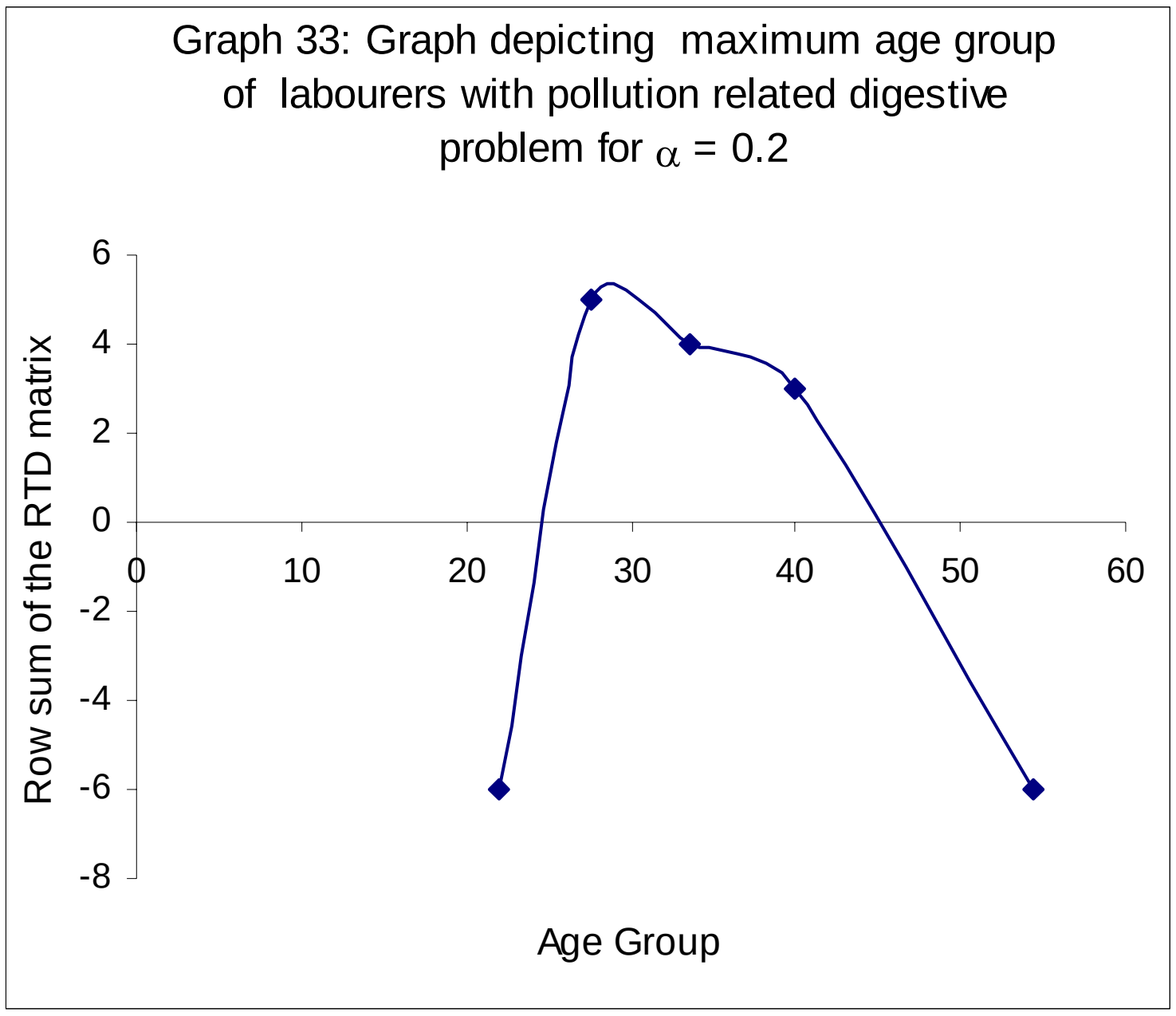
Graph 33: Graph depicting maximum age group of labourers with pollution related digestive problem for α = 0.2
Row sum of the RTD matrix
6
4
2
0
-2
-4
-6
-8
0
10
20
30
40
50
60
Age Group

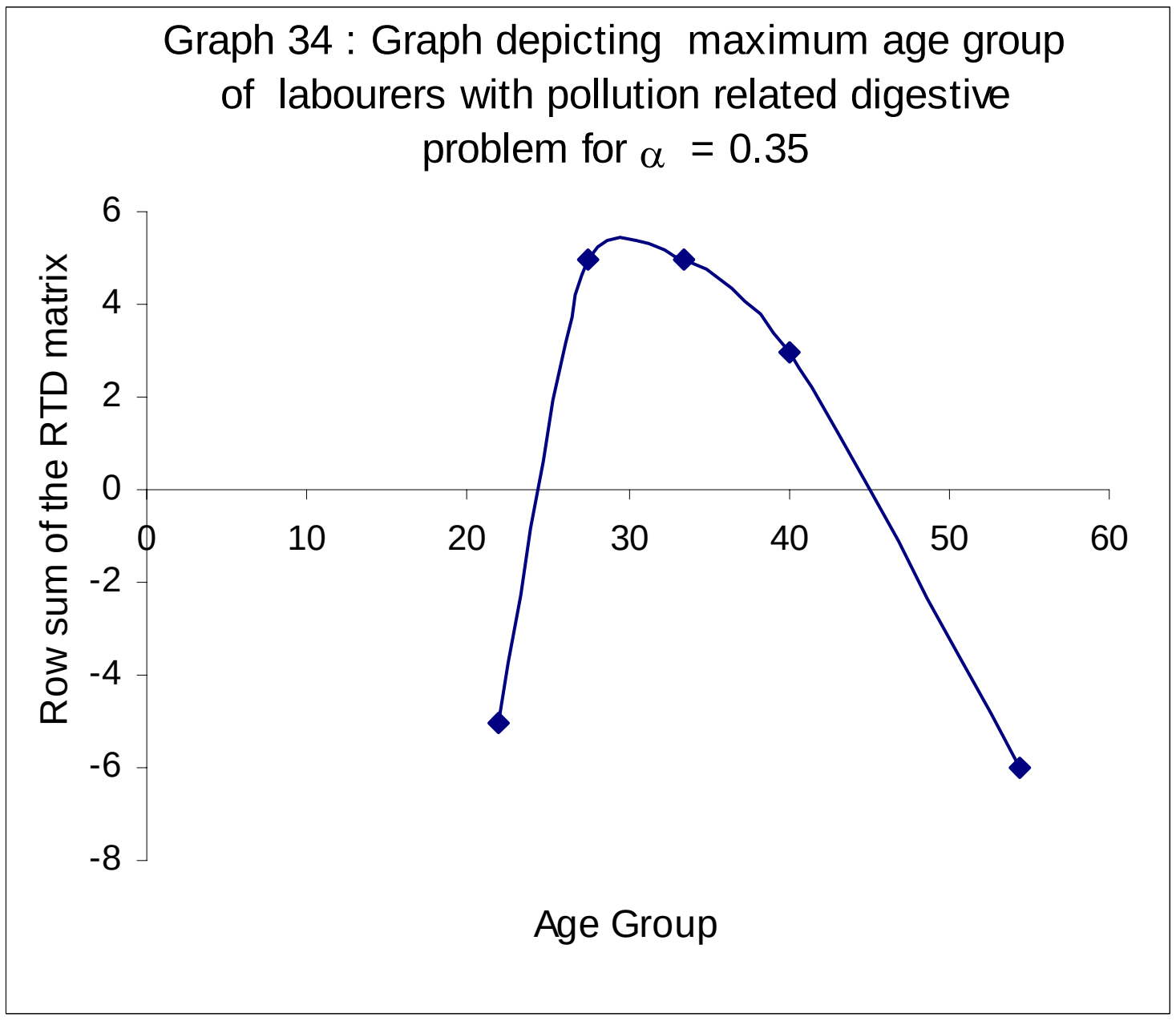
Graph 34 : Graph depicting maximum age group of labourers with pollution related digestive problem for α = 0.35
Row sum of the RTD matrix
6
4
2
0
-2
-4
-6
-8
0
10
20
30
40
50
60
Age Group

The comparative graph of the maximum age group of labourers who are affected by digestive problem is as follows:

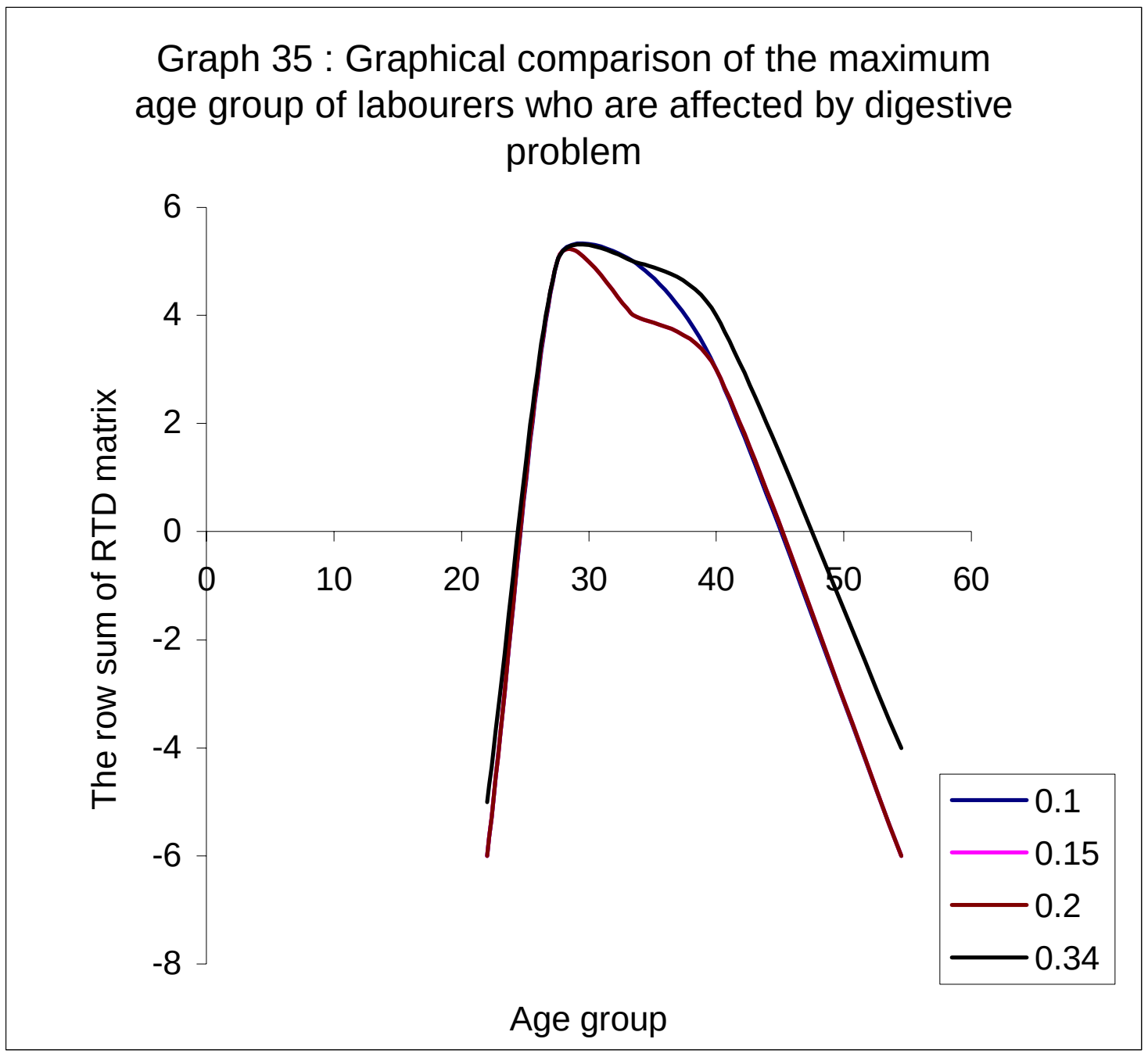

The CETD matrix and the respective row sum matrix is given below:

CETD matrix Row sum matrix

$$\begin{bmatrix} -4 & -4 & -4 & -3 & -4 & -4 \\ 4 & 0 & 4 & 4 & 4 & 4 \\ 4 & 4 & -2 & 4 & 4 & 4 \\ 4 & 4 & 4 & 0 & 0 & 0 \\ -4 & -4 & -4 & -4 & -4 & -4 \end{bmatrix} \qquad \begin{bmatrix} -23 \\ 20 \\ 18 \\ 12 \\ -24 \end{bmatrix}$$

The graph for the CETD matrix is given in the next page:

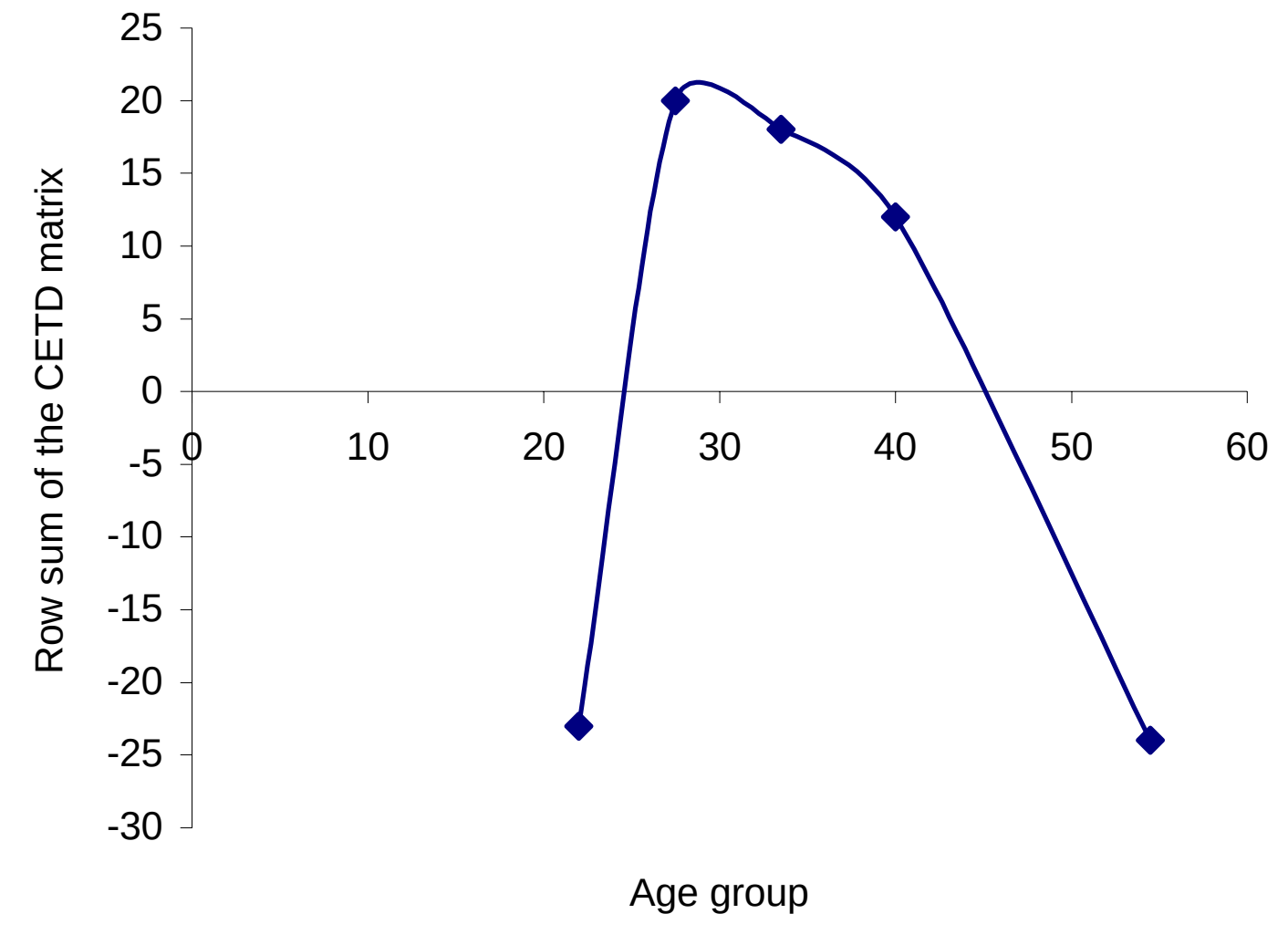

Graph 36: The graph depicting the maximum age group of the labourers who are affected by Digestive problem for CETD matrix

### 2.1.2.4 Conclusion

The general reason physicians give for digestive problem are food-fats, irregularity in food consumption and gastro-intestinal malfunction. The first two reasons are prevalent among urban areas. But in a rural community to have this problem between the age group of 20-30 with maximum level at this age of 27 as shown in the graph could be due to gastrointestinal malfunction. This age group generally shows a sound digestive capacity as their enzymatic secretion and intestinal motility are at its best. The findings have revealed that out of 110 interviewed, 68 of them showed loss of appetite; 58 showed mouth and stomach ulcer; 52 of them vomiting (nausea) and indigestion i.e. more than half the population interviewed complained of digestive problems. This contradicting note among the agricultural labourers makes one to conclude the proximity of pesticides as a strong reason. Thus at the stage it is pertinent to mention that only due to chemical pollutions these agriculture coolie are facing these health problems in these recent day. This conclusion has been completely supported by our analysis of the 110 agriculture laborers using fuzzy matrix.

### 2.1.3 Estimation of the Maximum Age Group of the Agriculture Labourer having Nervous Problem due to Chemical Pollution, using Matrices

Nervous problem is taken under eight symptom diseases viz.

$S_1$ - Headache,
$S_2$ - Fainting,
$S_3$ - Nausea,
$S_4$ - Dizziness,
$S_5$ - Getting angry,
$S_6$ - Getting irritated,
$S_7$ - Fits,
$S_8$ - Loss of memory.

By taking symptoms along the rows and age along the columns.

#### 2.1.3.1 Estimation of maximum age group of agriculture labourer with pollution related nervous problems using 3 × 8 matrix

Initial Raw Data Matrix of nervous problem of order 3 × 8

| Years | $S_1$ | $S_2$ | $S_3$ | $S_4$ | $S_5$ | $S_6$ | $S_7$ | $S_8$ |
|---|---|---|---|---|---|---|---|---|
| 20-30 | 28 | 27 | 15 | 27 | 25 | 20 | 13 | 16 |
| 31-43 | 38 | 42 | 27 | 39 | 34 | 34 | 33 | 23 |
| 44-65 | 23 | 23 | 17 | 22 | 20 | 24 | 19 | 15 |

The ATD Matrix of nervous problem of order 3 × 8

| Years | $S_1$ | $S_2$ | $S_3$ | $S_4$ | $S_5$ | $S_6$ | $S_7$ | $S_8$ |
|---|---|---|---|---|---|---|---|---|
| 20-30 | 2.55 | 2.45 | 1.36 | 2.45 | 2.27 | 1.82 | 1.18 | 1.45 |
| 30-43 | 2.92 | 3.23 | 2.08 | 3 | 2.62 | 2.62 | 2.54 | 1.77 |
| 44-65 | 1.05 | 1.05 | 0.77 | 1 | 0.91 | 1.09 | 0.86 | 0.86 |

The Average and Standard Deviation
of the above given ATD matrix

| | | | | | | | | |
|---|---|---|---|---|---|---|---|---|
| Average | 2.17 | 2.24 | 1.40 | 2.15 | 1.93 | 1.84 | 1.53 | 1.36 |
| Standard Deviation | 0.99 | 1.11 | 0.66 | 1.03 | 0.90 | 0.77 | 0.89 | 0.46 |

RTD matrix for $\alpha$ = 0.15 Row sum matrix

$$\begin{bmatrix} 1 & 1 & 0 & 1 & 1 & 0 & -1 & 1 \\ 1 & 1 & 1 & 1 & 1 & 1 & 1 & 1 \\ -1 & -1 & -1 & -1 & -1 & -1 & -1 & -1 \end{bmatrix} \qquad \begin{bmatrix} 4 \\ 8 \\ -8 \end{bmatrix}$$

RTD matrix for $\alpha$ = 0.35 Row sum matrix

$$\begin{bmatrix} 1 & 0 & 0 & 0 & 1 & 0 & -1 & 0 \\ 1 & 1 & 1 & 1 & 1 & 1 & 1 & 1 \\ -1 & -1 & -1 & -1 & -1 & -1 & -1 & -1 \end{bmatrix} \qquad \begin{bmatrix} 1 \\ 8 \\ -8 \end{bmatrix}$$

RTD matrix for $\alpha$ = 0.45 Row sum matrix

$$\begin{bmatrix} 0 & 0 & 0 & 0 & 0 & 0 & 0 & 0 \\ 1 & 1 & 1 & 1 & 1 & 1 & 1 & 1 \\ -1 & -1 & -1 & -1 & -1 & -1 & -1 & -1 \end{bmatrix} \qquad \begin{bmatrix} 0 \\ 8 \\ -8 \end{bmatrix}$$

RTD matrix for $\alpha$ = 0.75 Row sum matrix

$$\begin{bmatrix} 0 & 0 & 0 & 0 & 0 & 0 & 0 & 0 \\ 1 & 1 & 1 & 1 & 1 & 1 & 1 & 1 \\ -1 & -1 & -1 & -1 & -1 & -1 & -1 & -1 \end{bmatrix} \qquad \begin{bmatrix} 0 \\ 8 \\ -8 \end{bmatrix}$$

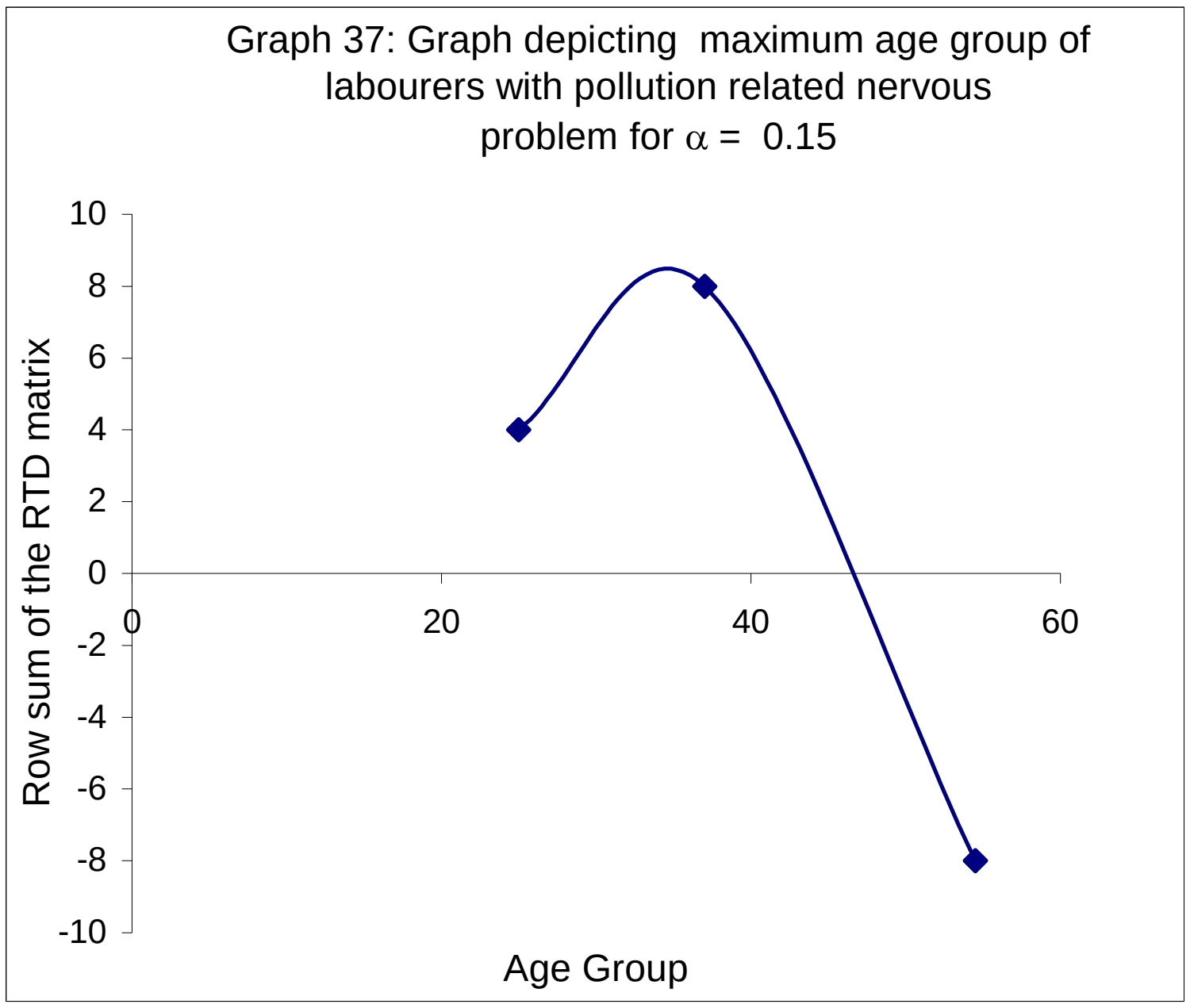
Graph 37: Graph depicting maximum age group of labourers with pollution related nervous problem for α = 0.15
10
8
6
4
2
0
-2
-4
-6
-8
-10
0
20
40
60
Row sum of the RTD matrix
Age Group

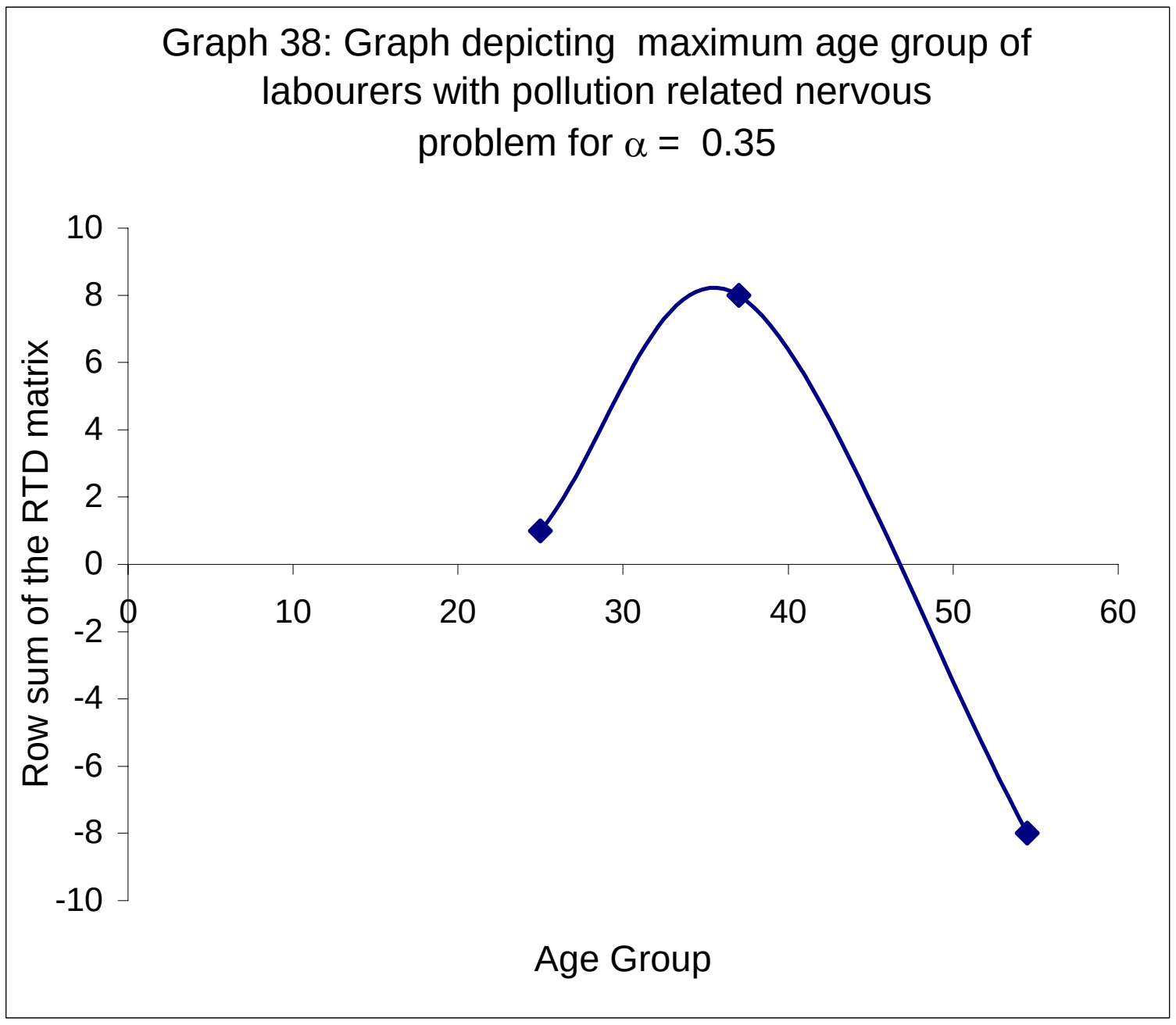
Graph 38: Graph depicting maximum age group of labourers with pollution related nervous problem for α = 0.35
10
8
6
4
2
0
-2
-4
-6
-8
-10
0
10
20
30
40
50
60
Row sum of the RTD matrix
Age Group

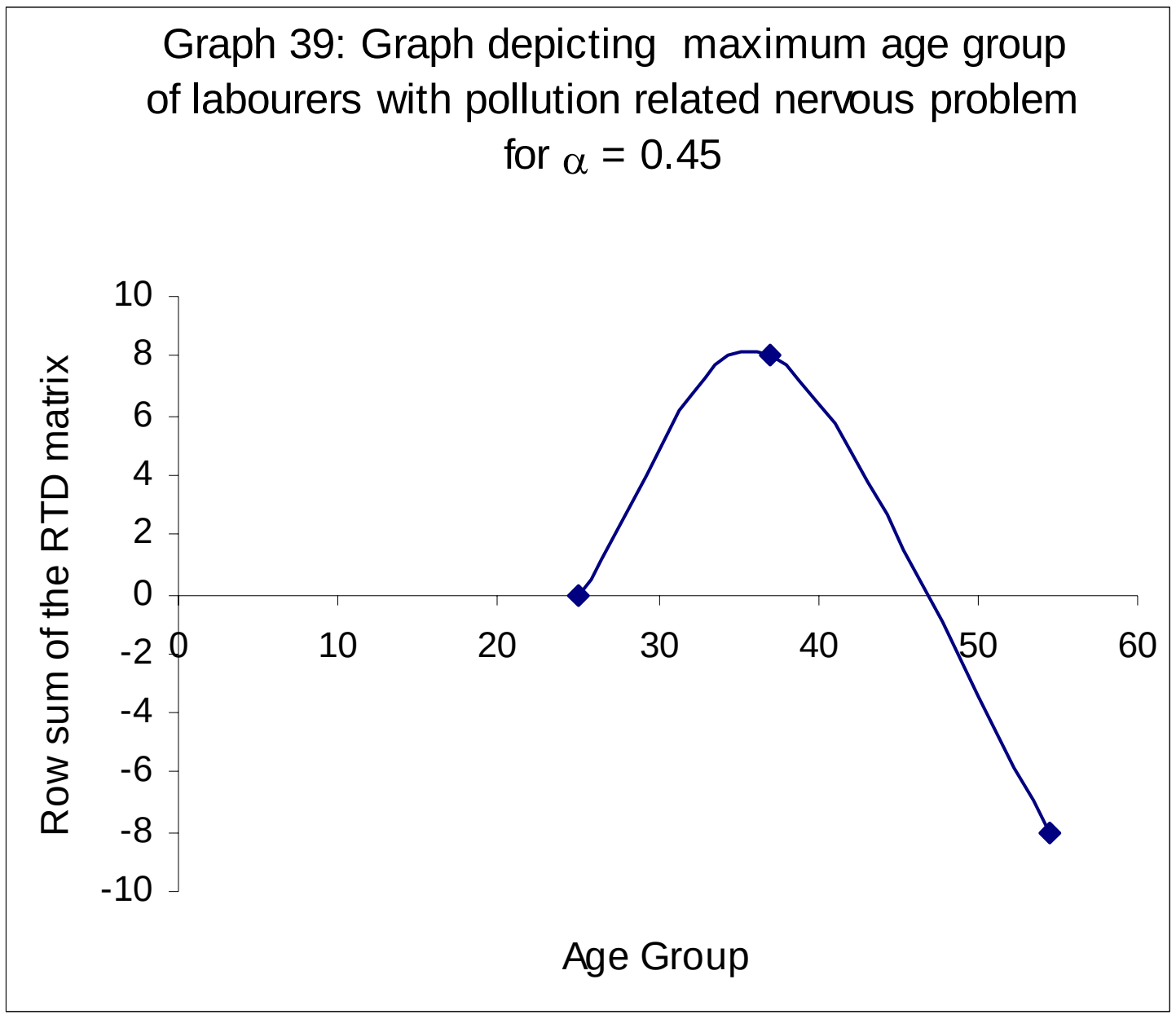
Graph 39: Graph depicting maximum age group of labourers with pollution related nervous problem for α = 0.45
Row sum of the RTD matrix
Age Group
10
8
6
4
2
0
-2
-4
-6
-8
-10
0
10
20
30
40
50
60

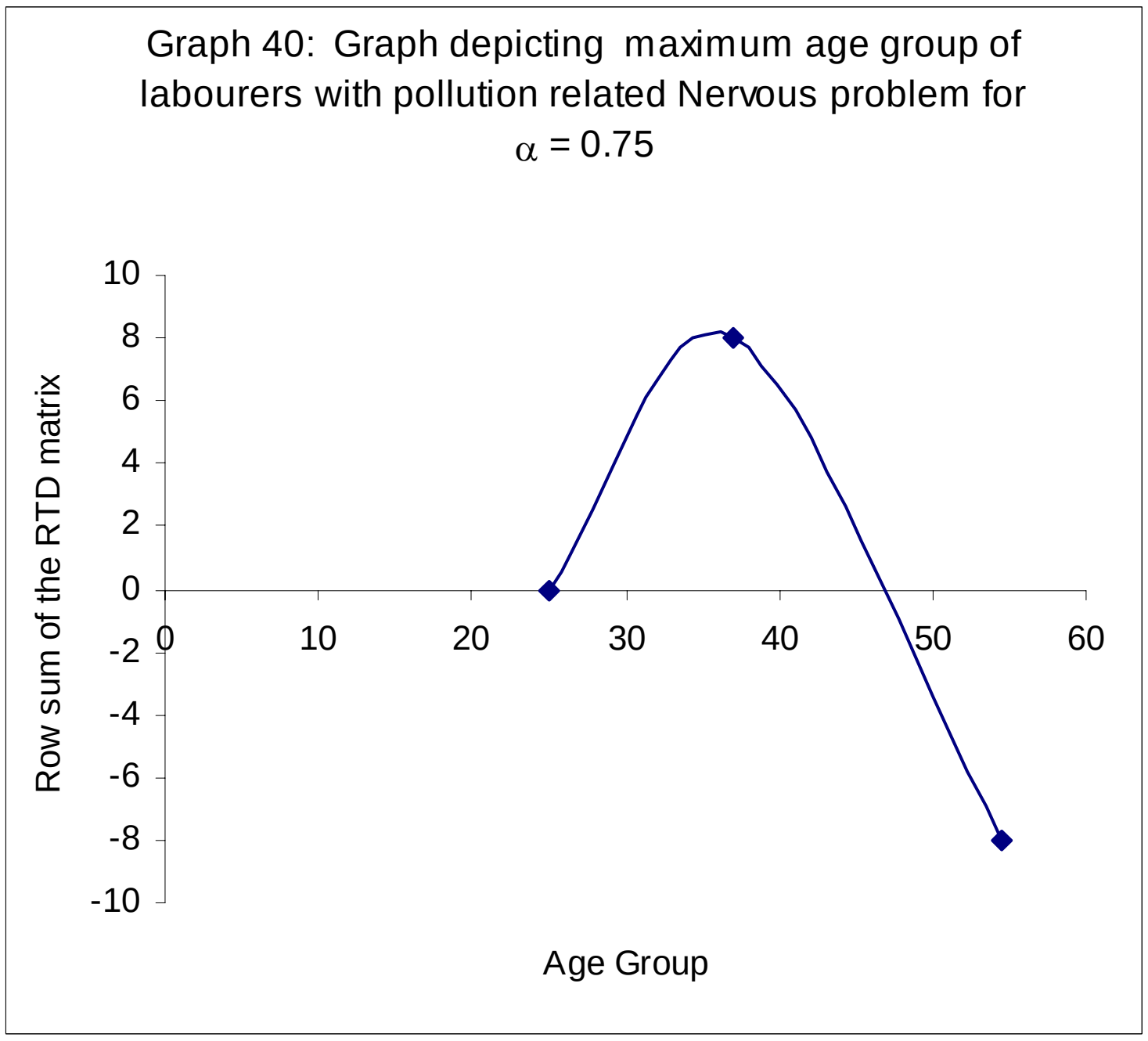
Graph 40: Graph depicting maximum age group of labourers with pollution related Nervous problem for α = 0.75
Row sum of the RTD matrix
Age Group
10
8
6
4
2
0
-2
-4
-6
-8
-10
0
10
20
30
40
50
60

Graph 41 gives the graphical comparison of the maximum age group of labourers affected by Nervous problem.

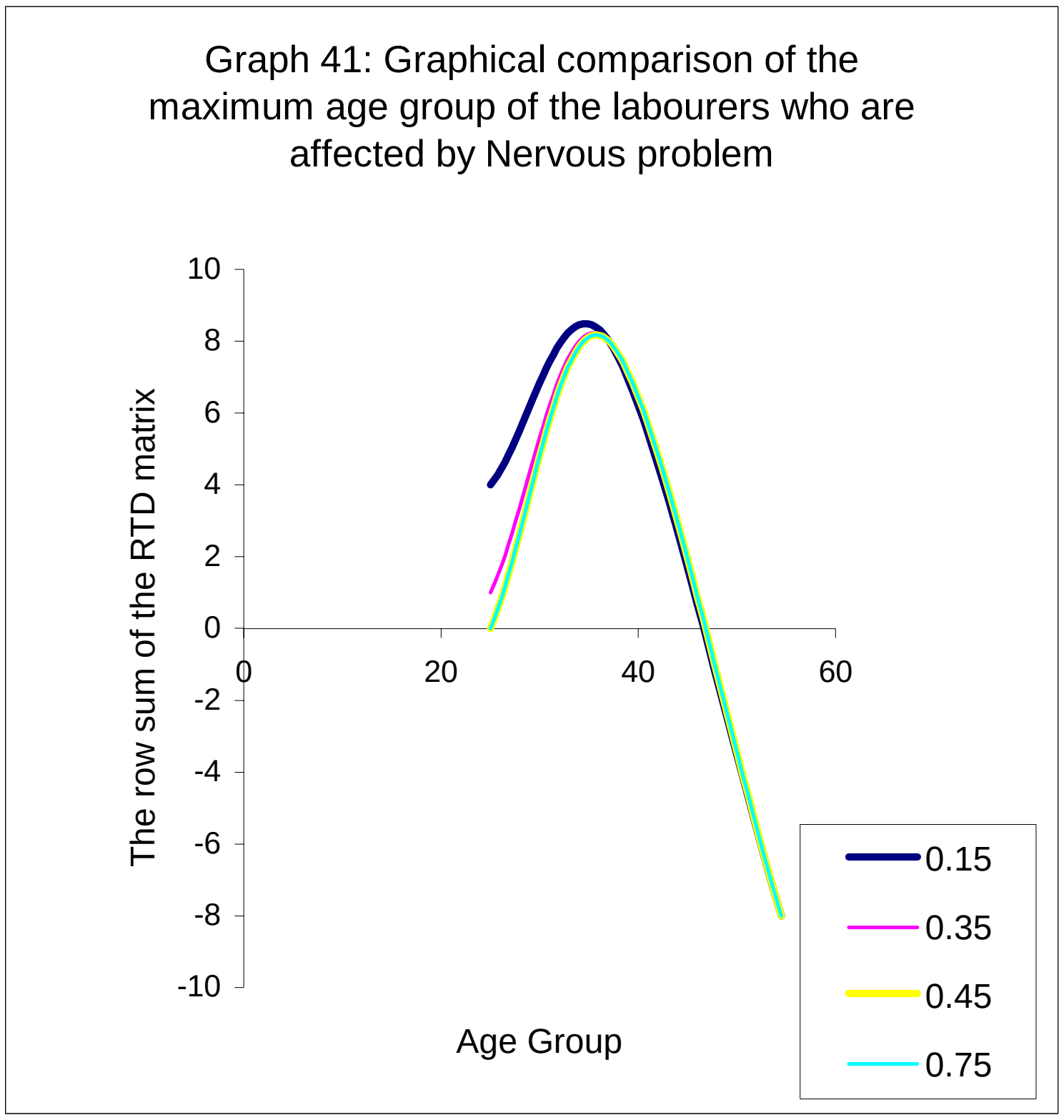

CETD matrix

$$\begin{bmatrix} 1 & 1 & -1 & 1 & 2 & 0 & -2 & 0 \\ 3 & 3 & 3 & 3 & 3 & 3 & 3 & 3 \\ -3 & -3 & -3 & -3 & -3 & -3 & -3 & -3 \end{bmatrix}$$

Row sum matrix

$$\begin{bmatrix} 2 \\ 24 \\ -24 \end{bmatrix}$$

The graph for the above given CETD matrix is given in the next page:

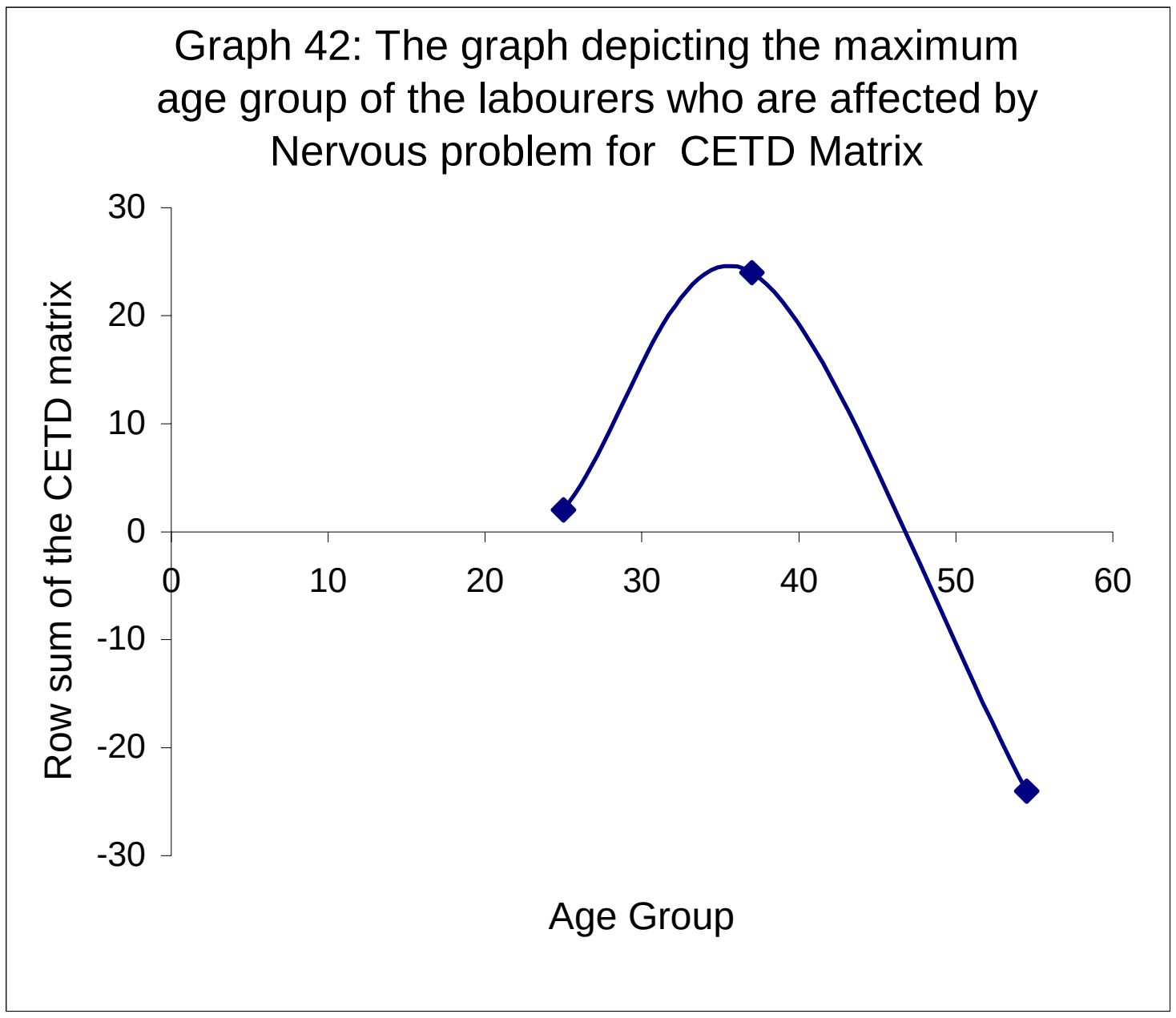

### 2.1.3.2 Estimation of maximum age group of agriculture labourer with pollution related Nervous problems using 4 × 8 matrix

Now to make the study more sensitive we increase the number of rows by 4 and see whether the decision arrived is more sensitive to the earlier one we have discussed. Thus we give the raw data of 4 × 8 matrix

Initial Raw Data Matrix of nervous problem of order 4 × 8

| Years | $S_1$ | $S_2$ | $S_3$ | $S_4$ | $S_5$ | $S_6$ | $S_7$ | $S_8$ |
|---|---|---|---|---|---|---|---|---|
| 20-30 | 28 | 27 | 15 | 27 | 25 | 20 | 13 | 16 |
| 31-36 | 22 | 22 | 12 | 21 | 18 | 19 | 10 | 15 |
| 37-43 | 16 | 20 | 15 | 18 | 16 | 15 | 13 | 8 |
| 44-65 | 23 | 23 | 17 | 22 | 20 | 24 | 19 | 15 |

The ATD Matrix of nervous problem of order 4 × 8

| Years | $S_1$ | $S_2$ | $S_3$ | $S_4$ | $S_5$ | $S_6$ | $S_7$ | $S_8$ |
|---|---|---|---|---|---|---|---|---|
| 20-30 | 2.55 | 2.45 | 1.36 | 2.45 | 2.27 | 1.82 | 1.18 | 1.45 |
| 31-36 | 3.67 | 3.67 | 2 | 3.5 | 3 | 3.17 | 1.67 | 2.5 |
| 37-43 | 2.29 | 2.86 | 2.14 | 2.57 | 2.29 | 2.14 | 1.86 | 1.14 |
| 44-65 | 1.05 | 1.05 | 0.77 | 1 | 0.91 | 1.09 | 0.86 | 0.86 |

The Average and Standard Deviation of the ATD matrix

| | | | | | | | | |
|---|---|---|---|---|---|---|---|---|
| Average | 2.39 | 2.51 | 1.57 | 2.38 | 2.12 | 2.06 | 1.39 | 1.49 |
| Standard Deviation | 1.08 | 1.1 | 0.63 | 1.03 | 0.87 | 0.86 | 0.46 | 0.72 |

RTD matrix for α = 0.15   Row sum matrix

$$\begin{bmatrix} 0 & 0 & -1 & 0 & 1 & -1 & -1 & 0 \\ 1 & 1 & 1 & 1 & 1 & 0 & 1 & -1 \\ -1 & 1 & 1 & 1 & 1 & 0 & 1 & -1 \\ -1 & -1 & -1 & -1 & -1 & -1 & -1 & -1 \end{bmatrix} \qquad \begin{bmatrix} -2 \\ 5 \\ 3 \\ -8 \end{bmatrix}$$

RTD matrix for α = 0.35   Row sum matrix

$$\begin{bmatrix} 0 & 0 & 0 & 0 & 0 & 0 & -1 & 0 \\ 1 & 1 & 1 & 1 & 1 & 1 & 1 & 1 \\ -1 & 0 & 1 & 0 & 0 & 0 & 1 & -1 \\ -1 & -1 & -1 & -1 & -1 & -1 & -1 & -1 \end{bmatrix} \qquad \begin{bmatrix} -1 \\ 8 \\ 0 \\ -8 \end{bmatrix}$$

RTD matrix for α = 0.45   Row sum matrix

$$\begin{bmatrix} 0 & 0 & 0 & 0 & 0 & 0 & 0 & 0 \\ 1 & 1 & 1 & 1 & 1 & 1 & 1 & 1 \\ -1 & 0 & 1 & 0 & 0 & 0 & 1 & -1 \\ -1 & -1 & -1 & -1 & -1 & -1 & -1 & -1 \end{bmatrix} \qquad \begin{bmatrix} 0 \\ 8 \\ 0 \\ -8 \end{bmatrix}$$

RTD matrix for α = 0.75

$$\begin{bmatrix} 0 & 0 & 0 & 0 & 0 & 0 & 0 & 0 \\ 1 & 1 & 0 & 1 & 1 & 1 & 0 & 1 \\ -1 & 0 & 1 & 0 & 0 & 0 & 1 & 0 \\ -1 & -1 & -1 & -1 & -1 & -1 & -1 & -1 \end{bmatrix}$$

Row sum matrix

$$\begin{bmatrix} 0 \\ 6 \\ 1 \\ -8 \end{bmatrix}$$

Graph 43: Graph depicting maximum age group of labourers with pollution related Nervous Problem for α = 0.15

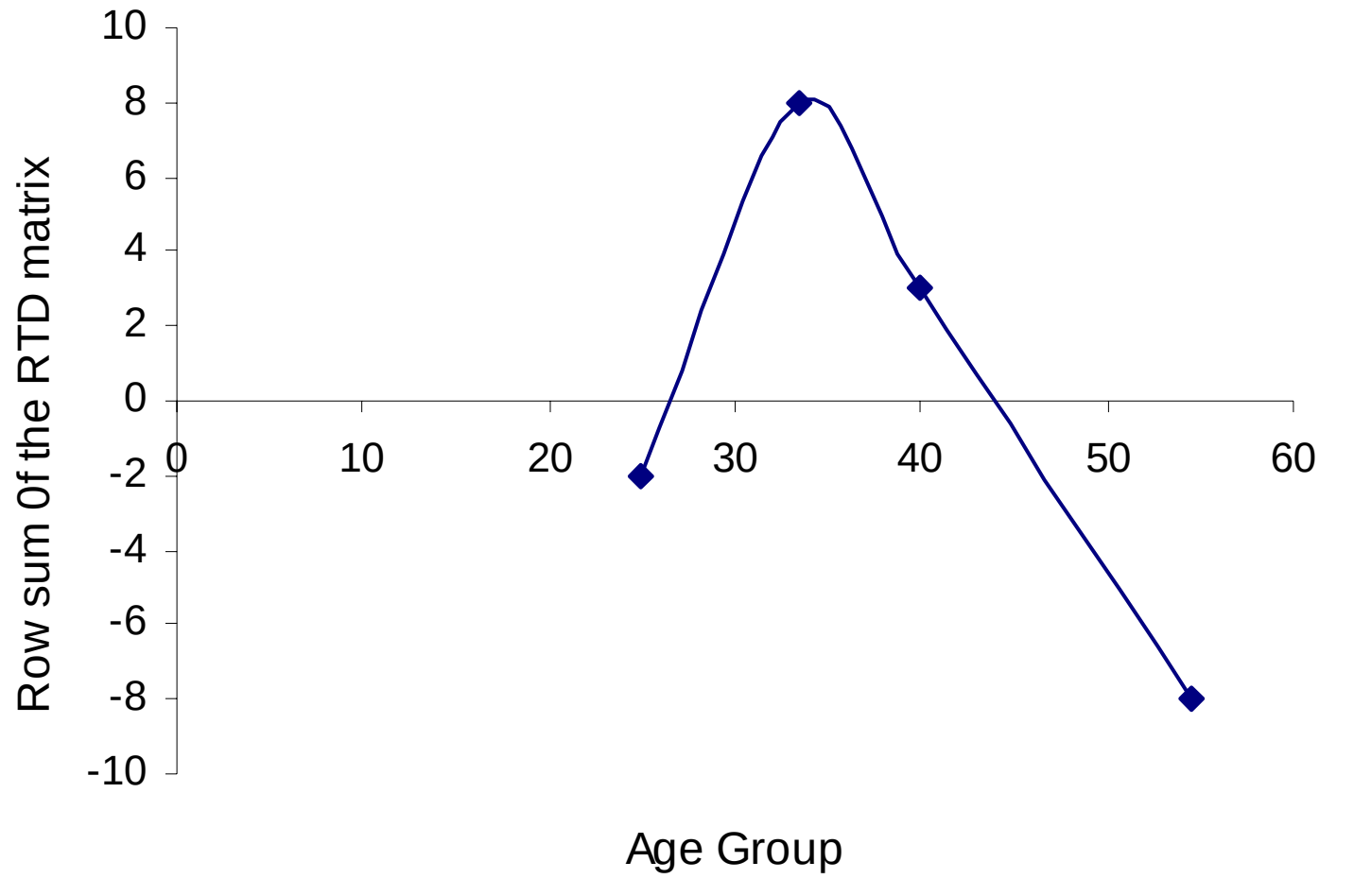

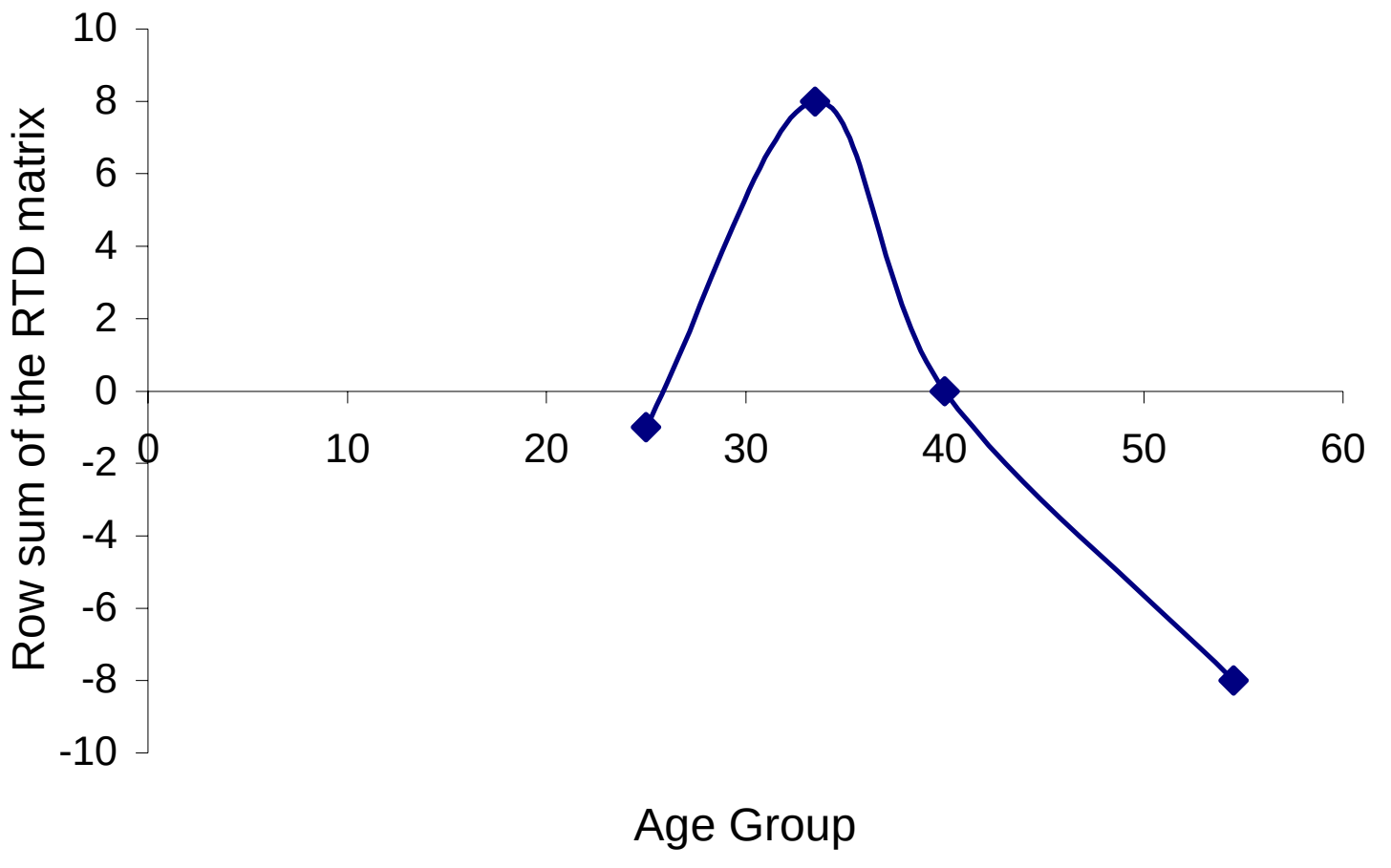

Graph 44: Graph depicting maximum age group of labourers with pollution related Nervous Problem for $\alpha = 0.35$

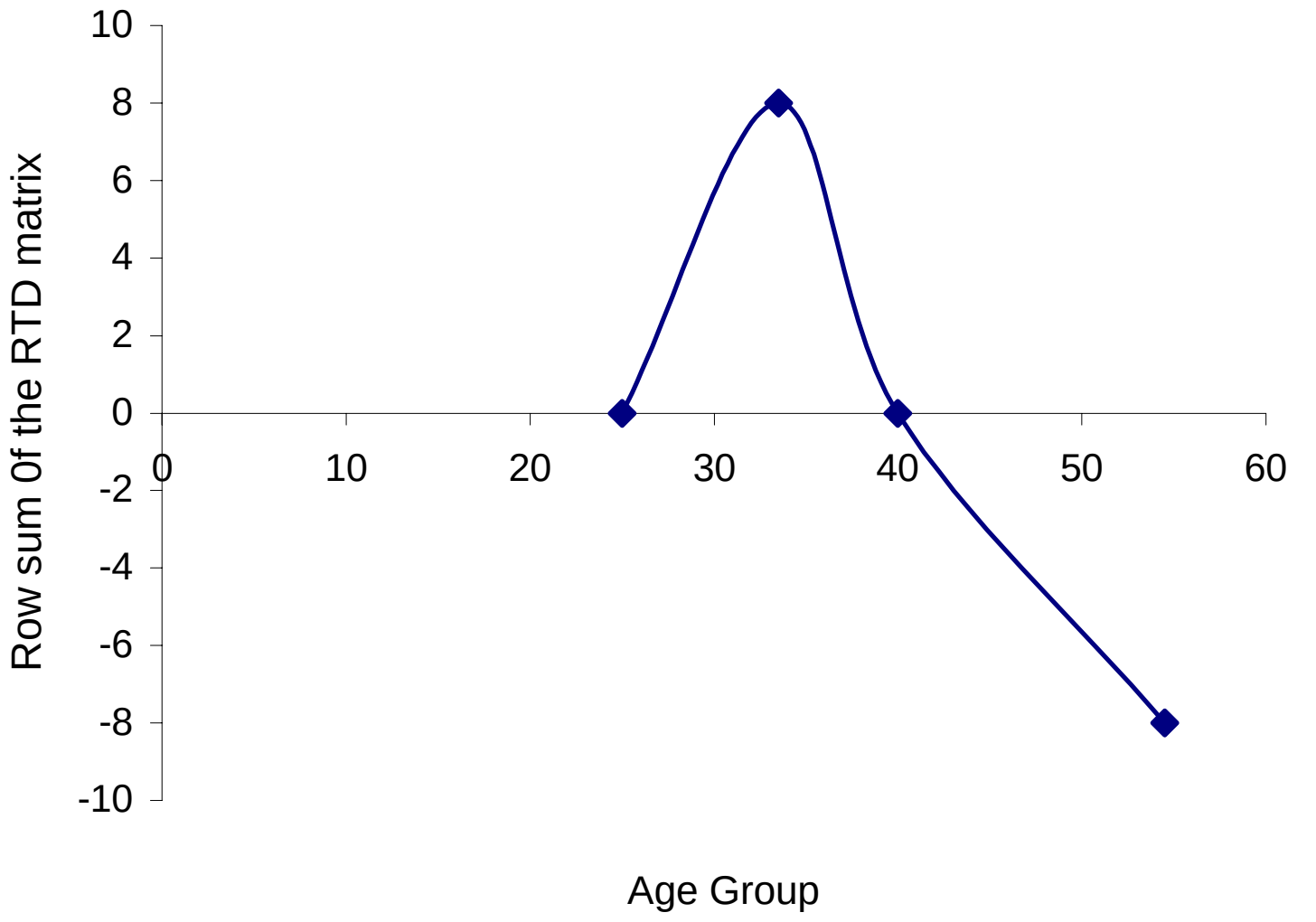

Graph 45: Graph depicting maximum age group of labourers with pollution related Nervous Problem for $\alpha = 0.45$

Graph 46: Graph depicting maximum age group of labourers with pollution related Nervous Problem for $\alpha = 0.75$

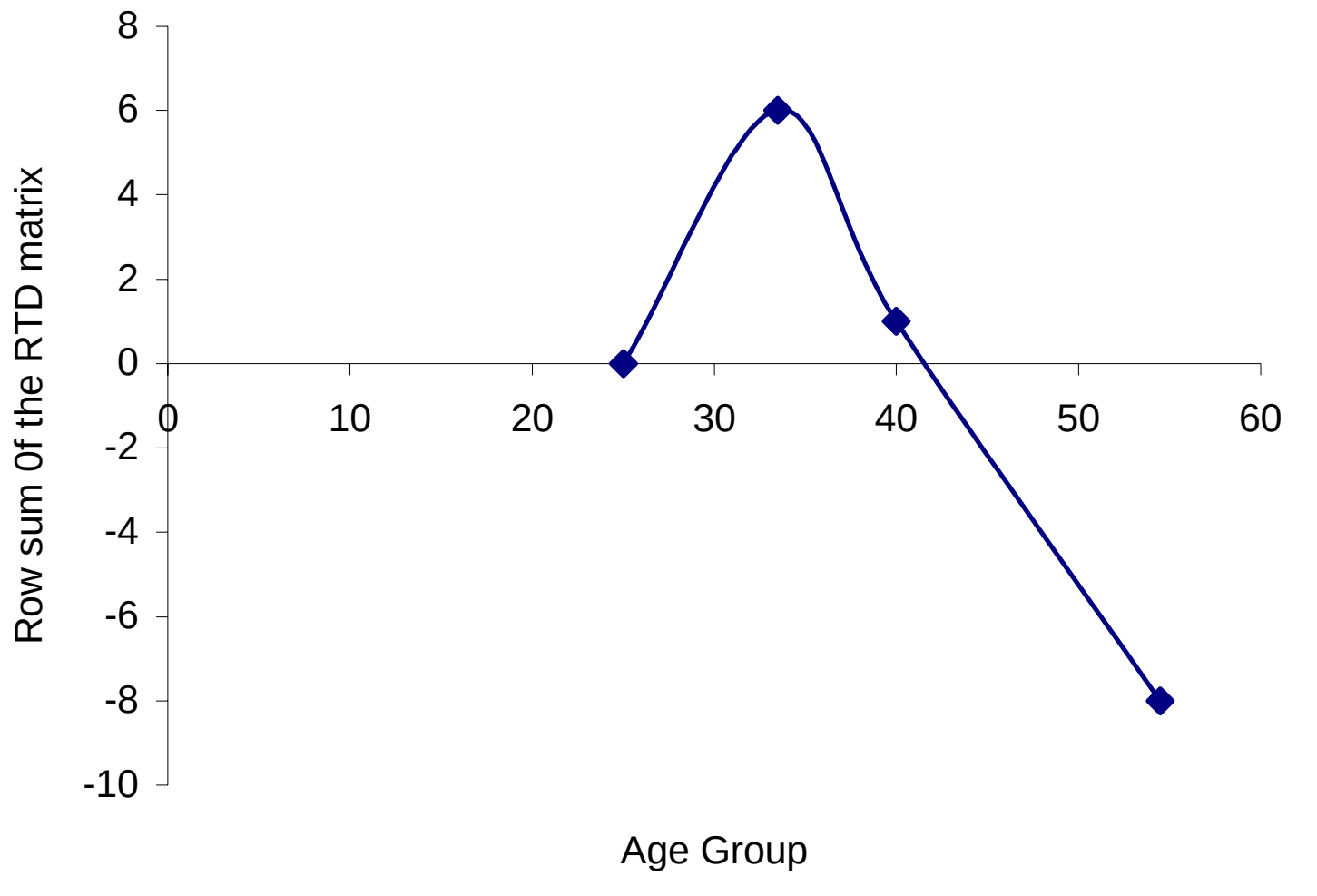

Graph 47: Graphical comparison of the maximum age group of the labourers who are affected by Nervous problem

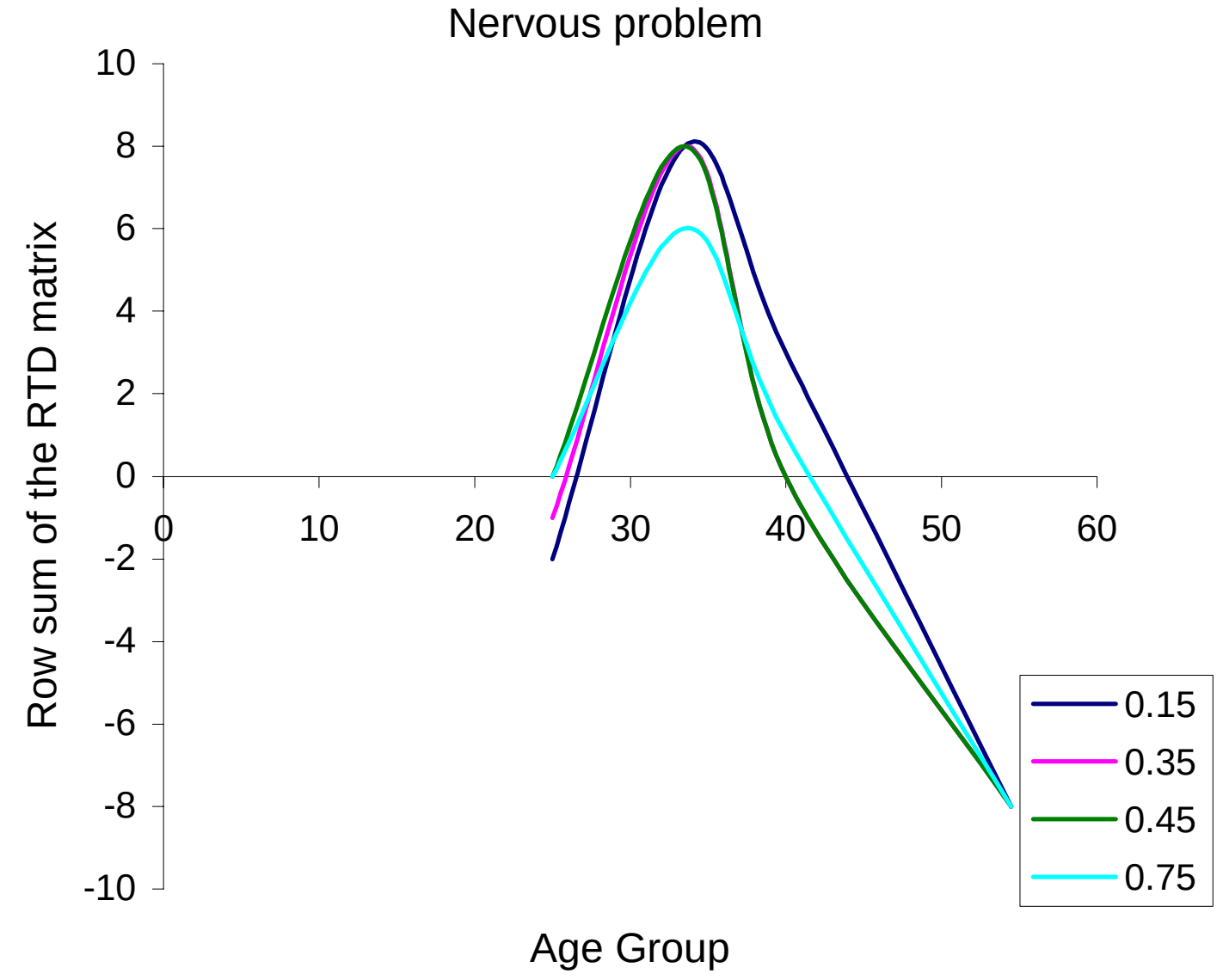

CETD matrix Row sum matrix

$$\begin{bmatrix} 0 & 0 & -1 & 0 & 1 & -1 & -2 & 0 \\ 4 & 4 & 3 & 4 & 4 & 4 & 3 & 4 \\ -4 & 1 & 4 & 1 & 1 & 0 & 4 & -3 \\ -4 & -4 & -4 & -4 & -4 & -4 & -4 & -4 \end{bmatrix} \qquad \begin{bmatrix} -3 \\ 30 \\ 4 \\ -32 \end{bmatrix}$$

The graph for the CETD matrix is given below :

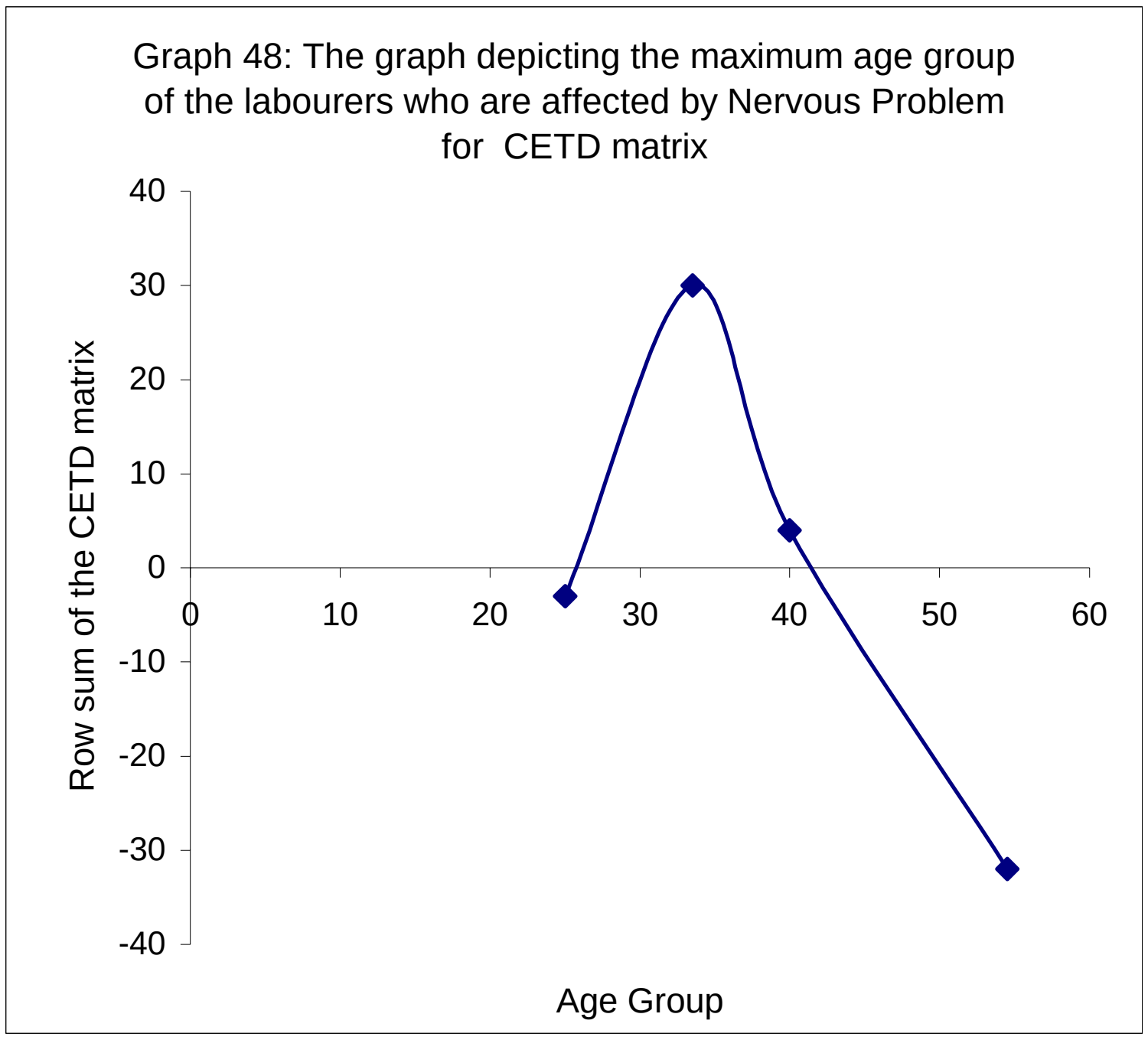

### 2.1.3.3 Estimation of maximum age group of agriculture labourers with pollution related Nervous problems using 5 × 8 matrix

Now to make study still more sensitive by increasing the number of rows by 5 and see whether the decision arrived is more sensitive to the earlier one we have discussed. Thus we give the raw data of 5 × 8 matrix.

The Initial Raw Data Matrix of
Nervous problem of order 5 × 8

| Years | $S_1$ | $S_2$ | $S_3$ | $S_4$ | $S_5$ | $S_6$ | $S_7$ | $S_8$ |
|---|---|---|---|---|---|---|---|---|
| 20-24 | 6 | 5 | 3 | 5 | 5 | 5 | 3 | 3 |
| 25-30 | 22 | 22 | 12 | 22 | 20 | 15 | 10 | 13 |
| 31-36 | 22 | 22 | 12 | 21 | 18 | 19 | 10 | 15 |
| 37-43 | 16 | 20 | 15 | 18 | 16 | 15 | 13 | 8 |
| 44-65 | 23 | 23 | 17 | 22 | 20 | 24 | 19 | 15 |

The ATD Matrix of
Nervous problem of order 5 × 8

| Years | $S_1$ | $S_2$ | $S_3$ | $S_4$ | $S_5$ | $S_6$ | $S_7$ | $S_8$ |
|---|---|---|---|---|---|---|---|---|
| 20-24 | 1.2 | 1 | 0.6 | 1 | 1 | 1 | 0.6 | 0.6 |
| 25-30 | 3.67 | 3.67 | 2 | 3.67 | 3.33 | 2.5 | 1.67 | 2.17 |
| 31-36 | 3.67 | 3.67 | 2 | 3.5 | 3 | 3.17 | 1.67 | 2.5 |
| 37-43 | 2.29 | 2.86 | 2.14 | 2.57 | 2.29 | 2.14 | 1.86 | 1.14 |
| 44-65 | 1.05 | 1.05 | 0.77 | 1 | 0.91 | 1.09 | 0.86 | 0.86 |

The Average and Standard Deviation
of the above given ATD matrix

| | | | | | | | | |
|---|---|---|---|---|---|---|---|---|
| Average | 2.38 | 2.45 | 1.50 | 2.35 | 2.11 | 1.98 | 1.33 | 1.45 |
| Standard Deviation | 1.28 | 1.34 | 0.75 | 1.3 | 1.12 | 0.93 | 0.56 | 0.83 |

RTD matrix for $\alpha = 0.1$     Row sum matrix

$$\begin{bmatrix} -1 & -1 & -1 & -1 & -1 & -1 & -1 & -1 \\ 1 & 1 & 1 & 1 & 1 & 1 & 1 & 1 \\ 1 & 1 & 1 & 1 & 1 & 1 & 1 & 1 \\ 0 & 1 & 1 & 1 & 1 & 1 & 1 & -1 \\ -1 & -1 & -1 & -1 & -1 & -1 & -1 & -1 \end{bmatrix} \qquad \begin{bmatrix} -8 \\ 8 \\ 8 \\ 5 \\ -8 \end{bmatrix}$$

RTD matrix for $\alpha = 0.15$     Row sum matrix

$$\begin{bmatrix} -1 & -1 & -1 & -1 & -1 & -1 & -1 & -1 \\ 1 & 1 & 1 & 1 & 1 & 1 & 1 & 1 \\ 1 & 1 & 1 & 1 & 1 & 1 & 1 & 1 \\ 0 & 1 & 1 & 1 & 1 & 1 & 1 & -1 \\ -1 & -1 & -1 & -1 & -1 & -1 & -1 & -1 \end{bmatrix} \qquad \begin{bmatrix} -8 \\ 8 \\ 8 \\ 5 \\ -8 \end{bmatrix}$$

RTD matrix for $\alpha = 0.2$     Row sum matrix

$$\begin{bmatrix} -1 & -1 & -1 & -1 & -1 & -1 & -1 & -1 \\ 1 & 1 & 1 & 1 & 1 & 1 & 1 & 1 \\ 1 & 1 & 1 & 1 & 1 & 1 & 1 & 1 \\ 0 & 1 & 1 & 0 & 0 & 0 & 1 & -1 \\ -1 & -1 & -1 & -1 & -1 & -1 & -1 & -1 \end{bmatrix} \qquad \begin{bmatrix} -8 \\ 8 \\ 8 \\ 2 \\ -8 \end{bmatrix}$$

RTD matrix for $\alpha = 0.35$     Row sum matrix

$$\begin{bmatrix} -1 & -1 & -1 & -1 & -1 & -1 & -1 & -1 \\ 1 & 1 & 1 & 1 & 1 & 1 & 1 & 1 \\ 1 & 1 & 1 & 1 & 1 & 1 & 1 & 1 \\ 0 & 0 & 1 & 0 & 0 & 0 & 1 & -1 \\ -1 & -1 & -1 & -1 & -1 & -1 & -1 & -1 \end{bmatrix} \qquad \begin{bmatrix} -8 \\ 8 \\ 8 \\ 1 \\ -8 \end{bmatrix}$$

Graph 49: Graph depicting maximum age group of labourers pollution related Nervous Problem for $\alpha$ = 0.1

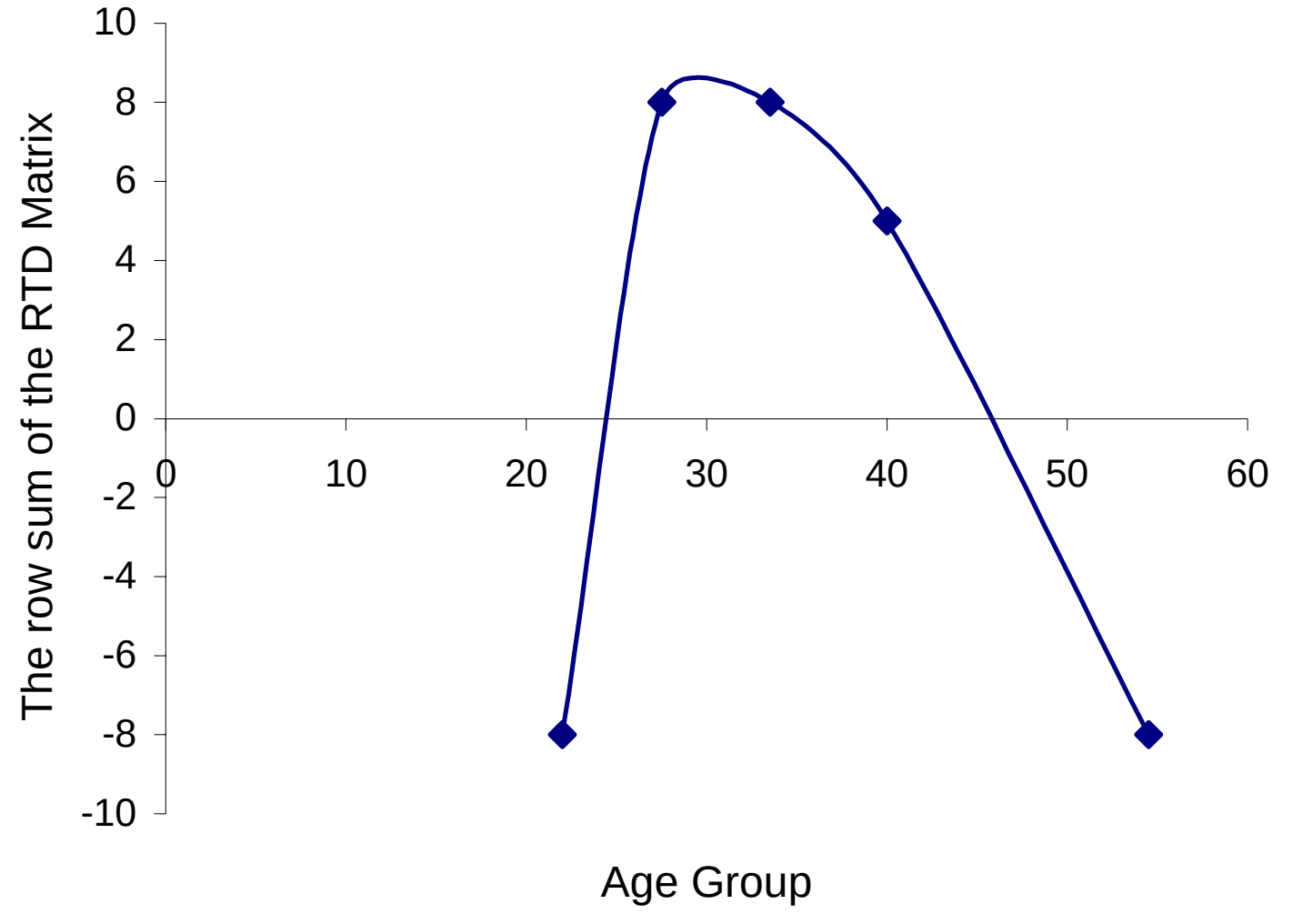

Graph 50: Graph depicting maximum age group of labourers pollution related Nervous Problem for $\alpha$ = 0.15

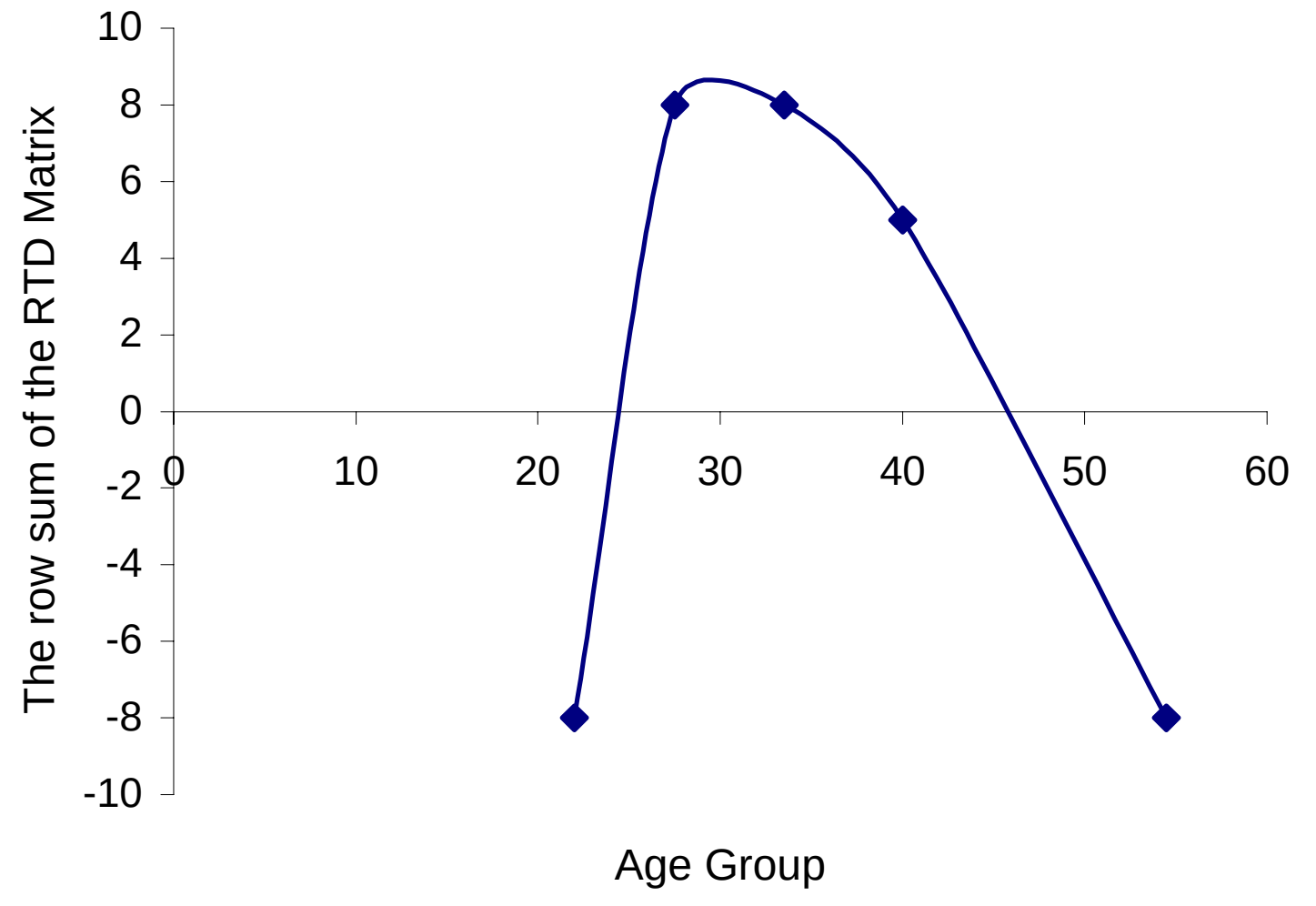

Graph 51: Graph depicting maximum age group of labourers pollution related Nervous Problem for $\alpha = 0.2$

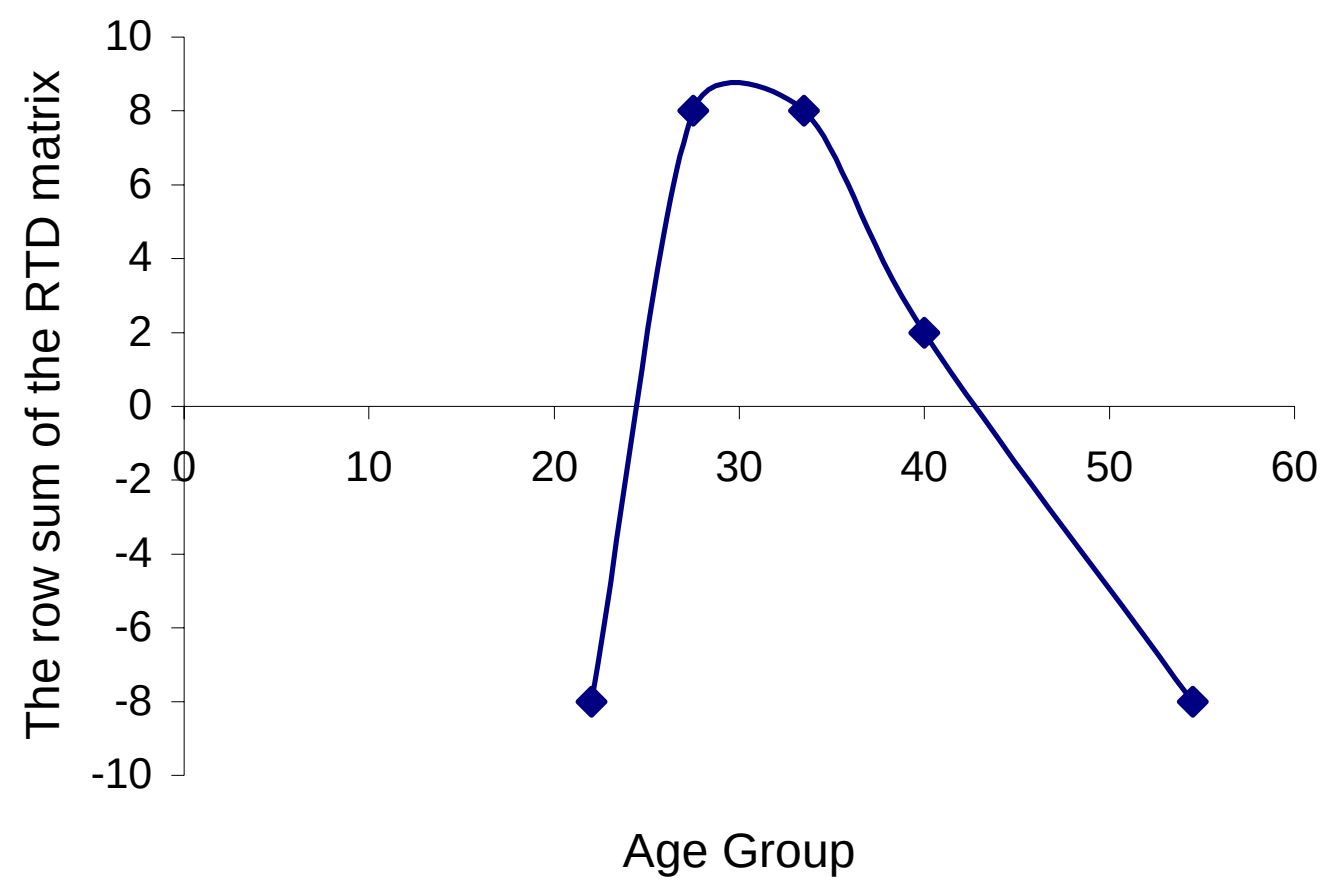

Graph 52: Graph depicting maximum age group of labourers pollution related Nervous Problem for $\alpha = 0.35$

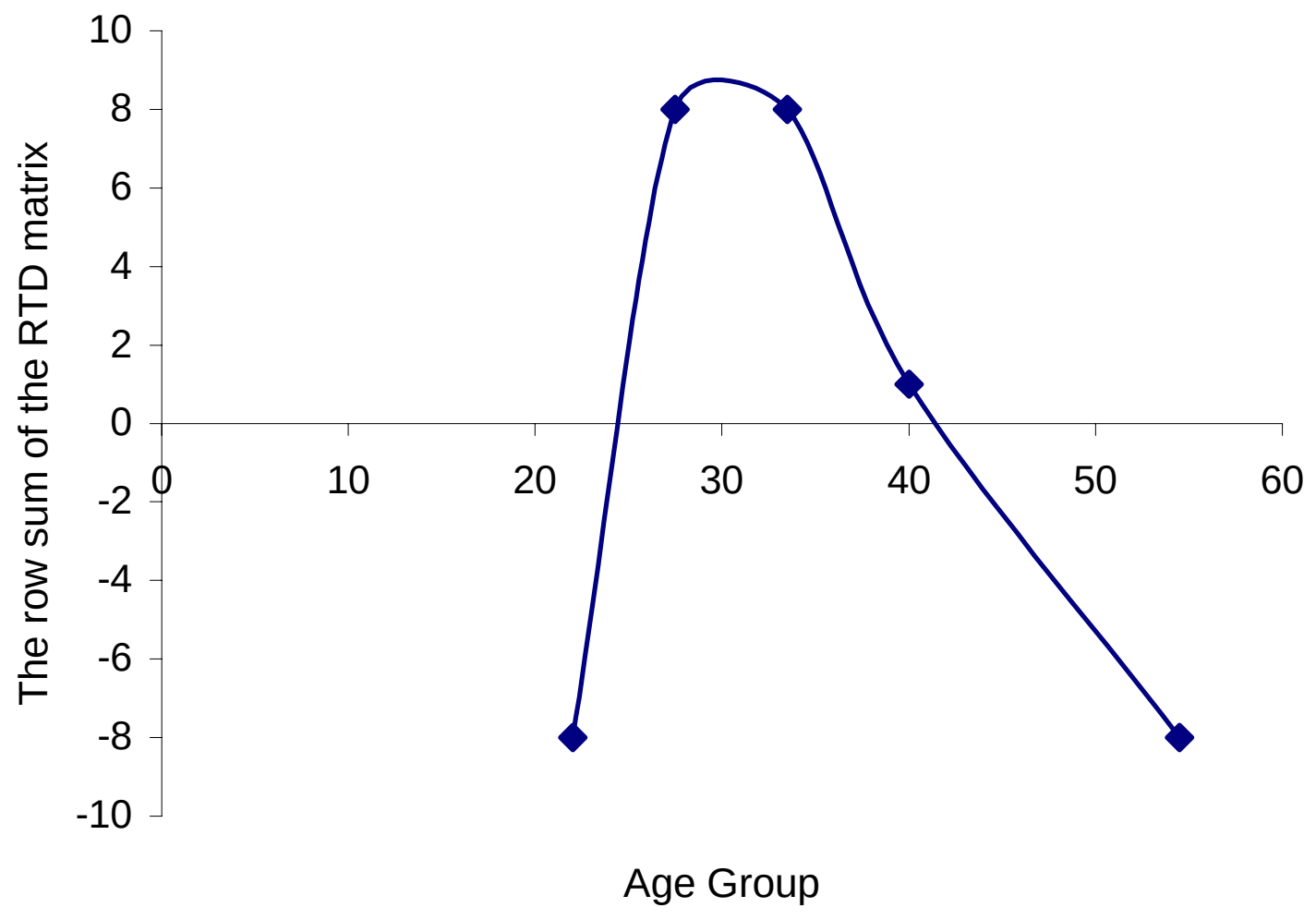

The comparative graph of the maximum age group of the labourers who are affected by nervous problem

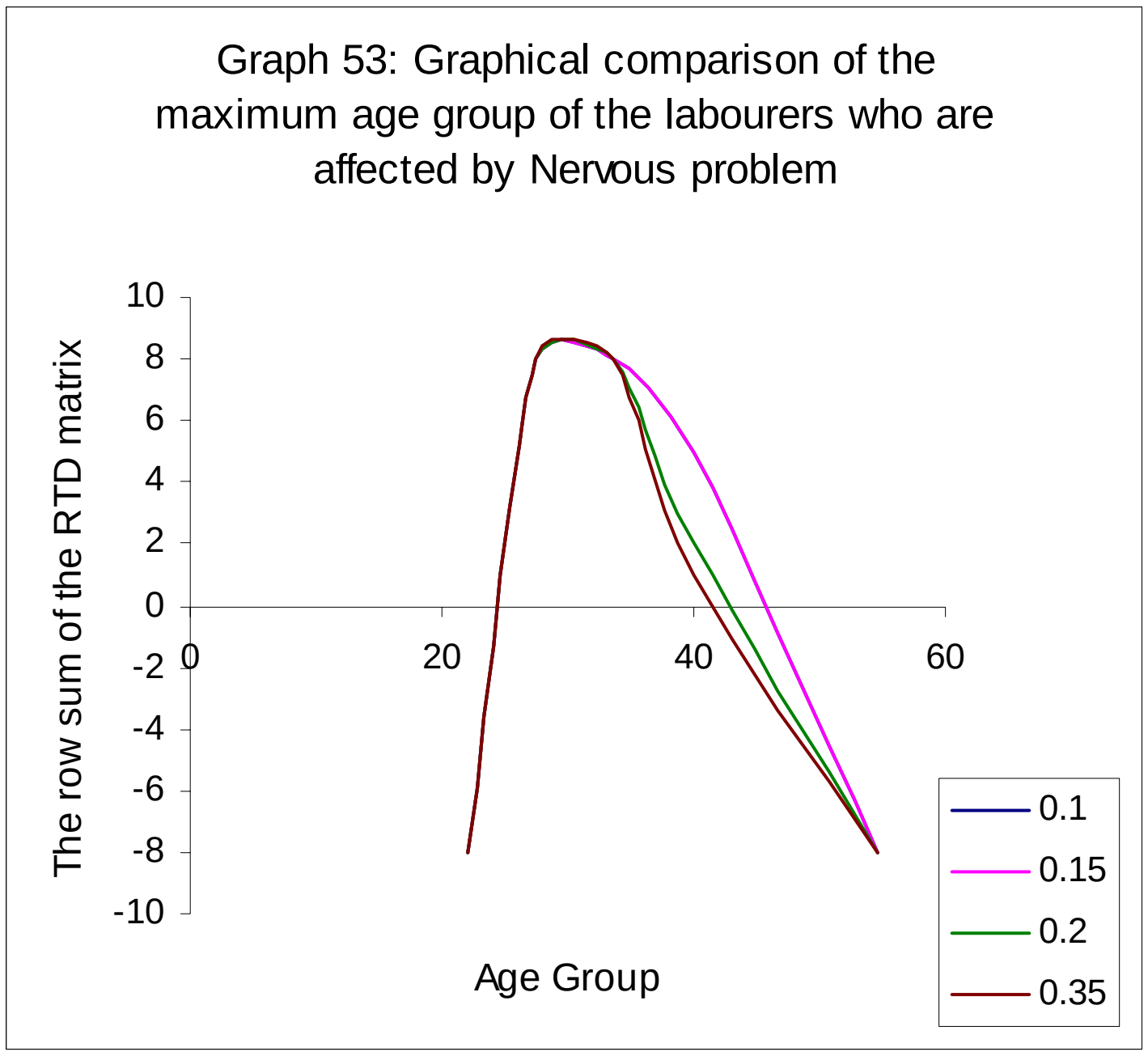

The CETD matrix and row sum matrix are given below:

The CETD matrix Row sum matrix

$$\begin{bmatrix} -4 & -4 & -4 & -4 & -4 & -4 & -4 & -4 \\ 4 & 4 & 4 & 4 & 4 & 4 & 4 & 4 \\ 4 & 4 & 4 & 4 & 4 & 4 & 4 & 4 \\ 0 & 3 & 4 & 2 & 2 & 2 & 4 & -4 \\ -4 & -4 & -4 & -4 & -4 & -4 & -4 & -4 \end{bmatrix} \qquad \begin{bmatrix} -32 \\ 32 \\ 32 \\ 13 \\ -32 \end{bmatrix}$$

The graph for the CETD matrix is given in the next page.

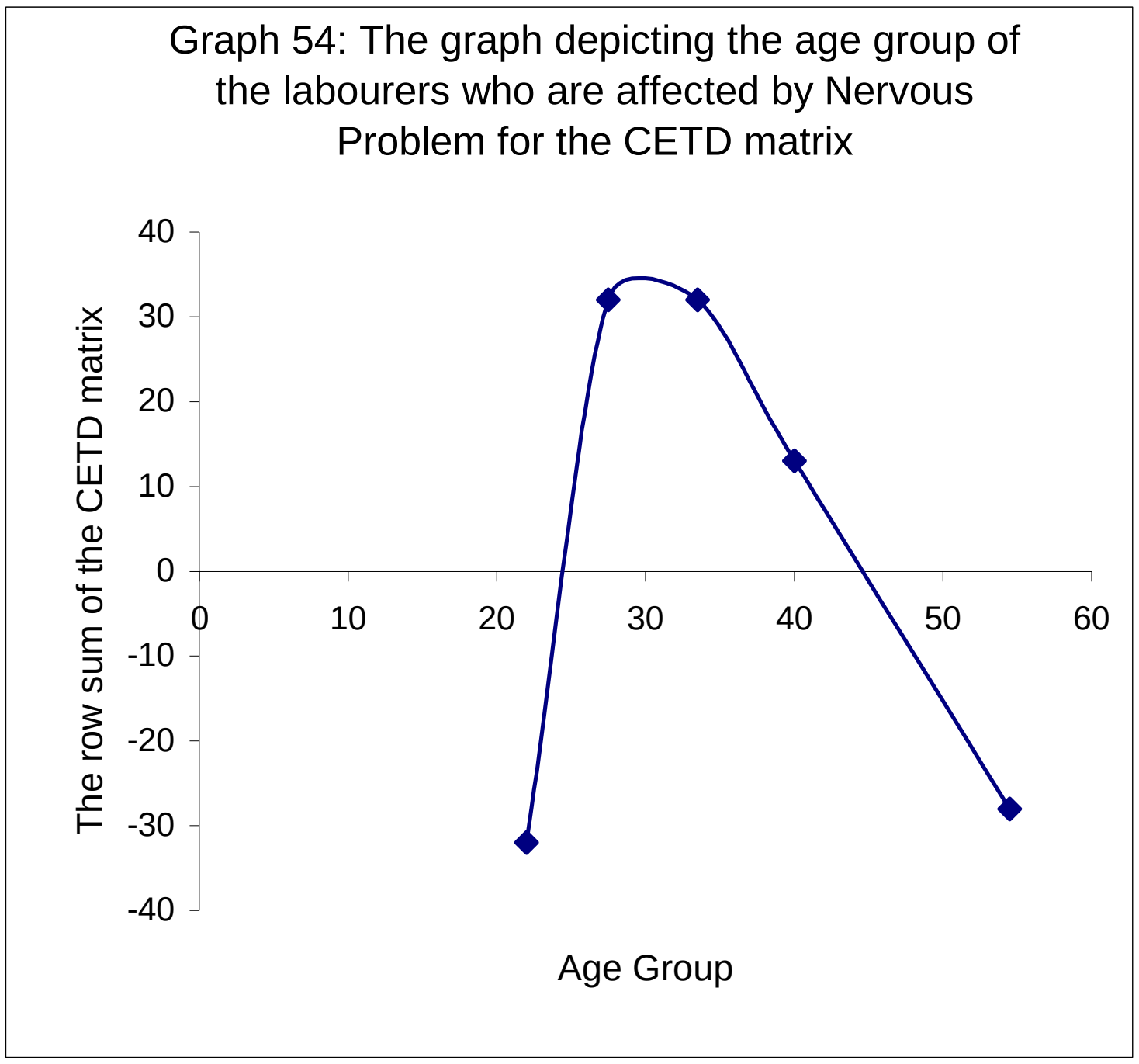

## 2.1.3.4 Conclusion

Strain and stress are often cited as the reason for various nervous problems. Surveys have thrown light into this problem, by explaining the fact that people like building contractors, musicians, technicians, surgeons, computer analysts, businessmen are those who face stress in their day-to-day life and hence suffer from various nervous problems. Nowhere in the existing literature are the farmers mentioned. But the graph picturises that age group between 31-43 have the neural symptom and maximum is reached at the age of 32. Out of 110 interviewed 89 showed headache and in that 89, 38 fall in the age group of 31-43. 92 of them showed fainting and in that 42 of them fall in the same age range. A similar majority is found in the same age range with other symptoms like dizziness, developing fits, loss of memory and vomiting. One can conclude from this graph that this age-range of people who are directly involved with spraying of pesticides remain in contact with the same over an extended period, ultimately they develop such nervous problem.

### 2.1.4 Estimation of The maximum age group of the Agricultural Labourers having Respiratory Problem Due To Chemical Pollution, Using Matrices

The Respiratory problem is taken under six symptom/ disease viz.

$S_1$ - Asthma,
$S_2$ - Allergy,
$S_3$ - Running nose,
$S_4$ - Sneezing,
$S_5$ - Coughing,
$S_6$ - Breathing difficulty .

Taking symptoms along with rows and age along the columns.

#### 2.1.4.1 Estimation of The Maximum Age Group Of The Agricultural Labourers Having Respiratory Problem due to Chemical Pollution, Using Matrices of order 3 × 6

The initial raw data matrix of respiratory problem of order 3 × 6

| Years | $S_1$ | $S_2$ | $S_3$ | $S_4$ | $S_5$ | $S_6$ |
|---|---|---|---|---|---|---|
| 20-30 | 5 | 15 | 18 | 26 | 24 | 21 |
| 31-43 | 13 | 26 | 30 | 35 | 33 | 35 |
| 44-65 | 3 | 14 | 17 | 22 | 23 | 23 |

The ATD Matrix of respiratory problem of order 3 × 6

| Years | $S_1$ | $S_2$ | $S_3$ | $S_4$ | $S_5$ | $S_6$ |
|---|---|---|---|---|---|---|
| 20-30 | 0.45 | 1.36 | 1.64 | 2.36 | 2.18 | 1.91 |
| 31-43 | 1 | 2 | 2.31 | 2.69 | 2.54 | 2.69 |
| 44-65 | 0.14 | 0.64 | 0.77 | 1 | 1.05 | 1.05 |

The Average and Standard Deviation
of the given ATD matrix

| | | | | | | |
|---|---|---|---|---|---|---|
| Average | 0.53 | 1.33 | 1.57 | 2.02 | 1.92 | 1.88 |
| Standard deviation | 0.44 | 0.68 | 0.77 | 0.9 | 0.78 | 0.82 |

RTD matrix for α = 0.15 Row sum matrix

$$\begin{bmatrix} 0 & 0 & 0 & 1 & 0 & 0 \\ 1 & 1 & 1 & 1 & 1 & 1 \\ -1 & -1 & -1 & -1 & -1 & -1 \end{bmatrix} \qquad \begin{bmatrix} 1 \\ 6 \\ -6 \end{bmatrix}$$

RTD matrix for α = 0.35 Row sum matrix

$$\begin{bmatrix} 0 & 0 & 0 & 1 & 0 & 0 \\ 1 & 1 & 1 & 1 & 1 & 1 \\ -1 & -1 & -1 & -1 & -1 & -1 \end{bmatrix} \qquad \begin{bmatrix} 1 \\ 6 \\ -6 \end{bmatrix}$$

RTD matrix for α = 0.45 Row sum matrix

$$\begin{bmatrix} 0 & 0 & 0 & 0 & 0 & 0 \\ 1 & 1 & 1 & 1 & 1 & 1 \\ -1 & -1 & -1 & -1 & -1 & -1 \end{bmatrix} \qquad \begin{bmatrix} 0 \\ 6 \\ -6 \end{bmatrix}$$

RTD matrix for α = 0.75 Row sum matrix

$$\begin{bmatrix} 0 & 0 & 0 & 0 & 0 & 0 \\ 1 & 1 & 1 & 0 & 0 & 1 \\ -1 & -1 & -1 & -1 & -1 & -1 \end{bmatrix} \qquad \begin{bmatrix} 0 \\ 4 \\ -6 \end{bmatrix}$$

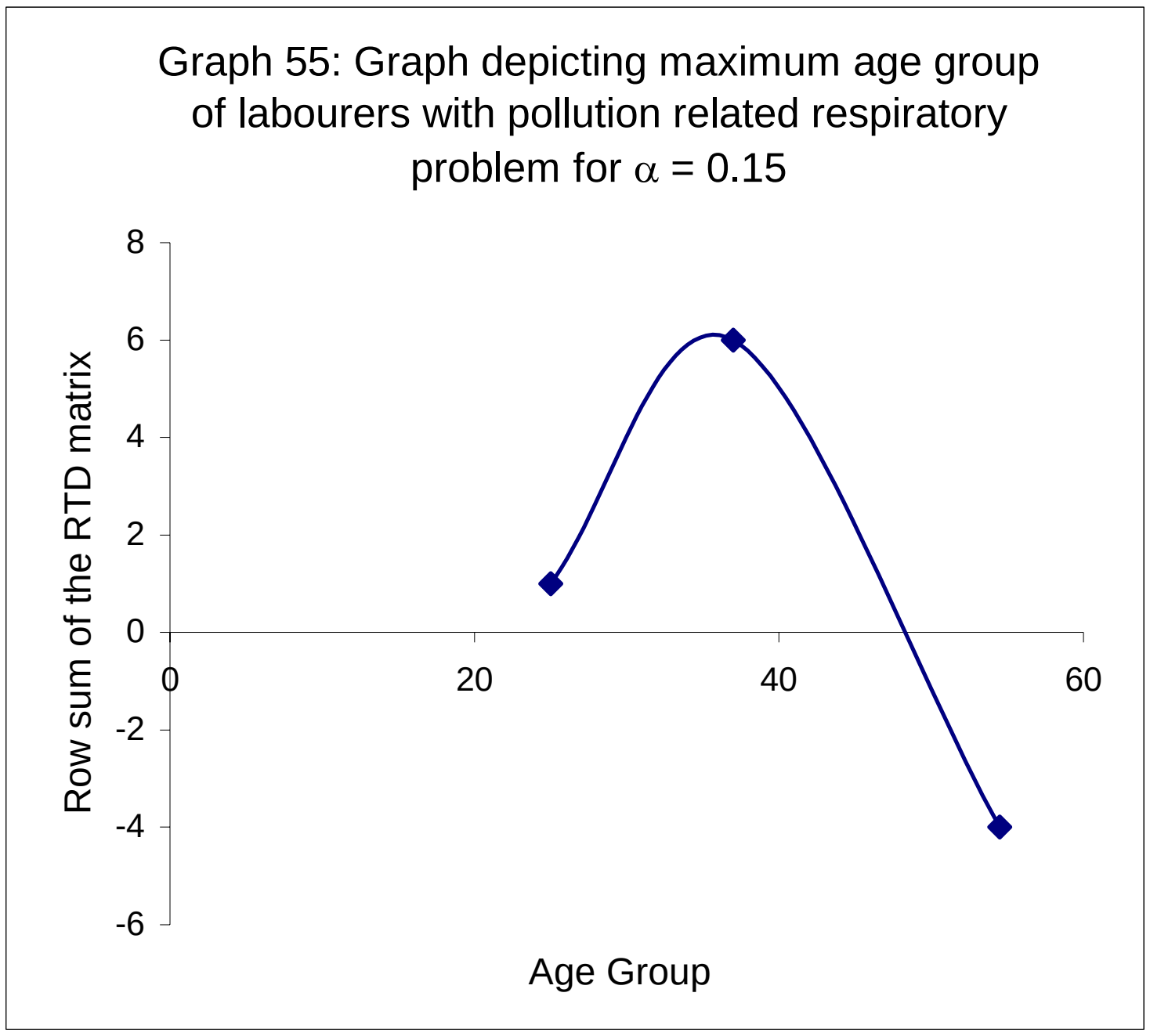
Graph 55: Graph depicting maximum age group of labourers with pollution related respiratory problem for α = 0.15
Row sum of the RTD matrix
Age Group
8
6
4
2
0
-2
-4
-6
0
20
40
60

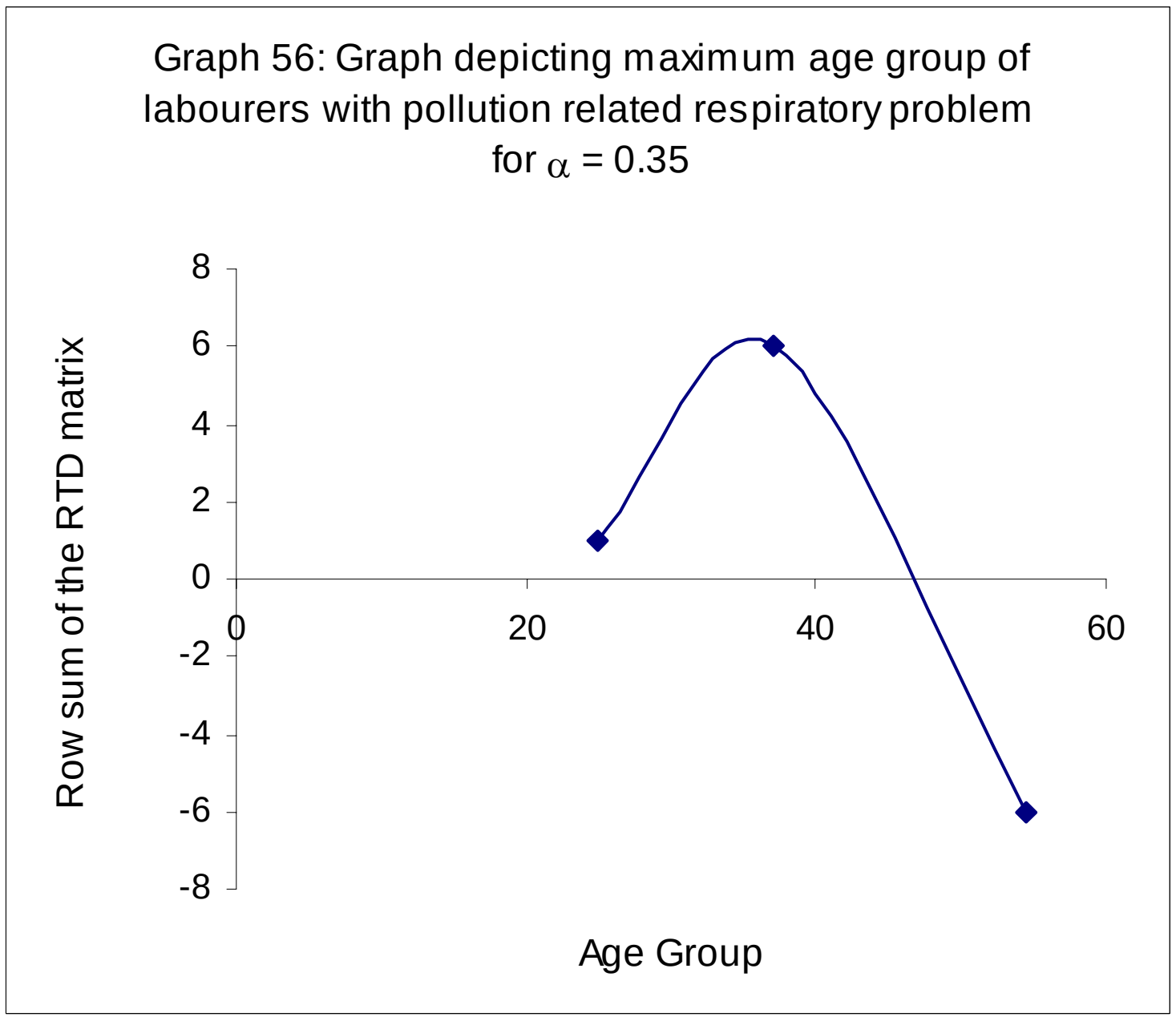
Graph 56: Graph depicting maximum age group of labourers with pollution related respiratory problem for α = 0.35
Row sum of the RTD matrix
Age Group
8
6
4
2
0
-2
-4
-6
-8
0
20
40
60

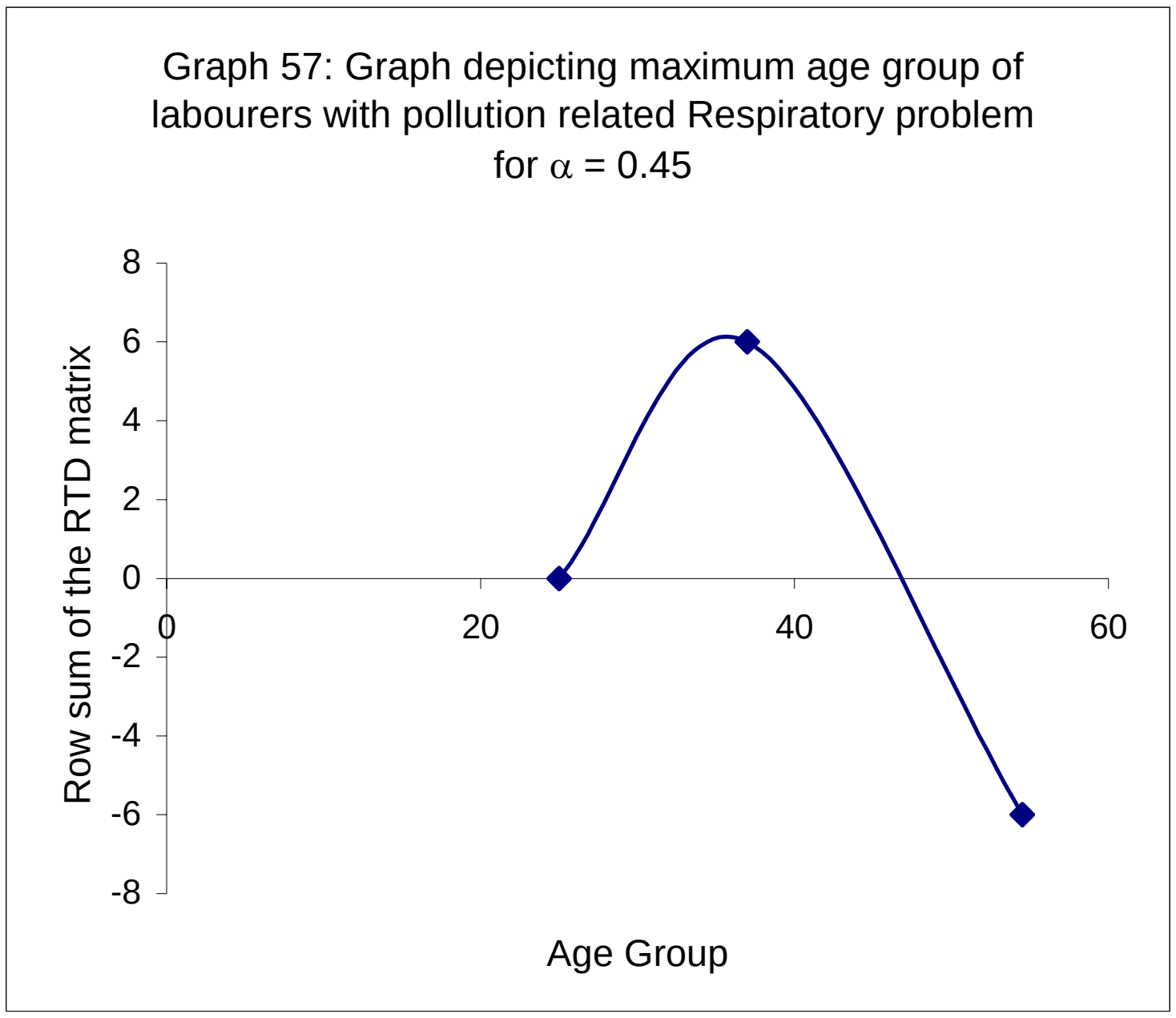
Graph 57: Graph depicting maximum age group of labourers with pollution related Respiratory problem for α = 0.45
Row sum of the RTD matrix
Age Group
8
6
4
2
0
-2
-4
-6
-8
0
20
40
60

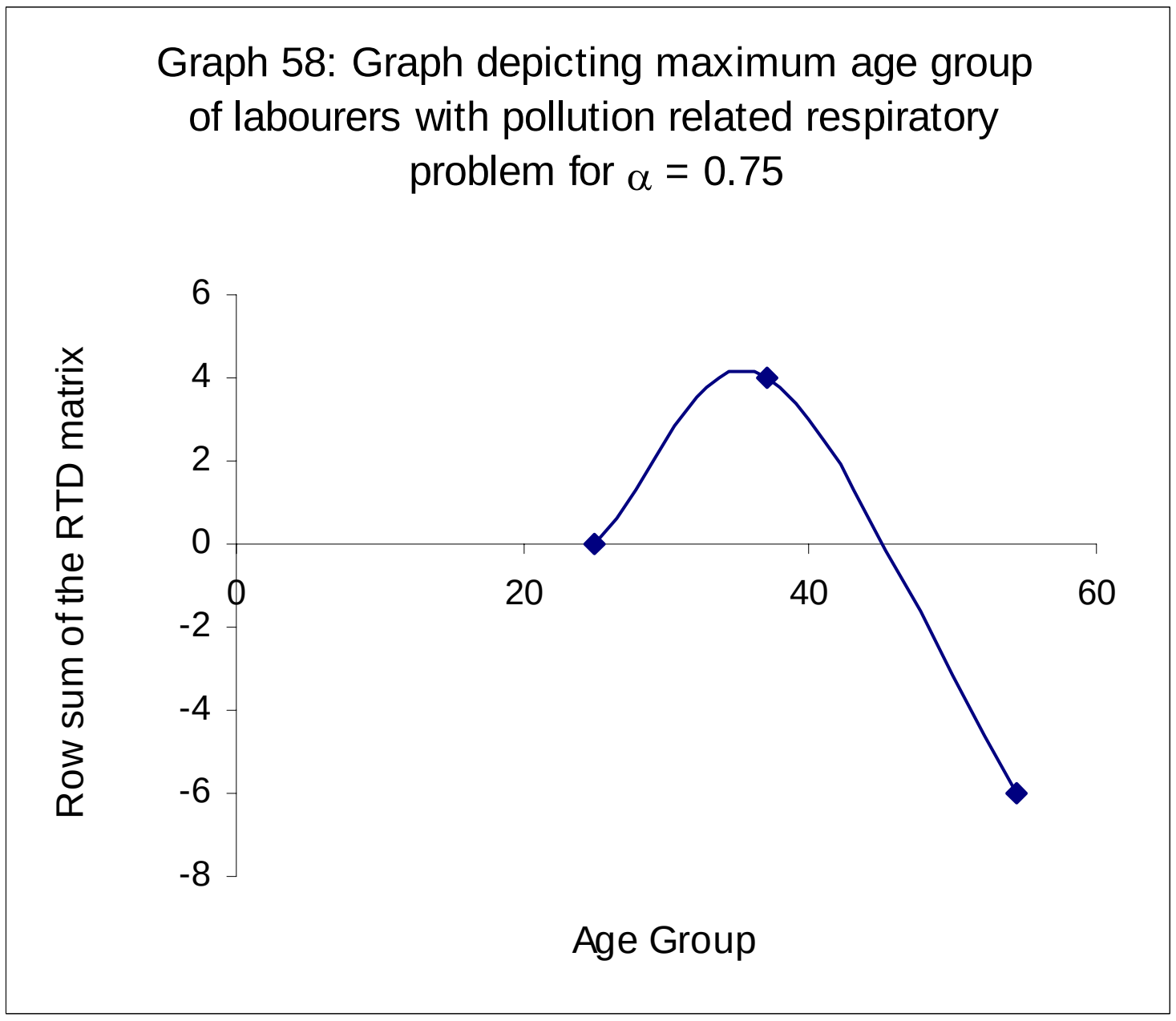
Graph 58: Graph depicting maximum age group of labourers with pollution related respiratory problem for α = 0.75
Row sum of the RTD matrix
Age Group
6
4
2
0
-2
-4
-6
-8
0
20
40
60

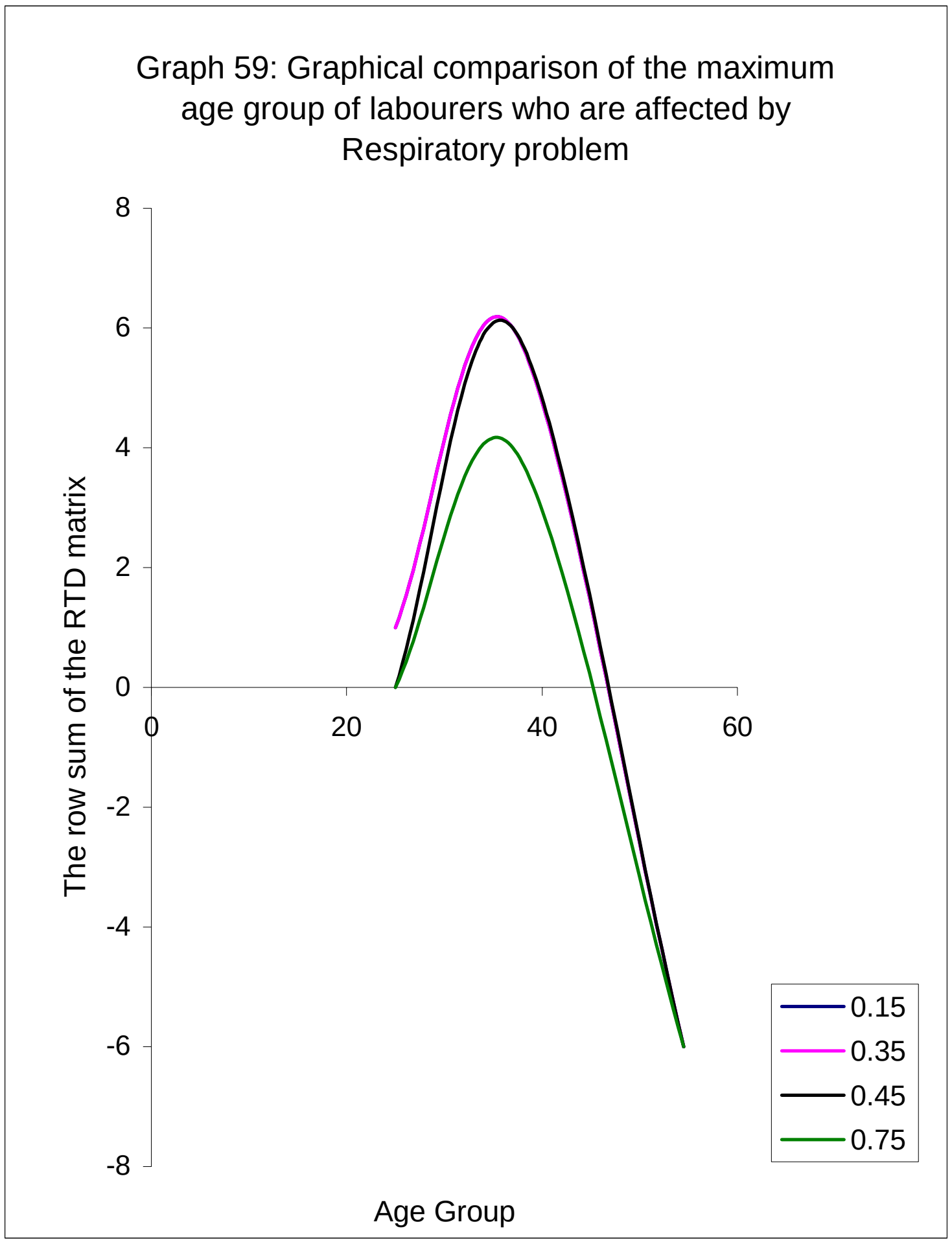

CETD matrix

$$\begin{bmatrix} 0 & 0 & 0 & 2 & 0 & 0 \\ 4 & 4 & 4 & 3 & 3 & 4 \\ -4 & -4 & -4 & -4 & -4 & -4 \end{bmatrix}$$

Row sum matrix

$$\begin{bmatrix} 2 \\ 22 \\ -24 \end{bmatrix}$$

The graph related to the above given CETD matrix is given in the next page:

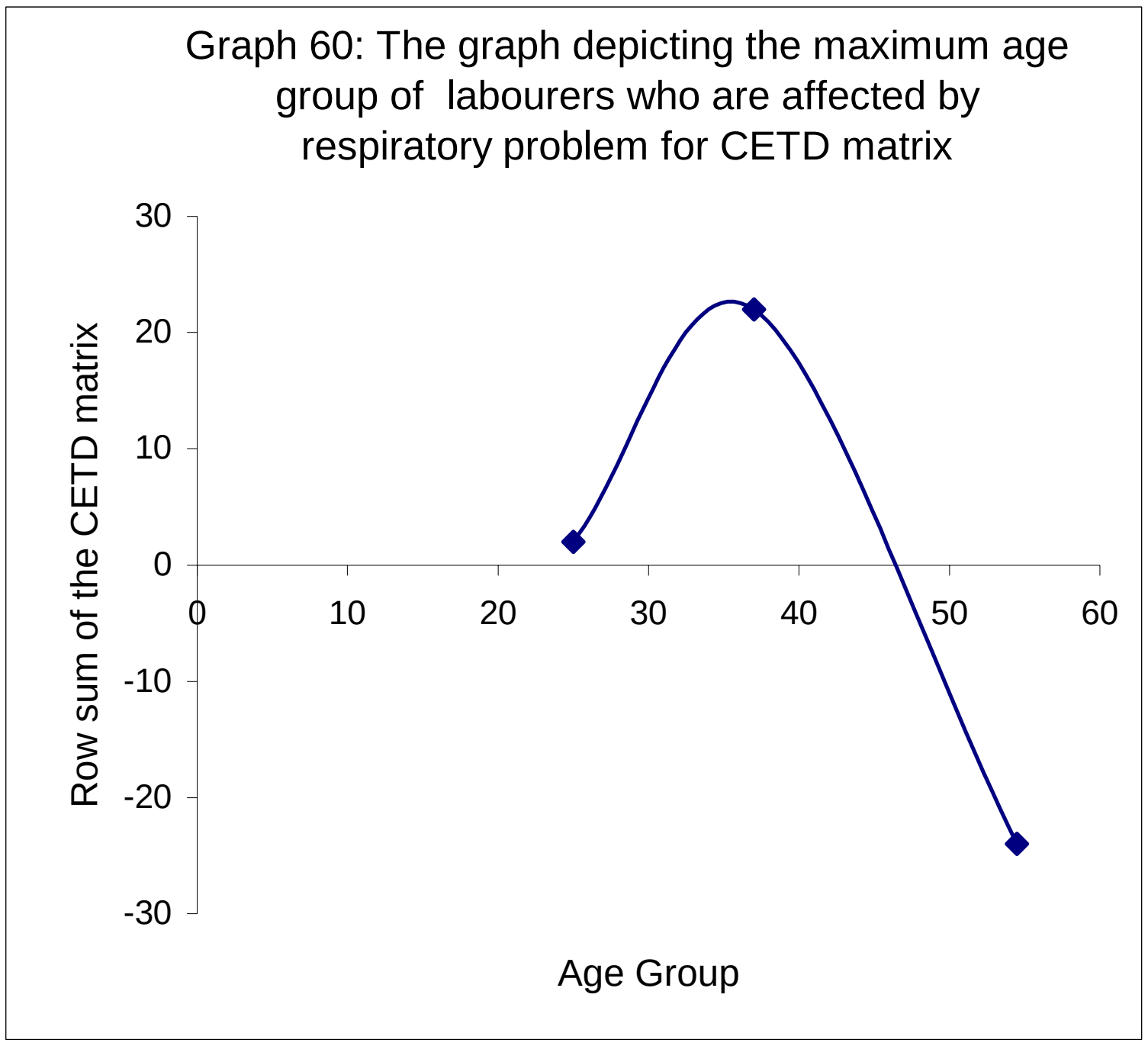

**2.1.4.2 Estimation of maximum age group of agricultural labourer with pollution related respiratory problem using 4 × 6 matrices**

Now to make the study more sensitive we increase the number of rows by 4 and see whether the decision arrived is more sensitive to the earlier one we have discussed. Thus we give the raw data of 4 × 6 matrix.

Initial Raw Data Matrix of respiratory problem of order 4 × 6

| Years | $S_1$ | $S_2$ | $S_3$ | $S_4$ | $S_5$ | $S_6$ |
|---|---|---|---|---|---|---|
| 20-30 | 5 | 15 | 18 | 26 | 24 | 21 |
| 31-36 | 7 | 16 | 13 | 17 | 15 | 16 |
| 37-43 | 6 | 10 | 17 | 18 | 18 | 19 |
| 44-65 | 3 | 14 | 17 | 22 | 23 | 23 |

The ATD Matrix of
Respiratory Problem of order 4 × 6

| Years | $S_1$ | $S_2$ | $S_3$ | $S_4$ | $S_5$ | $S_6$ |
|---|---|---|---|---|---|---|
| 20-30 | 0.45 | 1.36 | 1.64 | 2.36 | 2.18 | 1.91 |
| 31-36 | 1.17 | 2.67 | 2.17 | 2.83 | 2.5 | 2.67 |
| 37-43 | 0.86 | 1.43 | 2.43 | 2.57 | 2.57 | 2.71 |
| 44-65 | 0.14 | 0.64 | 0.77 | 1 | 1.05 | 1.05 |

The Average and Standard Deviation
of above given ATD matrix

| | | | | | | |
|---|---|---|---|---|---|---|
| Average | 0.66 | 1.53 | 1.75 | 2.19 | 2.08 | 2.09 |
| Standard deviation | 0.45 | 0.84 | 0.73 | 0.82 | 0.70 | 0.78 |

RTD matrix for $\alpha = 0.15$ Row sum matrix

$$\begin{bmatrix} -1 & -1 & 0 & 1 & 0 & -1 \\ 1 & 1 & 1 & 1 & 1 & 1 \\ 1 & 0 & 1 & 1 & 1 & 1 \\ -1 & -1 & -1 & -1 & -1 & -1 \end{bmatrix} \qquad \begin{bmatrix} -2 \\ 6 \\ 5 \\ -6 \end{bmatrix}$$

RTD matrix for $\alpha = 0.35$ Row sum matrix

$$\begin{bmatrix} -1 & -1 & 0 & 1 & 0 & 0 \\ 1 & 1 & 1 & 1 & 1 & 1 \\ 1 & 0 & 1 & 1 & 1 & 1 \\ -1 & -1 & -1 & -1 & -1 & -1 \end{bmatrix} \qquad \begin{bmatrix} -1 \\ 6 \\ 5 \\ -6 \end{bmatrix}$$

RTD matrix for $\alpha = 0.45$ Row sum matrix

$$\begin{bmatrix} -1 & 0 & 0 & 0 & 0 & 0 \\ 1 & 1 & 1 & 1 & 1 & 1 \\ 0 & 0 & 1 & 1 & 1 & 1 \\ -1 & -1 & -1 & -1 & -1 & -1 \end{bmatrix} \qquad \begin{bmatrix} -1 \\ 6 \\ 4 \\ -6 \end{bmatrix}$$

RTD matrix for $\alpha = 0.75$ Row sum matrix

$$\begin{bmatrix} 0 & 0 & 0 & 0 & 0 & 0 \\ 1 & 1 & 0 & 1 & 0 & 0 \\ 0 & 0 & 1 & 0 & 0 & 1 \\ -1 & -1 & -1 & -1 & -1 & -1 \end{bmatrix} \qquad \begin{bmatrix} 0 \\ 3 \\ 2 \\ -6 \end{bmatrix}$$

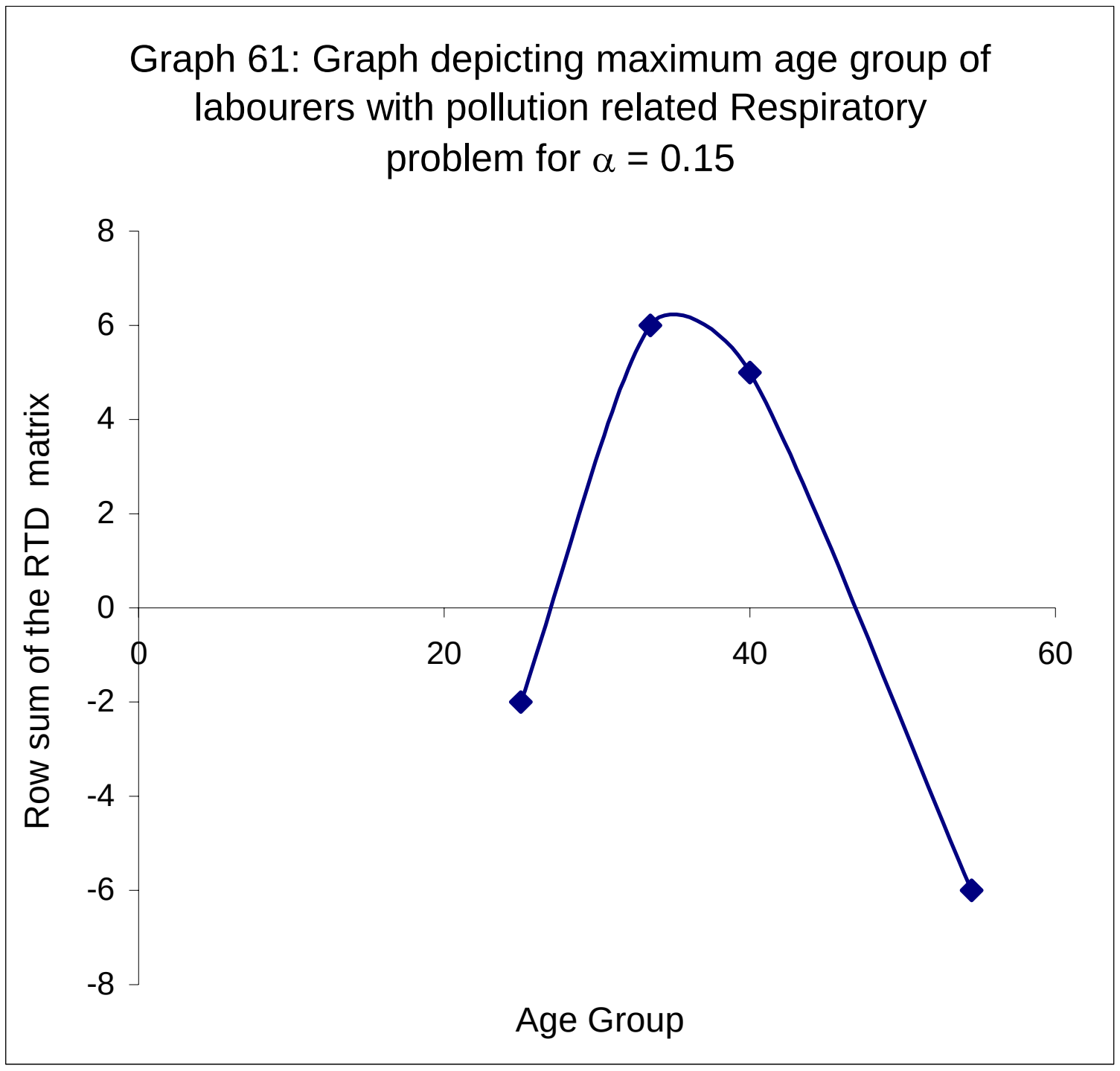

Graph 62: Graph depicting maximum age group of labourers with pollution related Respiratory problem for $\alpha$ = 0.35

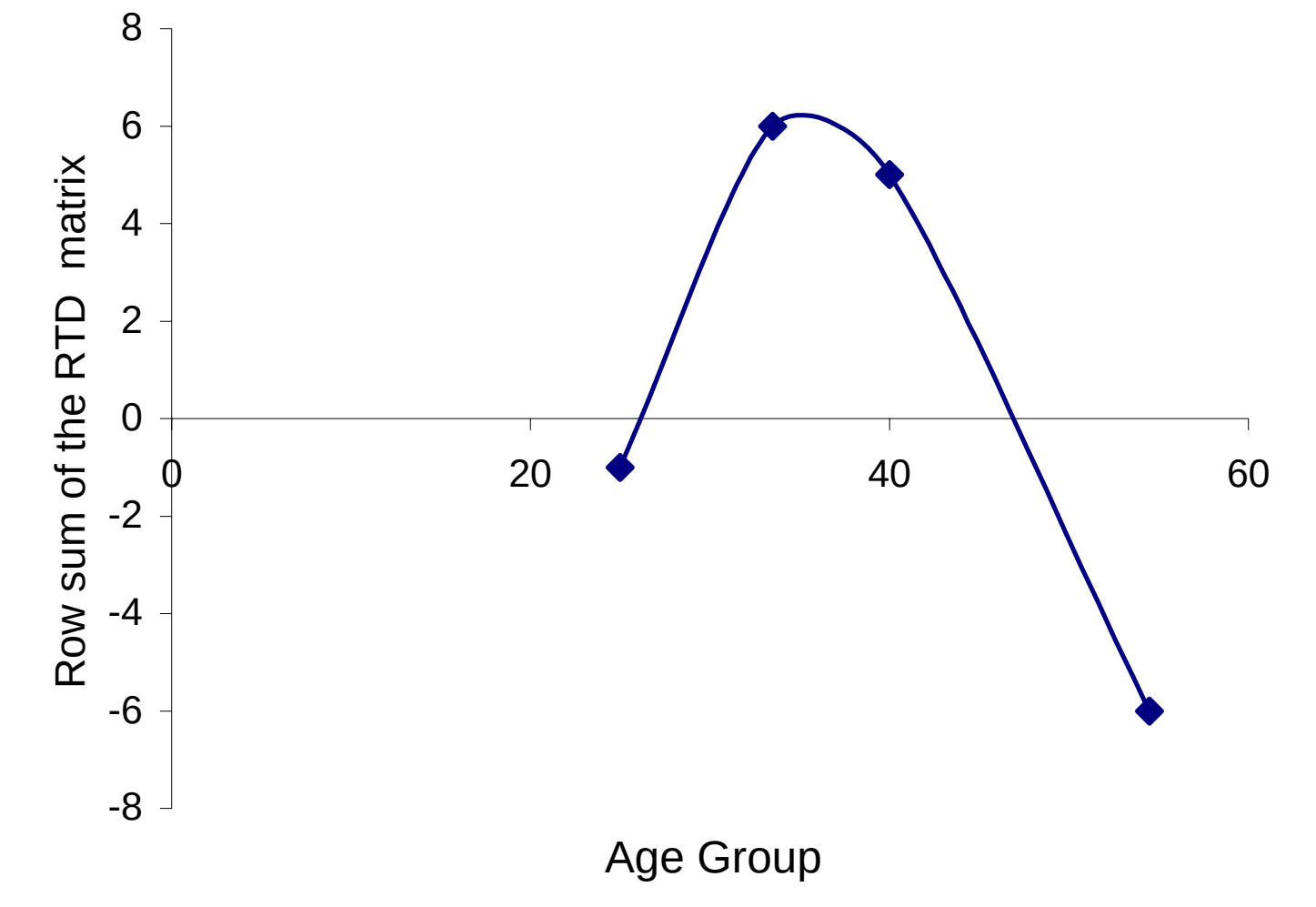

Graph 63: Graph depicting maximum age group of labourers with pollution related Respiratory problem for $\alpha$ = 0.45

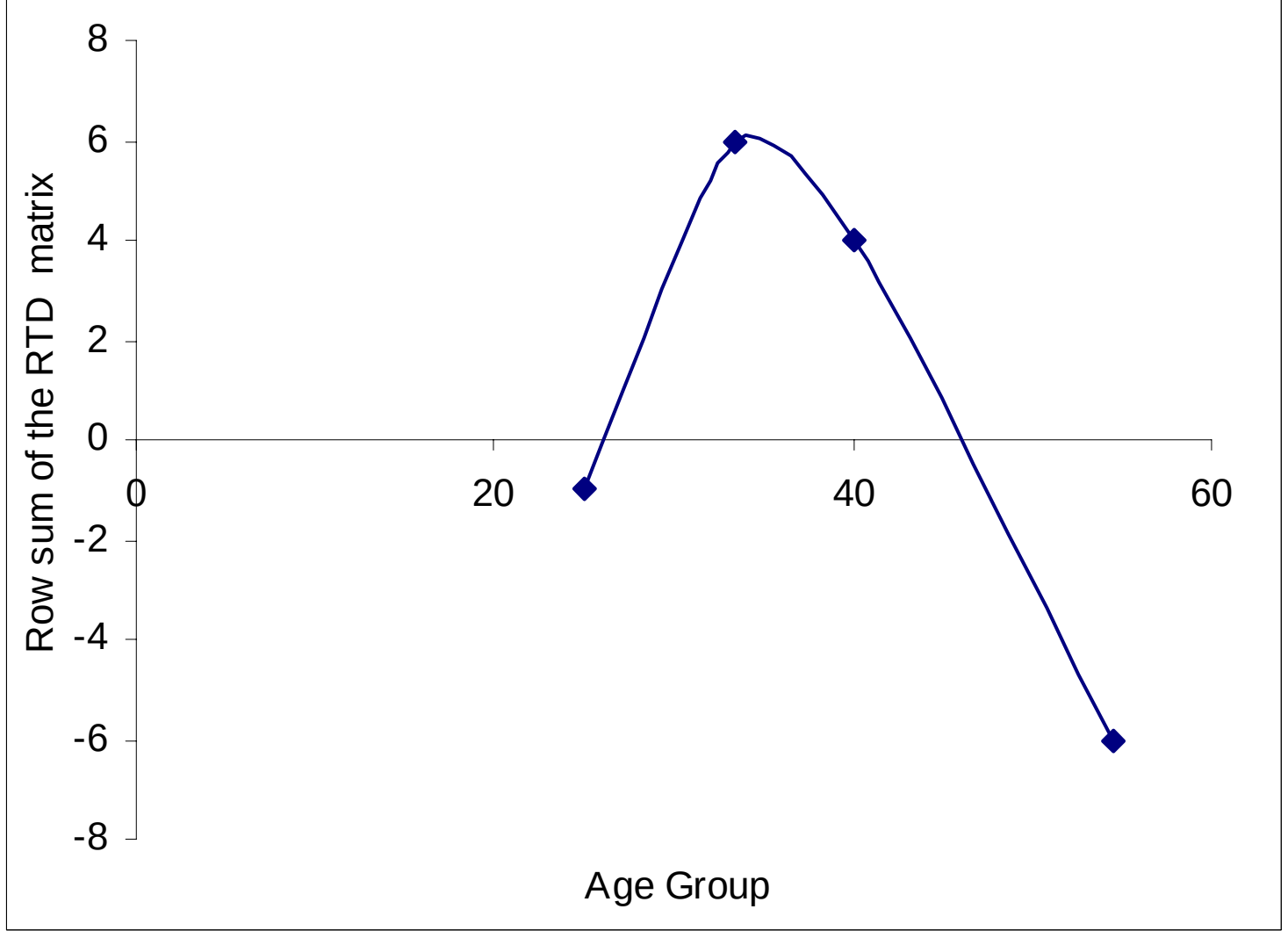

Graph 64: Graph depicting maximum age group of labourers with pollution related respiratory problem for $\alpha = 0.75$

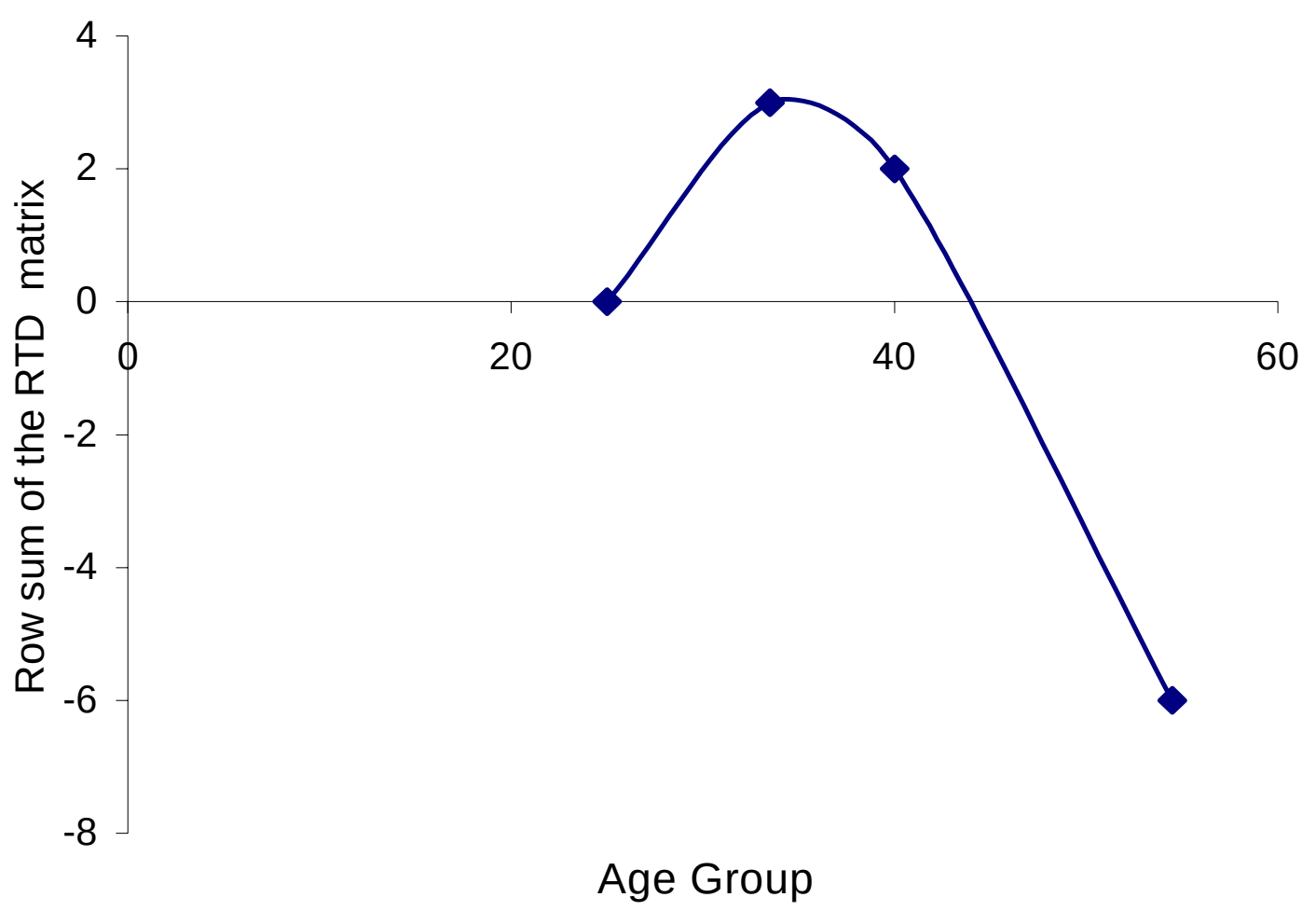

Graph 65: Graphical comparison of the maximum age group of the labourers who are affected due to Respiratory problem

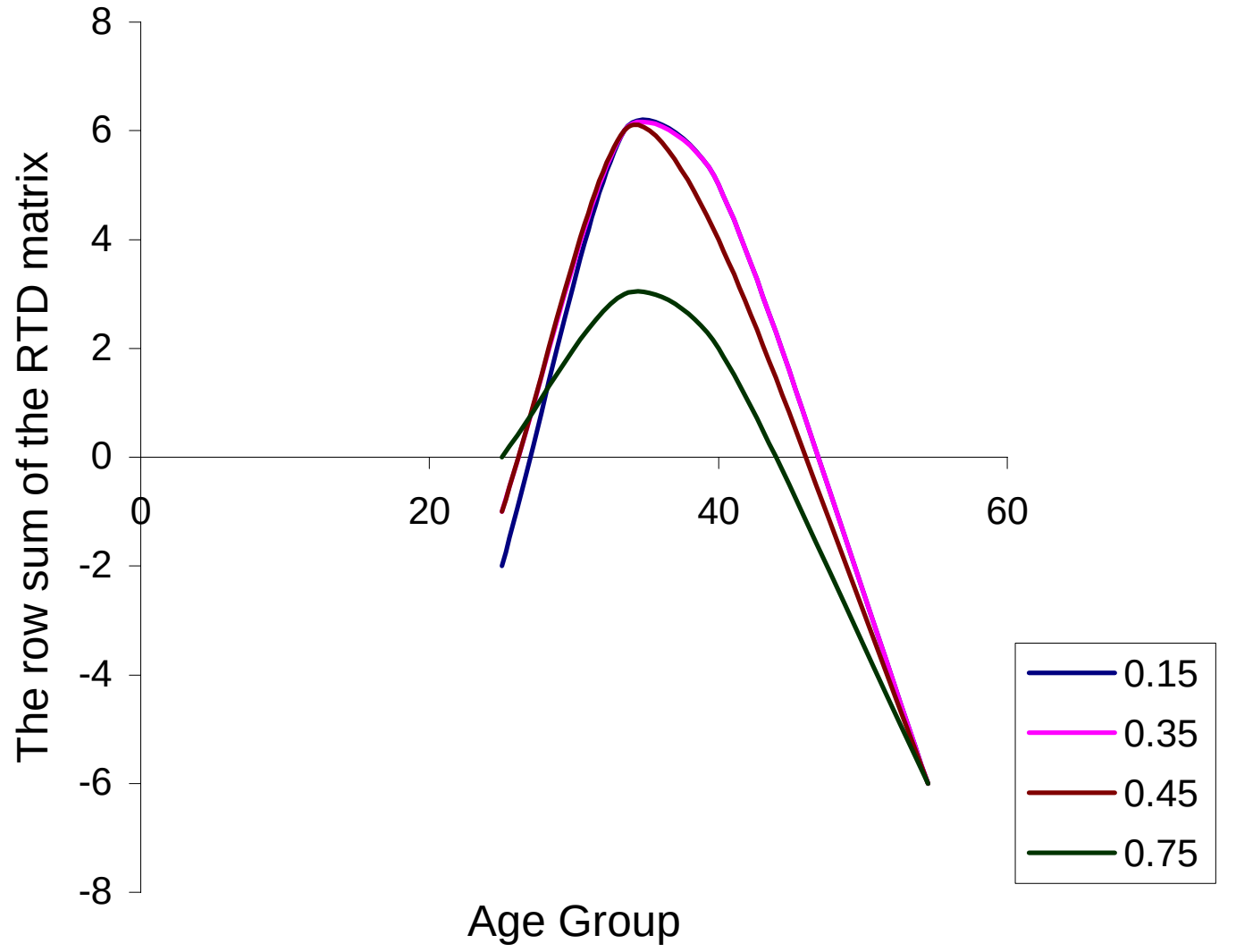

CETD matrix

$$\begin{bmatrix} -3 & -2 & 0 & 2 & 0 & -1 \\ 4 & 4 & 3 & 4 & 3 & 3 \\ 2 & 0 & 4 & 3 & 3 & 4 \\ -4 & -4 & -4 & -4 & -4 & -4 \end{bmatrix}$$

Row sum matrix

$$\begin{bmatrix} -4 \\ 21 \\ 16 \\ -24 \end{bmatrix}$$

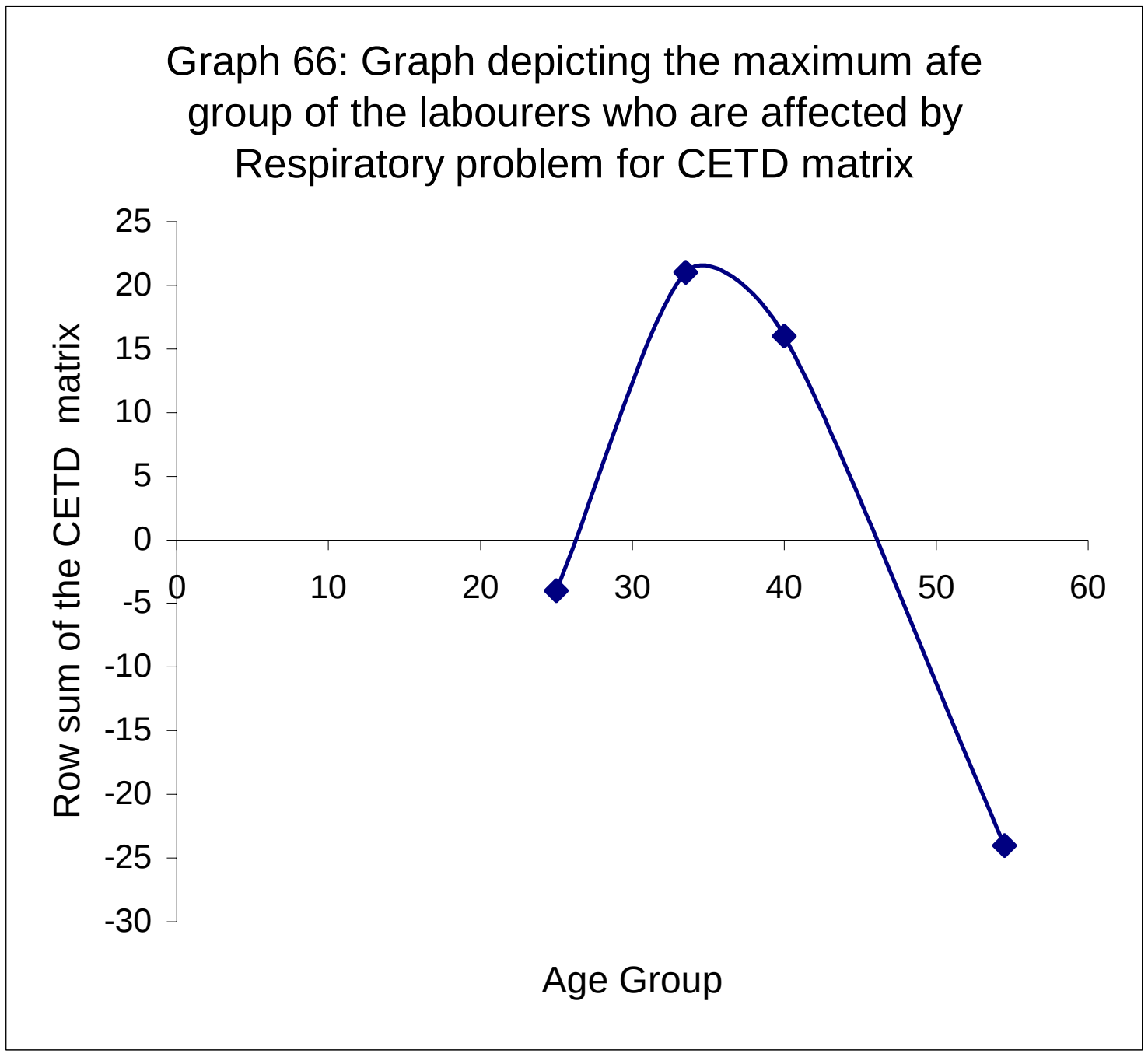

## 2.1.4.3 Estimation of maximum age group of agricultural labourers with pollution related respiratory problem using 5 × 8 matrices

Now to make the study more sensitive we increase the number of rows by 4 and see whether the decision arrived is more sensitive to the earlier one we have discussed. Thus we give the raw data of 5 × 6 matrix.

The Initial Raw Data Matrix
of Respiratory Problem of order 5 × 6

| Years | $S_1$ | $S_2$ | $S_3$ | $S_4$ | $S_5$ | $S_6$ |
|---|---|---|---|---|---|---|
| 20-24 | 1 | 3 | 4 | 6 | 4 | 5 |
| 25-30 | 4 | 12 | 14 | 20 | 20 | 16 |
| 31-36 | 7 | 16 | 13 | 17 | 15 | 16 |
| 37-43 | 6 | 10 | 17 | 18 | 18 | 19 |
| 44-65 | 3 | 14 | 17 | 22 | 23 | 23 |

ATD Matrix of
Respiratory problem of order 5 × 6

| Years | $S_1$ | $S_2$ | $S_3$ | $S_4$ | $S_5$ | $S_6$ |
|---|---|---|---|---|---|---|
| 20-24 | 0.2 | 0.6 | 0.8 | 1.2 | 0.8 | 1 |
| 25-30 | 0.67 | 2 | 2.33 | 3.33 | 3.33 | 2.67 |
| 31-36 | 1.17 | 2.67 | 2.17 | 2.83 | 2.5 | 2.67 |
| 37-43 | 0.86 | 1.43 | 2.43 | 2.57 | 2.57 | 2.71 |
| 44-65 | 0.14 | 0.64 | 0.77 | 1 | 1.05 | 1.05 |

The Average and Standard Deviation
of the given ATD matrix

| | | | | | | |
|---|---|---|---|---|---|---|
| Average | 0.61 | 1.47 | 1.7 | 2.19 | 2.05 | 2.02 |
| Standard deviation | 0.44 | 0.89 | 0.84 | 1.03 | 1.08 | 0.91 |

RTD matrix for $\alpha = 0.1$ Row sum matrix

$$\begin{bmatrix} -1 & -1 & -1 & -1 & -1 & -1 \\ 1 & 1 & 1 & 1 & 1 & 1 \\ 1 & 1 & 1 & 1 & 1 & 1 \\ 1 & 0 & 1 & 1 & 1 & 1 \\ -1 & -1 & -1 & -1 & -1 & -1 \end{bmatrix} \qquad \begin{bmatrix} -6 \\ 6 \\ 6 \\ 5 \\ -6 \end{bmatrix}$$

RTD matrix for $\alpha = 0.15$ Row sum matrix

$$\begin{bmatrix} -1 & -1 & -1 & -1 & -1 & -1 \\ 0 & 1 & 1 & 1 & 1 & 1 \\ 1 & 1 & 1 & 1 & 1 & 1 \\ 1 & 0 & 1 & 1 & 1 & 1 \\ -1 & -1 & -1 & -1 & -1 & -1 \end{bmatrix} \qquad \begin{bmatrix} -6 \\ 5 \\ 6 \\ 5 \\ -6 \end{bmatrix}$$

RTD matrix for $\alpha = 0.2$ Row sum matrix

$$\begin{bmatrix} -1 & -1 & -1 & -1 & -1 & -1 \\ 0 & 1 & 1 & 1 & 1 & 1 \\ 1 & 1 & 1 & 1 & 1 & 1 \\ 1 & 0 & 1 & 1 & 1 & 1 \\ -1 & -1 & -1 & -1 & -1 & -1 \end{bmatrix} \qquad \begin{bmatrix} -6 \\ 5 \\ 6 \\ 5 \\ -6 \end{bmatrix}$$

RTD matrix for $\alpha = 0.35$ Row sum matrix

$$\begin{bmatrix} -1 & -1 & -1 & -1 & -1 & -1 \\ 0 & 1 & 1 & 1 & 1 & 1 \\ 1 & 1 & 1 & 1 & 1 & 1 \\ 1 & 0 & 1 & 1 & 1 & 1 \\ -1 & -1 & -1 & -1 & -1 & -1 \end{bmatrix} \qquad \begin{bmatrix} -6 \\ 5 \\ 6 \\ 5 \\ -6 \end{bmatrix}$$

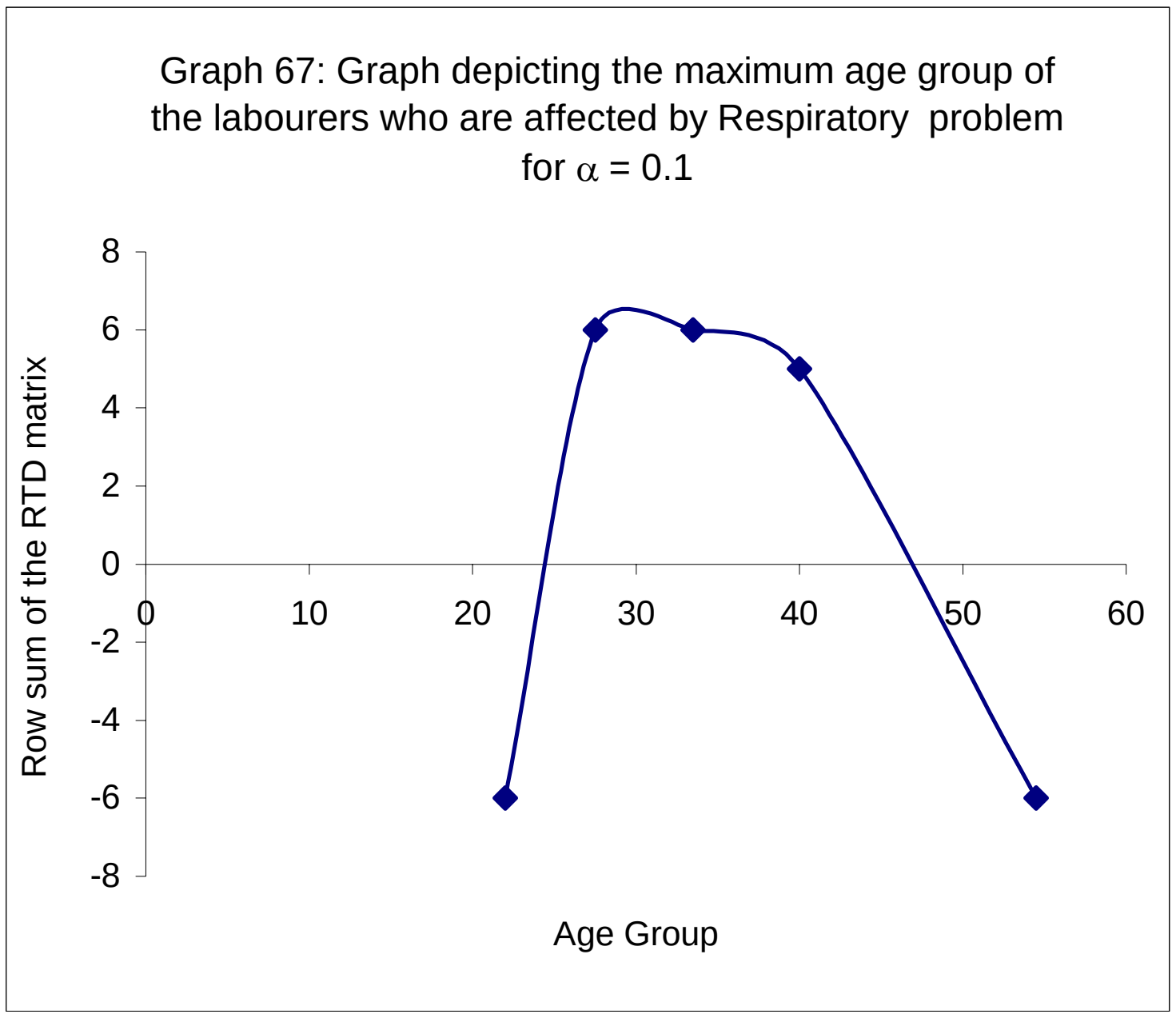
Graph 67: Graph depicting the maximum age group of the labourers who are affected by Respiratory problem for α = 0.1
Row sum of the RTD matrix
Age Group
8
6
4
2
0
-2
-4
-6
-8
0
10
20
30
40
50
60

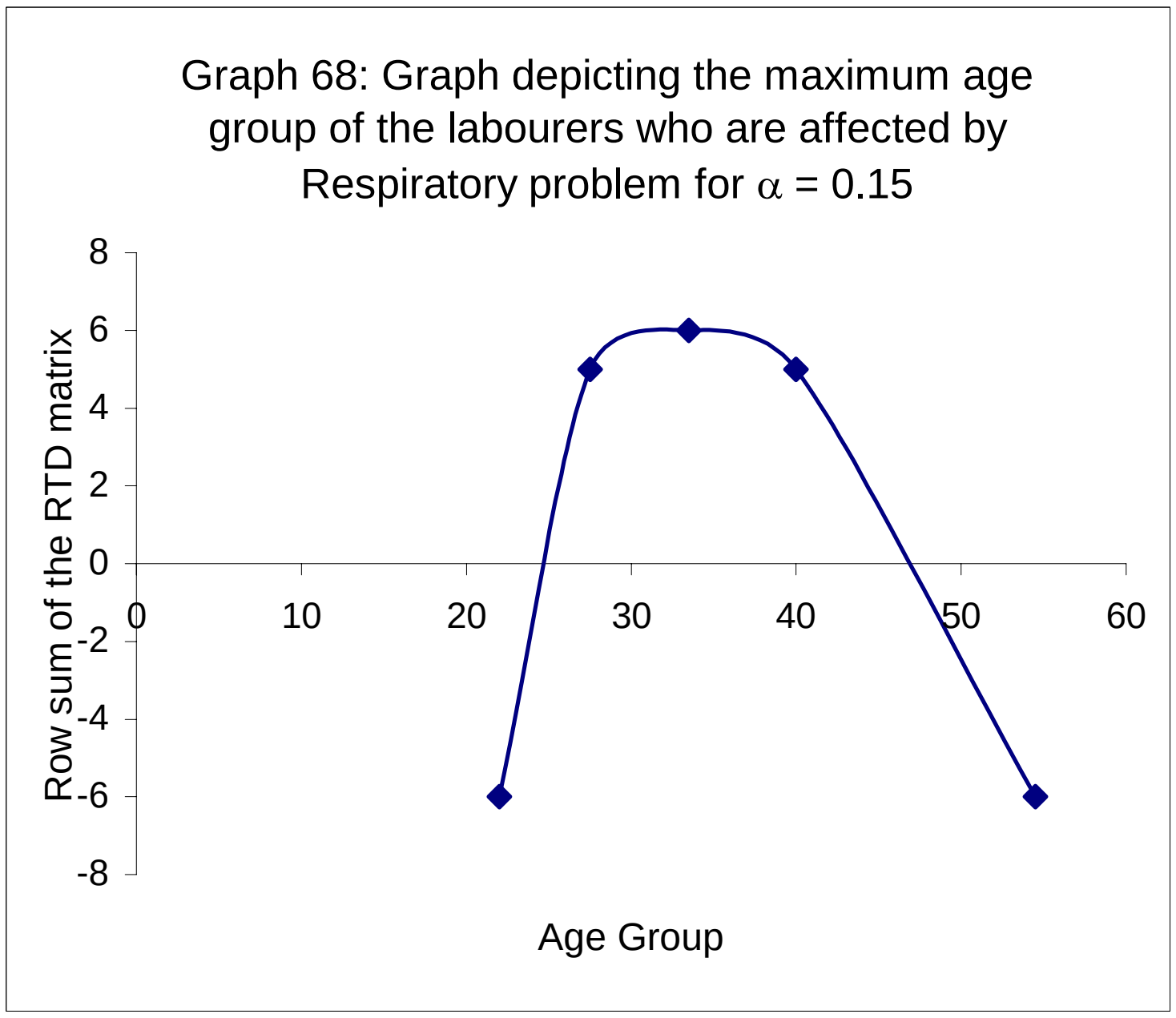
Graph 68: Graph depicting the maximum age group of the labourers who are affected by Respiratory problem for α = 0.15
Row sum of the RTD matrix
Age Group
8
6
4
2
0
-2
-4
-6
-8
0
10
20
30
40
50
60

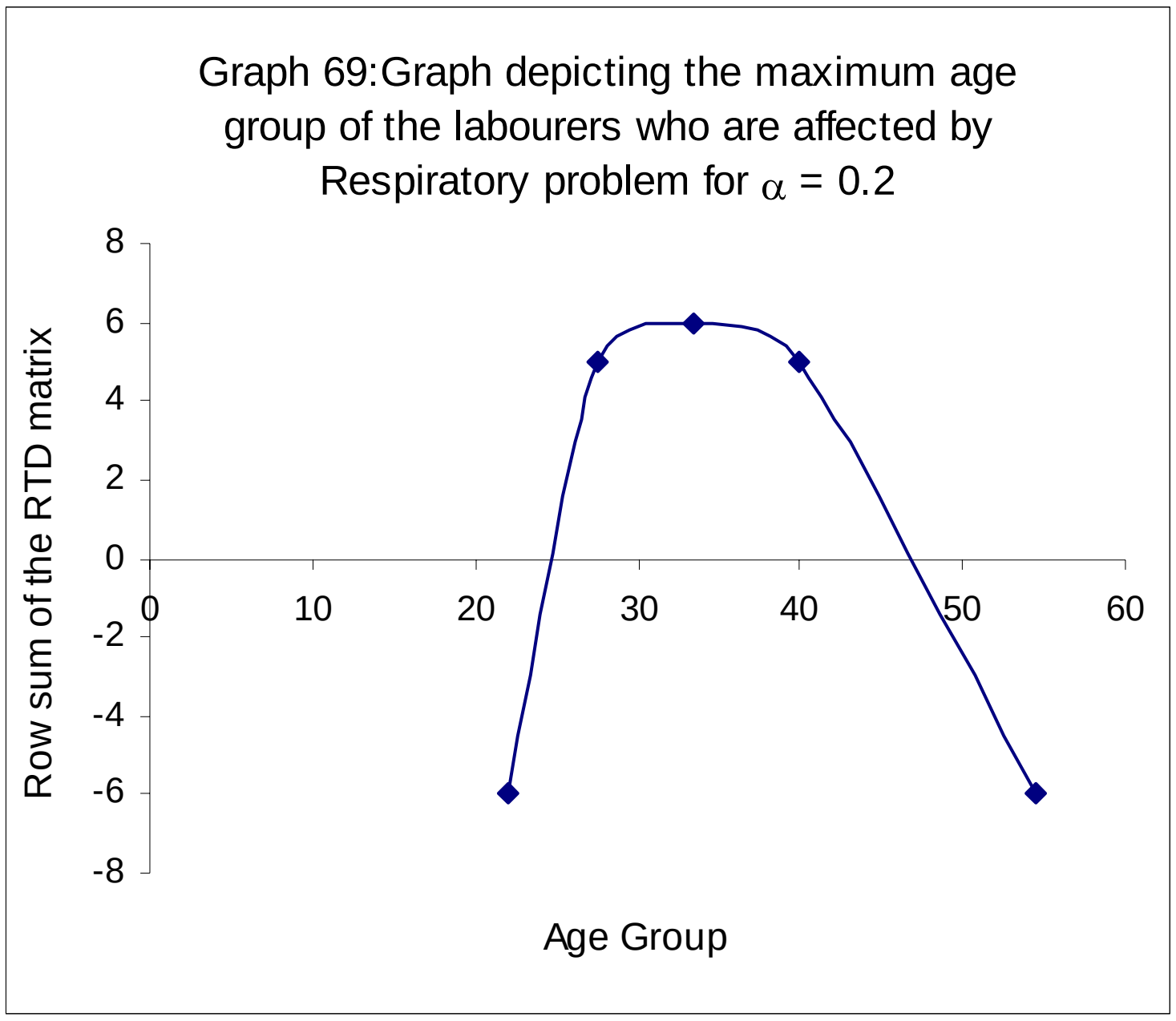
Graph 69:Graph depicting the maximum age group of the labourers who are affected by Respiratory problem for α = 0.2
Row sum of the RTD matrix
8
6
4
2
0
-2
-4
-6
-8
0
10
20
30
40
50
60
Age Group

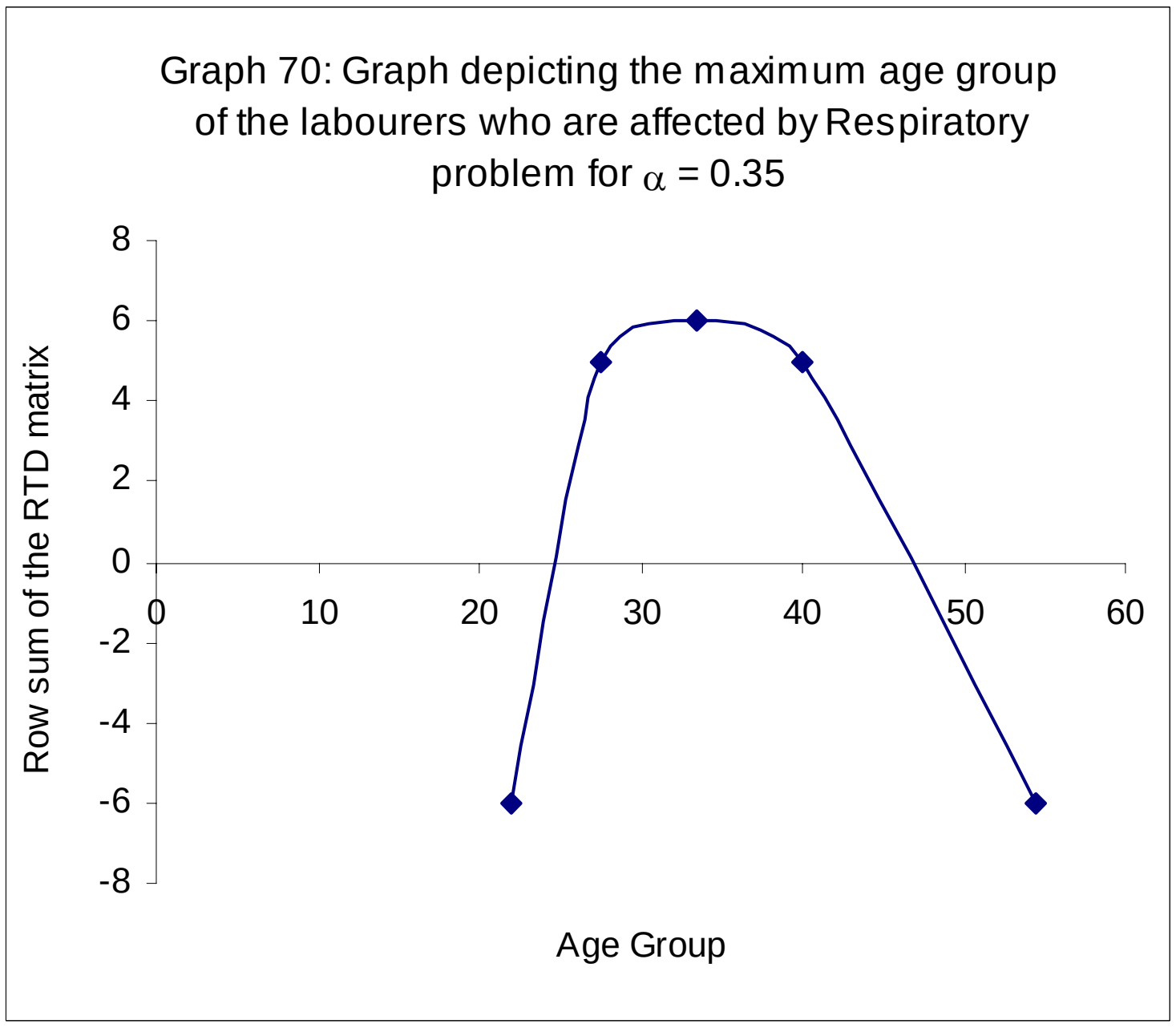
Graph 70: Graph depicting the maximum age group of the labourers who are affected by Respiratory problem for α = 0.35
Row sum of the RTD matrix
8
6
4
2
0
-2
-4
-6
-8
0
10
20
30
40
50
60
Age Group

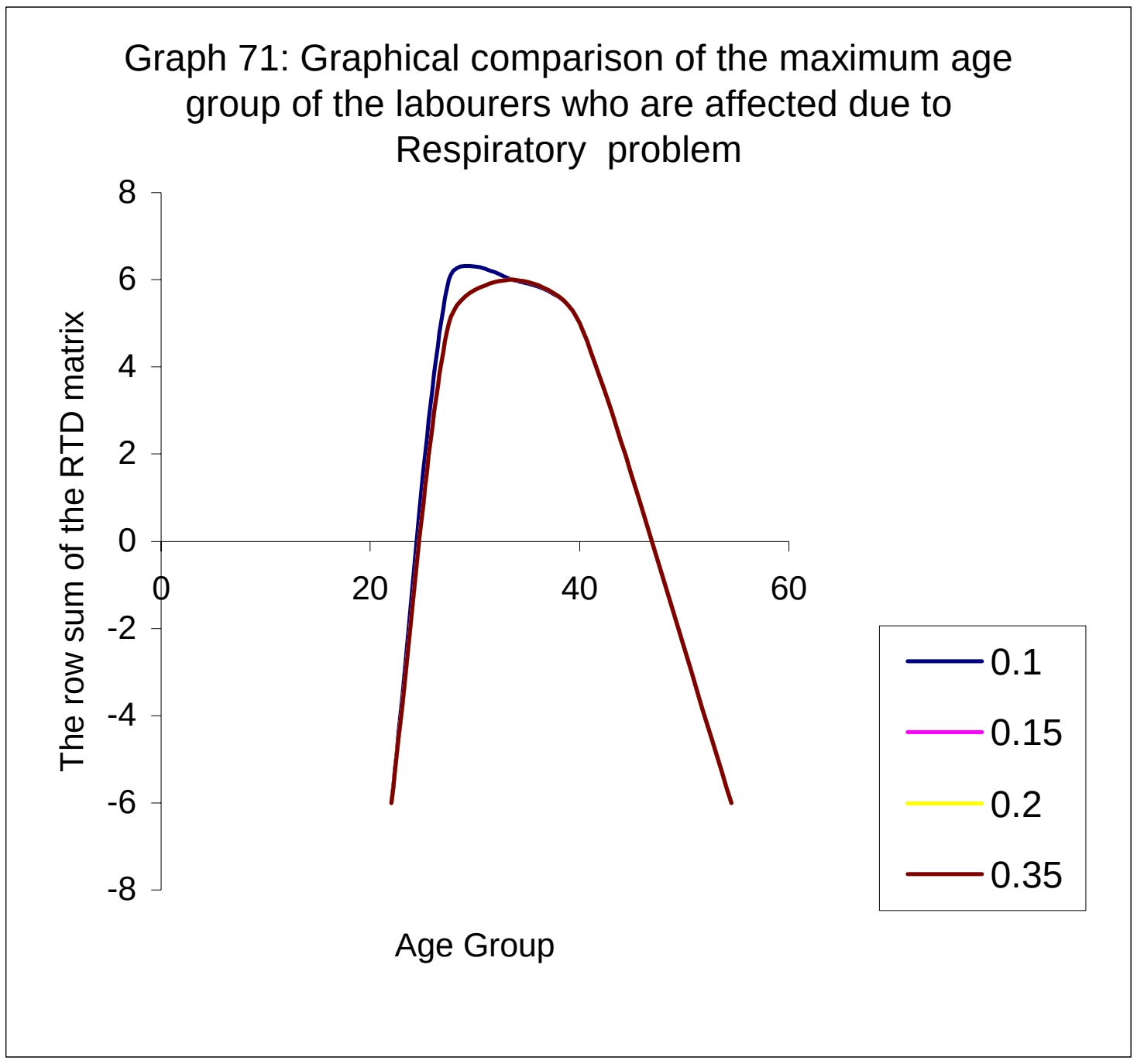

CETD matrix

$$\begin{bmatrix} -4 & -4 & -4 & -4 & -4 & -4 \\ 1 & 4 & 4 & 4 & 4 & 4 \\ 4 & 4 & 4 & 4 & 4 & 4 \\ 4 & 0 & 4 & 4 & 4 & 4 \\ -4 & -4 & -4 & -4 & -4 & -4 \end{bmatrix}$$

Row sum matrix

$$\begin{bmatrix} -24 \\ 21 \\ 24 \\ 20 \\ -24 \end{bmatrix}$$

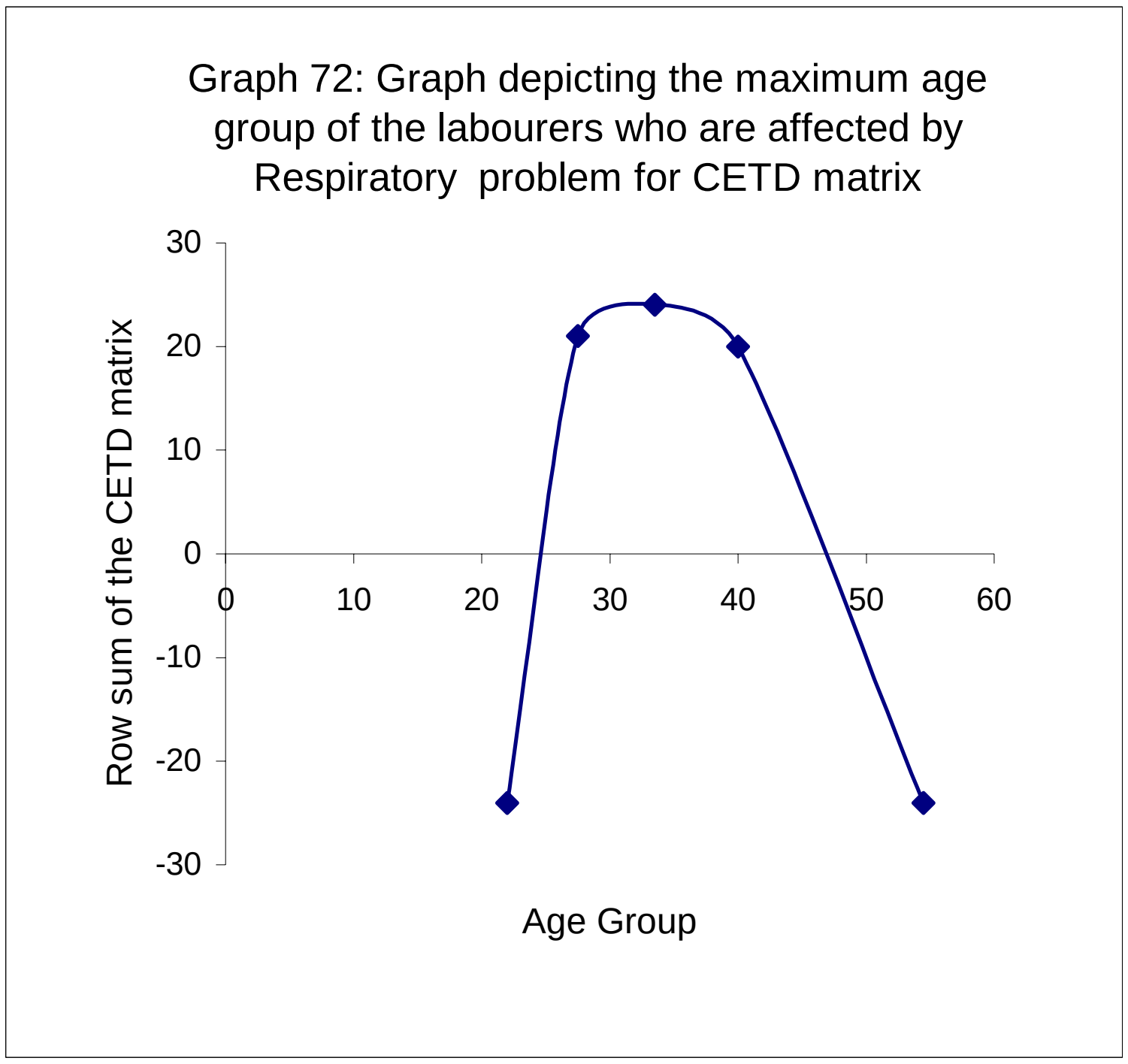

### 2.1.4.4 Conclusion

The respiratory symptoms like breathing with difficulty, allergy, asthma are the symptomatic indications of environmental pollution, which is true in urban community. But rural areas can surely boast about the pollutant free environment. In this graph, the above symptoms are at its peak at 29 yrs of age. Falling in the age range of 20-30 yrs. One can also suggest that alcoholism, drugs and smoke addiction as reasons for this. But medical history have put in the fact that such addiction has its symptomatic indication after the crossing of middle age. Out of 110 people interviewed 79 of them complained of difficulty in breathing and 80 complained of coughing and 83 of sneezing. These facts make us to point out that influence of the spraying of the pesticides using the modern technique like helicopter is one of the major reasons for the above problem in such a pollutant free environment. Thus the use of modern technique has only polluted the entire environment their by making the poor agriculture laborer as a passive victim of respiratory problem for which they are not even in a situation to get any type of medical aid for all the 11 villages in Chengalpet District in Tamil Nadu in India which we surveyed have no health centres.

## 2.2 Definition of Fuzzy Cognitive Maps with real world model representation

This section has five subsections. In subsection one we recall the definition and basic properties of Fuzzy Cognitive Maps (FCMs). In subsection two we give properties and models of FCMs and present some of its applications to problems such as the maximum utility of a route, Socio-economic problems and Symptom-disease model. In section three we use Combined Fuzzy Cognitive Maps to study the HIV/AIDS affected migrant labourers socio-economic problems. In section four we give Combined disjoint block FCM and its application to HIV/AIDS problem. In section five we use combined overlap block FCM to analyse the problem of HIV/AIDS affected migrant labourers.

### 2.2.1 Definition of Fuzzy Cognitive Maps

In this section we recall the notion of Fuzzy Cognitive Maps (FCMs), which was introduced by Bart Kosko [108] in the year 1986. We also give several of its interrelated definitions. FCMs have a major role to play mainly when the data concerned is an unsupervised one. Further this method is most simple and an effective one as it can analyse the data by directed graphs and connection matrices.

**DEFINITION 2.2.1.1:** *An FCM is a directed graph with concepts like policies, events etc. as nodes and causalities as edges. It represents causal relationship between concepts.*

***Example 2.2.1.1***: In Tamil Nadu (a southern state in India) in the last decade several new engineering colleges have been approved and started. The resultant increase in the production of engineering graduates in these years is disproportionate with the need of engineering graduates. This has resulted in thousands of unemployed and underemployed graduate engineers. Using an expert's opinion we study the effect of such unemployed people on the society. An expert spells out the five major concepts relating to the unemployed graduated engineers as

| | | |
|---|---|---|
| $E_1$ | – | Frustration |
| $E_2$ | – | Unemployment |
| $E_3$ | – | Increase of educated criminals |
| $E_4$ | – | Under employment |
| $E_5$ | – | Taking up drugs etc. |

The directed graph where $E_1, \ldots, E_5$ are taken as the nodes and causalities as edges as given by an expert is given in the following Figure 2.2.1.1:

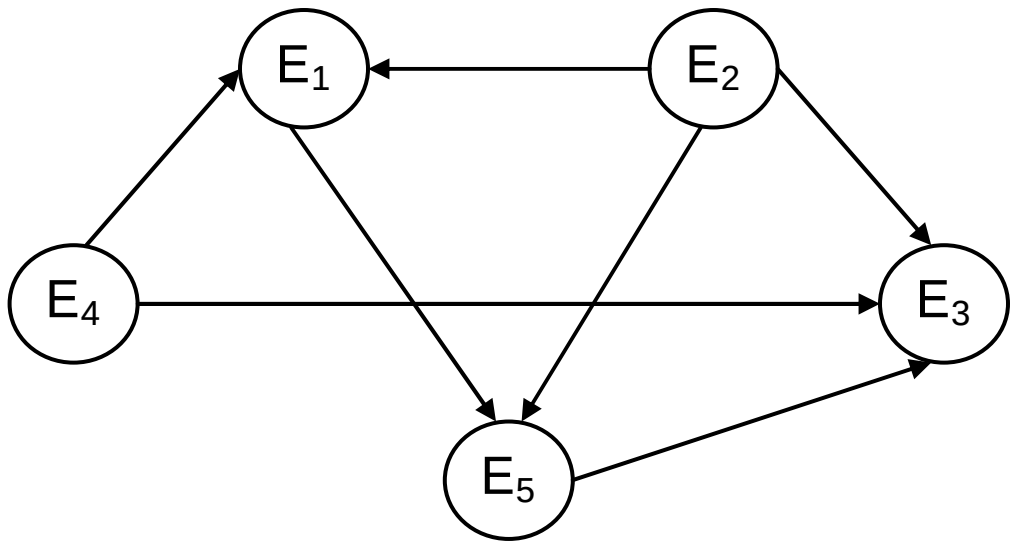

**FIGURE: 2.2.1.1**

According to this expert, increase in unemployment increases frustration. Increase in unemployment, increases the educated criminals. Frustration increases the graduates to take up to evils like drugs etc. Unemployment also leads to the increase in number of persons who take up to drugs, drinks etc. to forget their worries and unoccupied time. Under-employment forces then to do criminal acts like theft (leading to murder) for want of more money and so on. Thus one cannot actually get data for this but can use the expert's opinion for this unsupervised data to obtain some idea about the real plight of the situation. This is just an illustration to show how FCM is described by a directed graph.

{If increase (or decrease) in one concept leads to increase (or decrease) in another, then we give the value 1. If there exists no relation between two concepts the value 0 is given. If increase (or decrease) in one concept decreases (or increases) another, then we give the value –1. Thus FCMs are described in this way.}

**DEFINITION 2.2.1.2:** *When the nodes of the FCM are fuzzy sets then they are called as fuzzy nodes.*

**DEFINITION 2.2.1.3:** *FCMs with edge weights or causalities from the set {–1, 0, 1} are called simple FCMs.*

**DEFINITION 2.2.1.4:** *Consider the nodes / concepts $C_1, \ldots, C_n$ of the FCM. Suppose the directed graph is drawn using edge weight $e_{ij} \in \{0, 1, -1\}$. The matrix E be defined by $E = (e_{ij})$ where $e_{ij}$ is the weight of the directed edge $C_i C_j$. E is called the adjacency matrix of the FCM, also known as the connection matrix of the FCM.*

It is important to note that all matrices associated with an FCM are always square matrices with diagonal entries as zero.

**DEFINITION 2.2.1.5:** *Let $C_1, C_2, \ldots, C_n$ be the nodes of an FCM. $A = (a_1, a_2, \ldots, a_n)$ where $a_i \in \{0, 1\}$. A is called the instantaneous state vector and it denotes the on-off position of the node at an instant.*

$$a_i = 0 \text{ if } a_i \text{ is off and}$$
$$a_i = 1 \text{ if } a_i \text{ is on}$$

*for $i = 1, 2, \ldots, n$.*

**DEFINITION 2.2.1.6:** *Let $C_1, C_2, \ldots, C_n$ be the nodes of an FCM. Let $\overrightarrow{C_1C_2}, \overrightarrow{C_2C_3}, \overrightarrow{C_3C_4}, \ldots, \overrightarrow{C_iC_j}$ be the edges of the FCM ($i \neq j$). Then the edges form a directed cycle. An FCM is said to be cyclic if it possesses a directed cycle. An FCM is said to be acyclic if it does not possess any directed cycle.*

**DEFINITION 2.2.1.7:** *An FCM with cycles is said to have a feedback.*

**DEFINITION 2.2.1.8:** *When there is a feedback in an FCM, i.e., when the causal relations flow through a cycle in a revolutionary way, the FCM is called a dynamical system.*

**DEFINITION 2.2.1.9:** *Let* $\overline{C_1C_2}$, $\overline{C_2C_3}$, … , $\overrightarrow{C_{n-1}C_n}$ *be a cycle. When* $C_i$ *is switched on and if the causality flows through the edges of a cycle and if it again causes* $C_i$ *, we say that the dynamical system goes round and round. This is true for any node* $C_i$ *, for* $i = 1, 2, \ldots , n$*. The equilibrium state for this dynamical system is called the hidden pattern.*

**DEFINITION 2.2.1.10:** *If the equilibrium state of a dynamical system is a unique state vector, then it is called a fixed point.*

***Example 2.2.1.2:*** Consider a FCM with $C_1$, $C_2$, …, $C_n$ as nodes. For example let us start the dynamical system by switching on $C_1$. Let us assume that the FCM settles down with $C_1$ and $C_n$ on i.e. the state vector remains as (1, 0, 0, …, 0, 1) this state vector (1, 0, 0, …, 0, 1) is called the fixed point.

**DEFINITION 2.2.1.11:** *If the FCM settles down with a state vector repeating in the form*

$$A_1 \to A_2 \to \ldots \to A_i \to A_1$$

*then this equilibrium is called a limit cycle.*

Methods of finding the hidden pattern are discussed in the following Section 1.2.

**DEFINITION 2.2.1.12:** *Finite number of FCMs can be combined together to produce the joint effect of all the FCMs. Let* $E_1$, $E_2$, … , $E_p$ *be the adjacency matrices of the FCMs with nodes* $C_1$, $C_2$, …, $C_n$ *then the combined FCM is got by adding all the adjacency matrices* $E_1, E_2, \ldots, E_p$.

*We denote the combined FCM adjacency matrix by* $E = E_1 + E_2 + \ldots + E_p$.

**NOTATION:** Suppose A = ($a_1$, … , $a_n$) is a vector which is passed into a dynamical system E. Then AE = ($a'_1$, … , $a'_n$) after thresholding and updating the vector suppose we get ($b_1$, … , $b_n$) we denote that by

$$(a'_1, a'_2, \ldots , a'_n) \hookrightarrow (b_1, b_2, \ldots , b_n).$$

Thus the symbol ' $\hookrightarrow$ ' means the resultant vector has been thresholded and updated.

FCMs have several advantages as well as some disadvantages. The main advantage of this method it is simple. It functions on expert's opinion. When the data happens to be an unsupervised one the FCM comes handy. This is the only known fuzzy technique that gives the hidden pattern of the situation. As we have a very well known theory, which states that the strength of the data depends on, the number of experts' opinion we can use combined FCMs with several experts' opinions.

At the same time the disadvantage of the combined FCM is when the weightages are 1 and –1 for the same $C_i$ $C_j$, we have the sum adding to zero thus at all times the connection matrices $E_1$, … , $E_k$ may not be conformable for addition.

Combined conflicting opinions tend to cancel out and assisted by the strong law of large numbers, a consensus emerges as the sample opinion approximates the underlying population opinion. This problem will be easily overcome if the FCM entries are only 0 and 1.

We have just briefly recalled the definitions. For more about FCMs please refer Kosko [108-112].

### 2.2.2 Fuzzy Cognitive Maps – Properties and Models

Fuzzy Cognitive Maps (FCMs) are more applicable when the data in the first place is an unsupervised one. The FCMs work on the opinion of experts. FCMs model the world as a collection of classes and causal relations between classes.

FCMs are fuzzy signed directed graphs with feedback. The directed edge $e_{ij}$ from causal concept $C_i$ to concept $C_j$ measures how much $C_i$ causes $C_j$. The time varying concept function $C_i(t)$ measures the non negative occurrence of some fuzzy event, perhaps the strength of a political sentiment, historical trend or military objective.

FCMs are used to model several types of problems varying from gastric-appetite behavior, popular political developments etc. FCMs are also used to model in robotics like plant control.

The edges $e_{ij}$ take values in the fuzzy causal interval [–1, 1]. $e_{ij} = 0$ indicates no causality, $e_{ij} > 0$ indicates causal increase $C_j$ increases as $C_i$ increases (or $C_j$ decreases as $C_i$ decreases). $e_{ij} < 0$ indicates causal decrease or negative causality. $C_j$ decreases as $C_i$ increases (and or $C_j$ increases as $C_i$ decreases). Simple FCMs have edge values in {–1, 0, 1}. Then if causality occurs, it occurs to a maximal positive or negative degree. Simple FCMs provide a quick first approximation to an expert stand or printed causal knowledge.

We illustrate this by the following, which gives a simple FCM of a Socio-economic model. A Socio-economic model is constructed with Population, Crime, Economic condition, Poverty and Unemployment as nodes or concept. Here the simple trivalent directed graph is given by the following Figure 2.2.2.1 which is the experts opinion.

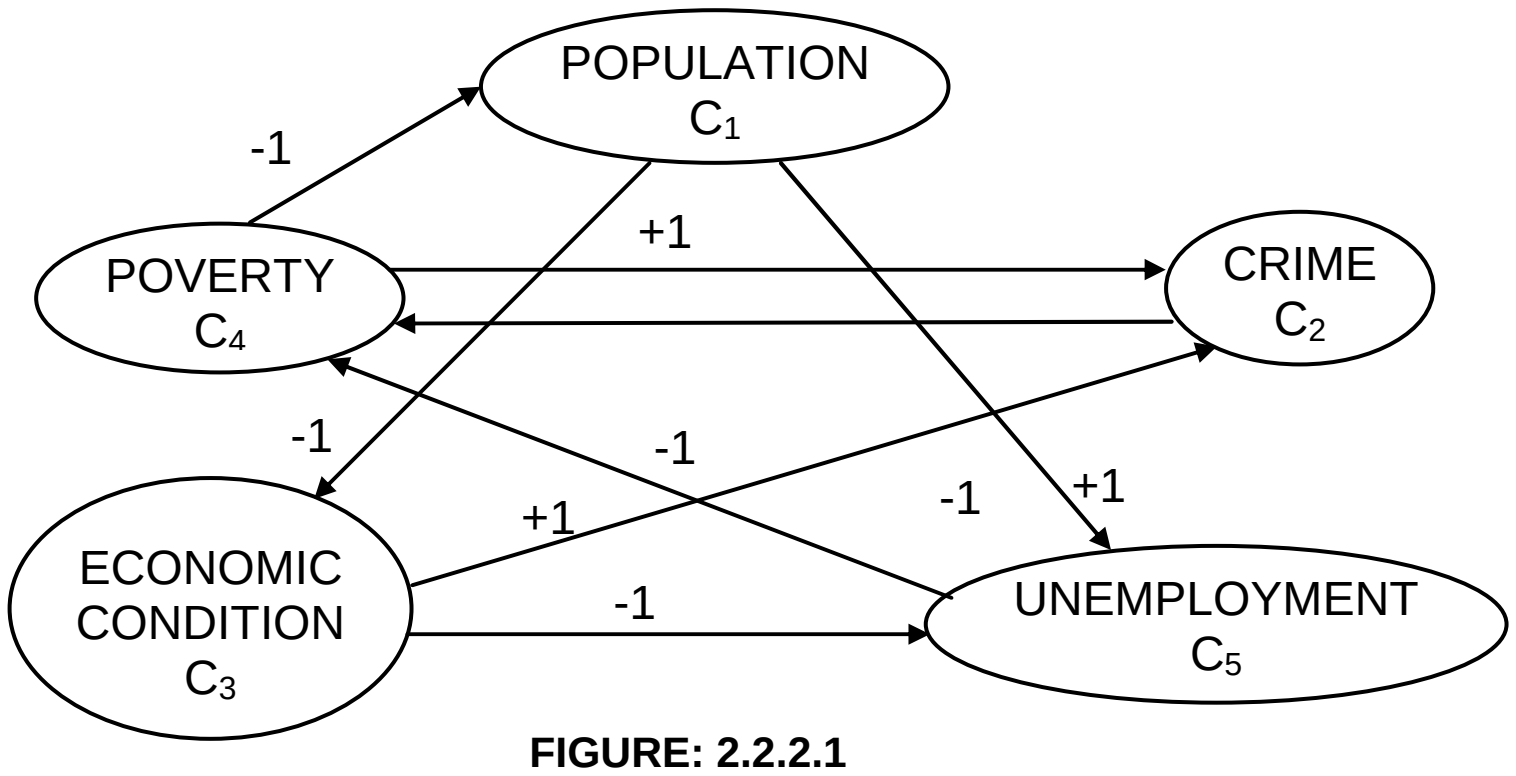

**FIGURE: 2.2.2.1**

Causal feedback loops abound in FCMs in thick tangles. Feedback precludes the graph-search techniques used in artificial-intelligence expert systems.

FCMs feedback allows experts to freely draw causal pictures of their problems and allows causal adaptation laws,

infer causal links from simple data. FCM feedback forces us to abandon graph search, forward and especially backward chaining. Instead we view the FCM as a dynamical system and take its equilibrium behavior as a forward-evolved inference. Synchronous FCMs behave as Temporal Associative Memories (TAM). We can always, in case of a model, add two or more FCMs to produce a new FCM. The strong law of large numbers ensures in some sense that knowledge reliability increases with expert sample size.

We reason with FCMs. We pass state vectors C repeatedly through the FCM connection matrix E, thresholding or non-linearly transforming the result after each pass. Independent of the FCMs size, it quickly settles down to a temporal associative memory limit cycle or fixed point which is the hidden pattern of the system for that state vector C. The limit cycle or fixed-point inference summarizes the joint effects of all the interacting fuzzy knowledge.

***Example 2.2.2.1:*** Consider the 5 × 5 causal connection matrix E that represents the socio economic model using FCM given in figure in Figure 2.2.2.1.

$$E = \begin{bmatrix} 0 & 0 & -1 & 0 & 1 \\ 0 & 0 & 0 & -1 & 0 \\ 0 & -1 & 0 & 0 & -1 \\ -1 & 1 & 0 & 0 & 0 \\ 0 & 0 & 0 & 1 & 0 \end{bmatrix}$$

Concept nodes can represent processes, events, values or policies. Consider the first node $C_1 = 1$. We hold or clamp $C_1$ on the temporal associative memories recall process. Threshold signal functions synchronously update each concept after each pass, through the connection matrix E. We start with the population $C_1 = (1\ 0\ 0\ 0\ 0)$. The arrow indicates the threshold operation,

$$C_1 E = (0\ 0\ {-1}\ 0\ 1) \hookrightarrow (1\ 0\ 0\ 0\ 1) = C_2$$
$$C_2 E = (0\ 0\ {-1}\ 1\ 1) \hookrightarrow (1\ 0\ 0\ 1\ 1) = C_3$$

$$C_3 E = (-1\ 1\ -1\ 1\ 1) \hookrightarrow (1\ 1\ 0\ 1\ 1) = C_4$$
$$C_4 E = (-1\ 1\ -1\ 0\ 1) \hookrightarrow (1\ 1\ 0\ 0\ 1) = C_5$$
$$C_5 E = (0\ 0\ -1\ 0\ 1) \hookrightarrow (1\ 0\ 0\ 0\ 1) = C_6 = C_2.$$

So the increase in population results in the unemployment problem, which is a limit cycle. For more about FCM refer Kosko [108] and for more about this socio economic model refer [235,251].

This example illustrates the strengths and weaknesses of FCM analysis. FCM allows experts to represent factual and evaluative concepts in an interactive framework. Experts can quickly draw FCM pictures or respond to questionnaires. Experts can consent or dissent to the local causal structure and perhaps the global equilibrium. The FCM knowledge representation and inferenceing structure reduces to simple vector-matrix operations, favors integrated circuit implementation and allows extension to neural statistical or dynamical systems techniques. Yet an FCM equally encodes the experts' knowledge or ignorance, wisdom or prejudice. Worse, different experts differ in how they assign causal strengths to edges and in which concepts they deem causally relevant. The FCM seems merely to encode its designers' biases and may not even encode them accurately.

FCM combination provides a partial solution to this problem. We can additively superimpose each experts FCM in associative memory fashion, even though the FCM connection matrices $E_1, \ldots, E_K$ may not be conformable for addition. Combined conflicting opinions tend to cancel out and assisted the strong law of large numbers a consensus emerges as the sample opinion approximates the underlying population opinion. FCM combination allows knowledge researchers to construct FCMs with iterative interviews or questionnaire mailings.

The laws of large numbers require that the random samples be independent identically distributed random variables with finite variance. Independence models each experts individually. Identical distribution models a particular domain focus.

We combine arbitrary FCM connection matrices $F_1, F_2, \ldots, F_K$ by adding augmented FCM matrices. $F_1, \ldots, F_K$. Each augmented matrix $F_i$ has n-rows and n-columns n equals the total number of distinct concepts used by the experts. We permute the rows and columns of the augmented matrices to bring them into mutual coincidence. Then we add the $F_i$ point wise to yield the combined FCM matrix F.

$$F = \sum_i F_i$$

We can then use F to construct the combined FCM directed graph.

Even if each expert gives trivalent description in {–1, 0, 1}, the combined (and normalized) FCM entry $f_{ij}$ tends to be in {–1, 1}. The strong law of large numbers ensures that $f_{ij}$ provides a rational approximation to the underlying unknown population opinion of how much $C_i$ affects $C_j$. We can normalize $f_{ij}$ by the number K of experts. Experts tend to give trivalent evaluations more readily and more accurately than they give weighted evaluations. When transcribing interviews or documents, a knowledge engineer can more reliably determine an edge's sign than its magnitude.

Some experts may be more credible than others. We can weight each expert with non-negative credibility weight, weighing the augmented FCM matrix.

$$F = \Sigma\, w_i\, F_i.$$

The weights need not be in [0, 1]; the only condition is they should be non-negative. Different weights may produce different equilibrium limit cycles or fixed points as hidden patterns. We can also weigh separately any submatrix of each experts augmented FCM matrix.

Augmented FCM matrices imply that every expert causally discusses every concept $C_1, \ldots, C_n$. If an expert does not include $C_j$ in his FCM model the expert implicitly say that $C_j$ is not causally relevant. So the $j^{th}$ row and the $j^{th}$ column of his augmented connection matrix contains only zeros.

***Example 2.2.2.2:*** In any nation the study of political situation i.e., the prediction of electoral winner or how people tend to prefer a particular politician and so on and so forth involves not only a lot of uncertainly for this, no data is available. They form an unsupervised data. Hence we study this model using FCM we see that while applying FCM to Indian politics the expert takes the following six nodes.

$x_1$ - Language
$x_2$ - Community
$x_3$ - Service to people, public figure configuration and Personality and nature
$x_4$ - Finance and Media
$x_5$ - Party's strength and opponent's strength
$x_6$ - Working members (Active volunteers) for the party.

We using these six nodes and obtain several experts opinion. Now let us consider a case where in three or four parties join up to make a coalition party and stand in the election. Now each party has its own opinion and when they join up as an united front we will see the result. So, let us take the opinion of 4 experts on any arbitrary four concepts of the given six. We form the directed graph of the four experts opinion, we obtain the corresponding relational matrices. Finally we will get the combined effect. The Figures 2.2.2.2 to 2.2.2.5 correspond to directed graphs of the expert opinion.

First expert's opinion in the form of the directed graph and its relational matrix is given below:

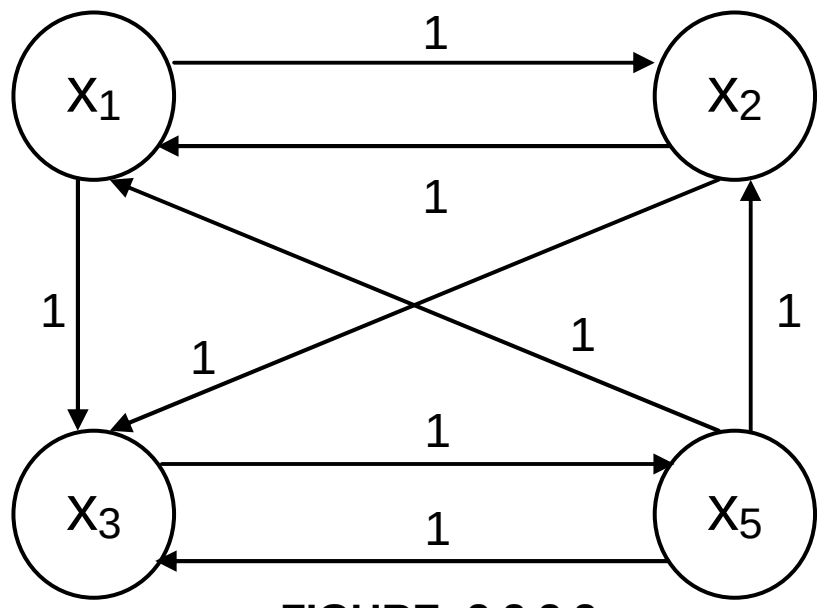

**FIGURE: 2.2.2.2**

$$E_1 = \begin{bmatrix} 0 & 1 & 1 & 0 & 0 & 0 \\ 1 & 0 & 1 & 0 & 0 & 0 \\ 0 & 0 & 0 & 0 & 1 & 0 \\ 0 & 0 & 0 & 0 & 0 & 0 \\ 1 & 1 & 1 & 0 & 0 & 0 \\ 0 & 0 & 0 & 0 & 0 & 0 \end{bmatrix}.$$

Second expert's opinion – directed graph is shown in Figure 2.2.2.3 and the relational matrix is given by $E_2$.

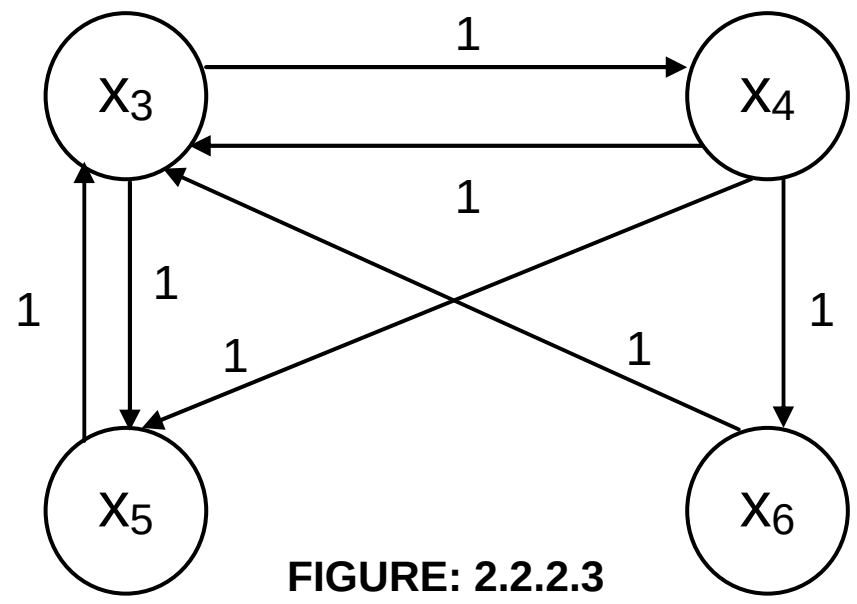

**FIGURE: 2.2.2.3**

$$E_2 = \begin{bmatrix} 0 & 0 & 0 & 0 & 0 & 0 \\ 0 & 0 & 0 & 0 & 0 & 0 \\ 0 & 0 & 0 & 1 & 1 & 0 \\ 0 & 0 & 1 & 0 & 1 & 1 \\ 0 & 0 & 1 & 0 & 0 & 0 \\ 0 & 0 & 1 & 0 & 0 & 0 \end{bmatrix}.$$

Third expert's opinion with the four conceptual nodes $x_2$, $x_3$, $x_4$ and $x_6$, its directed graph and the relational matrix are given in Figure 2.2.2.4 and its related matrix $E_3$.

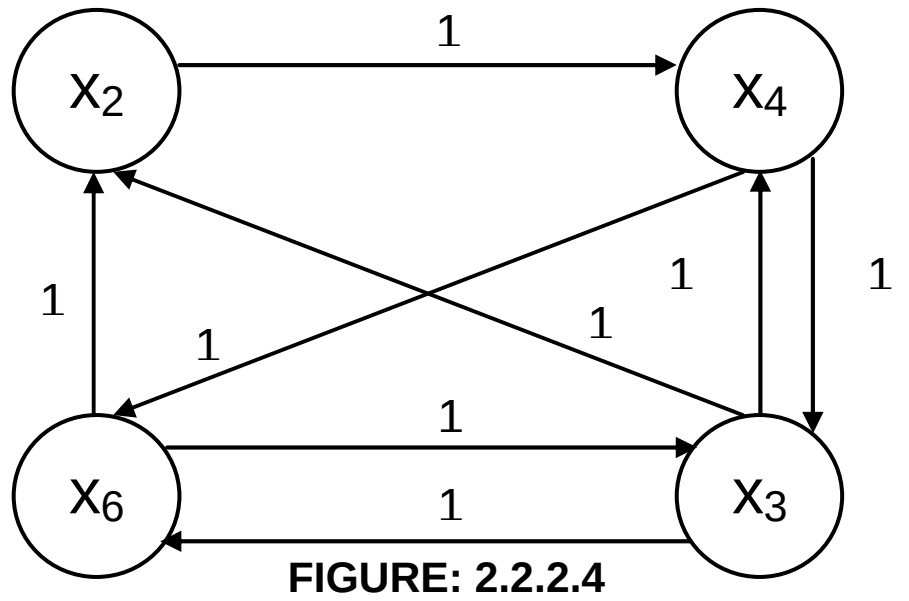

**FIGURE: 2.2.2.4**

$$E_3 = \begin{bmatrix} 0 & 0 & 0 & 0 & 0 & 0 \\ 0 & 0 & 0 & 1 & 0 & 0 \\ 0 & 1 & 0 & 1 & 0 & 1 \\ 0 & 0 & 1 & 0 & 0 & 1 \\ 0 & 0 & 0 & 0 & 0 & 0 \\ 0 & 1 & 1 & 0 & 0 & 0 \end{bmatrix}.$$

Directed graph and the relational matrix of a fourth expert using the concepts $x_2$, $x_5$, $x_4$ and $x_6$ is given in Figure 2.2.2.5 and the related matrix $E_4$ is given below:

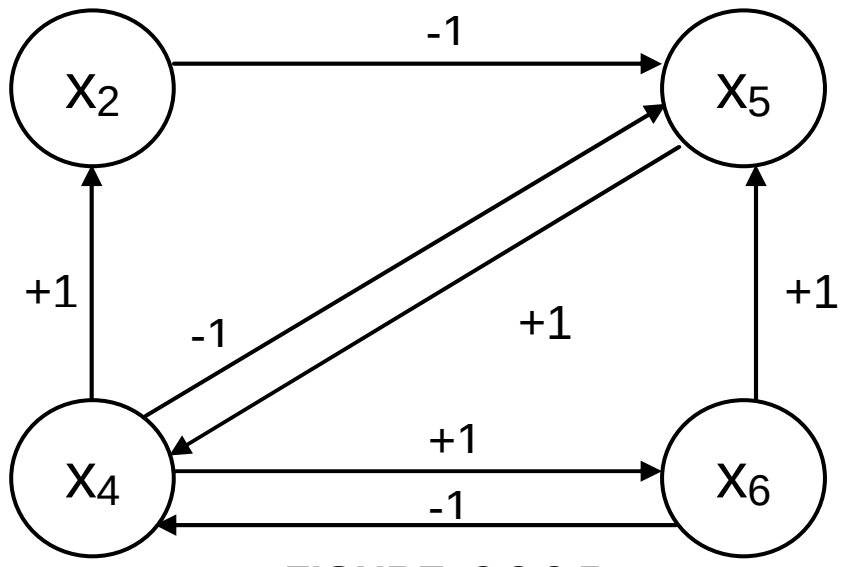

**FIGURE: 2.2.2.5**

$$E_4 = \begin{bmatrix} 0 & 0 & 0 & 0 & 0 & 0 \\ 0 & 0 & 0 & 0 & -1 & 0 \\ 0 & 0 & 0 & 0 & 0 & 0 \\ 0 & 1 & 0 & 0 & -1 & 1 \\ 0 & 0 & 0 & 1 & 0 & 0 \\ 0 & 0 & 0 & -1 & 1 & 0 \end{bmatrix}.$$

We note that each matrix contains two zero rows and two zero columns corresponding to the experts causally irrelevant concepts. We now combine the directed graph of the four experts and obtain the Figure 2.2.2.6.

The combined FCM matrix E, which is equal to the addition of the four matrices and its related directed graph, is as follows:

$$E = \begin{bmatrix} 0 & 1 & 1 & 0 & 0 & 0 \\ 1 & 0 & 1 & 1 & -1 & 0 \\ 0 & 1 & 0 & 2 & 2 & 1 \\ 0 & 1 & 2 & 0 & 0 & 3 \\ 1 & 1 & 2 & 1 & 0 & 0 \\ 0 & 1 & 2 & -1 & 1 & 0 \end{bmatrix}.$$

Combined directed graph of the combined FCM is given in Figure 2.2.2.6:

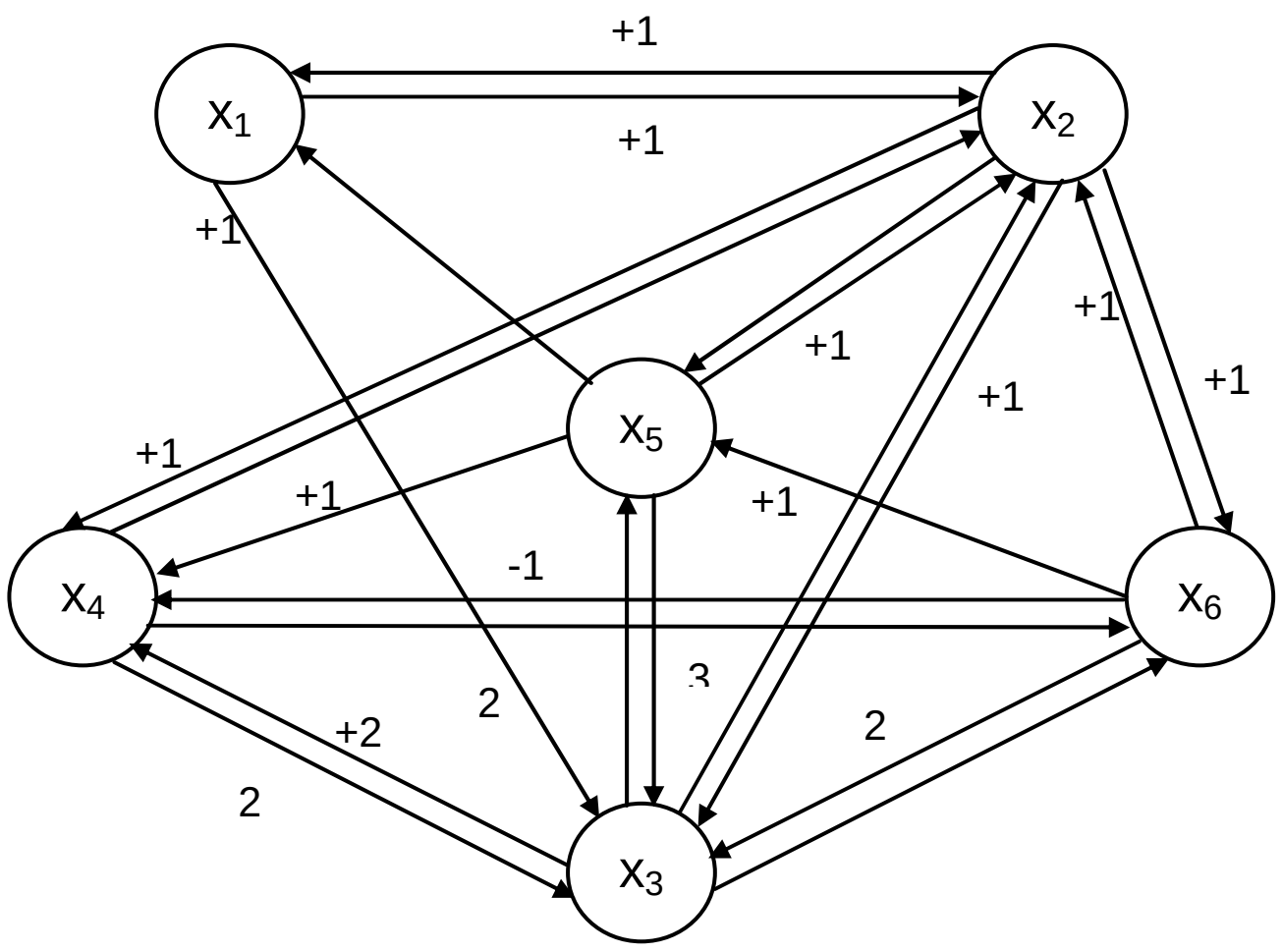

**FIGURE: 2.2.2.6**

We now consider a single state element through the matrix and multiply it with the above relational matrix and we see that nearly every state is active.

For example, let us consider the input with parties' strength and opponent's strength in the on state and rest of the coordinates as off state, vector matrix A = (0, 0, 0, 0, 1, 0). Now by passing on we get the matrix AE = (1, 1, 2, 1, 0, 0) $\hookrightarrow$ (1, 1, 2, 1, 1, 0).

Thus we see parties' strength and the opponent's strength is a very sensitive issue and it makes all other coordinate to be on except the state of the working members.

Likewise when the two concepts the community service to people and the public figure configuration is taken to be in the on state we see that on passing through the system all state becomes on.

Thus the FCMs gives us the hidden pattern. No other method can give these resultants that too with the unsupervised data. For more about this please refer [112].

***Example 2.2.2.3:*** The study of symptoms and its associations with disease in children happens to be very uncertain and difficult one.

For at times the doctor treats for the symptoms instead of treating for disease. So that when a ready-made model is made it may serve the better purpose for the doctor.

We have also adopted the FCM in case of Symptom-disease model in children [228]. To build the symptom-disease model for children we use the following 8 nodes of FCM, which are labeled as follows:

| | | |
|---|---|---|
| $C_1$ | - | Fever with cold |
| $C_2$ | - | Fever with vomiting (or) Loose Motions |
| $C_3$ | - | Fever with loss of appetite |
| $C_4$ | - | Fever with cough |
| $C_5$ | - | Respiratory diseases |
| $C_6$ | - | Gastroenteritis |
| $C_7$ | - | Jaundice |
| $C_8$ | - | Tuberculosis. |

The directed graph as given by the doctor who is taken as an expert is given in the Figure: 2.2.2.7, which is as follows:

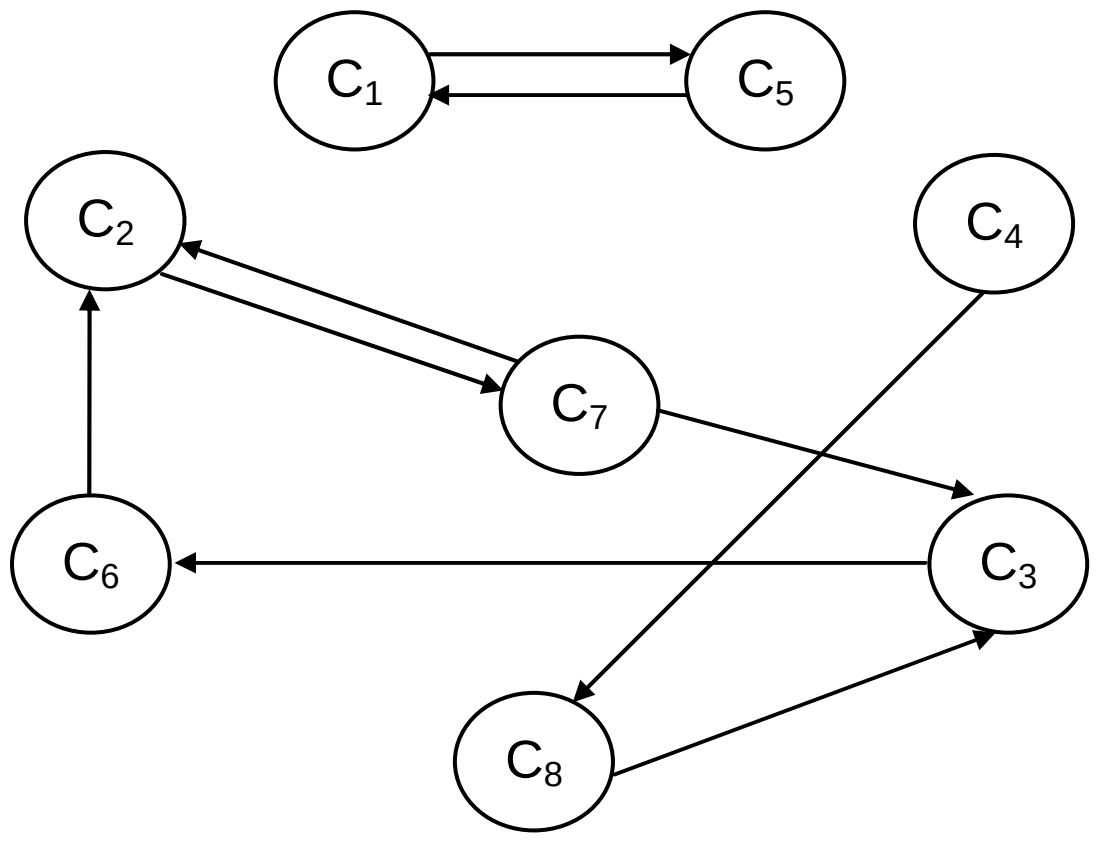

**FIGURE: 2.2.2.7**

The corresponding connection or adjacency matrix E is as follows:

$$E = \begin{bmatrix} 0 & 0 & 0 & 0 & 1 & 0 & 0 & 0 \\ 0 & 0 & 0 & 0 & 0 & 0 & 1 & 0 \\ 0 & 0 & 0 & 0 & 0 & 1 & 0 & 0 \\ 0 & 0 & 0 & 0 & 0 & 0 & 0 & 1 \\ 1 & 0 & 0 & 0 & 0 & 0 & 0 & 0 \\ 0 & 1 & 0 & 0 & 0 & 0 & 0 & 0 \\ 0 & 1 & 1 & 0 & 0 & 0 & 0 & 0 \\ 0 & 0 & 1 & 0 & 0 & 0 & 0 & 0 \end{bmatrix}.$$

Input the vector $A_1$ = (1 0 0 0 0 0 0 0)

$$A_1E \hookrightarrow (1\ 0\ 0\ 0\ 1\ 0\ 0\ 0) = A_2$$
$$A_2E \hookrightarrow (1\ 0\ 0\ 0\ 1\ 0\ 0\ 0) = A_3 = A_2.$$

According to this doctor's opinion, fever and cold induces the respiratory diseases.

Suppose the vector $D_1$ = (0 0 1 0 0 0 0 0) is taken as the input vector;

| | | | | |
|---|---|---|---|---|
| $D_1E$ | = | (0 0 0 0 0 1 0 0) | | |
| | ↪ | (0 0 1 0 0 1 0 0) | = | $D_2$ |
| $D_2E$ | ↪ | (0 1 1 0 0 1 0 0) | = | $D_3$ |
| $D_3E$ | ↪ | (0 1 1 0 0 1 1 0) | = | $D_4$ |
| $D_4E$ | = | $D_4$. | | |

Thus when a child suffers with the symptom fever and loss of appetite then the doctor suspects the child may develop fever with vomiting or loose motion leading to the sickness of gastroenteritis and jaundice. For further results in this direction please refer [228].

Use of FCMs in the study of the maximum utility of a bus route in Madras city (in South India) happens to be a difficult one for the concept deals with the many aspects of modern metropolitan public transportation. Now we just illustrate how we have applied FCM to the problem of determining the maximum utility of a route [225]. Here we not only give the FCM model to find the maximum utility of a route but also we have developed a Java program to study implications of the model [230].

***Example 2.2.2.4:*** The application of FCM in this case is done in a very ingenious way. Here for the first time FCMs are used in identifying the maximum utilization of a time period in a day, in predicting the overall utility of the routes and in the stability analysis.

We represent the different time periods of a day, viz morning, noon, evening etc by the total number of individual hours i.e. $\{H_1, H_2, \ldots, H_{24}\}$, where $H_i$ represents the $i^{th}$ hour ending of the day. In order to estimate the utility rate of a route and to identity the peak period of a route, we consider the total number of time periods in a day and the various attributes acting on these time periods as the conceptual nodes viz. $\{C_1, \ldots, C_{10}\}$. It is in the hands of the designer of the problem to assign the time periods from $\{H_1, H_2, \ldots, H_{24}\}$ to each value of $C_i$. Here we take n = 10 say $\{C_1, C_2, \ldots, C_{10}\}$ where $C_1$ corresponds to the early hours $\{H_6, H_7\}$, $C_2$ corresponds to the morning hours

{$H_8$, $H_9$, $H_{10}$}, $C_3$ corresponds to early noon hours {$H_{11}$, $H_{12}$, $H_{13}$}, $C_4$ refers to the evening {$H_{14}$, $H_{15}$, $H_{16}$}. $C_5$ corresponds to late evening hours {$H_{17}$, $H_{18}$, $H_{19}$}. $C_6$ corresponds to the night hours, {$H_{20}$, $H_{21}$, $H_{22}$}, $C_7$ indicates the number of passengers in each time period, $C_8$ corresponds to the total collection in each time period, $C_9$ denotes the number of trips made in each time period and $C_{10}$ corresponds to the hourly occupancy in each time period. We make an assumption i.e. an increase in the conceptual nodes depicting the time period say $C_i$ comprising of hours {$H_i$, $H_{(i+1)}$, …, $H_{(i+n)}$} will imply a graded rise in the hour of the day from hour $H_i$ to hour $H_{i+n}$.

To assess the interactions among the conceptual nodes, we collect expert opinion. The more number of experts the higher is the reliability in the knowledge base. The experts include the passengers traveling in the city transport service, the raw data obtained from the Pallavan Transport Corporation (currently renamed as Metropolitan Transport Corporation) (taken from the Madras city in India) denoting its route-wise loading pattern analysis.

Though we have taken several experts opinion, here we give only the expert opinion of a regular passenger traveling along the route 18B, the adjacency matrix of it is given below:

$$E = \begin{bmatrix} 0 & 0 & 0 & 0 & 0 & 0 & 1 & 1 & -1 & 1 \\ 0 & 0 & 0 & 0 & 0 & 0 & 1 & 1 & 1 & 0 \\ 0 & 0 & 0 & 0 & 0 & 0 & -1 & -1 & -1 & -1 \\ 0 & 0 & 0 & 0 & 0 & 0 & -1 & -1 & -1 & 0 \\ 0 & 0 & 0 & 0 & 0 & 0 & 1 & -1 & 1 & 0 \\ 0 & 0 & 0 & 0 & 0 & 0 & -1 & -1 & -1 & -1 \\ 0 & 0 & 0 & 0 & 0 & 0 & 0 & 1 & 1 & 0 \\ 0 & 0 & 0 & 0 & 0 & 0 & 1 & 0 & 1 & 1 \\ 0 & 0 & 0 & 0 & 0 & 0 & 1 & 0 & 0 & -1 \\ 0 & 0 & 0 & 0 & 0 & 0 & 1 & 1 & -1 & 0 \end{bmatrix}.$$

To study the hidden pattern first we clamp the concept with the first node in the on state and rest of the nodes in the off state

| | | | | |
|---|---|---|---|---|
| $A_1$ | = | (1 0 0 0 0 0 0 0 0 0) | | |
| $A_1E$ | = | (0 0 0 0 0 0 1 1 –1 1) | | |
| | $\hookrightarrow$ | (1 0 0 0 0 0 1 1 0 1) | = | $A_2$ |
| $A_2E$ | = | (0 0 0 0 0 0 2 3 0 3) | | |
| | $\hookrightarrow$ | (1 0 0 0 0 0 1 1 0 1) | = | $A_3$ |
| $A_3E$ | = | (0 0 0 0 0 0 3 2 0 2) | | |
| | $\hookrightarrow$ | (1 0 0 0 0 0 1 1 0 1) | = | $A_4 = A_3$. |

We observe that increase in time period $C_1$ increases the number of passengers, the total collection and the number of trips. Next consider the increase of the time period $C_2$. The input vector is

| | | | | |
|---|---|---|---|---|
| $B_1$ | = | (0 1 0 0 0 0 0 0 0 0) | | |
| $B_1E$ | = | (0 0 0 0 0 0 1 1 1 0) | | |
| | $\hookrightarrow$ | (0 1 0 0 0 0 1 1 1 0) | = | $B_2$ |
| $B_2E$ | = | (0 0 0 0 0 0 3 2 3 0) | | |
| | $\hookrightarrow$ | (0 1 0 0 0 0 1 1 1 0) | = | $B_2$. |

Thus we observe that an increase in the time period $C_2$ leads to an increase in number of passengers, the total collection and the increased occupancy.

Consider the combined effect of increasing $C_2$ and $C_{10}$ i.e. increase of trips in the time period $\{H_8, H_{10}, H_{10}\}$.

| | | | | |
|---|---|---|---|---|
| $W_1$ | = | (0 1 0 0 0 0 0 0 0 1) | | |
| $W_1E$ | = | (0 0 0 0 0 0 2 2 0 0) | | |
| | $\hookrightarrow$ | (0 1 0 0 0 0 1 1 0 1) | = | $W_2$ |
| $W_2E$ | = | (0 0 0 0 0 0 3 3 2 1) | | |
| | $\hookrightarrow$ | (0 1 0 0 0 0 1 1 1 1) | = | $W_3$ |
| $W_3E$ | = | (0 0 0 0 0 0 4 3 2 0) | | |
| | $\hookrightarrow$ | (0 1 0 0 0 0 1 1 1 1) | = | $W_4 = W_3$. |

Thus increase of time period $C_2$ and $C_{10}$ simultaneously increases the total number of passengers, the total collection and the occupancy rate.

## THE JAVA PSEUDO CODE TO EVALUATE THE FIXED POINTS FOR CLAMPED VECTORS

```
import java.applet.Applet;
import java.awt.*;
import java.lang.*;
import java.io.*;
import java.net.*;
import java.util.*;
public class cognitive extends java.applet.Applet
{
public void init()
   {
/* Through Text Label capture the Number of
Concepts used in this Study */
/* Define an Editable input window */
/* Define a non-Editable output window */
/* Add Button for Sample 1 to provide data for
sample matrix and the number of Concepts */
/* Add Evaluate Button: If clicked on this it
provides the values for Clamped vectors and the
corresponding Fixed point in the system */
/* Add a clear button to clear the values of
concepts, input window, output window */
/* Show all the values defined above */
}/* init() */

public boolean action(Event evt, Object arg)
{
String label = (String)arg;
if(evt.target instanceof Button)
{
if(label.equals("EVALUATE"))
{
/* Capture the number of rows in Matrix */
/* Call the function readandevaluate() to
evaluate the values */
/* show the calculated values */
}
else if (label.equals("Clear")
{
/*clear the rows, input window and output window
*/
}
else if (label.equals("Load Sample Matrix"))
{
/* set the number of values to 5 */
try
{
dataURL=new URL
("http://members.tripod.com/~mandalam/java1-
0/data3.txt");
```

```
try
{
/* get the sample data from the URL above and
populate the input window */
}
catch(IOException e)
{
System.err.println("Error:"+e);
}
}
catch(MalformedURLException e) {
return false;
}
return true;
}/* sample Matrix*/
  }
return false;
  }
void readandevaluate(String mystring)
  {
/*Tokenize the input window and get the initial
Adjacency matrix in A */
/* show the above read adjacency matrix A in
output window */
/* Determine the Hidden Pattern in B, the Clamped
Vector */
/* Determine the new vector C */
/* apply Threshold to new vector C */
/* check if C is same as B */
/* Determine the new vector C and compare it with
B */
/* if both are same this is the Fixed point in
the system, else continue by copying the values
of C into B */
  }/* readandevaluate() */
}/* end of class Cognitive */
```

Several new results in this direction can be had from [225]. Another type of problem about transportation is as follows.

The problem studied in this case is for a fixed source S, a fixed destination D and a unique route from the source to the destination, with the assumption that all the passengers travel in the same route, we identify the preferences in the regular services at the peak hour of a day.

We have considered only the peak-hour since the passenger demand is very high only during this time period, where the transport sector caters to the demands of the different groups of people like the school children, the office goers, the vendors etc.

We have taken a total of eight characteristic of the transit system, which includes the level of service and the convenience factors. We have the following elements, Frequency of the service, in-vehicle travel time, the travel fare along the route, the speed of the vehicle, the number of intermediate points, the waiting time, the number of transfers and the crowd in the bus or equivalently the congestion in the service.

Before defining the cognitive structure of the relationship, we give notations to the concepts involved in the analysis as below.

$C_1$ - Frequency of the vehicles along the route
$C_2$ - In-vehicle travel time along the route
$C_3$ - Travel fare along the route
$C_4$ - Speed of the vehicles along the route
$C_5$ - Number of intermediate points in the route
$C_6$ - Waiting time
$C_7$ - Number of transfers in the route
$C_8$ - Congestion in the vehicle.

The graphical representation of the inter-relationship between the nodes is given in the form of directed graph given in Figure: 2.2.2.8.

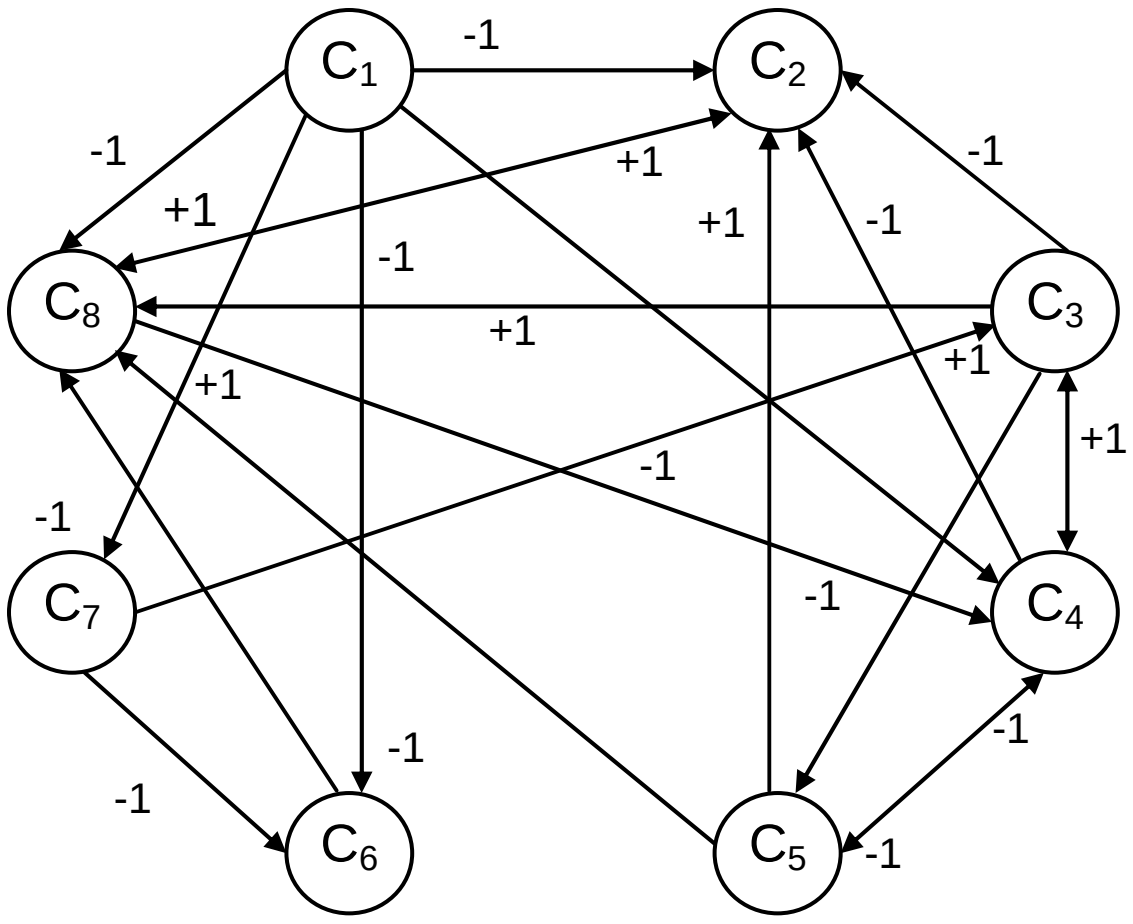

**FIGURE: 2.2.2.8**

From the above signed directed graph, we obtain a connection matrix E, since the number of concepts used here are eight, the connection matrix is a 8 × 8 matrix. Thus we have $E = [A_y]_{8\times 8}$

$$E = \begin{bmatrix} 0 & -1 & 0 & 1 & 0 & -1 & -1 & -1 \\ 0 & 0 & 0 & 0 & 0 & 0 & 0 & 1 \\ 0 & -1 & 0 & 1 & -1 & 0 & 0 & -1 \\ 0 & -1 & 1 & 0 & -1 & 0 & 0 & 0 \\ 0 & 1 & 0 & -1 & 0 & 0 & 0 & 1 \\ 0 & 0 & 0 & 0 & 0 & 0 & 0 & 1 \\ 0 & 0 & 1 & 0 & 0 & -1 & 0 & 0 \\ 0 & 1 & 0 & -1 & 0 & 0 & 0 & 0 \end{bmatrix}.$$

The model implications i.e. the predictions regarding changes in the state behaviour of a model are determined by activating the individual involved elements. This is achieved by 'Clamping' concepts of interest and including iterative operations on the matrix. We use the equation $I_{t+1} = O_{(t)} = I_{(t)} * E$ where $I_{(t)}$ is the input vector at the $t^{th}$ iteration, E is the connection matrix $O_{(t)}$ is the output vector at the $t^{th}$ iteration, used as the input for the $(t+1)^{th}$ iteration. Initially, we clamp the concept of interest. This is done by fixing the corresponding element in the input and the rest of the elements are given a value zero. Thus we have the input vector as a $1 \times n$ vector (where n is the number of concepts involved). Now we do matrix multiplication operation of $I_{1\times n}$ on the connection matrix $E_{(n\times n)}$. The output O is again a vector of order $1 \times n$. Now we implement the threshold function – a binary conversion of the output vector. Thus we have

$$O(x) = \begin{cases} 1 & \text{if } x > 1 \\ 0 & \text{otherwise} \end{cases}.$$

The output vector after the implementation of the threshold function is used as the input vector for the $(t+1)^{th}$ stage. This new input vector again is operated on the connection matrix. The working process described above can be expressed by this algorithm.

### 2.2.3 Combined FCM to study the HIV/AIDS affected migrants labourers socio-economic problem

The concept of Combined Fuzzy Cognitive Maps (CFCMs) was just defined in subsection 2.2.1. For more about CFCM please refer Kosko [108-112]. We analyze the problems of HIV/AIDS affected migrant labourers using CFCM.

***Example 2.2.3.1:*** Now we analyze the same problem using Combined Fuzzy Cognitive Maps (CFCM). Now we seek experts opinion about the seven attributes $A'_1, A'_2, \ldots, A'_7$ given below

- $A'_1$ - Easy Availability of Money
- $A'_2$ - Wrongful company, Addictive habits
- $A'_3$ - Visits CSWs
- $A'_4$ - Socially irresponsible/free
- $A'_5$ - Macho behaviour
- $A'_6$ - More leisure
- $A'_7$ - No awareness about the disease

The expert's opinion is given as graph and the related fuzzy relational matrix is given below:

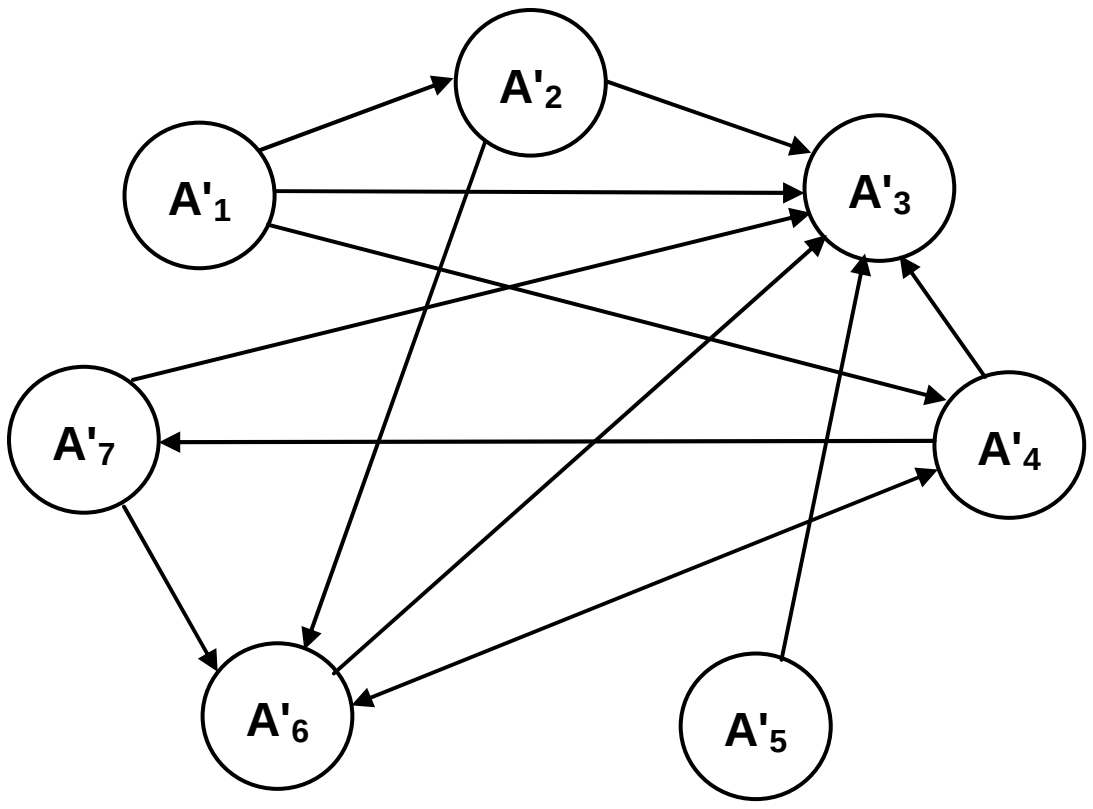

**FIGURE 2.2.3.1**

$$A = \begin{array}{c} \\ A'_1 \\ A'_2 \\ A'_3 \\ A'_4 \\ A'_5 \\ A'_6 \\ A'_7 \end{array} \begin{array}{c} \begin{array}{ccccccc} A'_1 & A'_2 & A'_3 & A'_4 & A'_5 & A'_6 & A'_7 \end{array} \\ \left[\begin{array}{ccccccc} 0 & 1 & 1 & 1 & 0 & 0 & 0 \\ 0 & 0 & 1 & 0 & 1 & 1 & 0 \\ 0 & 0 & 0 & 0 & 0 & 0 & 0 \\ 0 & 0 & 1 & 0 & 0 & 1 & 0 \\ 0 & 1 & 1 & 0 & 0 & 0 & 0 \\ 0 & 1 & 1 & 1 & 0 & 0 & 0 \\ 0 & 0 & 1 & 1 & 0 & 1 & 0 \end{array}\right] \end{array}$$

Suppose A denotes the connection matrix of the directed graph. Now to find the stability of the dynamical system or to be more precise the hidden pattern of the system which may be a fixed point or a limit cycle.

Consider the state vector / initial vector $X_1$ = (1 0 0 0 0 0 0) i.e. the only node 'easy availability of money which is dependent on the profession' i.e. $A'_1$ alone is in the on state, all other state vectors are in the off state. Now passing $X_1$ into the connection matrix A we get $X_1A \hookrightarrow$ (1 1 1 1 0 0 0) where '$\hookrightarrow$' denotes the resultant of the vector Xi after passing through A that has been updated and thresholded. Let

$$X_1A \hookrightarrow (1\ 1\ 1\ 1\ 0\ 0\ 0) = X_2$$
$$X_2A \hookrightarrow (1\ 1\ 1\ 1\ 1\ 1\ 0) = X_3$$

we see $X_3A = X_3$. Thus $X_i$, on passing through the dynamical system gives a fixed point so mathematically without any doubt we can say easy money based on profession (ie coolie, construction labourer, truck driver etc. easy money related to the individual who do not care for the family) leads to wrongful company, socially irresponsible and free, visiting of CSWs, male ego, more leisure. Only the node awareness of the disease remains in the off state, that is resulting in a fixed point.

Now let us take the node 'macho behaviour' to be in the on state i.e. Y = (0 0 0 0 1 0 0) i.e. all other states are in the off state, passing Y in the connection matrix YA $\hookrightarrow$ (0 1 1 0 1 0 0) of course after updating and thresholding the resultant vector YA, let

$$YA \hookrightarrow (0\ 1\ 1\ 0\ 1\ 0\ 0) = Y_1.$$
$$Y_1A \hookrightarrow (0\ 1\ 1\ 0\ 1\ 1\ 0) = Y_2.$$
$$Y_2A \hookrightarrow (0\ 1\ 1\ 1\ 1\ 1\ 0) = Y_3.$$

$Y_3$ which results in a fixed point. Hence we see that when a person displays a 'macho' and male-chauvinist behaviour he would opt for bad company and bad habits and would be socially irresponsible, visit the CSWs thus all the nodes are on except $A'_7$ and $A'_7$: 'Easy money' and 'No awareness about the disease.'

Now suppose we take the node $A'_7$ i.e. no awareness about the disease IS in the on state we will now find the hidden pattern with only $A'_7$ in the on state, and all other nodes to be in the off state. Let Z = (0 0 0 0 0 0 1), passing Z into the connection matrix A we get ZA $\hookrightarrow$ (0 0 1 1 0 1 1). Let $Z_1$ = (0 0 1 1 0 1 1), now

$$Z_1A \hookrightarrow (0\ 1\ 1\ 1\ 0\ 1\ 1) = Z_2 \text{ say}$$
$$Z_2A \hookrightarrow (0\ 1\ 1\ 1\ 1\ 1\ 1) = Z_3 \text{ say}$$
$$Z_3A \hookrightarrow (0\ 1\ 1\ 1\ 1\ 1\ 1) = Z_3.$$

Thus we see if the person is not aware of the disease, it implies that the person is generally a victim of wrongful company, that he is socially irresponsible, that he visits CSWs, has leisure and is full of macho behaviour and male chauvinism. However, only the concept of 'easy money and profession' has no relation with this. Thus we can in all cases obtain an hidden pattern which is never possible by any other model. Thus this model is well-suited to give the impact of each attribute in an HIV/AIDS patient and the effect or inter-relation between these attributes which no other mathematical model has given.

The directed graph given by the second expert

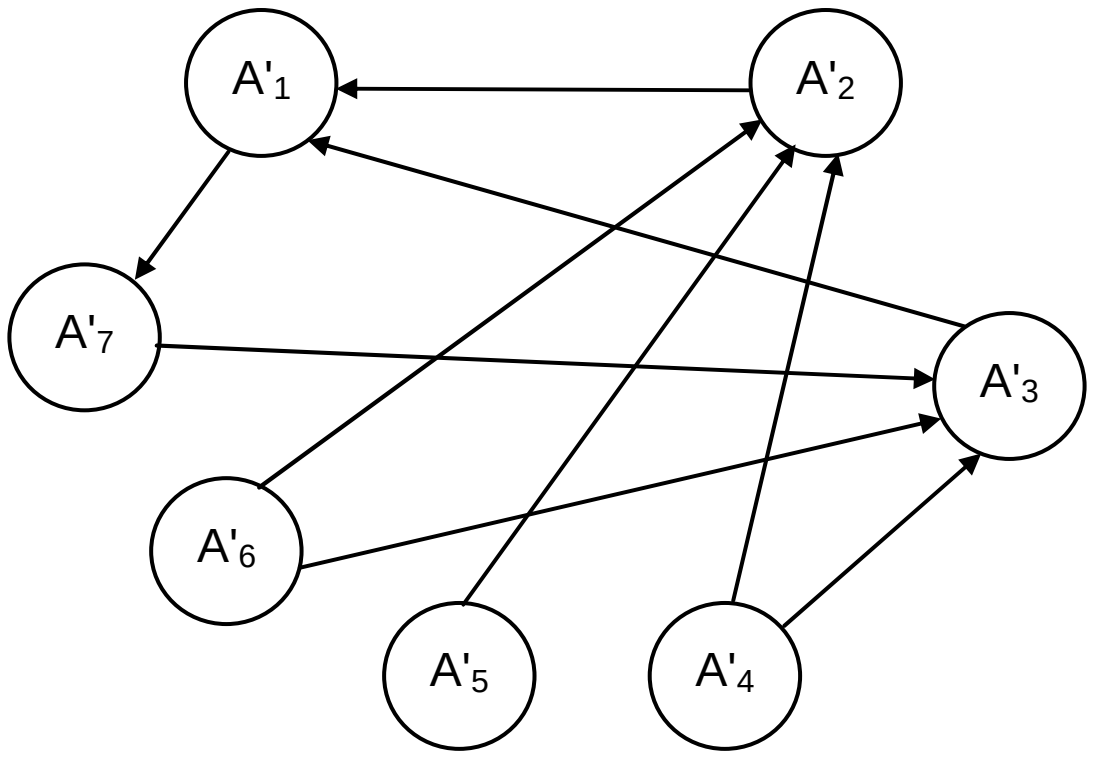

FIGURE 2.2.3.2

The related fuzzy relational matrix B is as follows:

$$B = \begin{array}{c} \\ A_1' \\ A_2' \\ A_3' \\ A_4' \\ A_5' \\ A_6' \\ A_7' \end{array} \begin{array}{c} \begin{array}{ccccccc} A_1' & A_2' & A_3' & A_4' & A_5' & A_6' & A_7' \end{array} \\ \begin{bmatrix} 0 & 0 & 0 & 0 & 0 & 0 & 1 \\ 1 & 0 & 0 & 0 & 0 & 0 & 0 \\ 1 & 0 & 0 & 0 & 0 & 0 & 0 \\ 0 & 1 & 1 & 0 & 0 & 0 & 0 \\ 0 & 1 & 0 & 0 & 0 & 0 & 0 \\ 0 & 1 & 1 & 0 & 0 & 0 & 0 \\ 0 & 0 & 1 & 0 & 0 & 0 & 0 \end{bmatrix} \end{array}$$

Suppose we consider the state vector X = (1 0 0 0 0 0 0) i.e., easy availability of money to be in the on state and all other nodes are in the off state. The effect of X on the dynamical system B is

$XB \hookrightarrow (1\ 0\ 0\ 0\ 0\ 0\ 1) = X_1$ (say)
$X_1B \hookrightarrow (1\ 0\ 1\ 0\ 0\ 0\ 1) = X_2$ (say)
$X_2B \hookrightarrow (1\ 0\ 1\ 0\ 0\ 0\ 1) = X_3 = X_2$

($X_2$ which a fixed point of the dynamical system.) i.e., easy money forces one to visit CSWs and he is also unaware of the disease and how it is communicated so he not only visits the CSWs but has unprotected sex with them there by infecting himself. As our motivation is to study the combined effect of the system we now for the same set of seven attributes seek the opinion of another expert. The directed graph given by the third expert is as follows:

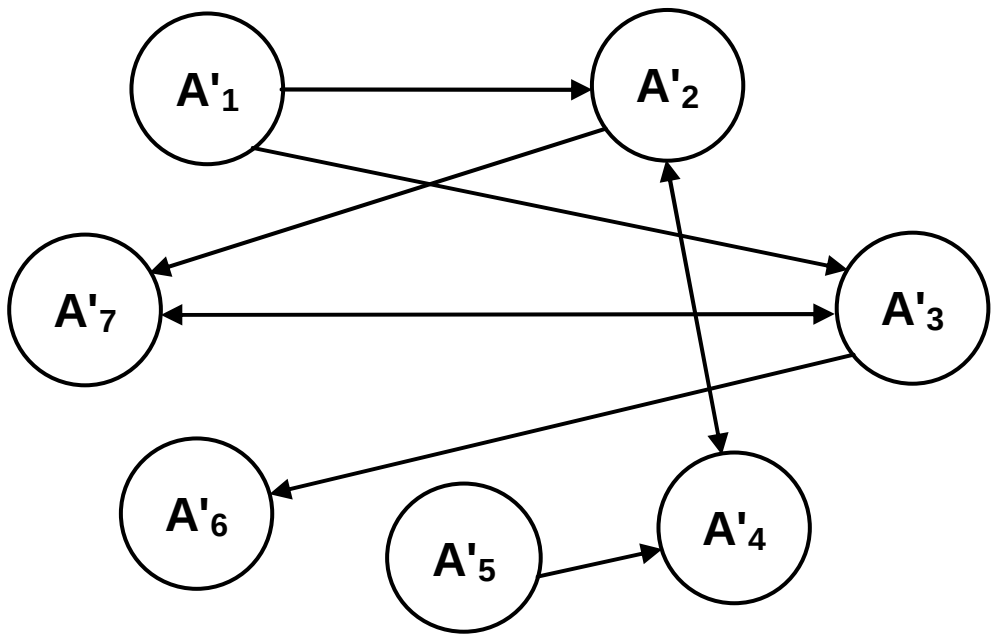

FIGURE 2.2.3.3

According to this expert when they have wrongful company and additive habits they will certainly become blind to awareness; for he feels as most of the migrant labourer who are HIV/AIDS patients acknowledge with full sense that they were in a complete drunken state when they visited CSWs and their statement they were not able to recollect whether they used *nirodh* (condom) is sufficient to explain when then have additive habits they fail to be aware of and to use protected sexual methods. Now we proceed on to give the related connection matrix C

$$
C = \begin{array}{c} \\ A_1^{'} \\ A_2^{'} \\ A_3^{'} \\ A_4^{'} \\ A_5^{'} \\ A_6^{'} \\ A_7^{'} \end{array}
\begin{array}{c} \begin{array}{ccccccc} A_1^{'} & A_2^{'} & A_3^{'} & A_4^{'} & A_5^{'} & A_6^{'} & A_7^{'} \end{array} \\
\left[\begin{array}{ccccccc}
0 & 1 & 1 & 0 & 0 & 0 & 0 \\
0 & 0 & 0 & 1 & 0 & 0 & 1 \\
0 & 0 & 0 & 0 & 0 & 1 & 1 \\
0 & 1 & 0 & 0 & 0 & 0 & 0 \\
0 & 0 & 0 & 1 & 0 & 0 & 0 \\
0 & 0 & 0 & 0 & 0 & 0 & 0 \\
0 & 0 & 1 & 0 & 0 & 0 & 0
\end{array}\right] \end{array}
$$

Suppose for this dynamical system C we consider the on state of the node visits CSWs to be in the on state and all other nodes are in the off state that is Y = (0 0 1 0 0 0 0); to analyze the effect of Y on C

| | | | | |
|---|---|---|---|---|
| YC | $\hookrightarrow$ | (0 0 1 0 0 1 1) | = | $Y_1$ (say) |
| $Y_1C$ | $\hookrightarrow$ | (0 0 1 0 0 1 1) | = | $Y_2$ = $Y_1$ (say). |

The hidden pattern is a fixed point which shows these migrant labourers who visit CSWs have more leisure and are not aware of how HIV/AIDS spreads.

Suppose we consider the on state of the vector easy availability of money to be in the on state the effect of this on the dynamical system C. Let P = (1 0 0 0 0 0 0). The effect of P on C is given by

| | | | | |
|---|---|---|---|---|
| PC | $\hookrightarrow$ | (1 1 1 0 0 0 0) | = | $P_1$ (say) |
| $P_1C$ | $\hookrightarrow$ | (1 1 1 1 0 1 1) | = | $P_2$ (say) |
| $P_2C$ | $\hookrightarrow$ | (1 1 1 1 0 1 1) | = | $P_3$ (a fixed point). |

Thus if they earn well they become easy victims of wrongful company and addictive habits, visits CSWs they are socially irresponsible, with more leisure and not knowing about the how the disease spreads.

Now we proceed on to take the fourth experts opinion. The directed graph given by him and the related matrix D is as follows:

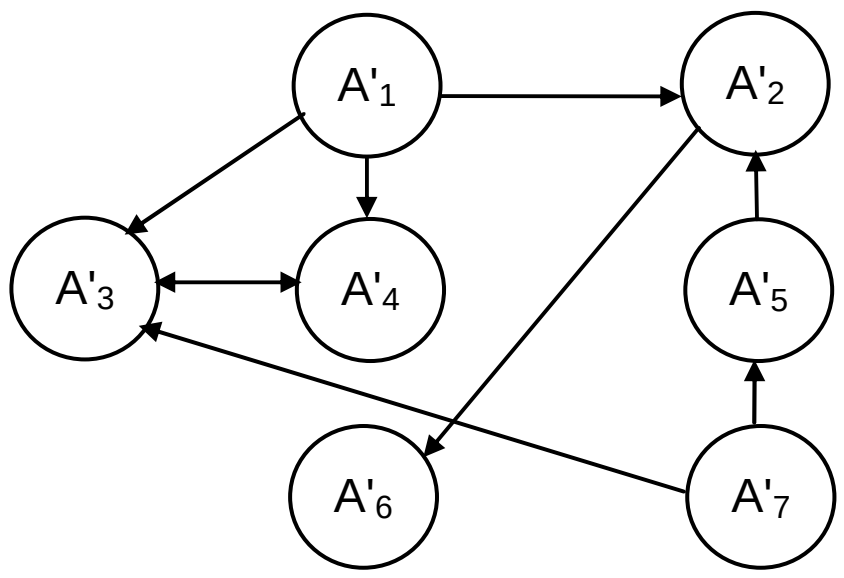

FIGURE 2.2.3.4

$$D = \begin{array}{c} \\ A_1' \\ A_2' \\ A_3' \\ A_4' \\ A_5' \\ A_6' \\ A_7' \end{array} \begin{array}{c} \begin{array}{ccccccc} A_1' & A_2' & A_3' & A_4' & A_5' & A_6' & A_7' \end{array} \\ \left[\begin{array}{ccccccc} 0 & 1 & 1 & 1 & 0 & 0 & 0 \\ 0 & 0 & 0 & 0 & 0 & 1 & 0 \\ 0 & 0 & 0 & 1 & 0 & 0 & 0 \\ 0 & 0 & 1 & 0 & 0 & 0 & 0 \\ 0 & 1 & 0 & 0 & 0 & 0 & 0 \\ 0 & 0 & 0 & 0 & 0 & 0 & 0 \\ 0 & 0 & 1 & 0 & 1 & 0 & 0 \end{array}\right] \end{array}$$

Suppose we consider the on state of the node socially irresponsible /free and all other nodes are in the off state. The effect of Z = (0 0 0 1 0 0 0) on the dynamical system D is given by

$$ZD \hookrightarrow (0\ 0\ 1\ 1\ 0\ 0\ 0) = Z_1 \text{ say}$$
$$Z_1 D \hookrightarrow (0\ 0\ 1\ 1\ 0\ 0\ 0) = Z_2 = Z_1$$

Thus the hidden pattern is a fixed point, leading to the conclusion socially irresponsible persons visit CSWs for they do not have binding on the family or even the responsibility that they would face health problem. Thus a socially free / irresponsible person is more prone to visit CSWs.

Now we take yet another experts opinion so that we can make conclusions on the combined FCMs. The directed graph given by the expert is in figure 2.2.3.5.

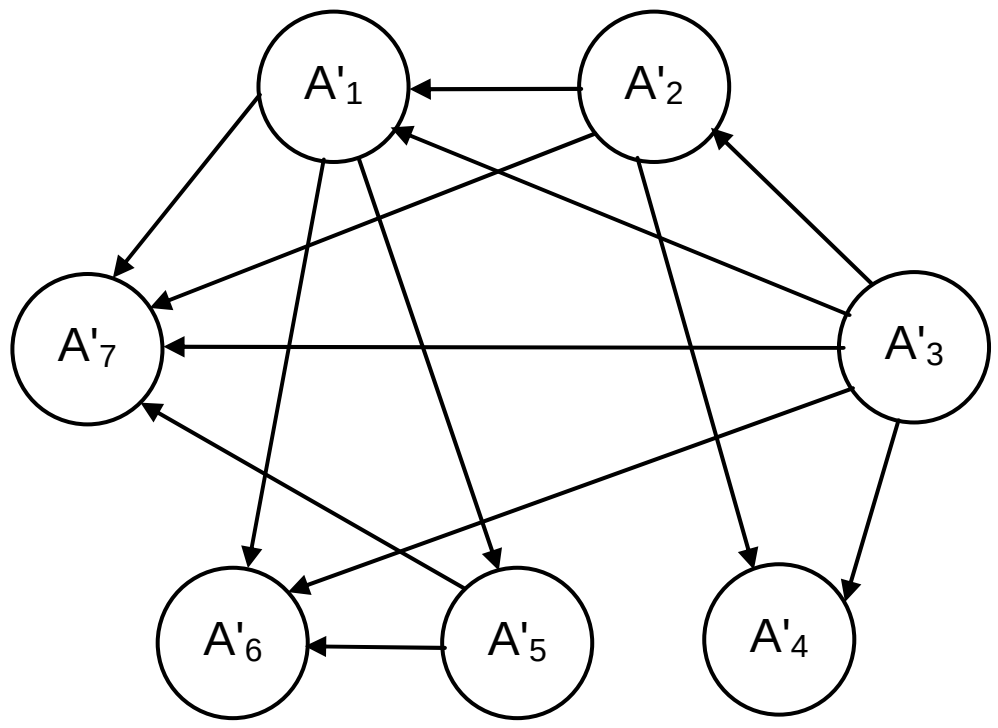

FIGURE 2.2.3.5

The related matrix E is given below

$$E = \begin{array}{c} \\ A_1' \\ A_2' \\ A_3' \\ A_4' \\ A_5' \\ A_6' \\ A_7' \end{array} \begin{array}{c} \begin{array}{ccccccc} A_1' & A_2' & A_3' & A_4' & A_5' & A_6' & A_7' \end{array} \\ \left[\begin{array}{ccccccc} 0 & 0 & 0 & 0 & 1 & 1 & 1 \\ 1 & 0 & 0 & 1 & 0 & 0 & 0 \\ 1 & 1 & 0 & 1 & 0 & 1 & 1 \\ 0 & 0 & 0 & 0 & 0 & 0 & 0 \\ 0 & 0 & 0 & 0 & 0 & 1 & 1 \\ 0 & 0 & 0 & 0 & 0 & 0 & 0 \\ 0 & 1 & 0 & 0 & 0 & 0 & 0 \end{array}\right] \end{array}$$

Suppose the node $A'_2$ i.e., wrongful company and addictive habits are in the on state.

The effect of

$$U = (0\ 1\ 0\ 0\ 0\ 0\ 0)$$

on the dynamical system E is given by

$$\begin{array}{llll} UE & \hookrightarrow & (1\ 1\ 0\ 1\ 0\ 0\ 0) & = U_1 \text{ (say)} \\ U_1E & \hookrightarrow & (1\ 1\ 0\ 1\ 1\ 1\ 1) & = U_2 \text{ (say)} \\ U_2E & \hookrightarrow & (1\ 1\ 0\ 1\ 1\ 1\ 1) & = U_3 = U_2. \end{array}$$

The hidden pattern is a fixed point. Thus according to this expert when a person is in a wrongful company and has additive habits it need not always imply that he will visit CSWs but for he has easy money, he is socially free /responsible, may have macho behaviour, he has more leisure and has no awareness of about the disease; as our data is an unsupervised one we have no means to change or modify the opinion of any expert we have to give his opinion as it is for otherwise the data would become biased.

Now we formulate the Combined Fuzzy Cognitive Maps using the opinion of the 5 experts.

Let S denote the combined connection matrix by S = A + B + C + D + E.

$$S = \begin{array}{c} \\ A_1' \\ A_2' \\ A_3' \\ A_4' \\ A_5' \\ A_6' \\ A_7' \end{array} \begin{array}{c} \begin{array}{ccccccc} A_1' & A_2' & A_3' & A_4' & A_5' & A_6' & A_7' \end{array} \\ \left[\begin{array}{ccccccc} 0 & 3 & 3 & 3 & 1 & 1 & 2 \\ 2 & 0 & 1 & 2 & 1 & 2 & 1 \\ 2 & 1 & 0 & 2 & 0 & 2 & 2 \\ 0 & 2 & 3 & 0 & 0 & 1 & 0 \\ 0 & 3 & 1 & 1 & 0 & 1 & 1 \\ 0 & 2 & 2 & 1 & 0 & 0 & 0 \\ 0 & 1 & 4 & 1 & 1 & 1 & 0 \end{array}\right] \end{array}$$

We threshold in a different way if a quality adds up to 2 or less than two we put 0 and if it is 3 or greater 3 than we put as 1.

Now we consider the effect of only easy availability of money to be in the on state.

The effect of

P = (1 0 0 0 0 0 0)

on the combined dynamical system S.

| | | | | |
|---|---|---|---|---|
| PS | = | (0 3 3 3 1 1 2) | | |
| | ↪ | (0 1 1 1 0 1 0) | = | $P_1$ |
| $P_1S$ | = | (4 5 6 5 1 4 3) | | |
| | ↪ | (1 1 1 1 1 1 1) | = | $P_2$ |
| $P_2S$ | = | (14 12 14 10 3 8 6) | | |
| | ↪ | (1 1 1 1 1 1 1) | = | $P_2$ |

The hidden pattern is a fixed point when they have easy availability of money they have all the states to be on. They visit CSW have bad habits, find more leisure (if no money they will do several types of work to earn) and as they are unaware of the disease they without any fear visit CSWs as they want to show their macho character they smoke and drink and so on. Thus the collective or combined FCM shows the on state of all vectors.

Suppose we consider the on state as the only attribute they have no awareness about the disease say L = (0 0 0 0 0 0 1) and all other attributes are off. We find the effect of L on the dynamical system E.

| | | | | |
|---|---|---|---|---|
| LE | = | (0 1 4 1 1 1 0) | | |
| | ↪ | (0 1 1 1 1 1 0) | = | $L_1$ say |
| $L_1E$ | = | (4 8 7 6 1 6 4) | | |
| | ↪ | (1 1 1 1 1 1 1) | = | $L_2$ say |
| $L_2E$ | ↪ | (1 1 1 1 1 1 1) | = | $L_3 = L_2$ say |

which is a fixed point of the hidden pattern. We can study and interpret the effect of each and every state vector using the C program given [220, 233, 235].

Also we have only illustrated using 7 attributes only with 5 experts, one can choose any number of attributes and use any desired number of experts and obtain the results using the C-program [220, 233, 235].

### 2.2.4. Combined Disjoint Block FCM and its Application

In this section we define for the first time the notion of Combined Disjoint Block FCM (CDB FCM) and apply them in the analysis of the socio economic problems of the HIV/AIDS affected migrant labourers.

Let $C_1, \ldots, C_n$ be n nodes / attributes related with some problem. The n may be very large and a non prime. Even though we have C- program to work finding the directed graph and the related connection matrix may be very unwieldy. In such cases we use the notion of combined disjoint block fuzzy cognitive maps.

We divide these n attributes into k equal classes and these k equal classes are viewed by k-experts or by the same expert and the corresponding directed graph and the connection matrices are got.

Now these connection matrices is made into a $n \times n$ matrix and using the C-program the results are derived. This type of Combined Disjoint Block FCM is know as Combined Disjoint Equal Block FCM. Now some times we may not be in a position to divide the 'n' under study into equal blocks in such cases we use the technique of dividing the n attributes say of some m blocks were each block may not have the same number of attributes, but it is essential that there n attributes are divided into disjoint classes. We use both these techniques in the analysis of the problem. This is the case when n happens to be a prime number.

We in the next stage proceed on to define the notion of combined overlap block Fuzzy cognitive maps, for this also we assume two types the overlap is regular and the overlap happens to be irregular. Both the cases are dealt in this book.

**DEFINITION 2.2.4.1:** *Let $C_1 \ldots C_n$ be n distinct attributes of a problem n very large and a non prime. If we divide n into k equal classes i.e., k/n and if n/k = t which are disjoint and if we find the directed graph of each of there k classes of attributes with t attributes each, then their corresponding connection matrices are formed and these connection matrices are joined as blocks to form a $n \times n$ matrix.*

*This $n \times n$ connection matrix forms the combined disjoint block FCM of equal classes. If the classes are not divided to have equal number of attributers but if they are disjoint classes we get a $n \times n$ connection matrix called the combined disjoint block FCM of unequal classes/ size.*

***Example 2.2.4.1:*** We illustrate first our models by a simple example. The attributes taken here are related with the HIV/AIDS affected Migrant Labourers. Suppose the 12 attributes $A_1, A_2, A_3, A_4, \ldots, A_{12}$ is given by the experts

| | | |
|---|---|---|
| $A_1$ | - | Easy availability of money |
| $A_2$ | - | Lack of education |
| $A_3$ | - | Visiting CSWs |
| $A_4$ | - | Nature of profession |
| $A_5$ | - | Wrong / Bad company |
| $A_6$ | - | Addiction to habit forming substances and visiting CSWs. |
| $A_7$ | - | Absence of social responsibility |
| $A_8$ | - | Socially free |
| $A_9$ | - | Economic status |
| $A_{10}$ | - | More leisure |

$A_{11}$ - Machismo / Exaggerated masculinity
$A_{12}$ - No awareness of the disease.

Using these attributes we give the combined disjoint block fuzzy cognitive map of equal size. We take 3 experts opinion on the 3 disjoint classes so that each class has four attributes / nodes. Let the disjoint classes be $C_1$, $C_2$ and $C_3$ be divided by the following:

$C_1 = \{A_1, A_6, A_7, A_{12}\}$,
$C_2 = \{A_2, A_3, A_4, A_{10}\}$ and
$C_3 = \{A_5, A_8, A_9, A_{11}\}$.

Now we collect the experts opinion on each of the classes $C_1$, $C_2$ and $C_3$.

The directed graph given by the expect on attributes $A_1$, $A_6$, $A_7$ and $A_{12}$ which forms the class $C_1$.

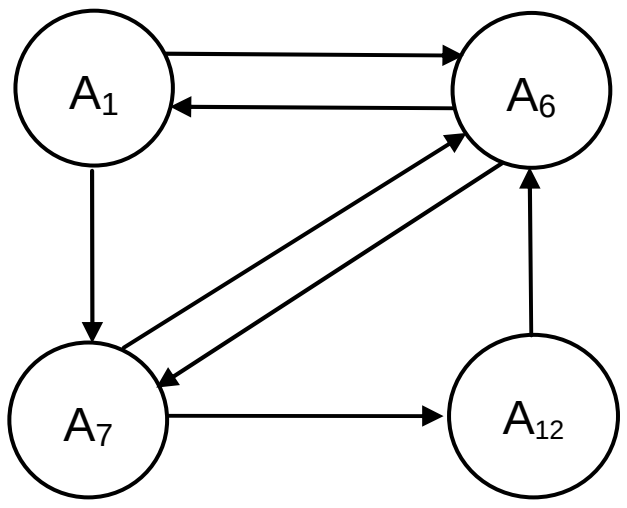

FIGURE 2.2.4.1

The related connection matrix $B_1$

$$B_1 = \begin{array}{c} \\ A_1 \\ A_6 \\ A_7 \\ A_{12} \end{array} \begin{array}{c} \begin{array}{cccc} A_1 & A_6 & A_7 & A_{12} \end{array} \\ \begin{bmatrix} 0 & 1 & 1 & 0 \\ 1 & 0 & 1 & 0 \\ 0 & 1 & 0 & 1 \\ 0 & 1 & 0 & 0 \end{bmatrix} \end{array}$$

The directed graph given by the expert on $A_2$, $A_3$, $A_4$, $A_{10}$ which forms the class $C_2$.

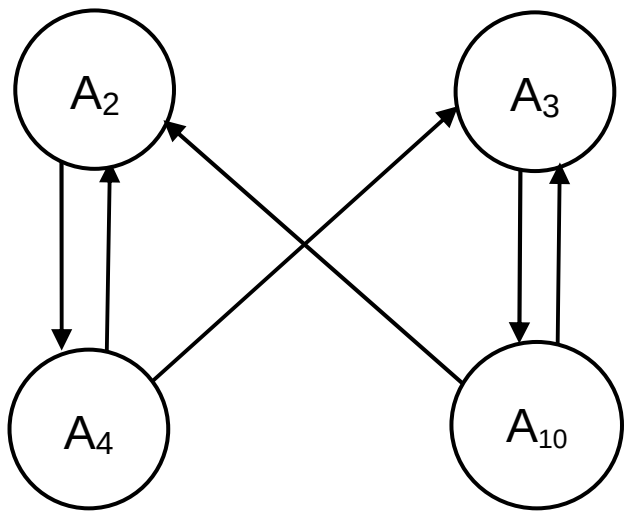

FIGURE 2.2.4.2

According to this expert the two nodes, lack of education and the nature of profession are very much interrelated. Also he feels most of the educated people have awareness about HIV/AIDS. The related connection matrix $B_2$ is given below:

$$B_2 = \begin{array}{c} \\ A_2 \\ A_3 \\ A_4 \\ A_{10} \end{array} \begin{array}{c} \begin{array}{cccc} A_2 & A_3 & A_4 & A_{10} \end{array} \\ \begin{bmatrix} 0 & 0 & 1 & 0 \\ 0 & 0 & 0 & 1 \\ 1 & 1 & 0 & 0 \\ 1 & 1 & 0 & 0 \end{bmatrix} \end{array}$$

Now we give the directed graph for the class $C_3$ as given by the expert $C_3 = \{A_5\ A_8\ A_9\ A_{11}\}$.

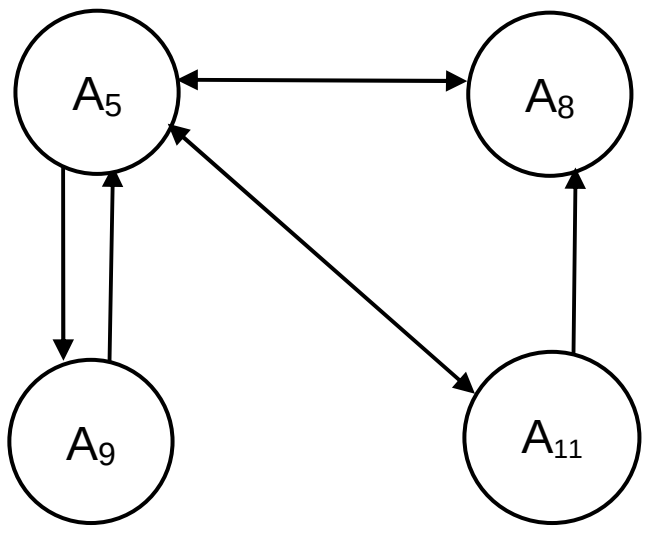

FIGURE 2.2.4.3

According to this expert bad company induces Macho behaviour. Also when a person has bad habit / bad company he is not bounded socially he without any fear engages himself in doing social evil, so he relates bad company with socially free.

If he has no food for the next day, i.e., with this economic status certainly he cannot enjoy or have bad company or bad habits. Usually person with bad company show more of macho behaviour. The related connection matrix

$$B_3 = \begin{array}{c} \\ A_5 \\ A_8 \\ A_9 \\ A_{11} \end{array} \begin{array}{c} \begin{array}{cccc} A_5 & A_8 & A_9 & A_{11} \end{array} \\ \begin{bmatrix} 0 & 1 & 1 & 1 \\ 1 & 0 & 0 & 0 \\ 1 & 0 & 0 & 0 \\ 0 & 1 & 0 & 0 \end{bmatrix} \end{array}$$

Now the combined disjoint block connection matrix of the FCM is given by B

$$\begin{array}{c} \\ A_1 \\ A_6 \\ A_7 \\ A_{12} \\ A_2 \\ A_3 \\ A_4 \\ A_{10} \\ A_5 \\ A_8 \\ A_9 \\ A_{11} \end{array} \begin{array}{c} \begin{array}{cccccccccccc} A_1 & A_6 & A_7 & A_{12} & A_2 & A_3 & A_4 & A_{10} & A_5 & A_8 & A_9 & A_{11} \end{array} \\ \begin{bmatrix} 0 & 1 & 1 & 0 & 0 & 0 & 0 & 0 & 0 & 0 & 0 & 0 \\ 1 & 0 & 1 & 0 & 0 & 0 & 0 & 0 & 0 & 0 & 0 & 0 \\ 0 & 1 & 0 & 1 & 0 & 0 & 0 & 0 & 0 & 0 & 0 & 0 \\ 0 & 1 & 0 & 0 & 0 & 0 & 0 & 0 & 0 & 0 & 0 & 0 \\ 0 & 0 & 0 & 0 & 0 & 0 & 1 & 0 & 0 & 0 & 0 & 0 \\ 0 & 0 & 0 & 0 & 0 & 0 & 0 & 1 & 0 & 0 & 0 & 0 \\ 0 & 0 & 0 & 0 & 1 & 1 & 0 & 0 & 0 & 0 & 0 & 0 \\ 0 & 0 & 0 & 0 & 1 & 1 & 0 & 0 & 0 & 0 & 0 & 0 \\ 0 & 0 & 0 & 0 & 0 & 0 & 0 & 0 & 0 & 1 & 1 & 1 \\ 0 & 0 & 0 & 0 & 0 & 0 & 0 & 0 & 1 & 0 & 0 & 0 \\ 0 & 0 & 0 & 0 & 0 & 0 & 0 & 0 & 1 & 0 & 0 & 0 \\ 0 & 0 & 0 & 0 & 0 & 0 & 0 & 0 & 0 & 1 & 0 & 0 \end{bmatrix} \end{array}$$

Suppose we consider the on state of the attribute easy availability of the money and all other states are off the effect of X = (1 0 0 0 0 0 0 0 0 0 0 0) on the CDB FCM is given by

$$\begin{array}{lllll} XB & \hookrightarrow & (1\ 1\ 1\ 0\ 0\ 0\ 0\ 0\ 0\ 0\ 0\ 0) & = & X_1 \text{ (Say)} \\ X_1B & \hookrightarrow & (1\ 1\ 1\ 1\ 0\ 0\ 0\ 0\ 0\ 0\ 0\ 0) & = & X_2 \text{ (Say)} \\ X_2B & \hookrightarrow & (1\ 1\ 1\ 1\ 0\ 0\ 0\ 0\ 0\ 0\ 0\ 0) & = & X_3 = X_2. \end{array}$$

($X_2$ is a fixed point of the dynamical system). Thus when one has easy money he has addiction to habit forming substances and visits CSWs he has no social responsibility and he is not aware of how the disease spreads.

Suppose we consider, the on state of the attributes visits CSWs and socially free to be in the on state and all other nodes are in the off state. Now we study the effect on the dynamical system B. Let T = (0 0 0 0 0 1 0 0 0 1 0 0) state vector depicting the on state of CSWs and socially free, passing the state vector T in to the dynamical system B.

$$TB \hookrightarrow (0\ 0\ 0\ 0\ 0\ 1\ 0\ 1\ 1\ 1\ 0\ 0) = T_1 \text{ (say)}$$
$$T_1B \hookrightarrow (0\ 0\ 0\ 0\ 1\ 1\ 1\ 1\ 1\ 1\ 1\ 1) = T_2 \text{ (say)}$$
$$T_2B \hookrightarrow (0\ 0\ 0\ 0\ 1\ 1\ 1\ 1\ 1\ 1\ 1\ 1) = T_3\ .$$

($T_2$ the fixed point of the dynamical system. Thus when he visits CSWs and is socially free one gets the lack of education, money is easily available to him, Absence of social responsibility and no awareness about the disease is in the off state and all other states become on. One can study this dynamical system using C-program given in [220, 233, 235].

***Example 2.2.4.2:*** Now we give another application of the Combined disjoint equal block fuzzy cognitive maps to the HIV/AIDS affected migrant labourers. Now we consider the following 15 conceptual nodes $S_1$, …, $S_{15}$ associated with the HIV/AIDS patients and migrancy

$S_1$ - No awareness of HIV/AIDS in migrant labourers
$S_2$ - Rural living with no education
$S_3$ - Migrancy as truck drivers / daily wagers
$S_4$ - Socially free and irresponsible
$S_5$ - Enjoy life jolly mood
$S_6$ - Away from family
$S_7$ - More leisure
$S_8$ - No association or union to protect and channelize their money / time
$S_9$ - Visits CSWs
$S_{10}$ - Unreachable by friends and relatives

$S_{11}$ - Easy victims of temptation
$S_{12}$ - Takes to all bad habit / bad company including CSWs
$S_{13}$ - Male chauvinism
$S_{14}$ - Women as inferior objects
$S_{15}$ - No family binding and respect for wife.

These 15 attributes are divided in 5 classes $C_1, C_2, \ldots, C_5$ with 3 in each class.

Let $C_1 = \{S_1, S_2, S_3\}$,
$C_2 = \{S_4, S_5, S_6\}$,
$C_3 = \{S_7, S_8, S_9\}$,
$C_4 = \{S_{10}, S_{11}, S_{12}\}$
and $C_5 = \{S_{13}, S_{14}, S_{15}\}$.

Now we take the experts opinion for each of these classes and take the matrix associated with the combined block disjoint FCMs.

Using the $15 \times 15$ matrix got using the combined connection matrix we derive our conclusions. The experts opinion for the class $C_1 = \{S_1, S_2, S_3\}$ in the form of the directed graph.

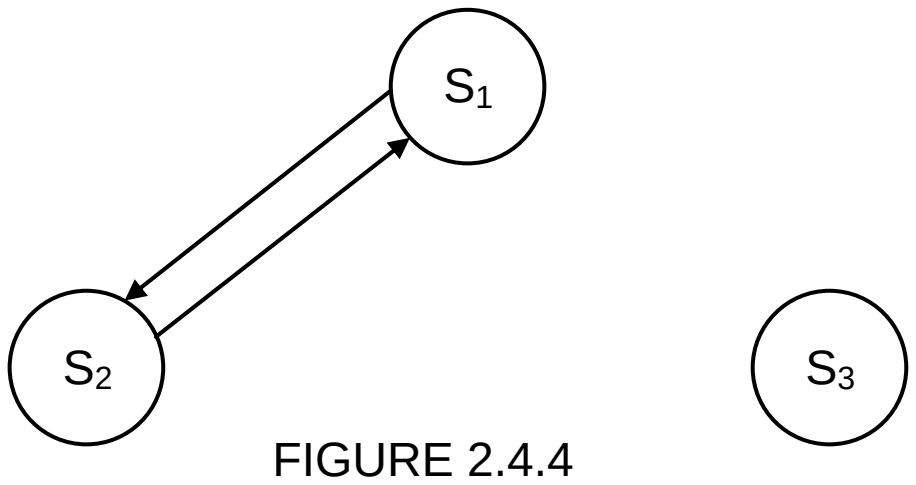

FIGURE 2.4.4

The related connection matrix is given by

$$\begin{array}{c} \\ S_1 \\ S_2 \\ S_3 \end{array}\begin{array}{c} \begin{array}{ccc} S_1 & S_2 & S_3 \end{array} \\ \begin{bmatrix} 0 & 1 & 0 \\ 1 & 0 & 0 \\ 0 & 0 & 0 \end{bmatrix} \end{array}.$$

The directed graph for the class $C_2 = \{S_4, S_5, S_6\}$ as given by the expert

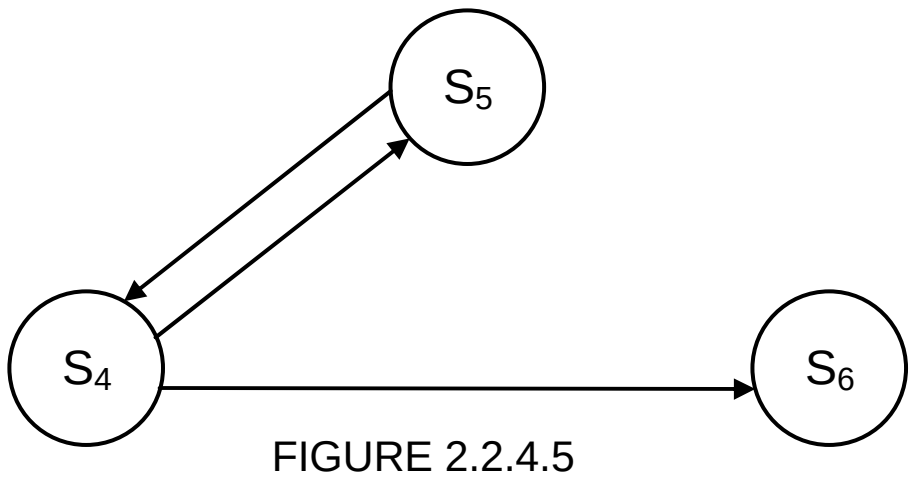

FIGURE 2.2.4.5

The related connection matrix is given as follows:

$$\begin{array}{c} \\ S_4 \\ S_5 \\ S_6 \end{array} \begin{array}{c} \begin{array}{ccc} S_4 & S_5 & S_6 \end{array} \\ \begin{bmatrix} 0 & 1 & 1 \\ 1 & 0 & 0 \\ 0 & 0 & 0 \end{bmatrix} \end{array}.$$

Now for the class $C_3 = \{S_7, S_8, S_9\}$ the directed graph is given below by the expert

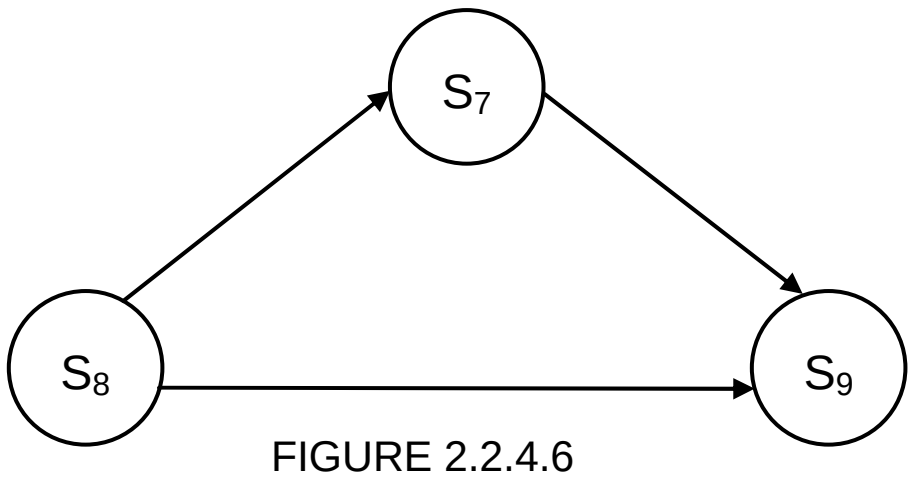

FIGURE 2.2.4.6

The related connection matrix

$$\begin{array}{c} \\ S_7 \\ S_8 \\ S_9 \end{array} \begin{array}{c} \begin{array}{ccc} S_7 & S_8 & S_9 \end{array} \\ \begin{bmatrix} 0 & 0 & 1 \\ 1 & 0 & 1 \\ 0 & 0 & 0 \end{bmatrix} \end{array}.$$

The directed graph given by the expert for the class $C_4 = \{S_{10}, S_{11}, S_{12}\}$ is given below:

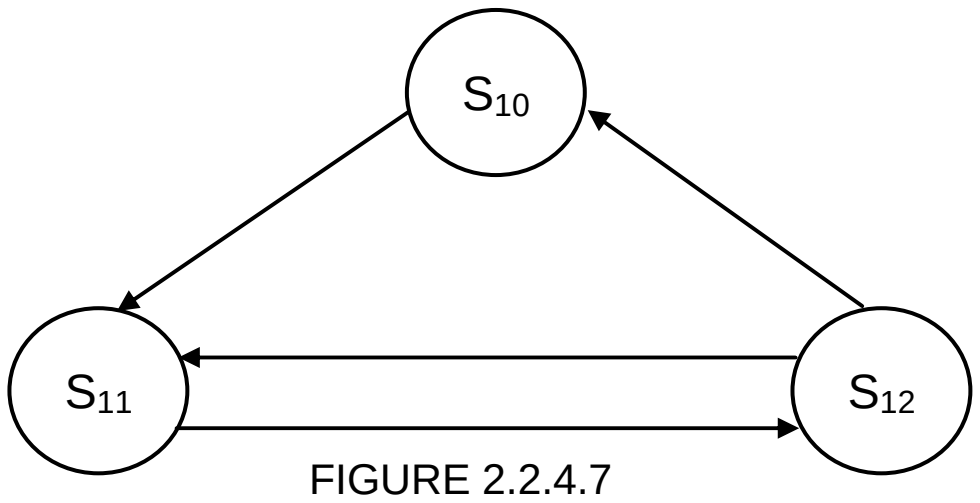

FIGURE 2.2.4.7

The related connection matrix

$$\begin{array}{c} \\ S_{10} \\ S_{11} \\ S_{12} \end{array}\begin{array}{c} \begin{array}{ccc} S_{10} & S_{11} & S_{12} \end{array} \\ \begin{bmatrix} 0 & 1 & 0 \\ 0 & 0 & 1 \\ 1 & 1 & 0 \end{bmatrix} \end{array}.$$

Now for the class $C_5 = \{S_{13}, S_{14}, S_{15}\}$, the directed graph is as follows:

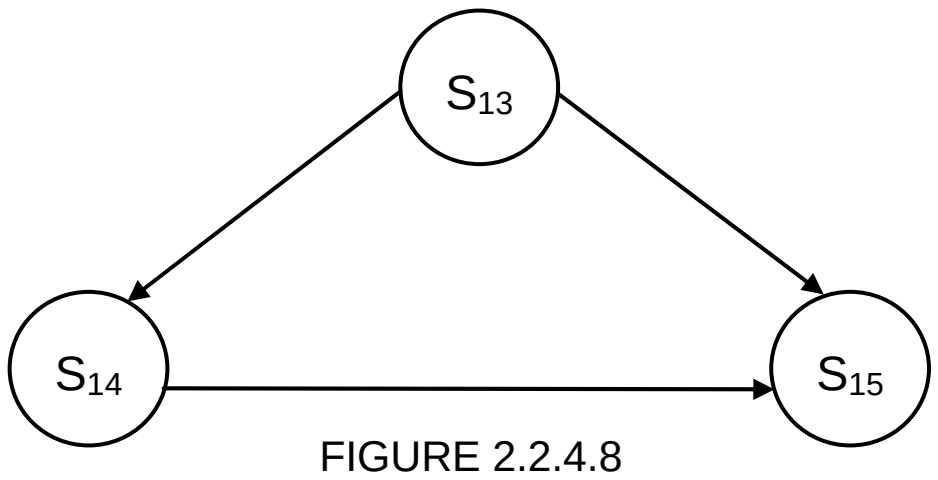

FIGURE 2.2.4.8

The related connection matrix

$$\begin{array}{c} \\ S_{13} \\ S_{14} \\ S_{15} \end{array}\begin{array}{c} \begin{array}{ccc} S_{13} & S_{14} & S_{15} \end{array} \\ \begin{bmatrix} 0 & 1 & 1 \\ 0 & 0 & 1 \\ 0 & 0 & 0 \end{bmatrix} \end{array}.$$

Now we give the combined block connection matrix related with the 15 attributes, let us denote it by S.

$$
S = \begin{array}{c} \\ S_1 \\ S_2 \\ S_3 \\ S_4 \\ S_5 \\ S_6 \\ S_7 \\ S_8 \\ S_9 \\ S_{10} \\ S_{11} \\ S_{12} \\ S_{13} \\ S_{14} \\ S_{15} \end{array}
\begin{array}{c}
\begin{array}{ccccccccccccccc} S_1 & S_2 & S_3 & S_4 & S_5 & S_6 & S_7 & S_8 & S_9 & S_{10} & S_{11} & S_{12} & S_{13} & S_{14} & S_{15} \end{array} \\
\left[\begin{array}{ccccccccccccccc}
0 & 1 & 0 & 0 & 0 & 0 & 0 & 0 & 0 & 0 & 0 & 0 & 0 & 0 & 0 \\
1 & 0 & 0 & 0 & 0 & 0 & 0 & 0 & 0 & 0 & 0 & 0 & 0 & 0 & 0 \\
0 & 0 & 0 & 0 & 0 & 0 & 0 & 0 & 0 & 0 & 0 & 0 & 0 & 0 & 0 \\
0 & 0 & 0 & 0 & 1 & 1 & 0 & 0 & 0 & 0 & 0 & 0 & 0 & 0 & 0 \\
0 & 0 & 0 & 1 & 0 & 0 & 0 & 0 & 0 & 0 & 0 & 0 & 0 & 0 & 0 \\
0 & 0 & 0 & 0 & 0 & 0 & 0 & 0 & 0 & 0 & 0 & 0 & 0 & 0 & 0 \\
0 & 0 & 0 & 0 & 0 & 0 & 0 & 0 & 1 & 0 & 0 & 0 & 0 & 0 & 0 \\
0 & 0 & 0 & 0 & 0 & 0 & 1 & 0 & 1 & 0 & 0 & 0 & 0 & 0 & 0 \\
0 & 0 & 0 & 0 & 0 & 0 & 0 & 0 & 0 & 0 & 0 & 0 & 0 & 0 & 0 \\
0 & 0 & 0 & 0 & 0 & 0 & 0 & 0 & 0 & 0 & 1 & 0 & 0 & 0 & 0 \\
0 & 0 & 0 & 0 & 0 & 0 & 0 & 0 & 0 & 0 & 0 & 1 & 0 & 0 & 0 \\
0 & 0 & 0 & 0 & 0 & 0 & 0 & 0 & 0 & 1 & 1 & 0 & 0 & 0 & 0 \\
0 & 0 & 0 & 0 & 0 & 0 & 0 & 0 & 0 & 0 & 0 & 0 & 0 & 1 & 1 \\
0 & 0 & 0 & 0 & 0 & 0 & 0 & 0 & 0 & 0 & 0 & 0 & 0 & 0 & 1 \\
0 & 0 & 0 & 0 & 0 & 0 & 0 & 0 & 0 & 0 & 0 & 0 & 0 & 0 & 0
\end{array}\right]
\end{array}.
$$

Now consider a state vector X = (0 1 0 0 0 1 0 0 0 1 0 0 0 1 0) where the nodes/concepts $S_2$, $S_6$, $S_{10}$ and $S_{14}$ are in the on state. The effect of X on the dynamical system S is given by

$$
\begin{array}{llll}
XS & \hookrightarrow & (1\ 1\ 0\ 0\ 0\ 1\ 0\ 0\ 0\ 1\ 1\ 0\ 0\ 1\ 1) & = X_1 \text{ say} \\
X_1S & \hookrightarrow & (1\ 1\ 0\ 0\ 0\ 1\ 0\ 0\ 0\ 1\ 1\ 1\ 0\ 1\ 1) & = X_2 \text{ say}
\end{array}
$$

and so on, one arrive at a fixed point.

### 2.2.5. Combined Overlap Block FCM and its Use in the Analysis of HIV/AIDS Affected Migrant Labourers

Next we give the model of the combined block overlap fuzzy cognitive maps and adapt it mainly on the HIV/AIDS affected migrant labourers relative to their socio economic problems.

***Example 2.2.5.1:*** We adapt it to the model $\{A_1, A_2, \ldots, A_{12}\}$ given in section 2.2.4 of this book. Let us consider the class

$C_1 = \{A_1\ A_2\ A_3\ A_4\}$, $C_2 = \{A_3\ A_4\ A_5\ A_6\}$,
$C_3 = \{A_5\ A_6\ A_7\ A_8\}$, $C_4 = \{A_7\ A_8\ A_9\ A_{10}\}$,
$C_5 = \{A_9\ A_{10}\ A_{11}\ A_{12}\}$ and $C_6 = \{A_{11}\ A_{12}\ A_1\ A_2\}$.

We give the directed graph for each of the classes of attributes $C_1$, $C_2$, …, $C_6$. The directed graph for the four attributes given by $C_1 = \{A_1\ A_2\ A_3\ A_4\}$ as given by the expert is as follows.

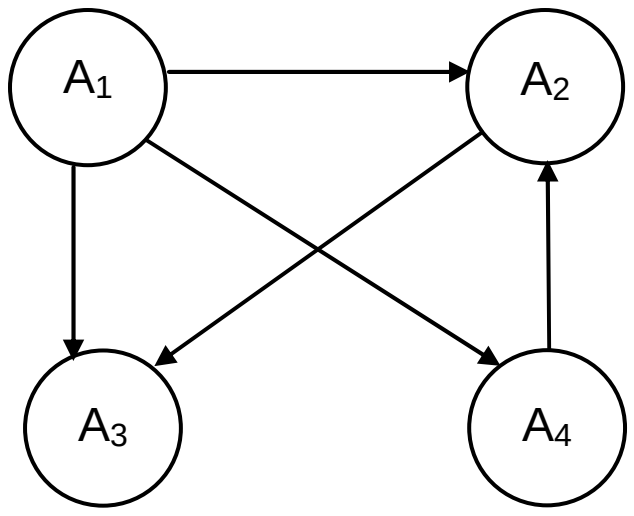

FIGURE 2.2.5.1

The related connection matrix

$$M_1 = \begin{matrix} & \begin{matrix} A_1 & A_2 & A_3 & A_4 \end{matrix} \\ \begin{matrix} A_1 \\ A_2 \\ A_3 \\ A_4 \end{matrix} & \begin{bmatrix} 0 & 1 & 1 & 1 \\ 0 & 0 & 1 & 0 \\ 0 & 0 & 0 & 0 \\ 0 & 1 & 0 & 0 \end{bmatrix} \end{matrix}.$$

The directed graph for the class of attributes $C_2 = \{A_3\ A_4\ A_5\ A_6\}$ as given by the expert is as follows:

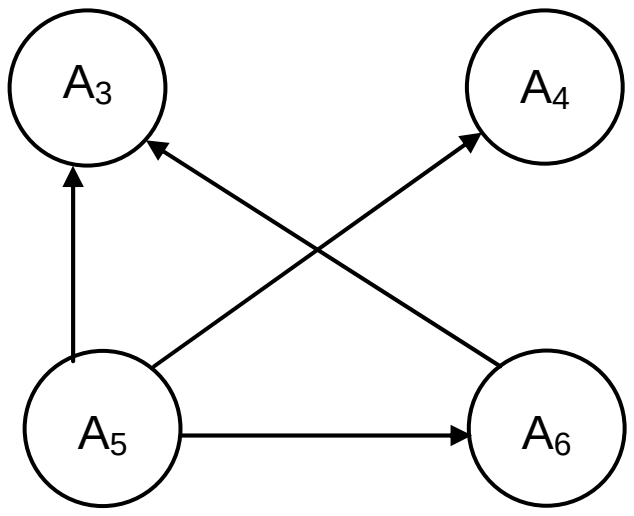

FIGURE 2.2.5.2

The related connection matrix

$$M_2 = \begin{array}{c} \\ A_3 \\ A_4 \\ A_5 \\ A_6 \end{array}\begin{array}{c} \begin{array}{cccc} A_3 & A_4 & A_5 & A_6 \end{array} \\ \begin{bmatrix} 0 & 0 & 0 & 0 \\ 0 & 0 & 0 & 0 \\ 1 & 1 & 0 & 1 \\ 1 & 0 & 0 & 0 \end{bmatrix} \end{array}.$$

The directed graph associated with the attributes $C_3 = \{A_5\ A_6\ A_7\ A_8\}$ as given by the expert.

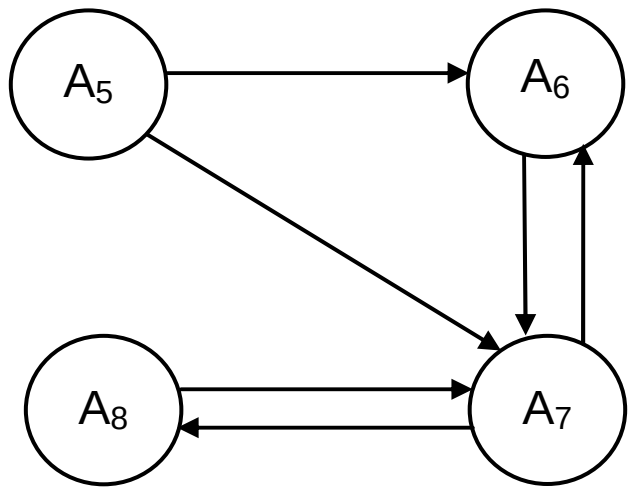

FIGURE 2.2.5.3

The related connection matrix

$$\begin{array}{c} \\ A_5 \\ A_6 \\ A_7 \\ A_8 \end{array}\begin{array}{c} \begin{array}{cccc} A_5 & A_6 & A_7 & A_8 \end{array} \\ \begin{bmatrix} 0 & 1 & 1 & 0 \\ 0 & 0 & 1 & 0 \\ 0 & 1 & 0 & 1 \\ 0 & 0 & 1 & 0 \end{bmatrix} \end{array}.$$

For the class of attributes $C_4 = \{A_7\ A_8\ A_9\ A_{10}\}$ the directed graph as given by the expert is as follows.

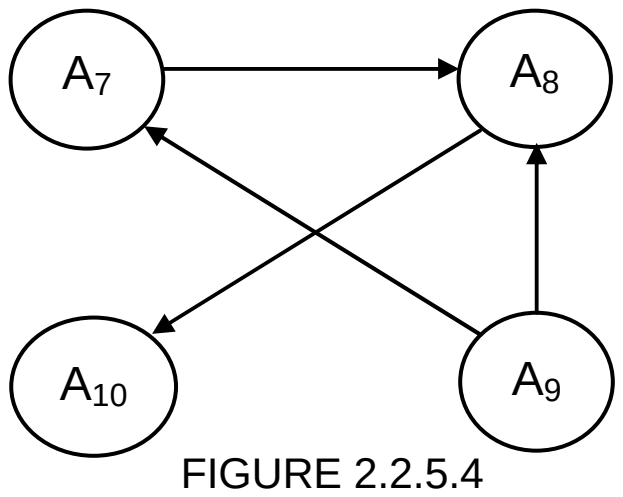

FIGURE 2.2.5.4

The connection matrix associated with the directed graph.

$$\begin{array}{c} \\ A_7 \\ A_8 \\ A_9 \\ A_{10} \end{array}\begin{array}{c} \begin{array}{cccc} A_7 & A_8 & A_9 & A_{10} \end{array} \\ \left[\begin{array}{cccc} 0 & 1 & 0 & 0 \\ 0 & 0 & 0 & 1 \\ 1 & 1 & 0 & 0 \\ 0 & 0 & 0 & 0 \end{array}\right]. \end{array}$$

The directed graph associated with the class C = {$A_9$ $A_{10}$ $A_{11}$ $A_{12}$} as given by the expert.

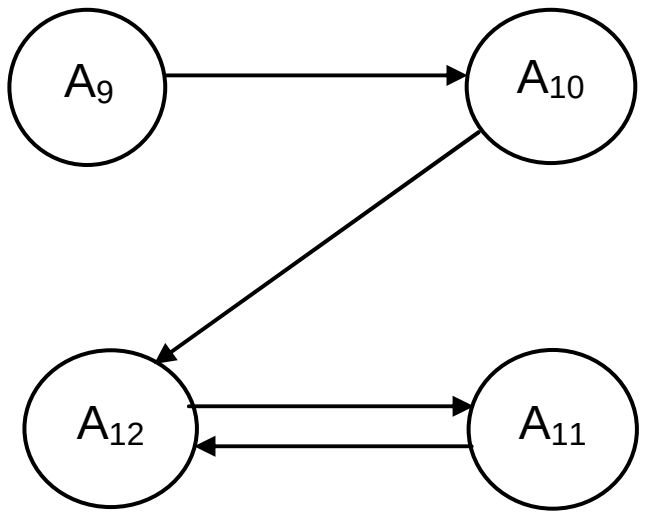

FIGURE 2.2.5.5

The related connection matrix is

$$\begin{array}{c} \\ A_9 \\ A_{10} \\ A_{11} \\ A_{12} \end{array}\begin{array}{c} \begin{array}{cccc} A_9 & A_{10} & A_{11} & A_{12} \end{array} \\ \left[\begin{array}{cccc} 0 & 1 & 0 & 0 \\ 0 & 0 & 0 & 1 \\ 0 & 0 & 0 & 1 \\ 0 & 0 & 1 & 0 \end{array}\right]. \end{array}$$

Now we give the directed graph of the last class given by the expert

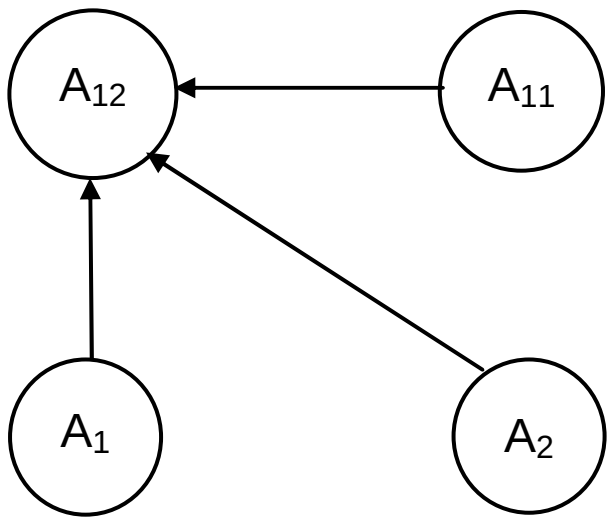

FIGURE 2.2.5.6

The related connection matrix

$$
\begin{array}{c} \\ A_{11} \\ A_{12} \\ A_1 \\ A_2 \end{array}
\begin{array}{c} \begin{array}{cccc} A_{11} & A_{12} & A_1 & A_2 \end{array} \\
\begin{bmatrix} 0 & 1 & 0 & 0 \\ 0 & 0 & 0 & 0 \\ 0 & 1 & 0 & 0 \\ 0 & 1 & 0 & 0 \end{bmatrix} \end{array}.
$$

The related Combined Block Overlap Fuzzy Cognitive Map (CBOFCM), matrix is as follows. We denote this matrix by W.

$$
\begin{array}{c} \\ A_1 \\ A_2 \\ A_3 \\ A_4 \\ A_5 \\ A_6 \\ A_7 \\ A_8 \\ A_9 \\ A_{10} \\ A_{11} \\ A_{12} \end{array}
\begin{array}{c} \begin{array}{cccccccccccc} A_1 & A_2 & A_3 & A_4 & A_5 & A_6 & A_7 & A_8 & A_9 & A_{10} & A_{11} & A_{12} \end{array} \\
\begin{bmatrix}
0 & 1 & 1 & 1 & 0 & 0 & 0 & 0 & 0 & 0 & 0 & 1 \\
0 & 0 & 1 & 0 & 0 & 0 & 0 & 0 & 0 & 0 & 0 & 1 \\
0 & 0 & 0 & 0 & 0 & 0 & 0 & 0 & 0 & 0 & 0 & 0 \\
0 & 1 & 0 & 0 & 0 & 0 & 0 & 0 & 0 & 0 & 0 & 0 \\
0 & 0 & 1 & 1 & 0 & 2 & 1 & 0 & 0 & 0 & 0 & 0 \\
0 & 0 & 0 & 1 & 0 & 0 & 1 & 0 & 0 & 0 & 0 & 0 \\
0 & 0 & 0 & 0 & 0 & 1 & 0 & 2 & 0 & 0 & 0 & 0 \\
0 & 0 & 0 & 0 & 0 & 0 & 1 & 0 & 0 & 1 & 0 & 0 \\
0 & 0 & 0 & 0 & 0 & 0 & 1 & 1 & 0 & 1 & 0 & 0 \\
0 & 0 & 0 & 0 & 0 & 0 & 0 & 0 & 0 & 0 & 0 & 1 \\
0 & 0 & 0 & 0 & 0 & 0 & 0 & 0 & 0 & 0 & 0 & 2 \\
0 & 0 & 0 & 0 & 0 & 0 & 0 & 0 & 0 & 0 & 1 & 0
\end{bmatrix} \end{array}.
$$

Let us consider the state vector X = (1 0 0 0 0 0 0 0 0 0 0 0) where easy availability of money alone is in the on state and all vectors are in the off state. The effect of X on the dynamical system W is give by

$$
\begin{array}{llll}
XW & \hookrightarrow & (1\ 1\ 1\ 1\ 0\ 0\ 0\ 0\ 0\ 0\ 0\ 1) & = X_1 \text{ (say)} \\
X_1W & \hookrightarrow & (1\ 1\ 1\ 1\ 0\ 0\ 0\ 0\ 0\ 0\ 1\ 1) & = X_2 \text{ say} \\
X_2\,W & \hookrightarrow & (1\ 1\ 1\ 1\ 0\ 0\ 0\ 0\ 0\ 0\ 1\ 1) & = X_3 \text{ (say)}
\end{array}
$$

Thus $X_2 = X_3$.

The hidden pattern of the CBOFCM is a fixed point.

***Example 2.2.5.2:*** We have described in this book 15 attributes denoted by $S_1, S_2, \ldots, S_{15}$ in page 182-3.

Now for some notational convenience take $S_i = P_i$, i = 1, 2, …, 15.

The attributes associated with it are $P_1, P_2, \ldots, P_{15}$. We format the following classes.

$C_1 = \{P_1\ P_2\ P_3\ P_4\ P_5\}$, $C_2 = \{P_3\ P_4\ P_5\ P_6\ P_7\}$
$C_3 = \{P_6\ P_7\ P_8\ P_9\ P_{10}\}$, $C_4 = \{P_8\ P_9\ P_{10}\ P_{11}\ P_{12}\}$
$C_5 = \{P_{11}\ P_{12}\ P_{13}\ P_{14}\ P_{15}\}$ and $C_6 = \{P_{13}\ P_{14}\ P_{15}\ P_1\ P_2\}$.

Using these 6 overlapping classes we give the directed graphs and their related connection matrices.

The directed graph related with $C_1 = \{P_1\ P_2\ P_3\ P_4\ P_5\}$ given by the expert is as follows:

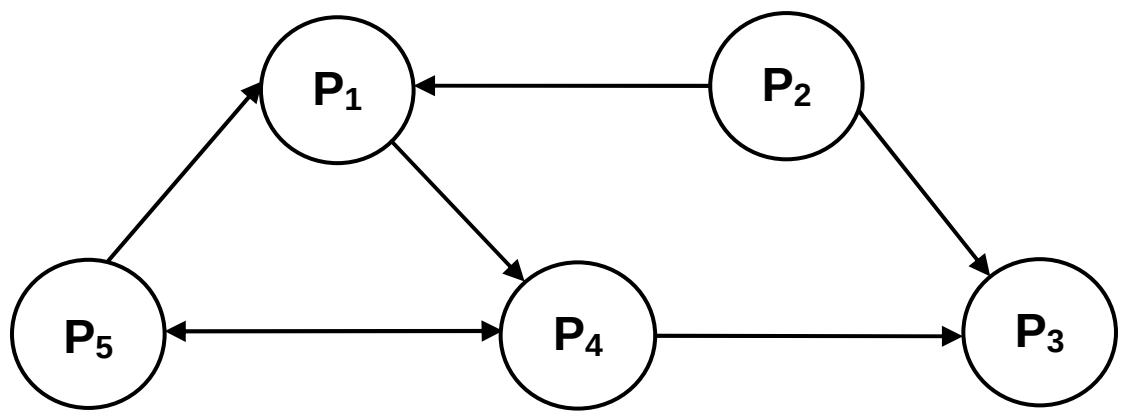

FIGURE 2.2.5.7

The related connection matrix is given below:

$$\begin{array}{c} \\ P_1 \\ P_2 \\ P_3 \\ P_4 \\ P_5 \end{array} \begin{array}{c} \begin{array}{ccccc} P_1 & P_2 & P_3 & P_4 & P_5 \end{array} \\ \left[\begin{array}{ccccc} 0 & 0 & 0 & 1 & 0 \\ 1 & 0 & 1 & 0 & 0 \\ 0 & 0 & 0 & 0 & 0 \\ 0 & 0 & 1 & 0 & 1 \\ 1 & 0 & 0 & 1 & 0 \end{array}\right] \end{array}.$$

Now we give the directed graph of the expert for the class

$$C_2 = \{P_3\ P_4\ P_5\ P_6\ P_7\}.$$

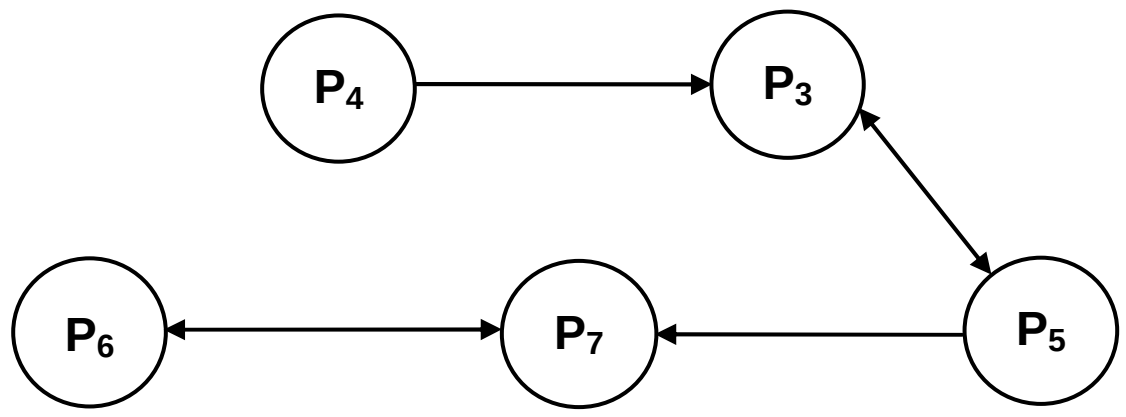

FIGURE 2.2.5.8

The related connection matrix is given below

$$\begin{array}{c} \\ P_3 \\ P_4 \\ P_5 \\ P_6 \\ P_7 \end{array}\begin{array}{c} \begin{array}{ccccc} P_3 & P_4 & P_5 & P_6 & P_7 \end{array} \\ \left[\begin{array}{ccccc} 0 & 0 & 1 & 0 & 0 \\ 1 & 0 & 0 & 0 & 0 \\ 1 & 0 & 0 & 0 & 1 \\ 0 & 0 & 0 & 0 & 1 \\ 0 & 0 & 0 & 1 & 0 \end{array}\right] \end{array}.$$

The directed graph related with the class $C_3 = \{P_6, P_7, P_8, P_9, P_{10}\}$ given by the expert

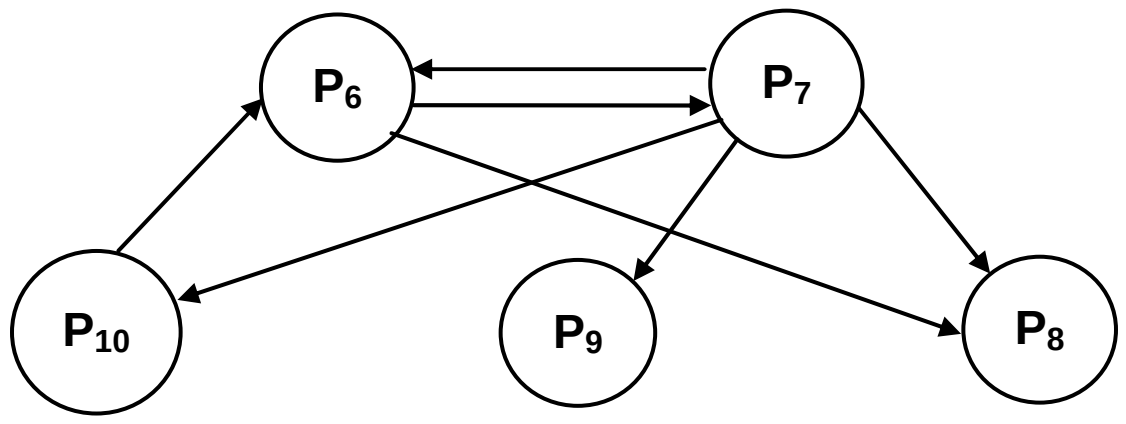

FIGURE 2.2.5.9

The connection matrix obtained from the directed graph is as follows:

$$\begin{array}{c} \\ P_6 \\ P_7 \\ P_8 \\ P_9 \\ P_{10} \end{array}\begin{array}{c} \begin{array}{ccccc} P_6 & P_7 & P_8 & P_9 & P_{10} \end{array} \\ \left[\begin{array}{ccccc} 0 & 1 & 1 & 0 & 0 \\ 1 & 0 & 1 & 1 & 1 \\ 0 & 0 & 0 & 0 & 0 \\ 0 & 0 & 0 & 0 & 0 \\ 1 & 0 & 0 & 0 & 0 \end{array}\right] \end{array}.$$

Now consider the class $C_4 = \{P_8, P_9, P_{10}\ P_{11}\ P_{12}\}$ the directed graph given by this expert is given below

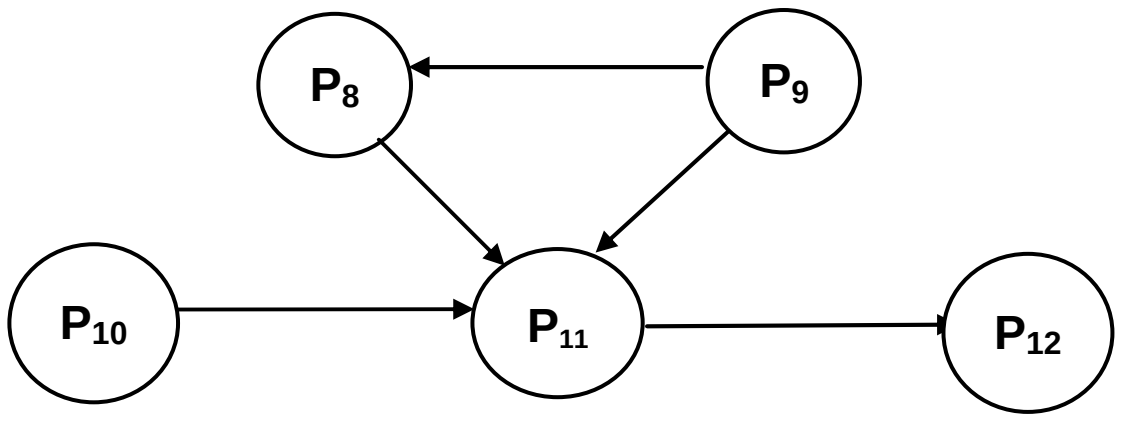

FIGURE 2.2.5.10

The related connection matrix is as follows

$$\begin{array}{c} \\ P_8 \\ P_9 \\ P_{10} \\ P_{11} \\ P_{12} \end{array}\begin{array}{c} \begin{array}{ccccc} P_8 & P_9 & P_{10} & P_{11} & P_{12} \end{array} \\ \left[\begin{array}{ccccc} 0 & 0 & 0 & 1 & 0 \\ 1 & 0 & 0 & 1 & 0 \\ 0 & 0 & 0 & 1 & 0 \\ 0 & 0 & 0 & 0 & 1 \\ 0 & 0 & 0 & 1 & 0 \end{array}\right] \end{array}.$$

The directed graph given by the expert related with the class $C_5 = \{P_{11}\ P_{12}\ P_{13}\ P_{14}\ P_{15}\}$ is given below:

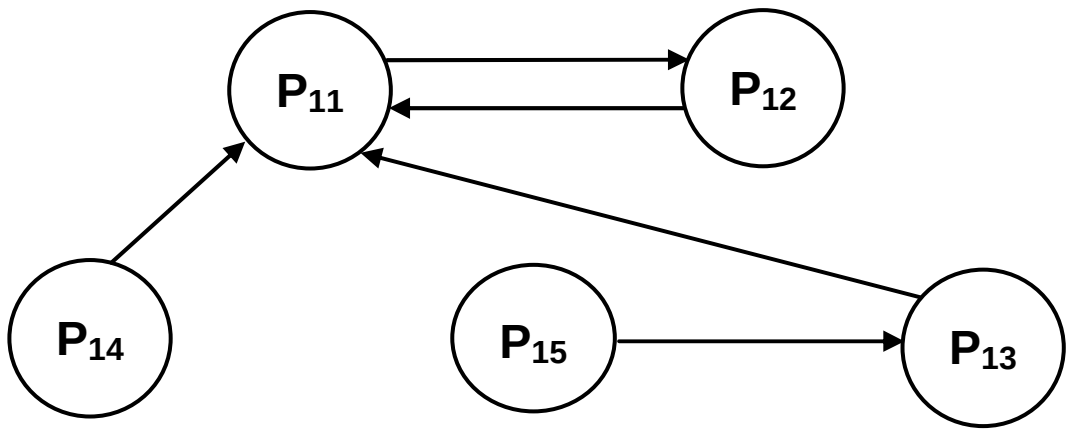

FIGURE 2.2.5.11

The related connection matrix is as follows:

$$\begin{array}{c} \\ P_{11} \\ P_{12} \\ P_{13} \\ P_{14} \\ P_{15} \end{array} \begin{array}{c} \begin{array}{ccccc} P_{11} & P_{12} & P_{13} & P_{14} & P_{15} \end{array} \\ \begin{bmatrix} 0 & 1 & 0 & 0 & 0 \\ 1 & 0 & 0 & 0 & 0 \\ 1 & 0 & 0 & 0 & 0 \\ 1 & 0 & 0 & 0 & 0 \\ 0 & 0 & 1 & 0 & 0 \end{bmatrix} \end{array}.$$

The directed graph given by the expert for the last class of attributes $C_6 = \{P_{13}\ P_{14}\ P_{15}\ P_1\ P_2\}$ is as follows

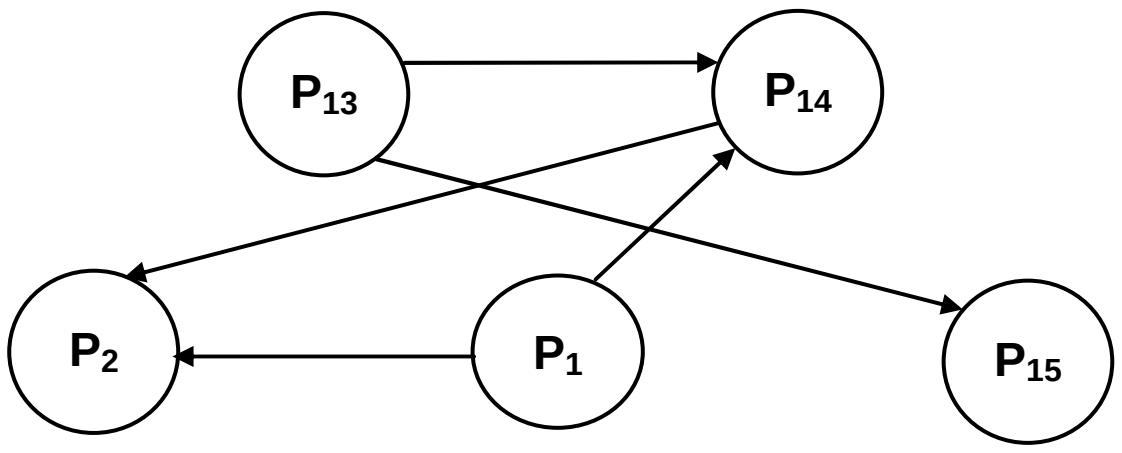

FIGURE 2.2.5.12

The related collection matrix

$$\begin{array}{c} \\ P_{13} \\ P_{14} \\ P_{15} \\ P_{1} \\ P_{2} \end{array} \begin{array}{c} \begin{array}{ccccc} P_{13} & P_{14} & P_{15} & P_{1} & P_{2} \end{array} \\ \begin{bmatrix} 0 & 1 & 1 & 0 & 0 \\ 0 & 0 & 0 & 0 & 1 \\ 0 & 0 & 0 & 0 & 0 \\ 0 & 1 & 0 & 0 & 1 \\ 0 & 0 & 0 & 0 & 0 \end{bmatrix} \end{array}.$$

Using the six classes of over lapping attributes we formulate the combined block overlap connection matrix V.

$$
\begin{array}{c}
 \\ P_1 \\ P_2 \\ P_3 \\ P_4 \\ P_5 \\ P_6 \\ P_7 \\ P_8 \\ P_9 \\ P_{10} \\ P_{11} \\ P_{12} \\ P_{13} \\ P_{14} \\ P_{15}
\end{array}
\begin{array}{c}
\begin{array}{ccccccccccccccc}
P_1 & P_2 & P_3 & P_4 & P_5 & P_6 & P_7 & P_8 & P_9 & P_{10} & P_{11} & P_{12} & P_{13} & P_{14} & P_{15}
\end{array} \\
\left[\begin{array}{ccccccccccccccc}
0 & 1 & 0 & 1 & 0 & 0 & 0 & 0 & 0 & 0 & 0 & 0 & 0 & 1 & 0 \\
1 & 0 & 1 & 0 & 0 & 0 & 0 & 0 & 0 & 0 & 0 & 0 & 0 & 0 & 0 \\
0 & 0 & 0 & 0 & 1 & 0 & 0 & 0 & 0 & 0 & 0 & 0 & 0 & 0 & 0 \\
0 & 0 & 2 & 0 & 1 & 0 & 0 & 0 & 0 & 0 & 0 & 0 & 0 & 0 & 0 \\
1 & 0 & 1 & 1 & 0 & 0 & 1 & 0 & 0 & 0 & 0 & 0 & 0 & 0 & 0 \\
0 & 0 & 1 & 0 & 0 & 0 & 2 & 1 & 0 & 0 & 0 & 0 & 0 & 0 & 0 \\
0 & 0 & 0 & 0 & 0 & 2 & 0 & 1 & 1 & 1 & 0 & 0 & 0 & 0 & 0 \\
0 & 0 & 0 & 0 & 0 & 0 & 0 & 0 & 0 & 0 & 1 & 0 & 0 & 0 & 0 \\
0 & 0 & 0 & 0 & 0 & 0 & 0 & 1 & 0 & 0 & 1 & 0 & 0 & 0 & 0 \\
0 & 0 & 0 & 0 & 0 & 1 & 0 & 0 & 0 & 0 & 1 & 0 & 0 & 0 & 0 \\
0 & 0 & 0 & 0 & 0 & 0 & 0 & 0 & 0 & 0 & 0 & 2 & 0 & 0 & 0 \\
0 & 0 & 0 & 0 & 0 & 0 & 0 & 0 & 0 & 0 & 2 & 0 & 0 & 0 & 0 \\
0 & 0 & 0 & 0 & 0 & 0 & 0 & 0 & 0 & 0 & 1 & 0 & 0 & 1 & 1 \\
0 & 1 & 0 & 0 & 0 & 0 & 0 & 0 & 0 & 0 & 1 & 0 & 0 & 0 & 0 \\
0 & 0 & 0 & 0 & 0 & 0 & 0 & 0 & 0 & 0 & 0 & 0 & 1 & 0 & 0
\end{array}\right]
\end{array}.
$$

Now V gives the 15 × 15 matrix of the dynamical system related with the combined block overlap fuzzy cognitive map. Suppose we are interested in studying the effect of state vectors.

Let us consider the state vector Y = (1 0 0 0 1 0 0 0 0 0 0 0 0 0 0) where only the two attributes, No binding with the family and the socially irresponsible co ordinates are in the on state and all other attributes are in the off state. To study the effect of Y on the dynamical system V;

$$
\begin{array}{llcl}
YV & \hookrightarrow (1\ 1\ 0\ 1\ 1\ 0\ 1\ 0\ 0\ 0\ 0\ 0\ 0\ 1\ 0) & = & Y_1 \text{ (say)} \\
Y_1V & \hookrightarrow (1\ 1\ 1\ 1\ 1\ 1\ 1\ 1\ 1\ 1\ 1\ 0\ 0\ 1\ 0) & = & Y_2 \text{ (say)} \\
Y_2V & \hookrightarrow (1\ 1\ 1\ 1\ 1\ 1\ 1\ 1\ 1\ 1\ 1\ 1\ 0\ 1\ 0) & = & Y_3 \text{ (say)} \\
Y_3\,V & \hookrightarrow (1\ 1\ 1\ 1\ 1\ 1\ 1\ 1\ 1\ 1\ 1\ 1\ 0\ 1\ 0) & = & Y_4 = Y_3
\end{array}
$$

(The hidden pattern of the system). Thus the hidden pattern of the system is a fixed point. When the migrant labourer has no binding with the family and is socially irresponsible we see he has all the other nodes to be working on his nature i.e., they all come to on state only the nodes unreachable by friends or relatives and failure of agriculture is in the off state.

Now we give the final model of this chapter namely blocks of varying size with varying overlap.

***Example 2.2.5.3:*** Consider the 12 attribute model studied relative to the problems of HIV/AIDS affected migrant labourers in page 178-9. Consider the 12 attributes $A_1$, $A_2$, …, $A_{12}$. Divide them into overlapping blocks as follows:

$C_1 = \{A_1\ A_2\ A_3\ A_4\ A_5\}$, $C_2 = \{A_4\ A_5\ A_6\ A_7\}$,
$C_3 = \{A_6\ A_7\ A_8\ A_9\}$, $C_4 = \{A_8\ A_9\ A_{10}\}$ and
$C_5 = \{A_{10}\ A_{11}\ A_{12}\ A_1\ A_2\ A_3\}$.

Now we analyze each of these classes of attributes $C_1$, $C_2$, $C_3$, $C_4$ and $C_5$.

The directed graph given by the expert for the attributes $\{A_1\ A_2\ A_3\ A_4\ A_5\}$.

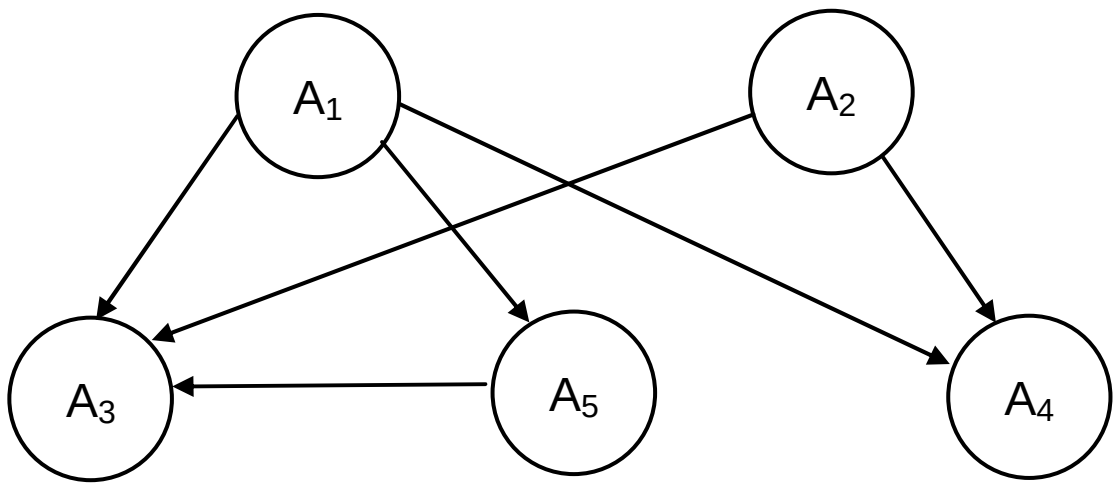

FIGURE 2.2.5.13

The related connection matrix.

$$\begin{array}{c} \\ A_1 \\ A_2 \\ A_3 \\ A_4 \\ A_5 \end{array} \begin{array}{c} \begin{array}{ccccc} A_1 & A_2 & A_3 & A_4 & A_5 \end{array} \\ \left[\begin{array}{ccccc} 0 & 0 & 1 & 1 & 1 \\ 0 & 0 & 1 & 1 & 0 \\ 0 & 0 & 0 & 0 & 0 \\ 0 & 0 & 0 & 0 & 0 \\ 0 & 0 & 1 & 0 & 0 \end{array}\right] \end{array}.$$

The directed graph related to the attributes {$A_4$, $A_5$, $A_6$, $A_7$} given by an expert.

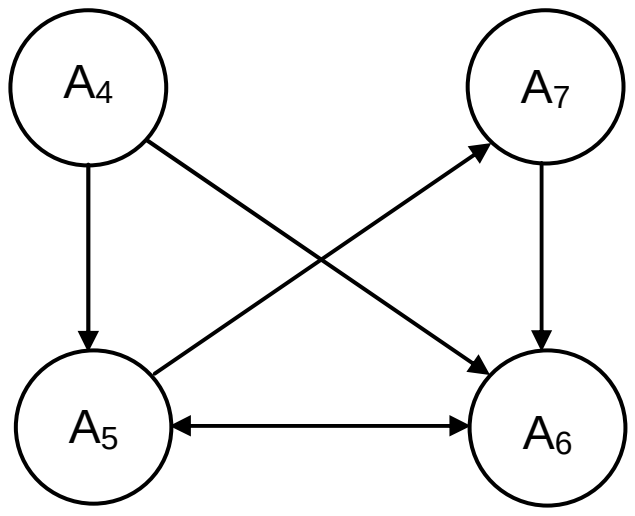

FIGURE 2.2.5.14

The related connection matrix is as follows.

$$\begin{array}{c} \\ A_4 \\ A_5 \\ A_6 \\ A_7 \end{array} \begin{array}{c} \begin{array}{cccc} A_4 & A_5 & A_6 & A_7 \end{array} \\ \begin{bmatrix} 0 & 1 & 1 & 0 \\ 0 & 0 & 1 & 1 \\ 0 & 1 & 0 & 0 \\ 0 & 0 & 1 & 0 \end{bmatrix} \end{array}.$$

The directed graph given by the expert related to the attributes $A_6$, $A_7$, $A_8$ and $A_9$.

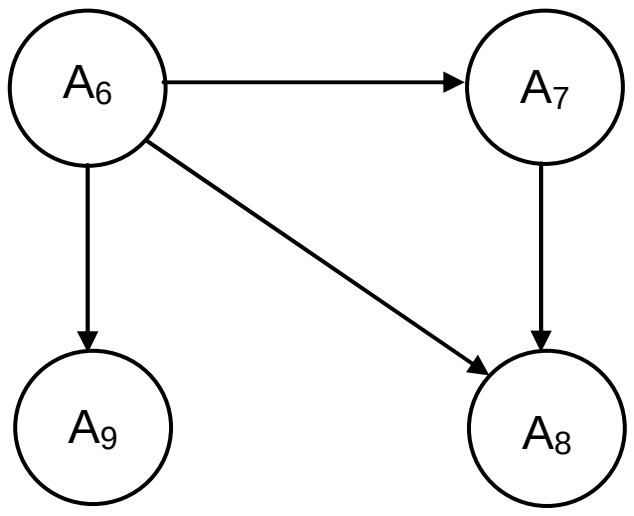

FIGURE 2.2.5.15

The related connection matrix is given in the next page:.

$$\begin{array}{c} \\ A_6 \\ A_7 \\ A_8 \\ A_9 \end{array}\begin{array}{c} \begin{array}{cccc} A_6 & A_7 & A_8 & A_9 \end{array} \\ \begin{bmatrix} 0 & 1 & 1 & 1 \\ 0 & 0 & 1 & 0 \\ 0 & 0 & 0 & 0 \\ 0 & 0 & 0 & 0 \end{bmatrix} \end{array}.$$

Now we give the directed graph given by an expert using the attributes $A_8$, $A_9$ and $A_{10}$.

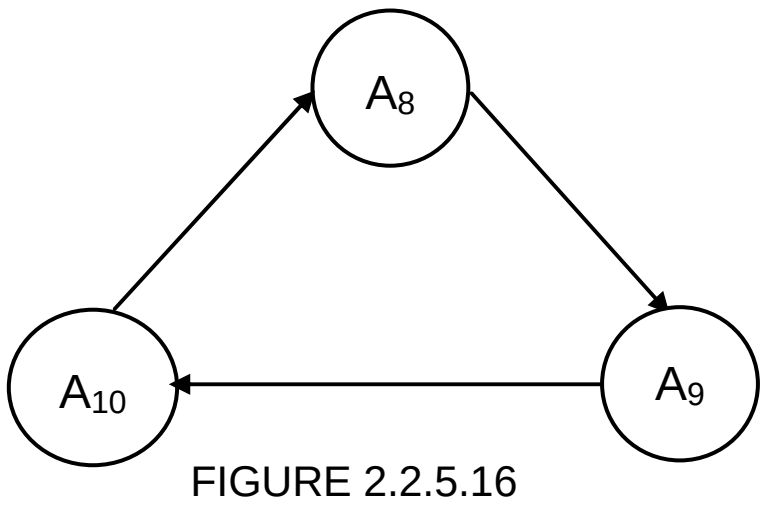

FIGURE 2.2.5.16

The connection matrix related with this directed graph.

$$\begin{array}{c} \\ A_8 \\ A_9 \\ A_{10} \end{array}\begin{array}{c} \begin{array}{ccc} A_8 & A_9 & A_{10} \end{array} \\ \begin{bmatrix} 0 & 1 & 0 \\ 0 & 0 & 1 \\ 1 & 0 & 0 \end{bmatrix} \end{array}.$$

Finally we obtain the connection matrix given by the expert related with the attributes $A_{10}$, $A_{11}$, $A_{12}$, $A_1$, $A_2$ and $A_3$.

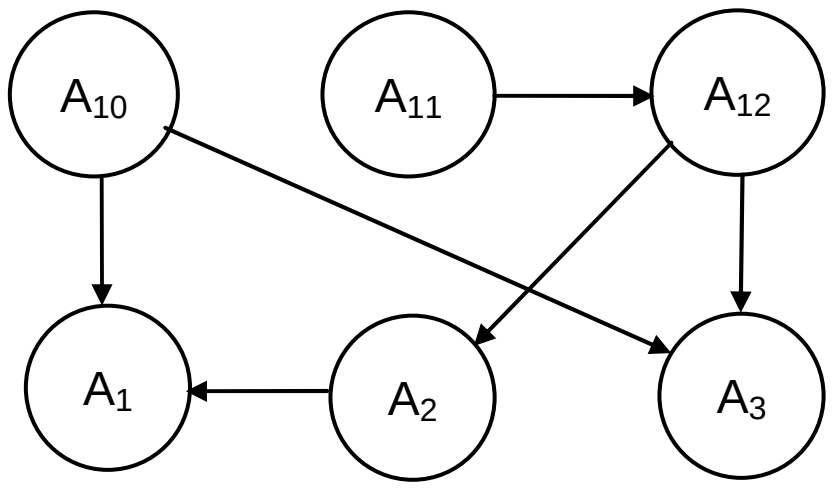

FIGURE 2.2.5.17

The related connection matrix.

$$
\begin{array}{c}
 \\ A_{10} \\ A_{11} \\ A_{12} \\ A_{1} \\ A_{2} \\ A_{3}
\end{array}
\begin{array}{c}
\begin{array}{cccccc} A_{10} & A_{11} & A_{12} & A_{1} & A_{2} & A_{3} \end{array} \\
\left[\begin{array}{cccccc}
0 & 0 & 0 & 1 & 0 & 1 \\
0 & 0 & 1 & 0 & 0 & 0 \\
0 & 0 & 0 & 0 & 1 & 1 \\
0 & 0 & 0 & 0 & 0 & 0 \\
1 & 0 & 0 & 0 & 0 & 0 \\
0 & 0 & 0 & 0 & 0 & 0
\end{array}\right]
\end{array}.
$$

Now we give the related connection combined block over lapping matrix N.

$$
\begin{array}{c}
 \\ A_{1} \\ A_{2} \\ A_{3} \\ A_{4} \\ A_{5} \\ A_{6} \\ A_{7} \\ A_{8} \\ A_{9} \\ A_{10} \\ A_{11} \\ A_{12}
\end{array}
\begin{array}{c}
\begin{array}{cccccccccccc} A_{1} & A_{2} & A_{3} & A_{4} & A_{5} & A_{6} & A_{7} & A_{8} & A_{9} & A_{10} & A_{11} & A_{12} \end{array} \\
\left[\begin{array}{cccccccccccc}
0 & 0 & 1 & 1 & 1 & 0 & 0 & 0 & 0 & 0 & 0 & 0 \\
0 & 0 & 1 & 1 & 0 & 0 & 0 & 0 & 0 & 1 & 0 & 0 \\
0 & 0 & 0 & 0 & 0 & 0 & 0 & 0 & 0 & 0 & 0 & 0 \\
0 & 0 & 0 & 0 & 1 & 1 & 0 & 0 & 0 & 0 & 0 & 0 \\
0 & 0 & 1 & 0 & 0 & 1 & 1 & 0 & 0 & 0 & 0 & 0 \\
0 & 0 & 0 & 0 & 1 & 0 & 1 & 1 & 1 & 0 & 0 & 0 \\
0 & 0 & 0 & 0 & 0 & 1 & 0 & 1 & 0 & 0 & 0 & 0 \\
0 & 0 & 0 & 0 & 0 & 0 & 0 & 0 & 1 & 0 & 0 & 0 \\
0 & 0 & 0 & 0 & 0 & 0 & 0 & 0 & 0 & 1 & 0 & 0 \\
1 & 0 & 1 & 0 & 0 & 0 & 0 & 1 & 0 & 0 & 0 & 0 \\
0 & 0 & 0 & 0 & 0 & 0 & 0 & 0 & 0 & 0 & 0 & 1 \\
0 & 1 & 1 & 0 & 0 & 0 & 0 & 0 & 0 & 0 & 0 & 0
\end{array}\right]
\end{array}.
$$

Now we can analyze the effect of any state vector on the dynamical system N.

Let

$$X = (0\ 0\ 0\ 1\ 0\ 0\ 0\ 0\ 0\ 0\ 0\ 0)$$

be the state vector with the only attribute profession in the on state the effect of X on N is given by

| | | | | |
|---|---|---|---|---|
| XN | $\hookrightarrow$ | (0 0 0 1 1 1 0 0 0 0 0 0) | = | $X_1$ (say) |
| $X_1N$ | $\hookrightarrow$ | (0 0 1 1 1 1 1 1 1 0 0 0) | = | $X_2$ (say) |
| $X_2N$ | $\hookrightarrow$ | (0 0 1 1 1 1 1 1 1 1 0 0) | = | $X_3$ (say) |
| $X_3N$ | $\hookrightarrow$ | (1 0 1 1 1 1 1 1 1 1 0 0) | = | $X_4$ (say) |
| $X_4N$ | $\hookrightarrow$ | (1 0 1 1 1 1 1 1 1 1 0 0) | = | $X_5$ (= $X_4$). |

The hidden pattern is a fixed point. Only the vectors $A_2$, $A_{11}$ and $A_{12}$ are in the off state and all other attributes come to the on state. We see there is difference between resultant vectors when using the dynamical system A and N.

Now we consider the state vector $R_1$ = (1 0 0 0 0 0 1 0 0 1 0 0) in which the attributes easy money, no social responsibility and more leisure i.e., $A_1$, $A_7$ and $A_{10}$ are in the on state and all other attributes are in the off state. The effect of $R_1$ on the dynamical system N is given by

| | | | | |
|---|---|---|---|---|
| $R_1N$ | $\hookrightarrow$ | (1 0 1 1 1 1 1 1 0 0 0 0) | = | $R_2$ (say) |
| $R_2N$ | $\hookrightarrow$ | (1 0 1 1 1 1 1 1 1 0 0 0) | = | $R_3$ (say) |
| $R_3N$ | $\hookrightarrow$ | (1 0 1 1 1 1 1 1 1 1 0 0) | = | $R_4$ (say) |
| $R_4N$ | $\hookrightarrow$ | (1 0 1 1 1 1 1 1 1 1 0 0) | = | $R_5 = R_4$ |

$R_4$ is a fixed point. The hidden pattern is a fixed point only the attributes $A_2$, $A_{11}$ and $A_{12}$ remain in the off state all other vectors became on.

***Example 2.2.5.4:*** Now we use the attributes $P_1, P_2, \ldots, P_{15}$ given in page 182-3 and model it using the combined block overlap FCM.

First we divide $P_1, P_2, \ldots, P_{15}$ into 5 classes

$C_1 = \{P_1\ P_2\ P_3\ P_4\ P_5\ P_6\}$,
$C_2 = \{P_4\ P_5\ P_6\ P_7\ P_8\ P_9\ P_{10}\}$,
$C_3 = \{P_7\ P_8\ P_9\ P_{10}\ P_{11}\}$,
$C_4 = \{P_{11}, P_{12}, P_{13}\}$
and $C_5 = \{P_{12}\ P_{13}\ P_{14}\ P_{15}\ P_1\ P_2\ P_3\}$.

The expert's directed graph using the attributes $P_1$, $P_2$, …, $P_6$.

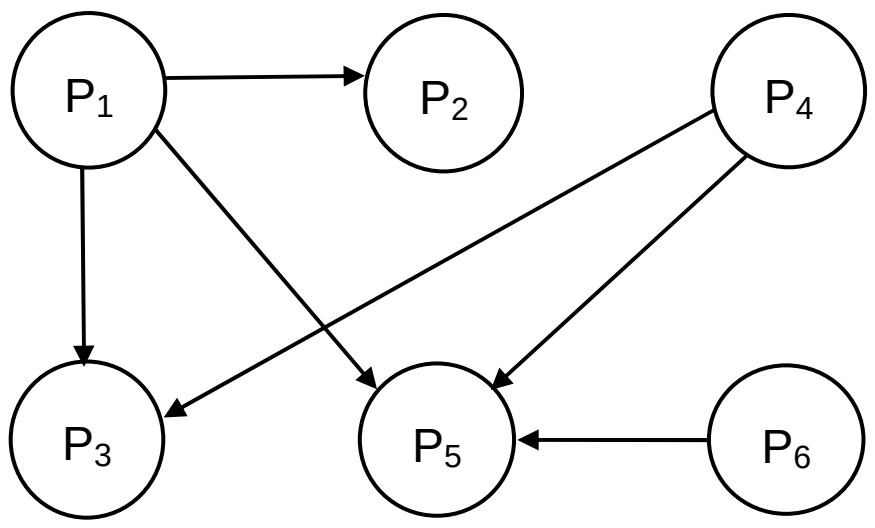

FIGURE 2.2.5.18

The resulting connection matrix using the directed graph is here.

$$
\begin{array}{c}
 \\ P_1 \\ P_2 \\ P_3 \\ P_4 \\ P_5 \\ P_6
\end{array}
\begin{array}{c}
\begin{matrix} P_1 & P_2 & P_3 & P_4 & P_5 & P_6 \end{matrix} \\
\begin{bmatrix}
0 & 1 & 1 & 0 & 1 & 0 \\
0 & 0 & 0 & 0 & 0 & 0 \\
0 & 0 & 0 & 0 & 0 & 0 \\
0 & 0 & 1 & 0 & 1 & 0 \\
0 & 0 & 0 & 0 & 0 & 0 \\
0 & 0 & 0 & 0 & 1 & 0
\end{bmatrix}
\end{array}.
$$

For the attributes $P_4$, $P_5$,…, $P_9$, $P_{10}$ using the experts opinion we have the following directed graph.

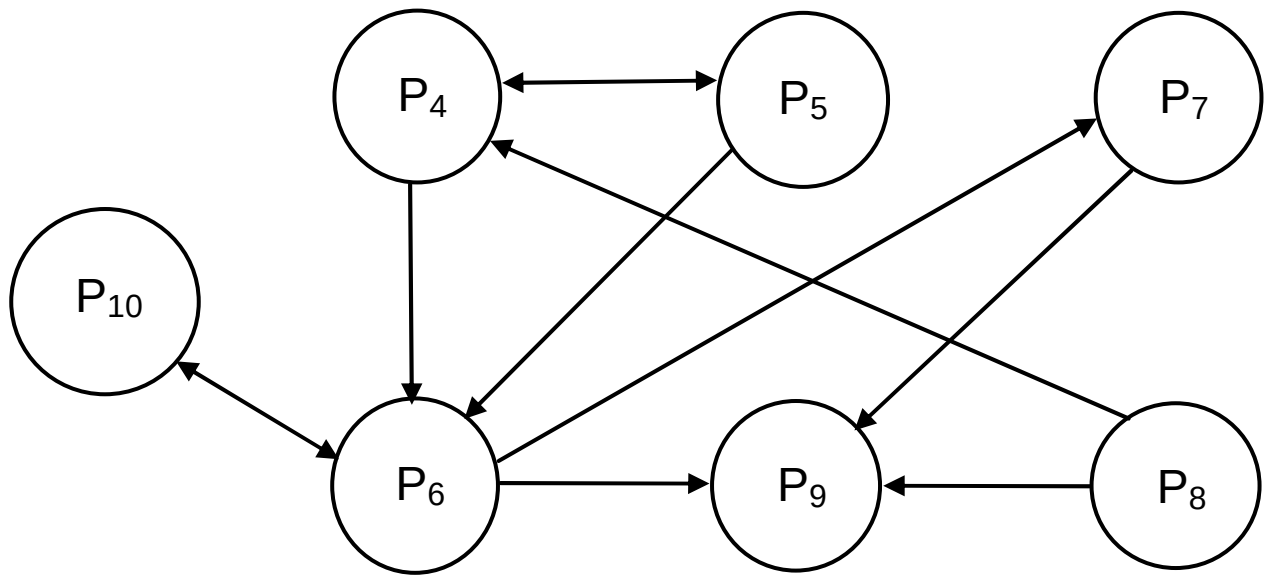

FIGURE 2.2.5.19

The related connection matrix is as follows:

$$
\begin{array}{c}
 \\ P_4 \\ P_5 \\ P_6 \\ P_7 \\ P_8 \\ P_9 \\ P_{10}
\end{array}
\begin{array}{c}
\begin{array}{ccccccc} P_4 & P_5 & P_6 & P_7 & P_8 & P_9 & P_{10} \end{array} \\
\begin{bmatrix}
0 & 1 & 1 & 0 & 0 & 0 & 0 \\
1 & 0 & 1 & 0 & 0 & 0 & 0 \\
0 & 0 & 0 & 1 & 0 & 1 & 1 \\
0 & 0 & 1 & 0 & 0 & 0 & 0 \\
1 & 0 & 0 & 0 & 0 & 1 & 0 \\
0 & 0 & 0 & 0 & 0 & 0 & 0 \\
0 & 0 & 1 & 0 & 0 & 0 & 0
\end{bmatrix}
\end{array}.
$$

The directed graph and the related connection matrix given by the expert for the next class is given below:

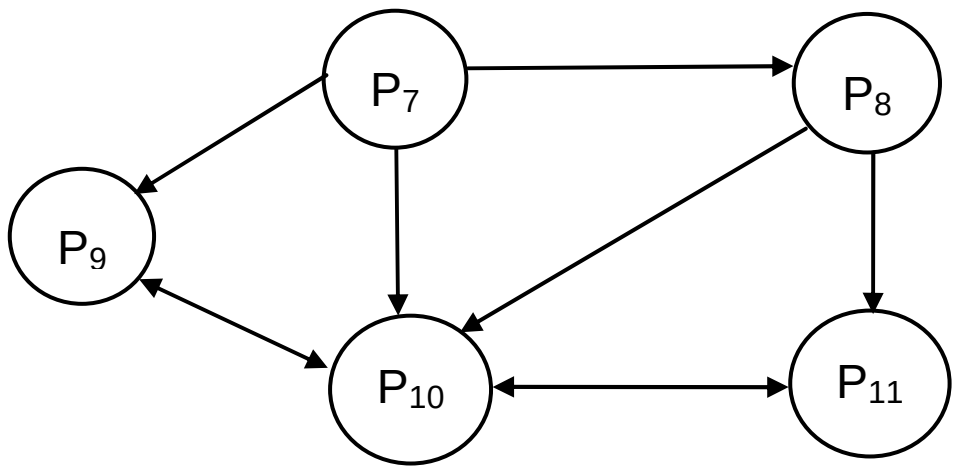

FIGURE 2.2.5.20

$$
\begin{array}{c}
 \\ P_7 \\ P_8 \\ P_9 \\ P_{10} \\ P_{11}
\end{array}
\begin{array}{c}
\begin{array}{ccccc} P_7 & P_8 & P_9 & P_{10} & P_{11} \end{array} \\
\begin{bmatrix}
0 & 1 & 1 & 1 & 0 \\
0 & 0 & 0 & 1 & 1 \\
0 & 0 & 0 & 1 & 0 \\
0 & 0 & 1 & 0 & 1 \\
0 & 0 & 0 & 1 & 0
\end{bmatrix}
\end{array}.
$$

Now the directed graph for the three attributes $P_{11}$ $P_{12}$ and $P_{13}$.

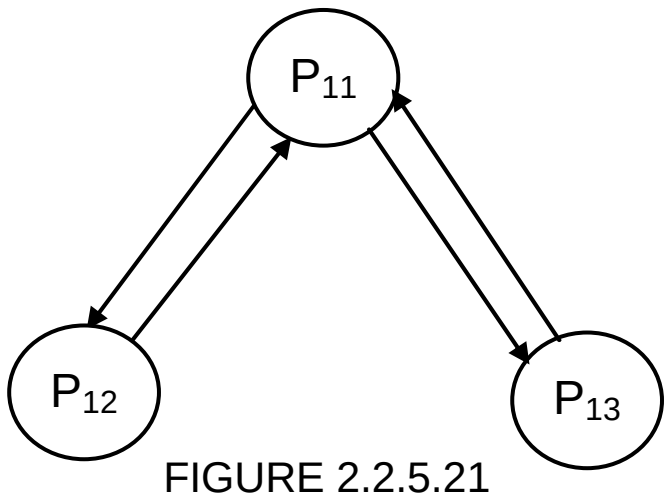

FIGURE 2.2.5.21

The related connection matrix

$$\begin{array}{c} \\ P_{11} \\ P_{12} \\ P_{13} \end{array} \begin{array}{c} \begin{array}{ccc} P_{11} & P_{12} & P_{13} \end{array} \\ \begin{bmatrix} 0 & 1 & 1 \\ 1 & 0 & 0 \\ 1 & 0 & 0 \end{bmatrix} \end{array}.$$

We give the expert's directed graph using $P_{12}$ $P_{13}$ $P_{14}$ $P_{15}$ $P_1$ $P_2$ $P_3$.

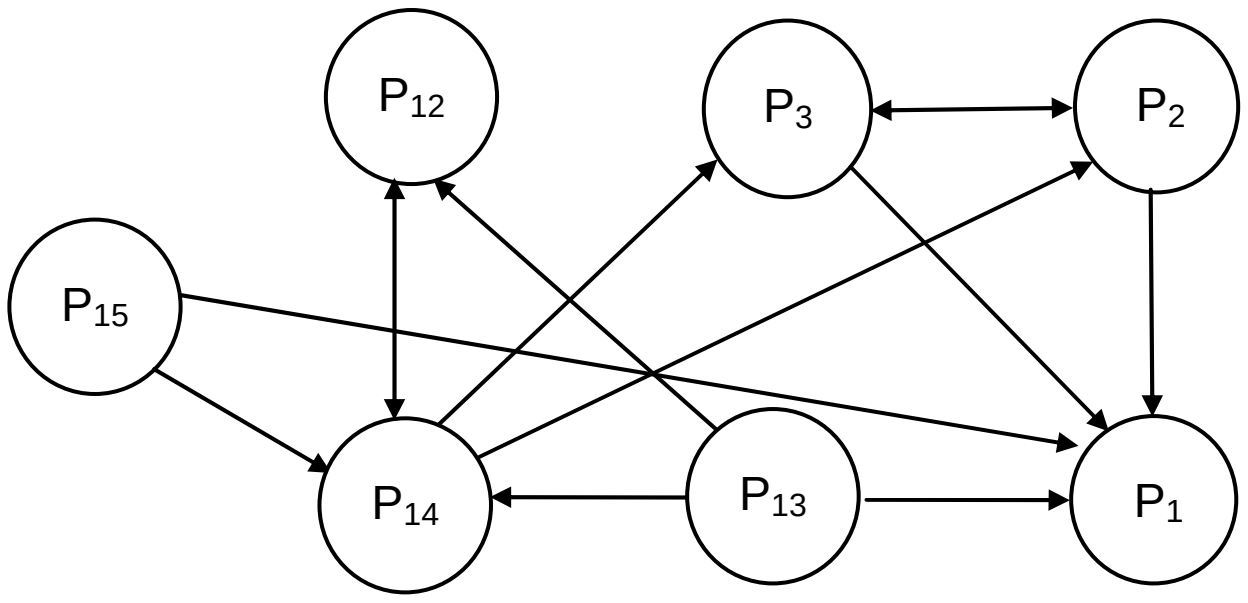

FIGURE 2.2.5.22

The related connection matrix is

$$\begin{array}{c} \\ P_{12} \\ P_{13} \\ P_{14} \\ P_{15} \\ P_1 \\ P_2 \\ P_3 \end{array} \begin{array}{c} \begin{array}{ccccccc} P_{12} & P_{13} & P_{14} & P_{15} & P_1 & P_2 & P_3 \end{array} \\ \begin{bmatrix} 0 & 0 & 1 & 0 & 0 & 0 & 0 \\ 1 & 0 & 1 & 0 & 1 & 0 & 0 \\ 1 & 0 & 0 & 0 & 0 & 0 & 1 \\ 0 & 0 & 0 & 0 & 1 & 0 & 0 \\ 0 & 0 & 0 & 0 & 0 & 0 & 0 \\ 0 & 0 & 0 & 0 & 1 & 0 & 1 \\ 0 & 0 & 0 & 0 & 1 & 1 & 0 \end{bmatrix} \end{array}.$$

Let O denote the connection matrix of the combined block overlap FCM.

$$
\begin{array}{c}
 \\ P_1 \\ P_2 \\ P_3 \\ P_4 \\ P_5 \\ P_6 \\ P_7 \\ P_8 \\ P_9 \\ P_{10} \\ P_{11} \\ P_{12} \\ P_{13} \\ P_{14} \\ P_{15}
\end{array}
\begin{array}{c}
\begin{array}{ccccccccccccccc} P_1 & P_2 & P_3 & P_4 & P_5 & P_6 & P_7 & P_8 & P_9 & P_{10} & P_{11} & P_{12} & P_{13} & P_{14} & P_{15} \end{array} \\
\left[\begin{array}{ccccccccccccccc}
0 & 1 & 1 & 0 & 1 & 0 & 0 & 0 & 0 & 0 & 0 & 0 & 0 & 0 & 0 \\
1 & 0 & 1 & 0 & 0 & 0 & 0 & 0 & 0 & 0 & 0 & 0 & 0 & 0 & 0 \\
1 & 1 & 0 & 0 & 0 & 0 & 0 & 0 & 0 & 0 & 0 & 0 & 0 & 0 & 0 \\
0 & 0 & 1 & 0 & 2 & 1 & 0 & 0 & 0 & 0 & 0 & 0 & 0 & 0 & 0 \\
0 & 0 & 0 & 1 & 0 & 1 & 0 & 0 & 0 & 0 & 0 & 0 & 0 & 0 & 0 \\
0 & 0 & 0 & 0 & 1 & 0 & 1 & 0 & 1 & 1 & 0 & 0 & 0 & 0 & 0 \\
0 & 0 & 0 & 0 & 0 & 1 & 0 & 1 & 1 & 1 & 0 & 0 & 0 & 0 & 0 \\
0 & 0 & 0 & 1 & 0 & 0 & 0 & 0 & 1 & 1 & 1 & 0 & 0 & 0 & 0 \\
0 & 0 & 0 & 0 & 0 & 0 & 0 & 0 & 0 & 1 & 0 & 0 & 0 & 0 & 0 \\
0 & 0 & 0 & 0 & 0 & 1 & 0 & 0 & 1 & 0 & 1 & 0 & 0 & 0 & 0 \\
0 & 0 & 0 & 0 & 0 & 0 & 0 & 0 & 0 & 1 & 0 & 1 & 1 & 0 & 0 \\
0 & 0 & 0 & 0 & 0 & 0 & 0 & 0 & 0 & 0 & 1 & 0 & 0 & 1 & 0 \\
1 & 0 & 0 & 0 & 0 & 0 & 0 & 0 & 0 & 0 & 1 & 1 & 0 & 1 & 0 \\
0 & 0 & 1 & 0 & 0 & 0 & 0 & 0 & 0 & 0 & 0 & 1 & 0 & 0 & 0 \\
1 & 0 & 0 & 0 & 0 & 0 & 0 & 0 & 0 & 0 & 0 & 0 & 0 & 0 & 0
\end{array}\right]
\end{array}.
$$

Now consider a state vector S = (0 1 0 0 0 0 0 0 0 0 0 0 0 0 0 ) the notion male ego alone is in the on state and all other attributes are in the off state.

The effect of S on the dynamical system O is given by

SO $\hookrightarrow$ (1 1 1 0 0 0 0 0 0 0 0 0 0 0 0) = $S_1$ (say)
$S_1$O $\hookrightarrow$ (1 1 1 0 1 0 0 0 0 0 0 0 0 0 0) = $S_2$ (say)
$S_2$O $\hookrightarrow$ (1 1 1 1 1 1 0 0 0 0 0 0 0 0 0) = $S_3$ (say)
$S_3$O $\hookrightarrow$ (1 1 1 1 1 1 1 0 1 1 0 0 0 0 0) = $S_4$
$S_4$O $\hookrightarrow$ (1 1 1 1 1 1 1 1 1 1 1 0 0 0 0) = $S_5$ (say)
$S_5$O $\hookrightarrow$ (1 1 1 1 1 1 1 1 1 1 1 1 1 0 0) = $S_6$ (say)
$S_6$O $\hookrightarrow$ (1 1 1 1 1 1 1 1 1 1 1 1 1 1 0) = $S_7$ ($S_6$).
$S_7$O $\hookrightarrow$ $S_8 = S_7$.

Thus the hidden pattern is a fixed point. All vector come to on state except the notion $P_{15}$. Using the C-program the on state of several vectors have been found out. They are used in deriving our conclusions.

## 2.3 Definition and Illustration of Fuzzy Relational Maps (FRMs)

This section has two subsections, subsection one recalls the notion of fuzzy relational maps. Subsection two gives illustrations of fuzzy relational maps.

### 2.3.1 Definition of Fuzzy Relational Maps

In this section, we introduce the notion of Fuzzy Relational Maps (FRMs); they are constructed analogous to FCMs described and discussed in the earlier sections. In FCMs we promote the correlations between causal associations among concurrently active units. But in FRMs we divide the very causal associations into two disjoint units, for example, the relation between a teacher and a student or relation between an employee or employer or a relation between doctor and patient and so on. Thus for us to define a FRM we need a domain space and a range space which are disjoint in the sense of concepts. We further assume no intermediate relation exists within the domain elements or node and the range spaces elements. The number of elements in the range space need not in general be equal to the number of elements in the domain space.

Thus throughout this section we assume the elements of the domain space are taken from the real vector space of dimension n and that of the range space are real vectors from the vector space of dimension m (m in general need not be equal to n). We denote by R the set of nodes $R_1,\ldots, R_m$ of the range space, where $R = \{(x_1,\ldots, x_m) \mid x_j = 0 \text{ or } 1\}$ for j = 1, 2, …, m. If $x_i = 1$ it means that the node $R_i$ is in the on state and if $x_i = 0$ it means that the node $R_i$ is in the off state. Similarly D denotes the nodes $D_1, D_2, \ldots, D_n$ of the domain space where $D = \{(x_1,\ldots, x_n) \mid x_j = 0 \text{ or } 1\}$ for i = 1, 2,…, n. If $x_i = 1$ it means that the node $D_i$ is in the on state and if $x_i = 0$ it means that the node $D_i$ is in the off state.

Now we proceed on to define a FRM.

**DEFINITION 2.3.1.1:** *A FRM is a directed graph or a map from D to R with concepts like policies or events etc, as nodes and causalities as edges. It represents causal relations between spaces D and R.*

*Let $D_i$ and $R_j$ denote that the two nodes of an FRM. The directed edge from $D_i$ to $R_j$ denotes the causality of $D_i$ on $R_j$ called relations. Every edge in the FRM is weighted with a number in the set $\{0, \pm 1\}$. Let $e_{ij}$ be the weight of the edge $D_iR_j$, $e_{ij} \in \{0, \pm 1\}$. The weight of the edge $D_i R_j$ is positive if increase in $D_i$ implies increase in $R_j$ or decrease in $D_i$ implies decrease in $R_j$ ie causality of $D_i$ on $R_j$ is 1. If $e_{ij} = 0$, then $D_i$ does not have any effect on $R_j$. We do not discuss the cases when increase in $D_i$ implies decrease in $R_j$ or decrease in $D_i$ implies increase in $R_j$.*

**DEFINITION 2.3.1.2:** *When the nodes of the FRM are fuzzy sets then they are called fuzzy nodes. FRMs with edge weights $\{0, \pm 1\}$ are called simple FRMs.*

**DEFINITION 2.3.1.3:** *Let $D_1, \ldots, D_n$ be the nodes of the domain space D of an FRM and $R_1, \ldots, R_m$ be the nodes of the range space R of an FRM. Let the matrix E be defined as $E = (e_{ij})$ where $e_{ij}$ is the weight of the directed edge $D_iR_j$ (or $R_jD_i$), E is called the relational matrix of the FRM.*

***Note:*** It is pertinent to mention here that unlike the FCMs the FRMs can be a rectangular matrix with rows corresponding to the domain space and columns corresponding to the range space. This is one of the marked difference between FRMs and FCMs.

**DEFINITION 2.3.1.4:** *Let $D_1, \ldots, D_n$ and $R_1, \ldots, R_m$ denote the nodes of the FRM. Let $A = (a_1, \ldots, a_n)$, $a_i \in \{0, 1\}$. A is called the instantaneous state vector of the domain space and it denotes the on-off position of the nodes at any instant. Similarly let $B = (b_1, \ldots, b_m)$ $b_i \in \{0, 1\}$. B is called instantaneous state vector of the range space and it denotes the on-off position of the nodes* *at any instant $a_i = 0$ if $a_i$ is off and $a_i = 1$ if $a_i$ is on for i= 1, 2,…, n Similarly, $b_i = 0$ if $b_i$ is off and $b_i = 1$ if $b_i$ is on, for i= 1, 2,…, m.*

**DEFINITION 2.3.1.5:** *Let $D_1$, …, $D_n$ and $R_1$,…, $R_m$ be the nodes of an FRM. Let $D_iR_j$ (or $R_j$ $D_i$) be the edges of an FRM, j = 1, 2,…, m and i= 1, 2,…, n. Let the edges form a directed cycle. An FRM is said to be a cycle if it posses a directed cycle. An FRM is said to be acyclic if it does not posses any directed cycle.*

**DEFINITION 2.3.1.6:** *An FRM with cycles is said to be an FRM with feedback.*

**DEFINITION 2.3.1.7:** *When there is a feedback in the FRM, i.e. when the causal relations flow through a cycle in a revolutionary manner, the FRM is called a dynamical system.*

**DEFINITION 2.3.1.8:** *Let $D_i$ $R_j$ (or $R_j$ $D_i$), $1 \leq j \leq m$, $1 \leq i \leq n$. When $R_i$ (or $D_j$) is switched on and if causality flows through edges of the cycle and if it again causes $R_i$ (or$D_j$), we say that the dynamical system goes round and round. This is true for any node $R_j$ (or $D_i$) for $1 \leq i \leq n$, (or $1 \leq j \leq m$). The equilibrium state of this dynamical system is called the hidden pattern.*

**DEFINITION 2.3.1.9:** *If the equilibrium state of a dynamical system is a unique state vector, then it is called a fixed point. Consider an FRM with $R_1$, $R_2$,…, $R_m$ and $D_1$, $D_2$,…, $D_n$ as nodes. For example, let us start the dynamical system by switching on $R_1$ (or $D_1$). Let us assume that the FRM settles down with $R_1$ and $R_m$ (or $D_1$ and $D_n$) on, i.e. the state vector remains as (1, 0, …, 0, 1) in R (or 1, 0, 0, … , 0, 1) in D), This state vector is called the fixed point.*

**DEFINITION 2.3.1.10:** *If the FRM settles down with a state vector repeating in the form*

$$A_1 \to A_2 \to A_3 \to \ldots \to A_i \to A_1 \text{ (or } B_1 \to B_2 \to \ldots \to B_i \to B_1\text{)}$$

*then this equilibrium is called a limit cycle.*

**METHODS OF DETERMINING THE HIDDEN PATTERN**

Let $R_1$, $R_2$, …, $R_m$ and $D_1$, $D_2$, …, $D_n$ be the nodes of a FRM with feedback. Let E be the relational matrix. Let us find a hidden pattern when $D_1$ is switched on i.e. when an input is given as vector $A_1$ = (1, 0, …, 0) in $D_1$, the data should pass through the relational matrix E. This is done by multiplying $A_1$ with the relational matrix E. Let $A_1E = (r_1, r_2, \ldots, r_m)$, after thresholding and updating the resultant vector we get $A_1 E \in R$. Now let $B = A_1E$ we pass on B into $E^T$ and obtain $BE^T$. We update and threshold the vector $BE^T$ so that $BE^T \in D$. This procedure is repeated till we get a limit cycle or a fixed point.

**DEFINITION 2.3.1.11:** *Finite number of FRMs can be combined together to produce the joint effect of all the FRMs. Let $E_1, \ldots, E_p$ be the relational matrices of the FRMs with nodes $R_1, R_2, \ldots, R_m$ and $D_1, D_2, \ldots, D_n$, then the combined FRM is represented by the relational matrix $E = E_1 + \ldots + E_p$.*

Now we give a simple illustration of a FRM, for more about FRMs please refer [237, 243-4].

***Example 2.3.1.1:*** Let us consider the relationship between the teacher and the student. Suppose we take the domain space as the concepts belonging to the teacher say $D_1, \ldots, D_5$ and the range space denote the concepts belonging to the student say $R_1$, $R_2$ and $R_3$.

We describe the nodes $D_1$, …, $D_5$ and $R_1$ , $R_2$ and $R_3$ as follows:

Domain Space

$D_1$ – Teaching is good
$D_2$ – Teaching is poor
$D_3$ – Teaching is mediocre
$D_4$ – Teacher is kind
$D_5$ – Teacher is harsh [or rude]

(We can have more concepts like teacher is non-reactive, unconcerned etc.)

Range Space

$R_1$ – Good Student
$R_2$ – Bad Student
$R_3$ – Average Student

The relational directed graph of the teacher-student model is given in Figure 2.3.1.1.

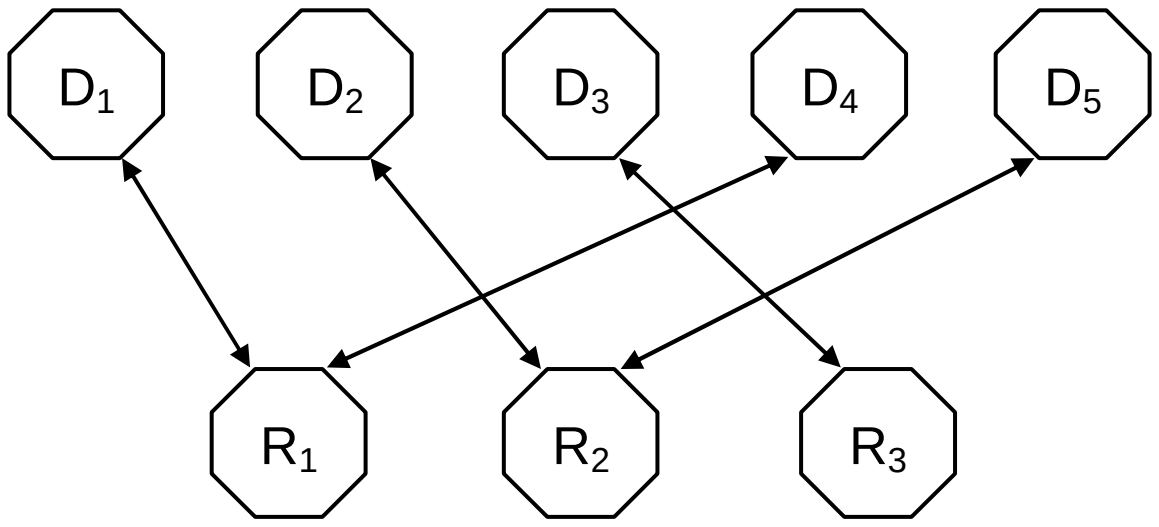

FIGURE: 2.3.1.1

The relational matrix E got from the above map is

$$E = \begin{bmatrix} 1 & 0 & 0 \\ 0 & 1 & 0 \\ 0 & 0 & 1 \\ 1 & 0 & 0 \\ 0 & 1 & 0 \end{bmatrix}.$$

If A = (1 0 0 0 0) is passed on in the relational matrix E, the instantaneous vector, AE = (1 0 0) implies that the student is a good student. Now let AE = B, $BE^T$ = (1 0 0 1 0) which implies that the teaching is good and he / she is a kind teacher. Let $BE^T$ = $A_1$, $A_1E$ = (2 0 0) after thresholding we get $A_1E$ = (1 0 0) which implies that the student is good, so on and so forth.

Now we have worked on the real world problem, which studies the Employee- Employer Relationship model in the following section. The reader is expected to develop several other model relationships and illustrations.

## 2.3.2 Models Illustrating FRM and combined FRMs

Now we in this section give two illustrations, which are data, collected from the real world problems analyzed using both FRMs and combined FRMs.

***Example 2.3.2.1:*** We use FRMs in Employee Employer Relationship model.

The employee-employer relationship is an intricate one. For, the employers expect to achieve performance in quality and production in order to earn profit, on the other hand employees need good pay with all possible allowances. Here we have taken three experts opinion in the study of Employee and Employer model. The three experts whose opinions are taken are the Industry Owner, Employees' Association Union Leader and an Employee. The data and the opinion are taken only from one industry. Using the opinion we obtain the hidden patterns.

The following concepts are taken as the nodes relative to the employee. We can have several more nodes and also several experts' opinions for it a clearly evident theory which professes that more the number of experts the better is the result.

We have taken as the concepts / nodes of domain only 8 notions which pertain to the employee.

$D_1$ – Pay with allowances and bonus to the employee
$D_2$ – Only pay to the employee
$D_3$ – Pay with allowances (or bonus) to the employee
$D_4$ – Best performance by the employee
$D_5$ – Average performance by the employee
$D_6$ – Poor performance by the employee
$D_7$ – Employee works for more number for hours
$D_8$ – Employee works for less number of hours.

$D_1, D_2, \ldots, D_8$ are elements related to the employee space which is taken as the domain space.

We have taken only 5 nodes / concepts related to the employer in this study.

These concepts form the range space which is listed below.

$R_1$ – Maximum profit to the employer
$R_2$ – Only profit to the employer
$R_3$ – Neither profit nor loss to the employer
$R_4$ – Loss to the employer
$R_5$ – Heavy loss to the employer

The directed graph as given by the employer is given in Figure 2.3.2.1. The associated relational matrix $E_1$ of the employer as given by following.

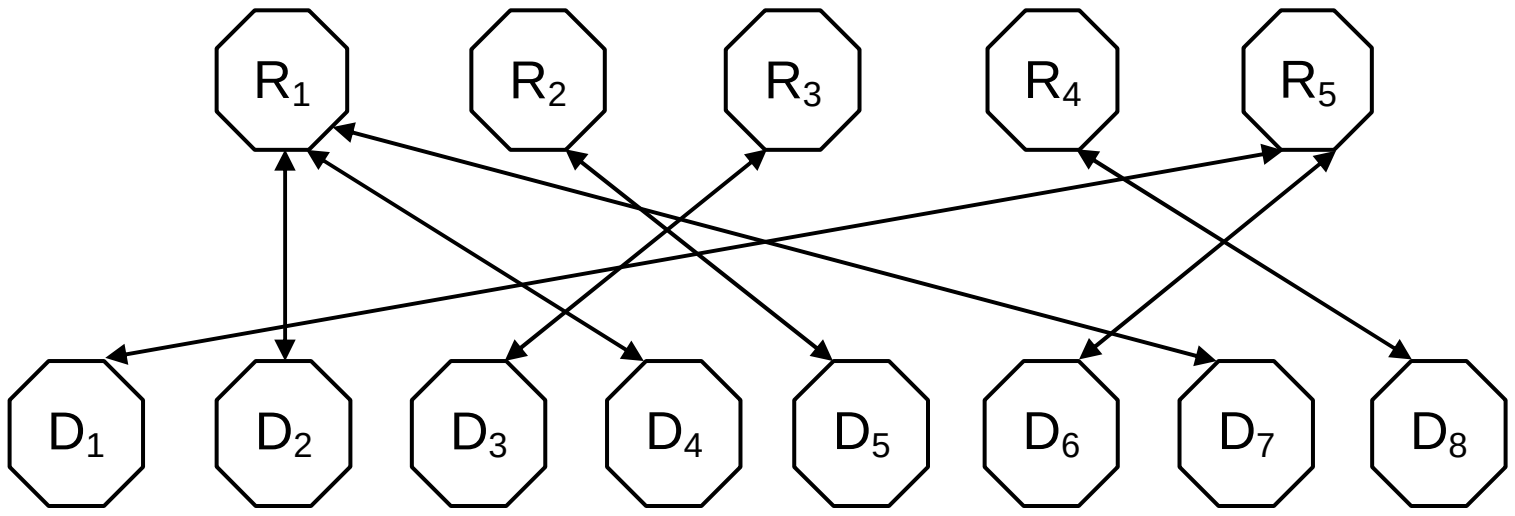

FIGURE: 2.3.2.1

$$E_1 = \begin{bmatrix} 0 & 0 & 0 & 0 & 1 \\ 1 & 0 & 0 & 0 & 0 \\ 0 & 0 & 1 & 0 & 0 \\ 1 & 0 & 0 & 0 & 0 \\ 0 & 1 & 0 & 0 & 0 \\ 0 & 0 & 0 & 0 & 1 \\ 1 & 0 & 0 & 0 & 0 \\ 0 & 0 & 0 & 1 & 0 \end{bmatrix}.$$

Suppose we consider the node $D_1$ to be in the on state and rest of the nodes in the off state. That is the employee is paid with allowance and bonus i.e. $G_1$ = (1 0 0 0 0 0 0 0). The effect of $G_1$ on the expert systems $E_1$ is $G_1E_1$ = (0 0 0 0 1), it belongs to the range space R.

Let $G_1E_1 = H_1$ = (0 0 0 0 1), $H_1E_1^T$ = (1 0 0 0 0 1 0 0) $\in$ D, after updating and thresholding the instantaneous vector at each stage we obtain the following chain

$$G_1 \to H_1 \to G_2 \to H_1$$

i.e., $G_1$ is a fixed point and according to the opinion of the employer who is taken as an expert we see if the employee is paid with allowances and bonus the company suffers a heavy loss due to the poor performance of the employee.

Suppose we input the vector $G_3$ = (0 0 0 1 0 0 0 0) which indicates that the node $D_4$ "best performance by the employee" is in the on state, Effect of $G_3$ on the system, $G_3$ $E_1$ = (1 0 0 0 0) = $H_3 \in$ R. $H_3E_1^T$ = (0 1 0 1 0 0 1 0) $\in$ D.

After updating and thresholding the instantaneous vector at each stage we obtain the following chain:

$$G_3 \to H_3 \to G_4 \to H_3$$

We see from the above the resultant is also a fixed point.

According to the first expert we see the company enjoys maximum profit by giving only pay to the employee in spite of his best performance and putting in more number of working hours.

Now to analyze the effect of employee on the employer let us input the vector $S_1$ = (1 0 0 0 0) indicating the on state of the node $R_1$ (maximum profit to the employer)

$$S_1E^T{}_1 = (0\ 1\ 0\ 1\ 0\ 0\ 1\ 0) = T_1 \in D.$$

After updating and thresholding the instantaneous vector at each stage we obtain the following chain

$$S_1 \to T_1 \to S_1.$$

From the above chain we see that when $S_1$ is passed on, it is a limit cycle.

The company enjoys maximum profit by getting best performance and more number of working hours from the employee and by giving only pay to the employee.

Suppose we input the vector $S_2 = (0\ 0\ 0\ 0\ 1)$ indicating the on state of the node $R_5$.

$$S_2 E_1^T = T_2 = (1\ 0\ 0\ 0\ 0\ 1\ 0\ 0) \in D.$$

After updating and thresholding the instantaneous vector at each stage we obtain the following chain

$$S_2 \to T_2 \to S_2.$$

Thus in his industry the employer feels that the employee performs poorly inspite of getting pay with allowances and bonus, company suffers a heavy loss.

The union leader of the same company was asked to give his opinion keeping the same nodes for the range space and the domain space i.e. as in case of the first expert. The directed graph given by the union leader of the company is given by the following diagram.

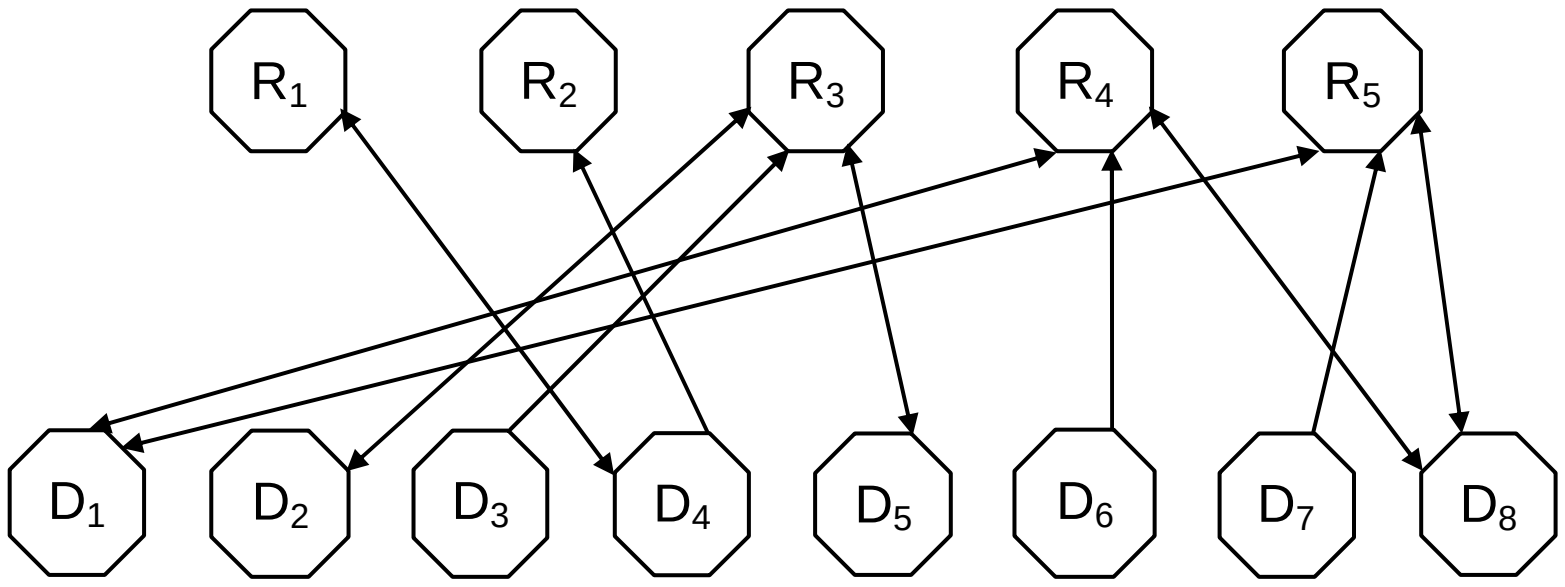

FIGURE: 2.3.2.2

The related matrix of the directed graph given by the second expert is given by $E_2$

$$E_2 = \begin{bmatrix} 0 & 0 & 0 & 1 & 1 \\ 0 & 0 & 1 & 0 & 0 \\ 0 & 0 & 1 & 0 & 0 \\ 1 & 1 & 0 & 0 & 0 \\ 0 & 0 & 1 & 0 & 0 \\ 0 & 0 & 0 & 1 & 0 \\ 0 & 0 & 0 & 0 & 1 \\ 0 & 0 & 0 & 1 & 1 \end{bmatrix}.$$

Let us input the vector $L_1 = (1\ 0\ 0\ 0\ 0\ 0\ 0\ 0)$ which indicates that the node $D_1$ viz. employee is paid with allowance and bonus is in the on state.

$$L_1E_2 \hookrightarrow (1\ 0\ 0\ 1\ 1) = N_1 \in R$$
$$N_1E_2^T \hookrightarrow (1\ 0\ 0\ 1\ 0\ 1\ 1\ 1) = L_2$$

$L_2$ is got only after updating and thresholding. Now $L_2E_2 \to N_2$ = (11011) which is got after updating and thresholding

$$N_2\, E_2^T \hookrightarrow (1\ 0\ 0\ 1\ 0\ 1\ 1\ 1) = L_3$$

$$L_1 \to N_1 \to L_2 \to N_2 \to L_3 = L_2 \to N_2$$

is a fixed point.

Thus the union leader's viewpoint is in a way very balanced though seemingly contradictory. When the company pays an employee with pay, allowances and bonus the company may have maximum profit or only profit. Or on the contrary if the employee does not put forth more number of working hours or puts in only the required number of working hours and if his/her performance is poor or average certainly the company will face loss and or heavy loss. Thus both can happen depending highly on the nature of the employees.

Now we proceed on to study the same industry using the third expert as an employee of the same industry.

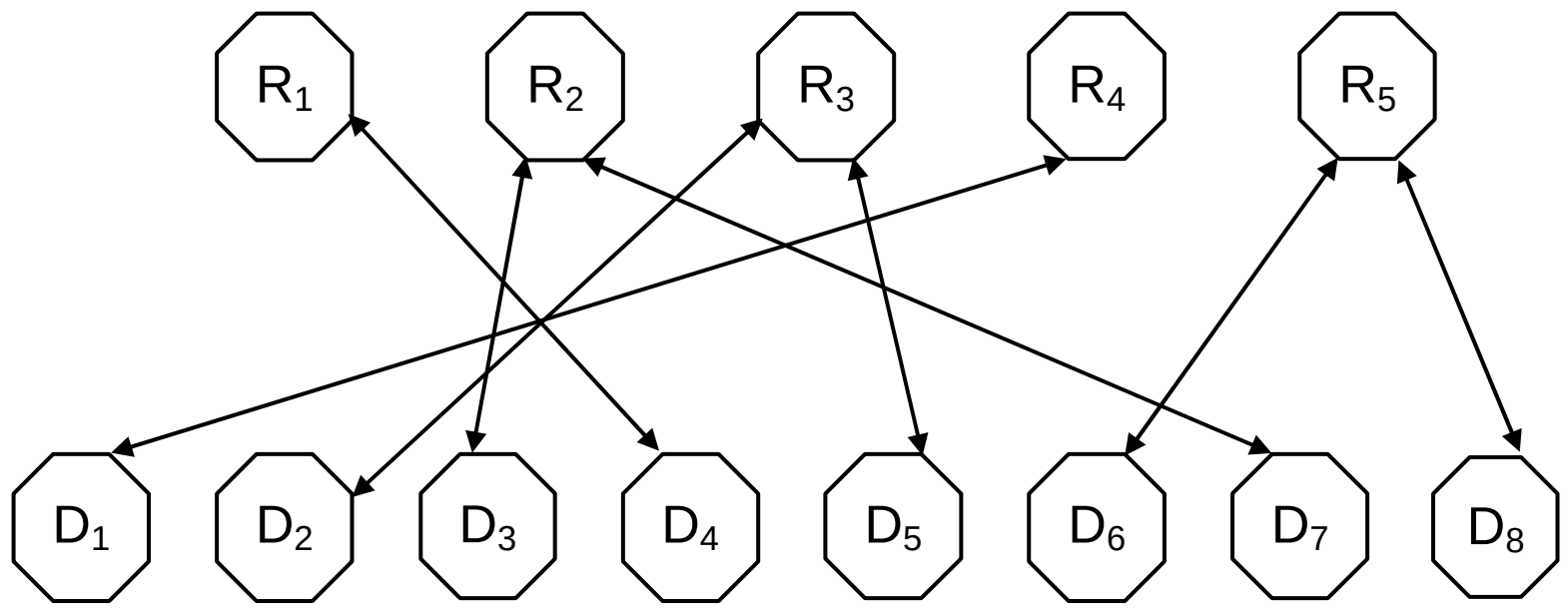

FIGURE: 2.3.2.3

The associated relational matrix $E_3$ of the third expert's opinion got from the directed graph given in Figure 2.3.2.3 is as follows:

$$E_3 = \begin{bmatrix} 0 & 0 & 0 & 1 & 0 \\ 0 & 0 & 1 & 0 & 0 \\ 0 & 1 & 0 & 0 & 0 \\ 1 & 0 & 0 & 0 & 0 \\ 0 & 0 & 1 & 0 & 0 \\ 0 & 0 & 0 & 0 & 1 \\ 0 & 1 & 0 & 0 & 0 \\ 0 & 0 & 0 & 0 & 1 \end{bmatrix}.$$

Let us input the vector $U_1$ = (1 0 0 0 0 0 0 0) indicating the on state of the note $D_1$, effect of $U_1$ on the system $E_3$ is

$$U_1 E_3 = V_1 = (0\ 0\ 0\ 1\ 0) \in R.$$

After thresholding and updating the vector. Now $V_1E^T_3$ = (1 0 0 0 0 0 0 0) = $U_1$. So

$$U_1 \to V_1 \to U_1.$$

Hence when $U_1$ is passed on it is a limit cycle. So according to the employee when the company gives pay, bonus/ perks and allowances the company suffers a loss.

Suppose we input vector $U_2$ = (0 0 0 1 0 0 0 0) which indicates the on state of the node $D_4$, employee puts in his best performance; effect of $U_2$ on the system $E_3$ is $U_2E_3$ = (1 0 0 0 0) $\in$ R.

After updating and thresholding the instantaneous vector at each state we get the following chain:

$$U_2 \to V_2 \to U_2$$

when $U_2$ is passed on, it is a limit cycle.

It is left for the reader to input any instantaneous state vector and obtain the resultant.

Now we proceed on to find the combined FRMs. We take the opinion of the three experts discussed above and find their opinions.

We first draw the directed graph of all the three experts, which is given by the Figure 2.3.2.4.

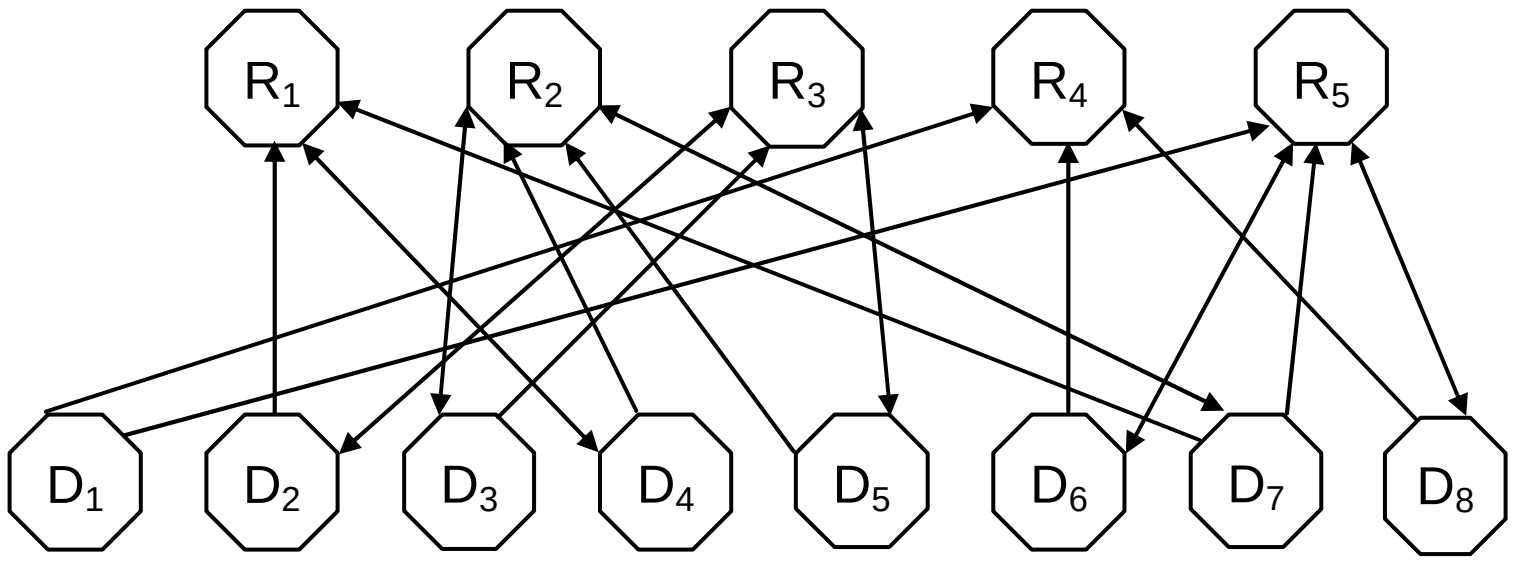

**FIGURE: 2.3.2.4**

The corresponding fuzzy relational matrix is given as the sum of the three fuzzy relational matrices $E_1$ , $E_2$ and $E_3$. Let

$$E = E_1 + E_2 + E_3 = \begin{bmatrix} 0 & 0 & 0 & 2 & 2 \\ 1 & 0 & 2 & 0 & 0 \\ 0 & 1 & 2 & 0 & 0 \\ 3 & 1 & 0 & 0 & 0 \\ 0 & 1 & 2 & 0 & 0 \\ 0 & 0 & 0 & 1 & 2 \\ 1 & 1 & 0 & 0 & 1 \\ 0 & 0 & 0 & 2 & 2 \end{bmatrix}.$$

Let us now input the state vector $X_1 = (1\ 0\ 0\ 0\ 0\ 0\ 0\ 0)$, indicating the on state of the node $D_1$, effect of $X_1$ on the combined system is

$$X_1 E = (0\ 0\ 0\ 2\ 2) \hookrightarrow (0\ 0\ 0\ 1\ 1) \in R.$$

Let $Y_1 = (0\ 0\ 0\ 1\ 1)$,

$$Y_1 E^T = (4\ 0\ 0\ 0\ 0\ 3\ 1\ 4) \hookrightarrow (1\ 0\ 0\ \ 0\ 0\ 1\ 1\ 1) \in D.$$

($\hookrightarrow$ denotes after updating and thresholding the vector).

Let $X_2 = (1\ 0\ 0\ 0\ 0\ 1\ 1\ 1)$;

$$X_2E = (1\ 1\ 0\ 5\ 7) \hookrightarrow (1\ 1\ 0\ 1\ 1) = Y_2 \in R.$$

$$Y_2\ E^T = (4\ 1\ 1\ 4\ 1\ 3\ 2\ 4) \hookrightarrow (1\ 1\ 1\ 1\ 1\ 1\ 1\ 1) \in D.$$

Let $X_3 = (1\ 1\ 1\ 1\ 1\ 1\ 1\ 1)$

$$X_3\ E = (5\ 4\ 6\ 5\ 7) \hookrightarrow (1\ 1\ 1\ 1\ 1) = Y_3.$$

$$X_1 \to Y_1 \to X_2 \to Y_2 \to X_3 - Y_3 \to X_3 \to Y_3$$

Thus $X_1$ is a fixed point.

In fact a lot of discussion can be made on such fixed points, for the opinion of the employee and the employer or the union leader happens to be very contradictory that is why it is visibly seen that all the nodes both in the domain space and the range space becomes on, at the very advent of seeing the effect of only one node of the domain space to be in the on state i.e. $D_1$ – pay, bonus and allowances. The interpretations are also to be carefully given. On one hand when the employee is paid with pay bonus and allowances, they may work so well to see that the company runs with maximum profit which will automatically turn on all other linguistic nodes like the notions of only profit to the employer, neither profit nor loss to the employer etc.

Likewise in case of the domain nodes all nodes come to the on state. Thus we are not able to distinguish or give a nice interpretation of the state vectors.

For more about this illustration refer [237, 243-4].

The fuzzy relational maps can also be used in the prediction using the past year data. So in this case it is not only the expert's opinion but the interpretation of the data and the opinion analyzed using other methods. We for the first time give such analysis in case of cement industries. We analyze how to maximize the production by giving the maximum satisfaction to the employees, we have used the data from Ramco Cement Industries for the years 1995-2001.

***Example 2.3.2.2:*** Maximizing Production by Giving Maximum Satisfaction to Employee using FRMs

At present India needs huge quantity of cement for the construction of dwellings, houses, apartments, dams, reservoirs, village roads, flyovers and cementing of the entire national highways. Hence maximizing cement production in each and every factory is essential which is directly dependent on the relationship between the employer and the employees. Employer and Employees congenial relationship (Industrial Harmony) is a most complicated one.

For example employer expects to achieve consistent production, quality product at optimum production to earn profit. However there may be profit, no loss and no profit, loss in business or heavy loss depending on various factors such as demand and supply, rent, electricity, raw materials transportation and consumables, safety, theft, medical aid, employees co-operation and enthusiasm. At the same time some of the expectations of the employees are pay, various allowances, bonus, welfare unanimity such as uniforms, hours of work and job satisfaction etc. Since it is a problem of both employer and employee we analyse taking into consideration the problem on both sides using fuzzy theory in general and FRMs in particular; as the very nature of the problem is one, which is dominated by uncertainties.

Some industries provide

1. Additional incentive for regularity in attendance
2. Additional days of vacation
3. Incentive linked with production level
4. Award for good suggestion for operation
5. Forming quality circles in various department for finding out new ways and means to carry out jobs economically, quickly and with safety
6. Provide hospital for employees and their family
7. Credit society and co-operative stores
8. Awards for highest sale made by salesman

9. Overtime wages
10. Loan facilities.

The employer's final achievement is to maximize the production level; so, when production linked incentive was introduced in a cement industry, the employees were so enthusiastic to run the kiln continuously and efficiently to ensure more than the targeted minimum production. Employees voluntarily offered their services in case of kiln stoppages and made the kiln to run again with a minimum shutdown. Small incentive to employees made the industry to earn considerable profit out of maximum production than the previous year. This study is significant because most of the cement industries have the common types of problems between employee and the employer.

Here in this analysis we are also trying to give a form of best relationship between the employees and employers, because the smooth relationship between the employers and employees is an important one for the cement industries to maximize the production. Probably, so far no one has approached the problem of finding a best form of relationship between employers and employees using FRMs. We use the connection fuzzy relational matrix to find the best form of relationship between the employees and employers and aim to maximize the production level.

Here, we approach the employee and employer problem of finding a best form of relationship between employee and employer. Thus this study tries to improve maximize the production level in cement industries. A good relationship between the employee and employer is very essential to run the industry successfully with least friction.

For this study the raw data is obtained from the cement industrialists and the employees; which is converted into a relational map. The relational map takes edge values as positive real numbers. If the number of concepts in the domain space is M and that of the range space is N we get an $M \times N$ causal relational matrix, which we call as the relational matrix.

Let $X_i$ and $Y_j$ denote the two nodes of the relational map. The nodes $X_i$ and $Y_j$ are taken as the attributes of the employee and employer respectively. The directed edge from $X_i$ to $Y_j$ and $Y_j$ to $X_i$ denotes the causality relations. Every edge in the relational map is weighted with a positive real number . Let $X_1$, $X_2,\ldots, X_M$ be the nodes of the domain space X and $Y_1, Y_2,\ldots, Y_N$ be the nodes of the range space Y of the relational map.

After obtaining a relational matrix, we find an average matrix; for simplification of calculations. The average matrix is then converted into a fuzzy matrix. Finally using the different parameter $\alpha$ (membership grades) we identify the best form of relationship between employee and employer in cement industries.

We approach the employee and the employer problems in cement industries using fuzzy relational matrix. Probably, so far no one has approached the employee and employer problems via a fuzzy relational matrix method. The raw data under investigation is classified under six broad heads viz. $Y_1$, $Y_2$, $Y_3$, $Y_4$, $Y_5$, $Y_6$ which are nodes associated with the employer. $X_1$, $X_2,\ldots, X_8$ are the eight broad heads, which deal with the attributes of the employees. It is pertinent to mention here one can choose any number of nodes associated with employees or employer. We prefer to choose $8 \times 6$ matrix for spaces X and Y respectively. Using the fuzzy matrix we find a best relationship between employer and employee in two-stages.

In the first stage, we convert the data into relational map, the relational matrix obtained from the relational map is then converted into an average matrix. In this matrix, take along the columns the data related to the employer and along the rows the data related to the employee.

In the second stage, we convert the average matrix into fuzzy matrix using different parameter $\alpha$ ($\alpha$-membership grade), and we give the graphical illustration to the fuzzy matrix, which enables one to give a better result.

We now describe the problem.

Using the data from (Ramco) cement industries, we analyse the data via a fuzzy relational matrix and obtain a best relationship between employee and employer in cement industries. Hence using these suggestions got from the study, cement industries can maximize their production. Here we find the best relationship between the employee and the employer.

Here we consider only eight attributes of the employee and its six effects on employers in production level. Hence, in a cement industry the attributes of the employee are described by the eight nodes $X_1$, $X_2$, $X_3$, $X_4$, $X_5$, $X_6$, $X_7$, $X_8$ and the six effects on the employers in production level are described by the nodes $Y_1$, $Y_2$, $Y_3$, $Y_4$, $Y_5$, $Y_6$.

In the first stage of the problem, the data obtained from the cement industry (1995-2001) are used in the illustration of the problem and is finally applied in the fuzzy matrix model to verify the validity of the method described. The relational map is obtained using the above nodes.

From the relational map we get the $M \times N$ relational matrix. Let M represent the six effects of employers. Let N represent the eight attributes of employees feelings. Production levels are treated as rows and the various attributes of employee feelings are treated as columns resulting in the $M \times N$ relational matrix into an average matrix.

In the second step we use mean, standard deviation of each of the columns of the $M \times N$ matrix and parameter $\alpha$ (membership grade $\alpha \in [0,1]$ to convert the average matrix $(a_{ij})$ into the fuzzy matrix $(b_{ij})$; where i represents the $i^{th}$ row and j represents the $j^{th}$ column. We calculate the mean $\mu_j$ and the standard deviation $\sigma_j$ for each attribute j, for j = 1, 2,…, n using the average matrix $(a_{ij})$.

For varying values of the parameter $\alpha$ where $\alpha \in [0, 1]$. We determine the values of the entry $b_{ij}$ in the average matrix using the following rule:

$$b_{ij} = \begin{cases} 0 & \text{if} \quad a_{ij} \le \mu_j - \alpha * \sigma_j \\ \dfrac{a_{ij} - (\mu_j - \alpha * \sigma_j)}{(\mu_j + \alpha * \sigma_j) - (\mu_j - \alpha * \sigma_j)} & \text{if } a_{ij} \in (\mu_j - \alpha * \sigma_j, \mu_j + \alpha * \sigma_j) \\ 1 & \text{if} \quad a_{ij} \ge \mu_j + \alpha * \sigma_j \end{cases}$$

Here ‘*’ denotes the usual multiplication.

Thus for different value of $\alpha$, we obtain different fuzzy matrices. Finally we add up the rows of each fuzzy matrix and we define the new fuzzy membership function and allocate a value between [0, 1] to each row sum. Here the highest membership grade gives the best form of relationship between employee and employer, which maximizes the production level.

STAGE 1

The cement industry is having eight types of attributes related to the employee, which are as follows:

- $X_1$ – Salaries and wages
- $X_2$ – Salaries with wages and bonus to the employee
- $X_3$ – Bonus to the employee
- $X_4$ – Provident fund (PF) to the employee
- $X_5$ – Employee welfare medical
- $X_6$ – Employee welfare LTA (Leave travel allowances)
- $X_7$ – Employee welfare others
- $X_8$ – Staff Training expenses.

Now we take the six effects on the employer when he proposes to pay and get the work done by the employee to maximize the production level, which are as follows:

- $Y_1$ – 1995 - 1996 Production level is 8,64,685 tonnes
- $Y_2$ – 1996 - 1997 Production level is 7,74,044 tonnes
- $Y_3$ – 1997 - 1998 Production level is 7,22,559 tonnes
- $Y_4$ – 1998 - 1999 Production level is 8,19,825 tonnes
- $Y_5$ – 1999 - 2000 Production level is 7,98,998 tonnes
- $Y_6$ – 2000 - 2001 Production level is 7,79,087 tonnes.

The problem now is to find to best relationship between the employee and employer, so we convert the $X_1$, $X_2$, $X_3$, $X_4$, $X_5$, $X_6$, $X_7$, $X_8$ and $Y_1$, $Y_2$, $Y_3$, $Y_4$, $Y_5$, $Y_6$ into the relational map. The relational map is shown in the Figure 2.3.2.5.

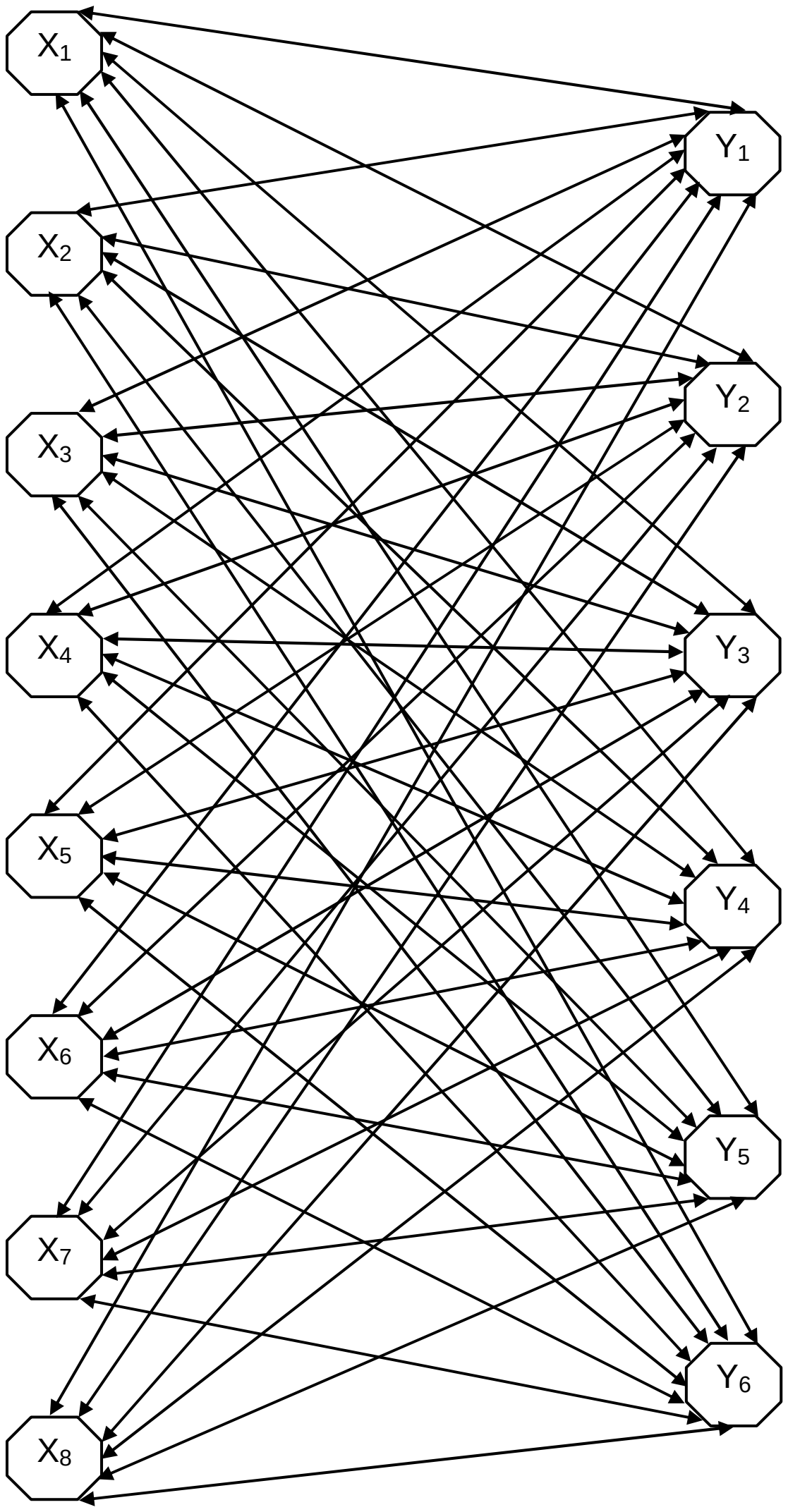

FIGURE: 2.3.2.5

We obtain the initial 8 × 6 relational matrix from the relational map using weights from the data.

$$\begin{bmatrix} 70.99 & 82.61 & 11.61 & 5.856 & 2.239 & 1.23 & 7.19 & 0.867 \\ 74.33 & 85.05 & 10.71 & 5.558 & 2.058 & 0.88 & 6.02 & 0.427 \\ 74.25 & 84.29 & 10.04 & 5.969 & 2.704 & 0.769 & 6.068 & 0.144 \\ 71.08 & 81.42 & 10.34 & 7.070 & 2.386 & 1.197 & 7.411 & 0.498 \\ 70.65 & 80.85 & 10.19 & 7.102 & 2.589 & 1.112 & 8.093 & 0.238 \\ 72.03 & 80.06 & 8.028 & 6.297 & 3.276 & 1.373 & 8.563 & 0.413 \end{bmatrix}$$

Now for simplification of calculation, we convert the relational matrix into the average matrix which is given by $(a_{ij})$.

$$\begin{bmatrix} 35.49 & 41.30 & 5.805 & 2.928 & 1.119 & 0.615 & 3.595 & 0.433 \\ 37.16 & 42.53 & 5.355 & 2.779 & 1.029 & 0.44 & 3.01 & 0.213 \\ 37.13 & 42.15 & 5.02 & 2.985 & 1.352 & 0.385 & 3.034 & 0.072 \\ 35.54 & 40.71 & 5.17 & 3.535 & 1.193 & 0.598 & 3.076 & 0.249 \\ 35.33 & 40.43 & 5.095 & 3.551 & 1.295 & 0.556 & 4.046 & 0.119 \\ 36.02 & 40.03 & 4.014 & 3.148 & 1.638 & 0.687 & 4.281 & 0.207 \end{bmatrix}$$

To find a best form of relationship between the employee and employer, we use mean, standard deviation and the parameter $\alpha \in [0, 1]$ and proceed to the second stage of the problem to convert the average matrix into a fuzzy matrix.

STAGE 2

In the second stage we use mean ($\mu$), standard deviation ($\sigma$) and the parameter $\alpha \in [0, 1]$ to find the best form of relationship between employer and employee. To convert the above average matrix into a fuzzy matrix $(b_{ij})$.

$$b_{ij} \in \left[\, 0,\ \ \frac{(a_{ij} - \mu_i - \alpha * \sigma_j)}{(\mu_j + \alpha * \sigma_j) - (\mu_j - \alpha * \sigma_j)},\ \ 1 \,\right]$$

where i represents the $i^{th}$ row and j represents the $j^{th}$ column. The value of the entry $b_{ij}$ corresponding to each intersection is determined from this interval. This interval is obtained strictly by using the average and standard deviation calculated from the raw data.

The calculations are as follows: First we calculate the $\mu_j$ corresponding to each column of the matrix. $\mu_1 = 36.1116$, $\mu_2 = 41.191$, $\mu_3 = 5.076$, $\mu_4 = 3.154$, $\mu_5 = 1.271$, $\mu_6 = 0.547$, $\mu_7 = 3.612$, $\mu_8 = 0.2155$, where $\mu_j$ are the average of each column respectively for j = 1, 2,…, 8.

Now the standard deviation $\sigma_j$ is calculated as follows.

**TABLE 1: MEAN AND STANDARD DEVIATION OF COLUMN 1**

When $\mu_1 = 36.1116$

| d | $d^2$ |
|---|---|
| 0.6216 | 0.386386 |
| 1.048 | 1.09830 |
| 1.0184 | 1.03713 |
| 0.5716 | 0.32672 |
| 0.7816 | 0.6108985 |
| 0.0916 | 0.008390 |
| $\Sigma d = 4.13280$ | $\Sigma d^2 = 3.4678245$ |

$$\sigma_1 = \sqrt{0.577970 - 0.4744454}$$

$$\sigma_1 = 0.32175.$$

In a similar way, the values of $\sigma_j$'s are as follows:

$\sigma_2 = 0.41368$ when $\mu_2 = 41.191$,
$\sigma_3 = 0.390285$ when $\mu_3 = 5.076$,
$\sigma_4 = 0.141763$ when $\mu_4 = 3.154$,
$\sigma_5 = 0.1162292$ when $\mu_5 = 1.271$,
$\sigma_6 = 0.05250$ when $\mu_6 = 0.547$,
$\sigma_7 = 0.253748$ when $\mu_7 = 3.612$,
$\sigma_8 = 0.07790$ when $\mu_8 = 0.2155$.

Now we calculate $\mu_j - \alpha * \sigma_j$ and $\mu_j + \alpha * \sigma_j$ when $\alpha = 0.1, 0.2, 0.3, 0.4, 0.5, 0.6, 0.7, 0.8, 0.9$ and 1.

**TABLE 2: FOR THE MEMBERSHIP GRADE $\alpha$ = 0.1**

| j | $\mu_j - \alpha * \sigma_j$ | $\mu_j + \alpha * \sigma_j$ |
|---|---|---|
| 1 | 36.0794250 | 36.1437750 |
| 2 | 41.1496320 | 41.2323680 |
| 3 | 5.03697150 | 5.11502850 |
| 4 | 3.13982370 | 3.16817630 |
| 5 | 1.25937700 | 1.28262200 |
| 6 | 0.54175000 | 0.55225000 |
| 7 | 3.58662500 | 3.63737400 |
| 8 | 0.20771000 | 0.22329000 |

From the table the corresponding fuzzy relational matrix for the parameter value $\alpha$ = 0.1 is as follows:

$$(b_{ij}) = \begin{array}{c} \\ Y_1 \\ Y_2 \\ Y_3 \\ Y_4 \\ Y_5 \\ Y_6 \end{array} \begin{array}{c} \begin{array}{cccccccc} X_1 & X_2 & X_3 & X_4 & X_5 & X_6 & X_7 & X_8 \end{array} \\ \begin{bmatrix} 0 & 1 & 1 & 0 & 0 & 1 & 0.165 & 1 \\ 1 & 1 & 1 & 0 & 0 & 0 & 0 & 0.339 \\ 1 & 1 & 0 & 0 & 1 & 0 & 0 & 0 \\ 0 & 0 & 1 & 1 & 0 & 1 & 1 & 1 \\ 0 & 0 & 0.743 & 1 & 1 & 1 & 1 & 0 \\ 0 & 0 & 0 & 0.288 & 1 & 1 & 1 & 0 \end{bmatrix} \end{array}.$$

The corresponding row sum of the above matrix is given by

$R_1$ = The first row sum = 4.165
$R_2$ = The second row sum = 3.3390
$R_3$ = The third row sum = 3
$R_4$ = The fourth row sum = 5
$R_5$ = The fifth row sum = 4.743
$R_6$ = The sixth row sum = 3.288

where $R_1$, $R_2 \ldots R_6$ are nothing but the corresponding years of $Y_1$, $Y_2, \ldots$ $Y_6$.

Now we define the new fuzzy membership function for graphical illustration, which converts the row sum to take values in the interval [0, 1].

We define the new fuzzy membership function as follows:

$$\mu_x(R_i) = \begin{cases} 1 & \text{if Row sum of max value} \\ \dfrac{R_i - \text{Rowsum of min value}}{\{\text{Rowsum of max value} - \text{Rowsum of min value}\}} & \text{if } \begin{array}{l}\text{Row sum of min value} \le R_i \\ \le \text{Row sum of max value}\end{array} \\ 0 & \text{if Row sum of min value} \end{cases}$$

Using the membership function the corresponding values of

$\mu_x$ $(R_i)$ is 0.5825, 0.1695, 0, 1, 0.87, 0.144

respectively where i = 1, 2, 3, 4, 5, 6.

Graphical illustration for $\alpha$ = 0.1 is shown in graph 2.3.2.1.

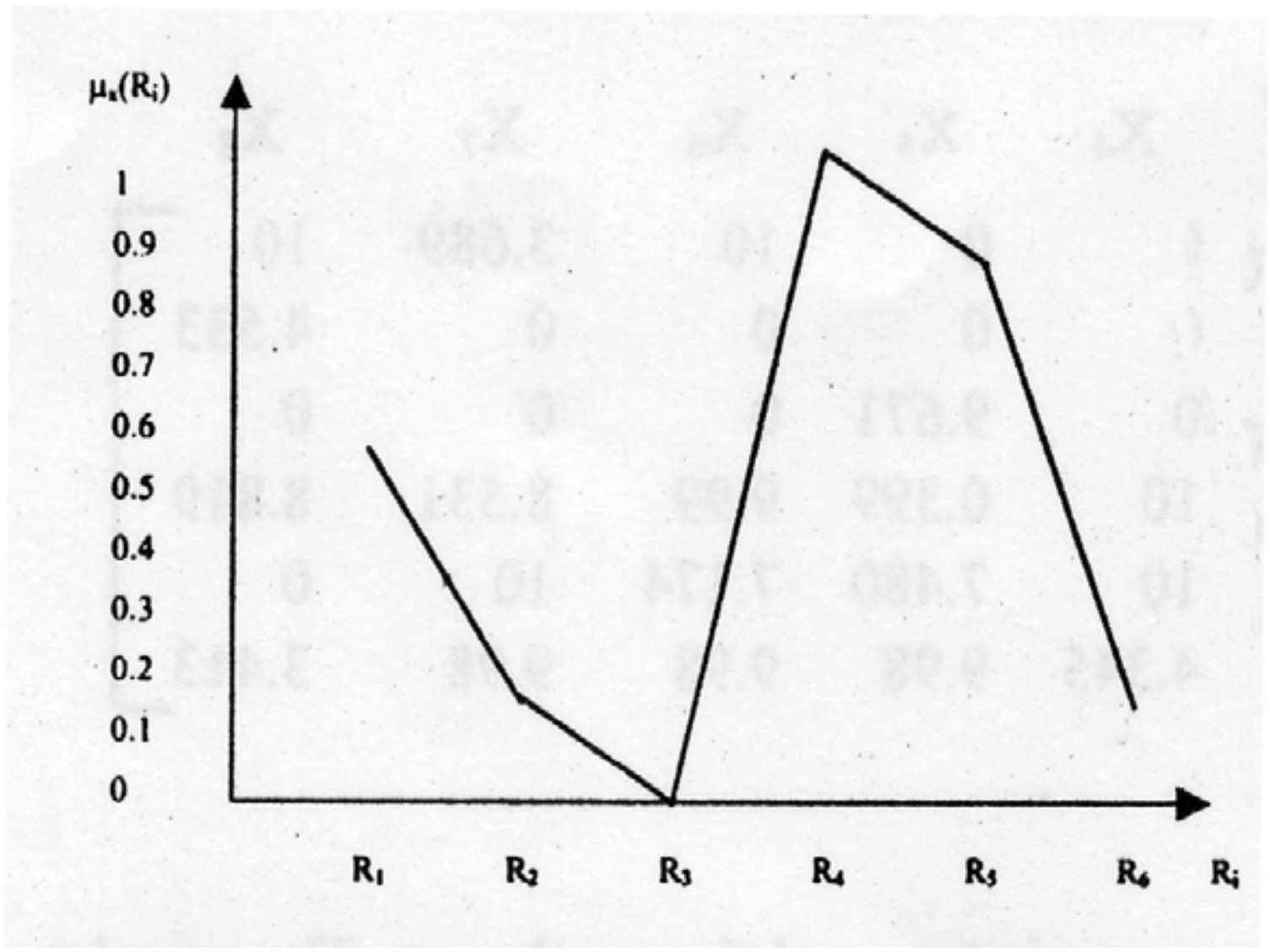

GRAPH 2.3.2.1: GRAPHICAL ILLUSTRATION FOR $\alpha$ = 0.1

Here the fourth row sum $R_4$ is getting the highest membership grade that is the membership grade is 1 and the fifth row sum $R_5$ is getting the next highest membership grade that is the membership grade is 0.87.

**TABLE 3: FOR THE MEMBERSHIP GRADE $\alpha$ = 0.2**

| J | $\mu_j - \alpha * \sigma_j$ | $\mu_j + \alpha * \sigma_j$ |
|---|---|---|
| 1 | 36.04720 | 36.175950 |
| 2 | 41.10826 | 41.273736 |
| 3 | 4.997900 | 5.1540500 |
| 4 | 3.125647 | 3.1823520 |
| 5 | 1.247754 | 1.2942450 |
| 6 | 0.536500 | 0.5575000 |
| 7 | 3.561250 | 3.6627496 |
| 8 | 0.199920 | 0.2310800 |

From the table the fuzzy matrix corresponding to the parameter value of $\alpha$ = 0.2 is as follows:

$$(b_{ij}) = \begin{array}{c} \\ Y_1 \\ Y_2 \\ Y_3 \\ Y_4 \\ Y_5 \\ Y_6 \end{array} \begin{array}{c} \begin{array}{cccccccc} X_1 & X_2 & X_3 & X_4 & X_5 & X_6 & X_7 & X_8 \end{array} \\ \left[\begin{array}{cccccccc} 0 & 1 & 1 & 0 & 0 & 1 & 0.333 & 1 \\ 1 & 1 & 1 & 0 & 0 & 0 & 0 & 0.419 \\ 1 & 1 & 0.1415 & 0 & 1 & 0 & 0 & 0 \\ 0 & 0 & 1 & 1 & 0 & 1 & 1 & 1 \\ 0 & 0 & 0.6218 & 1 & 1 & 0.9285 & 1 & 0 \\ 0 & 0 & 0 & 0.394 & 1 & 1 & 1 & 0.227 \end{array}\right] \end{array}$$

The corresponding row sum of the above matrix is given by

$R_1$ = The first row sum = 3.333
$R_2$ = The second row sum = 3.419
$R_3$ = The third row sum = 3.1415
$R_4$ = The fourth row sum = 5
$R_5$ = The fifth row sum = 4.5503
$R_6$ = The sixth row sum = 3.6210.

Using the membership function the corresponding values of $\mu_x(R_i)$ is 0.1030, 0.149, 0, 1, 0.7580, 0.258 respectively where i = 1, 2, 3, 4, 5, 6. Graphical illustration for $\alpha = 0.2$ is shown in graph 2.3.2.2.

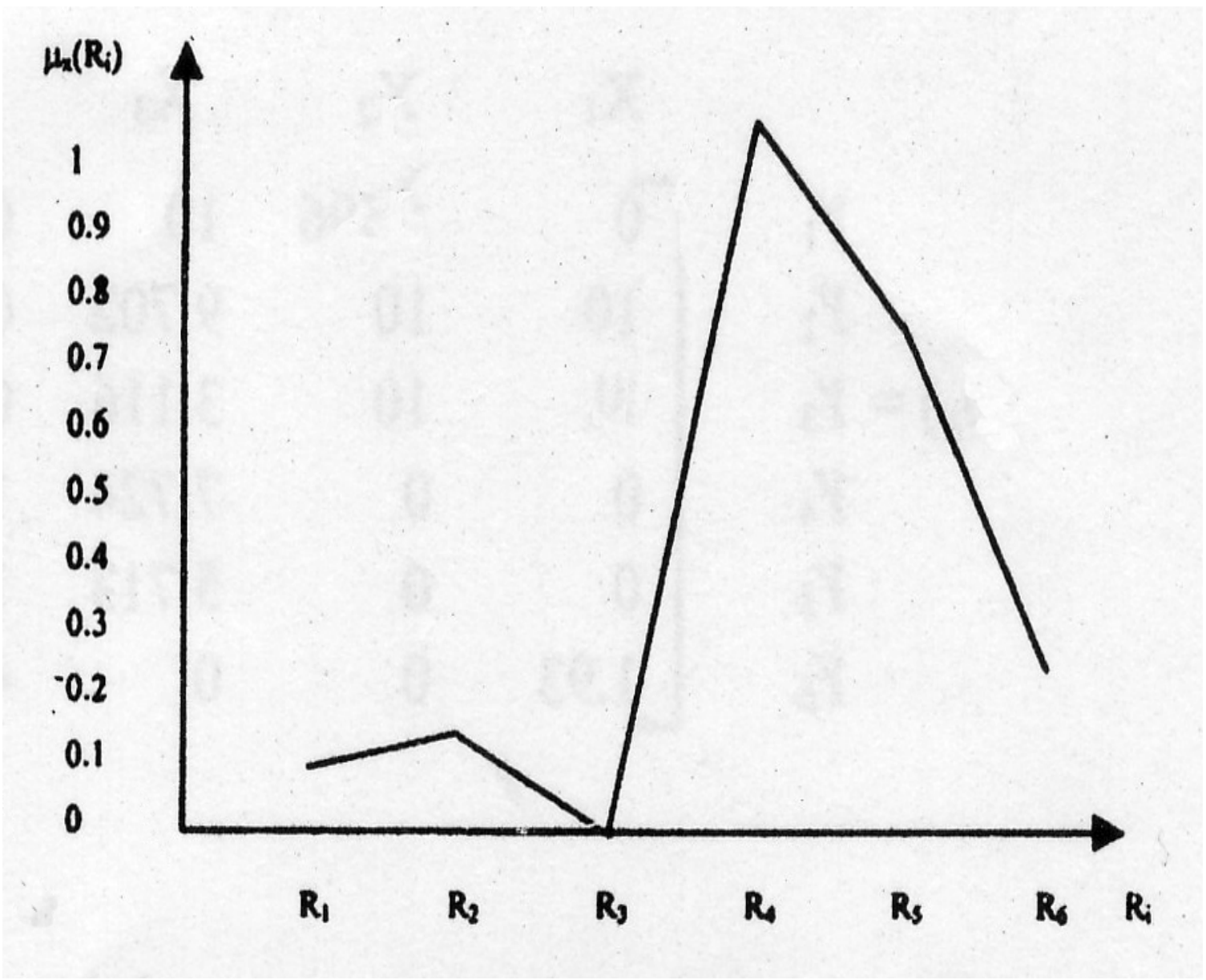

GRAPH 2.3.2.2: GRAPHICAL ILLUSTRATION FOR $\alpha = 0.2$

Here the fourth row sum $R_4$ is getting the highest membership grade that is the membership grade is 0.7580.

Using the same method described for the α cut values 0.1 and 0.2, we see when α takes the values 0.3, 0.5, 0.6 the corresponding row sum of the membership grades are as follows:

$$\mu_x(R_i)_{I=1,2,\ldots 6} = \begin{cases} 0.649, 0.113, 0, 1, 0.5791, 0.312 \ \text{ when } \ \alpha = 0.3 \\ 0.710, 0.094, 0, 1, 0.481, 0.597 \ \text{ when } \alpha = 0.5 \\ 0.79, 0.09, 0, 1, 0.492, 0.767 \ \text{ when } \ \alpha = 0.6 \end{cases} .$$

Here the fourth row sum $R_4$ is getting the highest membership grade that is the membership grade is 1 and the first row sum $R_1$ is getting the next highest membership grade.

Graphical illustration for α = 0.3, 0.5, 0.6 is shown in graph 2.3.2.3, graph 2.3.2.4 and graph 2.3.2.5.

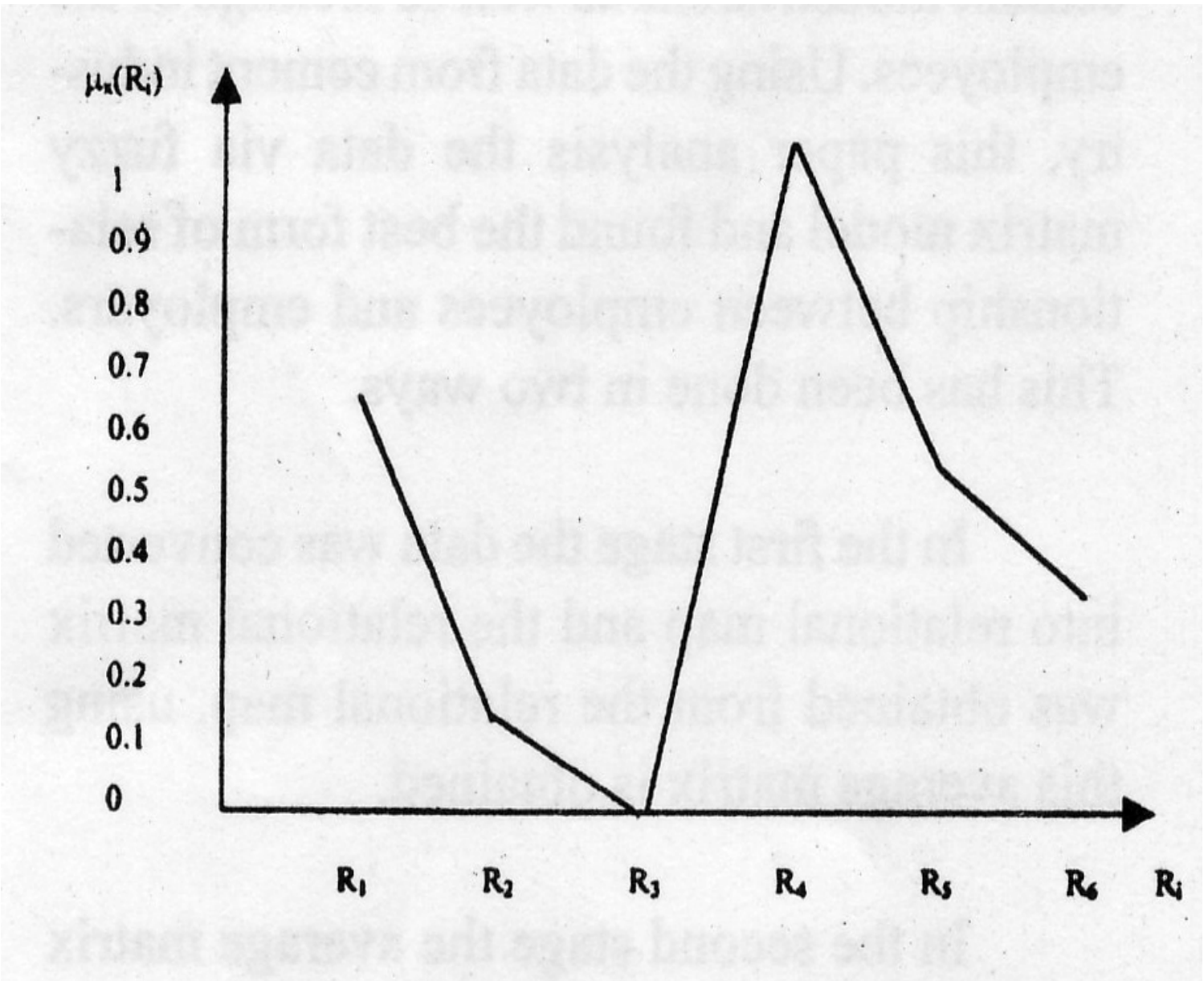

GRAPH 2.3.2.3: GRAPHICAL ILLUSTRATION FOR α = 0.3

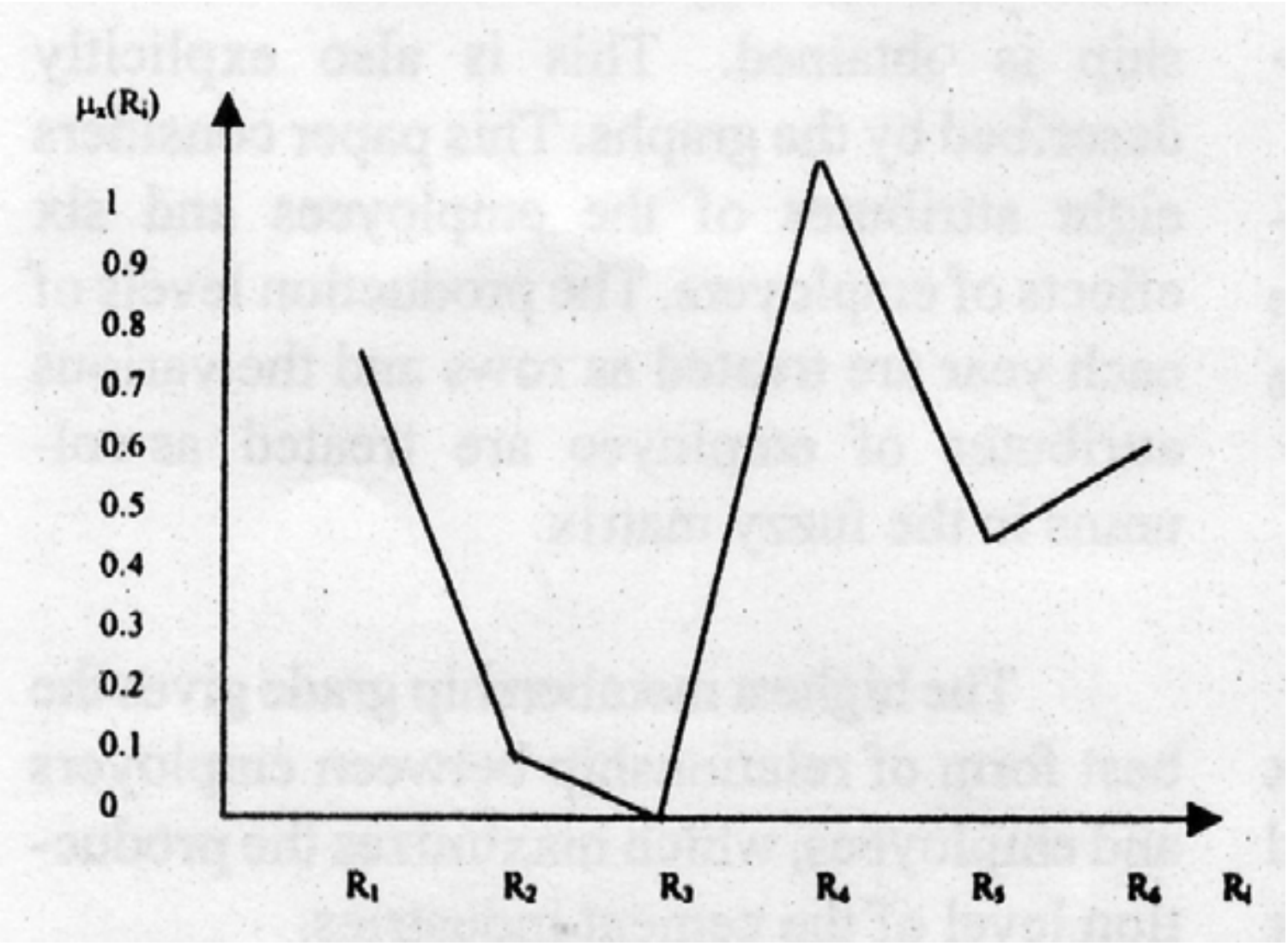

GRAPH 2.3.2.4: GRAPHICAL ILLUSTRATION FOR $\alpha = 0.5$

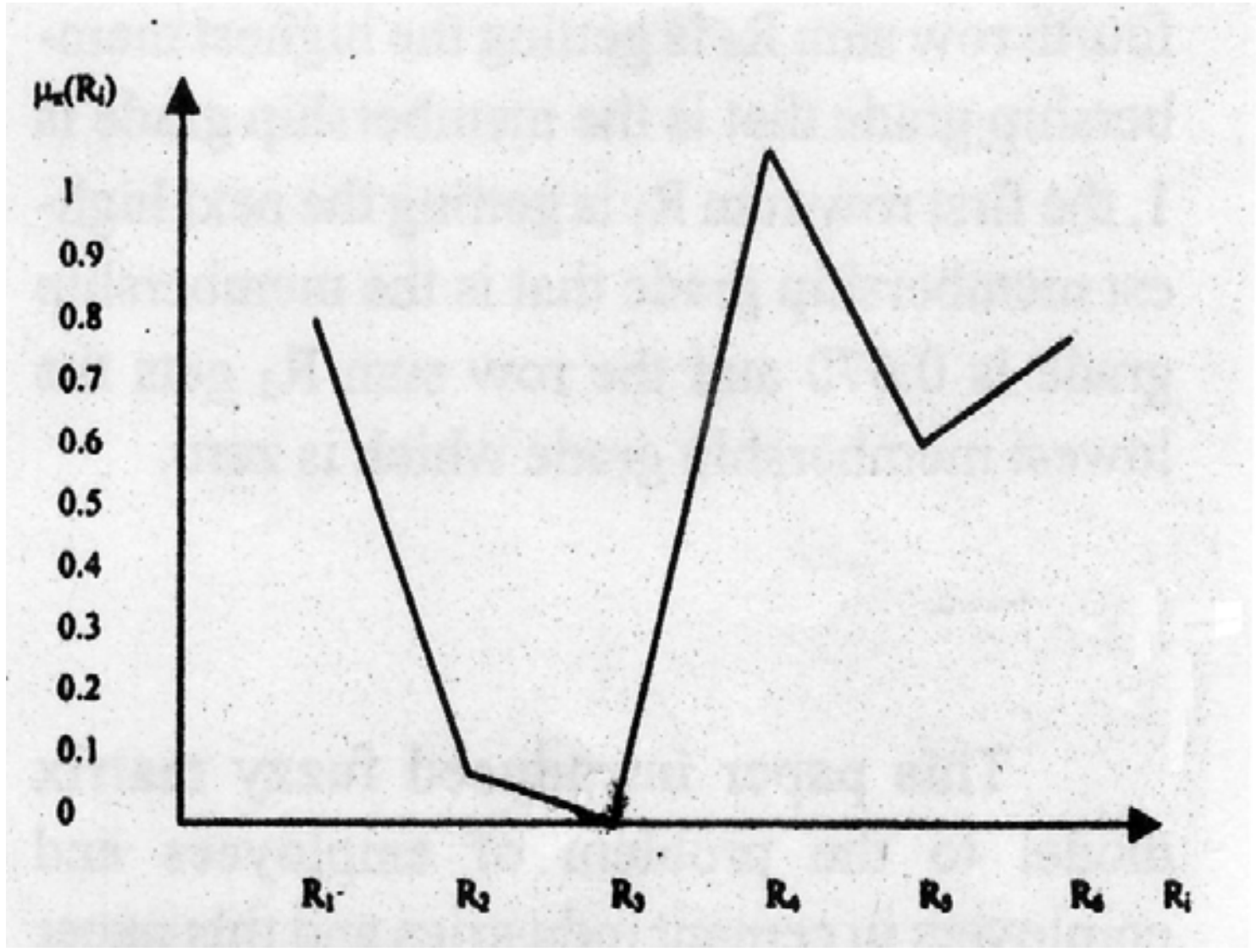

GRAPH 2.3.2.5: GRAPHICAL ILLUSTRATION FOR $\alpha = 0.6$

When $\alpha$ takes the value 0.4 the corresponding row sum of the membership grades are as follows:

$$\mu_x (R_i)_{i=1,2,\ldots 6} = \{0.411, 0.096, 0, 1, 0.493, 0.439 \text{ when } \alpha = 0.4.$$

Here the fourth row sum $R_4$ is getting the highest membership grade that is the membership grade is 1 and the fifth row sum $R_5$ is getting the next highest membership grade. Graphical illustration for $\alpha = 0.4$ is shown in graph 2.3.2.6.

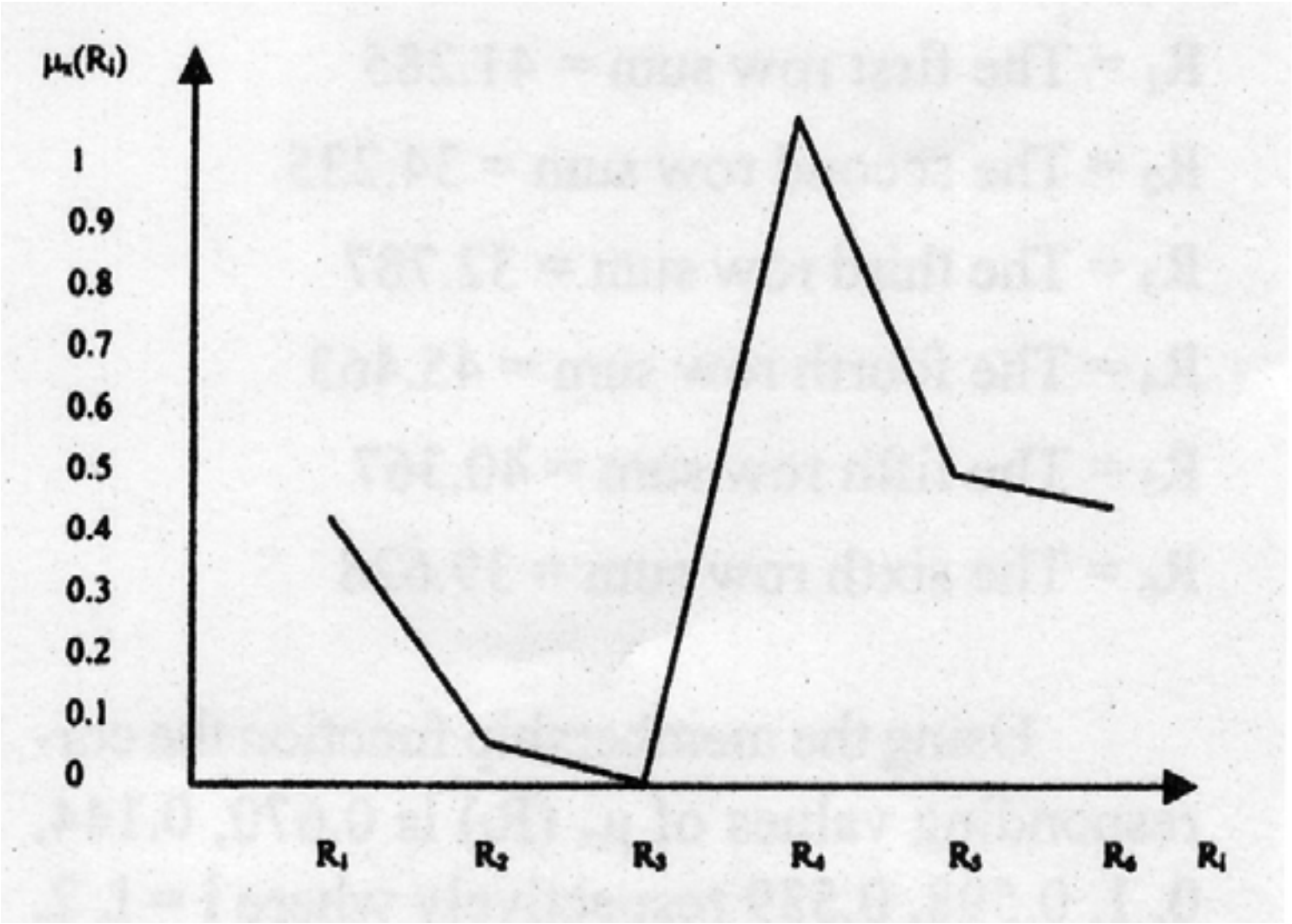

GRAPH 2.3.2.6: GRAPHICAL ILLUSTRATION FOR $\alpha$ = 0.4

When $\alpha$ takes the values 0.7, 0.8, 0.9, 1 the corresponding row sum of the membership grades are as follows:

$$\mu_x (R_i)_{i=1,2,\ldots 6} = \begin{cases} 0.521, 0.09, 0, 1, 0.478, 0.9 \text{ when } \alpha = 0.7 \\ 0.878, 0.092, 0, 1, 0.475, 0.994 \text{ when } \alpha = 0.8 \\ 0.892, 0.081, 0, 1, 0.483, 0.989 \text{ when } \alpha = 0.9 \\ 0.91, 0.07, 0, 1, 0.484, 0.99 \text{ when } \alpha = 1 \end{cases}.$$

Here the fourth row sum $R_4$ is getting the highest membership grade that is the membership grade is 1 and the sixth row sum

$R_6$ is getting the next highest membership grade. Graphical illustration for $\alpha$ = 0.7, 0.8 is shown in graphs given below:

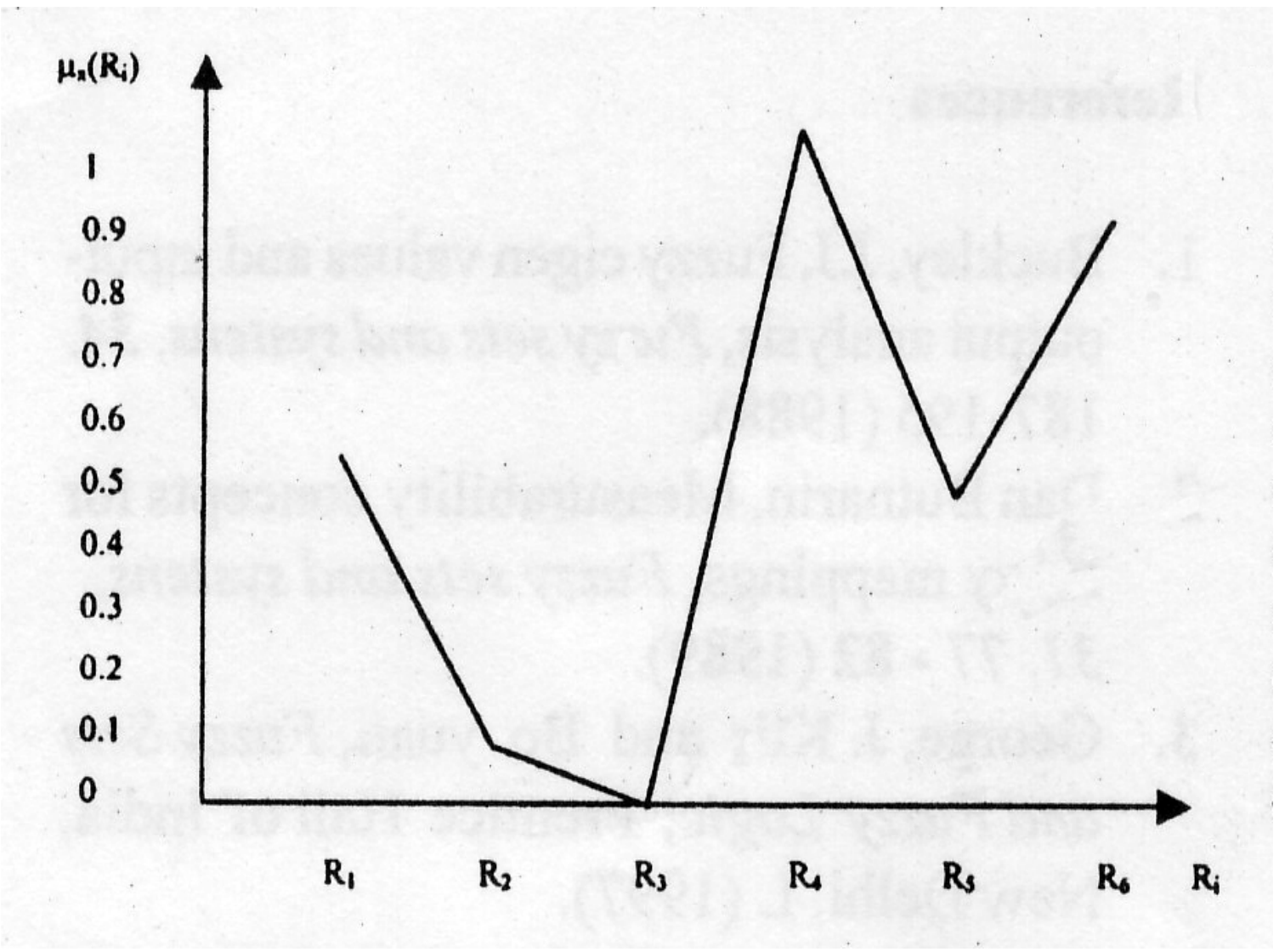

GRAPH 2.3.2.7: GRAPHICAL ILLUSTRATION FOR $\alpha$ = 0.7

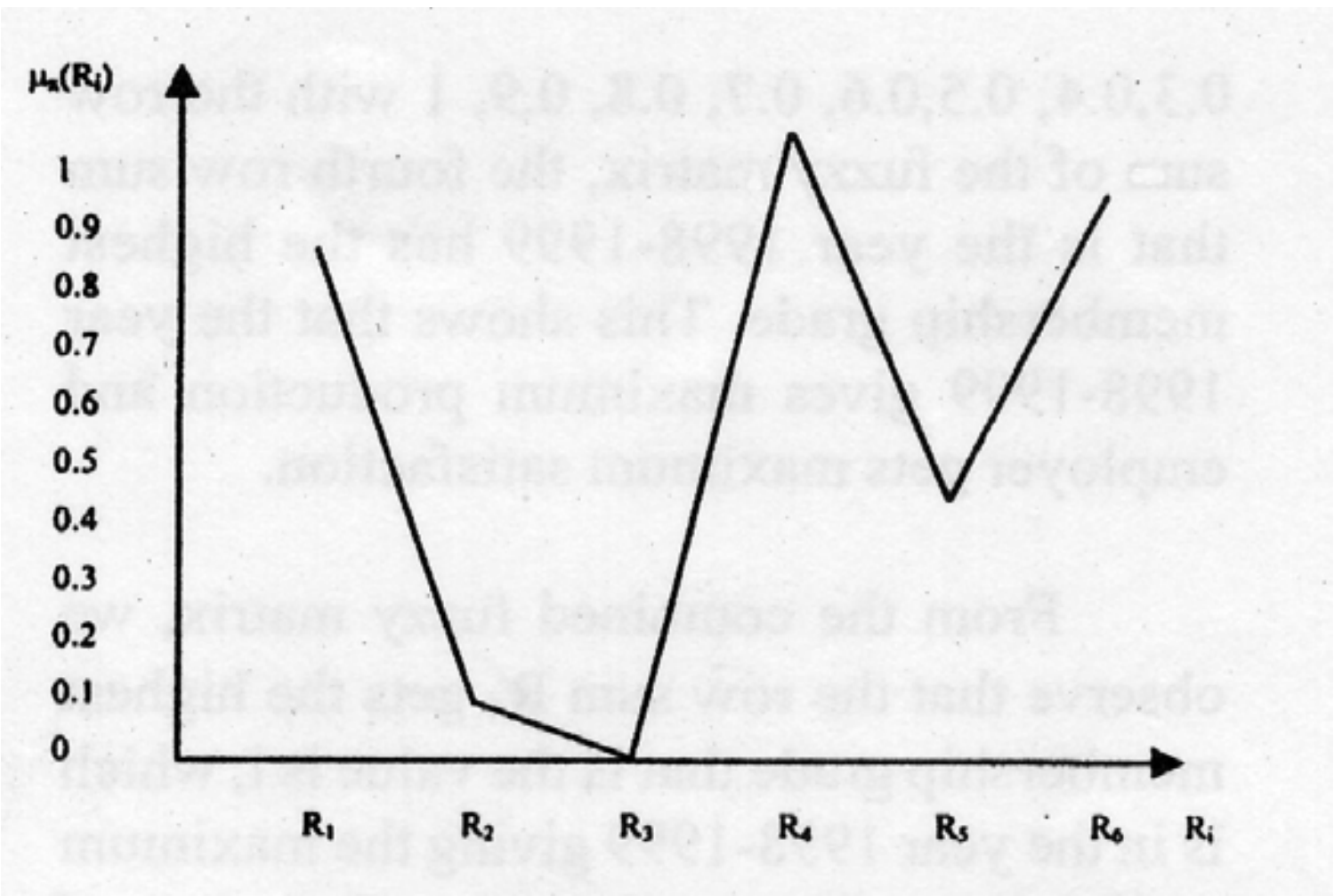

GRAPH 2.3.2.8: GRAPHICAL ILLUSTRATION FOR $\alpha$ = 0.8

Graphical illustration for $\alpha = 0.9$ and 1 is shown in graph 2.3.2.9 and graph 2.3.2.10 which is given below:

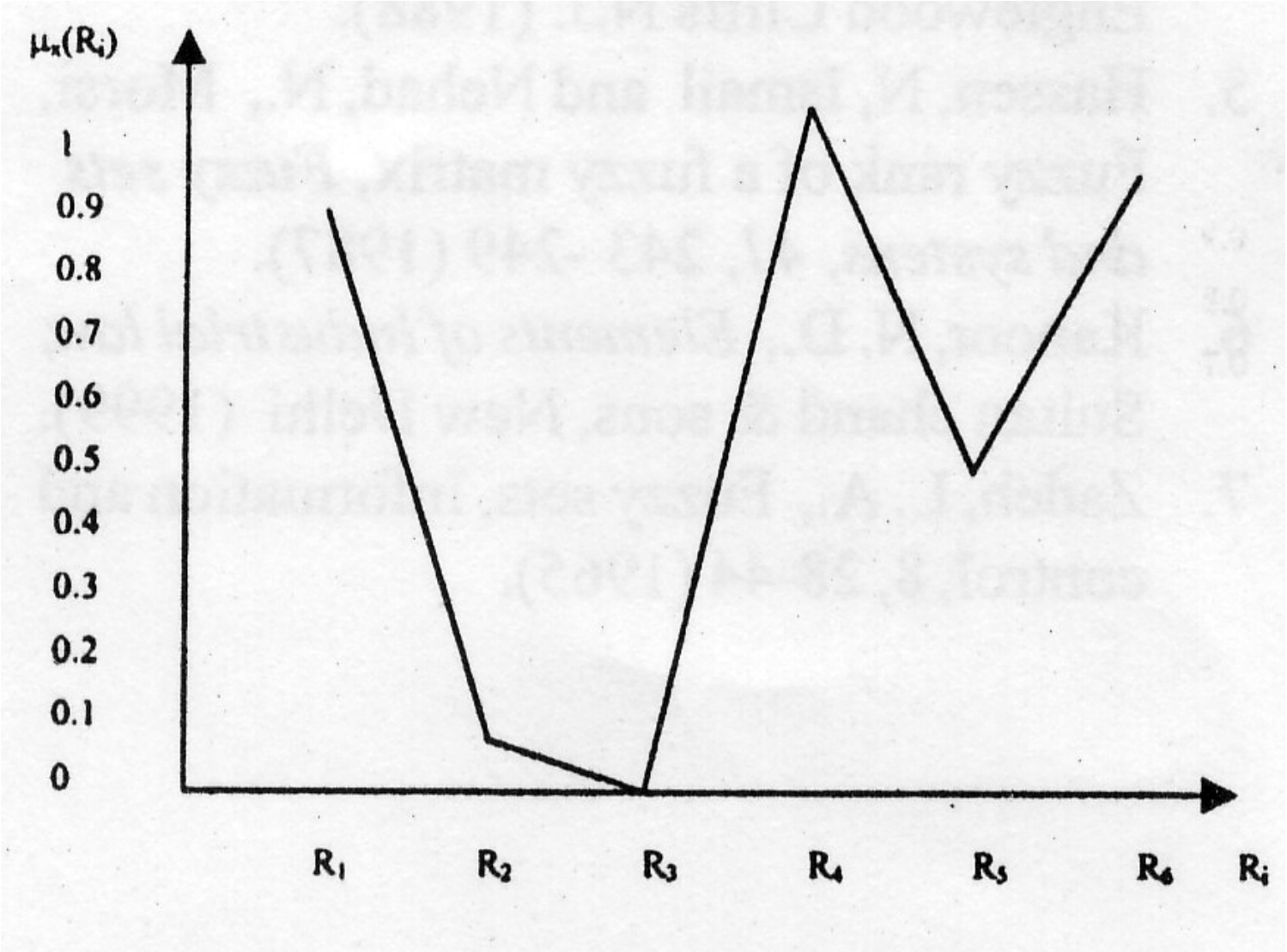

GRAPH 2.3.2.9: GRAPHICAL ILLUSTRATION FOR $\alpha$ = 0.9

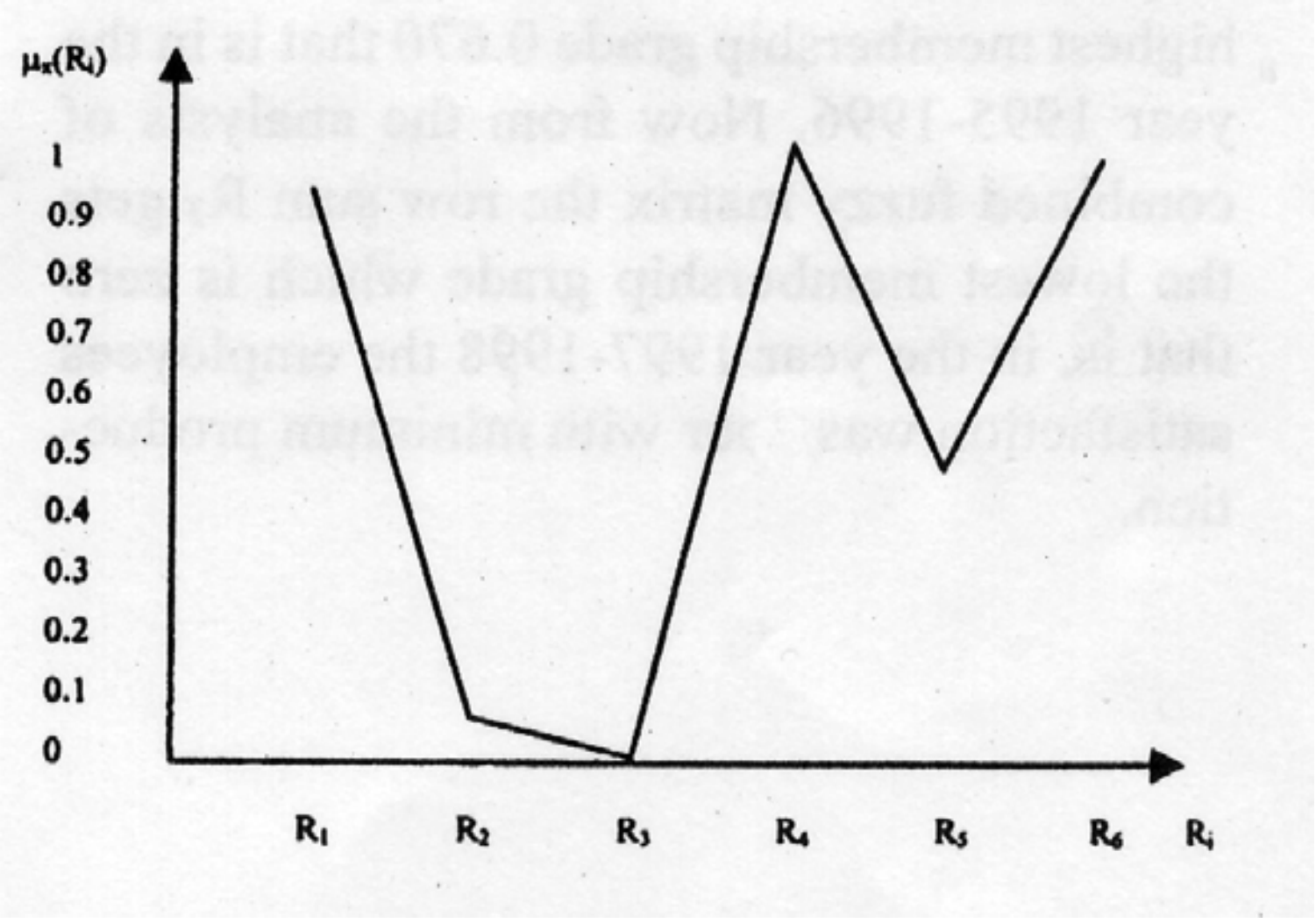

GRAPH 2.3.2.10: GRAPHICAL ILLUSTRATION FOR $\alpha$ =1

The combined fuzzy matrix for all the values of $\alpha \in [0, 1]$ is given below:

| | $X_1$ | $X_2$ | $X_3$ | $X_4$ | $X_5$ | $X_6$ | $X_7$ | $X_8$ |
|---|---|---|---|---|---|---|---|---|
| $Y_1$ | 0 | 7.596 | 10 | 0 | 0 | 10 | 3.689 | 10 |
| $Y_2$ | 10 | 10 | 9.702 | 0 | 0 | 0 | 0 | 4.533 |
| $Y_3$ | 10 | 10 | 3.116 | 0 | 9.671 | 0 | 0 | 0 |
| $Y_4$ | 0 | 0 | 7.724 | 10 | 0.399 | 9.99 | 8.531 | 8.819 |
| $Y_5$ | 0 | 0 | 5.713 | 10 | 7.480 | 7.174 | 10 | 0 |
| $Y_6$ | 1.93 | 0 | 0 | 4.345 | 9.98 | 9.98 | 9.98 | 3.413 |

The corresponding row sum of the above matrix is given by

| | | | | |
|---|---|---|---|---|
| $R_1$ | = | The first row sum | = | 41.285 |
| $R_2$ | = | The second row sum | = | 34.235 |
| $R_3$ | = | The third row sum | = | 32.787 |
| $R_4$ | = | The fourth row sum | = | 45.463 |
| $R_5$ | = | The fifth row sum | = | 40.367 |
| $R_6$ | = | The sixth row sum | = | 39.628. |

Using the membership function the corresponding values of $\mu_x(R_i)$ is 0.679, 0.144, 0, 1, 0.598, 0.539 respectively where i = 1, 2, 3, 4, 5, 6. Graphical illustration for the values of $\alpha \in [0, 1]$ is shown in graph 2.3.2.11.

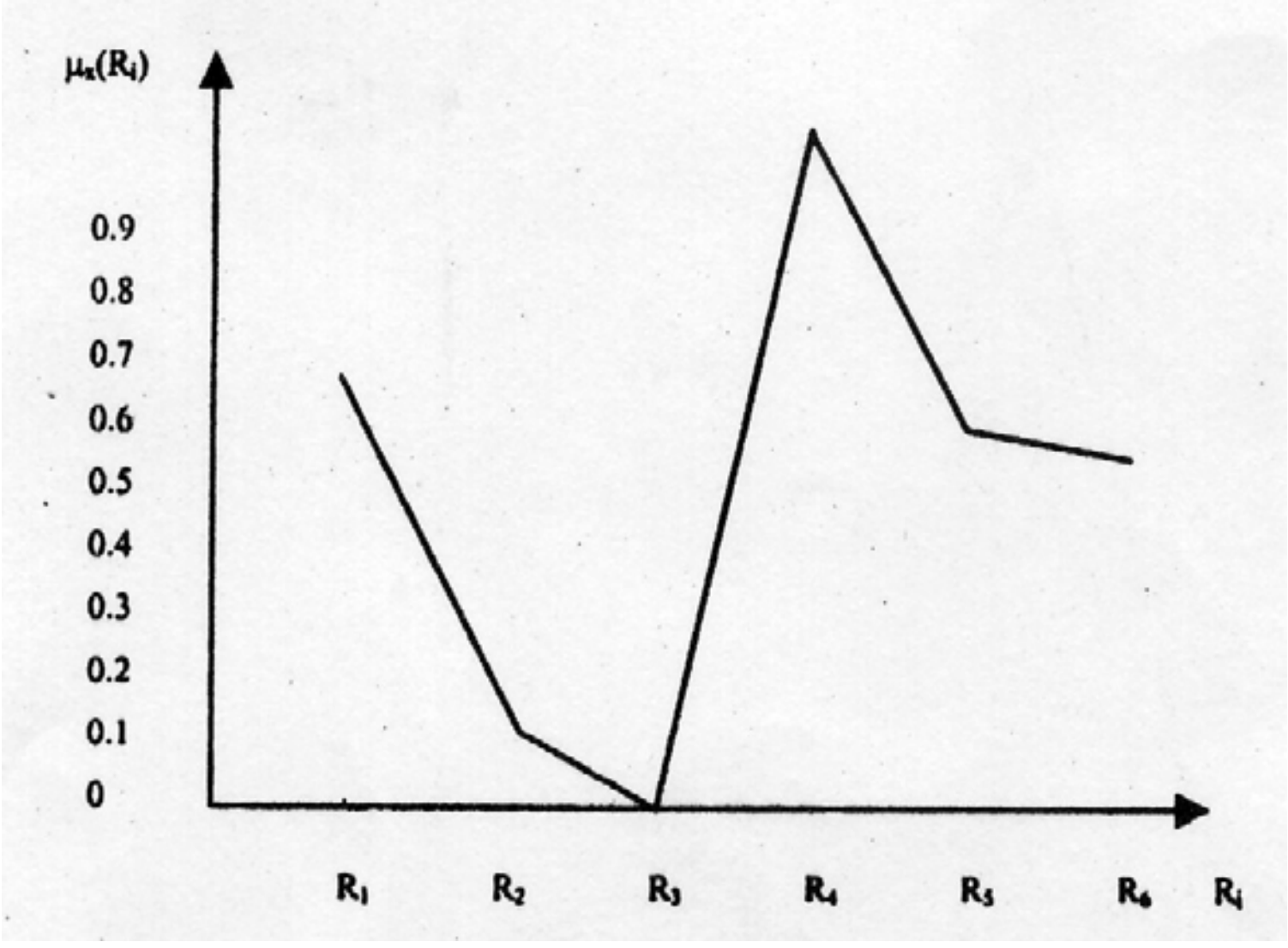

GRAPH 2.3.2.11: GRAPHICAL ILLUSTRATION FOR $\alpha \in [0, 1]$

Here the fourth row sum $R_4$ is getting the highest membership grade that is the membership grade is 1, first row sum $R_1$ is getting the next highest membership grade that is the membership grade is 0.670 and the row sum $R_3$ gets the lowest membership grade which is zero.

We now proceed onto to give the conclusions.

Here we introduce a fuzzy relational matrix model to the problem of employees and employers in cement industries and this finds best form of relationship between the employees and the employers. The raw data is obtained by the expert's opinion from the cement industrialists as well as feelings of the employees. Using the data from the cement industry, we carry out the analysis of the data via fuzzy relational matrix model and found the best form of relationship between employees and employers. This has been done in two ways.

In the first stage the data was converted into a relational map and the relational matrix was obtained form the relational map, and using this the average relational matrix is obtained.

In the second stage the average relational matrix is converted into fuzzy relational matrix and using different parameter $\alpha$, $\alpha \in [0, 1]$ the best form of relationship is obtained. This is also explicitly described by the graphs. Here we consider eight attributes of the employees and six effects of employers. The production levels of each year are treated as rows and the various attributes of employee are treated as columns in the fuzzy matrix.

The highest membership grade gives the best form of relationship between employers and employees, which maximizes the production level of the cement industries. Using the different parameter $\alpha$ = 0.1, 0.2, 0.3, 0.4, 0.5, 0.6, 0.7, 0.8, 0.9, 1 with the row sum of the fuzzy matrix, the fourth row sum that is the year 1998-1999 gives maximum production and employer gets maximum satisfaction.

From the combined fuzzy matrix, we observe that the row sum $R_4$ gets the highest membership grade that is the value is 1, which is in the year 1998-1999 giving the maximum production with maximum satisfaction of employees and the row sum $R_1$ gets the next highest membership grade 0.670 that is in the year 1995-1996. Now from the analysis of combined fuzzy matrix the row sum $R_3$ gets the lowest membership grade which is zero that is, in the year 1997-1998 the employees satisfaction was poor with minimum production.

For more please refer [237].

Thus we have seen how best the FRMs can be used when we do not have any data but only the opinions and when we have the data using the past experience how best we can give the predictions. Here we give some justifications to state why the use of FRMs are sometimes better than the FCMs.

The first marked difference between FCMs and FRMs is that FCMs cannot directly give the effect of one group on the other. But FRMs can give the effect of one group on the other group and vice versa .

FCMs cannot give any benefit when the nodes or causalities are mutually exclusive ones. But in the case of FRMs since we divide them into two groups and relational maps are sent from one group to other, it gives the maximum benefit.

FRMs give the direct effect of one node from space to other node or nodes of the other space very precisely. We see also in case of FRMs when more than one node is in the on state the hidden pattern ends in a limit cycle and when only one node is on the hidden pattern happens to be a fixed points.

Another positive point about the FRMs is when the data can be divided disjointly the size of the matrix is considerably and significantly reduced.

For if in FCMs say we have just 12 nodes then we have a 12 $\times$ 12 matrix with 144 entries. But if 12 nodes are divided into 7 and 5 we get only a 7 $\times$ 5 matrix with 35 entries. Likewise if it is a 8 $\times$ 4 matrix we deal only with 32 entries thus a three digit number is reduced to a two digit number.

## 2.4 Introduction to Bidirectional Associative Memories (BAM) Model and their Application

This section has two subsection. The first subsection just gives the functioning of the BAM model. Section two illustrates the functioning of the BAM model in the study of the cause vulnerability to HIV/AIDS and factors for migrations.

We live in a world of marvelous complexity and variety, a world where events never repeat exactly. Even though events are never exactly the same they are also not completely different. There is a threat of continuity, similarity and predictability that allows us to generalize often correctly from past experience to future events. Neural networks or neuro-computing or brain like computation is based on the wistful hope that we can reproduce at least some of the flexibility and power of the human brain by artificial brains. Neural networks consists of many simple computing elements generally simple non linear summing junctions connected together by connections of varying strength a gross abstraction of the brain which consists of very large number of far more complex neurons connected together with far more complex and far more structured couplings, neural networks architecture cover a wide range. In one sense every computer is a neural net, because we can view traditional digital logic as constructed from inter connected McCullouch-Pitts neurons. McCullouch-Pitts neurons were proposed in 1943 as models of biological neurons and arranged in networks for a specific purpose of computing logic functions. The problems where artificial neural networks have the most promise are those with a real-world flavor: medical research, signal processing, sociological problems etc.

Neural networks helps to solve these problems with natural mechanisms of generalizations. To over-simplify, suppose we represent an object in a network as a pattern of activation of several units. If a unit or two responds incorrectly the overall pattern stays pretty much of the same, and the network still respond correctly to stimuli when neural networks operate similar inputs naturally produce similar outputs. Most real world perceptual problems have this structure of input-output continuity. The prototype model provides a model for human categorization with a great deal of psychological support. The computational strategy leads to some curious human psychology. For instance in United States people imagine a prototype bird that looks somewhat like a sparrow or a robin. So they learn to judge ‘penguins’ or “ostriches” as “bad” birds because these birds do not resemble the prototype bird even though they are birds “Badness” shows up in a number of ways.

Neural networks naturally develop this kind of category structure. The problems that neural networks solved well and solved poorly were those where human showed comparable strengths and weaknesses in their cognitive computations”. For this reason until quite recently most of the study of neural networks has been carried out by psychologists and cognitive scientists who sought models of human cognitive function. Neural networks deal with uncertainty as humans do, not by deliberate design but as a by product of their parallel distributed structure. Because general statements about both human psychology and the structure of the world embed so deeply in both neural networks and fuzzy systems, it is very appropriate to study the psychological effects of HIV/AIDS patients and their influence on public using this theory. Like social customs these assumptions are obvious only if you grew up with them. Both neural networks and fuzzy systems break with the historical tradition prominent in western thought and we can precisely and unambiguously characterize the world, divide into two categories and then manipulate these descriptions according to precise and formal rules. Huang Po, a Buddhist teacher of the ninth century observed that “To make use of your minds to think conceptually is to leave the substance and attach yourself to form”, and “from discrimination between this and that a host of demons blazes forth”.

## 2.4.1 Some Basic Concepts of BAM

Now we go forth to describe the mathematical structure of the Bidirectional Associative Memories (BAM) model. Neural networks recognize ill defined problems without an explicit set of rules. Neurons behave like functions, neurons transduce an unbounded input activation x(t) at time t into a bounded output signal S(x(t)) i.e. Neuronal activations change with time.

Artificial neural networks consists of numerous simple processing units or neurons which can be trained to estimate sampled functions when we do not know the form of the functions. A group of neurons form a field. Neural networks contain many field of neurons. In our text $F_x$ will denote a neuron field, which contains n neurons, and $F_y$ denotes a neuron field, which contains p neurons. The neuronal dynamical system is described by a system of first order differential equations that govern the time-evolution of the neuronal activations or which can be called also as membrane potentials.

$$\dot{x}_i = g_i (X, Y, ...)$$
$$\dot{y}_j = h_j (X, Y, ...)$$

where $\dot{x}_i$ and $\dot{y}_j$ denote respectively the activation time function of the $i^{th}$ neuron in $F_X$ and the $j^{th}$ neuron in $F_Y$. The over dot denotes time differentiation, $g_i$ and $h_j$ are some functions of X, Y, ... where $X(t) = (x_1(t), ... , x_n(t))$ and $Y(t) = (y_1(t), ... , y_p(t))$ define the state of the neuronal dynamical system at time t. The passive decay model is the simplest activation model, where in the absence of the external stimuli, the activation decays in its resting value

$$\dot{x}_i = x_i$$

and $$\dot{y}_j = y_j$$

The passive decay rate $A_i > 0$ scales the rate of passive decay to the membranes resting potentials $\dot{x}_i = –A_ix_i$. The default rate is $A_i = 1$, i.e. $\dot{x}_i = –A_ix_i$. The membrane time constant $C_i > 0$ scales the time variables of the activation dynamical system. The default time constant is $C_i = 1$. Thus $C_i \dot{x}_i = –A_ix_i$.

The membrane resting potential $P_i$ is defined as the activation value to which the membrane potential equilibrates in the absence of external inputs. The resting potential is an additive constant and its default value is zero. It need not be positive.

$$P_i = C_i \dot{x}_i + A_i x_i$$
$$I_i = \dot{x}_i + x_i$$

is called the external input of the system. Neurons do not compute alone. Neurons modify their state activations with external input and with feed back from one another. Now, how do we transfer all these actions of neurons activated by inputs their resting potential etc. mathematically. We do this using what are called synaptic connection matrices.

Let us suppose that the field $F_X$ with n neurons is synaptically connected to the field $F_Y$ of p neurons. Let $m_{ij}$ be a synapse where the axon from the $i^{th}$ neuron in $F_X$ terminates. $M_{ij}$ can be positive, negative or zero. The synaptic matrix M is a n by p matrix of real numbers whose entries are the synaptic efficacies $m_{ij}$.

The matrix M describes the forward projections from the neuronal field $F_X$ to the neuronal field $F_Y$. Similarly a p by n synaptic matrix N describes the backward projections from $F_Y$ to $F_X$. Unidirectional networks occur when a neuron field synaptically intra connects to itself. The matrix M be a n by n square matrix. A Bidirectional network occur if $M = N^T$ and $N = M^T$. To describe this synaptic connection matrix more simply, suppose the n neurons in the field $F_X$ synaptically connect to the p-neurons in field $F_Y$. Imagine an axon from the $i^{th}$ neuron in $F_X$ that terminates in a synapse $m_{ij}$, that about the $j^{th}$ neuron in $F_Y$. We assume that the real number $m_{ij}$ summarizes the synapse and that $m_{ij}$ changes so slowly relative to activation fluctuations that is constant.

Thus we assume no learning if $m_{ij} = 0$ for all t. The synaptic value $m_{ij}$ might represent the average rate of release of a neuro-transmitter such as norepinephrine. So, as a rate, $m_{ij}$ can be positive, negative or zero.

When the activation dynamics of the neuronal fields $F_X$ and $F_Y$ lead to the overall stable behaviour the bidirectional networks are called as Bidirectional Associative Memories (BAM). As in the analysis of the HIV/AIDS patients relative to the migrancy we state that the BAM model studied presently and predicting the future after a span of 5 or 10 years may not be the same.

For the system would have reached stability and after the lapse of this time period the activation neurons under investigations and which are going to measure the model would be entirely different.

Thus from now onwards more than the uneducated poor the educated rich and the middle class will be the victims of HIV/AIDS. So for this study presently carried out can only give how migration has affected the life style of poor labourer and had led them to be victims of HIV/AIDS.

Further not only a Bidirectional network leads to BAM also a unidirectional network defines a BAM if M is symmetric i.e. $M = M^T$. We in our analysis mainly use BAM which are bidirectional networks. However we may also use unidirectional BAM chiefly depending on the problems under investigations. We briefly describe the BAM model more technically and mathematically.

An additive activation model is defined by a system of n + p coupled first order differential equations that inter connects the fields $F_X$ and $F_Y$ through the constant synaptic matrices M and N.

$$x_i = -A_i x_i + \sum_{j=1}^{p} S_j(y_j) n_{ji} + I_i \qquad (2.4.1.1)$$

$$y_i = -A_j y_j + \sum_{i=1}^{n} S_i(x_i) m_{ij} + J_j \qquad (2.4.1.2)$$

$S_i(x_i)$ and $S_j(y_j)$ denote respectively the signal function of the $i^{th}$ neuron in the field $F_X$ and the signal function of the $j^{th}$ neuron in the field $F_Y$.

Discrete additive activation models correspond to neurons with threshold signal functions.

The neurons can assume only two values ON and OFF. ON represents the signal +1, OFF represents 0 or – 1 (– 1 when the representation is bipolar). Additive bivalent models describe asynchronous and stochastic behaviour.

At each moment each neuron can randomly decide whether to change state or whether to emit a new signal given its current activation. The Bidirectional Associative Memory or BAM is a non adaptive additive bivalent neural network. In neural literature the discrete version of the equation (2.4.1.1) and (2.4.1.2) are often referred to as BAMs.

A discrete additive BAM with threshold signal functions arbitrary thresholds inputs an arbitrary but a constant synaptic connection matrix M and discrete time steps K are defined by the equations

$$x_i^{k+1} = \sum_{j=1}^{p} S_j(y_j^k) m_{ij} + I_i \qquad (2.4.1.3)$$

$$y_j^{k+1} = \sum_{i=1}^{n} S_i\left(x_i^k\right) m_{ij} + J_j \qquad (2.4.1.4)$$

where $m_{ij} \in M$ and $S_i$ and $S_j$ are signal functions. They represent binary or bipolar threshold functions. For arbitrary real valued thresholds $U = (U_1, ..., U_n)$ for $F_X$ neurons and $V = (V_1, ..., V_P)$ for $F_Y$ neurons the threshold binary signal functions corresponds to

$$S_i(x_i^k) = \begin{cases} 1 \quad \text{if} \qquad x_i^k > U_i \\ S_i(x_i^{k-1}) \text{ if } x_i^k = U_i \\ 0 \quad \text{if} \qquad x_i^k < U_i \end{cases} \qquad (2.4.1.5)$$

and

$$S_j(x_j^k) = \begin{cases} 1 \quad \text{if} \qquad y_j^k > V_j \\ S_j(y_j^{k-1}) \text{ if } y_j^k = V_j \\ 0 \quad \text{if} \qquad y_j^k < V_j \end{cases} \qquad (2.4.1.6)$$

The bipolar version of these equations yield the signal value –1 when $x_i < U_i$ or when $y_j < V_j$. The bivalent signal functions allow us to model complex asynchronous state change patterns. At any moment different neurons can decide whether to compare their activation to their threshold. At each moment any of the 2n subsets of $F_X$ neurons or 2p subsets of the $F_Y$ neurons can decide to change state. Each neuron may randomly decide whether to check the threshold conditions in the equations (2.4.1.5) and (2.4.1.6). At each moment each neuron defines a random variable that can assume the value ON(+1) or OFF(0 or -1). The network is often assumed to be deterministic and state changes are synchronous i.e. an entire field of neurons is updated at a time. In case of simple asynchrony only one neuron makes a state change decision at a time. When the subsets represent the entire fields $F_X$ and $F_Y$ synchronous state change results.

In a real life problem the entries of the constant synaptic matrix M depends upon the investigator's feelings. The synaptic matrix is given a weightage according to their feelings. If $x \in F_X$ and $y \in F_Y$ the forward projections from $F_X$ to $F_Y$ is defined by the matrix M. $\{F(x_i, y_j)\} = (m_{ij}) = M$, $1 \leq i \leq n$, $1 \leq j \leq p$.

The backward projections is defined by the matrix $M^T$. $\{F(y_i, x_i)\} = (m_{ji}) = M^T$, $1 \leq i \leq n$, $1 \leq j \leq p$. It is not always true that the backward projections from $F_Y$ to $F_X$ is defined by the matrix $M^T$.

Now we just recollect the notion of bidirectional stability. All BAM state changes lead to fixed point stability. The property holds for synchronous as well as asynchronous state changes. A BAM system $(F_X, F_Y, M)$ is bidirectionally stable if all inputs converge to fixed point equilibria. Bidirectional stability is a dynamic equilibrium. The same signal information flows back and forth in a bidirectional fixed point. Let us suppose that A denotes a binary n-vector and B denotes a binary p-vector. Let A be the initial input to the BAM system. Then the BAM equilibrates to a bidirectional fixed point $(A_f, B_f)$ as

$$\begin{array}{lllll} A & \rightarrow & M & \rightarrow & B \\ A' & \leftarrow & M^T & \leftarrow & B \\ A' & \rightarrow & M & \rightarrow & B' \\ A'' & \leftarrow & M^T & \leftarrow & B' \quad \text{etc.} \\ A_f & \rightarrow & M & \rightarrow & B_f \\ A_f & \leftarrow & M^T & \leftarrow & B_f \quad \text{etc.} \end{array}$$

where A', A'', ... and B', B'', ... represents intermediate or transient signal state vectors between respectively A and $A_f$ and B and $B_f$. The fixed point of a Bidirectional system is time dependent.

The fixed point for the initial input vectors can be attained at different times. Based on the synaptic matrix M which is developed by the investigators feelings the time at which bidirectional stability is attained also varies accordingly.

### 2.4.2 Use of BAM Model to Study the Cause of Vulnerability to HIV/AIDS and Factors for Migration

Now the object is to study the levels of knowledge and awareness relating to STD/HIV/AIDS existing among the migrant labourers in Tamil Nadu; and to understand the attitude, risk behavior and promiscuous sexual practice of migrant labourers.

***Example 2.4.2.1:*** This study was mainly motivated from the data collected by us of the 60 HIV/AIDS infected persons who belonged to the category that comes to be defined as migrant labourers. Almost all of them were natives of (remote) villages and had migrated to the city, typically, "in search of jobs", or because of caste and communal violence.

We have noticed how, starting from small villages with hopes and dreams these people had set out to the city, only to succumb to various temptations and finally all their dreams turned into horrid nightmares.

Our research includes probing into areas like: patterns and history of migration work, vulnerabilities and risk exposure in an alien surrounding, 'new' sexual practices or attitudes and discrimination, effect of displacement, coping mechanism etc. We also study the new economic policies of liberalization and globalization and how this has affected people to lose their traditional livelihood and sources of local employment, forcing them into migration.

Our study has been conducted among this informal sector mainly because migrant labourers are more vulnerable to HIV/AIDS infection, when compared to the local population for reasons which include easy money, poverty, powerlessness, inaccessibility to health services, unstable life-style such that insecurity, in jobs, lack of skills, alienation from hometown, lack of community bondage. Moreover, migrant labourers are also not organized into trade unions, as a result of which, they are made victims of horrendous exploitation: they are paid less than minimum wages, they don't receive legal protection, they are unaware of worker's rights issues and essentially lack stability.

Their work periods are rarely permanent; they work as short-term unskilled/ semi-skilled contract labourers or as daily wagers; but they earn well for a day and spend it badly without any social binding or savings or investing on their family members or children. A vast majority around 65% of those interviewed were essentially also part of the 'mobile' population, which was wrapped not only in a single migration from native village to metropolitan city, but also involved in jobs like driving trucks, taxis, etc. which gave them increased mobility. We have also analyzed the patients' feelings about the outreach and intervention programs related to HIV/AIDS and we have sought to comprehend the patterns of marginalization that has increased the predisposition of migrants to HIV and other infectious diseases.

A linguistic questionnaire which was drafted and interviews were conducted for 60 HIV/AIDS patients from the hospitals was the main data used in this analysis. Then the questionnaire was transformed into a Bidirectional Associative Memory (BAM) model [233].

Our sample group consisted of HIV infected migrant labourers whose age group ranged between 20-58 and they were involved in a variety of deregulated labour such as transport or truck drivers, construction labourers, daily wagers or employed in hotels or eateries. We have also investigated the feminization of migration and how women were vulnerable to HIV/AIDS only because of their partners.

Thus we have derived many notable conclusions and suggestion from our study of the socio-economic and psychological aspects of migrant labourers with reference to HIV/AIDS [233].

## DESCRIPTION OF THE PROBLEM

In view of the linguistic questionnaire we are analyzing the relation among

a. Causes for migrants' vulnerability to HIV/AIDS
b. Factors forcing migration
c. Role of the Government.

We take some subtitles for each of these three main titles.

For the sake of simplicity we are restricted to some major subtitles, which has primarily interested these experts. We use BAM model on the scale [–5, 5]. Here we mention that the analysis can be carried out on any other scale according to the whims and fancies of the investigator.

**A: Causes for Migrant Labourers Vulnerability to HIV/ AIDS**

$A_1$ - No awareness/ no education
$A_2$ - Social status
$A_3$ - No social responsibility and social freedom
$A_4$ - Bad company and addictive habits
$A_5$ - Types of profession
$A_6$ - Cheap availability of CSWs.

**F: Factors forcing people for migration**

$F_1$ - Lack of labour opportunities in their hometown
$F_2$ - Poverty/seeking better status of life
$F_3$ - Mobilization of labour contractors
$F_4$ - Infertility of lands due to implementation of wrong research methodologies/failure of monsoon.

**G: Role of the Government**

$G_1$ - Alternate job if the harvest fails there by stopping migration

| | | |
|---|---|---|
| $G_2$ | - | Awareness clubs in rural areas about HIV/AIDS |
| $G_3$ | - | Construction of hospitals in rural areas with HIV/AIDS counseling cell/ compulsory HIV/AIDS test before marriage |
| $G_4$ | - | Failed to stop the misled agricultural techniques followed recently by farmers |
| $G_5$ | - | No foresight for the government and no precautionary actions taken from the past occurrences. |

Now the experts opinion on the cause for vulnerability to HIV/AIDS and factors for migration are given in the following:

Taking the neuronal field $F_X$ as the attributes connected with the causes of vulnerability resulting in HIV/AIDS and the neuronal field $F_Y$ is taken as factors forcing migration.

The $6 \times 4$ matrix $M_1$ represents the forward synaptic projections from the neuronal field $F_X$ to the neuronal field $F_Y$.

The $4 \times 6$ matrix $M_1^T$ represents the backward synaptic projections $F_X$ to $F_Y$.

Now, taking $A_1, A_2, \ldots, A_6$ along the rows and $F_1, \ldots, F_4$ along the columns we get the synaptic connection matrix $M_1$ which is modeled on the scale [–5, 5]

$$M_1 = \begin{bmatrix} 5 & 2 & 4 & 4 \\ 4 & 3 & 5 & 3 \\ -1 & -2 & 4 & 0 \\ 0 & 4 & 2 & 0 \\ 2 & 4 & 3 & 3 \\ 0 & 1 & 2 & 0 \end{bmatrix}$$

Let $X_K$ be the input vector given as (3, 4, –1, –3, –2, 1) at the $K^{th}$ time period. The initial vector is given such that illiteracy, lack of awareness, social status and cheap availability of CSWs have stronger impact over migration. We suppose that all neuronal state change decisions are synchronous.

The binary signal vector

$$S(X_K) = (1\ 1\ 0\ 0\ 0\ 1).$$

From the activation equation

$$S(X_K)M_1 = (9, 6, 11, 7) = Y_{K+1}.$$

From the activation equation

$$S(Y_{K+1}) = (1\ 1\ 1\ 1).$$

Now

$$S(Y_{K+1})M_1^T = (15, 15, 1, 6, 12, 3) = X_{K+2}.$$

From the activation equation,

$$S(X_{K+2}) = (1\ 1\ 1\ 1\ 1\ 1),$$
$$S(X_{K+2})M_1 = (10, 12, 20, 10) = Y_{K+2.4.}$$
$$S(Y_{K+3}) = (1\ 1\ 1\ 1).$$

Thus the binary pair {(1 1 1 1 1 1), (1 1 1 1)} represents a fixed point of the dynamical system. Equilibrium of the system has occurred at the time K + 2, when the starting time was K. Thus this fixed point suggests that illiteracy with unawareness, social status and cheap availability of CSW lead to the patients remaining or becoming socially free with no social responsibility, having all addictive habits and bad company which directly depends on the types of profession they choose.

On the other hand, all the factors of migration also come to on state. Suppose we take only the on state that the availability of CSWs at very cheap rates is in the on state. Say at the $K^{th}$ time we have

$$
\begin{aligned}
P_K &= (0\ 0\ 0\ 0\ 0\ 4), \\
S(P_K) &= (0\ 0\ 0\ 0\ 0\ 1), \\
S(P_K)M_1 &= (0\ 1\ 2\ 0) \\
&= Q_{K+1}.
\end{aligned}
$$

$$
\begin{aligned}
S(Q_{K+1}) &= (0\ 1\ 1\ 0) \\
S(Q_{K+1})M_1^T &= (6, 8, 2, 6, 7, 3) \\
&= P_{K+2}
\end{aligned}
$$

$$
\begin{aligned}
S(P_{K+2}) &= (1\ 1\ 1\ 1\ 1\ 1) \\
S(P_{K+2})M_1 &= (10, 12, 20, 10) \\
&= P_{K+3} \\
S(P_{K+3}) &= (1\ 1\ 1\ 1).
\end{aligned}
$$

Thus the binary pair {(1 1 1 1 1 1), (1 1 1 1)} represents a fixed point. Thus in the dynamical system given by the expert even if only the cheap availability of the CSW is in the on state, all the other states become on, i.e. they are unaware of the disease, their type of profession, they have bad company and addictive habits, they have no social responsibility and no social fear.

Thus one of the major causes for the spread of HIV/AIDS is the cheap availability of CSWs which is mathematically confirmed from our study. Several other states of vectors have been worked by us for deriving the conclusions.

We now give the experts opinion on the role of government and causes for vulnerability of HIV/AIDS.

Taking the neuronal field $F_X$ as the role of Government and the neuronal field $F_Y$ as the attributes connected with the causes of vulnerabilities resulting in HIV/ AIDS.

The $6 \times 5$ matrix $M_2$ represents the forward synaptic projections from the neuronal field $F_X$ to the neuronal field $F_Y$.

The $5 \times 6$ matrix $M_2^T$ represents the backward synaptic projection $F_X$ to $F_Y$. Now, taking $G_1, G_2, \ldots, G_5$ along the rows and $A_1$, $A_2$, …, $A_6$ along the columns we get the synaptic connection matrix $M_2$ in the scale [–5, 5] is as follows:

$$M_2 = \begin{bmatrix} 3 & 4 & -2 & 0 & -1 & 5 \\ 5 & 4 & 3 & -1 & 0 & 4 \\ 1 & 3 & 0 & 1 & 4 & 2 \\ 2 & 3 & -2 & -3 & 0 & 3 \\ 3 & 2 & 0 & 3 & 1 & 4 \end{bmatrix}$$

Let $X_K$ be the input vector (–3, 4, –2, –1, 3) at the $K^{th}$ instant. The initial vector is given such that Awareness clubs in the rural villages and the Government's inability in foreseeing the conflicts have a stronger impact over the vulnerability of HIV/AIDS. We suppose that all neuronal state change decisions are synchronous.

The binary signal vector

$$S(X_K) = (0\ 1\ 0\ 0\ 1).$$

From the activation equation

$$S(X_K)M_2 = (8\ 6\ 3\ 2\ 1\ 8) = Y_{K+1}.$$

Now,

$$S(Y_{K+1}) = (1\ 1\ 1\ 1\ 1),$$
$$S(Y_{K+1})\ M_2^T = (9, 15, 11, 3, 13) = X_{K+2}.$$

Now,

$$S(X_{K+2}) = (1\ 1\ 1\ 1\ 1),$$
$$S(X_{K+2})M_2 = (14, 8, -1, 0, 4, 18) = Y_{K+2.4.}$$

Now,

$$S(Y_{K+3}) = (1\ 1\ 0\ 1\ 1\ 1),$$
$$S(Y_{K+3})\ M_2^T = (11, 12, 11, 5, 13) = X_{K+4}$$

Thus

$$S(X_{K+4}) = (1\ 1\ 1\ 1\ 1) = X_{K+2}$$

and

$$S(Y_{K+5}) = Y_{K+3}\ .$$

The binary pair ((1 1 1 1 1), (1 1 0 1 1 1)) represents a fixed point of the BAM dynamical system. Equilibrium for the system occurs at the time K + 4, when the starting time was K.

This fixed point reveals that the other three conditions cannot be ignored and have its consequences in spreading HIV/ AIDS.

Similarly by taking a vector $Y_K$ one can derive conclusions based upon the nature of $Y_K$.

We now give experts opinion on the factors of migration and the role of government.

Taking the neuronal field $F_X$ as the attributes connected with the factors of migration and the neuronal field $F_Y$ is the role of Government. The $4 \times 5$ matrix $M_3$ represents the forward synaptic projections from the neuronal field $F_X$ to the neuronal field $F_Y$.

The $5 \times 4$ matrix $M_3^T$ represents the backward synaptic projections $F_X$ to $F_Y$.

Now, taking $F_1, F_2, F_3, F_4$ along the rows and $G_1, G_2, \ldots, G_5$ along the columns we get the synaptic connection matrix $M_3$ as follows:

$$M_3 = \begin{bmatrix} 4 & 0 & 5 & 3 & 4 \\ 3 & -2 & -4 & 4 & 3 \\ 3 & 0 & 4 & -1 & -2 \\ 2 & 1 & 0 & 5 & 4 \end{bmatrix}$$

Let $X_K$ be the input vector (–2, 1, 4, –1) at the time K. The initial vector is given such that poverty and mobilization of labour contractors have a greater impact. We suppose that all neuronal state change decisions are synchronous.

The binary signal vector

$$S(X_K) = (0\ 1\ 1\ 0).$$

From the activation equation

$$S(X_K)\, M_3 = (6, -2, 0, 3, 1\ ) = Y_{K+1}.$$

Now

$$S(Y_{K+1}\ ) = (1\ 0\ 1\ 1\ 1),$$
$$S(Y_{K+1})\, M_3^T = (16, 6, 4, 11) = X_{K+2}.$$

Now

$$S(X_{K+2}) = (1\ 1\ 1\ 1\ ),$$
$$S(X_{K+2})M_3 = (12, -1, 5, 11, 9) = Y_{K+2.4.}$$

Now

$$S(Y_{K+3}\ ) = (1\ 0\ 1\ 1\ 1)\ ,$$
$$S(Y_{K+3})M_3^T = (16, 6, 9, 11) = X_{K+4}.$$

Thus

$$S(X_{K+4}) = (1\ 1\ 1\ 1) = X_{K+2},$$
$$S(Y_{K+5}) = Y_{K+3}\ .$$

The binary pair {((1 1 1 1), (1 0 1 1 1))} represents a fixed point of the BAM dynamical system. Equilibrium for the system occurs at the time K + 4, when the starting time was K. All the factors point out the failure of the Government in tackling HIV/AIDS. Similarly by taking a vector $Y_K$ one can derive conclusions based upon the nature of $Y_K$.

Thus these illustrations are given only for the sake of making the reader to understand the workings of the fuzzy models.

The calculations of the fixed point using in the BAM-model is given in Appendix 6 where C-program is used to make the computation easy.

***Example 2.4.2.2:*** Now we use the BAM model in the interval [–4 4] for analyzing the same linguistic questionnaire and we keep the main 3 heads as it is and make changes only the subtitles. These were the subtitles.

These were the subtitles proposed by the expert whose opinion is sought.

**A: Causes for migrant labourers vulnerability to HIV/AIDS.**

$A_1$ - No awareness / No education
$A_2$ - Away from the family for weeks
$A_3$ - Social status
$A_4$ - No social responsibility and social freedom
$A_5$ - Bad company and addictive habits
$A_6$ - Types of profession where they cannot easily escape from visit of CSWs
$A_7$ - Cheap availability of CSWs
$A_8$ - No union / support group to channelize their ways of living as saving for future and other monetary benefits from the employer
$A_9$ - No fear of being watched by friends or relatives.

**F: Factors forcing people for migration**

$F_1$ - Lack of labour opportunities in their hometown
$F_2$ - Poverty
$F_3$ - Unemployment
$F_4$ - Mobilization of contract labourers
$F_5$ - Infertility of land failure of agriculture

$F_6$ - Failed Government policies like advent of machinery, no value of small scale industries, weavers problem, match factory problem etc.

**G: Role of Government**

$G_1$ - Lack of awareness clubs in rural areas about HIV/AIDS
$G_2$ - No steps to help agricultural cooli only rich farmers are being helped i.e., landowners alone get the benefit from the government helped
$G_3$ - Failed to stop poor yield
$G_4$ - No alternative job provided for agricultural coolie / weavers
$G_5$ - No proper counselling centres for HIV/AIDS in villages
$G_6$ - Compulsory HIV/AIDS test before marriage
$G_7$ - Government does not question about the marriage age of women in rural areas. 90% of marriage is rural areas with no education takes place when women are just below 16 years in 20% of the cases even before they girls attain puberty they are married).

Now the expert's opinion on the cause of vulnerability to HIV/AIDS factors for migration is given below:

Taking the neuronal field $F_x$ as the attributes connected with the causes of vulnerability resulting in HIV/AIDS and the neuronal field $F_y$ is taken as factors forcing migration.

The $9 \times 6$ matrix $M_1$ represents the forward synaptic projections from the neuronal field $F_x$ to the neuronal field $F_Y$. The $6 \times 9$ matrix $M^T_1$ represents the backward synaptic projections $F_X$ to $F_Y$.

Now taking $V_2, \ldots, V_9$ along the rows and $F_1, \ldots, F_6$ along the columns we get the synaptic connection matrix $M_1$ which is modelled in the scale [– 4, 4].

$$M_1 = \begin{bmatrix} 2 & 0 & 0 & 0 & 0 & 0 \\ 3 & 2 & 2 & 2 & 1 & 3 \\ -2 & 3 & 0 & 2 & 0 & 0 \\ 0 & 0 & 0 & 0 & 0 & 0 \\ 0 & 0 & 0 & 1 & -2 & 0 \\ 4 & 3 & -2 & 3 & 2 & 0 \\ 0 & -1 & 0 & 0 & -2 & 0 \\ 0 & 0 & 0 & 2 & 0 & 0 \\ 0 & -3 & 0 & 0 & -2 & 0 \end{bmatrix}.$$

Let $X_K$ be the input vector given as (3 2 1 –1 0 –2 4 –2 1) at the $K^{th}$ time period. The initial vectors as $V_1$, $V_2$, …, $V_9$. We suppose that all neuronal state change decisions are synchronous.

The binary signal vector

$$S(X_K) = (1\ 1\ 1\ 0\ 0\ 0\ 1\ 0\ 1)$$

From the activation equation

$$\begin{aligned} S(X_K)M_1 &= (3\ 1\ 2\ 4\ {-3}\ 3) \\ &= Y_{K+1} \\ S(Y_{K+1}) &= (1\ 1\ 1\ 1\ 0\ 1) \end{aligned}$$

Now

$$\begin{aligned} S(Y_{K+1})\, M^T_1 &= (2\ 1\ 2\ 3\ 0\ 1\ 8\ {-1}\ 2\ {-3}) \\ &= X_{K+2} \end{aligned}$$

$$\begin{aligned} S(X_{K+2}) &= (1\ 1\ 1\ 0\ 1\ 1\ 0\ 1\ 0) \\ S(X_{K+2})\, M_1 &= (7, 8, 0\ , 10, 13) \\ &= Y_{K+3} \end{aligned}$$

$$\begin{aligned} S(Y_{K+3}) &= (1\ 1\ 0\ 1\ 1\ 1) \\ S(Y_{K+3})\, M^T_1 &= (2, 11, 3, 0, -1, 12, -3, 2, -5) \\ &= X_{K+4} \end{aligned}$$

$$S(X_{K+4}) = (1\ 1\ 1\ 0\ 0\ 1\ 0\ 1\ 0)$$
$$S(X_{K+4})\, M_1 = (7\ 8\ 0\ 9\ 3\ 3) = Y_{K+5}$$

$$S(Y_{K+5}) = (1\ 1\ 0\ 1\ 1\ 1)$$
$$S(Y_{K+5}) M^T_1 = X_{K+6} = X_{K+4.}$$

Thus the binary pair {(1 1 1 0 0 1 0 1 0), (1 1 0 1 1 1)} represents a fixed point of the dynamical system. Equilibrium of the system has occurred at the K + 6 time when the starting point was K.

Thus the fixed point suggests that when cheap availability of CSWs, with no proper awareness and no education, with when they are away from the family with poor social status with no fear of being watched lead the patients to become victims of HIV/AIDS due to their profession and no union/support group to channelize their ways of living. The factors that fuel this are lack of labour opportunity in their hometown, poverty, mobilization of contract labourers, infertility of land failure of agriculture and above all failed government policies like advent of machinery, no value of small scale industries, weavers problem, match factory problem etc.

Using the C program given in the book [233] we can find the equilibrium of the dynamical system for varying values.

Suppose we consider the input vector

$$Y_K = (4\ –1\ 0\ –3\ 0\ 3)$$

i.e. the lack of labour opportunities with failed government polices like advent of machinery which replaces human labour like weaver problem and match factor problems and no value for small scale industries to be vector with large positive values and other nodes being either 0 or negative value

$$S(Y_K) = (1\ 0\ 0\ 0\ 0\ 1)$$
$$S(Y_K)\, M^T_1 = (2\ 6\ –2\ 0\ 0\ 4\ 0\ 0\ 0) = X_{K+1}$$

$S(Y_{K+1}) = (1\ 1\ 0\ 0\ 0\ 1\ 0\ 0\ 0)$

$S(X_{K+1})\ M_1 = (9, 5, 0, 5, 3, 3)$

$= Y_{K+2}$

$S(Y_{K+2}) = (1\ 1\ 0\ 1\ 1\ 1)$

$S(Y_{K+2})\ M^T_1 = (2, 11, 3, 0, –1, 12, –3, 2, –5)$

$= X_{K+3}$

$S(X_{K+3}) = (1\ 1\ 1\ 0\ 0\ 1\ 0\ 1\ 0)$

$S(Y_{K+2})M_1 = (7, 8, 0, 9, 3, 3)$

$= Y_{K+3}$

$S(Y_{K+3}) = (1\ 1\ 0\ 1\ 1\ 1).$

The stability of the BAM is given by the fixed point, i.e., the binary pair {(1 1 1 0 0 1 1 0 1 0), (1 1 0 1 1 1)}. The conclusion as before is derived so if the labourers have no employment in their home town and failed government policies result in following problem n the on state of all vectors expect. $F_{2.4.}$

Now we proceed on to model using the experts opinion on the factors of migration and the role of government.

Taking the neuronal field $F_X$ as the attributes connected with the factors of migration and the neuronal field $F_Y$ is the role of government. The matrix $M_2$ represents the forward synaptic projections from the neuronal field $F_X$ to the neuronal field $F_Y$. The 7 × 6 matrix $M^T_2$ represents the backward synaptic projections $F_X$ to $F_Y$. Now taking the factors $F_1, F_2,\ldots, F_6$ along the rows and $G_1, G_2,\ldots, G_7$ along the columns we get the synaptic connection matrix $M_2$ as follows:

$$M_2 = \begin{array}{c} \\ F_1 \\ F_2 \\ F_3 \\ F_4 \\ F_5 \\ F_6 \end{array} \begin{array}{c} \begin{array}{ccccccc} G_1 & G_2 & G_3 & G_4 & G_5 & G_6 & G_7 \end{array} \\ \begin{bmatrix} 0 & 3 & 4 & 2 & 0 & 0 & 0 \\ 0 & 2 & 3 & 3 & 2 & 0 & 3 \\ 0 & 2 & 3 & 2 & 0 & 0 & 0 \\ 0 & 3 & 2 & 3 & 2 & 0 & 0 \\ 0 & 2 & 2 & 2 & 0 & 0 & 1 \\ 0 & 2 & 1 & 2 & 0 & 0 & 1 \end{bmatrix} \end{array}.$$

Here the expert feels one of the reasons for marrying off their daughter at a very young age is due to poverty. So they marry them off to old men and at times as second wife and so on.

Suppose we take the fit vector as $X_K$ = (-3 2 0 4 –1 –2) at the $K^{th}$ time period. We suppose that all neuronal state change decisions are synchronous.

The binary signal vector

$$S(X_K) = (0\ 1\ 0\ 1\ 0\ 0)$$

From the activation equation

$$S(X_K)\ M_2 = (0\ 5\ 5\ 6\ 4\ 0\ 3) = Y_{K+1}$$

$$S(Y_{K+1}) = (0\ 1\ 1\ 1\ 1\ 0\ 1)$$
$$S(Y_{K+1})\ M^T_2 = (9\ 13\ 7\ 10\ 7\ 6) = X_{K+2}$$

$$S(X_{K+2}) = (1\ 1\ 1\ 1\ 1\ 1)$$
$$S\ (X_{K+2})\ M_2 = (0\ 14\ 15\ 14\ 4\ 0\ 5) = Y_{K+3}$$

$$S\ (Y_{K+3}) = (0\ 1\ 1\ 1\ 1\ 0\ 1).$$

Thus the binary pair {(1 1 1 1 1 1), (0 1 1 1 1 0 1)} represents a fixed point of the dynamical system. Equilibrium of the system has occurred at the time K + 2 when the starting time was K. All attributes come to on state expect the factors lack of awareness clubs in rural areas about HIV/AIDS and compulsory HIV/AIDS lest before marriage are just zero, for they have no impact.

While the factor poverty and mobilization of the contract labourer are at the dominance. Let $Y_K$ = (3, -4, -2, -1 0 –1 2) be the input vector at the $K^{th}$ instant. The initial vector is given such that. Next we consider when in the state vector lack of awareness clubs in rural areas about HIV/AIDS is given the maximum priority followed by marriage age of women is not questioned by the government.

We study the impact of this on migration.

We suppose that all neuronal state change decisions are synchronous.

The binary signal vector.

$$S(Y_K) = (1\ 0\ 0\ 0\ 0\ 0\ 1)$$

From the activation equation

$$S(Y_K)\, M_2^T = (0\ 3\ 0\ 0\ 1\ 1) = X_{K+1}$$

$$S(X_{K+1}) = (0\ 1\ 0\ 0\ 1\ 1)$$
$$S(X_{K+1})\, M_2 = (0\ 6\ 6\ 7\ 2\ 0\ 5) = Y_{K+2}$$

$$S(Y_{K+2}) = (0\ 1\ 1\ 1\ 1\ 0\ 1)$$
$$S(Y_{K+2})\, M^T{}_2 = (9\ 13\ 7\ 10\ 7\ 6) = Y_{K+3}$$

$$S(X_{K+3}) = (1\ 1\ 1\ 1\ 1\ 1)$$
$$S(X_{K+3})\, M_2 = (0\ 14\ 15\ 14\ 4\ 0\ 5) = Y_{K+4}$$

$$S(Y_{K+4}) = (0\ 1\ 1\ 1\ 1\ 0\ 1).$$

Thus the binary pair {(0 1 1 1 1 0 1), (1 1 1 1 1 1)} represents a fixed point of dynamical system. Equilibrium for the system occurs at the time K + 3 when the starting time waste. All nodes except government failure to perform HIV/AIDS test as a compulsory one alone is unaffected and all factors $F_1, \ldots, F_6$ become on so these two nodes promote migration in all forms.

Next we study the experts opinion on the role of government and the causes for vulnerability of HIV/AIDs.

Now taking the neuronal field $F_X$ as the role of government. and the neuronal field $F_Y$ as the attributes connected with the causes of vulnerabilities resulting in HIV/AIDS.

The 7 × 9 matrix $M_3$ represents the forward synaptic projections from the neuronal field $F_X$ to the neuronal field $F_Y$ and the 9 × 7 matrix $M^T_3$ represents the backward synaptic projection $F_X$ to $F_Y$.

Now taking $G_1, \ldots, G_7$ along the rows and $A_1, \ldots, A_9$ along the columns we get the synaptic connection matrix $M_3$ in the scale [–4, 4] as follows:

$$M_3 = \begin{bmatrix} 3 & 0 & -2 & 0 & 0 & 4 & 3 & 4 & 2 \\ -1 & 0 & 3 & 0 & 0 & 0 & 2 & 0 & 0 \\ 0 & 3 & 1 & 0 & 1 & 0 & 0 & 0 & 0 \\ -1 & 2 & 3 & 0 & 0 & 0 & 0 & 1 & 0 \\ 3 & 0 & 1 & 0 & 0 & 1 & 1 & 0 & 0 \\ 0 & 0 & 0 & -1 & -1 & 0 & 0 & -1 & -1 \\ 3 & 0 & 1 & 0 & 0 & 0 & 0 & 0 & 0 \end{bmatrix}.$$

Let $X_K$ be the input vector (2, –1, 3, –2, –1, 0, –1) at the $K^{th}$ instant. The initial vector is given that failed to stop poor yield and lack of awareness clubs in rural area for the uneducated about HIV/AIDS, takes up the primary position. We shall assume all neuronal state change decisions are synchronous.

The binary signal vector

$$S(X_K) = (1\ 0\ 1\ 0\ 0\ 0\ 0)$$

From the activation equation

$$\begin{aligned} S(X_K)\, M_3 &= (3, 3\ –1\ 0\ 1\ 4\ 3\ 4\ 2) \\ &= Y_{K+1} \\ S\,(Y_{K+1}) &= (1\ 1\ 0\ 0\ 1\ 1\ 1\ 1\ 1) \\ S\,(Y_{K+1})\, M^T_3 &= (16, 1, 4, 2, 5\text{-}3\ 3) \\ &= X_{K+2} \end{aligned}$$

$S(X_{K+2})$ = (1 1 1 1 1 0 1)
$S\ (X_{K+2})\ M_3$ = (7, 5, 7 0 1 5 6 5 2)
= $Y_{K+3}$

$S\ (Y_{K+3})$ = (1 1 1 0 1 1 1 1 1)
$S\ (Y_{K+3})\ M_3$ = (5 3 4 –1 0 5 6 3 1)
= $X_{K+4}$

$S\ (X_{K+4})$ = (1 1 1 0 0 1 1 1 1)
$S\ (X_{K+4})\ M^T_3$ = (14 4 4 5 6 –2 4)
= $Y_{K+5}$

$S(Y_{K+5})$ = (1 1 1 1 1 0 1).

Thus the binary pair {(1 1 1 1 1 0 1), (1 1 10 1 1 1 1 1)} or {(1 1 1 1 1 0 1), (1 1 1 0 0 1 1 1 1)} where (1 1 1 1 1 0 1) is a fixed point and (1 1 1 0 1 1 1 1 1) or (1 1 1 0 0 1 1 1 1) occurs.

No social responsibility on the part of the part of the HIV/AIDS affected migrant labourers is 0 all other co-ordinates becomes on their by proving all other cases enlisted makes these labourers vulnerable to HIV/AIDS or the bad habits/bad company and no social responsibly on the parts of HIV/AIDS affected is 0, all other causes for vulnerability is 'on': all factors relating government becomes on except the node ($G_6$) – compulsory test for HIV/AIDS before marriage is off. Thus all factors $G_1$, $G_2$, …, $G_5$ and $G_7$ become on.

Now we consider the input vector $Y_K$ = (3, 2, –1, –2, 0, –3, 1, 0, –2) at the $K^{th}$ instant. The initial vector is given such that No awareness no education takes a priority followed by away from the family for weeks and cheap availability of CSWs.

The binary signal vector

$S\ (Y_K)$ = (1 1 0 0 0 0 1 0 0).

From the activation equation

$S(Y_K)M^T_3$ = (6 1 3 1 4 0 3)
= $X_{K+1}$

$$
\begin{aligned}
S(X_{K+1}) &= (1\ 1\ 1\ 1\ 1\ 0\ 1)\\
S(X_{K+1})\,M_3 &= (7\ 5\ 7\ 0\ 1\ 5\ 6\ 5\ 2)\\
&= Y_{K+2}\\[1em]
S(Y_{K+2}) &= (1\ 1\ 1\ 0\ 1\ 1\ 1\ 1\ 1)\\
S(Y_{K+2})\,M^T_3 &= (14, 4, 5, 5\ 6\ –3\ 4)\\
&= X_{K+3}\\[1em]
S(X_{K+3}) &= (1\ 1\ 1\ 1\ 1\ 0\ 1)\\
&= S(X_{K+1}).
\end{aligned}
$$

Thus the binary pair {(1 1 1 1 1 0 1), (1 1 1 0 1 1 1 1 1)} is a fixed point of the dynamical system. All factors become on in both the classes of attributes except $G_6$ and $A_4$ i.e., compulsory test for HIV/AIDS remains always off and no social responsibility and social freedom on the part of migrants in no other cause for vulnerability to HIV/AIDS.

Now we obtain the equation. The effect of government influence and the vulnerability of migrant labour using the product of the matrices $M_1$ and $M_2$ and taking the transpose.

$$
\left(M_1 \times M_2\right)^T = \begin{bmatrix}
0 & 6 & 8 & 4 & 0 & 0 & 0\\
0 & 31 & 33 & 33 & 8 & 0 & 10\\
0 & 6 & 8 & 11 & 10 & 0 & 9\\
0 & 0 & 0 & 0 & 0 & 0 & 0\\
0 & -1 & -2 & -1 & 0 & 0 & -2\\
0 & 27 & 29 & 26 & 12 & 0 & 11\\
0 & -6 & -7 & -7 & -2 & 0 & -5\\
0 & 6 & 4 & 6 & 4 & 0 & 0\\
0 & -10 & -13 & -13 & -6 & 0 & -13
\end{bmatrix}.
$$

We study this in the interval [–14, 14]. This gives the indirect relationship between government policies and the vulnerability of migrant labourers to HIV/AIDS.

The reader is expected to analyse this BAM model using C-program given in the book [233] and draw conclusions.

## 2.5 Description of Fuzzy Associative Memories (FAM) model and their illustrations

This section has three subsections. In subsection one, we introduce FAM and in subsection two FAM is used to study the socio economic problem of women with HIV/AIDS. In subsection three we study the environmental pollution due to dyeing industries using the FAM model.

### 2.5.1 Introduction to Fuzzy Associative Memories

A fuzzy set is a map $\mu : X \to [0, 1]$ where X is any set called the domain and [0, 1] the range i.e., μ is thought of as a membership function i.e., to every element $x \in X$ μ assigns membership value in the interval [0, 1]. But very few try to visualize the geometry of fuzzy sets. It is not only of interest but is meaningful to see the geometry of fuzzy sets when we discuss fuzziness. Till date researchers over looked such visualization [Kosko, 108-112], instead they have interpreted fuzzy sets as generalized indicator or membership functions mappings μ from domain X to range [0, 1]. But functions are hard to visualize. Fuzzy theorist often picture membership functions as two-dimensional graphs with the domain X represented as a one-dimensional axis.

The geometry of fuzzy sets involves both domain $X = (x_1,\ldots, x_n)$ and the range [0, 1] of mappings $\mu : X \to [0, 1]$. The geometry of fuzzy sets aids us when we describe fuzziness, define fuzzy concepts and prove fuzzy theorems. Visualizing this geometry may by itself provide the most powerful argument for fuzziness.

An odd question reveals the geometry of fuzzy sets. What does the fuzzy power set $F(2^X)$, the set of all fuzzy subsets of X, look like? It looks like a cube, What does a fuzzy set look like? A fuzzy subsets equals the unit hyper cube $I^n = [0, 1]^n$. The fuzzy set is a point in the cube $I^n$. Vertices of the cube $I^n$ define a non-fuzzy set. Now with in the unit hyper cube $I^n = [0, 1]^n$ we are interested in a distance between points, which led to measures of size and fuzziness of a fuzzy set and more fundamentally to a measure. Thus within cube theory directly extends to the continuous case when the space X is a subset of $R^n$. The next step is to consider mappings between fuzzy cubes. This level of abstraction provides a surprising and fruitful alternative to the prepositional and predicate calculus reasoning techniques used in artificial intelligence (AI) expert systems. It allows us to reason with sets instead of propositions. The fuzzy set framework is numerical and multidimensional. The AI framework is symbolic and is one dimensional with usually only bivalent expert rules or propositions allowed. Both frameworks can encode structured knowledge in linguistic form. But the fuzzy approach translates the structured knowledge into a flexible numerical framework and processes it in a manner that resembles neural network processing. The numerical framework also allows us to adaptively infer and modify fuzzy systems perhaps with neural or statistical techniques directly from problem domain sample data.

Between cube theory is fuzzy-systems theory. A fuzzy set defines a point in a cube. A fuzzy system defines a mapping between cubes. A fuzzy system S maps fuzzy sets to fuzzy sets. Thus a fuzzy system S is a transformation S: $I^n \to I^P$. The n-dimensional unit hyper cube $I^n$ houses all the fuzzy subsets of the domain space or input universe of discourse $X = \{x_1, \ldots, x_n\}$. $I^p$ houses all the fuzzy subsets of the range space or output universe of discourse, $Y = \{y_1, \ldots, y_p\}$. X and Y can also denote subsets of $R^n$ and $R^p$. Then the fuzzy power sets F $(2^X)$ and F $(2^Y)$ replace $I^n$ and $I^p$.

In general a fuzzy system S maps families of fuzzy sets to families of fuzzy sets thus S: $I^{n_1} \times \ldots \times I^{n_r} \to I^{p_1} \times \ldots \times I^{p_s}$ Here too we can extend the definition of a fuzzy system to allow arbitrary products or arbitrary mathematical spaces to serve as the domain or range spaces of the fuzzy sets. We shall focus on fuzzy systems S: $I^n \to I^P$ that map balls of fuzzy sets in $I^n$ to balls of fuzzy set in $I^p$. These continuous fuzzy systems behave as associative memories. The map close inputs to close outputs. We shall refer to them as Fuzzy Associative Maps or FAMs.

The simplest FAM encodes the FAM rule or association ($A_i$, $B_i$), which associates the p-dimensional fuzzy set $B_i$ with the n-dimensional fuzzy set $A_i$. These minimal FAMs essentially map one ball in $I^n$ to one ball in $I^p$. They are comparable to simple neural networks. But we need not adaptively train the minimal FAMs. As discussed below, we can directly encode structured knowledge of the form, "If traffic is heavy in this direction then keep the stop light green longer" is a Hebbian-style FAM correlation matrix. In practice we sidestep this large numerical matrix with a virtual representation scheme. In the place of the matrix the user encodes the fuzzy set association (Heavy, longer) as a single linguistic entry in a FAM bank linguistic matrix. In general a FAM system F: $I^n \rightarrow I^b$ encodes the processes in parallel a FAM bank of m FAM rules ($A_1$, $B_1$), …, ($A_m$ $B_m$). Each input A to the FAM system activates each stored FAM rule to different degree. The minimal FAM that stores ($A_i$, $B_i$) maps input A to $B_i$' a partly activated version of $B_i$. The more A resembles $A_i$, the more $B_i$' resembles $B_i$. The corresponding output fuzzy set B combines these partially activated fuzzy sets $B_1^1, B_2^1, \ldots, B_m^1$. B equals a weighted average of the partially activated sets B = $w_1 B_1^1 + \ldots + w_n B_m^1$ where $w_i$ reflects the credibility frequency or strength of fuzzy association ($A_i$, $B_i$). In practice we usually defuzzify the output waveform B to a single numerical value $y_j$ in Y by computing the fuzzy centroid of B with respect to the output universe of discourse Y.

More generally a FAM system encodes a bank of compound FAM rules that associate multiple output or consequent fuzzy sets $B^1{}_i$, …, $B_i^s$ with multiple input or antecedent fuzzy sets $A_i^1$, …, $A_i^r$. We can treat compound FAM rules as compound linguistic conditionals. This allows us to naturally and in many cases easily to obtain structural knowledge. We combine antecedent and consequent sets with logical conjunction, disjunction or negation. For instance, we could interpret the compound association ($A^1$, $A^2$, B), linguistically as the compound conditional "IF $X^1$ is $A^1$ AND $X^2$ is $A^2$, THEN Y is B" if the comma is the fuzzy association ($A^1$, $A^2$, B) denotes conjunction instead of say disjunction.

We specify in advance the numerical universe of discourse for fuzzy variables $X^1$, $X^2$ and Y. For each universe of discourse or fuzzy variable X, we specify an appropriate library of fuzzy set values $A_1^r$, …, $A_k^2$ Contiguous fuzzy sets in a library overlap. In principle a neural network can estimate these libraries of fuzzy sets. In practice this is usually unnecessary. The library sets represent a weighted though overlapping quantization of the input space X. They represent the fuzzy set values assumed by a fuzzy variable. A different library of fuzzy sets similarly quantizes the output space Y. Once we define the library of fuzzy sets we construct the FAM by choosing appropriate combinations of input and output fuzzy sets Adaptive techniques can make, assist or modify these choices.

An Adaptive FAM (AFAM) is a time varying FAM system. System parameters gradually change as the FAM system samples and processes data. Here we discuss how natural network algorithms can adaptively infer FAM rules from training data. In principle, learning can modify other FAM system components, such as the libraries of fuzzy sets or the FAM-rule weights $w_i$.

In the following subsection we propose and illustrate an unsupervised adaptive clustering scheme based on competitive learning to blindly generate and refine the bank of FAM rules. In some cases we can use supervised learning techniques if we have additional information to accurately generate error estimates. Thus Fuzzy Associative Memories (FAMs) are transformation. FAMs map fuzzy sets to fuzzy sets. They map unit cubes to unit cubes. In simplest case the FAM system consists of a single association. In general the FAM system consists of a bank of different FAM association. Each association corresponds to a different numerical FAM matrix or a different entry in a linguistic FAM-bank matrix. We do not combine these matrices as we combine or superimpose neural-network associative memory matrices. We store the matrices and access them in parallel. We begin with single association FAMs. We proceed on to adopt this model to the problem.

### 2.5.2 Use of FAM model to study the socio economic problem of women with HIV/AIDS

At the outset we first wish to state that these women are mainly from the rural areas, they are economically poor and uneducated and majority of them are infected by their husbands. Second we use FAM because only this will indicate the gradations of the causes, which is the major cause for women being affected by HIV/AIDS followed in order of gradation the causes. Further among the fuzzy tools FAM model alone can give such gradations so we use them in this analysis. Another reason for using FAM is they can be used with the same FRM that is using the attributes of the FRM, FAMs can also be formulated. Already FAMs are very briefly described in section one of this chapter. We now give the sketch of the analysis of this problem with a view that any reader with a high school education will be in a position to follow it.

***Example 2.5.2.1:*** We just illustrate two experts opinion though we have used several experts' opinion for this analysis. Using the problems of women affected with HIV/AIDS along the rows and the causes of it along the column we obtain the related fuzzy vector matrix M.

The following are taken as the attributes (concepts) related with women, which is taken along the rows.

$W_1$ - Child marriage / widower marriage / child married to men twice or thrice their age

$W_2$ - Causes of women being infected with HIV/AIDS

$W_3$ - Disease untreated till it is chronic or they are in last stages of life due to full-blown AIDS

$W_4$ - Women are not traditionally considered as breadwinners for the work they do and the money they earn to protect the family is never given any recognition

$W_5$ - Free of depression and despair in spite of being deserted by family when they have HIV/AIDS

$W_6$ - Faith in god / power of god will cure
$W_7$ - Married women have acquired the disease due to their husbands.

The concepts associated with society, men / husband is taken along the columns.

$R_1$ - Female child a burden so the sooner they get her married off, the better relief economically
$R_2$ - Poverty / not owners of property
$R_3$ - Bad habits of the men / husbands
$R_4$ - Infected women are left uncared by relatives, even by husbands
$R_5$ - No Guilt or fear of life
$R_6$ - Not changed religion and developed faith in god after the disease
$R_7$ - No moral responsibility on the part of husbands and they infect their wives willfully
$R_8$ - Frequent, natural abortion / death of born infants
$R_9$ - STD / VD infected husbands
$R_{10}$ - Husbands hide their disease from their family so the wife become HIV/AIDS infected.

The gradations are given in the form of the fuzzy vector matrix M which is as follows.

$$
\begin{array}{c} \\ W_1 \\ W_2 \\ W_3 \\ W_4 \\ W_5 \\ W_6 \\ W_7 \end{array}
\begin{array}{c}
\begin{array}{cccccccccc} R_1 & R_2 & R_3 & R_4 & R_5 & R_6 & R_7 & R_8 & R_9 & R_{10} \end{array} \\
\left[\begin{array}{cccccccccc}
0.9 & 0.8 & 0.7 & 0 & 0 & 0 & 0 & 0 & 0 & 0.7 \\
0.5 & 0.8 & 0.6 & 0 & 0 & 0 & 0 & 0 & 0 & 0 \\
0 & 0.3 & 0.6 & 0 & 0 & 0 & 0 & 0 & 0 & 0 \\
0 & 0 & 0 & 0.6 & 0 & 0 & 0 & 0 & 0 & 0 \\
0 & 0 & 0 & 0 & 0.9 & .6 & .7 & 0 & 0 & 0 \\
0 & 0 & 0 & 0 & 0 & 0.7 & 0.5 & 0 & 0 & 0 \\
0 & 0 & 0 & 0 & 0.6 & 0 & 0 & 0 & 0 & 0
\end{array}\right]
\end{array}
$$

Now using the expert's opinion we get the fit vectors. Suppose the fit vector B is given as B = (0 1 1 0 0 0 0 0 1 0). Using max-min in backward directional we get FAM described as

$$A = M \text{ o } B$$

that is $a_i = \max \min (m_{ij}, b_j) \quad 1 \leq j \leq 10.$

Thus A = (0.8, 0.8, 0.6 0 0 0 0) since 0.8 is the largest value in the fit vector in A and it is associated with the two nodes vide $W_1$ and $W_2$ child marriage / widower child marriage etc finds it first place also the vulnerability of rural uneducated women find the same states as that of $W_1$. Further the second place is given to $W_3$, disease untreated till it is chronic or they are in last stages, all other states are in the off state.

Suppose we consider the resultant vector

$$A = (0.8, 0.8, 0.6, 0\ 0\ 0\ 0)$$
$$A \text{ o } M = B. \text{ where}$$
$$B_j = \max (a_i, m_{ij}) \quad 1 \leq j \leq 10.$$

$$A \text{ o } M = (0.9, 0.8, 0.8, 0.6, 0.9, 0.7, 0.7, 0, 0, 0.7) = B_1.$$

Thus we see the major cause being infected by HIV/AIDS is $R_1$ and $R_5$ having the maximum value. viz. female child a burden so the sooner they get married off the better relief economically; and women suffer no guilt and fear for life. The next value being given by $R_2$ and $R_3$, poverty and bad habits of men are the causes of women becoming HIV/AIDS victims.

The next value being $R_6$, $R_7$ and $R_{10}$ taking the value 0.7, women have not changed religion and developed faith in god after the disease, lost faith in god after the disease. Husbands hide their disease from family so wife becomes HIV/AIDS infected. Several conclusions can be derived from the analysis of similar form.

Now we proceed on to describe one more model using FAM. Using the attributes related with women affected with HIV/AIDS along the rows and the causes of it along the column we obtain the fuzzy rector matrix $M_1$. The gradations are given in the form of the fuzzy vector matrix which is as follows.

$$\begin{array}{c} \\ W_1 \\ W_2 \\ W_3 \\ W_4 \\ W_5 \\ W_6 \\ W_7 \\ W_8 \\ W_9 \end{array} \begin{array}{c} \begin{array}{cccccccccc} R_1 & R_2 & R_3 & R_4 & R_5 & R_6 & R_7 & R_8 & R_9 & R_{10} \end{array} \\ \left[\begin{array}{cccccccccc} .9 & .8 & .6 & 0 & .5 & .7 & 0 & 0 & .4 & 0 \\ .8 & .7 & .8 & 0 & 0 & 0 & 0 & 0 & .7 & 0 \\ .7 & .6 & .5 & 0 & 0 & .8 & 0 & 0 & 0 & 0 \\ .9 & .7 & 0 & 0 & 0 & 0 & 0 & .8 & .7 & 0 \\ .8 & .6 & .7 & 0 & 0 & 0 & 0 & .7 & 0 & 0 \\ 0 & 0 & .4 & 0 & 0 & 0 & 0 & .6 & 0 & 0 \\ .4 & 0 & .6 & 0 & 0 & 0 & 0 & 0 & 0 & 0 \\ .6 & .5 & 0 & 0 & .4 & .3 & 0 & 0 & 0 & 0 \\ 0 & 0 & 0 & 0 & 0 & 0 & .8 & 0 & .4 & .6 \end{array}\right] \end{array}.$$

Now using the experts opinion we get the fit vectors B, given as B = (1 1 1 0 0 0 0 0 1), using max-min backward direction we get FAM described as

$$A = M_1 \text{ o } B.$$

That is

$$\begin{aligned} a_j &= \max \min (m_{ij}\, b_j),\ 1 \le j \le 10 \\ &= (0.9, 0.8, 0.8, 0, 0.5, 0.8, 0.8, 0, 0.7, 0.6). \end{aligned}$$

Since 0.9 is the largest value in the fit vector A and it is associated with $W_1$ and $W_4$, we see the root cause or the major reason of women becoming HIV/AIDS infected is due to no basic education, no recognition of their labour. One may be very much surprised to see why the co ordinate "no recognition of labour" find so much place the main reason for it if they (men i.e., the husband) had any recognition of their labour certainly he would not have the heart to infect his wife when he is fully aware of the fact that he is suffering from HIV/AIDS and by having unprotected sex with his wife certainly she too would soon be an HIV/AIDS patient. Thus $W_1$ and $W_4$ gets the highest value.

The next value is 0.8, which is got by the nodes $W_2$ and $W_5$. $W_2$ – No wealth and property and $W_5$ – No protection from their spouse. If wealth and property was in their hands (i.e., in the name of women) certainly they would fear her untimely death for the property may be taken by her relatives. Further these women cannot protect themselves from their husbands even if they are fully aware of the fact that their husbands are infected by HIV/AIDS.

Next higher value is .6 taken by $W_3$, $W_7$, $W_8$ and $W_9$ that is women are not empowered, they do not get nutritious food and they are not the decision makers or they do not have any voice in this traditional set up. Thus we have worked with several different nodes to arrive at the conclusions.

### 2.5.3 Use of FAM model to study the environmental pollution due to dyeing industries

We now give the description of the problem. Dyeing industries of Tamil Nadu, especially the ones situated in and around Tiruppur have become the concern of environmentalists, government, public and courts. For the safety, social, moral and ethical considerations of their functioning are questionable. The environmental pollution caused by these dyeing industries is unimaginable. It is reported on 2 September 2005 in the New Indian Express that 120 tonnes of dead fish have been removed from the dam site. It was further advised that people who consume these fish have the risk of being food-poisoned. A team lead by a retired PWD Chief Engineer to study the desilting of the Orathupalayam dam was to be taken. As per the orders of the Madras High Court the release of the entire volume of water into the Cauvery river resulted in the death of 120 tonnes of fish that could have yield profit in lakhs. Here, the labour situation and wages for the workers are in general bad and their issues are never taken up by anyone. Further no proper procedure is adopted in treating the wastes from the dyeing units. Thus very recently, on 27 August 2005 the Tamil Nadu government has planned with the help of Russian institutes to set up a center for treating the wastages from these dyeing industries for Rs.323 crores.

Not only fish, earlier it was reported on 28 August 2005 that several lakhs worth of betel leaves were damaged because the water had been contaminated by the dyeing unit's wastages. Over 100 acres of crop was damaged beyond use. The people in and around these units suffer several untold health problems. The politicians, public and the government are keeping quiet. The main reason attributed is the amount of revenue the nation gets because of these dyeing industries. But as our study is only about environmental pollution by these dyeing industries, so we do not indulge about the problems related with labour/wages/ hours of work.

***Example 2.5.3.1:*** Now we use FAM model to analyse the environmental pollution caused by dyeing industries

Since the environmental problems faced due to the wastages from these dyeing industries cannot be put in as a data for it pollutes, the living creatures of the water (dead fish due to poisoning came to 120 tonnes) including crabs and other aquatic plants and animals, it also destroyed plantations, crops and above all the health conditions of people at large cannot be comprehended. Hence we are justified in using fuzzy theory in general and FAM in particular to analyze this problem. As the analysis involves lots of uncertainty and unpredictability, we felt FAM will be the best tool and also it is a model in which the resultant vectors are graded. We had been to the neighbourhoods of these dyeing units and collected the data from people and also visited the spots of havoc. This data was obtained in the form of linguistic questionnaire. The expert's opinion was transformed into a FAM model. The attributes given by majority of the experts which are consolidated as follows:

Attributes related with the dyeing industries.

$D_1$ - No proper labour standard

$D_2$ - No proper means to treat the poisonous dyeing waste

$D_3$ - No moral responsibility of polluting the land and water resources etc

$D_4$ - Money making is the only motive
$D_5$ - Refusal to follow any industrial acts / trade union acts
$D_6$ - Not followed any norms of Tamil Nadu Pollution Control Board
$D_7$ - Legal remedies favouring them due to powerful economic strata
$D_8$ - No moral responsibility of any of the problems be it material / ethical/ human
$D_9$ - Very dangerous total environmental pollution
$D_{10}$ - Increases migrant labourers.

The attributes related to the environmental pollution are

$E_1$ - Death of over 120 tonnes of fish
$E_2$ - All living creations in water destroyed or slowly destroyed due to the mixing of wastages from these dyeing units with water resources in and around the industries
$E_3$ - Death of plantations
$E_4$ - Health hazards suffered by the labourers employed in these industries
$E_5$ - Health hazards suffered by the people living in and around these industries
$E_6$ - Pollution of soil
$E_7$ - Pollution of ground water.

Here it is important to mention that it is the liberty of one to work with more or less number of attributes. We have selected these attributes and we work with them using the FAM model. We illustrate here only one expert's opinion, though we have used several of the expert's opinions for the study. Using the attributes related with the industry along the rows and the environmental pollution by the dyeing industries along the columns we obtain the related fuzzy vector matrix $M^T$.

The gradations are given in the form of the fuzzy vector matrix $M^T$ which is as follows.

$$M^T = \begin{bmatrix} 0.9 & 0.8 & 0.7 & 0 & 0 & 0 & 0 & 0 & 0 & 0.7 \\ 0.5 & 0.8 & 0.6 & 0 & 0 & 0 & 0 & 0 & 0 & 0 \\ 0 & 0.3 & 0.6 & 0 & 0 & 0 & 0 & 0 & 0 & 0 \\ 0 & 0 & 0 & 0.6 & 0 & 0 & 0 & 0 & 0 & 0 \\ 0 & 0 & 0 & 0 & 0.9 & .6 & .7 & 0 & 0 & 0 \\ 0 & 0 & 0 & 0 & 0 & 0.7 & 0.5 & 0 & 0 & 0 \\ 0 & 0 & 0 & 0 & 0.6 & 0 & 0 & 0 & 0 & 0 \end{bmatrix}.$$

Now using the experts opinion we get the fit vectors.

Suppose the fit vector B is given as B = (0 1 1 0 0 0 0 0 1 0).

Using max-min in backward direction we get FAM described as

$$A = M^T \text{ o } B$$

$$a_i = \max_{1 \le j \le 10} \min (m_{ij}, b_j).$$

Thus A = (0.8, 0.8, 0.6, 0 0 0 0) since 0.8 is the largest value of the fit vector A and it is associated with the two nodes $E_1$ and $E_2$, Death of over 120 tonnes of fish and All living creations in water destroyed or slowly destroyed due to the mixing of wastages from these dyeing units with water resources in and around the industries finds its first place, also this is due to the fact local people have become very much affected by the death of 120 tonnes of fish and this is fresh in their memories, so they are also naturally agreeing to the fact that the natural water resources are highly polluted. Likewise the complete destruction of crops (betel leaves) has also affected them which are worth several lakhs. It is a pity they forget their own health hazards.

Suppose we consider the resultant vector A = (0.8, 0.8, 0.6, 0 0 0 0)

$$A \text{ o } M^T = B$$

where

$$B_j = \max_{1 \le j \le 10} (a_i, m_{ij})$$

$$A \text{ o } M^T = (0.9, 0.8, 0.8, 0.6, 0.9, 0.7, 0.7, 0, 0, 0.7) = B_i$$

Thus the major cause of being affected by these dyeing industries and $D_1$ and $D_5$ having the Maximum value viz. No proper labour standard and Refusal to follow any industrial acts / trade union acts. The next higher value being given by $D_2$ and $D_3$, i.e., No proper means to treat the poisonous Dyeing Waste and No moral responsibility of polluting the land and water resources etc. The next value being $D_6$, $D_7$ and $D_{10}$ taking the value 0.7, that is the industries have not followed any norms of Tamil Nadu Pollution Control Board, Legal Remedies favouring them due to powerful economic starta and increase of migrant labourers.

From our study we see that the major reason for the chaotic behaviour of these dyeing industries is that they do not follow any labour standards and their refusal to follow the industrial acts/ trade union acts has resulted in the dangerous pollution of the environment resulting in the death of fish and crop which is a major and only source of the livelihood of the natives of that place.

Secondly we see that they do not follow any method to treat the waste, we see that the two reasons which can be attributed are:

1. Treatment of these wastes is very costly.
2. They know that they can do anything and get away with it, if they are economically powerful.

Also the death of fish and crop are very strong in the memories of the public so whatever is spoken about the pollution ends in their description of the tragedy. It is high time the owners of the dyeing industries take up moral responsibilities and stop further pollution. They do not imagine the harm they do to the nation which is irreparable damage to the soil and water and can have disastrous consequences for people’s health. Unless very strong laws are made to punish the erring industries it would be impossible to save the nation from the environmental pollution.

## 2.6 Fuzzy Relational Equations and their Application

The notion of fuzzy relational equations based upon the max-min composition was first investigated by Sanchez [180]. He studied conditions and theoretical methods to resolve fuzzy relations on fuzzy sets defined as mappings from sets to [0, 1]. Some theorems for existence and determination of solutions of certain basic fuzzy relation equations were given by him. However the solution obtained by him is only the greatest element (or the maximum solution) derived from the max-min (or min-max) composition of fuzzy relations. [180]'s work has shed some light on this important subject. Since then many researchers have been trying to explore the problem and develop solution procedures [1, 4, 17, 37, 40-1, 46, 48, 50-4, 59, 61, 66, 75-6, 174, 237]. The max-min composition is commonly used when a system requires conservative solutions in the sense that the goodness of one value cannot compensate the badness of another value. In reality there are situations that allow compensatability among the values of a solution vector. In such cases the min operator is not the best choice for the intersection of fuzzy sets, but max-product composition, is preferred since it can yield better or at least equivalent result. Before we go into the discussion of these Fuzzy Relational Equations (FRE) and its properties, it uses and applications we just describe them. This section is divided into five subsections. Section one and two describe FREs and their properties. Fuzzy compatibility and composition of fuzzy relations are dealt with in section three. FRE is used in real world problems like Transportation and Globalizations in section four and five respectively.

### 2.6.1 Binary Fuzzy Relations and their properties

It is well known fact that binary relations are generalized mathematical functions. Contrary to functions from X to Y, binary relations R(X, Y) may assign to each element of X two or more elements of Y. Some basic operations on functions such as the inverse and composition are applicable to binary relations as well.

Given a fuzzy relation R(X, Y), its domain is a fuzzy set on X, dom R, whose membership function is defined by

$$\text{dom}\,(R(x)) = \max_{y \in Y} R\,(x, y)$$

for each $x \in X$. That is, each element of set X belongs to the domain of R to the degree equal to the strength of its strongest relation to any member of set Y. The range of R (X, Y) is a fuzzy relation on Y, ran R whose membership function is defined by

$$\text{ran}\,R(y) = \max_{x \in X} R\,(x, y)$$

for each $y \in Y$. That is, the strength of the strongest relation that each element of Y has to an element of X is equal to the degree of that elements membership in the range of R. In addition, the height of a fuzzy relation R(X,Y) is a number, h(R), defined by

$$h\,(R) = \max_{y \in Y}\ \max_{x \in X}\,(R\,(x, y).$$

That is h(R) is the largest membership grade attained by any pair (x, y) in R.

A convenient representation of binary relation R(X, Y) are membership matrices $R = [r_{xy}]$ where $r_{xy} = R(x, y)$. Another useful representation of binary relation is a sagittal diagram. Each of the sets X, Y is represented by a set of nodes in the diagram nodes corresponding to one set is distinguished from nodes representing the other set.

Elements of $X \times Y$ with non-zero membership grades in R(X, Y) are represented in the diagram by lines connecting the respective nodes.

We illustrate the sagittal diagram of a binary fuzzy relation R(X, Y) together with the corresponding membership matrix in Figure 2.6.1.1.

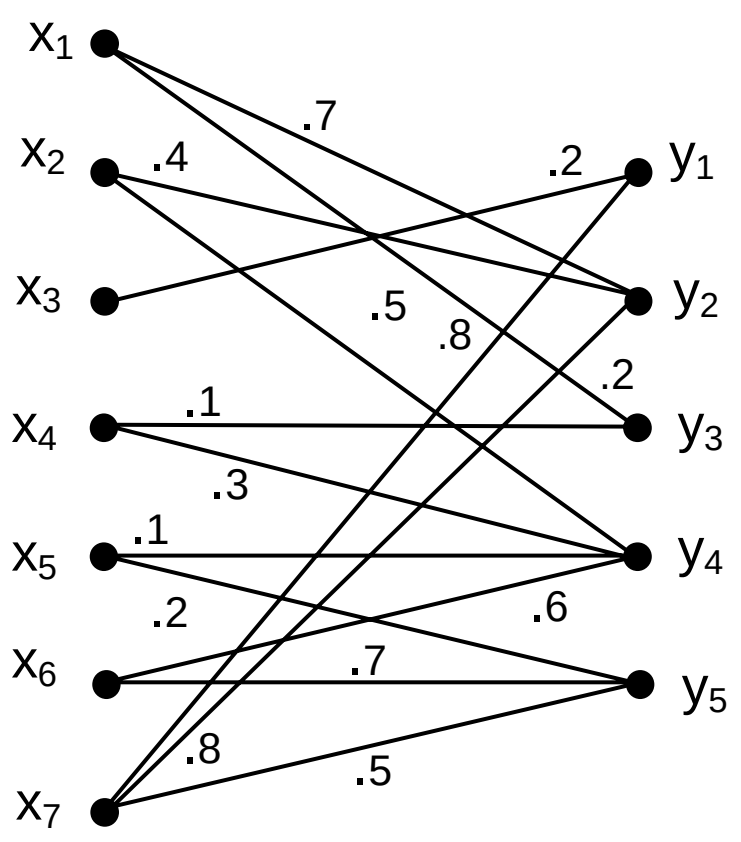

FIGURE: 2.6.1.1

The inverse of a fuzzy relation R(X, Y) denoted by $R^{-1}$(Y, X) is a relation on Y × X defined by $R^{-1}$ (y, x) = R (x, y) for all x ∈ X and for all y ∈ X. A membership matrix $R^{-1} = [r^{-1}_{yx}]$ representing $R^{-1}$ (Y, X) is the transpose of the matrix R for R (X, Y) which means that the rows of $R^{-1}$ equal the columns of R and the columns of $R^{-1}$ equal the rows of R.

Clearly $(R^{-1})^{-1} = R$ for any binary fuzzy relation. Thus a fuzzy binary relation can be represented by the sagittal diagram. The corresponding membership matrix:

$$
\begin{array}{c} \\ x_1 \\ x_2 \\ x_3 \\ x_4 \\ x_5 \\ x_6 \\ x_7 \end{array}
\begin{array}{c}
\begin{array}{ccccc} y_1 & y_2 & y_3 & y_4 & y_5 \end{array} \\
\begin{bmatrix}
0 & .7 & .5 & 0 & 0 \\
0 & .4 & 0 & .1 & 0 \\
.2 & 0 & 0 & 0 & 0 \\
0 & 0 & .1 & 1 & 0 \\
0 & 0 & 0 & .3 & .7 \\
0 & 0 & 0 & .6 & .7 \\
.2 & 0 & .8 & 0 & .5
\end{bmatrix}
\end{array}.
$$

R is the membership matrix.

Consider now two binary fuzzy relations P(X, Y) and Q(Y, Z) with a common set Y. The standard composition of these relations, which is denoted by P(X, Y) ∘ Q(Y, Z), produces a binary relation R (X, Z) on X × Z defined by

$$R(x, z) = [P \circ Q](x, z) = \max_{y \in Y} \min [P(x, y), Q(y, z)]$$

for all x ∈ X and all z ∈ Z. This composition, which is based on the standard t-norm and t-co-norm is often referred to as the max-min composition. It follows directly from the above equation that

$$[P(X, Y) \circ Q(Y, Z)]^{-1} = Q^{-1}(Z, Y) \circ P^{-1}(Y, X)$$
$$[P(X, Y) \circ Q(Y, Z)] \circ R(Z, W) = P(X, Y) \circ [Q(Y, Z) \circ R(Z, W)].$$

This is the standard (or max-min) composition, which is associative, and its inverse is equal to the reverse composition of the inverse relations.

However the standard composition is not commutative because Q(Y, Z) ∘ P(X, Y) is not even well defined when X ≠ Z. Even if X = Z and Q (Y, Z) ∘ P(X, Y) are well defined, we may have P(X, Y) ∘ Q (Y, Z) ≠ Q (Y, Z ) ∘ P(X, Y). Compositions of binary fuzzy relations can be performed conveniently in terms of membership matrices of the relations. Let P = $[p_{ik}]$, Q = $[q_{kj}]$ and R = $[r_{ij}]$ be the membership matrices of binary relations such that R = P ∘ Q. We can then write using this matrix notations

$$[r_{ij}] = [p_{ik}] \circ [q_{kj}]$$

where $r_{ij} = \max_{k} \min (p_{ik}\ q_{kj})$. Observe that the same elements of P and Q are used in the calculation of R as would be used in the regular multiplication of matrices, but the product and sum operations are here replaced with max and min operations respectively.

A similar operation on two binary relations, which differs from the composition in that it yields triples instead of pairs, is known as the relational join. For fuzzy relations P(X, Y) and Q(Y, Z), the relational join P * Q, corresponding to the standard max-min composition is a ternary relation R(X, Y, Z) defined by

$$R(x, y, z) = [P * Q](x, y, z) = \min [P(x, y), Q(y, z)]$$

for each $x \in X$, $y \in Y$ and $z \in Z$.

The fact that the relational join produces a ternary relation from two binary relations is a major difference from the composition, which results in another binary relation.

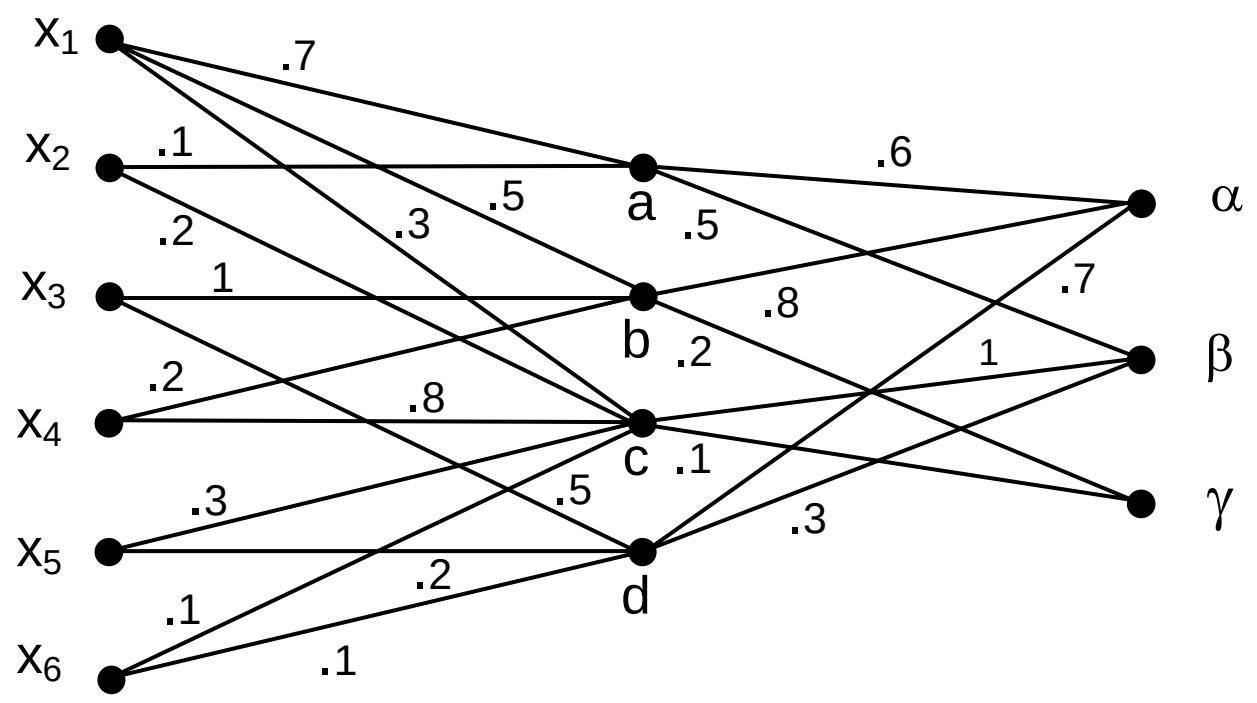

FIGURE: 2.6.1.2

Formally [P °Q] (x, z) = max [P * Q] (x, y, z) for each $x \in X$ and $z \in Z$. Now we just see what happens if the binary relation on a single set. Binary relation R (X, X) can be expressed by the same forms as general binary relations.

The following properties are to be observed:

i. Each element of the set X is represented as a single node in the diagram.
ii. Directed connection between nodes indicates pairs of elements of X, with the grade of membership in R is non-zero.

iii. Each connection in the diagram is labeled by the actual membership grade of the corresponding pair in R.

An example of this diagram for a relation R (X, X) defined on X = {$x_1$, $x_2$, $x_3$, $x_4$, $x_5$} is shown in figure 2.6.1.3. A crisp relation R (X, X) is reflexive if and only if (x, x) $\in$ R. for each x $\in$ R, that is if every element of X is related to itself, otherwise R(X, X) is called irreflexive. If (x, x) $\notin$ R for every x $\in$ X the relation is called anti reflexive.

A crisp relation R(X, X) is symmetric if and only if for every (x, y) $\in$ R, it is also the case that (y, x) $\in$ R where x, y $\in$ X. Thus whenever an element x is related to an element y through a symmetric relation, y is also related to x. If this is not the case for some x, y then the relation is called asymmetric. If both (x, y) $\in$ R and (y, x) $\in$ R implies x = y then the relation is called anti symmetric. If either (x, y) $\in$ R or (y, x) $\in$ R whenever x $\neq$ y, then the relation is called strictly anti symmetric.

$$\begin{array}{c} \\ x_1 \\ x_2 \\ x_3 \\ x_4 \\ x_5 \end{array} \begin{array}{c} \begin{array}{ccccc} x_1 & x_2 & x_3 & x_4 & x_5 \end{array} \\ \begin{bmatrix} .2 & 0 & .1 & 0 & .7 \\ 0 & .8 & .6 & 0 & 0 \\ .1 & 0 & .2 & 0 & 0 \\ .2 & 0 & 0 & .1 & .4 \\ .6 & .1 & 0 & 0 & .5 \end{bmatrix} \end{array}$$

The corresponding sagittal diagram is given in Figure 2.6.1.3:

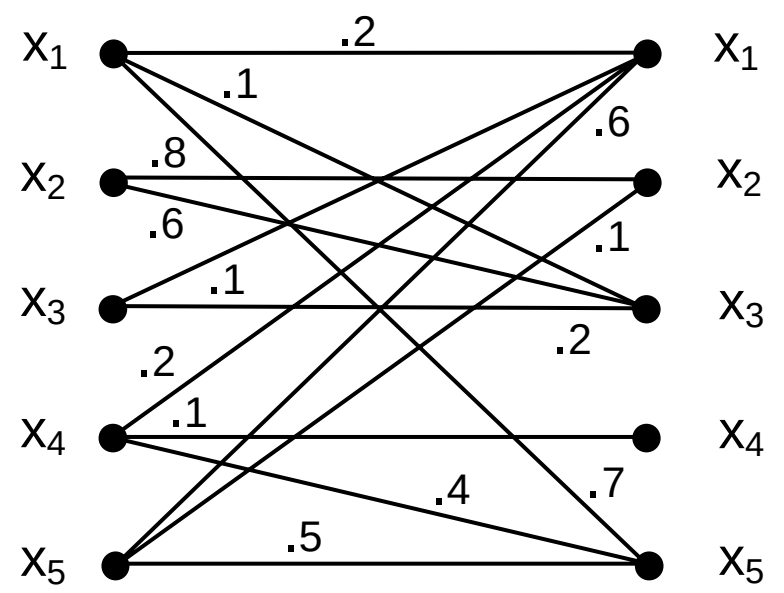

FIGURE: 2.6.1.3

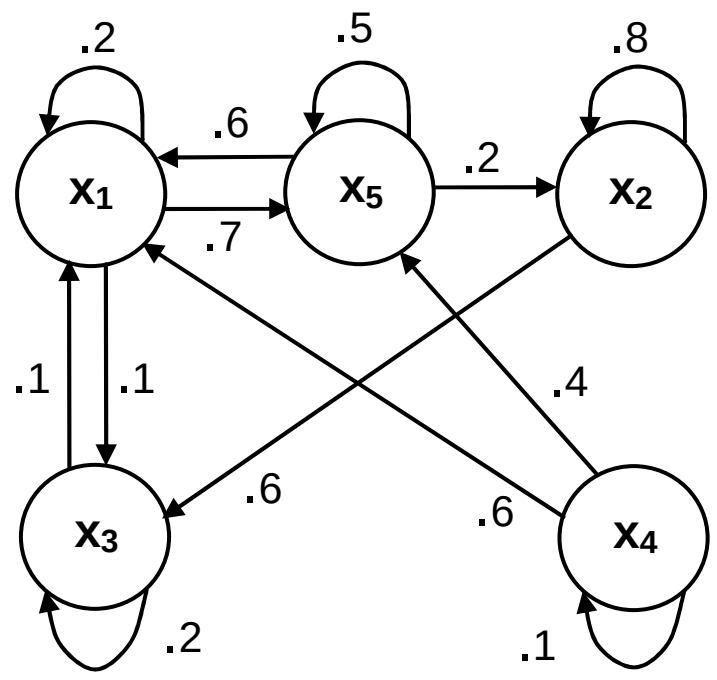

FIGURE: 2.6.1.4

Table

| x | y | R(x, y) |
|---|---|---|
| $x_1$ | $x_1$ | .2 |
| $x_1$ | $x_3$ | .1 |
| $x_1$ | $x_5$ | .7 |
| $x_2$ | $x_2$ | .8 |
| $x_2$ | $x_3$ | .6 |
| $x_3$ | $x_1$ | .1 |
| $x_3$ | $x_3$ | .2 |
| $x_4$ | $x_1$ | .2 |
| $x_4$ | $x_4$ | .1 |
| $x_4$ | $x_5$ | .4 |
| $x_5$ | $x_1$ | .6 |
| $x_5$ | $x_2$ | .1 |
| $x_5$ | $x_5$ | .5 |

A crisp relation R (X, Y) is called transitive if and only if (x, z) $\in$ R, whenever both (x, y) $\in$ R and (y, z) $\in$ R for at least one y $\in$ X.

In other words the relation of x to y and of y to z imply the relation x to z is a transitive relation. A relation that does not satisfy this property is called non-transitive. If (x, z) $\notin$ R whenever both (x, y) $\in$ R and (y, z) $\in$ R, then the relation is called anti transitive. The reflexivity, symmetry and transitivity is described by the following Figure 2.6.1.5:

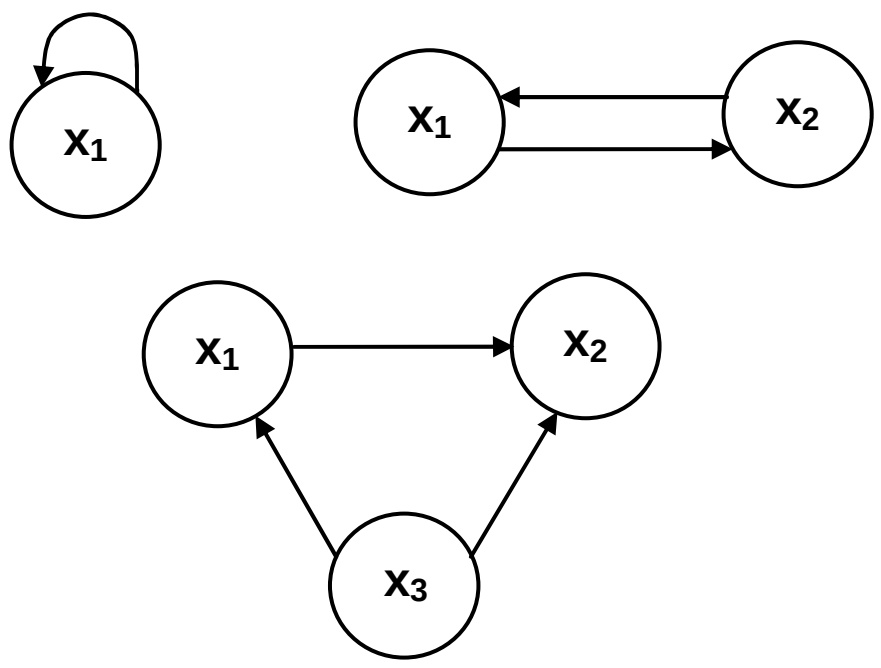

FIGURE: 2.6.1.5

A fuzzy relation R (X, X) is reflexive if and only if R(x,x) = 1 for all $x \in X$, if this is not the case for some $x \in X$, the relation is called irreflexive, if it is not satisfied for all $x \in X$, the relation is called anti reflexive. A weaker form of reflexivity referred to as $\in$ - reflexivity denoted by $R(x, x) \geq \in$ where $0 < \in < 1$. A fuzzy relation is symmetric if and only if

$$R(x, y) = R(y, x)$$

for all $x, y \in X$, if this relation is not true for some $x, y \in X$, the relation is called anti symmetric. Further more when $R(x, y) > 0$ and $R(y, x) > 0$ implies $x = y$ for all $x, y \in X$ the relation R is called anti symmetric.

A fuzzy relation R (X, X) is transitive if $R(x, z) \geq \max_{y \in Y} \min [R(x, y), R(y, z)]$ is satisfied for each pair $(x, z) \in X^2$. A relation failing to satisfy this inequality for some members of X is called non-transitive and if

$$R(x, z) < \max_{y \in Y} \min [R(x, y), R(y, z)]$$

for all $(x,z) \in X^2$, then the relation is called anti transitive.

### 2.6.2 Properties of Fuzzy Relations

In this section we just recollect the properties of fuzzy relations like, fuzzy equivalence relation, fuzzy compatibility relations, fuzzy ordering relations, fuzzy morphisms and sup-i-compositions of fuzzy relation. For more about these concepts please refer [106, 237].

Now we proceed on to define fuzzy equivalence relation. A crisp binary relation R(X, X) that is reflexive, symmetric and transitive is called an equivalence relation. For each element x in X, we can define a crisp set $A_x$, which contains all the elements of X that are related to x, by the equivalence relation.

$$A_x = \{y \mid (x, y) \in R (X, X)\}$$

$A_x$ is clearly a subset of X. The element x is itself contained in $A_x$ due to the reflexivity of R, because R is transitive and symmetric each member of $A_x$, is related to all the other members of $A_x$. Further no member of $A_x$, is related to any element of X not included in $A_x$. This set $A_x$ is referred to an as equivalence class of R (X, X) with respect to x. The members of each equivalence class can be considered equivalent to each other and only to each other under the relation R. The family of all such equivalence classes defined by the relation which is usually denoted by X / R, forms a partition on X.

A fuzzy binary relation that is reflexive, symmetric and transitive is known as a fuzzy equivalence relation or similarity relation. In the rest of this section let us use the latter term. While the max-min form of transitivity is assumed, in the following discussion on concepts; can be generalized to the alternative definition of fuzzy transitivity.

While an equivalence relation clearly groups elements that are equivalent under the relation into disjoint classes, the interpretation of a similarity relation can be approached in two different ways. First it can be considered to effectively group elements into crisp sets whose members are similar to each other to some specified degree. Obviously when this degree is equal to 1, the grouping is an equivalence class. Alternatively however we may wish to consider the degree of similarity that the elements of X have to some specified element $x \in X$. Thus for each $x \in X$, a similarity class can be defined as a fuzzy set in which the membership grade of any particular element represents the similarity of that element to the element x. If all the elements in the class are similar to x to the degree of 1 and similar to all elements outside the set to the degree of 0 then the grouping again becomes an equivalence class. We know every fuzzy relation R can be uniquely represented in terms of its α-cuts by the formula

$$R = \bigcup_{\alpha \in (0,1]} \alpha \bullet {}^{\alpha}R$$

It is easily verified that if R is a similarity relation then each α-cut, ${}^{\alpha}R$ is a crisp equivalence relation. Thus we may use any similarity relation R and by taking an α - cut ${}^{\alpha}R$ for any value $\alpha \in (0, 1]$, create a crisp equivalence relation that represents the presence of similarity between the elements to the degree α. Each of these equivalence relations form a partition of X. Let $\pi({}^{\alpha}R)$ denote the partition corresponding to the equivalence relation ${}^{\alpha}R$. Clearly any two elements x and y belong to the same block of this partition if and only if $R(x, y) \geq \alpha$. Each similarity relation is associated with the set $\pi(R) = \{\pi({}^{\alpha}R) \mid \alpha \in (0,1]\}$ of partition of X. These partitions are nested in the sense that $\pi({}^{\alpha}R)$ is a refinement of $\pi({}^{\beta}R)$ if and only $\alpha \geq \beta$.

The equivalence classes formed by the levels of refinement of a similarity relation can be interpreted as grouping elements that are similar to each other and only to each other to a degree not less than α.

Just as equivalences classes are defined by an equivalence relation, similarity classes are defined by a similarity relation. For a given similarity relation R(X, X) the similarity class for each $x \in X$ is a fuzzy set in which the membership grade of each element $y \in X$ is simply the strength of that elements relation to x or R(x, y). Thus the similarity class for an element x represents the degree to which all the other members of X are similar to x. Expect in the restricted case of equivalence classes themselves, similarity classes are fuzzy and therefore not generally disjoint.

Similarity relations are conveniently represented by membership matrices. Given a similarity relation R, the similarity class for each element is defined by the row of the membership matrix of R that corresponds to that element.

Fuzzy equivalence is a cutworthy property of binary relation R(X, X) since it is preserved in the classical sense in each $\alpha$-cut of R. This implies that the properties of fuzzy reflexivity, symmetry and max-min transitivity are also cutworthy. Binary relations are symmetric and transitive but not reflexive are usually referred to as quasi equivalence relations.

The notion of fuzzy equations is associated with the concept of compositions of binary relations. The composition of two fuzzy binary relations P (X, Y) and Q (Y, Z) can be defined, in general in terms of an operation on the membership matrices of P and Q that resembles matrix multiplication. This operation involves exactly the same combinations of matrix entries as in the regular matrix multiplication. However the multiplication and addition that are applied to these combinations in the matrix multiplication are replaced with other operations, these alternative operations represent in each given context the appropriate operations of fuzzy set intersections and union respectively. In the max-min composition for example, the multiplication and addition are replaced with the min and max operations respectively.

We shall give the notational conventions. Consider three fuzzy binary relations P (X, Y), Q (Y, Z) and R (X, Z) which are defined on the sets

$$X = \{x_i \mid i \in I\}$$
$$Y = \{y_j \mid j \in J\} \text{ and}$$
$$Z = \{z_k \mid k \in K\}$$

where we assume that $I = N_n$ $J = N_m$ and $K = N_s$. Let the membership matrices of P, Q and R be denoted by $P = [p_{ij}]$, $Q = [q_{ij}]$, $R = [r_{ik}]$ respectively, where $p_{ij} \in P(x_i, y_j)$, $q_{jk} \in Q(y_j, z_k)$ $r_{ij} \in R(x_i, z_k)$ for all $i \in I$ (= $N_n$), $j \in J$ (= $N_m$) and $k \in K$ (= $N_s$). This clearly implies that all entries in the matrices P, Q, and R are real numbers from the unit interval [0, 1]. Assume now that the three relations constrain each other in such a way that P ° Q = R where ° denotes max-min composition. This means that $\max_{j \in J} \min (p_{ij}, q_{jk}) = r_{ik}$ for all $i \in I$ and $k \in^{-} K$. That is the matrix equation P° Q = R encompasses $n \times s$ simultaneous equations of the form $\max_{j \in J} \min (p_{ij}, q_{jk}) = r_{ik}$. When two of the components in each of the equations are given and one is unknown these equations are referred to as fuzzy relation equations.

When matrices P and Q are given the matrix R is to determined using P ° Q = R. The problem is trivial. It is solved simply by performing the max-min multiplication – like operation on P and Q as defined by $\max_{j \in J} \min (p_{ij}, q_{jk}) = r_{ik}$. Clearly the solution in this case exists and is unique. The problem becomes far from trivial when one of the two matrices on the left hand side of P ° Q = R is unknown. In this case the solution is guaranteed neither to exist nor to be unique.

Since R in P ° Q = R is obtained by composing P and Q it is suggestive to view the problem of determining P (or alternatively Q ) from R to Q (or alternatively R and P) as a decomposition of R with respect to Q (or alternatively with respect to P). Since many problems in various contexts can be formulated as problems of decomposition, the utility of any method for solving P ° Q = R is quite high. The use of fuzzy relation equations in some applications is illustrated. Assume that we have a method for solving P ° Q = R only for the first decomposition problem (given Q and R).

Then we can directly utilize this method for solving the second decomposition problem as well. We simply write P ° Q = R in the form $Q^{-1}$ o $P^{-1} = R^{-1}$ employing transposed matrices. We can solve $Q^{-1}$ o $P^{-1} = R^{-1}$ for $Q^{-1}$ by our method and then obtain the solution of P ° Q = R by $(Q^{-1})^{-1} = Q$.

We study the problem of partitioning the equations P ° Q = R. We assume that a specific pair of matrices R and Q in the equations P ° Q = R is given. Let each particular matrix P that satisfies $P \circ Q = R$ is called its solution and let $S(Q, R) = \{P \mid P \circ Q = R\}$ denote the set of all solutions (the solution set).

It is easy to see this problem can be partitioned, without loss of generality into a set of simpler problems expressed by the matrix equations $p_i \circ Q = r_i$ for all $i \in I$ where

$$p_i = [p_{ij} \mid j \in J] \text{ and}$$
$$r_i = [r_{ik} \mid k \in K].$$

Indeed each of the equation in $\max_{j \in J} \min (p_{ij} q_{jk}) = r_{ik}$ contains unknown $p_{ij}$ identified only by one particular value of the index i, that is, the unknown $p_{ij}$ distinguished by different values of i do not appear together in any of the individual equations. Observe that $p_i$, Q, and $r_i$ in $p_i \circ Q = r_i$ represent respectively, a fuzzy set on Y, a fuzzy relation on $Y \times Z$ and a fuzzy set on Z. Let $S_i(Q, r_i) = [p_i \mid p_i \circ Q = r_i]$ denote, for each $i \in I$, the solution set of one of the simpler problem expressed by $p_i \circ Q = r_i$.

Thus the matrices P in $S(Q, R) = [P \mid P \circ Q = R]$ can be viewed as one column matrix

$$P = \begin{bmatrix} p_1 \\ p_2 \\ \vdots \\ p_n \end{bmatrix}$$

where $p_i \in S_i(Q, r_i)$ for all $i \in I = (=N_n)$. It follows immediately from $\max_{j \in J} \min (p_{ij} q_{jk}) = r_{ik}$. That if $\max_{j \in J} q_{jk} < r_{ik}$ for some $i \in I$ and some $k \in K$, then no values $p_{ij} \in [0, 1]$ exists ($j \in J$) that satisfy $P \circ Q = R$, therefore no matrix P exists that satisfies the matrix equation. This proposition can be stated more concisely as follows if

$$\max_{j \in J} q_{jk} < \max_{j \in J} r_{ik}$$

for some $k \in K$ then $S(Q, R) = \phi$. This proposition allows us in certain cases to determine quickly that $P \circ Q = R$ has no solutions its negation however is only a necessary not sufficient condition for the existence of a solution of $P \circ Q = R$ that is for $S(Q, R) \neq \phi$. Since $P \circ Q = R$ can be partitioned without loss of generality into a set of equations of the form $p_i \circ Q = r_i$ we need only methods for solving equations of the later form in order to arrive at a solution. We may therefore restrict our further discussion of matrix equations of the form $P \circ Q = R$ to matrix equation of the simpler form $P \circ Q = r$, where $p = [p_j \mid j \in J]$, $Q = [q_{jk} \mid j \in J, k \in K]$ and $r = \{r_k \mid k \in K]$.

We just recall the solution method as discussed by [106]. For the sake of consistency with our previous discussion, let us again assume that p, Q and r represent respectively a fuzzy set on Y, a fuzzy relation on $Y \times Z$ and a fuzzy set on Z. Moreover let $J = N_m$ and $K = N_s$ and let $S(Q, r) = \{p \mid p \circ Q = r\}$ denote the solution set of

$$p \circ Q = r.$$

In order to describe a method of solving $p \circ Q = r$ we need to introduce some additional concepts and convenient notation. First let $\wp$ denote the set of all possible vectors.

$$p = \{p_j \mid j \in J\}$$

such that $p_j \in [0, 1]$ for all $j \in J$, and let a partial ordering on $\wp$ be defined as follows for any pair ${}^1p, {}^2p \in \wp$ ${}^1p \leq {}^2p$ if and only if ${}^1p_j \leq {}^2p_j$ for all $j \in J$.

Given an arbitrary pair ${}^1p, {}^2p \in \wp$ such that ${}^1p \leq {}^2p$ let $[{}^1p, {}^2p] = \{p \in \wp \mid {}^1p \leq p < {}^2p\}$. For any pair ${}^1p, {}^2p \in \wp$ $(\{{}^1p, {}^2p\} \leq \}$ is a lattice.

Now we recall some of the properties of the solution set $S(Q, r)$. Employing the partial ordering on $\wp$, let an element $\hat{p}$ of $S(Q, r)$ be called a maximal solution of $p \circ Q = r$ if for all $p \in S(Q, r)$, $p \geq \hat{p}$ implies $p = \hat{p}$ if for all $p \in S(Q, r)$ $p < \tilde{p}$ then that is the maximum solution. Similar discussion can be made on the minimal solution of $p \circ Q = r$. The minimal solution is unique if $p \geq \hat{p}$ (i.e. $\hat{p}$ is unique).

It is well known when ever the solution set S (Q, r) is not empty it always contains a unique maximum solution $\hat{p}$ and it may contain several minimal solution.

Let $\breve{S}$ (Q, r) denote the set of all minimal solutions. It is known that the solution set S (Q, r) is fully characterized by the maximum and minimal solution in the sense that it consists exactly of the maximum solution $\hat{p}$ all the minimal solutions and all elements of $\wp$ that are between $\hat{p}$ and the numeral solution.

Thus S (Q, r) = $\bigcup_{p} [\tilde{p}, \hat{p}]$ where the union is taken for all $\tilde{p} \in \breve{S}$ (Q, r). When S (Q, r) ≠ ϕ, the maximum solution.

$\hat{p}$ = [ $\hat{p}_j$ | j ∈ J] of p ° Q = r is determined as follows:

$$\hat{p}_j = \min_{k \in K} \sigma\,(q_{ik}, r_k) \text{ where } \sigma\,(q_{jk}, r_k) = \begin{cases} r_k & \text{if } q_{jk} > r_k \\ 1 & \text{otherwise} \end{cases}$$

when $\hat{p}$ determined in this way does not satisfy p ° Q = r then S(Q, r) = ϕ. That is the existence of the maximum solution $\hat{p}$ as determined by $\hat{p}_j = \min_{k \in K} \sigma$ $(q_{jk}, r_k)$ is a necessary and sufficient condition for S (Q, r) ≠ ϕ. Once $\hat{p}$ is determined by $\hat{p}_j = \min_{k \in K} \sigma$ $(q_{jk}, r_k)$, we must check to see if it satisfies the given matrix equations p ° Q = r.

If it does not then the equation has no solution (S (Q, r) = ϕ), otherwise $\hat{p}$ in the maximum solution of the equation and we next determine the set $\tilde{S}$ (Q, r) of its minimal solutions.

### 2.6.3 Fuzzy compatibility relations and composition of fuzzy relations

In this subsection we recall the definition of fuzzy compatibility relations, fuzzy ordering relations, fuzzy morphisms, and sup and inf compositions of fuzzy relations.

**DEFINITION [106]:** *A binary relation R(X, X) that is reflexive and symmetric is usually called a compatibility relation or tolerance relation. When R(X, X) is a reflexive and symmetric fuzzy relation it is sometimes called proximity relation.*

An important concept associated with compatibility relations is compatibility classes. Given a crisp compatibility relation R(X, X), a compatibility class is a subset A of X such that (x, y) ∈ R for all x, y ∈ A.

A maximal compatibility class or maximal compatible is a compatibility class that is not properly contained with in any other compatibility class. The family consisting of all the maximal compatibles induced by R on X is called a complete cover of X with respect to R.

When R is a fuzzy compatibility relation, compatibility classes are defined in terms of a specified membership degree α. An α-compatibility class is a subset A of X such that R (x, y) ≥ α for all x, y ∈ A. Maximal α-compatibles and complete α-cover are obvious generalizations of the corresponding concepts for crisp compatibility relations.

Compatibility relations are often conveniently viewed as reflexive undirected graphs contrary to fuzzy cognitive maps that are directed graphs. In this context, reflexivity implies that each node of the graph has a loop connecting the node to itself the loops are usually omitted from the visual representations of the graph although they are assumed to be present. Connections between nodes as defined by the relation, are not directed, since the property of symmetry guarantees that all existing connections appear in both directions. Each connection is labeled with the value corresponding to the membership grade R (x, y) = R (y,x).

We illustrate this by the following example.

***Example 2.6.3.1:*** Consider a fuzzy relation R(X, X) defined on X = {$x_1$, $x_2$,…, $x_8$} by the following membership matrix:

$x_1$ $x_2$ $x_3$ $x_4$ $x_5$ $x_6$ $x_7$ $x_8$

$$\begin{array}{c} x_1 \\ x_2 \\ x_3 \\ x_4 \\ x_5 \\ x_6 \\ x_7 \\ x_8 \end{array}\begin{bmatrix} 1 & .3 & 0 & 0 & .4 & 0 & 0 & .6 \\ .3 & 1 & .5 & .3 & 0 & 0 & 0 & 0 \\ 0 & .5 & 1 & 0 & 0 & .7 & .6 & .8 \\ 0 & .3 & 0 & 1 & .2 & 0 & .7 & .5 \\ .4 & 0 & 0 & .2 & 1 & 0 & 0 & 0 \\ 0 & 0 & .7 & 0 & 0 & 1 & .2 & 0 \\ 0 & 0 & .6 & .7 & 0 & .2 & 1 & .8 \\ .6 & 0 & .8 & .5 & 0 & 0 & .8 & 1 \end{bmatrix}.$$

Since the matrix is symmetric and all entries on the main diagonal are equal to 1, the relation represented is reflexive, and symmetric therefore it is a compatibility relation. The graph of the relation is shown by the following figure 2.6.3.1, its complete $\alpha$-covers for $\alpha > 0$ and $\alpha \in \{0, .3, .1, .4, .6, .5, .2, .8, .7, 1\}$ is depicted.

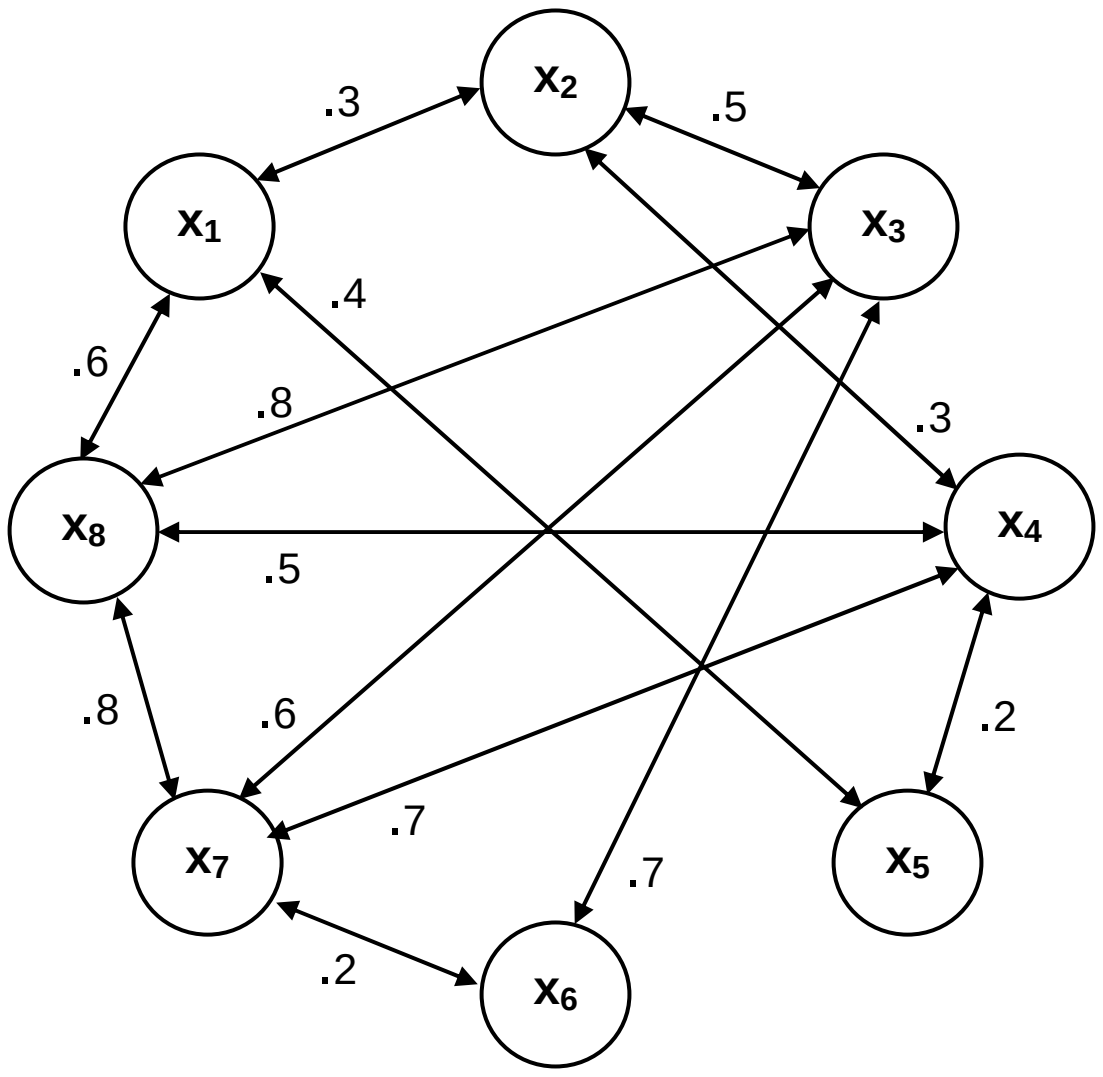

FIGURE: 2.6.3.1

Figure 2.6.3.1 is the graph of compatibility relation given in example 2.6.3.1 while similarity and compatibility relations are characterized by symmetry, ordering relations require asymmetry (or anti-symmetry) and transitivity. There are several types of ordering relations.

A crisp binary relation R(X, X) that is reflexive, anti symmetric and transitive is called a partial ordering. The common symbol ≤ is suggestive of the properties of this class of relations. Thus x ≤ y denotes (x, y) ∈ R and signifies that x precedes y. The inverse partial ordering $R^{-1}$ (X, X) is suggested by the symbol ≥.

If y ≥ x including that (y, x) ∈ $R^{-1}$ then we say that y succeeds x. When x ≤ y; x is also referred to as a predecessor of y while y is called a successor of x. When x ≤ y and there is no z such that x ≤ y and z ≤ y, x is called an immediate predecessor of y and y is called an immediate successor of x. If we need to distinguish several partial orderings, such as P, Q and R we use the symbol $\overset{P}{\leq}, \overset{Q}{\leq}$ and $\overset{R}{\leq}$ respectively.

Observe that a partial ordering '≤' on X does not guarantee that all pairs of elements x, y in X are comparable in the sense that either x ≤ y or y ≤ x. Thus, for some x, y ∈ X it is possible that x is neither a predecessor nor a successor of y. Such pairs are called non comparable with respect to ≤.

### 2.6.4 New FRE to estimate the peak hours of the day for transport system

In this section we just recall the notion of new fuzzy relational equations and study the estimation of the peak hour problem for transport systems using it we have also compared our results with the paper of W.B.Vasantha Kandasamy and V. Indra where FCMs are used [225]. We establish in our study which yields results, which are non over lapping and unique solutions.

***Example 2.6.4.1:*** Modified fuzzy relation equation has been derived for analyzing passenger preference for a particular hour in a day.

Since any transport or any private concern which plys the buses may not in general have only one peak hour in day, for; the peak hours are ones where there is the maximum number of passengers traveling in that hour. The passengers can be broadly classified as college students, school going children, office going people, vendors etc. Each category will choose a different hour according to their own convenience. For example the vendor group may go for buying good in the early morning hours and the school going children may prefer to travel from 7.00 a.m. to 8 a.m., college students may prefer to travel from 8 a.m. to 9 a.m. and the office going people may prefer to travel from 9.00 a.m. to 10.00 a.m. and the returning hours to correspond to the peak hours as the school going children may return home at about 3.00 p.m. to 4.00 p.m., college students may return home at about 2.00 p.m. to 3.30 p.m. and the office going people may return home at about 5.00 p.m. to 6.00 p.m. Thus the peak hours of a day cannot be achieved by solving the one equation $P \circ Q = R$. So we reformulate this fuzzy relation equation in what follows by partitioning $Q_i$'s. This in turn partition the number of preferences depending on the set Q which correspondingly partitions R also. Thus the fuzzy relation equation say $P \circ Q = R$ reduces to a set of fuzzy relations equations $P_1 \circ Q_1 = R_1$, $P_2 \circ Q_2 = R_2$, …, $P_s \circ Q_s = R_s$ where $Q = Q_1 \cup Q_2 \cup \ldots \cup Q_s$ such that $Q_i \cap Q_j = \phi$ for $i \neq j$. Hence by our method we get s preferences. This is important for we need at least 4 to 5 peak hours of a day. Here we give a new method by which we adopt the feed forward neural network to the transportation problem.

We briefly describe the modified or the new fuzzy relation equation used here. We know the fuzzy relation equation can be represented by neural network. We restrict our study to the form

$$P \circ Q = R \qquad (1)$$

where $\circ$ is the max-product composition; where $P = [p_{ij}]$, $Q = [q_{jk}]$ and $R = [r_{ik}]$, with $i \in N_n$, $j \in N_m$ and $k \in K_s$. We want to determine P. Equation (1) represents the set of equations.

$$\max_{j \in J_m} p_{ij}\, q_{jk} = \tau_{ik} \qquad (2)$$

for all $i \in N_n$ and $k \in N_s$.

To solve equation (2) for $P_{ij}$ ($i \in N_n$, $j \in N_m$), we use the feed forward neural network with m inputs and only one layer with n neurons.

First, the activation function employed by the neurons is not the sigmoid function, but the so-called linear activation function $f$ defined for all $a \in R$ by

$$f(a) = \begin{cases} 0 & \text{if } a < 0 \\ a & \text{if } a \in [0, 1] \\ 1 & \text{if } a > 1. \end{cases}$$

Second, the output $y_i$ of neuron i is defined by

$$y_i = f\left(\max_{j \in N} w_{ij} x_j\right) \quad \left(i \in N_n\right).$$

Given equation 1, the training set of columns $q_k$ of matrix Q as input ($x_j = q_{ik}$ for each $j \in N_m$, $k \in N_s$) and columns $r_k$ of matrix R as expected output ($y_i = r_{jk}$ for each $i \in N_n$ and $k \in N^s$). Applying this training set to the feed forward neutral network, we obtain a solution to equation 1, when the error function reaches zero. The solution is then expressed by the weight $w_{ij}$ as $p_{ij} = w_{ij}$ for all $i \in N_n$ and $j \in N_m$. Thus $p = (w_{ij})$ is a $n \times n$ matrix.

It is already known that the fuzzy relation equation is in the dominant stage and there is lot of scope in doing research in this area, further it is to be also tested in real data.

Here we are transforming the single equation $P \circ Q = R$ into a collection of equations. When the word preference is said, there should be many preferences. If only one choice is given or the equation results in one solution, it cannot be called as the preference. Further, when we do some experiment in the real data, we may have too many solutions of preferences. For a unique solution sought out cannot or may not be available in reality, so to satisfy all the conditions described above, we are forced to reformulate the equation $P \circ Q = R$. We partition the set Q into number of partition depending on the number preferences. When Q is partitioned, correspondingly R also gets partitioned, hence the one equation is transformed into the preferred number of equations.

Thus Q and R are given and P is to be determined. We partition Q into s sets, say $Q_1, Q_2, \ldots, Q_s$ such that $Q = Q_1 \cup Q_2 \cup \ldots \cup Q_s$, correspondingly R will be partitioned as $R = R_1 \cup R_2 \cup \ldots \cup R_s$. Now the resulting modified fuzzy equat6ions are $P_1 \circ Q_1 = R_1, P_2 \circ Q_2 = R_2, \ldots, P_s \circ Q_s = R_s$ respectively. Hence by our method, we obtain s preferences.

Since in reality it is difficult to make the error function $E_p$ to be exactly zero, we in our new fuzzy relation equation accept for the error function $E_p$ to be very close to zero. This is a deviation of the formula. Also we do not accept many stages in the arriving of the result. So once a proper guess is made even at the first stage we can get the desired solution by making $E_p$ very close to zero.

We are to find the passenger preference for a particular hour. The passenger preference problem for a particular hour reduces to finding the peak hours of the day (by peak hours of the day, we mean the number of passengers traveling in that hour is the maximum). Since the very term, preference by a passenger for a particular hour is an attribute, we felt it would be interesting if we adopt the modified fuzzy relation equation to this problem

So in our problem, we use the fuzzy relation equation $P \circ Q = R$, where P denotes the preference of a passenger to a particular hour, Q denotes the specific hour under study say $h_i$, $i = 1, 2, \ldots, 17$ where $h_i$ denotes the hour ending at 6 a.m., $h_2$ denotes the hour ending at 7 a.m., …, $h_{17}$ denotes the hour ending at 10 p.m. and $R_i$ denotes the number of passengers traveling during that particular hour $h_i$, for $i = 1, 2, \ldots, 17$.

Here we use the fuzzy relation equation to determine P. We formulate the problem as follows:

If $h_i$, for $i = 1, 2, \ldots, n$ are the n-hour endings, $R_i$, for $i = 1, 2, \ldots, n$ denote the number of passengers traveling during hour $h_i$, for $i = 1, 2, \ldots, n$. We denote by R the set $\{R_1, R_2, \ldots, R_n\}$ and $Q = \{h_1, h_2, \ldots, h_n\}$. To calculate the preference of a passenger to a particular hour we associative with each $R_i$, a weight $w_i$. Since $R_i$ correspond to the number of passenger traveling in that hour $h_i$, is a positive value and hence comparison between any two $R_i$ and $R_j$'s always exist. Therefore, if $R_i < R_j$, then we associate a weight $w_i$ to $R_i$ and $w_j$ to $R_j$ such that $w_i < w_j$, where $w_i$ and $w_j$ take values in the interval [0, 1].

Now we solve the matrix relation equation $P \circ Q = R$ and obtain the preference of the passenger P for a particular time period, which is nothing but the maximum number of passengers traveling in that hour.

If we wish to obtain only one peak hour for the day, we take all the n elements and form a matrix equation,

$$P_{max} \bullet \begin{bmatrix} h_1 \\ h_2 \\ \vdots \\ h_n \end{bmatrix}_{n\times 1} = \begin{bmatrix} R_1 \\ R_2 \\ \vdots \\ R_n \end{bmatrix}_{n\times 1}$$

and find the $n \times n$ matrix $P = (w_{ij})$ using the method described in the beginning. We choose in several steps the weight function $w_1, w_2, \ldots, w_n$ so that the error function $E_p$ reaches very near to zero. It is pertinent to mention here for our passengers preference problem we accept a value other than zero but which is very close to zero as limit which gives us the desired preference. If we wish to have two peaks hours, we partition Q into $Q_1$ and $Q_2$ so that correspondingly R gets partitioned in $R_1$ and $R_2$ and obtain the two peak hours using the two equations $P_1 \circ Q_1 = R_1$ and $P_2 \circ Q_2 = R_2$ respectively. Clearly $P_1 = \begin{pmatrix} 1 \\ w_{ij} \end{pmatrix}$ the weights associated with the set $R_1$ and $P_2 = \begin{pmatrix} 2 \\ w_{ij} \end{pmatrix}$ the weights associated with the set $R_2$.

If we wish to have a peak hours, $s < n$, then we partition $h_i$ for $i = 1, 2, \ldots, n$ into s disjoint sets and find the s peak hours of the day. This method of partitioning the fuzzy relation equation can be used to any real world data problem, though we have described in the context of the transportation problem.

We have tested our hypothesis in the real data got from Pallavan Transport Corporation.

Hour ending Q; 6, 7, 8, 9, 10, 11, 12, 13, 14, 15, 16, 17, 18, 19, 20, 21, 22.

Passengers per hour R: 96, 71, 222, 269, 300, 220, 241, 265, 249, 114, 381, 288, 356, 189, 376, 182, 67.

We have partitioned the 17 hours of the day Q.

i) by partitioning Q into three elements each so as to get five preferences,
ii) by partitioning Q into five elements each so as to get three preferences and
iii) by arbitrarily partitioning Q into for classes so as to get four preferences.

In all cases from these real data, our predicated value coincides with the real preference value.

Since all the concepts are to be realized as fuzzy concepts, we at the first state make the entries of Q and R to lie between 0 and 1. This is done by multiplying Q by $10^{-2}$ and R by $10^{-4}$ respectively.

We partition Q into three elements each by taking only the first 15 elements from the table. That is $Q = Q_1 \cup Q_2 \cup Q_3 \cup Q_4 \cup Q_5$ and the corresponding $R = R_1 \cup R_2 \cup R_3 \cup R_4 \cup R_5$.

| For | $Q_1 = x_i$ | $R_1 = r_{ik}$ |
|---|---|---|
| | 0.06 | 0.0096 |
| | 0.07 | 0.0071 . |
| | 0.08 | 0.0222 |

The fuzzy relation equation is

$$P1 = \begin{bmatrix} 0.06 \\ 0.07 \\ 0.08 \end{bmatrix} = \begin{bmatrix} 0.0096 \\ 0.0071 \\ 0.0222 \end{bmatrix}.$$

We employ the same method described earlier, where, the linear activation function *f* is defined by

$$f(a) = \left\{\begin{matrix} 0 & \text{if } a < 0 \\ a & \text{if } a \in [0, 1] \\ 1 & \text{if } a < 1 \end{matrix}\right\}$$

for all $a \in R$ and output $y_i$ of the neuron i is defined by

$$y_i = f\left(\max_{j \in N_m} w_{ij} \; x_j\right) \; (i \in N_n)$$

calculate $\max\limits_{j \in N_m} \; w_{ij} \quad x_j$ as follows:

(i) $w_{11}x_1 = 0.03 \times 0.06 = 0.0018$
$w_{12}x_2 = 0.0221875 \times 0.07 = 0.001553125$
$w_{13}x_3 = 0.069375 \times 0.08 = 0.00555$
$\therefore$ Max (0.0018, 0.001553125, 0.00555) = 0.00555

$f\left(\max\limits_{j \in N_m} w_{ij} x_j\right) = f(0.00555) = 0.00555$ (Since $0.00555 \in [0, 1]$)

$\therefore \; y_1 = 0.00555$.

(ii) $w_{21}x_1 = 0.06 \times 0.06 = 0.0036$
$w_{22}x_2 = 0.044375 \times 0.07 = 0.00310625$
$w_{23}x_3 = 0.13875 \times 0.08 = 0.0111$
$\therefore$ Max (0.0036, 0.00310625, 0.0111) = 0.0111

$f\left(\max\limits_{j \in N_m} w_{ij} x_j\right) = f(0.0111) = 0.0111$ (Since $0.0111 \in [0, 1]$)

$\therefore \; y_2 = 0.0111$.

(iii) $w_{31}x_1$ = $0.12 \times 0.06$ = 0.0072

$w_{32}x_2$ = $0.08875 \times 0.07$ = 0.0062125

$w_{33}x_3$ = $0.2775 \times 0.08$ = 0.0222

$\therefore$ Max (0.0072, 0.0062125, 0.0222) = 0.0222 (Since 0.0222 $\in$ [0, 1]

$\therefore$ $y_3$ = 0.0222.

$$f\left(\max_{j \in N_m} w_{ij} x_j\right) = f(0.0222) = 0.0222 \text{ (Since } 0.0222 \in [0, 1])$$

$\therefore$ $y_3$ = 0.0222.

Feed Forward Neural Network representing the solution is shown above.

$$\therefore P_1 = \begin{bmatrix} 0.03 & 0.06 & 0.12 \\ 0.0221875 & 0.044375 & 0.08875 \\ 0.069375 & 0.13875 & 0.2775 \end{bmatrix}.$$

Verification:

Consider, P $\circ$ Q = R

that is $\max_{j \in N_m} p_{ij}\, q_{jk} = r_{ik}$

$\therefore$ Max (0.0018, 0.0042, 0.0096) = 0.0096

Max (0.00133125, 0.00310625, 0.0071) = 0.0071

Max (0.0041625, 0.0097125, 0.0222) = 0.0222.

Similarly by adopting the above process, we have calculated the passenger preferences $P_2$, $P_3$, $P_4$ and $P_5$ for the pairs ($Q_2$, $R_2$), ($Q_3$, $R_3$), ($Q_4$, $R_4$) and ($Q_5$, $R_5$).

For

| $Q_2$ | $R_2$ |
|---|---|
| 0.09 | 0.0269 |
| 0.10 | 0.0300 |
| 0.11 | 0.0220 |

we have $P_2 = \begin{bmatrix} 0.1345 & 0.269 & 0.06725 \\ 0.15 & 0.3 & 0.075 \\ 0.11 & 0.22 & 0.00605 \end{bmatrix}$.

For

| $Q_3$ | $R_3$ |
|---|---|
| 0.12 | 0.0241 |
| 0.13 | 0.0265 |
| 0.14 | 0.0249 |

we have $P_3 = \begin{bmatrix} 0.2008 & 0.1004 & 0.0502 \\ 0.2208 & 0.1104 & 0.0552 \\ 0.2075 & 0.10375 & 0.051875 \end{bmatrix}$.

For

| $Q_4$ | $R_4$ |
|---|---|
| 0.15 | 0.0114 |
| 0.16 | 0.0381 |
| 0.17 | 0.0288 |

, $P_4 = \begin{bmatrix} 0.035625 & 0.07125 & 0.0178125 \\ 0.1190625 & 0.23125 & 0.05953125 \\ 0.09 & 0.18 & 0.045 \end{bmatrix}$,

and for

| $Q_5$ | $R_5$ |
|---|---|
| 0.18 | 0.0356 |
| 0.19 | 0.0189 |
| 0.20 | 0.0376 |

we have $P_5 = \begin{bmatrix} 0.0445 & 0.089 & 0.178 \\ 0.023625 & 0.04725 & 0.0945 \\ 0.047 & 0.094 & 0.188 \end{bmatrix}$.

On observing from the table, we see the preference $P_1$, $P_2$, $P_3$, $P_4$, and $P_5$ correspond to the peak hours of the day, $h_3$ that is 8 a.m. with 222 passengers, $h_5$ that is 10 a.m. with 300 passengers, $h_8$ that is 1 p.m. with 265 passengers, $h_{11}$ that is 4 p.m. with 381 passengers and $h_{15}$ that is 8 p.m. with 376 passengers. Thus this partition gives us five preferences with coincides with the real data as proved by the working.

max (0.003108, 0.00444, 0.00666, 0.0222, 0.01221) = 0.0222
max (0.003766, 0.00538, 0.0080694, 0.0269, 0.014795) = 0.0269
max (0.0042, 0.006, 0.009, 0.03, 0.0165) = 0.03
max (0.00308, 0.0044, 0.0065997, 0.0222, 0.0121) = 0.0222

Similarly we obtain the passenger preference P for the other entries using the above method.

For $Q_2$ $R_2$

| | |
|---|---|
| 0.12 | 0.0241 |
| 0.13 | 0.0265 |
| 0.14 | 0.0249 |
| 0.15 | 0.0114 |
| 0.16 | 0.0381 |

we have

$$P_2 = \begin{bmatrix} 0.030125 & 0.03765625 & 0.05020833 & 0.0753125 & 0.150625 \\ 0.033125 & 0.04140625 & 0.0552083 & 0.0828125 & 0.165625 \\ 0.031125 & 0.03890625 & 0.051875 & 0.0778125 & 0.155625 \\ 0.01425 & 0.0178125 & 0.02375 & 0.035625 & 0.07125 \\ 0.047625 & 0.05953 & 0.079375 & 0.1190625 & 0.238125 \end{bmatrix}$$

and for

| $Q_3$ | $R_3$ |
|---|---|
| 0.17 | 0.0288 |
| 0.18 | 0.0356 |
| 0.19 | 0.0189 |
| 0.20 | 0.0376 |
| 0.21 | 0.0182 we have |

$w_{54}\, x_4 \;=\; 0.15 \times 0.10 \;=\; 0.015$

$w_{55}\, x_5 \;=\; 0.11 \times 0.11 \;=\; 0.0121.$

$\therefore$ Max = (0.002485, 0.00888, 0.012105, 0.015, 0.0121) = 0.015.

$$f\left(\max_{j \in N_m} w_{ij}\; x_j\right) = f(0.015) = 0.015 \text{ (Since } 0.015 \in [0.1])$$

$\therefore y_5 \;=\; 0.015.$

Feed forward neural network representing the solution is shown above

$$\therefore P_1 = \begin{bmatrix} 0.0142 & 0.01775 & 0.02366 & 0.071 & 0.03555 \\ 0.0444 & 0.0555 & 0.074 & 0.222 & 0.111 \\ 0.0538 & 0.06725 & 0.08966 & 0.269 & 0.1345 \\ 0.06 & 0.075 & 0.1 & 0.3 & 0.15 \\ 0.044 & 0.055 & 0.07333 & 0.220 & 0.11 \end{bmatrix}.$$

Verification :

Consider, $P \circ Q = R$

that is $\max\limits_{j \in N_m} p_{ij}\, q_{jk} = r_{ik}$

max (0.000994, 0.00142, 0.0021294, 0.0071, 0.003905)
= 0.0071

$w_{24}\, x_4 = 0.075 \times 0.10 = 0.0075$

$w_{25}\, x_5 = 0.055 \times 0.11 = 0.00605$

$\therefore$ Max (0.0012425, 0.00444, 0.0060525, 0.0075, 0.00605)
= 0.0075

$f\left(\max\limits_{j \in N_m} w_{ij}\;\; x_j\right) = f\,(0.0075) = 0.0075$ (Since $0.0075 \in [0, 1]$)

$\therefore y_2 = 0.0075$

(iii)

$w_{31}\, x_1 = 0.02366 \times 0.07 = 0.016566$

$w_{32}\, x_2 = 0.074 \times 0.08 = 0.00592$

$w_{33}\, x_3 = 0.08966 \times 0.09 = 0.00807$

$w_{34}\, x_4 = 0.1 \times 0.10 = 0.010$

$w_{35}\, x_5 = 0.07333 \times 0.11 = 0.008066$

$\therefore$ Max (0.016566, 0.00592, 0.00807, 0.010, 0.008066)
= 0.016566

$f\left(\max\limits_{j \in N_m} w_{ij}\;\; x_j\right) = f\,(0.016566) = 0.016566$

(Since $0.016566 \in [0, 1]$)

$\therefore y_3 = 0.016566$

(iv)

$w_{41}\, x_1 = 0.071 \times 0.07 = 0.00497$

$w_{42}\, x_2 = 0.222 \times 0.08 = 0.01776$

$w_{43}\, x_3 = 0.269 \times 0.09 = 0.02421$

$w_{44}\, x_4 = 0.3 \times 0.10 = 0.03$

$w_{45}\, x_5 = 0.220 \times 0.11 = 0.0242$

$\therefore$ Max (0.00497, 0.01776, 0.02421, 0.03, 0.0242) = 0.03

$$f\left(\max_{j \in N_m} w_{ij}\ x_j\right) = f\ (0.03) = 0.03 \text{ (Since } 0.03 \in [0, 1])$$

$\therefore y_4 = 0.03$

(v)

| | | | | |
|---|---|---|---|---|
| $w_{51}\ x_1$ | = | $0.0355 \times 0.07$ | = | 0.002485 |
| $w_{52}\ x_2$ | = | $0.111 \times 0.08$ | = | 0.00888 |
| $w_{53}\ x_3$ | = | $0.1345 \times 0.09$ | = | 0.012105. |

Now, we partition Q into five elements each by leaving out the first and the last element from the table as $Q_1$ $Q_2$ and $Q_3$ and calculate $P_1$, $P_2$ and $P_3$ as in the earlier case:

| for | $Q_1$ | $R_1$ |
|---|---|---|
| | 0.07 | 0.0071 |
| | 0.08 | 0.0222 |
| | 0.09 | 0.0269 |
| | 0.10 | 0.0300 |
| | 0.11 | 0.0220. |

The fuzzy relation equation is

$$P_1 \circ \begin{bmatrix} 0.06 \\ 0.08 \\ 0.09 \\ 0.10 \\ 0.11 \end{bmatrix} = \begin{bmatrix} 0.0071 \\ 0.0222 \\ 0.0269 \\ 0.0300 \\ 0.0220 \end{bmatrix}.$$

Calculate $\max\limits_{j \in N_m} w_{ij}\ x_j$ as follows

(i)

| | | | | |
|---|---|---|---|---|
| $w_{11}\ x_1$ | = | $0.0142 \times 0.07$ | = | 0.000994 |
| $w_{12}\ x_2$ | = | $0.0444 \times 0.08$ | = | 0.003552 |
| $w_{13}\ x_3$ | = | $0.0538 \times 0.09$ | = | 0.004842 |
| $w_{14}\ x_4$ | = | $0.06 \times 0.10$ | = | 0.006 |
| $w_{15}\ x_5$ | = | $0.044 \times 0.11$ | = | 0.00484 |

$\therefore$ Max (0.000994, 0.003552, 0.004842, 0.006, 0.00484) = 0.006

$$f\left(\max_{j\in N_m} w_{ij}\ x_j\right) = f(0.006) = 0.006 \text{ (Since } 0.006 \in [0, 1])$$

$\therefore y_1 = 0.006$

(ii)

$w_{21}\ x_1 = 0.01775 \times 0.07 = 0.0012425$

$w_{22}\ x_2 = 0.0555 \times 0.08 = 0.00444$

$w_{23}\ x_3 = 0.06725 \times 0.09 = 0.0060525$

$$P_3 = \begin{bmatrix} 0.0288 & 0.036 & 0.048 & 0.149 & 0.072 \\ 0.0356 & 0.0445 & 0.05933 & 0.178 & 0.089 \\ 0.0189 & 0.023625 & 0.0315 & 0.0945 & 0.04725 \\ 0.0376 & 0.047 & 0.06266 & 0.0188 & 0.094 \\ 0.0182 & 0.02275 & 0.03033 & 0.091 & 0.0455 \end{bmatrix}.$$

On observing from the table, we see the preference $P_1$, $P_2$ and $P_3$ correspond to the peak hours of the day, $h_5$ that is 10 a.m. with 300 passengers, $h_{11}$ that is 4 p.m. with 381 passengers and $h_{15}$ that is 8 p.m. with 376 number of passengers. Thus this partition gives us three preferences, which coincides with the real data as proved by the working.

We now partition Q arbitrarily, that is the number of elements in each partition is not the same and by a adopting the above method we obtain the following results:

For

| $Q_1$ | $R_1$ |
|---|---|
| 0.06 | 0.0096 |
| 0.07 | 0.0071 |
| 0.08 | 0.0222 |

we have $P_1 = \begin{bmatrix} 0.03 & 0.06 & 0.12 \\ 0.0221875 & 0.044375 & 0.08875 \\ 0.069375 & 0.13875 & 0.2775 \end{bmatrix}.$

For $Q_2$ $R_2$

| | |
|---|---|
| 0.09 | 0.269 |
| 0.10 | 0.300 |
| 0.11 | 0.220 |
| 0.12 | 0.241 |
| 0.13 | 0.265 |

we have

$$P_2 = \begin{bmatrix} 0.1345 & 0.17933 & 0.1076 & 0.08966 & -.06725 \\ 0.15 & 0.2 & 0.12 & 0.1 & 0.075 \\ 0.11 & 0.14666 & 0.088 & 0.0733 & 0.055 \\ 0.1205 & 0.16066 & 0.0964 & 0.08033 & 0.06025 \\ 0.1325 & 0.17666 & 0.106 & 0.08833 & 0.06625 \end{bmatrix}.$$

For

| $Q_3$ | $R_3$ |
|---|---|
| 0.14 | 0.0249 |
| 0.15 | 0.0114 |

we have $P_3 = \begin{bmatrix} 0.083 & 0.166 \\ 0.038 & 0.076 \end{bmatrix}$

and for

| $Q_4$ | $R_4$ |
|---|---|
| 0.16 | 0.0381 |
| 0.17 | 0.0288 |
| 0.18 | 0.0356 |
| 0.19 | 0.0189 |
| 0.20 | 0.0376 |
| 0.21 | 0.0182 |

we have $P_4 =$

$$\begin{bmatrix} 0.09525 & 0.0635 & 0.047625 & 0.0381 & 0.1905 & 0.03175 \\ 0.072 & 0.048 & 0.036 & 0.0288 & 0.144 & 0.024 \\ 0.08233 & 0.05322 & 0.04111 & 0.03193 & 0.178 & 0.02411 \\ 0.4725 & 0.0312 & 0.023625 & 0.0189 & 0.0945 & 0.001575 \\ 0.0762 & 0.0508 & 0.0381 & 0.03048 & 0.188 & 0.0254 \\ 0.04275 & 0.0285 & 0.0213756 & 0.0171 & 0.091 & 0.001425 \end{bmatrix}$$

.

We obtain in the preferences $P_1$, $P_2$, $P_3$ and $P_4$ by partitioning the given data into a set of three elements, a set of five elements, a set of two elements and a set of six elements. On observing from the table, we see that these preferences correspond to the peak hours of the day, $h_3$ that is 8 a.m. with 222 passengers, $h_5$ that is 10 a.m. with 300 passengers, $h_9$ that is 2 p.m. with 249 number of passengers and $h_{11}$ that is 4 p.m. with 381 number of passengers. Thus this partition gives us four preferences which coincides with the real data as proved by the working.

Thus the Government sector can run more buses at the peak hours given and also at the same time restrain the number of buses in the non peak hours we derived the following conclusions:

1. The fuzzy relation equation described given by 1 can give only one preference function $P \circ Q = R$ but the partition method described by us in this paper can give many number of preferences or desired number of preferences.
2. Since lot of research is needed we feel some other modified techniques can be adopted in the FRE $P \circ Q = R$.
3. We have tested our method described in the real data taken from Pallavan Transport Corporation and our results coincides with the given data.
4. We see the number of preference is equal to the number of the partition of Q.
5. Instead of partitioning Q, if we arbitrarily take overlapping subsets of Q certainly we may get the same preference for two or more arbitrary sets.

We see that our method of the fuzzy relation equation can be applied to the peak hour problem in a very successful way. Thus only partitioning of Q can yield non-overlapping unique solution.

Finally, in our method we do not force the error function $E_p$ to become zero, by using many stages or intermittent steps. We accept a value very close to zero for $E_p$ as a preference solution.

For more please refer [230].

## 2.6.5 The effect of globalization on silk weavers who are bonded labourers using FRE

The strategies of globalization and the subsequent restructuring of economies, including the increased mechanization of labor has had stifling effects on the lives of the silk weavers in the famous Kancheepuram District in Tamil Nadu, India. Here, we study the effects of globalization, privatization and the mechanization of labor, and how this has directly affected (and ruined) the lives of thousands of silk weavers, who belong to a particular community whose tradition occupation is weaving. This research work is based on surveys carried out in the Ayyampettai village near Kancheepuram. The population of this village is around 200 families, and almost all of them are involved in the weaving of silk saris. They are skilled weavers who don't have knowledge of any other trade. Most of them are bond to labor without wages, predominantly because they had reeled into debt for sums ranging from Rs. 1000 upwards. They barely manage to have a square meal a day, and their work patterns are strenuous - they work from 6 a.m. to 7 p.m. on all days, expect the new moon day when they are forbidden from weaving.

Interestingly, their children are not sent to school, they are forced into joining the parental occupation, or into taking petty jobs in order to secure the livelihood. The villagers point to the advent of electric looms and reckon that their lives were much more bearable before this mechanization, at least they used to get better incomes. The wide scale introduction to electric looms / textile machines / power looms, has taken away a lot of their job opportunities. For instance, the machine can weave three silk saris which manually takes fifteen days to weave in just three hours. Also, machine woven silk saris are preferred to hand woven silk saris as a result of which their life is further shattered. Interviews with the weavers revealed the careless and negligent approach of the government to their problem.

***Example 2.6.5.1:*** Here, we study their problem and the effect of globalization on their lives using Fuzzy Relational Equations.

We have arrived at interesting conclusions to understand and assay this grave problem. We have made a sample survey of around 50 families out of the 200 families; who are bonded labourers living in Ayyampettai near Kancheepuram District in Tamil Nadu; have become bonded for Rs.1000 to Rs.2000. They all belong to Hindu community viz. weavers or they traditionally call themselves as Mudaliar caste. Most of the owners are also Mudaliars. They were interviewed using a linguistic questionnaire. Some of the notable facts about their lives are as follows:

1. They do not know any other trade or work but most of them like to learn some other work.
2. They are living now below the poverty line because of the advent of electrical or power looms which has drastically affected their income.
3. The whole family works for over 10 hours with only one day i.e. new moon day in a month being a holiday. On new moon day they don't weave and they are paid by their owners on that day.
4. Only one had completed his school finals. All others have no education for they have to learn the trade while very young.
5. They don't have even a square meal a day.
6. Becoming member of Government Society cannot be even dreamt for they have to pay Rs.3000/- to Rs.5000 to Government and 3 persons should give them surety. So out of the 200 families there was only one was a Government Society member. After the globalization government do not give them any work because marketers prefer machine woven saris to hand woven ones.
7. Owners of the bonded labourers are not able to give work to these labourers.
8. Observations shows that female infanticide must be prevalent in these families as over 80% of the children are only males.
9. The maximum salary a family of 3 to 4 members is around Rs. 2000/- 5% of them alone get this 90% of the

families get below Rs.2000 p.m.

10. Paying as rent, electricity, water, etc makes them live below poverty line.

The following attributes are taken as the main point for study:

$B_1$ – No knowledge of any other work has made them not only bonded but live in penury.

$B_2$ – Advent of power looms and globalization (modern textile machinery) has made them still poorer.

$B_3$ – Salary they earn in a month.

$B_4$ – No savings so they become more and more bonded by borrowing from the owners, they live in debts.

$B_5$ – Government interferes and frees them they don't have any work and Government does not give them any alternative job.

$B_6$ – Hours / days of work.

We have taken these six heads $B_1$, $B_2$, … , $B_6$ related to the bonded labourers as the rows of the fuzzy relational matrix.
The main attributes / heads $O_1$, $O_2$, $O_3$, $O_4$ related to the owners of the bonded labourers are :

$O_1$ – Globalization / introduction of modern textile machines

$O_2$ – Profit or no loss

$O_3$ – Availability of raw goods

$O_4$ – Demand for finished goods

Using these heads related to owners along columns the fuzzy relational equations are formed using the experts opinions.

The following are the limit sets using the questionnaire :

$B_1 \geq 0.5$ Means no knowledge of other work hence live in poverty.

$B_2 \geq 0.5$ Power loom / other modern textile machinery had made their condition from bad to worse.

| | |
|---|---|
| $B_3 \geq 0.5$ | Earning is mediocre. ($B_3 < 0.5$ implies the earning does not help them to meet both ends). |
| $B_4 \geq 0.4$ | No saving no debt. ($B_4 < 0.4$ implies they are in debt). |
| $B_5 \geq 0.5$ | Government interference has not helped. ($B_5 < 0.5$ implies Government Interference have helped). |
| $B_6 \geq 0.4$ | 10 hours of work with no holidays. ($B_6 < 0.4$ implies less than 10 hours of work). |
| | |
| $O_1 \geq 0.5$ | The globalizations / government has affected the owners of the bonded labourers drastically ($O_1 < 0.5$ implies has no impact on owners). |
| $O_2 \geq 0.5$ | Profit or no loss ($O_2 < 0.5$ implies total loss). |
| $O_3 \geq 0.6$ | Availability of raw materials. ($O_3 < 0.6$ implies shortage of raw material). |
| $O_4 \geq 0.5$ | Just they can meet both ends i.e. demand for finished goods and produced goods balance. ($O_4 < 0.5$ implies no demand for the finished product i.e. demand and supply do not balance). |

The opinion of the first expert who happens to be a bonded labor for the two generations aged in seventies is given vital importance and his opinion is transformed into the Fuzzy Relational Equation

$$P = \begin{array}{c} \\ B_1 \\ B_2 \\ B_3 \\ B_4 \\ B_5 \\ B_6 \end{array} \begin{array}{c} \begin{array}{cccc} O_1 & O_2 & O_3 & O_4 \end{array} \\ \begin{bmatrix} .8 & 0 & 0 & 0 \\ .8 & .3 & .3 & 0 \\ .1 & .2 & .3 & .4 \\ 0 & .1 & .1 & .1 \\ .8 & .1 & .2 & .4 \\ .2 & .4 & .4 & .9 \end{bmatrix} \end{array}.$$

By considering the profit suppose the owner gives values for Q where $Q^T = [.6, .5, .7, .5]$. Now P and Q are known in the fuzzy relational equation P o Q = R .

Using the max-min principle in the equation P o Q = R.

We get $R^T$ = [.6, .6, .4, .1, .6, .5] In the fuzzy relational equation P o Q = R, P corresponds to the weightages of the expert, Q is the profit the owner expects and R is the calculated or the resultant giving the status of the bonded labourers. When we assume the owners are badly affected by globalizations, but wants to carry out his business with no profit or loss, with moderate or good availability of the raw material and they have enough demand or demand and supply balance we obtain the following attributes related with the bonded labourers. The bonded labourers live in acute poverty as they have no other knowledge of any other work. The power loom has made their life from bad to worst, but the earning is medium with no savings and debts. They do not receive any help from the government, but they have to labor more than ten hours which is given by $[.6, .6, .4, .1, .6, .5]^T$. Using the same matrix P and taking the expected views of the bonded labourers R to be as $[.6, .4, .5, .4, .2, .6]^T$ .Using the equation $P^T$ o R = Q. We obtain Q = $[.6, .4, .4, .6]^T$.

The value of Q states the owners are affected by globalization. They have no profit but loss. They do not get enough raw materials to give to the bonded labor as the market prefers machine woven saris to hand made ones so the demand for the finished goods declines. Thus according to this expert, the main reason for their poverty is due to globalization i.e. the advent of power looms has not only affected them drastically as they do not have the knowledge of any other trade but it has also affected the lives of their owners.

A small owner who owns around ten bonded labor families opinion is taken as the second experts opinion. The weighted matrix P as given by the second expert is:

$$P = \begin{bmatrix} .7 & .1 & 0 & 0 \\ .9 & .2 & .3 & 0 \\ .0 & .1 & .2 & .3 \\ 0 & 0 & .1 & .1 \\ .9 & 0 & .1 & .4 \\ .1 & .2 & .4 & .7 \end{bmatrix}.$$

By considering the profit the owner expects i.e. taking Q = [.6, .5, .7, .5]$^T$

We calculate R using $P \text{ o } Q^T = R$

i.e. R = [.6, .6, .3, .1, .6, .5]$^T$.

We obtain the following attributes from R related with the bonded labourers.

They live in below poverty, as they have no other trade but the earning is medium with no savings and new debts. They do not get any help from the government, but they have to work more than 10 hours a day which is given by [.6, .6, .3, .1, .6, .5]$^T$ .

Using the same P i.e. the weightages we now find Q giving some satisfactory norms for the bonded labourers.

By taking R = [ .6, .4, .5, .4, .2, .6]$^T$ and using the equation $P^T \text{ o } R = Q$,

i.e. Q = [.6, .2, .4, .6]$^T$

which states the owners are badly affected by globalization. They have no profit but loss they do not get enough raw materials and the demand for the finished goods declines.

The third expert is a very poor bonded labor. The fuzzy relational matrix P given by him is

$$P = \begin{bmatrix} .9 & 0 & 0 & 0 \\ .5 & .3 & .4 & .1 \\ .2 & .2 & .2. & 3 \\ 0 & 0 & .1 & .2 \\ .7 & .2 & .2 & .4 \\ .2 & .3 & .3 & .8 \end{bmatrix}.$$

By considering the profit the owner expects i.e. taking Q = [.6, .5, .7, .5]$^T$ and using the relational equation P o Q = R, we calculate R;

R = [.6, .5, .3, .2, .6, .5]$^T$ .

We obtain the following attributes from R related with the bonded labourers.

This reveals that the bonded labourers standard of living is in a very pathetic condition. They do not have any other source of income or job. Their earning is bare minimum with no savings. Neither the government comes forward to help them nor redeem them from their bondage. In their work spot, they have to slog for 10 hours per day.

Using the same P i.e. the weightages we now find Q by giving some satisfactory norms for the bonded labourers.

By taking $R = [\ .6, .4, .5, .4, .2, .6]^T$ and using the equation $P^T \text{ o } R = Q$,

$$\text{i.e. } Q = [.6, .3, .4, .6]^T.$$

The value of Q states due to the impact of globalization (modern textile machinery), the owners are badly affected. They are not able to purchase enough raw materials and thus the out put from the industry declines. The owners do not get any profit but eventually end up in a great loss.

The following conclusions are not only derived from the three experts described here but all the fifty bonded labourers opinions are used and some of the owners whom we have interviewed are also ingrained in this analysis.

1. Bonded labourers are doubly affected people for the advent of globalization (modern textile machinery) has denied them small or paltry amount, which they are earning in peace as none of them have knowledge of any other trade.

2. The government has not taken any steps to give or train them on any trade or work or to be more precise they are least affected about their living conditions of them. Some of them expressed that government is functioning to protect the rich and see the rich do not loose anything but they do not even have any foresight about the bonded labourers or their petty owners by which they are making the poor more poorer.

3. Bonded labourers felt government has taken no steps to eradicate the unimaginable barrier to become members of the government society. They have to pay Rs.3000/- and also they should spell out and get the surety of 3 persons and the three persons demand more than Rs.3000/- each so only they are more comfortable in the hands of their owners i.e. as bonded labourers were at least they exist with some food, though not a square meal a day.

4. It is high time government takes steps to revive the life of the weavers who work as bonded labourers by training and giving them some job opportunities.

5. They felt government was killing the very talent of trained weavers by modernization as they have no knowledge of any other trade.

6. Child labor is at its best in these places as they cannot weave without the help of children. Also none of the children go to school and they strongly practice female infanticide.

7. Government before introducing these modern textile machineries should have analyzed the problem and should have taken steps to rehabilitate these weavers. Government has implemented textile machineries without foresight.

# INDEX

## C

## D

## E

## F

## G

## H

## I

## L

## M

## N

## O

## P

## R

## S

## T

## U

## V

## W

## Z

# ABOUT THE AUTHORS

**Dr.W.B.Vasantha Kandasamy** is an Associate Professor in the Department of Mathematics, Indian Institute of Technology Madras, Chennai. In the past decade she has guided 11 Ph.D. scholars in the different fields of non-associative algebras, algebraic coding theory, transportation theory, fuzzy groups, and applications of fuzzy theory of the problems faced in chemical industries and cement industries. Currently, four Ph.D. scholars are working under her guidance.

She has to her credit 636 research papers. She has guided over 51 M.Sc. and M.Tech. projects. She has worked in collaboration projects with the Indian Space Research Organization and with the Tamil Nadu State AIDS Control Society. This is her 30th book.

On India's 60th Independence Day, Dr.Vasantha was conferred the Kalpana Chawla Award for Courage and Daring Enterprise by the State Government of Tamil Nadu in recognition of her sustained fight for social justice in the Indian Institute of Technology (IIT) Madras and for her contribution to mathematics. (The award, instituted in the memory of Indian-American astronaut Kalpana Chawla who died aboard Space Shuttle Columbia). The award carried a cash prize of five lakh rupees (the highest prize-money for any Indian award) and a gold medal.
She can be contacted at vasanthakandasamy@gmail.com
You can visit her on the web at: http://mat.iitm.ac.in/~wbv or: http://www.vasantha.net

---

**Dr. Florentin Smarandache** is an Associate Professor of Mathematics at the University of New Mexico in USA. He published over 75 books and 100 articles and notes in mathematics, physics, philosophy, psychology, literature, rebus. In mathematics his research is in number theory, non-Euclidean geometry, synthetic geometry, algebraic structures, statistics, neutrosophic logic and set (generalizations of fuzzy logic and set respectively), neutrosophic probability (generalization of classical and imprecise probability). Also, small contributions to nuclear and particle physics, information fusion, neutrosophy (a generalization of dialectics), law of sensations and stimuli, etc.
He can be contacted at smarand@unm.edu

---

**K. Ilanthenral** is the editor of The Maths Tiger, Quarterly Journal of Maths. She can be contacted at ilanthenral@gmail.com